\documentclass[a4paper,11pt]{article}
\usepackage{color,xcolor,ucs}
\usepackage[top=0.1in, bottom=0.43in, left = 0.52in, right = 0.52in]{geometry}
\usepackage[linkcolor=black,colorlinks=true,urlcolor=blue]{hyperref}
\usepackage{mathtools}

\usepackage{color,xcolor,ucs}
\usepackage{mathtools}   \usepackage{tikz} 
\usepackage{ amssymb }
\usepackage{extarrows} 
\usepackage{pgf,tikz}
\usepackage{float}
\usetikzlibrary{positioning}
\usetikzlibrary{shapes.geometric}
\usetikzlibrary{shapes.misc}
\usetikzlibrary{arrows}
\usepackage{caption}
\usepackage{mathrsfs}
\usetikzlibrary{arrows,shapes,automata,backgrounds,petri,positioning}
\usetikzlibrary{decorations.pathmorphing}
\usetikzlibrary{decorations.shapes}
\usetikzlibrary{decorations.text}
\usetikzlibrary{decorations.fractals}
\usetikzlibrary{decorations.footprints}
\usetikzlibrary{shadows}
\usetikzlibrary{calc}
\usetikzlibrary{spy}
\usepackage{amsmath}
\usepackage{array}
\usepackage{ amssymb }
\usepackage{braket}
\usepackage{qcircuit}
\usepackage{soul}
\usepackage{braket}

\begin{document}
\title{{\huge From logarithmic delocalization of the six-vertex height function under sloped boundary conditions to weakened crossing probability estimates for the Ashkin-Teller, generalized random-cluster, and $(q_{\sigma},q_{\tau})$-cubic models}}
\author{{\large Pete Rigas}}
\date{}
\maketitle

\begin{abstract} 
    To obtain Russo-Seymour-Welsh estimates for the height function of the six-vertex model under sloped boundary conditions, which can be leveraged to demonstrate that the height function logarithmically delocalizes under a broader class of boundary conditions, we formulate crossing probability estimates in strips of the square lattice and in the cylinder, for parameters satisfying $a\equiv b$, $c \in [1,2]$, and $\mathrm{max} \{ a , b \} \leq c$, in which each of the first two conditions respectively relate to invariance under vertical and diagonal reflections enforced through the symmetry $\sigma \xi \geq -\xi$ for domains in strips of the square lattice, and satisfaction of the Fortuin-Kasteleyn-Ginibre property, for the height function and for its absolute value. To determine whether arguments for estimating crossing probabilities of the height function for flat boundary conditions from a recent work due to Duminil-Copin, Karila, Manolescu, and Oulamara remain applicable for boundary conditions that are not flat, from the set of all possible slopes of the height function given by the interior of the set of rational points from $[-1,1] \times [-1,1]$, we sample sloped Gibbs states in the square lattice, which do not have infinitely many disjointly oriented circuits. In comparison to Russo-Seymour-Welsh arguments for flat boundary conditions, Russo-Seymour-Welsh arguments for sloped boundary conditions present additional complications for both planar and cylindrical settings of the argument, in which crossing events that are considered in the strip, and then extended to the annulus and cylinder, must be maintained across rectangles of large aspect ratio, in spite of the fact that some positive proportion of faces within the strip freeze with positive probability. Once additional considerations have been raised for estimating crossing probabilities in the strip, which primarily relate to estimates of connectivity, and disconnectivity, events between frozen regions of the square lattice which are introduced as \textit{freezing clusters}, and in the cylinder, which primarily relate to the fraction of points about which cylindrical domains are situated that are deemed to satisfy the \textit{good} property after reversing the orientation of some loops to make the configuration balanced, we investigate additional models of interest for which weakened crossing probability estimates can be obtained, by making use of connections between the spin-representation of the six-vertex model and the Ashkin-Teller model, discussed in a $2019$ work due to Glazman and Peled, in addition to the random-cluster representation of the Ashkin-Teller model discussed in a $1997$ work due to Pfitzer and Velenik which shares connections with the anisotropic $N$-vector models. Besides being able to potentially  characterize various objects obtained throughout the argument in the future, some of which combinatorially relate to the shape of Young Tableaux which are in correspondence with directed, oriented loops obtained in the full packing limit for estimates of the free energy on the cylinder from balanced six-vertex configurations, by making use of properties of an FK-Ising representation, which provides an FK-Ising representation for the Ashkin-Teller model, estimates for crossing probabilities in the six-vertex model under sloped boundary conditions are propagated to obtain similar crossing estimates, under appropriate boundary conditions, for the Ashkin-Teller model. With a Proposition presented in the $1997$ work due to Pftizer and Velenik, crossing probabilities from the Ashkin-Teller model are shown to admit approximately equivalent representations as $\sigma$ connectivity events pushed forwards under the generalized random-cluster measure.  \footnote{\textbf{MSC Class}: 60K35; 82B02; 82D02; 82B30; 82B43; 82B44}
\end{abstract}

\section{Introduction}

\subsection{Overview}

The six-vertex model defined on the lattice was originally introduced by Pauling {\color{blue}[22]} to study thermodynamic properties of ice {\color{blue}[31]} and after having ferroelectric and antiferroelectric phases identified in {\color{blue}[20]} through a computation of the free energy, has subsequently been examined to explore conjectured universality results from fluctuations of solutions to the Kardar-Parisi-Zhang (KPZ) equation {\color{blue}[8]}, large scale simulations of Schramm-Loewner Evolution (SLE) curves {\color{blue}[19]}, connections to discontinuity of the phase transition for the random cluster model for $q > 4$ {\color{blue}[24]}, variants of the model through modifications to the ice rule with the 4, 5, 8, $\&$ 20-vertex models {\color{blue}[1,3,4,7,9]}, computational complexity of the partition function {\color{blue}[5,6]}, and connections to the corners process from random matrix theory {\color{blue}[10]}, in addition to other interpretations from the Yang-Baxter equation and Grassmannian {\color{blue}[17,18,26,30,31]}. 

To further study sloped boundary conditions for the six-vertex model, for establishing logarithmic delocalization over the torus, from previously mentioned constraints on parameters $a,b,c$, we obtain balanced six-vertex configurations that are measurable under the balanced probability measure $\textbf{P}^{Bal}_{\textbf{T}^N} \big[ \cdot \big]$ that is obtained from the measure $\textbf{P}_{\textbf{T}^N}\big[ \cdot \big]$ supported over unbalanced six-vertex configurations. Furthermore, to establish whether the height function logarithmically delocalizes for such boundary conditions rather than the complementary localization result for the height function for parameters $c>2$, we implement the following series of steps described in a summary of the paper organization below. Predominantly, arguments rely on providing Russo-Seymour-Welsh (RSW) estimates for crossing probabilities in strips of the square lattice and in the cylinder, which has recently attracted significant attention for other models, including the random cluster model {\color{blue}[13,15]}, whose phase transition has been characterized depending upon the number of colors that is captured with the cluster-weight parameter $q$ {\color{blue}[14,23]}, in addition to the dilute Potts model, an example of one model in Statistical Mechanics that is not self-dual at the critical point, for obtaining RSW estimates without making use of self-duality arguments {\color{blue}[25]} in opposition to classical arguments which are completely reliant upon duality  {\color{blue}[27,28]}. 

\subsection{Paper organization}

\noindent In \textit{1.3}, following statements of strict positivity of annulus crossing probabilities, lower and upper bounds on the variance of the height function, in addition to a lower bound for the annulus crossing probability expressed in terms of contributions from an exponential in the free energy $f_c$, we state corresponding results for sloped boundary conditions, parts of which the authors of {\color{blue}[11]} comment are expected to be similar. Following a statement of the (CBC), (FKG) and (SMP) properties of the six-vertex probability measure that are used for RSW arguments over the strip and cylinder, in \textit{2.2} we introduce objects used for our analyses of the six-vertex model in the strip first, including collections of faces that are frozen, denoted as \textit{freezing clusters}, as well as the \textit{inner} and \textit{outer diameters} of \textit{freezing clusters}, which are provided in order to demonstrate how macroscopic qualities of sloped Gibbs states in the six-vertex model impact crossing probability estimates. Such sloped configurations differ from their flat counterparts from the fact that the latter consist of infinitely many disjoint oriented circuits, as described in a characterization of different phases of the six-vertex model, namely the ferroelectric and disordered phase, as presented in {\color{blue}[16]}. To ensure that strips of the square lattice, in the first setting of RSW arguments, are conducive for \textit{sloped fence} and \textit{sloped ridge} events introduced in \textit{4.6} for the cylindrical RSW arguments, following \textbf{Lemma} \textit{4.4} we also introduce \textit{sloped-boundary domains} which are embedded within the strip, and are collections of subsets of faces bound within left, right, top, and bottom demarcations.

With these objects, we argue that modified RSW results in the strip, under sloped boundary conditions, hold, in which we introduce restrictions to domains within the strip over which long vertical bridgings occur with good probability, with results that are introduced for quantifying crossings across \textit{symmetric domains} in the strip that are obtained before the beginning of \textit{3}, in the case of non-intersecting, and intersecting, boundaries. By introducing a decomposition of boundary conditions, we are also capable of determining ways in which paths of faces can avoid \textit{freezing clusters} which can present compilications to RSW arguments in the strip, not only through preventing long vertical bridgings from occurring, but also through preventing \textit{cylindrical domains} from satisfying the \textit{tight} and \textit{good} properties that are introduced later with a list of $13$ properties in \textit{4.4} for a description of several characteristics of \textit{cylindrical domains}, which share in several characteristics of \textit{symmetric domains} initially considered in the strip. Before quantifying crossing probabilities over the cylinder, however, in \textit{2.5} we introduce \textbf{Proposition} \textit{2.7.5} before transferring estimates from the strip crossing probability to the annulus crossing probability at the beginning of \textit{3}. Such arguments allow us to obtain a lower bound for the annulus crossing probability that is provided in \textit{Theorem 6V 1} which is dependent upon a strictly positive parameter $\delta^{\prime}$ expressed in the connectivity event, which we take to be strictly smaller than a strictly positive parameter $\delta$ introduced in {\color{blue}[11]} for RSW results for the strip that are propagated to obtain RSW results for the cylinder. After having established the form of the power to which the exponential is raised in the annulus crossing probability lower bound, we proceed to the cylindrical setting of the RSW argument, after which results can be extended to simply connected domains, and finally, to the torus, for establishing logarithmic delocalization of the height-function, upon analyzing crossing probabilities over the cylinder from arguments presented throughout \textit{4}.

To demonstrate that upper and lower bounds for logarithmic delocalization of the height function hold under sloped boundary conditions, within the cylindrical setting we present modifications to the loop reversing map, which is denoted with $\mathcal{T}^{\prime}$, that was first introduced in {\color{blue}[11]} for obtaining balanced six-vertex configurations over the cylinder. In \textit{4}, properties that this map satisfies are further put to use, for establishing \textit{Theorem 6V 0}, and \textit{Theorem 6V 2}, in \textit{5}, after having introduced sloped analogues for the quantities $v_n$ and $w_n$, respectively with $v^{\xi^{\mathrm{Sloped}}}_n$ and $w^{\xi^{\mathrm{Sloped}}}_n$, which also were originally introduced in {\color{blue}[11]} under flat boundary conditions. The sloped counterparts of the two aforementioned quantities establish that logarithmic delocalization holds under sloped boundary conditions, provided in \textit{5.3}, following several inductive arguments which guarantee that annulus crossings, denoted with $\mathcal{O}$, and annulus $\mathrm{x}$-crossings, denoted with $\mathcal{O}^{\mathrm{x}}$, occur along sufficiently high level lines of either the height function, or of the level lines of the absolute value of the height function. Following the completion of arguments for establishing logarithmic delocalization of the height function under sloped boundary conditions, beginning in \textit{6} we introduce the Ashkin-Teller model, the spin-representation of the six-vertex model, and propagate results, when applicable, only for crossing probability estimates over strips of the square lattice. One significant difference in formulating crossing probability arguments between the six-vertex and its spin representation is that for the latter, the FKG inequality only holds on the marginal of spin configurations, which is related to the definition of the mixed Ashkin-Teller model on two dual lattices, namely the odd and even faces of the square lattice, over which connectivity events are probabilistically quantified. Most importantly, the positive association condition holding for the \textit{marginals} of mixed spin-configurations holds for the identical parameter regime as in the arguments that were obtained for RSW estimates under flat boundary conditions, namely that the weights satisfy the condition $a,b \leq c$, which is presented for RSW arguments of flat boundary conditions for the six-vertex model in {\color{blue}[11]}, in which the authors state that the parameter regime $\mathrm{max}\{a, b \} \leq c$ is required. The Ashkin-Teller model represents the special case $a=b\equiv 1$, which is also consistent with the requirement for RSW arguments for the six-vertex model under flat boundary conditions provided in {\color{blue}[11]}. Finally, for two more models of interest, from results provided by Pftizer and Velenik in {\color{blue}[23]}, in \textit{7} and \textit{8} we explore connections between the random-cluster representation of the Ashkin-Teller model, in addition to the $(q_{\sigma}, q_{\tau})$-cubic models which share connections with anisotropic $N$-vector models.

In other models of statistical mechanics, RSW results have been proposed in earlier works, {\color{blue}[15]}, in which the authors demonstrated that a suitable class of symmetric domains satisfying several desirable properties ensures that crossing probabilities occur with good probability over specifically chosen finite volumes. In the setting of the six-vertex model, and other closely related models, RSW results take place in both the strip of the square lattice, and also over the cylinder, in order to not only obtain estimates on the decay of the free energy over domains of aspect ratio modulo $6$, but also to make connections with sharp logarithmic bounds on the variance of the height function. For sloped boundary conditions in the six-vertex model which are chosen suitably so that complete freezing within the domain does not occur (ie, that the slopes of the Gibbs configuration do not belong to the boundary of the rational points of $\big[ -1 , 1 \big] \times \big[ -1 , 1 \big]$), symmetric domains with several additional properties must be considered. Throughout the paper, the modulo $6$ condition of symmetric domains over the cylinder refers to a statement originally provided in {\color{blue}[11]}, which is restated in the next section as \textbf{Theorem} \textit{1.5}, that the probability of a certain annulus event occuring over the cylinder can be lower bounded by the exponential of a suitable free energy function for the six-vertex model, given factors of $6$ which appear in the dimensions of the annulus event in the upper bound that are pushed forward under the six-vertex probability measure $\textbf{P} \big[ \cdot \big]$.

The properties that such domains satisfy, building upon previous RSW arguments for demonstrating that the height function of the six-vertex model delocalizes under sufficiently flat boundary conditions, {\color{blue}[11]}, partially include: constructing symmetric domains in finite volumes of the strip, excluding collection of frozen faces which are denoted with a countable collection of \textit{freezing clusters} $\mathscr{F}\mathscr{C}$; obtaining constants for which vertical crossing events across the strip occur, with sufficiently good probability, that are dependent upon the occurrence of events in which multiple \textit{freezing clusters} in the strip are disconnected, due to the fact that countably many \textit{freezing clusters} in the strip being connected to one another could form a blocking surface preventing crossings with the strip from occurring with sufficiently good probability in the infinite volume limit; incorporating additional properties rather than the tightness, good, and related properties, of cylindrical domains, which, in opposition to cylindrical domains originally introduced in {\color{blue}[11]}, do not respect the modulo $6$ condition which appears in the lower bound for the exponential of the free energy decay of the six-vertex model under sufficiently flat boundary conditions. Regardless of these differences in the symmetric domains for sloped Gibbs states over both the strip and cylinder, we achieve a logarithmic delocalization result which parallels that provided in {\color{blue}[11]}, with appropriately adjusted constants.

The notions that we expand upon in this paper for the six-vertex model under sloped boundary conditions naturally allow for extensions to other models, in which we also explore the Ashkin-Teller, generalized random-cluster model, and $\big( q_{\sigma} , q_{\tau} \big)$-spin model. In comparison to the six-vertex model under flat and sloped boundary conditions alike, these other models satisfy very similar properties of the Spatial Markov Property, Comparison between Boundary Conditions, and Fortuin-Kestelyn-Ginibre lattice condition through the FKG inequality. However, especially in the case of the Ashkin-Teller model, the correspondence between the six-vertex model and the Ashkin-Teller model only exists along the self-dual line of parameters, in which the exponential of one Potts coupling constant can be expressed in terms of the hyperbolic sine of another Potts coupling constant. More signficantly, the three main ingredients for obtaining RSW results for the six-vertex model {\color{blue}[11]}, random-cluster model {\color{blue}[15]}, and dilute Potts model {\color{blue}[25]}, each hold for all faces over the square lattice in opposition to only holding over collections of even and odd faces of the square lattice for the Ashkin-Teller model, which poses significant constraints on the random geometry of configurations that can be drawn from the sample space for which crossings across strips of the square lattice occur with sufficiently good probability. For the generalized random-cluster and $\big( q_{\sigma} , q_{\tau} \big)$-spin models, the main ingredient allowing for an extension of RSW results from the Ashkin-Teller model consists of a result from {\color{blue}[23]}, in which, from a random-cluster representation of the Ashkin-Teller model, the authors provide a correspondence between spin-spin correlation type quantities under one probability measure to connectivity properties under another probability measure. By leveraging this approximate correspondence, it is possible to rephrase arguments for weakened crossing probability estimates in the Ashkin-Teller model over strips of the square lattice to arguments for weakened crossing probability estimates for the generalized random-cluster, and $\big( q_{\sigma} , q_{\tau} \big)$, models, over the same environment.

Beyond strips of the square lattice for which RSW results are first obtained, several components of RSW arguments for relating crossing probability estimates from the cylinder back to the torus for obtaining the logarithmic delocalization result under sloped boundary conditions parallel arguments first presented in {\color{blue}[11]} for establishing that logarithmic delocalization holds under sufficiently flat boundary conditions. For other models of interest besides the six-vertex model which are considered in forthcoming arguments, encoding boundary conditions in the Ashkin-Teller, generalized random-cluster, and $\big( q_{\sigma} , q_{\tau} \big)$, models largely differs. While encoding boundary conditions in the six-vertex model is dependent upon the slope of the height function, which has several alternate characterizations, boundary conditions for the Ashkin-Teller model largely differ, from the fact that there are only $++$ or $--$ sets of boundary conditions that can be enforced over the boundary of finite volumes so that the positive association (FKG) inequality holds. From a suitable restriction of parameters to the self-dual line of the Ashkin-Teller model which makes it possible for one to relate the six-vertex and Ashkin-Teller models to each other, one loses the interpretation of the slope of the height function for the six-vertex model, which does not make it immediately clear that RSW results for the Askin-Teller model can be obtained as RSW results for the six-vertex model are, from characteristics of the slope of the height function. Nevertheless, it is possible to reinterpret RSW results over strips of the square lattice under $++$ and $--$ boundary conditions of the Ashkin-Teller probability measure, given the caveat that crossing probability events over strips of the square lattice are argued to hold over collections of even and odd faces over the square lattice, hence constituting a weakening of crossing probability estimates for the six-vertex model when restricting to the self-dual line of parameters for the Ashkin-Teller model. Given this correspondence, and lack of information from the slope of the height function that is encoded in boundary conditions for the Ashkin-Teller model, it is possible to further rephrase results for crossing probability estimates over the strip of the Ashkin-Teller to crossing probability estimates over the strip for the generalized random-cluster and $\big( q_{\sigma} , q_{\tau} \big)$-spin model, which are weakened for the same fundamental reason - the fact that the positive association (ie, the FKG inequality) does not hold over the entire square lattice but rather only over the even and odd square sublattices.

The results obtained from analyses of crossing probabilities for the six-vertex model under sloped boundary conditions hold promise for future studies, mostly notably in being able to characterize the behaviors of several models together simultaneously. As presented in several recent research articles in the areas of Statistical Mechanics, Probability, and Mathematical Physics, it is of upmost importance, when determining whether analogues, or reinterpretations, of one result can be applied from one model to another, the manner in which boundary conditions, correlations, and related quantities are encoded in the probability measure of the model.

\subsection{Statements of previous six-vertex results}

We define the cylindrical six-vertex model and from such configurations, distinguish between localization and delocalization results. As mentioned previously, in comparison to the sample space the sample space $\Omega$ of all six-vertex configurations, cylindrical six-vertex configurations are not balanced at each level. Before defining the probability measure for configurations supported over the cylinder, we recall that some six-vertex configuration $\omega$ on $\textbf{T}_N = (V(\textbf{T}_N) , E(\textbf{T}_N))$, for $N$ even, is weighed according to the following product over each possible vertex type,

\begin{figure}
\begin{align*}
\includegraphics[width=0.53\columnwidth]{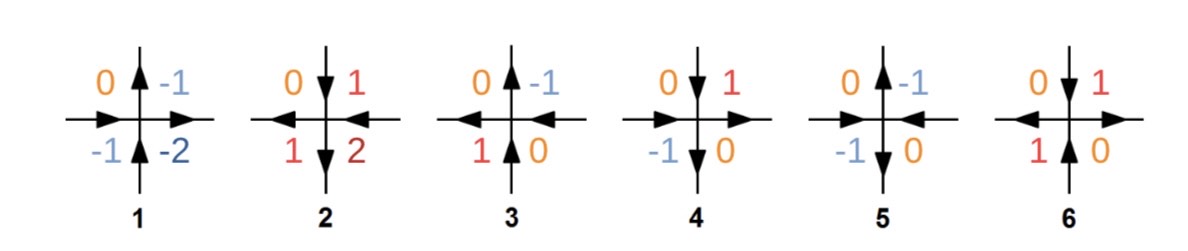}
\end{align*}
\caption{Each possible vertex for the six-vertex model, adapted from ${\color{blue}[11]}$.}
\end{figure}

\begin{align*}
  w_{6V}(\omega) \equiv w(\omega) =        a_1^{n_1} a_2^{n_2} b_1^{n_3} b_2^{n_4} c_1^{n_5} c_2^{n_6}  \underset{c_1 \equiv c_2 \equiv c}{\underset{b_1 \equiv b_2 \equiv b}{\underset{a_1 \equiv a_2 \equiv a}{\overset{\mathrm{Isotropic}}{\Longleftrightarrow}} }}          a^{n_1+n_2}  b^{n_3+n_4}  c^{n_5+n_6}    \text{ , }   \\
\end{align*}

\noindent for $a_1,a_2,b_1,b_2,c_1,c_2 \geq 0$, and each $n_i$ indicates the number of vertices of type $i$, from which the corresponding probability measure takes the form,

\begin{align*}
  \textbf{P}_{\textbf{T}_N}[      \omega         ]   \equiv \textbf{P}[   \omega     ]     =  \frac{w(\omega)}{Z_{\textbf{T}_N}}  \text{ , }  \\
\end{align*}

\noindent where the $Z_{\textbf{T}_N} \equiv Z = \sum_{\omega \in \Omega} w(\omega)$ is the partition function (see \textit{Figure 1}). We denote $\textbf{E}_{\textbf{T}_N} \equiv \textbf{E}$ as the accompanying expectation of the probability measure $\textbf{P}$. For balanced 6V configurations that cross each level of the torus the same number of times in the upwards and downwards directions, the measure of such configurations takes the form,

\begin{align*}
         \textbf{P}^{Bal}_{\textbf{T}_N}[\text{ }  \cdot \text{ } ]  = \textbf{P}^{\mathrm{B}}[ \text{ }  \cdot \text{ }  | \text{ }  \omega \in \Omega^{\mathrm{Bal}} ]         \text{ , } \\
    \end{align*}
    
    \noindent where the sample space $\Omega^{\mathrm{Bal}} \subset \Omega$, and $\textbf{E}^{\mathrm{Bal}}_{\textbf{T}_N} \equiv \textbf{E}^{\mathrm{B}}$ is the accompanying expectation. Moving forwards, in the weight $w$ for each six-vertex configuration we set $a \equiv a_1 \equiv a_2$, $b \equiv b_1 \equiv b_2$, and $c \equiv c_1 \equiv c_2$. Finally, we denote $h: F(\textbf{T}_N) \rightarrow \textbf{Z}$ as the six-vertex height function, which as a graph homomorphism is oriented along incoming edges of vertices on the square lattice so that $h$ is greater to the left of the incoming edge rather than to the right (see \textit{Figure 1} above). Under the balanced expectation, the square of the difference of the height function for $x,y \in F(\textbf{T}_N)$ is equivalent to,
    
    \begin{align*}
      \textbf{E}^{\mathrm{B}}[    (h(x) - h(y))^2      ]  = \mathrm{Var}^{\mathrm{B}}[     h(x) - h(y)  ] \text{ , } \\
    \end{align*}

   \noindent  where $\mathrm{Var}^{\mathrm{B}}$ denotes the balanced variance under $\textbf{P}^{\mathrm{B}}$. As mentioned in the previous section, to study the delocalization, or localization, behaviors of the height function, we introduce the following previously established results.

    \bigskip \noindent \textbf{Theorem} $\textit{1.1}$ (\textit{delocalization phase of the height function with a tight logarithmic bound on the balanced expectation} {\color{blue}[11]}, Theorem $1.1$): For $c \in [1,2]$, there exits constants $c,C > 0$ such that the following expectation of balanced six-vertex configurations satisfies,
    
    \begin{align*}
    c \text{ } \mathrm{log}(  d(x,y))   \leq  \textbf{E}^{\mathrm{B}}[  (h(x) - h(y))^2      ]  \leq  C \text{ }  \mathrm{log}(d(x,y)) \text{ , } \\
    \end{align*}
    
    \noindent where $N$ is even, $x,y \in F(\textbf{T}_N)$, and $d(x,y) \geq 2$, for the graph distance of the dual graph $\textbf{T}^{*}_N$ of $\textbf{T}_N$.

    \bigskip
    
    \noindent \textbf{Theorem} $\textit{1.2}$ (\textit{localiztion phase of the height function} {\color{blue}[11]}, Theorem $1.2$): For $c >2$, there exists some $C \in (0,\infty)$ such that the six-vertex height function under the balanced expectation satisfies,
    
    \begin{align*}
         \textbf{E}^{\mathrm{B}}[    (h(x) - h(y))^2    ]     \leq C \text{ , } \\
    \end{align*}
    
    \noindent for $N$ even and, $x,y \in F(\textbf{T}_N)$.

    \bigskip
    
    \noindent For planar domains, to construct symmetric domains similar to arguments provided in {\color{blue}[25]} that are applicable to other models, we introduce annuli crossing events. For $n \in \textbf{Z}$ with $0 < n < N$, for the box $\Lambda_n \equiv (-n,n)^2$ we denote the annulus $A(n,N) \equiv \Lambda_N \backslash \Lambda_n$ across which crossing events occur as $\mathcal{O}_{h \geq k}(n,N)$; the crossing event demands that $h \geq k$ for the circuit of adjacent faces across which crossings occur with positive probability. From $\textbf{P}$, we denote $\textbf{P}^{\xi}$ as the six-vertex measure with boundary conditions $\xi$. We gather two more related results concerning a uniformly positive lower bound for circuit probabilities, in addition to logarithmic variance of $h$.

    \bigskip

    \noindent \textbf{Theorem} $\textit{1.3}$ (\textit{lower bound on circuit crossing probabilities} {\color{blue}[11]}, Theorem $1.4$): From the annuli crossing event $\mathcal{O}_{h \geq k}(n,N)$, for $c \in [1,2]$ and every $k,l > 0$, there exists $c \equiv c^{\prime} \equiv c^{\prime}(c,l,k) > 0$, under admissible boundary conditions $\xi$ on $\partial D$ (or a boundary subset $\mathcal{D} \subset \partial D$), with $|\xi| \leq l$, the following circuit probability is uniformly bounded below by $\textbf{c}$, 
    
    \begin{align*}
            \textbf{P}^{\xi}_D[     \mathcal{O}_{h\geq k}(n,2n)    ] \geq \textbf{c}   \text{ , } \\
    \end{align*}

    \noindent where $\textbf{P}$ is supported over $\Lambda_{2n} \subset D$. 
    
    \bigskip

    \noindent Across a restricted scale in the annulus, we provide a similar statement under sloped boundary conditions.

    \bigskip
    
    \noindent \textbf{Theorem} \textit{6V 0} (\textit{strict positivity of the annulus crossing probability under sloped boundary conditions across a restricted scale}). Across the restricted annulus $A( n, n^{\prime} n ) \subset A(n , 2n)$, for $n^{\prime} < 2$,

    \begin{align*}
        \textbf{P}^{\xi^{\mathrm{Sloped}}}_D [   \mathcal{O}_{h \geq k}(n,n^{\prime} n)     ] \text{ } \geq \text{ }   \textbf{c}^{\prime\prime}     \text{ } \text{ , } 
    \end{align*}

    \noindent for some $\textbf{c}^{\prime\prime}$ so that $\textbf{c}^{\prime\prime} \neq \textbf{c}$.

    \bigskip
    
    \noindent \textbf{Corollary} $\textit{1}$ (\textit{tight logarithmic bound on the variance of the height function} {\color{blue}[11]}, Corollary 1.5): For $c \in [1,2]$, there exists constants $c,C >0$ such that for a discrete domain $D$, for any admissible boundary condition $\xi$ on $\partial D$ and any face $x \in D \backslash \partial D$, the variance under $\textbf{P}$ with respect to $\xi$ satisfies,
    
    \begin{align*}
         c \text{ } \mathrm{log}  ( d(x,\partial D)) - 4 l^2 \leq \mathrm{Var}^{\xi}_D ( h(x)) \leq C \text{ } \mathrm{log}( d(x , \partial D) ) + 4 l^2     \text{ , } \\
    \end{align*}

    \noindent where in the upper and lower bounds, $l = \mathrm{max}_{y \in \partial D} |\xi(y)|$.

    \bigskip
    
    \noindent \textbf{Corollary} \textit{6V 1} (\textit{accompanying tight logarithmic bound on the sloped variance of the height function from \textbf{Corollary} \textit{1} above}). Under the same assumptions of \textbf{Corollary} \textit{1}, the sloped variance satisfies, under sloped boundary conditions,
    
    \begin{align*}
       c^{\prime} \mathrm{log} ( d(x,\partial D))       - 4 |\xi^{\mathrm{Sloped}}|^2     \text{ }  \leq \text{ }  \mathrm{Var}^{\mathrm{Sloped}}   \big[                       h(x)      \big] \text{ }   \leq            c^{\prime} \mathrm{log} ( d(x,\partial D))         + 4 |\xi^{\mathrm{Sloped}}|^2      \text{ } \text{ , } 
    \end{align*}
    
    \noindent for parameters $\xi^{\mathrm{Sloped}} \sim \textbf{B}\textbf{C}^{\mathrm{Sloped}}$ and $c^{\prime} \neq c$.

    \bigskip
    
    \noindent In comparison to the measures $\textbf{P}$ and $\textbf{P}^{\mathrm{B}}$, to apply RSW arguments for crossing probabilities over the cylindrical square lattice $O_{N,M}$ whose height consists of $M$ faces, and whose perimeter consists of $N$ faces (also for $N$ even), we introduce the partition function,
    
    \begin{align*}
      Z^{s}_{N,M} \equiv Z^{s} = \sum_{\omega \in \Omega^s} w(\omega)  \text{ , } \\
    \end{align*}

    \noindent corresponding to six-vertex configurations defined on the cylinder for which every row of $N$ faces is crossing $2 \lceil s \rceil $ more times in the upwards direction than in the downwards direction. To conclude the section, we define the cylindrical 6V free energy and the lower bound estimate for circuit probabilities.

    \bigskip
    
    \noindent \textbf{Theorem} $\textit{1.4}$ (\textit{cylindrical six-vertex free energy} {\color{blue}[11]}, Theorem $1.6$): For $c > 0$, for the unbalance parameter $\alpha$ the cylindrical free energy is of the form, for $M$, and then $N \longrightarrow +\infty$,
    
    \begin{align*}
      \underset{N \longrightarrow +\infty}{\underset{{N \text{ even}}}{\mathrm{lim \text{ }  inf}}} \underset{{M \longrightarrow +\infty}}{\mathrm{lim \text{ } inf}} \big( \text{ } NM \text{ } \big)^{-1} \text{ } \mathrm{log} \text{ } \big(   Z^{\alpha}_{N,M} \big) \text{ } =\text{ }  f_c(\alpha)  \text{ }  \text{ , } \\
    \end{align*}
    
    \noindent where the function $f_c : (- \frac{1}{2} , \frac{1}{2}) \longrightarrow \textbf{R}^{+}$. For another regime of $c$, namely $0 < c \leq 2$, for every $\alpha \in (- \frac{1}{2}, \frac{1}{2})$ the free energy for the six-vertex model satisfies,
    
    \begin{align*}
            f_c(\alpha) \leq f_c(0) - C \alpha^2   \text{ , } \\
    \end{align*}

    \noindent where $C \equiv C(c) > 0$.

    \bigskip
    
    \noindent \textbf{Theorem} $\textit{1.5}$ (\textit{lower bound for height function circuit probabilities} {\color{blue}[11]} Theorem $1.7$): For all integers $k \leq r$, given $k$ sufficiently large and $c \geq 1$, there exists $\eta,c , C > 0$ so that,
    
    \begin{align*}
       \textbf{P}^{0,1}_{\Lambda_{12r}} [      \mathcal{O}_{h\geq ck}(12r,6r)    ] \geq  c \text{ } \mathrm{exp}\big[   C r^2 \big(     f_c ( \frac{k }{\eta  r} ) - f_c(0)    \big)   \big] \text{ , } \\
    \end{align*}
    
    \noindent where superscript $0,1$ on the probability measure indicates mixed boundary conditions on $\partial \Lambda_{12r}$ take values of either $0$ or $1$.
    
    \bigskip 
    
    \noindent We introduce an analogue of \textbf{Theorem} \textit{1.5} above in \textit{4.3} for sloped boundary conditions. In the following properties listed below, consider boundary conditions $\xi$, $\xi^{\prime}$ with $\xi^{\prime}\geq \xi$. Besides the previously mentioned results, the six-vertex model enjoys the two following properties, the first of which provides a lower bound for the following,
    
    \begin{align*}
 \textbf{P}^{\xi}[ F(h) \text{ } G(h) ] \text{ } \geq  \textbf{P}^{\xi}[ F(h) ] \text{ } \textbf{P}^{\xi} [ G(h) ]   \Longleftrightarrow  \textbf{E}^{\xi}[ F(h) \text{ } G(h) ] \text{ } \geq  \textbf{E}^{\xi}[ F(h) ] \text{ } \textbf{E}^{\xi} [ G(h) ]  \tag{$\mathrm{FKG}$}  \text{ , }\\
   \textbf{P}^{\xi^{\prime}}[ F(h)  ]    \geq \textbf{P}^{\xi}[ F(h) ]   \Longleftrightarrow  \textbf{E}^{\xi^{\prime}}[ F(h)  ]    \geq \textbf{E}^{\xi}[ F(h) ]  \tag{$\mathrm{CBC}$} \\ \Updownarrow \\  \textbf{P}^{\xi}[ F(|h|) \text{ } G(|h|) ] \text{ } \geq  \textbf{P}^{\xi}[ F(|h|) ] \text{ } \textbf{P}^{\xi} [ G(|h|) ]   \Longleftrightarrow  \textbf{E}^{\xi}[ F(|h|) \text{ } G(|h|) ] \text{ } \geq  \textbf{E}^{\xi}[ F(|h|) ] \text{ } \textbf{E}^{\xi} [ G(|h|) ]  \tag{$\mathrm{FKG-|h|}$}  \text{ , }\\
   \textbf{P}^{\xi^{\prime}}[ F(|h|)  ]    \geq \textbf{P}^{\xi}[ F(|h|) ]   \Longleftrightarrow  \textbf{E}^{\xi^{\prime}}[ F(|h|)  ]    \geq \textbf{E}^{\xi}[ F(|h|) ]  \tag{$\mathrm{CBC-|h|}$}    \text{ , }\\
    \end{align*}
    
    \noindent for increasing functions $F,G$ defined over the set of all admissible height functions, in addition to, for $A,B \nearrow$,

    \begin{align*}
     \textbf{P}^{\xi} [  A \cap B ] \geq     \textbf{P}^{\xi} [   A ]  \textbf{P}^{\xi} [   B ] \tag{$\mathrm{FKG}$} \text{ } \text{ , } \text{ increasing } A, B    \\ \textbf{P}^{\xi}  [ \text{ }  \cdot \text{ }  | \text{ }   h = \xi \text{ } \mathrm{ on } \text{ }  D^c \cup \partial D      ] =     \textbf{P}^{\xi_{\partial D}}[  \text{ }   \cdot   \text{ }   ]                             \tag{$\mathrm{SMP}$}   \text{ . } \\
    \end{align*}

    \noindent For subsets of the plane over which $h$ is frozen, from each possible six-vertex weight the (CBC) inequality for pairs of boundary conditions takes the form, for increasing events $\mathcal{A}$ in finite volume $\Lambda$, with associated vertex set $V( \text{ } \Lambda \text{ })$ and boundary vertices $\partial V ( \text{ } \Lambda \text{ })$, boundary conditions $\xi \leq \xi^{\prime}$, and $i > 0$, 
  
    \begin{align*}
    \mathcal{S} \text{ }    \bigg[    \sum_{x \in V( \text{ } \Lambda\text{ } ) }   \text{ } \Delta^{\xi,\xi^{\prime}}_{x}  +  \sum_{x \in \partial V( \text{ } \Lambda \text{ } ) } \text{ } \big( h^{\xi^{\prime}}|_{\partial \Lambda} - h^{\xi}|_{\partial  \Lambda \text{ } } \big) \text{ }       \Delta^{\xi,\xi^{\prime}}_{x}   \bigg]   \text{ }    \leq  \frac{\textbf{P}^{\xi}_{\Lambda}[ \mathcal{A}  ] }{ \textbf{P}^{\xi^{\prime}}_{\Lambda}[  \mathcal{A}  ]  }  \text{ }    \text{ , } \\ 
    \end{align*}
    
    \noindent for $\textbf{P}^{\xi^{\prime}}_{\Lambda}[ \text{ } \mathcal{A} \text{ }  ] \neq 0$, where, for a collection of faces $\big\{ \mathcal{F}_k \big\}_{k \in \textbf{N}}$ with $\cup_k \mathcal{F}_k = \Lambda$,

    \begin{align*}
        \mathcal{S} \equiv    {\underset{ \forall \xi^{\prime} \text{ } , \text{ }  \mathcal{F}^{\xi^{\prime}}_k \text{ frozen }: \text{ } \mathcal{F}^{\xi^{\prime}}_k \subsetneq \Lambda  }{\underset{k \in \textbf{N} \text{ } }{\mathrm{sup}}}}  \text{ }          \frac{ \big|  \{ v : v \in \mathcal{F}^{\xi^{\prime}}_k\} \big| }{\big| \{ v : v \in V( \text{ } \Lambda\text{ } )  \} \big| }      \text{ } \text{ , }   \\
    \end{align*}

    \noindent in addition to the following quantity,
    
    \begin{align*}
      \text{ }   \Delta^{\xi,\xi^{\prime}}_{x} \equiv \text{ }   \frac{a_x^{(n^{\xi^{\prime}}_1 + n^{\xi^{\prime}}_2) \textbf{1}_{(v_x \equiv 1) \text{ or }  (v_x \equiv 2)} }   \text{ } b_x^{(n^{\xi^{\prime}}_3 + n^{\xi^{\prime}}_4) \textbf{1}_{(v_x \equiv 3) \text{ or } (v_x \equiv 4)}} \text{ } c_x^{(n^{\xi^{\prime}}_5 + n^{\xi^{\prime}}_6) \textbf{1}_{(v_x \equiv 5) \text{ or }  (v_x \equiv 6)} }  }{a_x^{(n^{\xi}_1 + n^{\xi}_2) \textbf{1}_{(v_x \equiv 1) \text{ or }  (v_x \equiv 2)}} \text{ } b_x^{(n^{\xi}_3 + n^{\xi}_4) \textbf{1}_{(v_x \equiv 3) \text{ or }  (v_x \equiv 4)}}  \text{ } c_x^{(n^{\xi}_5 + n^{\xi}_6) \textbf{1}_{(v_x \equiv 5) \text{ or }  (v_x \equiv 6)}} }      \text{ }   \text{ , } \\
    \end{align*}
    
    \noindent for $v_x \in V(\textbf{Z}^2)$. In the special case of the (CBC) inequality above, collections of frozen faces under boundary condition $\xi$ are denoted with $\mathcal{F}^{\xi}_i$. Besides the multiplicative prefactor in the supremum taken over a countable number of freezing domains, the summation over lattice sites of the square lattice arises from the Isotropic parameter choice, and accounts for the manner in which changes of $h$ at the boundary of the finite volume impact the values of $h$ in the bulk.

    In comparison to the dependence of boundary conditions for other models in Statistical Physics, sloped boundary conditions for the six-vertex model imply, for some finite volume $\Lambda$ of $\textbf{Z}^2$, that,

\[
\left\{\!\begin{array}{ll@{}>{{}}l} \forall F \in F \big( \Lambda \big) \text{ } \mathrm{boundary \text{ } conditions\text{ }} \xi^{\mathrm{Sloped}} \text{ } \mathrm{sloped\text{ } if}  & \xi^{\mathrm{Sloped}} \in  \textbf{Q}  \cap \big( \big( -1 , 1 \big) \times \big( - 1 , 1 \big) \big) \Longleftrightarrow   \\
  & \textbf{P}^{\xi^{\mathrm{Sloped}}}_{F ( \Lambda )}  \big[ \mathrm{all \text{ }vertices\text{ } are \text{ } of \text{ } the \text{ } same \text{ } type} \big]\equiv 0    \text{ , } \\
     \forall F \in F \big( \Lambda \big) \text{ } \mathrm{boundary \text{ } conditions\text{ }} \xi^{\mathrm{Sloped}} \text{ } \mathrm{sloped\text{ } if}  &   \xi^{\mathrm{Sloped}} \in \textbf{Q} \cap \big( \big[ -1 , 1 \big] \times \big[ - 1 , 1 \big] \big) \Longleftrightarrow     \text{ }  \\     &  \textbf{P}^{\xi^{\mathrm{Sloped}}}_{F ( \Lambda )} \big[ \mathrm{all \text{ }vertices\text{ } are \text{ } of \text{ } the \text{ } same \text{ } type} \big] \equiv 1      \text{ }  \\
\end{array}\right.
\]

    Within the parameter regime with $c \geq 1$, the FKG inequality is not satisfied. Otherwise, independently of choices for $a$ and $b$ the inequalities provide conditions for positive association amongst crossing events on the lattice, in addition to comparison between boundary conditions, conditionally upon the absolute value of the height function for some real $k$,
    
    \begin{align*}
                  \textbf{E}^{\xi} [ F(|h|)   G(|h|)] \geq \textbf{E}^{\xi} [F(|h|) ] \textbf{E}^{\xi}[G(|h|)]  \tag{$\mathrm{FKG-|h|}$}   \text{ , }   \\    \textbf{E}^{\xi^{\prime}} [ F(|h|)  ]   \geq \textbf{E}^{\xi}[F(|h|)  ]  \tag{$\mathrm{CBC-|h|}$}  \text{ . }   \\
    \end{align*}

\noindent Each pair of inequalities will be extensively used in the forthcoming construction of symmetric domains with sloped boundary conditions. The pair of inequalities above also hold conditionally upon  the value of $h$, or upon the value of $|h|$. When analyzing the Ashkin-Teller, generalized random-cluster, and $\big( q_{\sigma} , q_{\tau} \big)$-spin models, there are no such families of inequalities due to a lack of analogy for the absolute value of the height function. In spite of such differences between these models of Statistical Mechanics for which positive association properties hold, crossing probability estimates along strips of the square lattice can still be obtained, albeit with different consequences. Out of the four models considered in this work, the six-vertex model is the most rich in the sense that it allows for analyses of the sloped free energy functional over the cylinder. Due to there being no direct analogy for the slope of boundary conditions, in addition to the existence of $|h|$ for the height function of the six-vertex model, extending weakened crossing probability estimates obtained for the Ashkin-Teller, generalized random-cluster, and $\big( q_{\sigma}, q_{\\tau} \big)$-spin models for interpretations of the free energy does not remain entirely clear. 

\subsection{Beyond the minimal height function for sloped boundary conditions}

\noindent We denote the set of sloped height functions over an arbitrary finite volume $\Lambda$ as $\mathcal{H}^{\text{ } \mathrm{Sloped}}_{\Lambda}$, and the set of all sloped boundary conditions at $\textbf{BC}^{\mathrm{Sloped}}$. Under such boundary conditions, denote the probability measure $\textbf{P}^{\mathrm{Sloped}}_{\Lambda} [ \cdot ]$, for some finite volume $\Lambda \subset \textbf{Z}^2$, that is normalized by the sloped partition function,

\begin{align*}
Z^{\mathrm{Sloped}}_{M,N} \equiv    Z^{\mathrm{S}}_{M,N}   \equiv \sum_{\omega \in \Omega^{\mathrm{Sloped}}} \text{ } w ( \omega )   \equiv \sum_{\omega \in \Omega^{\mathrm{Sloped}}} \text{ } \mathrm{exp} \big[   g_c(\omega)  (    1      + o(1) ) \big]    \equiv \sum_{\omega \in \Omega^{\mathrm{Sloped}}} \mathrm{exp} \big[   g_c(\omega)  (    1      +        \mathcal{O} \big)  \big]   \text{ } \text{ } \text{ , } 
\end{align*}        

\noindent which will be further analyzed in \textit{4.3} when defining the \textit{sloped} free energy function.

\bigskip

 \noindent \textbf{Proposition} $\textit{1.1}$ (\textit{Comparison between boundary conditions under sloped boundary conditions, in which boundary conditions of the lower bound probability are taken along a finite volume boundary}): For a pair of sloped boundary conditions $\xi \leq \xi^{\prime}$ and $h \in \mathcal{H}^{\text{ } \mathrm{Sloped}}_{\Lambda}$, there exists a strictly positive constant $\mathscr{C}$ such that the probability of an increasing connectivity event $\mathcal{E}_k \equiv \{    \mathscr{F}_i           \overset{ h \geq k }{\longleftrightarrow}      \mathscr{F}_j    \}$ in cluster of faces $\mathcal{C}$ taken under the restriction of boundary conditions for any $\Lambda^{\prime} \subset \Lambda $  satisfies,
    
    \begin{align*}
       \text{ }   \textbf{P}_{\Lambda}^{\xi|_{\partial \Lambda^{\prime}}} [    \mathcal{E}_k              ]      \text{ } \leq \text{ } \mathscr{C} \text{ }  \textbf{P}_{\Lambda}^{\xi^{\prime}} [  \mathcal{E}_k   ]          \text{ }   \text{ , } \\
    \end{align*}

    \noindent where $\mathscr{F}_i, \mathscr{F}_j \in \mathcal{C}$ are faces, with $\mathscr{F}_i \cap \mathscr{F}_j = \emptyset$.
    
    \bigskip

\noindent \textit{Proof of Proposition 1.1}. It suffices to upper bound the ratio of probabilities,

\begin{align*}
      \text{ }   \frac{\textbf{P}_{\Lambda}^{\xi^{\prime}}[    \mathcal{E}_k   ]}{\textbf{P}_{\Lambda}^{\xi}[  \mathcal{E}_k   ]}         \text{ }  \text{ , } \\
\end{align*}

\noindent for the increasing event $\mathcal{E}_k$ defined in the statement, and $\textbf{P}^{\xi}_{\Lambda}[ \mathcal{E}_k] \neq 0$. The (CBC) condition implies the product of observables,

\begin{align*}
    \text{ }  \mathcal{S}^{-1} \text{ }     \bigg[   \sum_{x \in V( \text{ } \Lambda \text{ } )}   \text{ } \Delta^{\xi,\xi^{\prime}}_{x}  +  \sum_{x \in  \partial V( \text{ } \Lambda \text{ } ) }  \text{ } \big( h^{\xi^{\prime}}|_{\partial \Lambda} - h^{\xi}|_{\partial \Lambda}  \big) \text{ }       \Delta^{\xi,\xi^{\prime}}_{x}   \bigg]^{-1}            \text{ } \text{ , } \text{ }    \\
\end{align*}

\noindent in which the reciprocal of a strictly positive prefactor is strictly positive, as there can be at most countably many $\mathcal{F}_k^{\xi^{\prime}}$ with $\mathcal{F}_k^{\xi^{\prime}} \subsetneq \Lambda$ for each $k$ under $\xi^{\prime}$.

\bigskip 

\noindent Besides these contributions in the upper bound, changing the slope of boundary conditions from the values assigned by $\xi$ to those assigned by $\xi^{\prime}$ can combinatorially result in  the following probabilistic changes for the occurrence of crossing events, in which under an isotropic choice of parameters ($a\equiv b$ and $c \geq 1$) contributions from the bulk in (CBC) take the form,

\begin{align*}
    \sum_{x\in V( \text{ } \Lambda \text{ } ) } \frac{a_x^{(n^{\xi^{\prime}}_1 + n^{\xi^{\prime}}_2 - n^{\xi}_1 - n^{\xi}_2 ) \text{ }  \textbf{1}_{(v_x \equiv 1) \text{ or }  (v_x \equiv 2)}}}{b_x^{(n^{\xi}_3 + n^{\xi}_4 - n^{\xi^{\prime}}_3 - n^{\xi^{\prime}}_4 ) \text{ }  \textbf{1}_{(v_x \equiv 3) \text{ or }  (v_x \equiv 4)}}} c_x^{(n^{\xi^{\prime}}_5 + n^{\xi^{\prime}}_6 - n^{\xi}_5 - n^{\xi}_6 ) \text{ }  \textbf{1}_{(v_x \equiv 5) \text{ or }  (v_x \equiv 6)}} \text{ } \text{ . } \end{align*}

    \noindent The expression above can be rewritten as,
    
    \begin{align*}
    \text{ }     \sum_{x \in V( \text{ }  \Lambda \text{ } ) } a_x^{(n^{\xi^{\prime}}_1 + n^{\xi^{\prime}}_2 - n^{\xi}_1 - n^{\xi}_2 ) \text{ }  \textbf{1}_{(v_x \equiv 1) \text{ or }  (v_x \equiv 2)}- (n^{\xi}_3 + n^{\xi}_4 -n^{\xi^{\prime}}_3 - n^{\xi^{\prime}}_4 )  \text{ }  \textbf{1}_{(v_x \equiv 3) \text{ or }  (v_x \equiv 4)}}  \text{ }  \text{ }  c_x^{( n^{\xi^{\prime}}_5 + n^{\xi^{\prime}}_6 - n^{\xi}_5 - n^{\xi}_6 \text{ } ) \textbf{1}_{(v_x \equiv 5) \text{ or }  (v_x \equiv 6)}}       \text{ }  \text{ , }  \\
\end{align*}

\noindent while contributions from the boundary dependent summation in (CBC) take the form,

\begin{align*}
            \text{ }     \sum_{x \in \partial V( \text{ } \Lambda \text{ } ) } \text{ }  \big( h^{\xi^{\prime}}|_{\partial \Lambda} - h^{\xi}|_{\partial \Lambda}  \big)  a_x^{(n^{\xi^{\prime}}_1 + n^{\xi^{\prime}}_2 - n^{\xi}_1 - n^{\xi}_2 ) \text{ }  \textbf{1}_{(v_x \equiv 1) \text{ or }  (v_x \equiv 2)} - (n^{\xi}_3 + n^{\xi}_4 - n^{\xi^{\prime}}_3 - n^{\xi^{\prime}}_4 ) \text{ }  \textbf{1}_{(v_x \equiv 3) \text{ or }  (v_x \equiv 4)}} \\ \times   c_x^{(n^{\xi^{\prime}}_5 + n^{\xi^{\prime}}_6 - n^{\xi}_5 - n^{\xi}_6 )\text{ }  \textbf{1}_{(v_x \equiv 5) \text{ or }  (v_x \equiv 6)}}       \text{ , } \text{ }     \\
\end{align*}

\noindent where $n_1^{\xi^{\prime}} + n_2^{\xi^{\prime}} - n_1^{\xi} - n_2^{\xi}$, $n_3^{\xi} + n_4^{\xi} - n_3^{\xi^{\prime}} - n_4^{\xi^{\prime}}$, $ n_5^{\xi^{\prime}} + n_6^{\xi^{\prime}} - n_5^{\xi} - n_6^{\xi} \text{ } > 0  $, in which vertices of Type $1$ and $2$ in finite volume increase.

\noindent Otherwise, if vertices of type $5$ and $6$ simultaneously decrease under $\xi^{\prime}$ with vertices of type $3$ and $4$,  $n_5^{\xi^{\prime}} + n_6^{\xi^{\prime}} - n_5^{\xi} - n_6^{\xi} \text{ } < 0 $, and contributions from the bulk in (CBC) instead take the form,

\begin{align*}
            \sum_{x\in V( \text{ } \Lambda \text{ } )}    \frac{  a_x^{(n^{\xi^{\prime}}_1 + n^{\xi^{\prime}}_2 - n^{\xi}_1 - n^{\xi}_2) \text{ }  \textbf{1}_{(v_x \equiv 1) \text{ or }  (v_x \equiv 2)}  - (n^{\xi}_3 + n^{\xi}_4 - n^{\xi^{\prime}}_3 - n^{\xi^{\prime}}_4 )  \text{ }  \textbf{1}_{(v_x \equiv 3) \text{ or }  (v_x \equiv 4)} }   }{c_x^{(n_5^{\xi} + n_6^{\xi} - n_5^{\xi^{\prime}} - n_6^{\xi^{\prime}})  \text{ }  \textbf{1}_{(v_x \equiv 5) \text{ or }  (v_x \equiv 6)}}}    \text{ }     \text{ , } \text{ }  \\
\end{align*}

\noindent while contributions from the boundary in (CBC) take the form,

\begin{align*}
        \sum_{x \in \partial V( \text{ } \Lambda \text{ } ) } \big( h^{\xi^{\prime}}|_{\partial \Lambda} - h^{\xi}|_{\partial \Lambda}  \big)  \text{ }  \frac{  a_x^{(n^{\xi^{\prime}}_1 + n^{\xi^{\prime}}_2 - n^{\xi}_1 - n^{\xi}_2 )  \text{ }  \textbf{1}_{(v_x \equiv 1) \text{ or }  (v_x \equiv 2)} - (n^{\xi}_3 + n^{\xi}_4 - n^{\xi^{\prime}}_3 - n^{\xi^{\prime}}_4 )  \text{ }  \textbf{1}_{(v_x \equiv 3) \text{ or }  (v_x \equiv 4)} }   }{c_x^{( n_5^{\xi} + n_6^{\xi} - n_5^{\xi^{\prime}} - n_6^{\xi^{\prime}}) \text{ }  \textbf{1}_{(v_x \equiv 5) \text{ or }  (v_x \equiv 6)}}}      \text{ }   \text{ , } \text{ } \\
\end{align*}

\noindent In remaining cases, there can be an asymmetry between the number of vertices of Type $3$ and $4$ under $\xi$ and $\xi^{\prime}$, in which contributions from the bulk in (CBC) take the same form, but are instead subject to the constraints in which vertices of Type $3$ and $4$, or of Type $5$ and $6$, with $n_3^{\xi^{\prime}} +n_4^{\xi^{\prime}} - n_3^{\xi} - n_4^{\xi}- n_3^{\xi} - n_4^{\xi} - n_1^{\xi} - n_2^{\xi} +n_3^{\xi^{\prime}} + n_4^{\xi^{\prime}}$, $n_5^{\xi^{\prime}} + n_6^{\xi^{\prime}} - n_5^{\xi} - n_6^{\xi} > 0$ or with $n_3^{\xi^{\prime}} +n_4^{\xi^{\prime}} - n_3^{\xi} - n_4^{\xi}- n_3^{\xi} - n_4^{\xi} - n_1^{\xi} - n_2^{\xi} +n_3^{\xi^{\prime}} + n_4^{\xi^{\prime}}>0$, $n_5^{\xi^{\prime}} + n_6^{\xi^{\prime}} - n_5^{\xi} - n_6^{\xi} < 0$. Contributions from the bulk and boundary of finite volume in each case are strictly positive, and remain strictly positive when multiplied by $\frac{1}{\mathcal{S}}$. \boxed{}

\begin{figure}
\begin{align*}
\includegraphics[width=0.32\columnwidth]{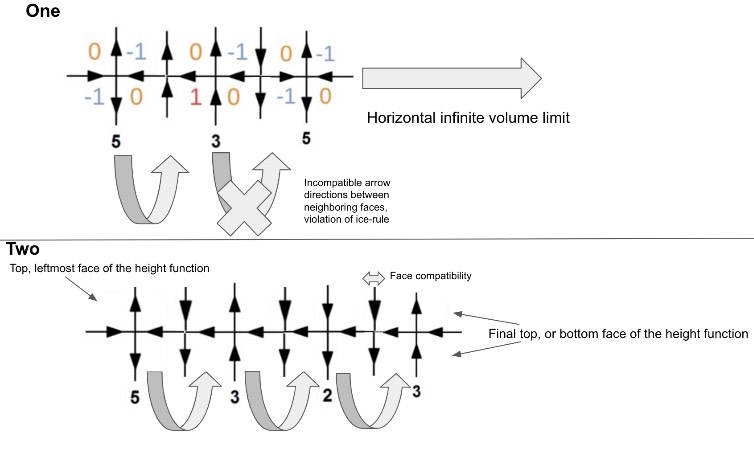}\\
\includegraphics[width=0.32\columnwidth]{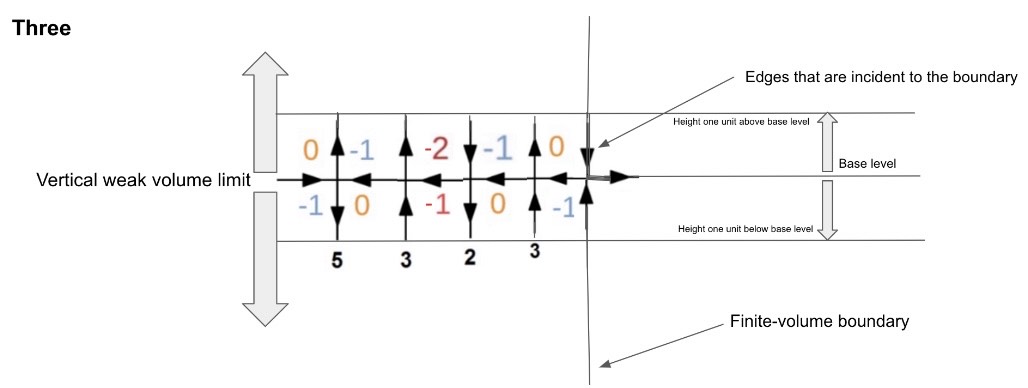}\\
\includegraphics[width=0.32\columnwidth]{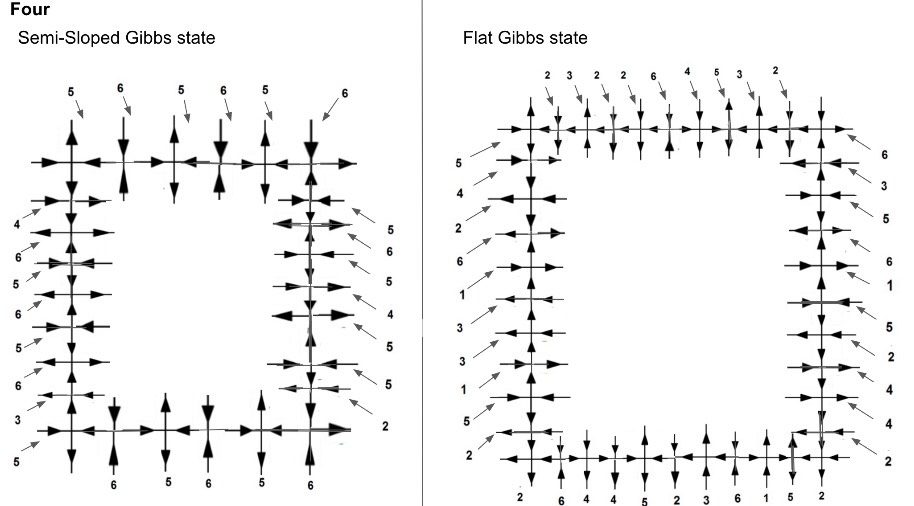}\\
\end{align*}
\caption{\textit{Distinguishing between flat and sloped Gibbs states in the six-vertex model}. In \textit{panel one}, six-vertex configurations on the square lattice in the horizontal weak volume limit must be individually composed of vertices satisfying the ice rule. In \textit{panel two}, a sequence of compatible faces is constructed, and exhibits partial freezing of the height function with repeating vertices of type $2$. In \textit{panel three}, a vertical weak volume limit taken along the finite volume boundary exhibits an interface consisting of edges oriented inwards or outwards. In \textit{panel four}, a Flatt Gibbs state is provided in comparison to the Semi-Sloped Gibbs state, which has vertices partially freezed along the boundary which are predominantly either of Type $5$ or $6$.}
\end{figure}

\bigskip

\noindent \textbf{Corollary} $\textit{1.2}$ (\textit{Applying (CBC) and (SMP) properties of the six-vertex probability measure for an upper bound of the probability of height-function crossing events at higher levels with the probability of height-function crossing events at lower levels}): For sloped boundary conditions $\xi \leq \xi^{\prime}$ and some $j>0$, and finite volumes $\Lambda^{\prime} \subset \Lambda$ such that $\xi$ and $\xi^{\prime}$ are equal on $\partial \Lambda \cap \partial \Lambda^{\prime}$, the probability of a height-function crossing event $\mathcal{E}_{k} \equiv \{    \mathscr{F}_i            \text{ }  \overset{ h \geq k }{\longleftrightarrow} \text{ }          \mathscr{F}_j       \} $ occurring satisfies the upper bound,

\begin{align*}
  \frac{1}{\mathscr{C}} \text{ }  \textbf{P}^{\xi^{\prime}}_{\Lambda^{\prime}}[ \mathcal{E}_{k+j}   ] \text{ } \leq \text{ }  \textbf{P}^{\xi}_{\Lambda}[   \mathcal{E}_k    ] \text{ }              \text{ } \text{ , } \text{ } \\
\end{align*}

\noindent where the probability measure in the lower bound is supported over the finite volume subdomain, and the probability measure in the upper bound is supported over the entire finite volume. $\mathcal{E}_k$ is the height-function crossing event defined in \textit{Proposition 1.1}, and $\mathscr{C}$ is also defined in \textit{Proposition 1.1} from the upper bound.

\bigskip

\noindent \textit{Proof of Corollary 1.2}. Applying (CBC) to the probability of $\mathcal{E}_{k+j}$ occurring in the lower bound of the \textbf{Corollary} yields, for $\xi \geq \xi^{\prime}$,

\begin{align*}
         \textbf{P}^{\xi}_{\Lambda^{\prime}}[ \mathcal{E}_{k+j}  ] \text{ } \leq \text{ } \textbf{P}^{\xi^{\prime}}_{\Lambda^{\prime}} [       \mathcal{E}_{k+j} ]        \text{ } \text{ , } \text{ }    \\
\end{align*}

\noindent which yields the lower bound,

\begin{align*}
         \textbf{P}^{\xi^{\prime}}_{\Lambda^{\prime}} [            \mathcal{E}_{k+j}        ]  \text{ } \overset{(\textbf{Proposition}\text{ }  \textit{1.1})}{\geq} \text{ } \mathscr{C} \text{ }  \textbf{P}^{\xi|_{\partial \Lambda^{\prime}}}_{\Lambda^{\prime}} [   \mathcal{E}_{k+j}   ]  \text{ } \overset{\mathrm{(SMP)}}{\geq} \text{ }    \mathscr{C} \text{ }       \textbf{P}^{\xi|_{\partial \Lambda}}_{\Lambda}   [   \mathcal{E}_{k+j}   ]   \text{ } \text{ . } \text{ } \\
\end{align*}

\noindent We conclude with the following lower bound for boundary conditions along $\partial \Lambda^{\prime}$, as collecting previous estimates yields,

\begin{align*}
      \textbf{P}^{\xi|_{\partial \Lambda}}_{\Lambda} [  \mathcal{E}_{k+j}  ]       \text{ }  \overset{\mathrm{(CBC)}}{\leq} \text{ }      \frac{1}{\mathscr{C}} \text{ } \textbf{P}^{\xi}_{\Lambda}    [   \mathcal{E}_{k+j}  ]  \text{ } \equiv \text{ }   \frac{1}{\mathscr{C}} \text{ }  \textbf{P}^{\xi -j}_{\Lambda} [ \mathcal{E}_k     ]              \text{ } \Leftrightarrow \text{ }    \textbf{P}^{\xi^{\prime}}_{\Lambda^{\prime}} [            \mathcal{E}_{k+j}        ]  \text{ } \leq   \text{ }         \frac{1}{\mathscr{C}} \text{ }  \textbf{P}^{\xi -j}_{\Lambda} [  \mathcal{E}_k  ]  \text{ } \equiv \text{ } \frac{1}{\mathscr{C}} \text{ }  \textbf{P}^{\xi}_{\Lambda} [  \mathcal{E}_{k}    ]         \text{ } \text{ , } \text{ } \\
\end{align*}

\noindent hence concluding the argument. \boxed{}

\bigskip

\noindent To proceed, we characterize the crossing probability across the annulus and variance of the height function, each of which are dependent upon the slope of the boundary conditions. For boundary conditions in which $h$ is not sufficiently flat, it is expected that $\mathrm{Var}^{\xi}(h(x))$ should also contain a lower order term that scales with $\frac{d(x,y)}{N}$,

\begin{align*}
    \mathrm{Var}^{\xi} ( h(x) ) =   \text{ }   \mathrm{log}(  d( x , y  )  )    \text{ } + \text{ } \frac{d(x,y)}{N} \text{ }   \text{ , }  \\
\end{align*}

\noindent where $x,y$ belong to disjoint faces, and $N$ is one degree of freedom about which the infinite weak volume limit is taken (recall the infinite limit volume limit given in \textit{Theorem 4}).

In the next section, we further detail estimates relating to various distributions of vertex weights which are first studied in domains in the strip, afterwards allowing us to obtain approximations of the cylindrical free energy over finite volume.

\subsection{Heuristics for the change in annuli crossing probabilities between flat and sloped boundary conditions}

To quantify the change in probability of obtaining crossings across annuli of sloped boundary conditions of the height function, we study crossing probabilities of events supported over strips of the square lattice. When taking the infinite volume limit, in the presence of sloped boundary conditions, the height function can freeze uniformly to any one of the six possible configurations within a finite volume if any of the slopes of the height function are the rational points of $\big[ -1 , 1 \big] \times \big[ -1 , 1 \big]$. Otherwise, the height function with positive probability some proportion of the faces of the height function within the finite volume will freeze. When constructing sequences of rational points for the values of the height function within a finite volume, if the final value of the height function is a rational number lying on the boundary of $\big[ -1 , 1 \big] \times \big[ -1 , 1 \big]$, then with positive probability sampling a configuration of the height function for which a crossing event occurs in a long strip decreases as there will be a nonzero number of frozen over which the height function is freezes. To counteract the possibility of blocking interfaces forming within the strip from frozen faces of the height function, crossing estimates in long strips are obtained by allowing for the possibility that countably many blocking surfaces of frozen faces - \textit{freezing clusters} - exist with positive probability that intersect the boundary of the strip. Moreover, the cardinality of such a collection of faces that are frozen in finite volume is dependent upon the value of boundary conditions of the height function along the bottom and top of the strip.

\bigskip

Define the polarization of a finite volume $\Lambda_1 \subset \textbf{Z}^2$ as the fluctuation of the height function, namely the difference $m_{\Lambda} = \mathrm{max}_{\mathcal{F}_1,\mathcal{F}_2 \in \mathcal{F}_{\textbf{Z}^2}}  h_{\mathcal{F}_1} - h_{\mathcal{F}_2}$, for faces $\mathcal{F}_1$ and $\mathcal{F}_2$ over $\Lambda$. Dependent upon the choice of finite volume, the maximum order of fluctuation of $h$ can give rise to the \textit{freezing} phenomena, in which the the random function $f$ undergoes the following change of behavior,

\[
f \equiv f(h) \text{ } \equiv
\left\{\!\begin{array}{ll@{}>{{}}l}
> 0  & \text{ for }  & \text{ }  c(k+\beta) \geq h \geq ck \text{ } \text{ , } \\
\ 0 &  \text{ for } & \text{ } ck < h  \text{ } \text{ , } \\
\end{array}\right.
\]

\noindent for connectivity event of the height function across faces. Specifically, the random $f$ is given by the probability,

\begin{align*}
      \textbf{P}^{\xi_1}_{\Lambda_1}[  A \overset{ c(k+\beta) \geq  h \geq ck}{\longleftrightarrow} B      ]      \text{ , }     \\
\end{align*}

\noindent for $\beta$ sufficiently small and boundary conditions $\xi_1$, in which is a positive probability that the faces $A,B \subset \mathcal{F}_{\textbf{Z}^2}$ are connected along a path of height for a restricted regime of parameters of the height function. Correspondingly, if one encounters another subset of the square lattice over which the maximum fluctuation of the height function does not concentrate within $c(k-\beta) \leq  h \leq ck$ on $\Lambda_2$ with $\Lambda_1 \cap \Lambda_2 = \emptyset$, the height function connectivity event,

\begin{align*}
       \textbf{P}^{\xi_2}_{\Lambda_2}   [              A \overset{ h \geq ck}{\longleftrightarrow} B   ]          \text{ , }     \\
\end{align*}

\noindent holds with positive probability, for boundary conditions $\xi_2$.

\bigskip

\noindent Under sloped boundary conditions, if the value of the slope along each face of the interface belongs to the boundary of the set of rational points of $[-1,1] \times [-1,1]$, the following transition parameterized in the height function is expected to occur, in which,

\[
\mathscr{T}_1 \equiv \left\{\!\begin{array}{ll@{}>{{}}l}
        \textbf{P}^{\xi_{\mathrm{Supcritical}}}_{\Lambda_1}[ \text{ }  A \overset{  h \geq c_{\mathrm{Supcritical}}k}{\longleftrightarrow} B     \text{ } ]       & \text{ with } h_{\xi_{\mathrm{Supcritical}}} \equiv c_{\mathrm{Supcritical}}k + \epsilon_1 \text{ } \text{ , } \\
    \textbf{P}^{\xi_{\mathrm{Near-critical}}}_{\Lambda_1}[ \text{ }  A \overset{  h \geq c_{\mathrm{Near-critical}}k}{\longleftrightarrow} B     \text{ } ]           & \text{ with } h_{\xi_{\mathrm{Supcritical}}} \equiv c_{\mathrm{Near-critical}} k + \epsilon_2 \text{ } \text{ , }\\
    \textbf{P}^{\xi_{\mathrm{Subcritical}}}_{\Lambda_1}[ \text{ }  A \overset{    h \geq c_{\mathrm{Subcritical}}k}{\longleftrightarrow} B     \text{ } ]           & \text{ with } h_{\xi_{\mathrm{Subcritical}}} \equiv c_{\mathrm{Subcritical}} k + \epsilon_3 \text{ , } \\
\end{array}\right.
\]

\noindent occurs, where $c_{\mathrm{Supcritical}} \geq c_{\mathrm{Near-critical}} \geq  c_{\mathrm{Subcritical}} \geq 0$, and $\epsilon_1,\epsilon_2,\epsilon_3 > 0$. To empirically realize distribution of six-vertex weights for which the transition above occurs, a sequencing approximating each $\epsilon_1,\epsilon_2,\epsilon_3$, which after being set to be arbitrarily large, would induce a freezing of the height function in the finite volume bulk. As a discontinuous phase transition, with respect to stipulations of the slope of $h$ on the boundary, RSW arguments can be appropriately fashioned by determining the critical height at which the connected components of the connectivity event initially form. The superscripts of $\textbf{P}$, $\xi_{\mathrm{Supcritical}}$, $\xi_{\mathrm{Near-critical}}$, $\xi_{\mathrm{Subcritical}}$, respectively correspond to boundary conditions of the subcritical, near-critical, and supercritical phases, under which there emerges a positive probability of crossing an annulus of a $\mathrm{Near-critical}$ height $c_{\mathrm{Near-critical}} k$.

\noindent The transition,

\[
\mathscr{T}_2 \equiv \left\{\!\begin{array}{ll@{}>{{}}l}
        \textbf{P}^{\xi_{\mathrm{Supcritical}_1}}_{\Lambda_1}[ \text{ }  A \overset{  h \geq c_{\mathrm{Supcritical}}k}{\longleftrightarrow} B     \text{ } ]  -   \textbf{P}^{\xi_{\mathrm{Supcritical}_2}}_{\Lambda_1}[ \text{ }  A \overset{  h \geq c_{\mathrm{Supcritical}}k}{\longleftrightarrow} B     \text{ } ]   \approx 0    & \text{, } \xi_{\mathrm{Supcritical}_1} \geq \xi_{\mathrm{Supcritical}_2} \text{ , }  \\
    \textbf{P}^{\xi_{\mathrm{Near-critical}_1}}_{\Lambda_1}[ \text{ }  A \overset{  h \geq c_{\mathrm{Near-critical}}k}{\longleftrightarrow} B     \text{ } ]  - \textbf{P}^{\xi_{\mathrm{Near-critical}_2}}_{\Lambda_1}[ \text{ }  A \overset{  h \geq c_{\mathrm{Near-critical}}k}{\longleftrightarrow} B     \text{ } ]    > c        & \text{, } \xi_{\mathrm{Near-critical}_1} \geq \xi_{\mathrm{Near-critical}_2} \text{ , } \\
    \textbf{P}^{\xi_{\mathrm{Subcritical}_1}}_{\Lambda_1}[ \text{ }  A \overset{    h \geq c_{\mathrm{Subcritical}}k}{\longleftrightarrow} B     \text{ } ]   - \textbf{P}^{\xi_{\mathrm{Subcritical}_2}}_{\Lambda_1}[ \text{ }  A \overset{    h \geq c_{\mathrm{Subcritical}}k}{\longleftrightarrow} B     \text{ } ]    \approx 1     & \text{, } \xi_{\mathrm{Subcritical}_1} \geq \xi_{\mathrm{Subcritical}_2} \text{ , }  \\
\end{array}\right.
\]

\noindent encapsulates whether there is positive probability of the height function exhibiting \textit{frozen} versus \textit{liquid} behavior in the finite volume across which the crossing probability is quantified, where $0 < c < 1$. In the Near-critical phase, the mixing time with respect to the two pairs of boundary conditions under which the crossing probability across an annulus is bounded below by a positive constant. Similar distributions of six-vertex weights, as those provided under the transition $\mathscr{T}_2$, will be analyzed in future sections. 

\subsection{Annuli crossing probabilities across spacing that is not uniform}

As a byproduct of introducing sloped boundary conditions of the height function that are not as flat as those considered in {\color{blue}[11]}, one must consider annuli of unequal spacing instead of always considering the annulus comprised of a larger box of length $2n$, enclosed within a smaller box of length $n$. In addition to the fact that the probability of a horizontal crossing event of the height function occurring in long, near infinite volume, strips decreases as the number of frozen faces increases, the probability that crossings of the height function can still be sustained results from estimates on crossings of the height function across faces that are not frozen in finite volume. Across faces of the height function that are  not frozen despite imposing boundary conditions that the slope of the height function lies within the interior of the set of rational points of $\big[ -1 , 1 \big] \times \big[ -1 , 1 \big]$, long horizontal crossings of the height function can be obtained in the presence of sloped boundary conditions. For such long horizontal crossings, sloped boundary conditions within the interior of the set of rational of $\big[ -1 , 1 \big] \times \big[ -1 , 1 \big]$ dictate that crossings of the height function could 'escape' finite volumes contained within the strip by passing through the bottom and/or top. However, given the encoding of sloped boundary conditions for the six-vertex model, long horizontal crossings occur with sufficiently good probability if the height function encounters as few frozen faces as possible. In the optimal case in which the frozen faces are dispersed nearly homogeneously throughout finite volume, with positive probability there exists crossings of the height function which traverse in the long horizontal direction. Under less homogeneous conditions, one must consider crossing events for which the height function completely avoids collections of frozen faces altogether. If the distribution of frozen faces over a finite volume of the height function under sloped boundary conditions is significantly larger than the number of faces in finite volume that are not frozen, then there exists a slightly larger finite volume for which the number of frozen and unfrozen faces is equal, allowing for a similar use of the following arguments. 

 \bigskip

    \noindent \textit{Heuristic} (\textit{annuli crossing probabilities under flat and sloped boundary conditions are not uniformly positive}): For sloped boundary conditions $\xi$, crossings of $h$ are not uniformly positive for all annuli aspect ratios, in which the following spatially-dependent change in the following probabilities occurs, for $\mathcal{A}_{n,N} \equiv \Lambda_N \backslash \Lambda_n$,
    
    \[
\mathcal{S} \mathcal{T} \equiv \left\{\!\begin{array}{ll@{}>{{}}l}
          \textbf{P}^{\xi}_{\Lambda_n} [      \mathcal{O}_{h\geq k }(n,N)  ]    >  c      \text{ }  \text{ for }      n \leq N + 1      \text{ , } \text{ with } c > 0  \text{ , }  \\
       \textbf{P}^{\xi}_{\Lambda_n} [   \mathcal{O}_{h\geq k}(n,N)      ] \approx 0             \text{ }  \text{ for }    N + \delta > N + 1         \text{ , } \text{ with } \delta > 0 \text{ , }  \\
    \textbf{P}^{\xi}_{\Lambda_n}[    \mathcal{O}_{h \geq k}(n, N )              ]   \equiv 0      \text{ } \text{ for }  n << N   \text{ , } \\
\end{array}\right.
\]

\noindent in which the three phases above indicate different scales over which the crossing $\mathcal{C}_{h\geq k}(\mathcal{A}_{n,N}) \equiv \mathcal{O}_{h \geq k}(n,N)$ can occur. We denote the scales over which the annulus crossing probability is either strictly positive or zero, a spatial transition, abbreviated $\mathcal{S} \mathcal{T}$.

\bigskip

\noindent Besides the \textit{Heuristic} above, we further study how RSW arguments can be formulated for sloped boundary conditions, which produce six-vertex configurations in comparison to those of the flat Gibbs states which are defined below. Over finite-volume, we denote the total number of vertices of each type, $k$, as $\mathcal{T}_k$.

\bigskip

\noindent \textbf{Definition} $\textit{1}$ (\textit{flat Gibbs states in the six-vertex model}). A \textit{flat Gibbs state} over a finite-volume $\Lambda$ is a configuration sampled from the six-vertex probability measure satisfying, for $ v^k_k , v^k_{i+1} \in \mathcal{T}_k$, and the Euclidean metric $d$,

\begin{align*}
  \text{ }    \textbf{P}^{\xi_{\mathrm{Flat}}}_{\Lambda} \big[ x_{v^k_i} , x_{v^k_{i+1}}  \in V(  {\textbf{Z}^2} ) : \text{ } d( x_{v^k_i} , x_{v^k_{i+1}}   ) \text{ } \geq 2       \big]   >  0      \text{ }   \text{ , } 
\end{align*}

\noindent  where $\text{ } k \in \textbf{Z}_6 \text{ } ,\text{ }  i \in \textbf{Z}_{\geq 0}$, in which under sufficiently flat boundary conditions $\xi_{\mathrm{Flat}}$, the $i$th vertex of type $k$ is at least within Euclidean distance of two of the $(i+1)$th vertex of type $k$ over $\Lambda$, making the probability of encountering the $(i+1)$th vertex of type $k$ within distance $1$ of the $i$th vertex vanish.

\bigskip

\noindent By contrast, six-vertex configurations taken under sloped boundary conditions can be composed of multiple faces, in succession, over which are all nearest-neighbor vertices are of the same type, which raise immediate implications for the level of the height function required for crossings across symmetric domains, in addition for other components of the RSW argument.

\bigskip

\noindent \textbf{Definition} $\textit{2}$ (\textit{sloped Gibbs states in the six-vertex model}). A \textit{sloped Gibbs state} over a finite volume $\Lambda$ is a configuration sampled from the six-vertex probability measure satisfying, for $ v^k_k , v^k_{i+1} \in \mathcal{T}_k$, with respect to the Euclidean metric $d$,

\begin{align*}
     \textbf{P}^{\xi_{\mathrm{Sloped}}}_{\Lambda}\big[  x_{v^k_i} , x_{v^k_{i+1}} \in V(  {\textbf{Z}^2}  ) : \text{ }       d(    x_{v^k_i} , x_{v^k_{i+1}}    )  \geq 1    \big] \text{ } > \text{ } 0      \text{ , } \\
\end{align*}

\noindent where $\text{ } k \in \textbf{Z}_6 \text{ } ,\text{ }  i \in \textbf{Z}_{\geq 0}$, in which under sloped boundary conditions $\xi_{\mathrm{Sloped}}$, the probability of obtaining vertices of the same type is strictly positive.

\bigskip

\noindent Differences between six-vertex configurations that can be obtained from conditions placed on $d$ in \textit{Definition 1} and in \textit{Definition 2} are provided in \textit{panel four} of \textit{Figure 2}.

\section{Sloped Gibbs states for crossings across the strip}

\subsection{Overview}

\noindent Depending upon the distribution of weights, the height function can interpolate between regions of the plane over which it is increasing or decreasing. To account for the influence of sloped boundary conditions on crossings across symmetric domains, and consequently for RSW results, modifications to the boundary conditions arising from the slope of the height function at the boundary of finite volumes is first addressed for height function crossings across the strip. With negligible probability, the change of sampling a sloped Gibbs state in the strip for which there does not exist a long horizontal crossing is demonstrated to hold, as it does from arguments developed in {\color{blue}[11]} for sufficiently flat boundary conditions. Besides this commonality, one must additionally stipulate that \textit{sloped} symmetric domains in strips of the square lattice avoid, in the optimal case, all collections of frozen faces of the height function.

\bigskip

\noindent We denote $\mathcal{I}_j = [j \lfloor{\delta n}\rfloor{} , ( j + 1 ) \lfloor{\delta n }\rfloor{}] \times \{ 0 \} $, and $\widetilde{\mathcal{I}_j} = [j \lfloor{\delta n} \rfloor{} , ( j + 1 ) \lfloor{\delta n } \rfloor{ } ]  \times \{ n \} $ for $\delta , n \in \textbf{R}$ and $j \in \textbf{N}$, with the corresponding connectivity event $\widetilde{E}_{\geq k} \equiv \text{ } \{    \mathcal{I}_j            \text{ }  \overset{ h \geq k }{\longleftrightarrow} \text{ }      \widetilde{\mathcal{I}_j}      \} $. In comparison to previous notation used for sloped boundary conditions, in all forthcoming results sloped boundary conditions are denoted with $\xi_{\mathrm{Sloped}}$. We begin by addressing scaling of horizontal crossings in the strip, in which we examine the weak infinite volume limit of strip probability measures as $m \longrightarrow \infty$ under sloped boundary conditions, as $\textbf{P}^{\xi_{\mathrm{Sloped}}}_{[-m,m] \times [0,\delta n]} \longrightarrow \textbf{P}^{\xi_{\mathrm{Sloped}}}_{\textbf{Z} \times [0,\delta n]}$.

\subsection{Six-vertex objects}

\noindent \textbf{Definition} $\textit{3}$ (\textit{freezing cluster}). A \textit{freezing cluster} is collection of faces $\mathscr{F}$ over $\textbf{Z}^2$ for which each vertex is of the same type, which is given by the countable collection, for an index set $\mathcal{I}_{\Lambda}$ corresponding to the faces in the finite volume interior,

\begin{align*} 
 \mathscr{F} \mathscr{C} \text{ } = \bigcup_{ i \text{ }  \in \text{ }  \mathcal{I}_{\Lambda} } \text{ } \big\{  \forall \text{ }    \mathscr{F}_i  \in F\big(\textbf{Z}^2  \cap \Lambda  \big) \text{ , }  \exists \text{ }  j \neq i \in \mathcal{I}_{\Lambda} :    \text{ } \textbf{1}_{ \{ v_{\mathscr{F}_i} = i \} }\text{ } =  \text{ }   \textbf{1}_{\{ v_{\mathscr{F}_j} = i \} }  \big\}         \text{ } \text{ , } \text{ } 
\end{align*}

\noindent where vertices are checked for compatibility of the same type.

\bigskip

\noindent We denote the set of all such \textit{freezing clusters} as $\mathcal{FC}$.

\bigskip

\noindent \textbf{Definition} $\textit{4}$ (\textit{freezing cluster inner and outer diameters}). The \textit{inner diameter of a freezing cluster} sampled from $\mathcal{FC}$ is the maximum length for which a horizontal or vertical connectivity event within the \textit{freezing cluster} occurs between $D$ faces, of height at least $ck$,

\begin{align*}
\mathcal{D}(D,I) \text{ } \equiv \text{ }   \mathcal{D}_I  \text{ } = \text{ }  \underset{ \mathscr{F}_i \in \mathscr{F} \mathscr{C} }{\mathrm{max}}  \text{ }   \big\{     \forall \text{ } \mathscr{F}_i \text{ } \exists \text{ }   \big\{  \mathscr{F}_j \text{ } \big\}_{i+1 \leq j \leq D}         :          \text{ }  \textbf{P}_{\Lambda}^{\xi_{\mathrm{Sloped}}} [ \text{ }  \mathscr{F}_i  \underset{h \geq ck}{\overset{  \mathscr{F} \mathscr{C} }{\longleftrightarrow}}  \mathscr{F}_{D}  \text{ }  ] \text{ } > 0   \big\}               \text{ } \text{ , } \text{ } \\
\end{align*}

\noindent consisting of $|\mathcal{D}_I| = D$ faces, while the outer diameter is the maximum length for which the number of faces for which a horizontal or vertical connectivity event within the the \textit{freezing cluster} occurs between fewer than $D-l$ faces, for positive $l^{\prime} \neq i$, of height at least $ck$,

\begin{align*}
     \mathcal{D} (D-L , O)     \text{ }  \equiv \text{ }   \mathcal{D}_O \text{ } = \text{ }  \underset{ \mathscr{F}_{l^{\prime}} \in \mathscr{F} \mathscr{C} }{\mathrm{max}}  \text{ }    \big\{      \forall \text{ } \mathscr{F}_{l^{\prime}} \text{ } \exists \text{ }   \big\{  \mathscr{F}_j \text{ } \big\}_{l^{\prime}+1 \leq j < D - l }         :        \text{ }  \textbf{P}_{\Lambda}^{\xi_{\mathrm{Sloped}}} [ \text{ }  \mathscr{F}_{l^{\prime}}  \underset{h\geq ck}{\overset{  \mathscr{F} \mathscr{C} }{\longleftrightarrow}}  \mathscr{F}_{D-l-1}  \text{ }  ] \text{ } > 0  \big\}              \text{ } \text{ , } \text{ } \\
\end{align*}

\noindent consisting of $| \mathcal{D}_O| = D-l$ faces. See \textit{Figure 3} for an illustration.

\bigskip

\noindent \textbf{Definition} $\textit{5}$ (\textit{sloped-boundary domains}). A \textit{domain} $\Gamma$ is a face-contour satisfying, for $n^{\prime} > 1$ and $\mathcal{I}_L , \widetilde{\mathcal{I}_L}, \mathcal{I}_R , \widetilde{\mathcal{I}_R} \subset [-m,m] \times [0,n^{\prime} N ]$,

\[
F(\text{ } \textbf{Z}^2 \text{ } )\text{ }  \supsetneq    \text{ } \Gamma \equiv \Gamma\big(\text{ } \gamma_L , \gamma_R , [a_L , \widetilde{a_L}] , [a_R , \widetilde{a_R} ]  \text{ } \big)  \text{ } \Longleftrightarrow \text{ } 
\left\{\!\begin{array}{ll@{}>{{}}l} \gamma_L \subsetneq F( \text{ } \textbf{Z}^2 \text{ } )  & \text{ : }  & \text{ } \textbf{P}^{\xi_{\mathrm{Sloped}}}_{[-m,m] \times [0,n^{\prime} N]} [ \text{ }   \mathcal{I}_L \overset{  h \geq c k}{\longleftrightarrow}  \widetilde{\mathcal{I}_L} \text{ } ] \text{ }  >  \text{ } 0  \text{ } \text{ , } \\
 \gamma_R \subsetneq F( \text{ } \textbf{Z}^2 \text{ } )  & \text{ : }  & \text{ } \textbf{P}^{\xi_{\mathrm{Sloped}}}_{[-m,m] \times [0,n^{\prime} N]} [ \text{ }   \mathcal{I}_R \overset{  h \geq c k}{\longleftrightarrow}  \widetilde{\mathcal{I}_R} \text{ } ] \text{ }  >  \text{ } 0  \text{ } \text{ , } \\
\end{array}\right.
\]

\noindent where $\gamma_L , \gamma_R$ denote the left and right boundaries comprising two sides of $\Gamma$ which are also face-contours, and $a_L \in \mathcal{I}_L$, $\widetilde{a_L} \in  \widetilde{\mathcal{I}_L}$, $a_R \in \mathcal{I}_R$, $\widetilde{a_R} \in \widetilde{\mathcal{I}_R}$ denote points demarcating the top and bottom sides of the domain along the boundary of the strip.

 \subsection{RSW scaling for frozen square sublattices}

To obtain desired estimates for crossing probabilities in the strip, we characterize properties of strip symmetric domains. If the boundary conditions of the height function along finite volume are sloped, then, with positive probability, if the aspect dimensions of the supporting finite volume on the probability measure are taken to be too narrow, the probability of a long vertical crossing vanishes. To obtain such desired crossing estimates for horizontal traversals in infinite volume, we introduce properties of the strip environment, ranging from observations that long horizontal crossings in the strip occur with sufficiently good probability in infinite volume if, with positive probability, there exists paths which can avoid at least a positive proportion of the faces that are frozen so that long horizontal crossings occur with sufficient probability; in finite volume, the crossing probability of a horizontal crossing occurring between two segments of the lower boundary of the strip occurs with negligible probability; the occurrence of countably many collections of frozen faces in the strip does not prevent long horizontal crossings from occurring with sufficiently good probability.

\bigskip

    \noindent \textbf{Lemma} \textit{2.1} (\textit{Scaling of the width of the square lattice strip for long horizontal crossings of the height-function under sloped boundary conditions}): For any \textit{freezing cluster}, there exists $n \in\textbf{R}_{>1}$ such that the probability of $\widetilde{E}_{\ngeq k}$ occurring is proportional to the conditional probability of a crossing across height $k$ not occurring, where any one of the \textit{three} conditional connectivity events over $\textbf{Z}^2$ holds:

\begin{itemize}
    \item[$\bullet$] \textit{Blocking inner diameter of the freezing cluster}

    The conditional crossing event is between the boundary of the inner diameter, one face of the \textit{freezing cluster} and the boundary of the inner diameter, in which, 
    
      \begin{align*}
 \textbf{P}^{\xi_{\mathrm{Sloped}}}_{\textbf{Z} \times [0, n N]} \bigg[   \widetilde{E}_{\ngeq k} \text{ }  | \text{ }                 \big\{ \text{ }   \exists \text{ }   \mathscr{F}\mathscr{C}  \text{ } \in \mathcal{FC}  \text{ } \text{ with } | \mathcal{D}_{I}  |\text{ } < \infty \text{ }  \big\} \text{ } \cap  \big\{ \text{ }     \mathcal{I}_j  \overset{h \geq ck}{\not\longleftrightarrow} \text{ }     \partial \mathcal{D}_I \text{ } \big\}  \cap \big\{ \text{ }  \mathscr{F} \mathscr{C} \overset{h \geq ck}{\not\longleftrightarrow} \text{ }   \mathcal{I}_j \text{ }  \big\}  \bigg]     \text{ } \text{ . } \text{ } \\
    \end{align*}

        \item[$\bullet$] \textit{Blocking outer diameter of the freezing cluster}
        
          The conditional crossing event is between the boundary of the outer diameter, one face of the \textit{freezing cluster} and the boundary of the outer diameter, in which, 
        
        \begin{align*}
          \textbf{P}^{\xi_{\mathrm{Sloped}}}_{\textbf{Z} \times [0, n N]} \bigg[   \widetilde{E}_{\ngeq k} \text{ }  | \text{ }             \big\{ \text{ }        \exists \text{ }   \mathscr{F}\mathscr{C}  \text{ } \in \mathcal{FC}  \text{ } \text{ with }  |\mathcal{D}_{O} |\text{ } < \infty \text{ }   \big\} \cap \big\{ \text{ }  \mathcal{I}_j  \overset{h\geq ck}{\not\longleftrightarrow} \text{ }     \partial \mathcal{D}_O \text{ } \big\} \cap \big\{ \text{ }  \mathscr{F} \mathscr{C} \overset{h\geq ck}{\not\longleftrightarrow} \text{ }    \mathcal{I}_j  \text{ } \big\}    \bigg]         \text{ } \text{ . } \text{ }   \\
        \end{align*}
            \item[$\bullet$] \textit{Simultaneous blocking of inner and outer diameters of the freezing cluster}
            
                The conditional crossing event is between the boundary of the outer diameter, one face of the \textit{freezing cluster} and the boundary of the outer diameter, in which, 
        
        \begin{align*}
        \textbf{P}^{\xi_{\mathrm{Sloped}}}_{\textbf{Z} \times [0, n N]} \bigg[ \text{ }  \widetilde{E}_{\ngeq k} \text{ }  | \text{ }              \big\{ \text{ }      \exists \text{ }   \mathscr{F}\mathscr{C}  \text{ } \in \mathcal{FC}  \text{ } \text{ with } | \mathcal{D}_I \cup  \mathcal{D}_{O} | \text{ } < \infty \text{ } \big\} \cap \big\{ \text{ }    \mathcal{I}_j  \overset{h\geq ck}{\not\longleftrightarrow} \text{ }   \big(  \partial \mathcal{D}_O \text{ } \cup \text{ }  \partial \mathcal{D}_I  \big) \text{ } \big\} \cap \big\{ \text{ }   \mathscr{F} \mathscr{C} \overset{h\geq ck}{\not\longleftrightarrow} \text{ }    \mathcal{I}_j \text{ } \big\} \bigg]    \text{ . } \\
        \end{align*}

\end{itemize}

\noindent From configurations drawn from the six-vertex sample space, it is possible to encounter any one of the three possible configurations above with positive probability, namely the first configuration in which, given the failure of the height function to maintain a level $\geq k$ for all faces in the strip, the existence of a \textit{freezing cluster} for which $\textbf{P}^{\xi_{\mathrm{Sloped}}}_{\textbf{Z} \times [ 0 , nN]} \big[     \mathcal{I}_j \overset{h \geq ck}{\longleftrightarrow}  \partial \mathcal{D}_I     \big] =0$ and $\textbf{P}^{\xi_{\mathrm{Sloped}}}_{\textbf{Z} \times [ 0 , nN]} \big[      \mathscr{F}\mathscr{C} \overset{ h \geq ck}{\longleftrightarrow}    \mathcal{I}_j \big] =0$, consisting of finitely many frozen faces; the second configuration in which, given the failure of the height function to maintain a level $\geq c k$ for all faces in the strip, the existence of a \textit{freezing cluster} for which $\textbf{P}^{\xi_{\mathrm{Sloped}}}_{\textbf{Z} \times [ 0 , nN]} \big[     \mathcal{I}_j   \overset{h \geq ck}{\longleftrightarrow} \partial \mathcal{D}_O    \big] =0$ and $\textbf{P}^{\xi_{\mathrm{Sloped}}}_{\textbf{Z} \times [ 0 , nN]} \big[      \mathscr{F} \mathscr{C}          \overset{h \geq ck}{\longleftrightarrow}     \mathcal{I}_j            \big] =0$, consisting of finitely many frozen faces; the third configuration in which, given the failure of the height function to maintain a level $\geq k$ for all faces in the strip, the existence of a \textit{freezing cluster} for which $\textbf{P}^{\xi_{\mathrm{Sloped}}}_{\textbf{Z} \times [ 0 , nN]} \big[       \big\{  \mathcal{I}_j          \overset{h \geq ck}{\longleftrightarrow}          \partial \mathcal{D}_O  \big\}  , \big\{ \mathcal{I}_j \overset{h \geq ck}{\longleftrightarrow} \partial \mathcal{D}_I \big\}  \big] =0$ and $\textbf{P}^{\xi_{\mathrm{Sloped}}}_{\textbf{Z} \times [ 0 , nN]} \big[        \mathscr{F} \mathscr{C}        \overset{h \geq ck}{\longleftrightarrow}      \mathcal{I}_j           \big] =0$, consisting of finitely many frozen faces.

\bigskip

\noindent \textit{Proof of Lemma 2.1}. We characterize the crossing probability in the first case with the \textit{blocking inner diameter of the freezing cluster}. As a matter of notation, we write

\begin{align*}
  \partial \big([-m,m] \times [0,nN]\big) \text{ } = \big( [-m,m] \times  \{ 0 \} \big) \cup \big( [-m,m] \times \{ n N \} \big)  \text{ } \text{ . } \text{ } \\
\end{align*}

\noindent The conditional probability under sloped boundary conditions is equivalent to the product of crossing probabilities, in finite volume,

\begin{align*}
   \textbf{P}^{\xi_{\mathrm{Sloped}}}_{[-m,m] \times [0,nN]} [ \text{ } \widetilde{E}_{\ngeq k} \text{ } ] \text{ }     \textbf{P}^{\xi_{\mathrm{Sloped}}}_{[-m,m] \times [0, n N]} \bigg[              \big\{    \exists \text{ }   \mathscr{F}\mathscr{C}  \text{ } \in \mathcal{FC}  \text{ } \text{ with } | \mathcal{D}_{I}  |\text{ } < \infty \big\} \text{ } \cap  \big\{     \mathcal{I}_j  \overset{h\geq ck}{\not\longleftrightarrow}      \partial \mathcal{D}_I \text{ } \big\}  \cap \big\{ \mathscr{F} \mathscr{C} \overset{h\geq ck}{\not\longleftrightarrow}   \mathcal{I}_j  \big\}  \bigg]          \text{ }     \text{ . } \\
\end{align*}
  
    \noindent To lower bound the second term of the above product, observe, by (FKG), 
    
    \begin{align*}
    \text{ }    \textbf{P}^{\xi_{\mathrm{Sloped}}}_{[-m,m] \times [0, n N]} [  \text{ }       \exists \text{ }   \mathscr{F}\mathscr{C}   \in \mathcal{FC}  \text{ with } | \mathcal{D}_{I}  |\text{ } < \infty     \text{ } ] \text{ }   \textbf{P}^{\xi_{\mathrm{Sloped}}}_{[-m,m] \times [0, n N]} [     \mathcal{I}_j  \overset{h\geq ck}{\not\longleftrightarrow} \text{ }     \partial \mathcal{D}_I ] \text{ }  \textbf{P}^{\xi_{\mathrm{Sloped}}}_{[-m,m] \times [0, n N]} [  \text{ }       \mathscr{F} \mathscr{C} \overset{h\geq ck}{\not\longleftrightarrow} \text{ }   \mathcal{I}_j \text{ }       \text{ } ]        \text{ }   \text{ , } \text{ }    \\
    \end{align*}    
    
    \noindent from which the disconnective probabilities of $ \{  \mathcal{I}_j  \overset{h\geq ck}{\not\longleftrightarrow} \text{ }     \partial \mathcal{D}_I \}$, and of $ \{ \mathscr{F} \mathscr{C} \overset{h \geq ck}{\not\longleftrightarrow} \text{ }   \mathcal{I}_j \} $, occurring can each respectively be expressed with,

    \begin{align*}
     \text{ }     \textbf{P}^{\xi_{\mathrm{Sloped}}}_{[-m,m] \times [0, n N]} [       \mathcal{I}_j   \overset{h< ck}{\longleftrightarrow}    \partial \mathcal{D}_I           ]  \text{ } \text{ , } \text{ } \text{ and }   \text{ } \textbf{P}^{\xi_{\mathrm{Sloped}}}_{[-m,m] \times [0, n N]} [  \mathscr{F} \mathscr{C}   \overset{h< ck}{\longleftrightarrow}     \mathcal{I}_j          ]      \text{ } \text{ , } \text{ } \\
    \end{align*}  
   
   \noindent from properties of x-crossings introduced in {\color{blue}[11]} for crossings of the height-function in which the distance between neighboring faces can be two. For the first probability, under $\xi_{\mathrm{Sloped}}$ and $\chi_{\mathrm{Sloped}} \neq \xi_{\mathrm{Sloped}}$, there exists a positive integer $n^{\prime} > n \geq 1$ for which, 
   
   \begin{align*}
    [ - m , m ] \times [0 , n N ] \subset [ -m , m ] \times [0 , n^{\prime} N ] \text{ } \text{ , }   \text{ } \\
   \end{align*}
   
   \noindent implying, by (SMP), 
   
        \begin{align*}
         \textbf{P}^{\xi_{\mathrm{Sloped}}}_{[-m,m] \times [0,n^{\prime}N]}   [    \mathcal{I}_j   \overset{h< ck}{\longleftrightarrow}    \partial \mathcal{D}_I    |       \mathcal{C}           ] \text{ }  \equiv \text{ }   \textbf{P}^{\chi_{\mathrm{Sloped}}|_{\partial ([-m,m] \times [0,nN])}}_{[-m,m] \times [0,nN]}  [      \mathcal{I}_j   \overset{h< ck}{\longleftrightarrow}    \partial \mathcal{D}_I         ]    \text{ }  \text{ , } \text{ } \\
    \end{align*} 
   
         \begin{figure}
\begin{align*}
\includegraphics[width=1.03\columnwidth]{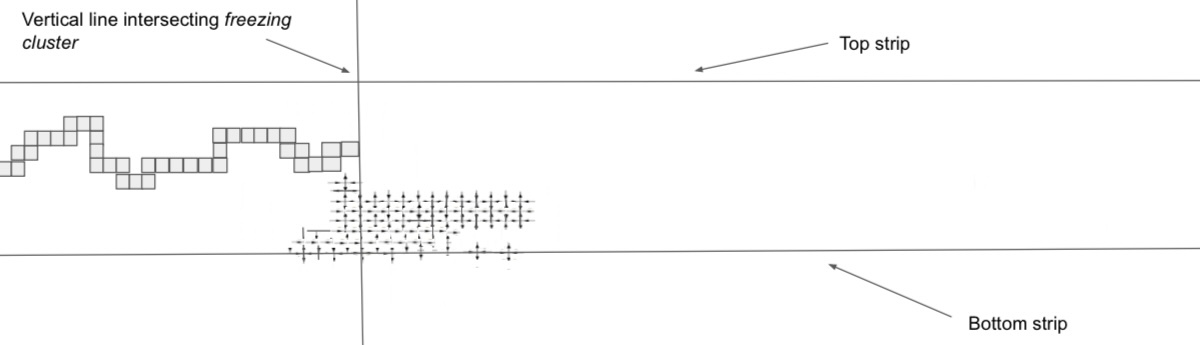}\\
\end{align*}
\caption{\textit{One realization of a crossing event of height $<k$ occurring in the strip, conditionally upon the existence of a freezing cluster composed of finitely many faces intersecting the bottom boundary of the strip}. One vertical line intersecting the freezing cluster, in addition to the final face of the path, is also depicted.}
\end{figure}

   \noindent where the conditioning is over a sub-lattice of the strip given by the union below,

   \begin{align*}
   \mathcal{C} =  \big\{  h =   \text{ }    \chi_{\mathrm{Sloped}} \text{ on } \big([ -m , m ] \times [n N , n ^{\prime} N ] \big) \cup  \text{ }      \big( \text{ }    \partial ([-m,m] \times [0,nN])     \text{ }  \big)              \text{ } \big\}  \text{ , } \text{ }  \\
   \end{align*}

   \noindent in addition to the existence of $  \xi^{\prime}_{\mathrm{Sloped}} \geq \xi_{\mathrm{Sloped}}$ for which, by (CBC),

        \begin{align*}
      \textbf{P}^{\xi_{\mathrm{Sloped}}}_{[-m,m] \times [0, n N]}  [    \mathcal{I}_j   \overset{h< ck}{\longleftrightarrow}    \partial \mathcal{D}_I        ]  \leq  \textbf{P}^{\xi^{\prime}_{\mathrm{Sloped}}}_{[-m,m] \times [0, n N]}  [      \mathcal{I}_j   \overset{h< ck}{\longleftrightarrow}    \partial \mathcal{D}_I     ]   \text{ }  \text{ . } \text{ } \\
    \end{align*} 
   
   \noindent We relate these spatial properties by taking $m \longrightarrow \infty$ and varying the slope of the boundary conditions, from which we obtain the lower bound estimate for the conditional crossing event,
   
   \begin{align*}
          \text{ }       \textbf{P}^{\xi_{\mathrm{Sloped}}}_{[-m,m] \times [0,n^{\prime}N]}   [      \mathcal{I}_j   \overset{h< ck}{\longleftrightarrow}    \partial \mathcal{D}_I    ] \text{ }          = \text{ }        \frac{\textbf{P}^{\chi_{\mathrm{Sloped}}|_{\partial ([-m,m] \times [0,nN])}}_{[-m,m] \times [0,nN]}  [       \mathcal{I}_j   \overset{h< ck}{\longleftrightarrow}    \partial \mathcal{D}_I          ]     }{\textbf{P}^{\xi_{\mathrm{Sloped}}}_{[-m,mn] \times [0,n^{\prime} N]}[ \mathcal{C}  ] } \text{ }   \\ \text{ }  > \text{ }   \mathcal{C}^{-1}_{\xi_{\mathrm{Sloped}}}         \text{ }    \textbf{P}^{ \chi_{\mathrm{Sloped}}|_{\partial ([-m,m] \times [0,nN])}}_{[-m,m] \times [0,nN]}  [      \mathcal{I}_j   \overset{h< ck}{\longleftrightarrow}    \partial \mathcal{D}_I       ]         \text{ } \text{ , } \text{ }   \\
   \end{align*}

   \noindent where the probability in the uppermost bound is dependent upon $h$ achieving the boundary conditions stipulated by $\xi_{\mathrm{Sloped}}$, which occurs with positive probability because over all faces $\underset{\mathscr{F}_i \neq \mathscr{F}_j \in F(\textbf{Z}^2)}{\mathrm{max}} \big\{ \text{ }    h_{\mathscr{F}_i}^{\chi_{\mathrm{Sloped}}} - h_{\mathscr{F}_j}^{\xi_{\mathrm{Sloped}}}             \text{ }   \big\} \text{ } < \pm  \text{ }  \infty$; hence this term can be bounded below with strictly positive $\mathcal{C}_{\xi_{\mathrm{Sloped}}}$.

   \bigskip
   
   \noindent For the numerator term, under $\chi_{\mathrm{Sloped}}|_{\partial ([-m,m] \times [0,nN])}$, 
   
   \begin{align*}
          \textbf{P}^{\chi_{\mathrm{Sloped}}|_{\partial ([-m,m] \times [0,nN])}}_{[-m,m] \times [0,nN]}  [          \mathcal{I}_j   \overset{h< ck}{\longleftrightarrow}    \partial \mathcal{D}_I    ] \text{ } \geq \text{ }     \textbf{P}^{\chi_{\mathrm{Sloped}}|_{\partial ([-m,m] \times [0,nN])}}_{[-m,m] \times [0,nN]}  \bigg[   \text{ } \underset{\forall i^{\prime} ,  e \in E(\textbf{Z}^2) \exists \# \{ \mathscr{F} \in F(\textbf{Z}^2) \}    = n  \text{ } : \text{ }    \mathscr{F} \cap \mathcal{I}_j = e              }{\underset{i^{\prime}=j+1}{\overset{n}{\bigcap}}} \text{ } \big\{\text{ }  \mathcal{I}_j   \overset{h< ck}{\longleftrightarrow}  \mathcal{I}_{i^{\prime}} \text{ } \big\}  \text{ } \bigg] \text{ } \\ \text{ } \overset{(\mathrm{FKG})}{\geq} \text{ } \prod_{i^{\prime}=1}^n \text{ }   \textbf{P}^{\chi_{\mathrm{Sloped}}|_{\partial ([-m,m] \times [0,nN])}}_{[-m,m] \times [0,nN]}  \big[   \text{ }                \mathcal{I}_j   \overset{h< ck}{\longleftrightarrow}  \mathcal{I}_{i^{\prime}}   \text{ } \big] \\ \text{ }           > \text{ }      \prod_{i^{\prime}=1}^n \text{ }   \textbf{P}^{\mathrm{min} \text{ } \chi_{\mathrm{Sloped}}|_{\partial ([-m,m] \times [0,nN])}}_{[-m,m] \times [0,nN]}  \big[  \text{ }                        \mathcal{I}_j   \overset{h< ck}{\longleftrightarrow}  \mathcal{I}_{i^{\prime}}                            \text{ } \big]  \\        \text{ }   \geq \text{ }      c_{(\mathrm{min} \text{ } \chi_{\mathrm{Sloped}})}^n        \text{ }  \text{ , } \text{ } \\
   \end{align*}

   \noindent is the infimum over segment connectivity events dependent upon $i^{\prime}$. We characterize the sloped boundary dependence on the conditional event $\mathscr{E}$ stated in the \textit{Lemma} by collecting previous estimates, yielding
   
      \begin{figure}
\begin{align*}
\includegraphics[width=0.8\columnwidth]{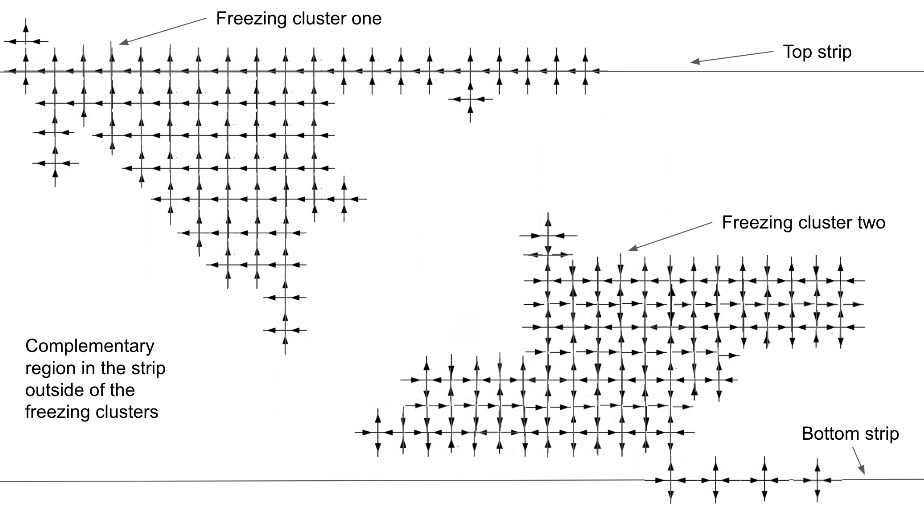}\\
\end{align*}
\caption{\textit{Representative freezing clusters which macroscopically penetrate the finite volume bulk from the boundary of the strip}. Under sloped boundary conditions, the height function can become frozen, inducing a cluster of vertices that are either all of the same type, or are of alternating type between pairs of vertices (\textit{Freezing cluster one} is a six-vertex configuration sample for which the height function is frozen on vertices of type $3$, while \textit{Freezing cluster two} is a six-vertex configuration sample for which the height function is frozen on a combination of vertices of type $5$ and $6$). In the complementary region of the strip outside of the \textit{freezing clusters} on the upper and lower portion of the strip, standard RSW arguments for flat Gibbs remain applicable.}
\end{figure}

   \noindent where in the first inequality, a lower bound is provided for $\mathcal{I}_{i^{\prime}} \subset [-m,m] \times [0,nN]$ with $\mathcal{I}_j \cap \mathcal{I}_{i^{\prime}} \neq \emptyset$ for $i^{\prime}=j+1$ and $\mathcal{I}_{j+2} \cap \mathcal{I}_{j+3} \neq \emptyset$ for each $i^{\prime} > j + 1$; in the second inequality, a lower bound is provided by (FKG) for crossings across the set of $\mathcal{I}_{i^{\prime}}$; in the third inequality, a lower bound is provided from minimal boundary conditions $\mathrm{min} \text{ } \xi_{\mathrm{Sloped}}|_{\partial ([-m,m] \times [0,nN])}$; in the fourth inequality, there exists a uniformly, strictly positive constant $ c_{(\mathrm{min} \text{ } \chi_{\mathrm{Sloped}})}^n$ for the connectivity event under minimal sloped boundary conditions occurs, where,

   \begin{align*}
c_{(\mathrm{min} \text{ } \chi_{\mathrm{Sloped}})}    \text{ } = \text{ } \underset{i^{\prime} \in \textbf{Z}}{\mathrm{inf}} \text{ } \big\{  \textbf{P}^{\mathrm{min} \text{ } \chi_{\mathrm{Sloped}}|_{\partial ([-m,m] \times [0,nN])}}_{[-m,m] \times [0,nN]}  \big[                        \mathcal{I}_j   \overset{h< ck}{\longleftrightarrow}  \mathcal{I}_{i^{\prime}}                      \big] \big\}  \text{ } \in ( 0 , 1 ]   \text{ , } \text{ }     \\
   \end{align*}
   
\noindent which appears as a prefactor to the probability in the lower bound,

   \begin{align*}
            \textbf{P}^{\xi_{\mathrm{Sloped}}}_{[-m,m] \times [0,nN]} \text{ } \big[     \mathscr{E}  \big] \text{ }       > \text{ }  \frac{c^n_{(\mathrm{min} \text{ } \chi_{\mathrm{Sloped}})} }{\mathcal{C}_{\xi_{\mathrm{Sloped}}}} \text{ }   \text{ }      \textbf{P}^{\xi_{\mathrm{Sloped}}}_{[-m,m] \times [0,nN]} \big[  \widetilde{E}_{\ngeq k}  \big] \text{ }     \textbf{P}^{\xi_{\mathrm{Sloped}}}_{[-m,m] \times [0, n N]} \big[             \big\{ \text{ }   \exists    \mathscr{F}\mathscr{C}   \in \mathcal{FC}   \text{ with } | \mathcal{D}_{I}  | < \infty  \big\} \big]     \text{ } \text{ , } \text{ } \\
   \end{align*}
   
   \noindent which reads, in infinite volume, 
   
     \begin{align*}
            \textbf{P}^{\xi_{\mathrm{Sloped}}}_{\textbf{Z} \times [0,nN]} \text{ } \big[     \mathscr{E}  \big] \text{ }       > \text{ }  \frac{c^n_{(\mathrm{min} \text{ } \chi_{\mathrm{Sloped}})} }{\mathcal{C}_{\xi_{\mathrm{Sloped}}}} \text{ }   \text{ }      \textbf{P}^{\xi_{\mathrm{Sloped}}}_{\textbf{Z} \times [0,nN]} \big[  \widetilde{E}_{\ngeq k} \big] \text{ }     \textbf{P}^{\xi_{\mathrm{Sloped}}}_{\textbf{Z} \times [0, n N]} \big[                  \big\{   \exists    \mathscr{F}\mathscr{C}  \in \mathcal{FC}  \text{ } \text{ with } | \mathcal{D}_{I}  | < \infty   \big\} \big]     \text{ } \text{ , } \text{ } \\
   \end{align*}
   
   \noindent To conclude the argument, observe that the second probability in the inequality above occurs with positive probability as the number of frozen faces of the height function in a \textit{freezing cluster} over $[-m,m] \times [0,nN]$ is $8 m n N$, and that the last probability can be similarly bounded from below by following the same argument as given for the previous two crossing events (ie, express the disconnectivity event in terms of a connectivity event through complementary x-crossings, and apply (FKG) to $n^{\prime}$ crossings). Altogether,

   \begin{align*}
       \frac{c^n_{(\mathrm{min} \text{ } \chi_{\mathrm{Sloped}})} \text{ } (c^{\prime})^{n^{\prime}} c_{3}}{\mathcal{C}_{\xi_{\mathrm{Sloped}}}} \text{ } >  c_0 > 0 \text{ , } \\
   \end{align*}
   
   \noindent for suitable constants, $(c^{\prime})^{n^{\prime}}$, $c^n_{(\mathrm{min} \text{ } \chi_{\mathrm{Sloped}}}$, and,  
   
   \begin{align*}
   c_3 \text{ }  \equiv \text{ } c_3 \big(   |\mathcal{D}_I|    \big)   \text{ }  \text{ , } \text{ } \\
   \end{align*}
   
   \noindent given $|\mathcal{D}_I| < \infty$ and arbitrary $c_0$. 
  
   \bigskip
   
   \noindent The remaining cases follow by executing the same argument, in which the conditioning is dependent upon the boundary faces incident to $\partial \mathcal{D}_0$, or incident to $\partial \mathcal{D}_O \cup \partial \mathcal{D}_I$. Hence there exists such an $n^{\prime}$ for each case. \boxed{} 
   
   \bigskip
   
   \noindent To proceed, we analyze how placement of the left and right boundaries of the \textit{symmetric domain} from strip crossings placed relative to the occurrence of \textit{freezing clusters} provides  the following upper bound.
   
   \bigskip
    
    \noindent \textbf{Proposition} $\textit{1.2}$ (\textit{sloped-boundary crossing event upper bound}): Fix suitable $\delta , \delta^{\prime\prime}, n >0$ such that $\lfloor{}\delta n\rfloor{}> \lfloor{}\delta^{\prime\prime} n\rfloor{} \in \textbf{Z}$. For any $k \geq \frac{1}{c}$ and $i  \in \textbf{Z}$, the crossing probability in the strip, 
   
    \begin{align*}
          \textbf{P}^{\xi_{\mathrm{Sloped}}}_{\textbf{Z} \times [0,n^{\prime} N]}  \big[    [0 , \delta^{\prime\prime} n]                             \text{ }   \times \{ 0 \}      \text{ } \overset{ h \geq (1-c) k }{ \longleftrightarrow } \text{ } [ i , i + \delta^{\prime\prime} n ]      \text{ }  \times \{ n \}                      \big] \text{ } < \text{ } 1 - c  \text{ } \text{ , } \text{ }   \\
    \end{align*}

    \noindent admits an upper bound dependent on $c$.
    
    \bigskip
    
    \noindent The boundary conditions along the strip for the study of \textit{sloped Gibbs states }completely differ from those imposed for \textit{flat Gibbs states} in {\color{blue}[11]}. In the following arguments, the changes primarily emerge from the fact that long horizontal crossings in infinite volume under sloped boundary conditions

 \bigskip
 
 \noindent  \textit{Proof of Proposition 1.2}. We establish the following.
 \bigskip
    
    \noindent \textbf{Lemma} $\textit{2.2}$ (\textit{maintaining sufficient distance between left and right boundaries of the domain with freezing clusters as the slope of boundary conditions increases}): Suppose that there is a countable number of freezing clusters $N$ distributed over the finite-volume strip boundary, in which,
    
    \begin{align*}
 \big\{   \forall \text{ } 1 \leq i_{\mathscr{F} \mathscr{C}} \equiv i  \leq \text{ } N , \text{ } \exists \text{ }  \mathscr{F} \mathscr{C}_i  \in \mathcal{F} \mathcal{C} \text{ } : \text{ } \mathscr{F} \mathscr{C}_i \cap  \big(   \partial \big([-m,m] \times [0,nN]\big)   \big) \text{ } \neq \emptyset \big\}   \text{ } \text{ . } \text{ } \\
    \end{align*}
    
   \noindent  For $v \in V(\textbf{Z}^2)$, $e \in E(\textbf{Z}^2)$, and $\mathscr{F}_v  \in F( \text{ } \textbf{Z}^2 \text{ } ) $ with $\mathscr{F}_v \cap \gamma_L = e$, the difference is bounded below by,

    \begin{align*}
         \text{ }  \bigg|  \textbf{P}^{\xi^2_{\mathrm{Sloped}}}_{[-m,m] \times [0,n^{\prime}N]}     \big[ \gamma_L \overset{h \geq ck}{\longleftrightarrow} \mathscr{F}\mathscr{C}_1    \big]           \text{ } -  \text{ }     \textbf{P}^{\xi^1_{\mathrm{Sloped}}}_{[-m,m] \times [0,n^{\prime}N]}  \big[   \gamma_L  \overset{h \geq ck}{\longleftrightarrow}    \mathscr{F}\mathscr{C}_1      \big]     \bigg|         \text{ } \geq \text{ }               { \bigg|  F (  \gamma_L     \cap   \mathcal{L}^1_{\mathscr{F}\mathscr{C}_1}   ) \bigg|  \bigg[ \frac{\delta^2_L}{\delta^1_L } \bigg]   }  \text{ } \text{ , } \text{ } \\ \tag{\textit{2.2 left boundary symmetric domain lower bound}}
    \end{align*}

    \noindent between crossing events respectively taken under boundary conditions $\xi^2_{\mathrm{Sloped}} \equiv \xi^1_{\mathrm{Sloped}} + \text{ }    \tau_k         \text{ }  \geq \xi^1_{\mathrm{Sloped}} \text{ }$ for arbitrary $0 < \delta^1_L \text{ }  \neq \text{ } \delta^2_L \neq 0$, $\mathcal{L}^1_{\mathscr{F} \mathscr{C}}$ is a line intersecting the strip at an arbitrary distance to the left of $\mathscr{F} \mathscr{C}_1$, and $\tau_k$ denotes a horizontal translation of boundary conditions rightwards along the boundary of the strip by $k$ faces.
    
    \bigskip
    
    \noindent Similarly, given some $\mathscr{F}\mathscr{C}_N$ for which the below infimum is achieved,

    \begin{align*}
  \mathscr{F}\mathscr{C}_N \text{ } \equiv \text{ }    \underset{i \in \textbf{Z}}{\mathrm{inf}} \text{ } \big\{   d  (     \mathscr{F}\mathscr{C}_i     ,      \gamma_R   )        \big\}           \text{ , } 
    \end{align*}
    
    \noindent one similarly has,
    
    \begin{align*}
\text{ }   \bigg|     \textbf{P}^{(\xi^2_{\mathrm{Sloped}})^{\prime}}_{[-m,m] \times [0,n^{\prime}N]}     \big[ \mathscr{F}\mathscr{C}_N   \overset{h \geq ck}{\longleftrightarrow}      \gamma_R  \big]           \text{ } -  \text{ }     \textbf{P}^{(\xi^1_{\mathrm{Sloped}})^{\prime}}_{[-m,m] \times [0,n^{\prime}N]}  \big[    \mathscr{F}\mathscr{C}_N \overset{h \geq ck}{\longleftrightarrow}      \gamma_R              \big]      \bigg|   \geq          { \text{ } \bigg|  F (  \gamma_R    \cap   \mathcal{L}^2_{\mathscr{F}\mathscr{C}_N}   ) \bigg| \bigg[ \frac{\delta^2_R}{\delta^1_R} \bigg]  }       \text{ } \text{ , } \text{ } \\
    \end{align*}

    \noindent between crossing events respectively taken under boundary conditions $(\xi^2_{\mathrm{Sloped}})^{\prime} \equiv (\xi^1_{\mathrm{Sloped}})^{\prime} - \text{ }    \tau_k         \text{ }  \geq (\xi^1_{\mathrm{Sloped}})^{\prime} \text{ }$ for $0 < \delta^1_R \neq \delta^2_R \neq 0$, and $\mathcal{L}^2_{\mathscr{F}\mathscr{C}_N}$ is a line intersecting the strip at an arbitrary distance to the right of $\mathscr{F}\mathscr{C}_N$.
    
\bigskip

\noindent \textit{Proof of Lemma 2.2}. Performing straightforward rearrangement to obtain the first relation implies,

\begin{align*}
\text{ }  \bigg|    \textbf{P}^{\xi^2_{\mathrm{Sloped}}}_{[-m,m] \times [0,n^{\prime}N]}     \big[   \gamma_L \overset{h \geq ck}{\longleftrightarrow} \mathscr{F}\mathscr{C}     \big]           \text{ } -  \text{ }     \textbf{P}^{\xi^1_{\mathrm{Sloped}}}_{[-m,m] \times [0,n^{\prime}N]}  \big[   \gamma_L \overset{h \geq ck}{\longleftrightarrow} \mathscr{F}\mathscr{C}        \big]  \bigg|             \text{ } \\ \Updownarrow  \text{                } \text{ (SMP) applied to the second term} \\  \text{ }        \bigg|       \textbf{P}^{\xi^2_{\mathrm{Sloped}}}_{[-m,m] \times [0,n^{\prime}N]}     \big[  \gamma_L \overset{h \geq ck}{\longleftrightarrow} \mathscr{F}\mathscr{C}     \big]          -   \textbf{P}^{\xi^2_{\mathrm{Sloped}}}_{[-m,m] \times [0,n^{\prime}N]}     \big[   \gamma_L \overset{h \geq ck}{\longleftrightarrow} \mathscr{F}\mathscr{C}   \big] \text{ }  \\ \times    \textbf{P}^{\xi^2_{\mathrm{Sloped}}}_{[-m,m] \times [0,n^{\prime}N]}     \big[     h  = \text{ } \xi^1_{\mathrm{Sloped}} \text{ }                      \text{on }  \partial \big(  [-m,m] \times  [n N , n^{\prime} N ]\big)  \big] \bigg| \text{ } \\ \text{ } \Updownarrow \text{ } \\ \text{ }  \bigg| \text{ } \textbf{P}^{\xi^2_{\mathrm{Sloped}}}_{[-m,m] \times [0,n^{\prime}N]}     \big[  \gamma_L \overset{h \geq ck}{\longleftrightarrow} \mathscr{F}\mathscr{C}     \big]  \bigg[     1 -    \textbf{P}^{\xi^2_{\mathrm{Sloped}}}_{[-m,m] \times [0,n^{\prime}N]}     \big[    h  = \text{ } \xi^1_{\mathrm{Sloped}} \text{ }                      \text{on }  \partial \big( [-m,m] \times  [n N , n^{\prime} N ]  \big)   \big]   \bigg]  \bigg| \text{ }     \text{ . }       \text{ } 
\end{align*}

    \noindent As a result of the existence of $\xi^{3}_{\mathrm{Sloped}} \geq \xi^{2}_{\mathrm{Sloped}}$ for which,

    \begin{align*}
              1 - \textbf{P}^{\xi^2_{\mathrm{Sloped}}}_{[-m,m] \times [0,n^{\prime} N]} \big[     \gamma_L \overset{h \geq ck}{\longleftrightarrow} \mathscr{F}\mathscr{C}       \big]    \geq   1 - \textbf{P}^{\xi^3_{\mathrm{Sloped}}}_{[-m,m] \times [0,n^{\prime} N]} \big[     \gamma_L \overset{h \geq ck}{\longleftrightarrow} \mathscr{F}\mathscr{C}     \big]   \text{ }  \text{ , } \text{ } \\
    \end{align*}
    
    \noindent one also has that, 
    
    \begin{align*}
\cdots  \overset{\mathrm{(CBC)}}{\geq}   \text{ }  
\bigg| \textbf{P}^{\xi^2_{\mathrm{Sloped}}}_{[-m,m] \times [0,n^{\prime}N]}     \big[   \gamma_L \overset{h \geq ck}{\longleftrightarrow} \mathscr{F}\mathscr{C}   \big]   \bigg[     1 -    \textbf{P}^{\xi^3_{\mathrm{Sloped}}}_{[-m,m] \times [0,n^{\prime}N]}     \big[    h  = \text{ } \xi^1_{\mathrm{Sloped}} \text{ }                      \text{on }  \partial \big(  [-m,m] \times  [n N , n^{\prime} N ]  \big)   \big]  \bigg]  \bigg|    \\ \text{ }     \geq \text{ }               \bigg|  \textbf{P}^{\xi^2_{\mathrm{Sloped}}}_{[-m,m] \times [0,n^{\prime}N]}     \big[  \gamma_L \overset{h \geq ck}{\longleftrightarrow} \mathscr{F}\mathscr{C}         \big]   \bigg|  \bigg|                     1 -    \textbf{P}^{\xi^3_{\mathrm{Sloped}}}_{[-m,m] \times [0,n^{\prime}N]}     \big[    h  = \text{ } \xi^1_{\mathrm{Sloped}} \text{ }                      \text{on }  \partial \big(  [-m,m] \times  [n N , n^{\prime} N ]  \big)    \big]   \bigg|                      \text{ }  \\ \equiv     \bigg|    \frac{  \big\{ \mathscr{F}  \in   F ( \gamma_L     \cap     \mathcal{L}^1_{\mathscr{F}\mathscr{C}_1       }  ) :  \mathscr{F}   \overset{h \geq ck}{\leftrightarrow} \mathscr{F}\mathscr{C}_1      \big\}    }{   \big\{ \mathscr{F} \in     F ( \gamma_L     \cap     \mathcal{L}^1_{\mathscr{F}\mathscr{C}_1}   )   :    \mathscr{F}   \overset{h < ck}{\leftrightarrow} \mathscr{F}\mathscr{C}_1  \big\}       }  \bigg|  \bigg|    1 -    \textbf{P}^{\xi^3_{\mathrm{Sloped}}}_{[-m,m] \times [0,n^{\prime}N]}     \big[   h  = \text{ } \xi^1_{\mathrm{Sloped}} \text{ }                      \text{on }  \partial \big(  [-m,m] \times  [n N , n^{\prime} N ] \big)   \big]  \bigg|        \\ 
   \text{ }     \geq \text{ }     \frac{ \big| \big\{ \mathscr{F}  \in   F ( \gamma_L     \cap     \mathcal{L}^1_{\mathscr{F}\mathscr{C}_1       }  ) :  \mathscr{F}   \overset{h \geq ck}{\leftrightarrow} \mathscr{F}\mathscr{C}_1      \big\}   \big|   }{ \big|  \big\{ \mathscr{F} \in     F ( \gamma_L     \cap     \mathcal{L}^1_{\mathscr{F}\mathscr{C}_1}   )   :    \mathscr{F}   \overset{h < ck}{\leftrightarrow} \mathscr{F}\mathscr{C}_1  \big\}  \big|      }   \bigg|    1 -    \textbf{P}^{\xi^3_{\mathrm{Sloped}}}_{[-m,m] \times [0,n^{\prime}N]}     \big[   h  = \text{ } \xi^1_{\mathrm{Sloped}} \text{ }                      \text{on }  \partial \big(  [-m,m] \times  [n N , n^{\prime} N ] \big)   \big]  \bigg|   \\ 
     \text{ } \geq \text{ }  \frac{N_1}{N_2}      \bigg|    1 -    \textbf{P}^{\xi^3_{\mathrm{Sloped}}}_{[-m,m] \times [0,n^{\prime}N]}     \big[   h  = \text{ } \xi^1_{\mathrm{Sloped}} \text{ }                      \text{on }  \partial \big(  [-m,m] \times  [n N , n^{\prime} N ]  \big)    \big]  \bigg|  \text{ }        \\                 \text{ } \geq \text{ } \frac{N_1}{\delta^1_L} \text{ }       \big(   1 -    \epsilon_L         \big) \\ \text{ } \geq N_3  \text{ }  \big( 1 - \epsilon_L \big)   \text{ }    \geq \text{ } N_3 \text{ }  \delta^2_L     \text{ } \text{ , } \text{ }  
    \end{align*}
    
    \noindent from a lower bound on $\textbf{P}^{\xi^2_{\mathrm{Sloped}}}_{[-m,m] \times [0,n^{\prime}N]}     \big[  \gamma_L \overset{h \geq ck}{\longleftrightarrow} \mathscr{F}\mathscr{C}_1  \big]  $ dependent upon the number of $\mathscr{F}$ for which the desired crossing occurs, suitable $N_1, N_2, N_3$, and $\epsilon_L$, which can be bound below with arbitrary $\delta^2_L$. For $\gamma_R$, following the previously described argument, instead for the connectivity event between $\mathscr{F}\mathscr{C}_n$ and $\gamma_R$ provides the accompanying lower bound for the other case. \boxed{}
    
    \bigskip

     \begin{figure}
\begin{align*}
\includegraphics[width=0.73\columnwidth]{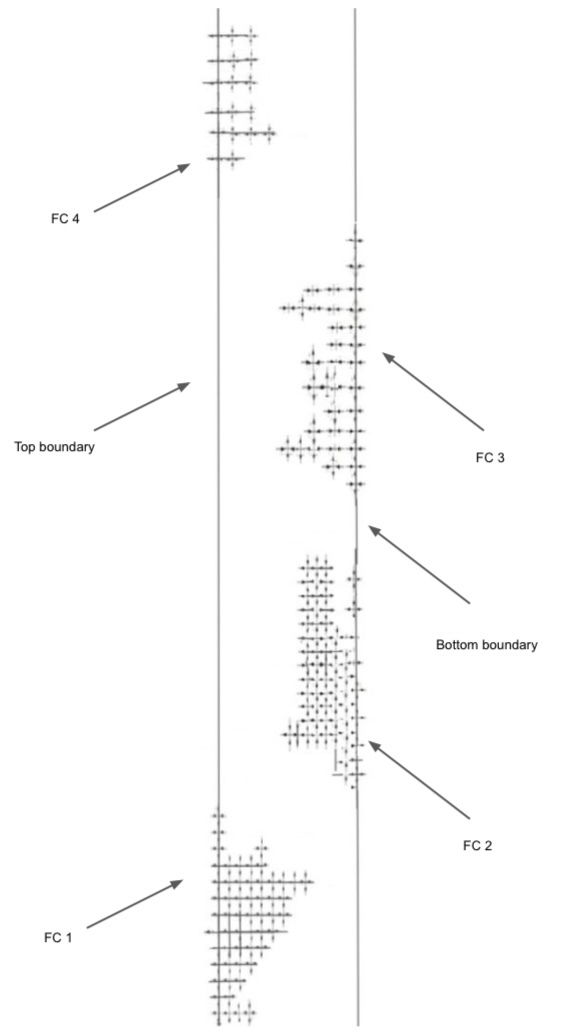}\\
\end{align*}
\caption{\textit{$\frac{\pi}{2}$ rotation of a depiction of four freezing clusters intersecting the strip, in which the first and fourth freezing cluster intersect the top strip boundary, while the second and third freezing cluster intersect the bottom strip boundary}. Within the complementary space to the four \textit{freezing clusters} in the strip, the probability of a macroscopic horizontal crossing is quantified.}
\end{figure}

    \noindent 
    
    \bigskip
    
    \noindent Besides crossing events defined in the vicinity of left and right boundaries, with positive probability freezing clusters can also be embedded within the interior of strip domains, as mentioned from previous remarks and objects. To maintain long horizontal crossings in the infinite volume strip, in spite of the fact that with positive probability there exists countably many \textit{freezing clusters} intersecting the finite, and infinite, volume strip, one can quantify crossing probabilities across smaller regions of the finite volume strip, which occur with good probability in infinite volume. That is, there exists other connectivity events within smaller parts of the strip that occur with good probability between faces within small enough neighborhoods of frozen faces.
    
    \bigskip
    
     \noindent \textbf{Lemma} $\textit{2.3}$ (\textit{strip crossing events with obstructions from a countable number of freezing clusters}): Given a countable number of \textit{freezing clusters} sampled from $\mathcal{F}\mathcal{C}$, the probability of a connectivity even between each pair of neighboring \textit{freezing clusters},

     \begin{align*}
        \textbf{P}^{\xi_{\mathrm{Sloped}}}_{[-m,m] \times [0,n^{\prime}N]} \big[   \mathscr{F}\mathscr{C}_2     \overset{h < ck }{\longleftrightarrow }   \mathscr{F}\mathscr{C}_{N-1}      \big] \text{ } \text{ , } \text{ } \\
     \end{align*}
     
     \noindent has lower bound $c_{\mathrm{FC \text{ } scale}} > 0$.
 
 \bigskip
 
 \noindent \textit{Proof of Lemma 2.3}. We analyze the complementary $h$ $\mathrm{x}$-crossing event,

 \begin{align*}
  \textbf{P}^{\xi_{\mathrm{Sloped}}}_{[-m,m] \times [0,n^{\prime} N]} \big[     \mathscr{F}\mathscr{C}_2     \overset{h < ck }{\longleftrightarrow }   \mathscr{F}\mathscr{C}_{N-1}     \big] \text{ } = \text{ }  \textbf{P}^{\xi_{\mathrm{Sloped}}}_{[-m,m] \times [0,n^{\prime} N]} \big[      \mathscr{F}\mathscr{C}_2     \overset{h \geq ck }{\longleftrightarrow }_{\mathrm{x}}   \mathscr{F}\mathscr{C}_{N-1}       \big]   \text{ } \text{ , } \text{ } \\
 \end{align*}
 
 \noindent and subsequently apply (FKG),
 
 \begin{align*}
         \text{ }       \textbf{P}^{\xi_{\mathrm{Sloped}}}_{[-m,m] \times [0,n^{\prime} N]} \big[       \mathscr{F}\mathscr{C}_2     \overset{h \geq ck }{\longleftrightarrow }_{\mathrm{x}}   \mathscr{F}\mathscr{C}_{N-1}       \big]   \text{ } \equiv   \textbf{P}^{\xi_{\mathrm{Sloped}}}_{[-m,m] \times [0,n^{\prime} N]} \big[ \bigcap_{i=2}^{N-1} \big\{  \mathscr{F}\mathscr{C}_i     \overset{h \geq ck }{\longleftrightarrow }_{\mathrm{x}}   \mathscr{F}\mathscr{C}_{i+1}  \big\}      \big] \text{ }   \\ \text{ }    \geq \text{ }  \prod_{i=2}^{N-1} \text{ }  \textbf{P}^{\xi_{\mathrm{Sloped}}}_{[-m,m] \times [0,n^{\prime} N]} \big[  \mathscr{F}\mathscr{C}_i     \overset{h \geq ck }{\longleftrightarrow }_{\mathrm{x}}   \mathscr{F}\mathscr{C}_{i+1}  \big]  \text{ } \\ \text{ }            \geq \text{ }     \prod_{i=2}^{N-1} \text{ }    c_i  \\ \text{ } \geq \text{ }   c_{\mathrm{FC \text{ } scale}} \text{ } \text{ , } \text{ } \\
 \end{align*}

 \noindent for a horizontal crossing event across all \textit{freezing clusters} simultaneously from $N-2$ vertical crossings and $c_i > 0 \text{ } \forall \text{ } i$ where,
 
 \begin{align*}
        \text{ }  \sqrt[N-2]{c_{\mathrm{FC \text{ } scale}}} \text{ } = \text{ }     \underset{i < N , i , N \in \textbf{Z}: 2 \leq i \leq N-1}{\mathrm{inf}} \text{ }   c_i                   \text{ } \text{ , } \text{ } \\
 \end{align*}
 
 \noindent appears in the penultimate lower bound. \boxed{}
 
 \bigskip
 
 \noindent In arguments relating to crossings across domains, we will also need to make use of the following property of \textit{freezing clusters} pertaining to their separation with respect to boundary conditions. In the following, we consider the number of connected components across \textit{freezing clusters} as a sequence of probability measures with increasing slope along the boundaries.
 
\subsection{Frozen, and unfrozen, faces of the height function in the striip environment}

\bigskip
 
 \noindent \textbf{Definition} $\textit{6A}$ (\textit{partitioning faces of the height function along the strip finite volume boundary}). For a collection of faces over $\textbf{Z}^2$, denote

 \begin{align*}
        \partial^{+}_{\Lambda}  \mathscr{F}   \equiv \big\{     \mathscr{F}^{+} \in F(  \textbf{Z}^2  \cap \Lambda^c ) :  \mathscr{F}^{+} \cap \mathscr{F}(\partial \Lambda ) \neq \emptyset       \big\}                 \text{ } \text{ , } \text{ } \\
 \end{align*}

 \noindent as the collection of faces contained within a subset of $\Lambda^c$, and,

 \begin{align*}
        \partial^{-}_{\Lambda}       \mathscr{F} \text{ }   \equiv  \big\{    \mathscr{F}^{-} \in F(  \Lambda  )  :  \mathscr{F}^{-}       \cap  \mathscr{F}(\partial \Lambda ) \neq \emptyset            \big\}            \text{ } \text{ , } \text{ } \\
 \end{align*}
 
 \noindent as the collection of faces contained within a subset of $\Lambda$, where,

\[
   \mathscr{F} (\partial \Lambda)\text{ }  \equiv    \text{ } \left\{\!\begin{array}{ll@{}>{{}}l} \mathscr{F}^{\partial} = \{ e^{\partial}_1 , e^{\partial}_2 , e^{\partial}_3 , e^{\partial}_4 \}    \in F(     \textbf{Z}^2   \cap \Lambda_{\mathrm{Restriction}}   ) \text{: }      \text{edges }                \text{ } e^{\partial}_i \text{ and } e^{\partial}_j \text{ }         \text{lie incident to } F(  \partial \Lambda           )                      \text{ } \text{ , } \\ \forall\text{ } 
    \mathscr{F}^{\mathrm{in}}\text{ } \in F (  \textbf{Z}^2  \cap   \Lambda    )         \text{ } , \text{ }  \mathscr{F}^{\mathrm{out}} \text{ } \in   F (         \textbf{Z}^2 \cap \Lambda_{\text{ } \mathrm{Restriction}}^c   )\text{ } \text{ , }  \exists \text{ } 1 \leq j \leq 4  \text{: } \text{ }     \mathscr{F}^{\mathrm{in}} \cap \mathscr{F}^{\mathrm{out}} = e^{\partial}_j       \text{ , } \text{ } \\
\end{array}\right.
\]

 \noindent denotes the collection of faces incident to $\partial \Lambda$, and $\Lambda_{\mathrm{Restriction}}$ is the restriction of the faces enclosed within $\Lambda$ that are within Euclidean distance $1$ of any boundary faces, in which any face incident to the boundary is spanned by the four edges $e^{\partial}_1 , \cdots , e^{\partial}_4$, any two of which overlap with intersecting faces from $ F (\text{ }  \textbf{Z}^2 \text{ } \cap \text{ }  \Lambda  \text{ }  ) $, and from $ F ( \text{ }         \textbf{Z}^2 \cap \Lambda_{\text{ } \mathrm{Restriction}}^c    \text{ } )$. Examples of incident edges to the finite volume boundary are provided in the fourth panel of $\textit{Figure 2}$.

 \bigskip

 \noindent With the definition below, we construct sequences of sloped boundary conditions. Such sequences of boundary conditions are significant not only for further discussion of differences in the strip, but also for obtaining desirable properties of slopped symmetric domains in the strip.

 \bigskip
 
 \noindent \textbf{Definition} $\textit{6B}$ (\textit{sequences of sloped-boundary conditions, from the irreducibility property of the height function for the six-vertex model}). Over the finite strip $[-m,m] \times [0,n^{\prime} N]$, a \textit{sequence of sloped boundary-conditions},

 \begin{align*}
\{   \xi^{\mathrm{Sloped}}_k  \big(   F \big( \partial \big( [-m,m] \times [0,n^{\prime} N] \big) \big)    \big)     \}_{k \in \textbf{N}} \equiv \{   \xi^{\mathrm{Sloped}}_k    \}_{k \in \textbf{N}} \text{ } \text{ , } 
\end{align*}

\noindent with a corresponding sequence of \textit{irreducible} height functions,

\begin{align*}
\{  h^{\mathrm{Sloped}}_k \big(   F \big( [-m,m] \times [0,n^{\prime} N]  \big)    \big)  \}_{k \in \textbf{N}} \equiv \{  h^{\mathrm{Sloped}}_k   \}_{k \in \textbf{N}} \text{ } \text{ , } 
\end{align*}

\noindent satisfies the condition,
 
 \begin{align*}
       \bigg[ \underset{\mathscr{F} \in  F \big( [-m,m] \times [0,n^{\prime} N]  \big) }{\bigcup} h^{\mathrm{Sloped}}_{k-1} \big(   \mathscr{F}  \big)        \bigg]      \cap  \bigg[         \underset{\mathscr{F} \in  F \big( [-m,m] \times [0,n^{\prime} N]  \big) }{\bigcup} h^{\mathrm{Sloped}}_k \big(  \mathscr{F}      \big)  \bigg]   \equiv   \mathscr{F} \text{ } \text{ , } 
\end{align*}
 
\noindent for $0 \leq k - 1 < k < \big|  F \big( [-m,m] \times [0,n^{\prime} N]  \big) \big|$, and some $\mathscr{F}$.
    
     \begin{figure}
\begin{align*}
\includegraphics[width=0.63\columnwidth]{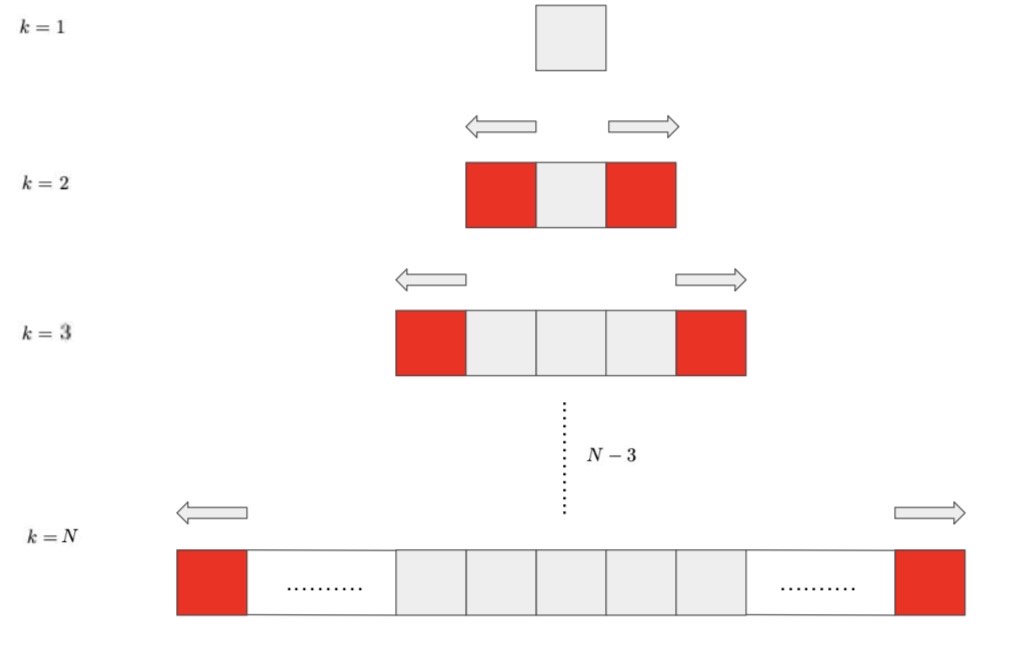}\\
\end{align*}
\caption{\textit{A depiction of the recursive procedure for defining the boundary conditions from the height function along the bottom and top boundaries of strips of the square lattice}. For a single face given by $k \equiv 1$, the interface created from sequence of \textit{irreducible} height functions begins by enforcing the boundary condition of the height function, through the image of the graph homomorphism, on a single face. In all subsequent steps, for an arbitrary, strictly positive $N$, faces are appended simultaneously on the left and right sides of the top and bottom boundaries.}
\end{figure}

 \bigskip
 
 \noindent \textbf{Definition} $\textit{7}$ (\textit{collections of faces incident to the strip over which h achieves its minima, maxima and all intermediate height values}). Define,

 \begin{align*}
      \partial  \mathscr{F}_{ \mathrm{inf} \text{ } h } = \big\{    \mathscr{F}_1,   \mathscr{F}_2 \in F(  \textbf{Z}^2  ) :       h ( \mathscr{F}_1 )\equiv \underset{\mathscr{F}_2}{\mathrm{inf}} \text{ } h (  \mathscr{F}_2 )          \big\}   \text{ } \text{ , } \text{ }    \\
 \end{align*}

 \noindent corresponding to the collection of faces over the boundary of the strip for which $h$ achieves its infimum. Similarly,
 
 \begin{align*}
    \partial  \mathscr{F}_{\text{ } h < ck \text{ } }    \text{ }     = \text{ }   \big\{      \mathscr{F}^{\prime}_1 , \mathscr{F}_2 \in F(  \textbf{Z}^2  )  : \underset{\mathscr{F}_2}{\mathrm{inf}} \text{ } h (  \mathscr{F}_2  )       <     h (  \mathscr{F}^{\prime}_1         ) <   \underset{\mathscr{F}_2}{\mathrm{sup}} \text{ } h ( \mathscr{F}_2  )    \big\}     \text{ } \text{ , }  \\
 \end{align*}
 
 \noindent corresponding to the collection of faces over the boundary for which $h$ achieves a height above the global minimum, and maximum, of $h$, and finally,

 \begin{align*}
    \partial  \mathscr{F}_{\text{ } h \geq ck \text{ } }    \text{ }  = \text{ }   \big\{       \mathscr{F}^{\prime}_2 , \mathscr{F}_2 \in F( \textbf{Z}^2 ) \text{ } , \mathscr{F}_3 \in \text{ }  \partial \mathscr{F}_{\text{ } h < ck\text{ } }     :    h ( \mathscr{F}^{\prime}_2  ) <    \underset{\mathscr{F}_2}{\mathrm{sup} }  \text{ } h (  \mathscr{F}_2   )   \text{ } , \text{ }  \underset{\mathscr{F}_3}{\mathrm{sup}} \text{ }  h (\mathscr{F}_3  )  <   h (\mathscr{F}^{\prime}_2  )     \big\}      \text{ } \text{ , } \text{ } \\
 \end{align*}

 \noindent corresponding to the collection of faces over the boundary for which $h$ achieves a height exceeding $ck$.

 
 \bigskip

\noindent With the following, we demonstrate that the crossing event that one would like to occur with sufficiently good probability over long strips is dependent upon two factors, the first of which captures effects from the (SMP) property as boundary conditions of new faces along the strip are appended to the interface in infinite volume, and the second of which captures effects from the number of frozen faces relative to unfrozen faces that appear when taking the weak volume limit. For estimates on other crossing probabilities later in the section before transitioning to the cylinder, a combination of effects from both (SMP) and (CBC) arise for similar reasons, in addition to the fact that sloped Gibbs states over the strip can, with positive probability, contain finitely many connected components of frozen faces, in comparison to flat Gibbs states which have no frozen faces. Regardless of the difference in slope of boundary conditions, logarithmic delocalization results for either case for the six-vertex model can be established, by further analyzing how flat symmetric domains change in the presence of rational slopes. 

\bigskip

\noindent \textbf{Proposition} $\textit{2.3}$ (\textit{decomposition of the the lower bound constant from Proposition 2.2}): The constant in the lower bound of the ratio of crossing probabilities given in \textit{Proposition 2.2} can be decomposed as,

\begin{align*}
    \mathcal{C}_{\mathscr{I}} \equiv \text{ } \mathcal{C}_{1,2} \equiv     \mathcal{C}_1 +  \mathcal{C}_2                    \text{ } \text{ , } \text{ } \\
\end{align*}

\noindent where $\mathcal{C}_{1,2} \equiv \mathcal{C}_{1,2}\big(  [-m,m] \times [0,n^{\prime} N ] , a_x , b_x , c_x  \big)$, $\mathcal{C}_1 \equiv \mathcal{C}_1\big(  [-m,m] \times [0,n^{\prime} N ] , a_x , b_x , c_x \big)$, and $\mathcal{C}_2 \equiv \mathcal{C}_2\big( [-m,m] \times [0,n^{\prime} N ] , a_x , b_x , c_x \big)$.

\bigskip

\noindent \textit{Proof of Proposition 2.3}. From arguments used previously for \textit{Proposition 2.2}, one reads, from contributions in the lower bound, that

\begin{align*}
  \text{ }       \mathcal{C}_{\mathscr{I}} \equiv \text{ }   \mathcal{S}^{\prime\prime} \bigg[                  \mathscr{I}_x^{| V(  [-m,m] \times [0 , n^{\prime} N]  )       |}    \text{ } + \text{ }       2   \sum_{x \in \partial  V( [-m,m] \times [0,n^{\prime} N ] ) } \text{ }        h^{\xi^{\mathrm{Sloped}^{\prime\prime\prime\prime}}} \Delta^{\xi^1 , \xi^2}_x          \bigg]   \text{ } \text{ , } 
  \end{align*}

\noindent from which it suffices to show that,

\begin{align*}
   \mathcal{C}_1 \propto  \text{ } \mathcal{S}^{\prime\prime}   \mathscr{I}_x^{|  V(  [-m,m] \times [0 , n^{\prime} N]  )   |}   \text{ } \text{ , } \text{ } \\
\end{align*}
  
  \noindent   and also that,
  
  \begin{align*}
    \mathcal{C}_2 \text{ }  \propto  \text{ } 2  \text{ }  \mathcal{S}^{\prime\prime}      \sum_{x \in \partial  V( \text{ } [-m,m] \times [0,n^{\prime} N ]\text{ } ) } \text{ }        h^{\xi^{\mathrm{Sloped}^{\prime\prime\prime\prime}}} \Delta^{\xi^1 , \xi^2}_x                  \text{ } \text{ . } \text{ }   \\
  \end{align*}
  
  \noindent To accomplish this, observe, for the first term,

  \begin{align*}
    \mathcal{S}^{\prime\prime}   \mathscr{I}_x^{|\text{ }  V( \text{ } [-m,m] \times [0 , n^{\prime} N] \text{ } )      \text{ } |}   \propto \big(  V \big(  [-m,m] \times [0,n^{\prime} N]  \big) \big)    \big(   w_{6V}      \big)    \text{ }  \propto \mathcal{C}_1 \big(  [-m,m] \times [0,n^{\prime} N ] , a_x , b_x , c_x \big) \equiv \mathcal{C}_1   \text{ }  \text{ , } \text{ } \\
  \end{align*}

  \noindent where $w_{6V}$ is the weight introduced in \textit{1.3}, while similarly, for the second term,

  \begin{align*}
       \text{ }   2 \text{ }   \mathcal{S}^{\prime\prime}      \sum_{x \in \partial  V( \text{ } [-m,m] \times [0,n^{\prime} N ]\text{ } ) } \text{ }        h^{\xi^{\mathrm{Sloped}^{\prime\prime\prime\prime}}} \Delta^{\xi^1 , \xi^2}_x     \propto \text{ }  \big(   V( [-m,m] \times [0 , n^{\prime} N] ) \big) \big(     w_{6V}       \big) \\    \propto \text{ }   \mathcal{C}_2  \bigg[  [-m,m] \times [0,n^{\prime} N ] , a_x , b_x , c_x \bigg]  \equiv \mathcal{C}_2 \text{ } \text{ . } \text{ }  \\
  \end{align*}

  \noindent As a result, 
  
  \begin{align*}
       \mathcal{C}_{1,2} \propto \mathcal{C}_1 + \mathcal{C}_2      \equiv \text{ } \mathcal{C}_1 \bigg[   [-m,m] \times [0,n^{\prime} N ] , a_x , b_x , c_x \bigg] +  \mathcal{C}_2 \bigg[  [-m,m] \times [0,n^{\prime} N ] , a_x , b_x , c_x \bigg]              \text{ } \text{ , }  \\
\end{align*}

\noindent from which we conclude the argument. \boxed{}

\bigskip

\noindent Besides the result above, there should also exist a constant - which also exhibits the difficulty that sloped boundary conditions impose on horizontal crossing probabilities - which we denote as $\mathcal{C}^{\prime}_{\mathrm{Sloped}}$. Given the existence of a horizontal crossing event $\mathcal{H}^{\prime}$, the probability of $\mathcal{H}^{\prime}$ occurring under the six-vertex model with sufficiently flat boundary conditions is greater than the same event occurring under the six-vertex model with sloped boundary conditions. Quantifying the absolute value distance between the ratio of these crossing probabilities for different sequences of sloped boundary conditions,  $ \frac{\textbf{P}^{\xi^{\mathrm{Sloped}}_{h_k}}_{[-m,m]\times[0,n^{\prime}N]} [    \mathcal{H}^{\prime}     ]}{\textbf{P}^{\xi^{\mathrm{Sloped}}_{h_{k-1}}}_{[-m,m]\times[0,n^{\prime}N]} [   \mathcal{H}^{\prime} ]  } $, with $1$, yields the constant. By taking the width of the strip large enough, with positive probability there exists finite volumes of the strip which satisfy any one of the possible configurations, described in \textbf{Lemma} \textit{2.1}, to counterparts in infinite volume satisfying the same property.

\bigskip
 
\noindent \textbf{Lemma} $\textit{2.4}$ (\textit{freezing clusters are 'well separated' with respect to values of the height function encoded through boundary conditions of the probability measure over the strip}): For a sequence of boundary conditions as provided in \textbf{Definition} $6$B, the difference in a horizontal crossing across the strip is bounded below by the following quantity, 

\begin{align*}
    \bigg|   \frac{\textbf{P}^{\xi^{\mathrm{Sloped}}_{h_k}}_{[-m,m]\times[0,n^{\prime}N]} \big[    \mathcal{H}^{\prime}   \big]}{\textbf{P}^{\xi^{\mathrm{Sloped}}_{h_{k-1}}}_{[-m,m]\times[0,n^{\prime}N]} \big[  \mathcal{H}^{\prime}    \big]}  -  1   \bigg|   \text{ }  \geq \text{ }    \mathcal{C}^{\prime}_{\mathrm{Sloped}}   \text{ }   \text{ , } \text{ }    \\
\end{align*}


\noindent for $j > 0$, $\mathcal{H}^{\prime} \equiv \mathcal{H}^{\prime}\big([-m,m]\times[0,n^{\prime}N]\big)$ and suitable, strictly positive, $\mathcal{C}^{\prime}_{\mathrm{Sloped}}$, with,

\begin{align*}
 \text{ }       \mathcal{C}^{\prime}_{\mathrm{Sloped}} \text{ } \equiv \text{ }  \mathcal{C}^{\prime}_{\mathrm{Sloped}} \bigg[   \big[-\frac{m}{2},\frac{m}{2}\big] \times [0,n^{\prime} N ]   ,       a_x , b_x , c_x  \bigg]                        \text{ } \text{ . } \text{ } \\
\end{align*}

\bigskip

\noindent By making use of properties such as the one above for sequences of sloped boundary conditions, there exists a third constant, for which the distance between \textit{freezing clusters} which possibly intersect the bottom and top boundaries of the strip, can be taken sufficiently small. This sufficiently small constant is dependent upon the two previously obtained constants depending on the sloped boundary conditions.

\bigskip

\noindent \textit{Proof of Lemma 2.4}. By direct computation, from the LHS given in the statement of \textbf{Lemma} $\textit{2.4}$,

\begin{align*}
        \bigg| \frac{\textbf{P}^{\xi^{\mathrm{Sloped}}_{h_k}}_{[-m,m]\times[0,n^{\prime}N]} \big[   \mathcal{H}^{\prime}  \big]}{\textbf{P}^{\xi^{\mathrm{Sloped}}_{h_{k-1}}}_{[-m,m]\times[0,n^{\prime}N]} \big[  \mathcal{H}^{\prime}   \big]} -  1   \bigg| \text{ }  \overset{\mathrm{(\textit{*} )}}{\geq} \text{ }    \big| \mathcal{C}_{\mathscr{I}} -  1  \big| \text{ }   \overset{\mathrm{(\textbf{Proposition} \text{ } \textit{2.3} )}}{=} \text{ } \big|  \mathcal{C}_{1,2} \text{ }  - \text{ }  1      \big| 
        \end{align*}
        
        \noindent in which case we further lower bound the last estimate with,
        
        \begin{align*}
        \text{ }      \big|  \mathcal{C}^{\prime}_{\mathrm{Sloped}}  -  1  \big|           \text{ } \text{ , }  \\
\end{align*}

   \noindent after having set $\mathscr{H}\mathscr{C} \equiv \mathcal{H}^{\prime}$ from applying the bound given in Proposition \textit{2.2}, given a suitable constant in the lower bound for which,
   
   \begin{align*}
    \mathcal{C}_{1,2} \equiv    \mathcal{C}_1 + \mathcal{C}_2  \geq  \mathcal{C}^{\prime}_1 + \mathcal{C}^{\prime}_2 \text{ , }      
    \end{align*}
    
    \noindent where, 
    
    \begin{align*}
  \mathcal{C}^{\prime}_1 + \mathcal{C}^{\prime}_2 \equiv \mathcal{C}^{\prime}_1 \bigg[   \big[-\frac{m}{2} , \frac{m}{2}\big] \times [0, n^{\prime} N ]     ,    a_x , b_x , c_x     \bigg]  + \mathcal{C}^{\prime}_2 \bigg[      \big[-\frac{m}{2},\frac{m}{2}\big] \times [0,n^{\prime} N ]         ,  a_x , b_x ,c_x  \bigg] \text{ }  \propto \text{ }   \mathcal{C}^{\prime}_{\mathrm{Sloped}}          \text{ } \text{ , } \text{ } \\
   \end{align*}
   
   \noindent in addition to the fact that, for $(\textit{*})$,

   \begin{align*}
        \text{ }           \frac{ \textbf{P}^{\xi^{\mathrm{Sloped}}_{h_k}}_{[-m,m]\times[0,n^{\prime} N]} \big[       \mathscr{H}\mathscr{C}            \big]  \text{ }}{ \textbf{P}^{\xi^{\mathrm{Sloped}}_{h_{k-1}}}_{[-m,m]\times[0,n^{\prime} N]} \big[                     \mathscr{H}\mathscr{C}        \big]    \text{ }    }               \text{ } > \text{ } \mathcal{C}_{\mathscr{I}}          \text{ } \text{ , } 
   \end{align*}

   \noindent for the intersection,
   
    \begin{align*}
       \text{ }     \partial^{k^{\prime}}_{< ck}  \cap  \partial^{k^{\prime}}_{\text{ } \mathrm{inf} \text{ } }   =    \big\{  \mathscr{F}^{\prime\prime}_0     , \mathscr{F}^{\prime\prime}_1 \in \partial^{k^{\prime}}_{< ck} \text{ } , \text{ }   \mathscr{F}^{\prime\prime}_2  \in \text{ }  \partial^{k^{\prime}}_{\mathrm{inf}} :                       \underset{\mathscr{F}^{\prime\prime}_2}{\mathrm{inf}} \text{ }   h ( \ \mathscr{F}^{\prime\prime}_2  )    \text{ }  \leq  h (      \mathscr{F}^{\prime\prime}_0 )  <\underset{\mathscr{F}^{\prime\prime}_1}{\mathrm{sup}} \text{ }   h (     \mathscr{F}^{\prime\prime}_1      )        \big\}           \text{ } \text{ , } \text{ } \\
 \end{align*}

  \noindent for all $k^{\prime} >0$, where $\partial \mathscr{F}^{k^{\prime}}_{\text{ } \mathrm{inf} \text{ } h     \text{ } } \equiv \partial^{k^{\prime}}_{\mathrm{inf}}$, $\partial \mathscr{F}^{k^{\prime}}_{h < ck }\equiv \partial^{k^{\prime}}_{< ck}$, and $\partial  \mathscr{F}^{k^{\prime}}_{h \geq ck}  \equiv \partial^{k^{\prime}}_{\geq ck}$, $\mathscr{F}\mathscr{C}^{\xi_{k+j}}_k$ denotes the \textit{freezing cluster} under $\xi^{\mathrm{Sloped}}_{k+j}$, $\mathscr{F}\mathscr{C}^{\xi_{k}}_k$ denotes the \textit{freezing cluster} under $\xi^{\mathrm{Sloped}}_k$, and $\mathscr{H}\mathscr{C}$ a horizontal crossing, and,

 \begin{align*}
\mathcal{C}_{\mathscr{I}} \propto \text{ }  \bigg|   \big\{      \text{connected components of }  \mathscr{F}\mathscr{C}^{\xi_{k+j}}_k                                \big\} \cap \big\{      \text{connected components of }  \mathscr{F}\mathscr{C}^{\xi_{k}}_k    \big\}   \bigg|   \text{ , }  \\ 
  \end{align*}

  \noindent denotes a suitable constant satisfying, which provides the desired estimate. \boxed{}

 \bigskip
 
 \noindent With this result, we examine the following bound involving height-function crossings in domains. The inequality is expected to hold as for sufficiently flat boundary conditions, however with an adjusted constant for sloped boundary conditions. 
 
 \bigskip
    
    \noindent \textbf{Lemma} $\textit{2.5}$ (\textit{upper bound for the probability of the segment connectivity event between} $\mathcal{I}_0$ and $\widetilde{\mathcal{I}_0}$, and between $\mathcal{I}_j$ and $\widetilde{\mathcal{I}_j}$, \textit{under sloped boundary conditions}). Fix $\xi^{\prime}_{\mathrm{Sloped}} \leq \xi_{\mathrm{Sloped}}$. With regards to connectivity events between two segments $\mathcal{I}_0$ and $\widetilde{\mathcal{I}_0}$, which are respectively on the top and bottom boundaries of the finite volume strip, 
    
          \begin{figure}
\begin{align*}
\includegraphics[width=0.75\columnwidth]{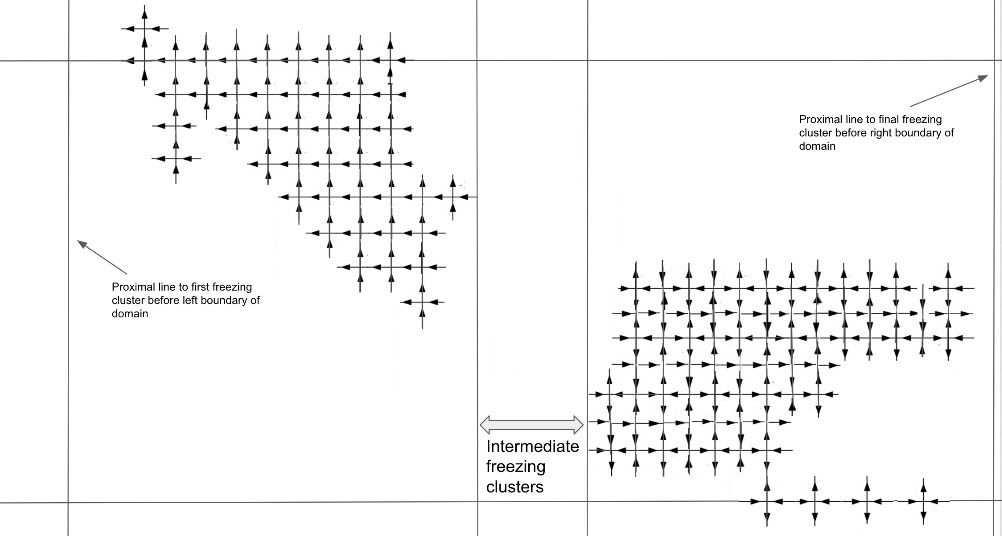}\\
\end{align*}
\caption{\textit{Lines placed after right and left boundaries of strip domains}. For a countable number of \textit{freezing clusters} along the top and/or bottom boundaries of the strip, the probability of a horizontal crossing is dependent upon the values of the height function surrounding each \textit{freezing cluster}.}
\end{figure}

    \begin{align*}
  \textbf{P}^{\xi_{\mathrm{Sloped}}}_{\textbf{Z} \times [ 0 , n ^{\prime} N ]}  \big[    \mathcal{I}_0  \overset{h \geq (1-c) k }{\longleftrightarrow}  \widetilde{\mathcal{I}_0 }                         \big] \leq 1 -   \widetilde{\mathscr{R}} \text{ }  \textbf{P}^{\xi^{\prime}_{\mathrm{Sloped}}}_{\textbf{Z} \times [ 0 , n ^{\prime} N ]}  \bigg[     \big\{   \mathcal{I}_{-\mathscr{I}}  \overset{h \leq ck }{\longleftrightarrow}  \widetilde{\mathcal{I}_{-\mathscr{I}}}     \big\} \cap \big\{     \mathcal{I}_{-\mathscr{J}}  \overset{h \leq ck }{\longleftrightarrow}  \widetilde{\mathcal{I}_{-\mathscr{J}}}         \big\}  \cap \big\{     \mathcal{I}_{-\mathscr{K}}  \overset{h \leq ck }{\longleftrightarrow}  \widetilde{\mathcal{I}_{-\mathscr{K}}}         \big\}   \bigg]                      \text{ } \text{ , } \text{ } \\
    \end{align*}
    
    \noindent for $ 0 \neq \mathscr{I} < \mathscr{J} < \mathscr{K}$ with a strictly positive multiplicative factor $\widetilde{\mathscr{R}}$ for the intersection of crossing probabilities appearing in the upper bound, which is dependent on the slope of the boundary conditions, and the width $n^{\prime}$ of the strip,

    \begin{align*}
  \widetilde{\mathscr{R}} \equiv \widetilde{\mathscr{R}}_{\{ \mathcal{L}_{\mathcal{I}_{-\mathscr{I}}} , \mathcal{L}_{\widetilde{\mathcal{I}_{-\mathscr{I}}}}  \} } +  \widetilde{\mathscr{R}}_{\{ \text{ } \mathcal{L}_{\mathcal{I}_{-\mathscr{J}}} , \mathcal{L}_{\widetilde{\mathcal{I}_{-\mathscr{J}}}} \} } + \widetilde{\mathscr{R}}_{\{ \text{ } \mathcal{L}_{\mathcal{I}_{-\mathscr{K}}} , \mathcal{L}_{\widetilde{\mathcal{I}_{-\mathscr{K}}}} \text{ } \} } \equiv   \widetilde{\mathscr{R}}_{\{ \text{ } \mathcal{L}_{\mathcal{I}_{-\mathscr{J}}} , \cdots , \text{ }  \mathcal{L}_{\widetilde{\mathcal{I}_{-\mathscr{H}}}}\text{ } \}}\text{ }  \equiv \widetilde{\mathscr{R}}_{\{ \text{ } \mathcal{L}_{\mathcal{I}_{-\mathscr{J}}} , \cdots , \text{ }  \mathcal{L}_{\widetilde{\mathcal{I}_{-\mathscr{H}}}}\text{ } \}}(\xi^{\prime}_{\mathrm{Sloped}} , n^{\prime})  \text{ , }
        \end{align*}

    \noindent for $\mathcal{L}_{\mathcal{I}_{-k}}$ denotes the line centered about some point in $\mathcal{I}_{-k}$. 
    
    \bigskip

  \noindent  \textit{Proof of Lemma 2.5}. We adopt notation from previous arguments, in which we begin with flat boundary conditions $\xi_{\mathrm{Flat}}$ and translate the vertical crossing event along the strip before \textit{freezing clusters} occur with positive probability as the height of the boundary conditions increase. To begin, fix $i=0$, in which case the RHS of the inequality provided in the statement above can be lower bounded with

  \begin{align*}
      \textbf{P}^{\xi^{\prime}_{\mathrm{Sloped}}}_{[-m,m] \times [ 0 , n ^{\prime} N ]}  \big[   \mathcal{I}_{-\mathscr{I}}  \overset{h \leq ck }{\longleftrightarrow}  \widetilde{\mathcal{I}_{-\mathscr{I}}}  ]        \textbf{P}^{\xi^{\prime}_{\mathrm{Sloped}}}_{[-m,m] \times [ 0 , n ^{\prime} N ]}  \big[  \mathcal{I}_{-\mathscr{J}}  \overset{h \leq ck }{\longleftrightarrow}  \widetilde{\mathcal{I}_{-\mathscr{J}}} ]   \textbf{P}^{\xi^{\prime}_{\mathrm{Sloped}}}_{[-m,m] \times [ 0 , n ^{\prime} N ]}  \big[     \mathcal{I}_{-\mathscr{K}}  \overset{h \leq ck }{\longleftrightarrow}  \widetilde{\mathcal{I}_{-\mathscr{K}}}     \big]                      \text{ }                          \text{ } \text{ , } \text{ } \\
  \end{align*}

 \noindent by (FKG) from the intersection of the three crossing events occurring. Next, for each probability, given the existence of flat boundary conditions for which the crossing event occurs, there exists a lower bound,
    
    \begin{align*}
   \underset{ \text{ } k   \in  \{  \mathscr{I} ,   \mathscr{J}, \mathscr{K}  \}}{\mathrm{inf}} \big\{    \textbf{P}^{\xi^{\prime}_{\mathrm{Sloped}}}_{[-m,m] \times [ 0 , n ^{\prime} N ]}  \big[   \mathcal{I}_{-k}  \overset{h \leq ck }{\longleftrightarrow}  \widetilde{\mathcal{I}_{-k}}  ]      \big\} \text{ }  \geq \text{ }  \textbf{P}^{\xi_{\mathrm{Flat}}}_{[-m,m] \times [0,n^{\prime} N]} \big[   \mathcal{I}_0  \overset{h \geq (1-c) k }{\longleftrightarrow}  \widetilde{\mathcal{I}_0 }       \big]   \text{ } \text{ , } \text{ } \\
    \end{align*}

    \noindent dependent upon the segment connectivity event for $i=0$, implying,
    
    \begin{align*}
          \underset{ k   \in   \{  \mathscr{I} ,   \mathscr{J}, \mathscr{K}  \}}{\prod} \textbf{P}^{\xi^{\prime}_{\mathrm{Sloped}}}_{[-m,m] \times [ 0 , n ^{\prime} N ]}  \big[ \mathcal{I}_{-k}  \overset{h \leq ck }{\longleftrightarrow}  \widetilde{\mathcal{I}_{-k}}  ]     \text{ }            \geq \text{ }   \big( \textbf{P}^{\xi_{\mathrm{Flat}}}_{[-m,m] \times [0,n^{\prime} N]} \big[   \mathcal{I}_0  \overset{h \geq (1-c) k }{\longleftrightarrow}  \widetilde{\mathcal{I}_0 }       \big]   \big)^3  \text{ } \text{ . }  \end{align*}

       \noindent To ensure that the probability in the lower bound of \textbf{Lemma} \textit{2.5} is strictly less than $1$,

          \begin{align*}
          \textbf{P}^{\xi_{\mathrm{Sloped}}}_{\textbf{Z} \times [ 0 , n ^{\prime} N ]}  \big[     \mathcal{I}_0  \overset{h \geq (1-c) k }{\longleftrightarrow}  \widetilde{\mathcal{I}_0 }                      \big] < 1 \\    \Updownarrow \\ \text{ }     \textbf{P}^{\xi_{\mathrm{Sloped}}}_{\textbf{Z} \times [ 0 , n ^{\prime} N ]}  \big[     \mathcal{I}_0  \overset{h \geq (1-c) k }{\longleftrightarrow}  \widetilde{\mathcal{I}_0 }                       \big] \text{ } \overset{(\mathrm{FKG})}{\leq} 1 - \text{ }   \widetilde{\mathscr{R}} \bigg[         \underset{  k   \in  \{  \mathscr{I} , \mathscr{J}, \mathscr{K}  \}}{\prod} \textbf{P}^{\xi^{\prime}_{\mathrm{Sloped}}}_{[-m,m] \times [ 0 , n ^{\prime} N ]}  \big[ \mathcal{I}_{-k}  \overset{h \leq ck }{\longleftrightarrow}  \widetilde{\mathcal{I}_{-k}}  ]          \bigg]       \\ \leq 1 -   \widetilde{\mathscr{R}}   \big(     \textbf{P}^{\xi_{\mathrm{Sloped}}}_{[-m,m] \times [0,n^{\prime} N]} \big[   \mathcal{I}_0  \overset{h \geq (1-c) k }{\longleftrightarrow}  \widetilde{\mathcal{I}_0 }       \big]        \big)^3  \text{ } < \text{ } 1 \\ \Updownarrow \\                                     1 -  \widetilde{\mathscr{R}}  \big(      \textbf{P}^{\xi_{\mathrm{Sloped}}}_{[-m,m] \times [0,n^{\prime} N]} \big[   \mathcal{I}_0  \overset{h \geq (1-c) k }{\longleftrightarrow}  \widetilde{\mathcal{I}_0 }       \big]           \big)^3  < 1 \\  \Updownarrow \\ 
          \big(    \textbf{P}^{\xi_{\mathrm{Sloped}}}_{[-m,m] \times [0,n^{\prime} N]} \big[   \mathcal{I}_0  \overset{h \geq (1-c) k }{\longleftrightarrow}  \widetilde{\mathcal{I}_0 }       \big]          \big)^{-1}      \text{ } > \text{ }      \sqrt[3]{\widetilde{\mathscr{R}}} \text{ } > \text{ } 0       \text{ }     \text{ , }  \tag{\textit{Constant for upper bound on segment connectivity event}}
    \end{align*}

    \noindent exhibiting that the claim in the \textbf{Lemma} above holds for $i \neq 0$, given $\widetilde{\mathscr{R}}$ sufficiently small satisfying the conditions above. Besides $i \neq 0$, for other $i$ so that $i \in [-2 \delta^{\prime\prime} n , 2 \delta^{\prime\prime}n]$, observe, given $\delta^{\prime\prime} \equiv \frac{\delta}{3}$,
    
    \begin{align*}
     \textbf{P}^{\xi^{\mathrm{Sloped}}}_{[-m,m] \times [0,n^{\prime} N]} \big[  [ 0 , \lfloor \delta^{\prime\prime}  n \rfloor     \times \{ 0 \}   \overset{h \geq (1-c)k }{\longleftrightarrow} [  i   , i +  \lfloor \delta^{\prime\prime}  n \rfloor  ]  \times \{n \}       \big]   \\  \overset{(\mathrm{CBC})}{\leq}        \textbf{P}^{(\xi^{\mathrm{Sloped}})_2}_{[-m,m] \times [0,n^{\prime} N]} \big[  [ 0 , \lfloor \delta^{\prime\prime}  n \rfloor     \times \{ 0 \}   \overset{h \geq (1-c)k }{\longleftrightarrow} [  i   , i +  \lfloor \delta^{\prime\prime}  n \rfloor  ]  \times \{n \}       \big]    \leq 1-c                 \text{ } \text{ . } 
    \end{align*}
    
    \noindent Finally, for $i \notin [-2 \delta^{\prime\prime} n , 2 \delta^{\prime\prime}n]$, namely $i > 2 \delta^{\prime\prime} n$, observe,
    
   \begin{align*}
    \textbf{P}^{\xi^{\mathrm{Sloped}}}_{[-m,m]\times[0,n^{\prime}N]} \big[ [ \lfloor \delta^{\prime\prime} n \rfloor , 2  \lfloor \delta^{\prime\prime} n \rfloor ] \times \{ 0 \}  \overset{ h \leq ck}{\longleftrightarrow }         [ - i + \lfloor \delta^{\prime\prime} n \rfloor , - i + 2  \lfloor \delta^{\prime\prime} n \rfloor ]     \times \{ n \} \big] \\ \Updownarrow \\   \textbf{P}^{\xi^{\mathrm{Sloped}}-k^{\prime\prime}}_{[-m,m]\times[0,n^{\prime}N]} \big[ [ \lfloor \delta^{\prime\prime} n \rfloor , 2  \lfloor \delta^{\prime\prime} n \rfloor ] \times \{ 0 \}  \overset{ h \leq (1-c) k }{\longleftrightarrow }         [ - i + \lfloor \delta^{\prime\prime} n \rfloor , - i + 2  \lfloor \delta^{\prime\prime} n \rfloor ]     \times \{ n \} \big]         \text{ } \text{ , }
   \end{align*} 
    
    \noindent for suitable $k^{\prime\prime} > 0$, from which the final probability obtained under boundary conditions $\xi^{\mathrm{Sloped}}-k^{\prime\prime}$ obtained above can itself be lower bounded with,
    
       \begin{align*}
            \textbf{P}^{\xi^{\mathrm{Sloped}}}_{[-m,m]\times[0,n^{\prime}N]} \big[ [ 0 ,   \lfloor \delta^{\prime\prime} n \rfloor ] \times \{ 0 \}  \overset{ h \leq ck}{\longleftrightarrow }         [  i  , i +  \lfloor \delta^{\prime\prime} n \rfloor ]     \times \{ n \} \big]   \text{ } \text{ , } 
    \end{align*}
    
    \noindent implying,
    
    \begin{align*}
         \textbf{P}^{\xi^{\mathrm{Sloped}}}_{[-m,m]\times[0,n^{\prime}N]} \big[ [ 0 ,   \lfloor \delta^{\prime\prime} n \rfloor ] \times \{ 0 \}  \overset{ h \leq ck}{\longleftrightarrow }         [  i  , i +  \lfloor \delta^{\prime\prime} n \rfloor ]     \times \{ n \} \big]  <  \frac{1}{2}   \text{ } \text{ , } 
    \end{align*}

    \noindent from which we conclude the argument, for a suitable constant $< \frac{1}{2}$. \boxed{}

\bigskip
    
    \noindent With the claim below, we characterize vertical crossing probabilities across domains.
    
    \bigskip
    
   \noindent \textbf{Lemma} $\textit{2.6}$ (\textit{vertical crossing probabilities across domains}).  Lower bound estimates for vertical crossing events $\mathcal{V}^{h \leq (1-c)k}(\mathscr{D}) \equiv \mathcal{V}^{h \leq (1-c)k}$ under $\xi^{\mathrm{Sloped}}$ across domains $\mathscr{D}$ take the form, with suitable constants $\mathfrak{C}_{\mathscr{D}}^{\prime,\prime\prime,\prime\prime\prime}  $ and $C^{\emptyset}_{\mathcal{V}}$,

   \begin{align*}
    \mathfrak{C}_{\mathscr{D}}^{\prime,\prime\prime,\prime\prime\prime}      \leq \textbf{P}^{\xi^{\mathrm{Sloped}}}_{\mathscr{D}}  [   \mathcal{V}^{h \leq (1-c)k}    ]  \leq          C^{\emptyset}_{\mathcal{V}}    \text{ } \text{ , } \text{ } \\
   \end{align*}
   
   \noindent for non-intersecting left and right boundaries of $\mathscr{D}$, while the vertical crossing estimate takes the form,
   
     \begin{align*}
        \textbf{P}^{\xi^{\mathrm{Sloped}}}_{\mathscr{D}}  [   \mathcal{V}^{h \leq (1-c)k}       ]  \leq     C_{\mathcal{V}}       \text{ } \text{ , } \text{ } \\
   \end{align*}
   
   \noindent for intersecting left and right boundaries of $\mathscr{D}$.
    
    \bigskip
    
  \noindent \textit{Proof of Lemma 2.6}. Abbreviate $ \mathcal{V}^{h \leq (1-c)k}_{[-m,m] \times [0,n^{\prime} N]}( \mathscr{F}_1 ,\cdots, \mathscr{F}_n, FV) \equiv \mathcal{V}^{h \leq (1-c)k}_{[-m,m] \times [0,n^{\prime} N]}(FV) \equiv \mathcal{V}^{h \leq (1-c)k}$ for a vertical crossing that is a function of the height function over $n$ faces over a finite strip volume $FV$. 
    
    \bigskip
    
    \noindent \textbf{Case one} (\textit{non-intersecting left and right boundaries of the strip domain}). We have,
    
    \begin{align*}
   \text{ }      \textbf{P}^{\xi_{\mathrm{Sloped}}}_{[-m,m]\times[0,n^{\prime}N]} \big[      \mathcal{V}^{h \leq (1-c)k}      \big]   =    \text{ }    \textbf{P}^{\xi_{\mathrm{Sloped}}}_{[-m,m]\times[0,n^{\prime}N]} \bigg[   \big\{    \gamma_L \overset{h\geq ck}{\leftrightarrow} \mathscr{F}\mathscr{C}_1       \big\}     \cap      \big\{       \mathscr{F}\mathscr{C}_2      \overset{h\geq ck}{\not\leftrightarrow}    \mathscr{F}\mathscr{C}_{N-1}        \big\}  \cap   \big\{     \mathscr{F}\mathscr{C}_{N}               \overset{h\geq ck}{\leftrightarrow}    \gamma_R      \big\}   \bigg]  \\ \text{ }  \overset{(\mathrm{FKG})}{\geq}           \textbf{P}^{\xi_{\mathrm{Sloped}}}_{[-m,m]\times[0,n^{\prime}N]} \big[     \gamma_L \overset{h\geq ck}{\leftrightarrow} \mathscr{F}\mathscr{C}_1  \big]  \textbf{P}^{\xi_{\mathrm{Sloped}}}_{[-m,m]\times[0,n^{\prime}N]} \big[ \mathscr{F}\mathscr{C}_2      \overset{h\geq ck}{\not\leftrightarrow}    \mathscr{F}\mathscr{C}_{N-1} \big]  \textbf{P}^{\xi_{\mathrm{Sloped}}}_{[-m,m]\times[0,n^{\prime}N]} \big[  \mathscr{F}\mathscr{C}_{N}               \overset{h\geq ck}{\leftrightarrow}    \gamma_R    \big]   \text{ } \text{ , } \text{ } 
   \end{align*}

   \noindent in which the second term can be bounded below by another application of (FKG) over countably many horizontal, nonintersecting crossings with \textit{freezing clusters},
   
   \begin{align*}
   \textbf{P}^{\xi_{\mathrm{Sloped}}}_{[-m,m]\times[0,n^{\prime}N]} \big[ \mathscr{F}\mathscr{C}_2      \overset{h\geq ck}{\not\leftrightarrow}    \mathscr{F}\mathscr{C}_{N-1} \big]  \overset{\mathrm{(FKG)}}{\geq}     \prod_{2 \leq i < i + 1 \leq N-1} \text{ } \textbf{P}^{\xi_{\mathrm{Sloped}}}_{[-m,m] \times [0,n^{\prime} N]} \big[   \mathscr{F}\mathscr{C}_{i} \overset{h \geq ck}{\not\leftrightarrow}   \mathscr{F}\mathscr{C}_{i+1}   \big]   \text{ } \end{align*}
   
   \noindent as the lower bound probability is equivalent to, for a subdomain $\mathscr{S} \equiv \text{ }      [ x_{\mathscr{L}_1} ,     x_{\mathscr{L}_2} ]  \times   [0,n^{\prime} N]    \text{ } \subset \text{ }      \text{ } [-m,m]   \times[0,n^{\prime} N]$ of the strip,
   
   \begin{align*}
   \prod_{2 \leq i < i + 1 \leq N-1} \text{ }     \textbf{P}^{\xi_{\mathrm{Sloped}}}_{[-m,m] \times [0,n^{\prime} N]} \big[        \exists \text{ adjacent}  \mathscr{F}\mathscr{C}_{i}, \mathscr{F}\mathscr{C}_{i+1}  \in \mathcal{F}\mathcal{C}  ,  F_1 , \cdots , F_n \in       \mathcal{C}_{\mathrm{env}} :                                            F_1    \underset{\mathscr{S}}{\overset{h \leq ck }{\longleftrightarrow}}     F_n   ,         \mathscr{F}\mathscr{C}_{i} \overset{h \geq ck}{\not\leftrightarrow}   \mathscr{F}\mathscr{C}_{i+1}                      \big]     \\ \geq \text{ }  \prod_{2 \leq i < i + 1 \leq N-1} \text{ } \textbf{P}^{\xi_{\mathrm{Sloped}}}_{[-m,m] \times [0,n^{\prime} N]}   \bigg[    \big\{       F_1    \underset{\mathscr{S}}{\overset{h \leq ck }{\longleftrightarrow}}     F_n     \big\} \cap \big\{     \mathscr{F}\mathscr{C}_{i} \overset{h \geq ck}{\not\leftrightarrow}   \mathscr{F}\mathscr{C}_{i+1}    \big\}    \bigg]  
  \\ \overset{(*)}{\geq}          \text{ }  \prod_{2 \leq i < i + 1 \leq N-1} \text{ } \textbf{P}^{\xi_{\mathrm{Sloped}}}_{[-m,m] \times [0,n^{\prime} N]}   \bigg[   \big\{         \mathscr{L}_1                    \overset{h \leq ck}{\longleftrightarrow}    F_{\mathrm{crit}}               \big\}       \cap     \big\{      F_{\mathrm{supcrit}}                        \overset{h \leq ck }{\longleftrightarrow}            \mathscr{L}_2       \big\}    \cap \big\{\mathscr{F}\mathscr{C}_{i} \overset{h \geq ck}{\not\leftrightarrow}   \mathscr{F}\mathscr{C}_{i+1}   \big\}     \bigg]            \text{ }  \\ \text{ } \overset{\mathrm{(FKG)}}{\geq} \text{ }     \prod_{2 \leq i < i + 1 \leq N-1} \text{ } \textbf{P}^{\xi_{\mathrm{Sloped}}}_{[-m,m] \times [0,n^{\prime} N]}   \big[   \mathscr{L}_1   \overset{h \leq ck}{\longleftrightarrow}              F_{\mathrm{crit}}     \big]  \text{ } \textbf{P}^{\xi_{\mathrm{Sloped}}}_{[-m,m] \times [0,n^{\prime} N]}   \big[         F_{\mathrm{supcrit}} \overset{h \leq ck}{\longleftrightarrow}      \mathscr{L}_2              \big]  \text{ } \\ \times  \textbf{P}^{\xi_{\mathrm{Sloped}}}_{[-m,m] \times [0,n^{\prime} N]}   \big[  \mathscr{F}\mathscr{C}_i   \overset{h \geq ck}{\not\leftrightarrow}      \mathscr{F}\mathscr{C}_{i+1}    \big]       \text{ }\text{ , } \tag{\textit{2.6.1}} \end{align*}
   
   \noindent in which in the sequence of inequalities above, we decompose the product of crossing events dependent on $i$ as two events, the first of which depends upon $\mathcal{C}_{\mathrm{env}}$, defined to be the values of $h$ on all faces outside of the adjacent \textit{freezing clusters}, while for the second event we examine the lack of existence of crossings of the height function of level at least $ck$. From the second to third inequality annotated with (*), we lower bound the probability in the previous inequality, which is dependent upon the intersection of two events occurring simultaneously, by the intersection of three events occurring simultaneously. 
   
   \bigskip

   \noindent We lower bound each portion of \textit{(2.6.1)} with the following series of results. The following estimates will be incorporated to obtain the desired bound for other crossing probabilities taken under $\textbf{P}^{\xi^{\mathrm{Sloped}}}_{[-m,m] \times [0 , n^{\prime}N]} \big[     \cdot   \big]$.
   
   \bigskip

   \noindent \textbf{Lemma} $\textit{2.6.1.1}$ (\textit{a lower bound for the first crossing probability of \textit{(2.6.1)}}). For the line $\mathscr{L}_1$ appearing before the first \textit{freezing cluster} in the finite-volume strip, one has the strictly positive lower bound,
   
   \begin{align*}
          \textbf{P}^{\xi^{\mathrm{Sloped}}}_{[-m,m]\times[0,n^{\prime}N]} \big[   \mathscr{L}_1   \overset{h \leq ck}{\longleftrightarrow}               F_{\mathrm{crit}}     \big]  \text{ } \geq \text{ }   \mathscr{C}_1        \text{ } \text{. }
   \end{align*}

   \bigskip
   
   \noindent \textit{Proof of Lemma 2.6.1.1}. The probability of the connectivity event between the first line and the critical \textit{freezing cluster} occurring will be analyzed through the following series of inequalities. Beginning with the probability given on the LHS of \textit{2.6.1.1}, for some index set $\mathcal{I}^{\prime\prime}$,
   
   \begin{align*}
     \textbf{P}^{\xi^{\mathrm{Sloped}}}_{[-m,m]\times[0,n^{\prime}N]} \big[   \mathscr{L}_1   \overset{h \leq ck}{\longleftrightarrow}               F_{\mathrm{crit}}   \big]  \geq   \textbf{P}^{\xi^{\mathrm{Sloped}}}_{[-m,m]\times[0,n^{\prime}N]} \big[   \underset{{i \in \mathcal{I}^{\prime\prime}}}{\bigcap} \text{ }\big\{ \mathscr{L}_i     \overset{h \leq ck}{\longleftrightarrow} \mathscr{L}_{i+1} \big\}    \big]  \\  \text{ } \overset{(\mathrm{FKG})}{\geq} \prod_{i \in \mathcal{I}^{\prime\prime}} \text{ } \textbf{P}^{\xi^{\mathrm{Sloped}}}_{[-m,m] \times [0,n^{\prime}N]} \big[    \mathscr{L}_i     \overset{h \leq ck}{\longleftrightarrow} \mathscr{L}_{i+1}        \big] \text{ } \\  \geq \text{ }  \mathscr{C}_1\big(|\mathcal{I}^{\prime\prime}|\big) \equiv \mathscr{C}_1  \text{ } \text{ , } \\
   \end{align*}
   
   \noindent where in the second inequality after the LHS probability, we introduce a series of connectivity events between lines $\mathscr{L}_i$ and $\mathscr{L}_{i+1}$, where the final line in the series intersects $F_{\mathrm{crit}}$. Hence we obtain the desired lower bound for the first probability. \boxed{}
   
   \bigskip
   
      \noindent \textbf{Lemma} $\textit{2.6.1.2}$ (\textit{a lower bound for the second crossing probability of \textit{(2.6.1)}}). For the line $\mathscr{L}_1$ appearing before the first \textit{freezing cluster} in the finite-volume strip, one has the strictly positive lower bound,
   
   \begin{align*}
        \textbf{P}^{\xi_{\mathrm{Sloped}}}_{[-m,m] \times [0,n^{\prime} N]}   \big[       F_{\mathrm{supcrit}} \overset{h \leq ck}{\longleftrightarrow}      \mathscr{L}_2            \big] \geq            \mathscr{C}_2     \text{ } \text{ , }  \\
   \end{align*}
   
   \bigskip
   
   \noindent \textit{Proof of Lemma 2.6.1.2}. Beginning with the probability on the LHS of \textit{2.6.1.2} we obtain the lower bound, for some index set $\mathcal{I}^{\prime\prime\prime}$,

   \begin{align*}
    \textbf{P}^{\xi_{\mathrm{Sloped}}}_{[-m,m] \times [0,n^{\prime} N]}   \big[     F_{\mathrm{supcrit}} \overset{h \leq ck}{\longleftrightarrow}      \mathscr{L}_2            \big]  \geq    \textbf{P}^{\xi_{\mathrm{Sloped}}}_{[-m,m] \times [0,n^{\prime} N]}   \bigg[ \bigcap_{i \in \mathcal{I}^{\prime\prime\prime}}  \bigg\{  \big\{      \mathscr{L}^i_{\mathrm{supcrit}}   \overset{h \leq ck}{\longleftrightarrow}      F_{\mathrm{supcrit}}  \big\}                 \cap      \big\{   F_{\mathrm{supcrit}}       \overset{h \leq ck}{\longleftrightarrow}        \mathscr{L}^i_2   \big\}    \bigg\}    \bigg]   \text{ , } \text{ } \\
   \end{align*}
   
   \noindent where $\mathscr{L}^i_{\mathrm{supcrit}} \equiv \mathscr{L}^i$ is a series of lines, indexed by $i$, each of which appears before $F_{\mathrm{supcrit}}$, and $\mathscr{L}^i_2$ is a series of lines appearing before $\mathscr{L}_2$. The lower bound provided above strictly dominates the probability, 
   
   \begin{align*}
   \text{ }   \prod_{i \in \mathcal{I}^{\prime\prime\prime}} \text{ } \textbf{P}^{\xi^{\mathrm{Sloped}}}_{[-m,m]\times[0,n^{\prime}N]} \bigg[                \big\{      \mathscr{L}^i_{\mathrm{supcrit}}   \overset{h \leq ck}{\longleftrightarrow}      F_{\mathrm{supcrit}}    \big\}             \cap       \big\{ F_{\mathrm{supcrit}}       \overset{h \leq ck}{\longleftrightarrow}        \mathscr{L}^i_2    \big\}         \bigg]      \text{ } \text{ , } \text{ } \\
   \end{align*}
   
   \noindent by (FKG), from which we observe that the probability above strictly dominates,

         \begin{figure}
\begin{align*}
\includegraphics[width=0.85\columnwidth]{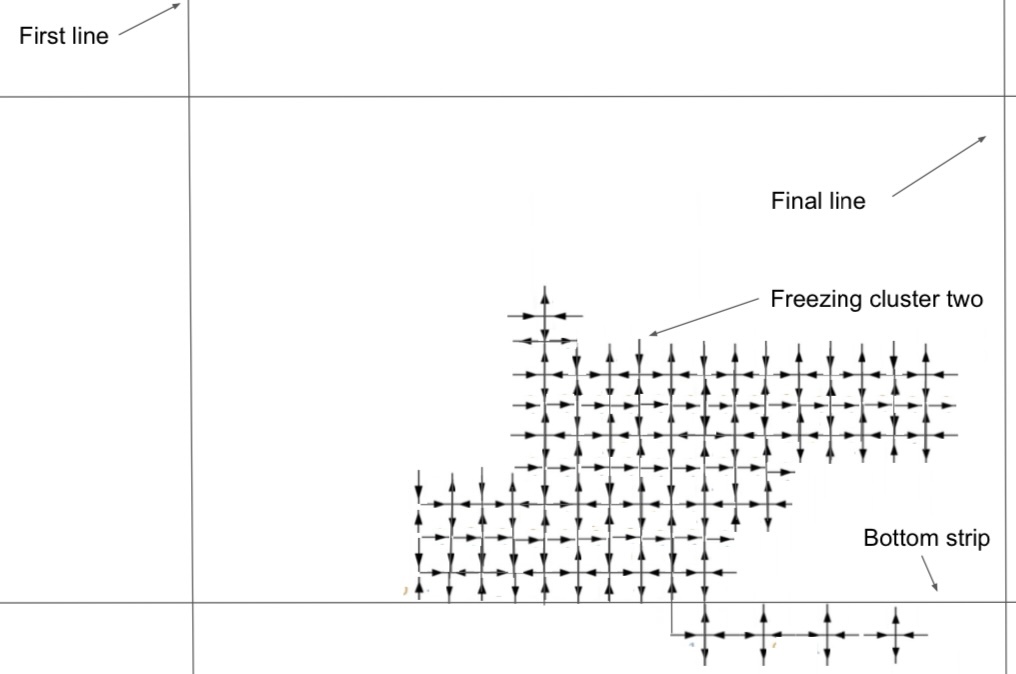}\\
\end{align*}
\caption{\textit{One sample of the random environment in the strip, with one freezing cluster depicted}. By situating vertical lines that intersect the top and bottom of the strip, arguments to quantify the crossing probability between series of lines, such as the first, and final, lines depicted above, are provided.}
\end{figure}

   \begin{align*}
         \prod_{i \in \mathcal{I}^{\prime\prime}} \text{ } \textbf{P}^{\xi^{\mathrm{Sloped}}}_{[-m,m] \times [0,n^{\prime} N]} \big[    \mathscr{L}^i_{\mathrm{supcrit}}   \overset{h \leq ck}{\longleftrightarrow}      F_{\mathrm{supcrit}}       \big]\textbf{P}^{\xi^{\mathrm{Sloped}}}_{[-m,m] \times [0,n^{\prime} N]} \big[    F_{\mathrm{supcrit}}       \overset{h \leq ck}{\longleftrightarrow}        \mathscr{L}^i_2      \big]   \text{ } \text{ , } \text{ } \\
   \end{align*}
   
   \noindent also by (FKG). Hence it suffices to bound, from below,

   \begin{align*}
     \textbf{P}^{\xi^{\mathrm{Sloped}}}_{[-m,m] \times [0,n^{\prime} N]} \big[     \mathscr{L}^i_{\mathrm{supcrit}}   \overset{h \leq ck}{\longleftrightarrow}      F_{\mathrm{supcrit}}     \big]                  \text{ } \text{ , }  \\
   \end{align*}

   \noindent and,

   \begin{align*}
            \textbf{P}^{\xi^{\mathrm{Sloped}}}_{[-m,m] \times [0,n^{\prime} N]} \big[     F_{\mathrm{supcrit}}       \overset{h \leq ck}{\longleftrightarrow}        \mathscr{L}^i_2      \big]            \text{ } \text{ . } \text{ }  \\
   \end{align*}

   \noindent For the first term, we observe that the crossing probability can be rearranged as, for $\mathscr{L}_{F_{\mathrm{supcrit}}} \equiv \mathscr{L}_F$,

   \begin{align*}
        \textbf{P}^{\xi^{\mathrm{Sloped}}}_{[-m,m] \times [0,n^{\prime} N]} \big[     \mathscr{L}^i_{\mathrm{supcrit}}   \overset{h \leq ck}{\longleftrightarrow}      F_{\mathrm{supcrit}}     \big]   \text{ }  \geq \text{ }   \textbf{P}^{\xi^{\mathrm{Sloped}}}_{[-m,m] \times [0,n^{\prime} N]} \bigg[   \bigcap_{i \in \mathcal{I}^{\prime\prime\prime}} \big\{ \mathscr{L}^i_{\mathrm{supcrit}}         \overset{h \leq ck}{\longleftrightarrow}       \mathscr{L}^i_F     \big\}      \bigg]\text{ }    \\ \overset{\mathrm{(FKG)}}{\geq} \prod_{i \in \mathcal{I}^{\prime\prime\prime}} \textbf{P}^{\xi^{\mathrm{Sloped}}}_{[-m,m] \times [0,n^{\prime} N]} \big[   \mathscr{L}^i_{\mathrm{supcrit}}         \overset{h \leq ck}{\longleftrightarrow}       \mathscr{L}^i_F                                 \big]                      \text{ } \text{ , } \text{ }  \\
   \end{align*}

   \noindent where $\mathscr{L}^i_F$ are a series of lines appearing in the strip before $ F_{\mathrm{supcrit}}$. Furthermore, the following lower bound from $\mathrm{(FKG)}$,
   
   \begin{align*}
         \textbf{P}^{\xi^{\mathrm{Sloped}}}_{[-m,m] \times [0,n^{\prime} N]} \big[ \mathscr{L}^i_{\mathrm{supcrit}}         \overset{h \leq ck}{\longleftrightarrow}       \mathscr{L}^i_F                             \big]   \geq          \textbf{P}^{\xi^{\mathrm{Sloped}}}_{[-m,m] \times [0,n^{\prime} N]} \big[    F_{\mathrm{supcrit}}          \overset{h \leq ck}{\longleftrightarrow}     \mathscr{L}^i_F          \big] \text{ }  \\  \geq  \textbf{P}^{\xi^{\mathrm{Sloped}}}_{[-m,m] \times [0,n^{\prime} N]} \bigg[  \big\{    F_{\mathrm{supcrit}}          \overset{h \leq ck}{\longleftrightarrow}     \mathscr{L}^i_F          \big\}    \cap  \big\{       F_{\mathrm{supcrit}} \overset{h \leq ck}{\not\longleftrightarrow}     \mathscr{L}^i_{\mathrm{supcrit}} \big\}  \cap  \big\{       \mathscr{L}^i_F      \overset{h \leq ck}{\longleftrightarrow}       \mathscr{L}^i_{\mathrm{supcrit}}   \big\}  \bigg] \text{ } \\ \overset{\mathrm{(FKG)}}{\geq} \text{ } \textbf{P}^{\xi^{\mathrm{Sloped}}}_{[-m,m] \times [0,n^{\prime} N]} \big[  F_{\mathrm{supcrit}}          \overset{h \leq ck}{\longleftrightarrow}     \mathscr{L}^i_F         \big] \text{ } \textbf{P}^{\xi^{\mathrm{Sloped}}}_{[-m,m] \times [0,n^{\prime} N]} \big[           F_{\mathrm{supcrit}} \overset{h \leq ck}{\not\longleftrightarrow}     \mathscr{L}^i_{\mathrm{supcrit}}                \big]  \\ \times  \textbf{P}^{\xi^{\mathrm{Sloped}}}_{[-m,m] \times [0,n^{\prime} N]} \big[          \mathscr{L}^i_F    \overset{h \leq ck}{\longleftrightarrow}   \mathscr{L}^i_{\mathrm{supcrit}}         \big]  \text{ } \end{align*}

         \noindent yields the lower bound,

         \begin{align*} \text{ }   \mathcal{C}_{\mathrm{supcrit, F}}  \mathcal{C}_{F,\mathrm{supcrit}}  \text{ } \textbf{P}^{\xi^{\mathrm{Sloped}}}_{[-m,m] \times [0,n^{\prime} N]} \big[            F_{\mathrm{supcrit}} \overset{h \leq ck}{\not\longleftrightarrow}     \mathscr{L}^i_{\mathrm{supcrit}}           \big]  \text{ }  \text{ , }   \tag{\textit{2.6.1.2 A}} \text{ }  \\
   \end{align*}
   
   \noindent where the two constants appearing in the lower bound, which are respectively given by $\mathcal{C}_{\mathrm{supcrit, F}}$, and $\mathcal{C}_{F,\mathrm{supcrit}}$, lower bound the probabilities of the connectivity events $\{ F_{\mathrm{supcrit}}          \overset{h \leq ck}{\longleftrightarrow}     \mathscr{L}^i_F   \}$, and $\{      \mathscr{L}^i_F    \overset{h \leq ck}{\longleftrightarrow}   \mathscr{L}^i_{\mathrm{supcrit}}     \}$, of occurring.
   
   \bigskip

   \noindent For the disconnective probability term, given an index set $\mathcal{I}^{\prime\prime\prime\prime}$,
   
   \begin{align*}
    \textbf{P}^{\xi^{\mathrm{Sloped}}}_{[-m,m] \times [0,n^{\prime} N]} \big[            F_{\mathrm{supcrit}} \overset{h \leq ck}{\not\longleftrightarrow}     \mathscr{L}^i_{\mathrm{supcrit}}                \big] \geq    \textbf{P}^{\xi^{\mathrm{Sloped}}}_{[-m,m] \times [0,n^{\prime} N]} \bigg[    \bigcap_{j \in \mathcal{I}^{\prime\prime\prime\prime}}         \big\{    F_{\mathrm{supcrit}}           \overset{h \leq ck}{\not\longleftrightarrow}      \mathscr{L}^j_{\mathrm{supcrit}}         \big\}    \bigg] \text{ } \\ \overset{\mathrm{(FKG)}}{\geq} \text{ }    \prod_{i \in \mathcal{I}^{\prime\prime\prime\prime}} \text{ } \textbf{P}^{\xi^{\mathrm{Sloped}}}_{[-m,m]\times[0,n^{\prime}N]}  \big[  F_{\mathrm{supcrit}}    \overset{h \leq ck}{\not\longleftrightarrow}       \mathscr{L}^j_{\mathrm{supcrit}}        \big]             \text{ } \text{ . } \text{ }   
    \end{align*}

    \noindent The last term in the inequality above can be bound below by a strictly positive constant, $c^{\prime}_i$, as,
    
    \begin{align*}
   \textbf{P}^{\xi^{\mathrm{Sloped}}}_{[-m,m]\times[0,n^{\prime}N]} \text{ } \big[  F_{\mathrm{supcrit}}    \overset{h \leq ck}{\not\longleftrightarrow}       \mathscr{L}^j_{\mathrm{supcrit}}        \big]  =             1  -  \textbf{P}   ^{\xi^{\mathrm{Sloped}}}_{[-m,m]\times[0,n^{\prime}N]} \text{ } \big[        F_{\mathrm{supcrit}}    \overset{h > ck}{\not\longleftrightarrow}_x       \mathscr{L}^j_{\mathrm{supcrit}}  ] \\ \geq 1 -      c^{\prime}_i       \text{ } \text{ , }    \tag{\textit{2.6.1.2 B}} \\
   \end{align*}
   
   \noindent in which the x-crossing occurs with strictly positive probability by finite energy, giving,
   
   \begin{align*}
    \prod_{i \in \mathcal{I}^{\prime\prime\prime\prime}}
   \textbf{P}^{\xi^{\mathrm{Sloped}}}_{[-m,m]\times[0,n^{\prime}N]} \text{ } \big[  F_{\mathrm{supcrit}}    \overset{h \leq ck}{\not\longleftrightarrow}       \mathscr{L}^j_{\mathrm{supcrit}}        \big]   \geq   \prod_{i \in \mathcal{I}^{\prime\prime\prime\prime}}  \big(  1 - c^{\prime}_i   \big)  \geq  \big(    1 -     c^{\prime}_{\mathrm{min},i}              \big)^{|\mathcal{I}^{\prime\prime\prime\prime}|}   \equiv    C \big(    |\mathcal{I}^{\prime\prime\prime\prime}|  ,  c^{\prime}_i  \text{ } \big) \equiv C
   \text{ } \text{ , } \text{ } \\
   \end{align*}
   
   \noindent for, $c^{\prime}_{\mathrm{min},i}  \text{ } \equiv \text{ }     \underset{i \in \mathcal{I}^{\prime\prime\prime}}{\mathrm{inf}}  \text{ }  c^{\prime}_i$, and a series of lines appearing before $F_{\mathrm{supcrit}}$ in the strip, which demonstrates that the below estimate holds,

   \begin{align*}
           \textbf{P}^{\xi^{\mathrm{Sloped}}}_{[-m,m] \times [0,n^{\prime} N]} \big[           \mathscr{L}^i_{\mathrm{supcrit}}   \overset{h \leq ck}{\longleftrightarrow}      F_{\mathrm{supcrit}}       \big]  \geq \mathcal{C}_{\mathrm{supcrit, F}}\mathcal{C}_{F,\mathrm{supcrit}}   \textbf{P}^{\xi^{\mathrm{Sloped}}}_{[-m,m] \times [0,n^{\prime} N]} \big[   F_{\mathrm{supcrit}} \overset{h \leq ck}{\not\longleftrightarrow}     \mathscr{L}^i_{\mathrm{supcrit}}          \big]                    \text{ } \text{ , } \text{ } \\
   \end{align*}

   \noindent as a result of $(\textit{2.6.1.2 A})$, and,
   
   \begin{align*}
   \textbf{P}^{\xi^{\mathrm{Sloped}}}_{[-m,m] \times [0,n^{\prime} N]} \big[            \mathscr{L}^i_{\mathrm{supcrit}}   \overset{h \leq ck}{\longleftrightarrow}      F_{\mathrm{supcrit}}      \big]  \geq      \mathcal{C}_{\mathrm{supcrit, F}} \mathcal{C}_{F,\mathrm{supcrit}}   \text{ }                   C \equiv   \mathscr{C}^{\prime}_3\big(  \mathcal{C}_{\mathrm{supcrit, F}} , \mathcal{C}_{F,\mathrm{supcrit}}  , C    \big) \equiv \mathscr{C}^{\prime}_3  \text{ } \text{ , } \text{ } \\
   \end{align*}
   
   \noindent as a result of $(\textit{2.6.1.2 B})$, for strictly positive $\mathscr{C}_3$. For the remaining term, observe, for an index set $\mathcal{I}^{\prime\prime\prime\prime\prime}$,
   
      \begin{align*}
          \textbf{P}^{\xi^{\mathrm{Sloped}}}_{[-m,m] \times [0,n^{\prime} N]} \big[    F_{\mathrm{supcrit}}       \overset{h \leq ck}{\longleftrightarrow}        \mathscr{L}^i_2     \big]      \geq      \textbf{P}^{\xi^{\mathrm{Sloped}}}_{[-m,m] \times [0,n^{\prime} N]} \bigg[      \bigcap_{l \in\mathcal{I}^{\prime\prime\prime\prime\prime}}    \big\{     F_{\mathrm{supcrit}}               \overset{h \leq ck}{\longleftrightarrow}      \mathscr{L}^l_2               \big\}                   \bigg]             \text{ }             \\ \text{ } \overset{\mathrm{(FKG)}}{\geq} \text{ } \prod_{l \in \mathcal{I}^{\prime\prime\prime\prime\prime}} \text{ } \underset{\equiv \text{ }   \sqrt[n^{\prime\prime\prime}]{ \mathscr{C}^{\prime\prime}_3 }   \text{ } }{\underset{ \text{ }  \underbrace{\geq \big(  \mathscr{C}^{\prime\prime}_3 \big(  | \mathcal{I}^{\prime\prime\prime\prime\prime}| \big) \big)^{\frac{1}{n^{\prime\prime\prime}}}} \text{ }             \text{ } }{\underset{     \underbrace{\geq ( c^{\prime\prime}_{\mathrm{min},l} )^{\frac{|\mathcal{I}^{\prime\prime\prime\prime\prime}|}{n^{\prime\prime\prime}}}}       \text{ } }{\underset{  }{ \underbrace{{\textbf{P}^{\xi^{\mathrm{Sloped}}}_{[-m,m] \times [0,n^{\prime} N]} \big[     F_{\mathrm{supcrit}}               \overset{h \leq ck}{\longleftrightarrow}      \mathscr{L}^l_2           \big]  }}}}}} \text{ }   \text{ , } \text{ } \\
   \end{align*}

   \noindent for a strictly positive number $n^{\prime\prime\prime}$ crossings, in turn yielding the estimate,
   
         \begin{align*}
        \prod_{l \in \mathcal{I}^{\prime\prime\prime\prime\prime}} \text{ }       c^{\prime\prime}_{\mathrm{min},l}            \geq    \big( c^{\prime\prime}_{\mathrm{min},l}  )^{|\mathcal{I}^{\prime\prime\prime\prime\prime}|}  \equiv     \mathscr{C}^{\prime\prime}_3\big(     |\mathcal{I}^{\prime\prime\prime\prime\prime}|    \big) \equiv \mathscr{C}^{\prime\prime}_3          \text{ } \text{ , } \text{ } \\
   \end{align*}

   \noindent for lines $\mathscr{L}^l_2$, indexed in $l$, appearing in the strip before $\mathscr{L}^i_2$, and,
   
  \begin{align*}
     c^{\prime\prime}_{\mathrm{min},l}   \text{ } = \text{ } \underset{l \in \mathcal{I}^{\prime\prime\prime\prime\prime}}{\mathrm{inf} }    c^{\prime\prime}_l      \text{ } \text{ , } \text{ } \\
  \end{align*} 
  
  \noindent for,
  
  \begin{align*}
 c^{\prime\prime}_l \text{ } =      \text{ }      \textbf{P}^{\xi^{\mathrm{Sloped}}}_{[-m,m] \times [0,n^{\prime} N]} \big[      F_{\mathrm{supcrit}}               \overset{h \leq ck}{\longleftrightarrow}      \mathscr{L}^l_2             \big]  \text{ }  \text{ , } \text{ } \\
  \end{align*}
  
\noindent   so that, $\sqrt[n]{c^{\prime\prime}_{\mathrm{min},l}}  \text{ } \text{ } \leq \text{ } c^{\prime\prime}_l \text{ }$, for $l \in \mathcal{I}^{\prime\prime\prime\prime\prime}$, and for a strictly positive $N$ crossings between $F_{\mathrm{supcrit}}$ and $\mathscr{L}^i_2$. Finally, in light of previously obtained estimates,
   
   \begin{align*}
     \prod_{i \in \mathcal{I}^{\prime\prime}} \text{ }   \textbf{P}^{\xi^{\mathrm{Sloped}}}_{[-m,m] \times [0,n^{\prime} N]} \big[     \mathscr{L}^i_{\mathrm{supcrit}}   \overset{h \leq ck}{\longleftrightarrow}      F_{\mathrm{supcrit}}       \big] \text{ }  \textbf{P}^{\xi^{\mathrm{Sloped}}}_{[-m,m] \times [0,n^{\prime} N]} \big[      F_{\mathrm{supcrit}}       \overset{h \leq ck}{\longleftrightarrow}        \mathscr{L}^i_2       \big]       \geq           \prod_{k,k^{\prime} \in \mathcal{I}^{\prime\prime}} \text{ }  \big( \mathscr{C}^{\prime}_3   \big)_{k}    \big(   \mathscr{C}^{\prime\prime}_3      \big)_{k^{\prime}}       \text{ } \\  \equiv   \bigg[  \prod_{k \in \mathcal{I}^{\prime\prime}}   \big(  \mathscr{C}^{\prime}_3    \big)_{k}   \bigg] \bigg[   \prod_{k^{\prime} \in \mathcal{I}^{\prime\prime}}   \big( \mathscr{C}^{\prime\prime}_3       \big)_{k^{\prime}}    \bigg] \text{ }   \\ \geq        \bigg[ \big(    \mathscr{C}^{\prime}_3        \big)_{k_m}  \big(   \mathscr{C}^{\prime\prime}_3 \big)_{k^{\prime}_m} \bigg]^{|
     \mathcal{I}^{\prime\prime}|}\text{ }  \\  \equiv  \bigg[  \big(   \mathscr{C}^{\prime}_{3,k_m}      \big)  \big(   \mathscr{C}^{\prime\prime}_{3,k^{\prime}_m}    \big) \bigg]^{|
     \mathcal{I}^{\prime\prime}|} \\ \equiv \bigg[ \mathscr{C}^{\prime,\prime\prime}_{3,k_m,k^{\prime}_m}               \bigg]^{|
     \mathcal{I}^{\prime\prime}|}     \text{ } \text{ , } \text{ } \tag{\textit{2.6.1.2 C}} \\
   \end{align*}
   
   \noindent where,
   
   \begin{align*}
  \big(    \mathscr{C}^{\prime}_3  \big)_{k_m} \equiv      \mathscr{C}^{\prime}_{3,k_m}  = \underset{{k \in    \mathcal{I}^{\prime\prime\prime\prime}      }\subset \mathcal{I}^{\prime\prime} }{\mathrm{inf}} \text{ } \big(   \mathscr{C}^{\prime}_3  \big)_k  \equiv    \underset{{k \in    \mathcal{I}^{\prime\prime\prime\prime}      }\subset \mathcal{I}^{\prime\prime} }{\mathrm{inf}}  \text{ }  \mathscr{C}^{\prime}_{3,k}             \text{ , } \text{ } \\
   \end{align*}
   
   \noindent and,
   
   \begin{align*}
   \big(     \mathscr{C}^{\prime\prime}_3 \big)_{k^{\prime}_m}   \equiv       \mathscr{C}^{\prime\prime}_{3,k^{\prime}_m}      = \underset{\text{ }    k^{\prime} \in         \mathcal{I}^{\prime\prime\prime\prime\prime}        \subset \mathcal{I}^{\prime\prime}}{\mathrm{inf}} \text{ }       \big(  \mathscr{C}^{\prime\prime}_3  \big)_{k^{\prime}}    \equiv        \underset{\text{ }    k^{\prime} \in        \mathcal{I}^{\prime\prime\prime\prime\prime}        \subset \mathcal{I}^{\prime\prime}}{\mathrm{inf}} \text{ }        \mathscr{C}^{\prime\prime}_{3,{k^{\prime}}}          \text{ } \text{ . } \\
   \end{align*}
   
   \noindent Hence we conclude the argument, setting $\mathscr{C}^{\emptyset}_2 = \mathscr{C}^{\prime,\prime\prime}_{3,k_m}$ from $(\textit{2.6.1.2 C})$. \boxed{}

\bigskip

\noindent From \textit{2.6.1.1} and \textit{2.6.1.2}, we obtain the final estimate with \textit{2.6.1.3} below.

   \bigskip
   
      \noindent \textbf{Lemma} $\textit{2.6.1.3}$ (\textit{a lower bound for the third crossing probability of \textit{(2.6.1)}}). For the line $\mathscr{L}_1$ appearing before the first \textit{freezing cluster} in the finite-volume strip, one has the strictly positive lower bound,
    
   \begin{align*}
        \textbf{P}^{\xi_{\mathrm{Sloped}}}_{[-m,m] \times [0,n^{\prime} N]}   \big[   \mathscr{F}\mathscr{C}_i   \overset{h \geq ck}{\not\leftrightarrow}      \mathscr{F}\mathscr{C}_{i+1}      \big]   \geq \mathscr{C}_3    \text{ } \text{ . } \text{ } \\
   \end{align*}

   \bigskip
   
   \noindent \textit{Proof of Lemma 2.6.1.3}. Beginning with the probability on the LHS of \textit{2.6.1.3} to obtain the lower bound, write,
   
   \begin{align*}
     \textbf{P}^{\xi_{\mathrm{Sloped}}}_{[-m,m] \times [0,n^{\prime} N]}   \big[   \mathscr{F}\mathscr{C}_i   \overset{h \geq ck}{\not\leftrightarrow}      \mathscr{F}\mathscr{C}_{i+1}   \big]  =  1 - \textbf{P}^{\xi_{\mathrm{Sloped}}}_{[-m,m] \times [0,n^{\prime} N]}   \big[    \mathscr{F}\mathscr{C}_i   \overset{h < ck}{\not\leftrightarrow}_x      \mathscr{F}\mathscr{C}_{i+1}     \big]    \text{ } \text{ , } \text{ } \\
   \end{align*}
   
   \noindent allowing one to conclude that,
   
   \begin{align*}
    \textbf{P}^{\xi_{\mathrm{Sloped}}}_{[-m,m] \times [0,n^{\prime} N]}   \big[    \mathscr{F}\mathscr{C}_i   \overset{h < ck}{\not\leftrightarrow}_x      \mathscr{F}\mathscr{C}_{i+1}      \big]   \geq \textbf{P}^{\xi^{\mathrm{Sloped}}}_{[-m,m] \times [0,n^{\prime} N]} \big[     \bigcap_{l \in \textbf{N}: l \in \mathcal{J}} \text{ } \big\{         \mathscr{F}\mathscr{C}_i  \overset{h < ck}{\not\longleftrightarrow}_x          \mathscr{L}^{i+1}_l    \big\}      \big] \text{ }  \\ \overset{\mathrm{(FKG)}}{\geq} \text{ } \prod_{l \in \textbf{N}: l \in \mathcal{J}} \textbf{P}^{\xi^{\mathrm{Sloped}}}_{[-m,m] \times [0,n^{\prime} N]} \big[  \mathscr{F}\mathscr{C}_i   \overset{h < ck}{\not\longleftrightarrow}_x \mathscr{L}^{i+1}_l       \big] \\ \geq   \prod_{l \in \textbf{N}: l \in \mathcal{J}} \text{ } c^{i,i+1}_l         \\  \geq       \big( \underset{l \in \textbf{N}: l \in \mathcal{J}}{\mathrm{inf}} \text{ } c^{i,i+1}_l     \big)^{|\mathcal{J}|} \text{ } \\ \equiv   C^{i,i+1}  \text{ } \text{ , } \text{ } \tag{\textit{2.6.1.3.1}}\\
   \end{align*}
   
   \noindent for some index set $\mathcal{J}$, and a series of lines $\mathscr{L}^{i+1}_l$ in the strip appearing before $\mathscr{F}\mathscr{C}_{i+1}$. For the lower bound of each crossing event, which takes the form,
   
   \begin{align*}
     \textbf{P}^{\xi^{\mathrm{Sloped}}}_{[-m,m] \times [0,n^{\prime} N]} \big[   \mathscr{F}\mathscr{C}_i   \overset{h < ck}{\not\longleftrightarrow}_x   \mathscr{L}^{i+1}_l            \big]   \text{ } \text{ , } \text{ } \\
   \end{align*}
   
   \noindent each occurs by finite energy with strictly positive probability for every $l$ in the index set; each crossing probability is bound below by $c^{i,i+1}_l$. Hence the desired constant $\mathscr{C}^{\emptyset}_3 \equiv C^{\text{ } i,i+1}$, appearing in the expression for the lower bound from $(\textit{2.6.1.3.1})$. \boxed{}

         \begin{figure}
\begin{align*}
\includegraphics[width=0.82\columnwidth]{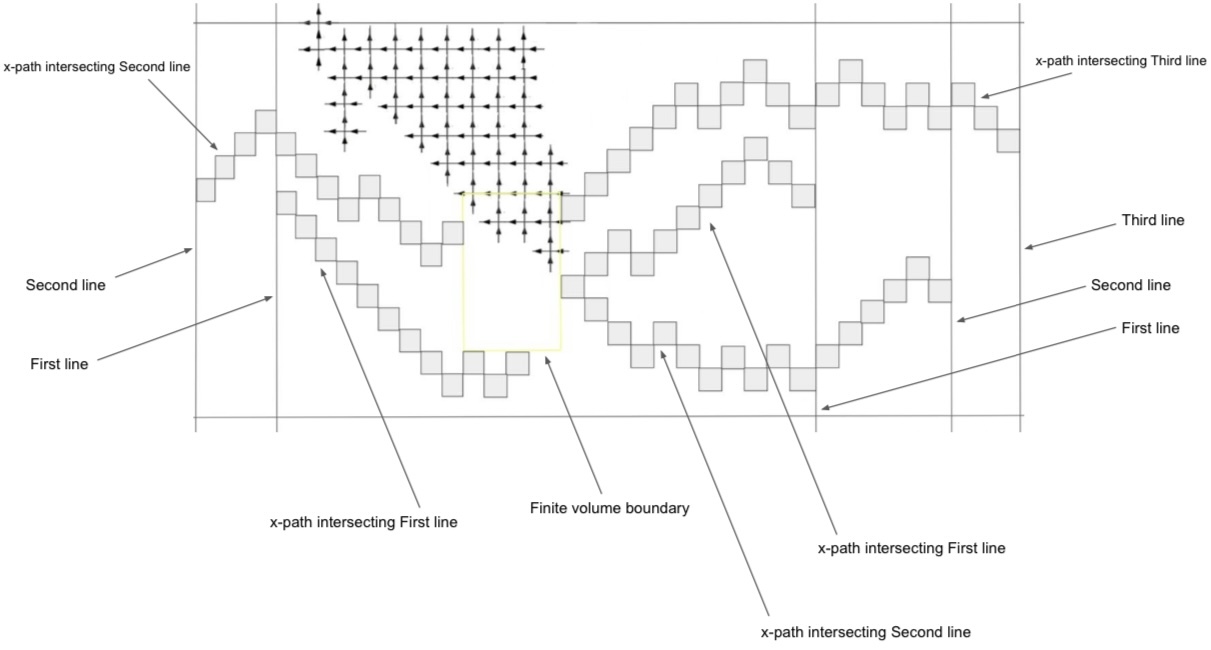}\\
\end{align*}
\caption{\textit{A finite volume configuration with yellow finite subvolume surrounding a boundary of the freezing cluster in the strip}. Relating to arguments in the previous \textbf{Lemma}, two lines are positioned to the left of the freezing cluster, with horizontal $\mathrm{x}$-paths intersectingi the yellow finite volume boundary and each of the two lines. On the other side, there are three $\mathrm{x}$-paths, each of which respectively intersect three vertically placed lines to the right of the freezing cluster.}
\end{figure}
   
   \bigskip
   
   \noindent \textbf{Lemma }$\textit{2.6.1.4}$  (\textit{comparison of disconnectivity between adjacent freezing clusters with connectivity across the complement of freezing clusters in the strip}). The disconnective probability between adjacent \textit{freezing clusters} admits the lower bound,
   
   \begin{align*}
        \textbf{P}^{\xi^{\mathrm{Sloped}}}_{[-m,m] \times [0,n^{\prime}N]} \big[       \mathscr{F}\mathscr{C}_i \overset{h \geq ck}{\not\longleftrightarrow} \mathscr{F}\mathscr{C}_{i+1}   \big]     \geq    \textbf{P}^{\xi^{\mathrm{Sloped}}}_{[-m,m] \times [0,n^{\prime}N]} \bigg[  \big\{         \mathscr{L}_i   \overset{h \leq ck}{\longleftrightarrow}               (F_{\mathrm{crit}})_i            \big\}   \cap     \big\{     (F_{\mathrm{supcrit}})_i \overset{h \leq ck}{\longleftrightarrow}      \mathscr{L}_2             \big\}   \\ \cap    \big\{  \mathscr{F}\mathscr{C}_i     \underset{\textbf{Z}^2 \cap ( \mathscr{F} \mathscr{C})^c_i}{\overset{h \leq ck}{\longleftrightarrow}}    \mathscr{F}\mathscr{C}_{i+1}               \big\}         \bigg]   \\ \geq   \big( \mathscr{C}_4  \big)_i     \text{ } \text{ , } \text{ }  \\
   \end{align*}
   
   \noindent for every $i$, $\big( \text{ } \mathscr{C}_4 \text{ }  \big)_i > 0$, and critical and supercritical faces surrounding in the the strip between $\mathscr{F}\mathscr{C}_{i-1}$ and $\mathscr{F}\mathscr{C}_i$, respectively given by $(F_{\mathrm{crit}})_i$ and $(F_{\mathrm{supcrit}})_i $.

   \bigskip
   
           \begin{figure}
\begin{align*}
\includegraphics[width=0.88\columnwidth]{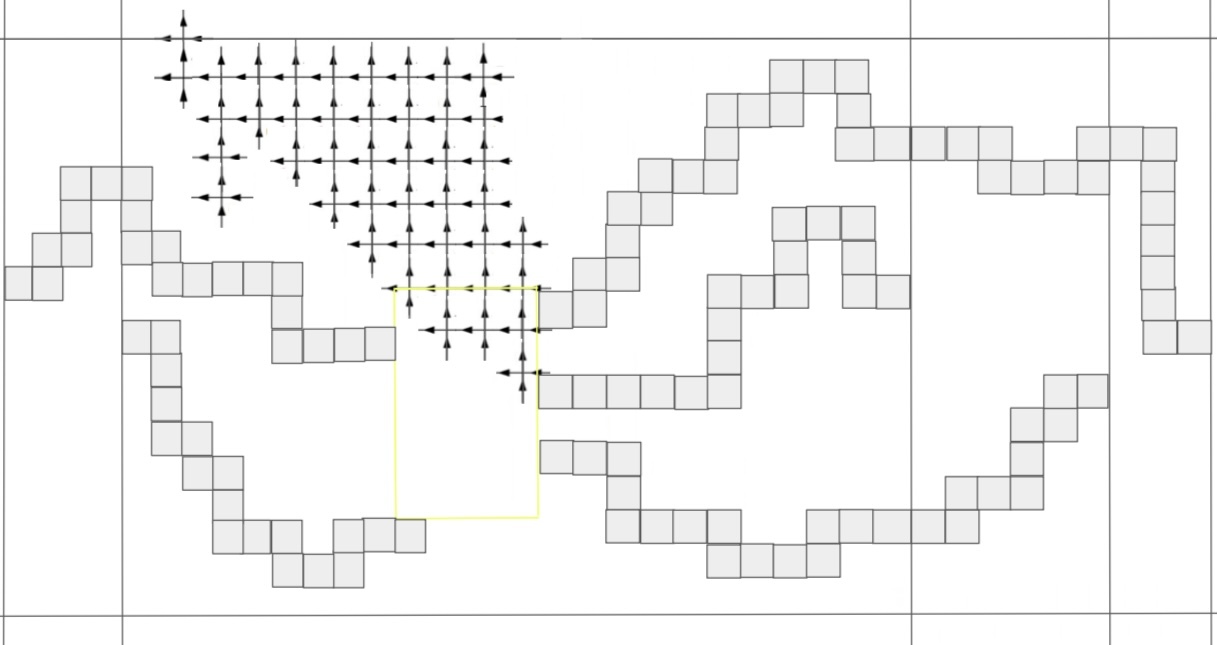}\\
\end{align*}
\caption{\textit{A depiction of crossings, in comparison to \text{x}-crossings, provided in the previous figure}. For regular crossings, the requirement that the distance between neighboring faces of the path be $1$, instead of $2$, can lead to differences in the number of faces contained in paths which intersect lines to the left, or to the right, of the \textit{freezing cluster}.}
\end{figure}

   \noindent \textit{Proof of Lemma 2.6.1.4}. We decompose the crossing event and produce a lower bound estimate of the disconnectivity event between adjacent \textit{freezing clusters},
   
   \begin{align*}
            \textbf{P}^{\xi^{\mathrm{Sloped}}}_{[-m,m] \times [0,n^{\prime}N]}\big[      \mathscr{F}\mathscr{C}_i \overset{h \geq ck}{\not\longleftrightarrow} \mathscr{F}\mathscr{C}_{i+1}     \big] \geq \textbf{P}^{\xi^{\mathrm{Sloped}}}_{[-m,m] \times [0,n^{\prime}N]} \bigg[     \underset{i \in \mathcal{I}}{\bigcap}  \bigg\{ \big\{ \text{ }        \mathscr{F}\mathscr{C}_i       \overset{h \geq ck}{\longleftrightarrow}      \mathscr{L}^{\prime}  \big\}  \cap  \big\{   \mathscr{L}^{\prime\prime} \overset{h \geq ck}{\not\longleftrightarrow }   \mathscr{F}\mathscr{C}_{i+1}    \big\}   \bigg\}       \bigg] \\ \overset{\mathrm{(FKG)}}{\geq}  \textbf{P}^{\xi^{\mathrm{Sloped}}}_{[-m,m] \times [0,n^{\prime}N]}  \big[  \mathscr{F}\mathscr{C}_i       \overset{h \geq ck}{\longleftrightarrow}      \mathscr{L}^{\prime} \big]   \textbf{P}^{\xi^{\mathrm{Sloped}}}_{[-m,m] \times [0,n^{\prime}N]}  \big[        \mathscr{L}^{\prime\prime} \overset{h \geq ck}{\not\longleftrightarrow }   \mathscr{F}\mathscr{C}_{i+1}       \big] \text{ , } \text{ }  \\
   \end{align*}
   
     \noindent implying that it suffices to provide a lower bound for,

     \begin{align*}
         \textbf{P}^{\xi^{\mathrm{Sloped}}}_{[-m,m] \times [0,n^{\prime}N]} \big[     \mathscr{F}\mathscr{C}_i       \overset{h \geq ck}{\longleftrightarrow}      \mathscr{L}^{\prime}       \big]       \text{ } \text{ , } \text{ }  \\
   \end{align*}

   \noindent with,
   
   \begin{align*}
   \textbf{P}^{\xi^{\mathrm{Sloped}}}_{[-m,m] \times [0,n^{\prime}N]} \big[          \mathscr{L}_i   \overset{h \leq ck}{\longleftrightarrow}               (F_{\mathrm{crit}})_i      \big]     \text{ } \text{ , } \text{ } \\
   \end{align*}

   \noindent and also, 
   
   \begin{align*}
   \textbf{P}^{\xi^{\mathrm{Sloped}}}_{[-m,m] \times [0,n^{\prime}N]} \big[           \mathscr{L}^{\prime\prime} \overset{h \geq ck}{\not\longleftrightarrow }   \mathscr{F}\mathscr{C}_{i+1}      \big]     \text{ } \text{ , } \text{ } \\
   \end{align*}

   \noindent with,
   
   \begin{align*}
   \textbf{P}^{\xi^{\mathrm{Sloped}}}_{[-m,m] \times [0,n^{\prime}N]} \bigg[ \big\{    (F_{\mathrm{supcrit}})_i \overset{h \leq ck}{\longleftrightarrow}      \mathscr{L}_2     \big\}    \cap  \big\{   \mathscr{F}\mathscr{C}_i     \underset{\textbf{Z}^2 \cap ( \mathscr{F} \mathscr{C})^c_i}{\overset{h \leq ck}{\longleftrightarrow}}    \mathscr{F}\mathscr{C}_{i+1}             \big\}     \bigg] \text{ } \text{ }  \text{ , } \text{ } \\
   \end{align*}
   
   \noindent because, by $\mathrm{(FKG)}$, the lower bound given in the statement of \textbf{Lemma} $\textit{2.6.1.4}$,
  
  \begin{align*}
     \text{ }  \textbf{P}^{\xi^{\mathrm{Sloped}}}_{[-m,m] \times [0,n^{\prime}N]} \bigg[ \big\{   (F_{\mathrm{supcrit}})_i \overset{h \leq ck}{\longleftrightarrow}      \mathscr{L}_2       \text{ } \big\}    \cap    \big\{ \mathscr{F}\mathscr{C}_i     \underset{\textbf{Z}^2 \cap ( \mathscr{F} \mathscr{C})^c_i}{\overset{h \leq ck}{\longleftrightarrow}}    \mathscr{F}\mathscr{C}_{i+1}        \text{ }        \big\}   \bigg] \text{ } \text{ }  \text{ . } \text{ } \\
   \end{align*}
  
  \noindent For the first comparison, we begin with the crossing,

  \begin{align*}
 \textbf{P}^{\xi^{\mathrm{Sloped}}}_{[-m,m] \times [0,n^{\prime} N ]} \big[  \mathscr{L}^{\prime}        \overset{h \geq ck}{\longleftrightarrow}         \mathscr{F}\mathscr{C}_i       \big]  \text{ } \text{ , } \text{ } \\
  \end{align*}


  \noindent which is compared to the lower bound of the connectivity event between $\mathscr{L}_i$ and $(F_{\mathrm{crit}})_i$, which from $\mathrm{(FKG)}$ can be expressed as,

  \begin{align*}
        \textbf{P}^{\xi^{\mathrm{Sloped}}}_{[-m,m] \times [0,n^{\prime}N]} \big[           \mathscr{L}_i   \overset{h \leq ck}{\longleftrightarrow}               (F_{\mathrm{crit}})_i      \big]    \geq     \textbf{P}^{\xi^{\mathrm{Sloped}}}_{[-m,m] \times [0,n^{\prime}N]} \bigg[  \bigg\{   \bigcap_{o \in \mathcal{I}} \big\{   \mathscr{L}_i   \overset{h \leq ck}{\longleftrightarrow}          \mathscr{L}^{\prime}_o  \big\}  \bigg\}  \cap    \big\{       \mathscr{L}_i   \overset{h \leq ck}{\longleftrightarrow}               (F_{\mathrm{crit}})_i \big\}     \bigg]  \\ \overset{\mathrm{(FKG)}}{\geq}  \textbf{P}^{\xi^{\mathrm{Sloped}}}_{[-m,m] \times [0,n^{\prime}N]} \big[   \bigcap_{o \in \mathcal{I}} \big\{    \mathscr{L}_i   \overset{h \leq ck}{\longleftrightarrow}          \mathscr{L}^{\prime}_o  \big\}                                \big]   \textbf{P}^{\xi^{\mathrm{Sloped}}}_{[-m,m] \times [0,n^{\prime}N]}   \big[     \mathscr{L}_i   \overset{h \leq ck}{\longleftrightarrow}               (F_{\mathrm{crit}})_i                    \big]   \\ \overset{\mathrm{(FKG)}}{\geq}  \bigg[   \prod_{o  \in \mathcal{I} : o \cap i \neq \emptyset } \textbf{P}^{\xi^{\mathrm{Sloped}}}_{[-m,m] \times [0,n^{\prime}N]}  \big[     \mathscr{L}_i   \overset{h \leq ck}{\longleftrightarrow}          \mathscr{L}^{\prime}_o  \big]   \bigg]                    \textbf{P}^{\xi^{\mathrm{Sloped}}}_{[-m,m] \times [0,n^{\prime}N]}   \big[      \mathscr{L}_i   \overset{h \leq ck}{\longleftrightarrow}               (F_{\mathrm{crit}})_i                       \big]              \text{ } \text{ . } \text{ } \\
  \end{align*}
  
  \noindent In turn, part of the comparison for the first term follows from the estimate,
  
  \begin{align*}
       \textbf{P}^{\xi^{\mathrm{Sloped}}}_{[-m,m] \times [0,n^{\prime}]} \big[ \mathscr{L}^{\prime}        \overset{h \geq ck}{\longleftrightarrow}         \mathscr{F}\mathscr{C}_i                            \big] \geq      \textbf{P}^{\xi^{\mathrm{Sloped}}}_{[-m,m] \times [0,n^{\prime}]}  \big[    \bigcap_{i^{\prime\prime}=1}^N    \big\{ \mathscr{L}^{\prime}_{i^{\prime}}                        \overset{h \geq ck}{\longleftrightarrow}       \mathscr{L}^{\prime}_{i^{\prime}+1}              \big\}                                    \big]        \text{ } \text{ , } \text{ } \\
  \end{align*}
  
  \noindent where $\mathscr{L}^{\prime}_N \cap \mathscr{F}\mathscr{C}_i \neq \emptyset$ for $N$ large enough, because, by $\mathrm{(FKG)}$,

  \begin{align*}
     \prod_{i^{\prime}=1}^N \text{ } \textbf{P}^{\xi^{\mathrm{Sloped}}}_{[-m,m] \times [0,n^{\prime}N]} \big[   \mathscr{L}^{\prime}_{i^{\prime}}         \overset{h \geq ck}{\longleftrightarrow}   \mathscr{L}^{\prime}_{i^{\prime}+1}         \big]   \text{ } \text{ . } \text{ } \\
  \end{align*}
  
  \noindent The above countable product of probabilities admits a lower bound,
  
   \begin{align*}
  C^{\prime}_{i,i+1}  \equiv    \big(          c^{\prime}_{i,i+1}           \big)^N        \text{ } \text{ , } \text{ } \\
  \end{align*}

  \noindent further demonstrating that the remaining term for the first crossing probability can be analyzed with the following,

  \begin{align*}
     \textbf{P}^{\xi^{\mathrm{Sloped}}}_{[-m,m] \times [0,n^{\prime}N]}   \big[     \mathscr{L}_i   \overset{h \leq ck}{\longleftrightarrow}               (F_{\mathrm{crit}})_i                   \big]      \text{ } \overset{\mathrm{(FKG)}}{\geq} \text{ }    \prod_{i^{\prime\prime}=1}^{N^{\prime}} \textbf{P}^{\xi^{\mathrm{Sloped}}}_{[-m,m] \times [0,n^{\prime}N]}   \big[                       \mathscr{L}_{i^{\prime\prime}}           \overset{h \leq ck}{\longleftrightarrow}   \mathscr{L}_{i^{\prime\prime}+1}      \big] \\ {\geq}   \prod_{i^{\prime\prime}=1}^{N^{\prime}} \textbf{P}^{\xi^{\mathrm{Sloped}}}_{[-m,m] \times [0,n^{\prime}N]}   \bigg[ \big\{    \mathscr{L}_{i^{\prime\prime}}           \overset{h \leq ck}{\longleftrightarrow}   \mathscr{L}_{i^{\prime\prime}+1}          \big\}   \cap    \big\{   \mathscr{L}_{i^{\prime\prime}+1} \text{ }             \overset{h \leq ck}{\longleftrightarrow}    ( F_{\mathrm{crit}})_i             \big\}         \bigg]  \text{ } \\ \overset{\mathrm{(FKG)}}{\geq} \text{ }    \prod_{i^{\prime\prime}=1}^{N^{\prime}} \text{ }                \textbf{P}^{\xi^{\mathrm{Sloped}}}_{[-m,m] \times [0,n^{\prime}N]}   \big[      \mathscr{L}_{i^{\prime\prime}}           \overset{h \leq ck}{\longleftrightarrow}   \mathscr{L}_{i^{\prime\prime}+1}        \big]     \textbf{P}^{\xi^{\mathrm{Sloped}}}_{[-m,m] \times [0,n^{\prime}N]}   \big[   \mathscr{L}_{i^{\prime\prime}+1}        \overset{h \leq ck}{\longleftrightarrow}    ( F_{\mathrm{crit}})_i           \big]      \text{ } \text{ , } \text{ } \\
  \end{align*}

  \noindent for $N^{\prime}$ large enough, in addition to two estimates for each probability in the product, including,
  
  \begin{align*}
             \textbf{P}^{\xi^{\mathrm{Sloped}}}_{[-m,m] \times [0,n^{\prime}N]}   \big[       \mathscr{L}_{i^{\prime\prime}}           \overset{h \leq ck}{\longleftrightarrow}   \mathscr{L}_{i^{\prime\prime}+1}            \big] \geq    \textbf{P}^{\xi^{\mathrm{Sloped}}}_{[-m,m] \times [0,n^{\prime}N]}   \big[      \mathscr{L}_{i^{\prime\prime}}           \overset{h \leq ck}{\longleftrightarrow}   \mathscr{L}_{i^{\prime\prime\prime}+1}           \big]    \text{ } \overset{\mathrm{(FKG)}}{\geq}       C^1_{i^{\prime\prime},i^{\prime\prime\prime}}   \equiv   C^1_{\prime\prime,\prime\prime\prime}       \text{ } \text{ , }   \tag{\textit{2.6.1.4.1}}        \text{ } \\
  \end{align*}
  
  \noindent for $i^{\prime\prime\prime} > i^{\prime\prime}$, with $i^{\prime\prime\prime} >0$ in the first probability, which is further rearranged from a countable intersection of crossing events to obtain the final lower bound $C_{\prime\prime,\prime\prime\prime}$, as,

  \begin{align*}
           C^1_{\prime\prime,\prime\prime\prime}  \leq        \big(    c^1_{\prime\prime,\prime\prime\prime}  \big)^{|\mathcal{I}_{\prime\prime,\prime\prime\prime}|}        \leq  \prod_{k\in\mathcal{I}_{\prime\prime,\prime\prime\prime}}     \textbf{P}^{\xi^{\mathrm{Sloped}}}_{[-m,m] \times [0,n^{\prime}N]}    \big[                             \mathscr{L}_{i^{\prime\prime}}     \overset{h \leq ck}{\longleftrightarrow}      \mathscr{L}_k        \big]                             \text{ } \text{ , } \text{ } \\
  \end{align*}

  \noindent for some index set $\mathcal{I}_{\prime\prime,\prime\prime\prime}$, yielding the estimate,

  \begin{align*}
  \textbf{P}^{\xi^{\mathrm{Sloped}}}_{[-m,m] \times [0,n^{\prime}N]}   \big[    \mathscr{L}_{i^{\prime\prime}}           \overset{h \leq ck}{\longleftrightarrow}   \mathscr{L}_{i^{\prime\prime}+1}             \big]   \geq     C^1_{\prime\prime,\prime\prime\prime}                \text{ }\text{ , } \\
  \end{align*}

  \noindent  implying,

  \begin{align*}
        \textbf{P}^{\xi^{\mathrm{Sloped}}}_{[-m,m] \times [0,n^{\prime}N]}   \big[   \mathscr{L}_{i^{\prime\prime}+1} \text{ }             \overset{h \leq ck}{\longleftrightarrow}    ( F_{\mathrm{crit}})_i         \big]    \overset{\mathrm{(FKG)}}{\geq}  \prod_{i^{\prime\prime} \in \mathcal{I}}   \textbf{P}^{\xi^{\mathrm{Sloped}}}_{[-m,m] \times [0,n^{\prime}N]}   \big[       \mathscr{L}_{i^{\prime\prime}+1}                     \overset{h \leq ck}{\longleftrightarrow}        \mathscr{L}^{\prime}_{i^{\prime\prime}+1}                            \big]    \\  \geq    \prod_{i^{\prime\prime} \in \mathcal{I}}                 C^1_{\prime\prime,\prime\prime\prime}  \geq       c^1_{\prime\prime,\prime\prime\prime}  \big(   |   \mathcal{I} |     \big)            \text{ } \text{ , } \text{ }    \\
  \end{align*}
  
  \noindent with the above lower bound explicitly taking the form,

  \begin{align*}
       c^1_{\prime\prime,\prime\prime\prime} \big(  |   \mathcal{I} |  \big)   \equiv             \big(   C^1_{\prime\prime,\prime\prime\prime}  \big)^{|\mathcal{I}|}                 \text{ }    \tag{\textit{2.6.1.4.2}} \text{ , } \\
  \end{align*}
  
  \noindent in the second probability, where the lines satisfy the condition in which there are finitely many of which have nonempty intersection with $( F_{\mathrm{crit}})_i$,
  
  \begin{align*}
      \mathscr{L}^{\prime}_{i^{\prime\prime}+1}        \cap  ( F_{\mathrm{crit}})_i     \neq  \emptyset \text{ , } \text{ } \\
  \end{align*}

  \noindent yielding the estimate,

  \begin{align*}
           \textbf{P}^{\xi^{\mathrm{Sloped}}}_{[-m,m] \times [0,n^{\prime}N]}   \big[  \mathscr{L}_{i^{\prime\prime}+1} \text{ }             \overset{h \leq ck}{\longleftrightarrow}    ( F_{\mathrm{crit}})_i         \big]     \geq    C^2_{\prime\prime,\prime\prime\prime}      \text{ } \text{ , } \text{ } \tag{\textit{2.6.1.4.3}}\\
  \end{align*}
  
  \noindent where as for $C^1_{\prime\prime,\prime\prime\prime}$, $C^2_{\prime\prime,\prime\prime\prime}$ is obtained from the countable intersection of crossing events given below, 
  
  \begin{align*}
        C^2_{\prime\prime,\prime\prime\prime} \text{ } \leq  \big(  c^2_{\prime\prime,\prime\prime\prime}  \big)^{ | \mathcal{I}|
         }  \leq        \prod_{k^{\prime} \in \mathcal{I}} \textbf{P}^{\xi^{\mathrm{Sloped}}}_{[-m,m] \times [0,n^{\prime}N]} \big[   \mathscr{L}_{i^{\prime\prime}+1}                     \overset{h \leq ck}{\longleftrightarrow}        \mathscr{L}^{\prime}_{i^{\prime\prime}+1}      \big]                   \text{ } \text{ . } \text{ } \\
  \end{align*}


    \noindent Altogether, the constant from the two terms provides the lower bound,
    
    \begin{align*}
      \textbf{P}^{\xi^{\mathrm{Sloped}}}_{[-m,m] \times [0,n^{\prime}N]}   \big[      \mathscr{L}_{i^{\prime\prime}}           \overset{h \leq ck}{\longleftrightarrow}   \mathscr{L}_{i^{\prime\prime}+1}          \big]      \textbf{P}^{\xi^{\mathrm{Sloped}}}_{[-m,m] \times [0,n^{\prime}N]}   \big[   \mathscr{L}_{i^{\prime\prime}+1} \text{ }             \overset{h \leq ck}{\longleftrightarrow}    ( F_{\mathrm{crit}})_i          \big]    \geq  C^1_{\prime\prime,\prime\prime\prime} C^2_{\prime\prime,\prime\prime\prime}       \text{ } \text{ . } \text{ } \tag{\textit{2.6.1.4.4}}
    \end{align*}

         \begin{figure}
\begin{align*}
\includegraphics[width=1.08\columnwidth]{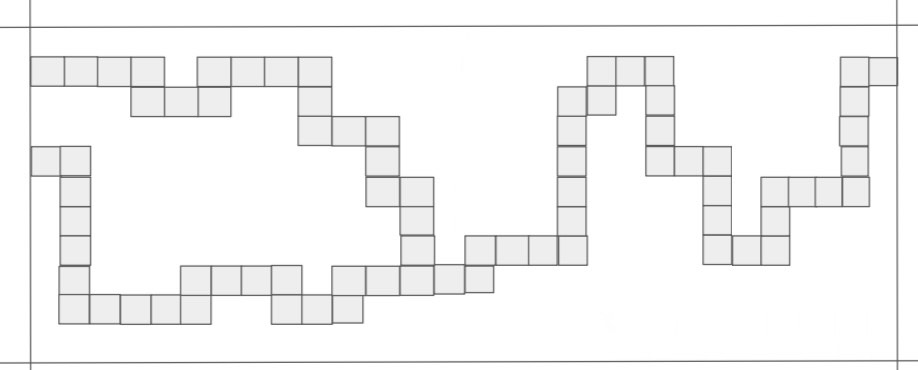}\\
\end{align*}
\caption{\textit{A macroscopic presence of faces vertically crossing the strip between two demarcated lines, with two connected components}. }
\end{figure}

  \noindent Similarly, for the second comparison, to bound the disconnective probability,
  
  \begin{align*}
        \textbf{P}^{\xi^{\mathrm{Sloped}}}_{[-m,m] \times [0,n^{\prime}N]} \big[         \mathscr{L}^{\prime\prime} \overset{h \geq ck}{\not\longleftrightarrow }   \mathscr{F}\mathscr{C}_{i+1}        \big]                  \text{ } \text{ , } \text{ }      \\
  \end{align*}
  
  \noindent we decompose the crossing event into a countable intersection, in which the terms in the following sequence of inequalities each term of which provides a lower bound,

  \begin{align*}
       \textbf{P}^{\xi^{\mathrm{Sloped}}}_{[-m,m] \times [0,n^{\prime}N]} \big[  \bigcap_{k \in \mathcal{I}} \big\{           \mathscr{L}^{\prime\prime} \overset{h \geq ck}{\not\longleftrightarrow }   \mathscr{L}^{\prime\prime}_k   \big\}    \big]    \overset{\mathrm{(FKG)}}{\geq}  \prod_{k \in \mathcal{I}}   \textbf{P}^{\xi^{\mathrm{Sloped}}}_{[-m,m] \times [0,n^{\prime}N]} \big[    \mathscr{L}^{\prime\prime} \overset{h \geq ck}{\not\longleftrightarrow }   \mathscr{L}^{\prime\prime}_k                          \big]             \text{ }   \text{ , } \text{ } \\
  \end{align*}

  \noindent from which each term in the product can be estimated with,

  \begin{align*}
     \textbf{P}^{\xi^{\mathrm{Sloped}}}_{[-m,m] \times [0,n^{\prime}N]} \big[   \mathscr{L}^{\prime\prime} \overset{h \geq ck}{\not\longleftrightarrow }   \mathscr{L}^{\prime\prime}_k                              \big]   \geq   \prod_{k^{\prime} \in\mathcal{I}}   \text{ }          \textbf{P}^{\xi^{\mathrm{Sloped}}}_{[-m,m] \times [0,n^{\prime}N]} \big[ \mathscr{L}^{\prime\prime} \overset{h \geq ck}{\not\longleftrightarrow }   \mathscr{L}^{\prime\prime}_{k^{\prime}}   \big]          \text{ }    \text{ , } \text{ } \\
  \end{align*}
  
  \noindent for some $\mathcal{I}^{\prime} \subset \mathcal{I}$, given the existence of finitely many lines for which,

  \begin{align*}
    \mathscr{L}^{\prime\prime}_k \cap    \mathscr{F}\mathscr{C}_{i+1}   \neq \emptyset \text{ } \text{ , } \text{ }    \\
  \end{align*}
  
  \noindent implying that the product lower bound holds, 
  
  \begin{align*}
       \prod_{k \in \mathcal{I}} \text{ }    \textbf{P}^{\xi^{\mathrm{Sloped}}}_{[-m,m] \times [0,n^{\prime}N]} \big[    \mathscr{L}^{\prime\prime} \overset{h \geq ck}{\not\longleftrightarrow }   \mathscr{L}^{\prime\prime}_k                               \big]   \geq       \prod_{k \in \mathcal{I}} \text{ }   \bigg[   \prod_{k^{\prime} \in \mathcal{I}^{\prime}}   \textbf{P}^{\xi^{\mathrm{Sloped}}}_{[-m,m] \times [0,n^{\prime}N]} \big[   \mathscr{L}^{\prime\prime} \overset{h \geq ck}{\not\longleftrightarrow }   \mathscr{L}^{\prime\prime}_{k^{\prime}}     \big]    \bigg]     \text{ } \\ \text{ }  \equiv \underset{k\in\mathcal{I}}{\prod_{k^{\prime} \in \mathcal{I}^{\prime}}} \text{ }   \textbf{P}^{\xi^{\mathrm{Sloped}}}_{[-m,m] \times [0,n^{\prime}N]} \big[   \mathscr{L}^{\prime\prime} \overset{h \geq ck}{\not\longleftrightarrow }   \mathscr{L}^{\prime\prime}_{k^{\prime}} \big]    \\ \geq    \underset{k\in\mathcal{I}}{\prod_{k^{\prime} \in \mathcal{I}^{\prime}}}     \mathscr{C}_{k,k^{\prime}}       \geq  \big( \mathscr{C}_{k,k^{\prime}}        \big)^{| \mathcal{I}^{\prime} \cap   \mathcal{I}|} \text{ } \text{ , } \text{ }   \tag{\textit{2.6.1.4.5}} 
  \end{align*}
  
  \noindent with,
  
  \begin{align*}
   \mathscr{C}_{k,k^{\prime}} \leq  \underset{k^{\prime}\in\mathcal{I}^{\prime}}{\mathrm{inf}} \text{ } \big\{   \textbf{P}^{\xi^{\mathrm{Sloped}}}_{[-m,m] \times [0,n^{\prime}N]} \big[    \mathscr{L}^{\prime\prime} \overset{h \geq ck}{\not\longleftrightarrow }   \mathscr{L}^{\prime\prime}_{k^{\prime}}   \big]         \big\}                      \text{ } \text{ . } \text{ } \\
  \end{align*}
  
  \noindent Upon examination of crossings of $\{ h \geq ck \}$ with crossings of $\{ h \leq ck\}$, one has the equality,
  
  \begin{align*}
       \textbf{P}^{\xi^{\mathrm{Sloped}}}_{[-m,m] \times [0,n^{\prime}N]} \big[    \mathscr{L}^{\prime\prime} \overset{h \geq ck}{\not\longleftrightarrow }   \mathscr{L}^{\prime\prime}_k              \big]  =   \textbf{P}^{\xi^{\mathrm{Sloped}}}_{[-m,m] \times [0,n^{\prime}N]} \big[   \mathscr{L}^{\prime\prime} \overset{h \leq ck}{\longleftrightarrow }   \mathscr{L}^{\prime\prime}_k             \big]   \text{ } \text{ , } \text{ } \\
  \end{align*}
  
  \noindent implying that it suffices to provide the same lower bound for,
  
  \begin{align*}
            \textbf{P}^{\xi^{\mathrm{Sloped}}}_{[-m,m] \times [0,n^{\prime}N]} \big[   \mathscr{L}^{\prime\prime} \overset{h \leq ck}{\longleftrightarrow }   \mathscr{L}^{\prime\prime}_k            \big]            \text{ } \text{ , } \text{ } \\
  \end{align*}
  
  \noindent which is given by the product lower bound below,
  
  \begin{align*}
  \prod_{k \in \mathcal{I}} \text{ }      \textbf{P}^{\xi^{\mathrm{Sloped}}}_{[-m,m] \times [0,n^{\prime}N]} \big[  \mathscr{L}^{\prime\prime} \overset{h \leq ck}{\longleftrightarrow }   \mathscr{L}^{\prime\prime}_k                  \big]    \geq  \mathcal{C}_k     \text{ } \text{ , } \text{ }     \tag{\textit{2.6.1.4.6}}            \\
  \end{align*}

    \noindent which is achieved with the two estimates,

  \begin{align*}
        \text{ } \big(  \mathcal{C}^{\prime}_k  \big)^{|\mathcal{I}|}  \geq  \mathcal{C}_k \text{ }    \tag{\textit{2.6.1.4.7}}  \text{ , } \\
  \end{align*}

  \noindent for $\mathcal{C}_k$ small enough, with the constant in the upper bound of \textit{(2.6.1.4.3)} taking the form,
  
  \begin{align*}
        \mathcal{C}^{\prime}_k =  \underset{k \in \mathcal{I}}{\mathrm{inf}}  \big\{   1 -    \textbf{P}^{\xi^{\mathrm{Sloped}}}_{[-m,m] \times [0,n^{\prime}N]} \big[   \mathscr{L}^{\prime\prime} \overset{h \leq ck}{\longleftrightarrow }   \mathscr{L}^{\prime\prime}_k   \big]    \big\}  \text{ } \text{ , } \text{ } \\
  \end{align*}

 \noindent  by finite energy. From such estimates, the remaining probability appearing in the intersection of crossing events for,

 \begin{align*}
    \textbf{P}^{\xi^{\mathrm{Sloped}}}_{[-m,m] \times [0,n^{\prime}N]} \big[         \mathscr{L}^{\prime\prime} \overset{h \geq ck}{\not\longleftrightarrow }   \mathscr{F}\mathscr{C}_{i+1}         \big]    \text{ , } \\
 \end{align*}

 \noindent is,
 
     \begin{align*}
          \textbf{P}^{\xi^{\mathrm{Sloped}}}_{[-m,m] \times [0,n^{\prime}N]} \big[       \mathscr{F}\mathscr{C}_i     \underset{\textbf{Z}^2 \cap ( \mathscr{F} \mathscr{C})^c_i}{\overset{h \leq ck}{\longleftrightarrow}}    \mathscr{F}\mathscr{C}_{i+1}    \big]       \text{ } \text{ , } \text{ }   \end{align*}

  \bigskip

         \begin{figure}
\begin{align*}
\includegraphics[width=0.88\columnwidth]{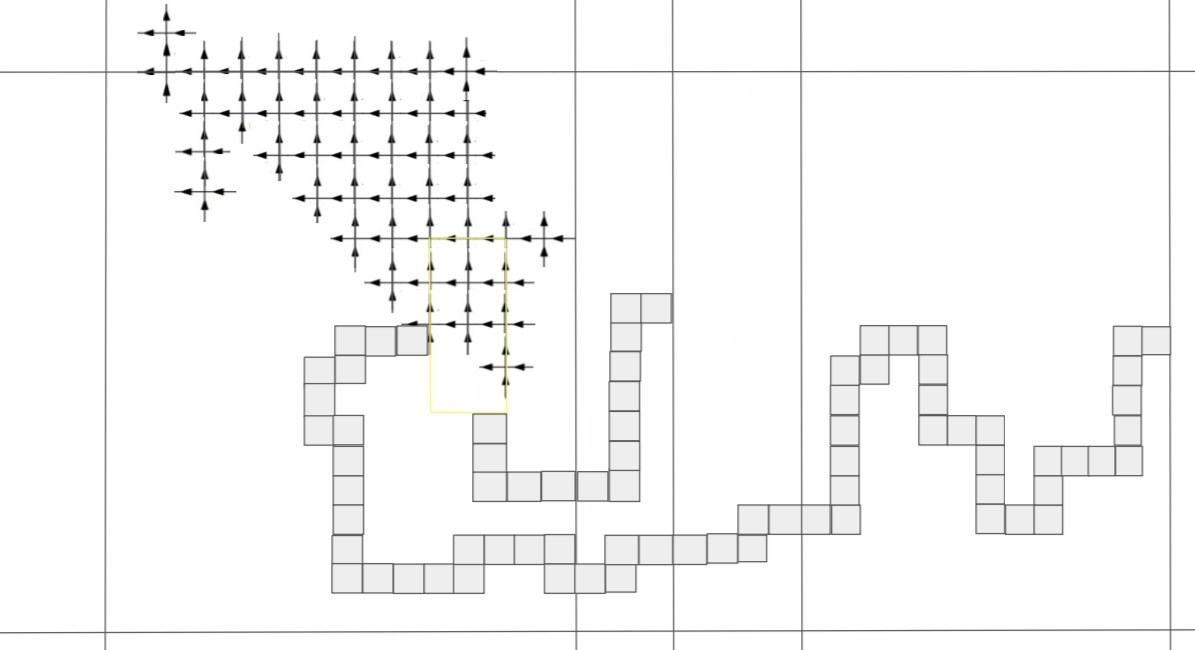}\\
\end{align*}
\caption{\textit{A finite volume configuration with yellow finite subvolume surrounding a boundary of the freezing cluster in the strip}. Vertical crossings in the complement of a \textit{freezing cluster} in the strip are depicted above, which relates to several lower bounds obtained throughout the section. One such crossing intersects the first and second lines to the right of the subset of frozen faces, while another crossing consisting of more faces also intersects the third and fourth lines to the right of the subset of frozen faces.}
\end{figure}

  \noindent which, upon rearrangement, produces another series of estimates, satisfying,

  \begin{align*}
    \textbf{P}^{\xi^{\mathrm{Sloped}}}_{[-m,m] \times [0,n^{\prime}N]} \big[       \mathscr{F}\mathscr{C}_i     \underset{\textbf{Z}^2 \cap ( \mathscr{F} \mathscr{C})^c_i}{\overset{h \leq ck}{\longleftrightarrow}}    \mathscr{F}\mathscr{C}_{i+1}     \big]  \geq    \textbf{P}^{\xi^{\mathrm{Sloped}}}_{[-m,m] \times [0,n^{\prime}N]} \big[  \bigcap_{l^{\prime\prime} \in \mathcal{I}}  \big\{            \mathscr{F}\mathscr{C}_i      \underset{\textbf{Z}^2 \cap ( \mathscr{F} \mathscr{C})^c_i}{\overset{h \leq ck}{\longleftrightarrow}}     \mathscr{L}_{l^{\prime\prime}}      \big\}      \big]                  \\ \overset{\mathrm{(FKG)}}{\geq} \text{ } \prod_{l^{\prime\prime} \in \mathcal{I}} \text{ } \textbf{P}^{\xi^{\mathrm{Sloped}}}_{[-m,m] \times [0,n^{\prime}N]} \big[    \mathscr{F}\mathscr{C}_i      \underset{\textbf{Z}^2 \cap ( \mathscr{F} \mathscr{C})^c_i}{\overset{h \leq ck}{\longleftrightarrow}}     \mathscr{L}_{l^{\prime\prime}}       \big]            \text{ } \text{ , } \text{ } \\
  \end{align*}
  
  \noindent which, upon further rearrangement, implies the following lower bound for the connectivity event connecting the \textit{freezing cluster} to $ \mathscr{L}_{l^{\prime\prime}}$, in $\textbf{Z}^2 \cap ( \mathscr{F} \mathscr{C})^c_i$,
  
  \begin{align*}
         \textbf{P}^{\xi^{\mathrm{Sloped}}}_{[-m,m] \times [0,n^{\prime}N]} \big[       \mathscr{F}\mathscr{C}_i      \underset{\textbf{Z}^2 \cap ( \mathscr{F} \mathscr{C})^c_i}{\overset{h \leq ck}{\longleftrightarrow}}     \mathscr{L}_{l^{\prime\prime}}       \big]   \text{ } \geq     \mathfrak{C}   ( l^{\prime\prime})  \equiv \mathfrak{C}_{l^{\prime\prime}}      \text{ } \text{ , } \text{ } \\
  \end{align*}
  
  \noindent by finite energy, with the corresponding product lower bound,
  
  \begin{align*}
      \prod_{l^{\prime\prime} \in \mathcal{I}}   \text{ }  \mathfrak{C}_{l^{\prime\prime}}  \geq  \mathfrak{C}^{|\mathcal{I}|}_{l^{\prime\prime}}  \equiv  \mathfrak{C}^{\prime}_{l^{\prime\prime}}     \text{ } \text{ , } \text{ }  \tag{\textit{2.6.1.4.8}}  \\
  \end{align*}
  
  \noindent from which we obtain the cumulative estimate, respectively taken over natural indices $k$ and $l^{\prime\prime}$,
  
  \begin{align*}
 \text{ }  \bigg[     \underset{l^{\prime\prime} \in \mathcal{I}}{\prod_{k \in \mathcal{I}} }  \text{ }      \textbf{P}^{\xi^{\mathrm{Sloped}}}_{[-m,m] \times [0,n^{\prime}N]} \big[    \mathscr{L}^{\prime\prime} \overset{h \leq ck}{\longleftrightarrow }   \mathscr{L}^{\prime\prime}_k               \big]              \bigg]  \bigg[       \textbf{P}^{\xi^{\mathrm{Sloped}}}_{[-m,m] \times [0,n^{\prime}N]} \big[     \mathscr{F}\mathscr{C}_i      \underset{\textbf{Z}^2 \cap ( \mathscr{F} \mathscr{C})^c_i}{\overset{h \leq ck}{\longleftrightarrow}}     \mathscr{L}_{l^{\prime\prime}}       \big]             \bigg]                     \text{ } \text{ , }  \text{ }
  \end{align*}
  
  \noindent can be rearranged insofar as to obtain the following equivalent product expansion,
  
  \begin{align*}
     \bigg[   \prod_{k \in \mathcal{I}} \text{ }      \textbf{P}^{\xi^{\mathrm{Sloped}}}_{[-m,m] \times [0,n^{\prime}N]} \big[   \mathscr{L}^{\prime\prime} \overset{h \leq ck}{\longleftrightarrow }   \mathscr{L}^{\prime\prime}_k                  \big] \bigg]   \bigg[     \prod_{l^{\prime\prime} \in \mathcal{I}} \text{ } \textbf{P}^{\xi^{\mathrm{Sloped}}}_{[-m,m] \times [0,n^{\prime}N]} \big[       \mathscr{F}\mathscr{C}_i      \underset{\textbf{Z}^2 \cap ( \mathscr{F} \mathscr{C})^c_i}{\overset{h \leq ck}{\longleftrightarrow}}     \mathscr{L}_{l^{\prime\prime}}     \big]    \bigg]            \text{ } \text{ , } \text{ } \\
  \end{align*}

  \noindent which from respectively applying the estimates of \textit{(2.6.1.4.4)} and \textit{(2.6.1.4.6)}, imply,
  
  \begin{align*}
          \bigg[   \prod_{k \in \mathcal{I}} \text{ }                \mathcal{C}_k      \bigg]  \bigg[     \prod_{l^{\prime\prime} \in \mathcal{I}}                           \textbf{P}^{\xi^{\mathrm{Sloped}}}_{[-m,m] \times [0,n^{\prime}N]} \big[       \mathscr{F}\mathscr{C}_i      \underset{\textbf{Z}^2 \cap ( \mathscr{F} \mathscr{C})^c_i}{\overset{h \leq ck}{\longleftrightarrow}}     \mathscr{L}_{l^{\prime\prime}}        \big]                     \bigg]\geq  \bigg[  \prod_{k \in \mathcal{I}} \text{ }                \mathcal{C}_k   \bigg]    \bigg[   \prod_{l^{\prime\prime} \in \mathcal{I}}  \text{ }   \mathfrak{C}^{\prime}_{l^{\prime\prime}}                        \bigg]         \\ \geq \underset{l^{\prime\prime} \in \mathcal{I}}{\prod_{k \in \mathcal{I}}} \text{ } \mathcal{C}_k  \mathfrak{C}^{\prime}_{l^{\prime\prime}}   \text{ } \\  \geq      \mathcal{C}^{|\mathcal{I}|}_k  \big( \mathfrak{C}^{\prime}_{l^{\prime\prime}}  \big)^{|\mathcal{I}|}  \equiv   \big(  \mathcal{C}_k  \mathfrak{C}^{\prime}_{l^{\prime\prime}}    \big)^{|\mathcal{I}|}  \equiv      \mathfrak{C}^{\prime}_{k, l^{\prime\prime}}   \big(   |\mathcal{I}|     \big)      \text{ } \text{ . } \text{ }     \tag{\textit{2.6.1.4.9}}                         \\
  \end{align*}

  \noindent From arguments provided for the lower bound of the two probabilities, the superposition,

  \begin{align*}
     \textbf{P}^{\xi^{\mathrm{Sloped}}}_{[-m,m] \times [0,n^{\prime}N]} \big[    \mathscr{F}\mathscr{C}_i       \overset{h \geq ck}{\longleftrightarrow}      \mathscr{L}^{\prime}       \big]               \textbf{P}^{\xi^{\mathrm{Sloped}}}_{[-m,m] \times [0,n^{\prime}N]} \big[           \mathscr{L}^{\prime\prime} \overset{h \geq ck}{\not\longleftrightarrow }   \mathscr{F}\mathscr{C}_{i+1}      \big]  \text{ } \text{ , }
     \end{align*}
     
     \noindent can be lower bounded with,
       \begin{align*}
     \text{ }           \textbf{P}^{\xi^{\mathrm{Sloped}}}_{[-m,m] \times [0,n^{\prime}N]} \big[           \mathscr{L}_i   \overset{h \leq ck}{\longleftrightarrow}               (F_{\mathrm{crit}})_i     \big]       \textbf{P}^{\xi^{\mathrm{Sloped}}}_{[-m,m] \times [0,n^{\prime}N]} \bigg[ \big\{   (F_{\mathrm{supcrit}})_i \overset{h \leq ck}{\longleftrightarrow}      \mathscr{L}_2      \big\}    \cap   \big\{ \mathscr{F}\mathscr{C}_i     \underset{\textbf{Z}^2 \cap ( \mathscr{F} \mathscr{C})^c_i}{\overset{h \leq ck}{\longleftrightarrow}}    \mathscr{F}\mathscr{C}_{i+1}         \big\}      \bigg] 
  \end{align*}
  
  \noindent which yields the lower bound, from $\textit{(2.6.1.4.2)}$, $\textit{(2.6.1.4.3)}$, and $\textit{(2.6.1.4.4)}$,
  
   \begin{align*}
          C^1_{\prime\prime,\prime\prime\prime}  C^2_{\prime\prime,\prime\prime\prime}     \textbf{P}^{\xi^{\mathrm{Sloped}}}_{[-m,m] \times [0,n^{\prime}N]} \bigg[  \big\{    (F_{\mathrm{supcrit}})_i \overset{h \leq ck}{\longleftrightarrow}      \mathscr{L}_2    \big\}   \cap \big\{  \mathscr{F}\mathscr{C}_i     \underset{\textbf{Z}^2 \cap ( \mathscr{F} \mathscr{C})^c_i}{\overset{h \leq ck}{\longleftrightarrow}}    \mathscr{F}\mathscr{C}_{i+1}              \big\}   \bigg]          \text{ } \text{ , } \text{ } 
 \end{align*}
  
 \noindent corresponding to one estimate, while the other estimate yields, from \textit{(2.6.1.4.8)} and \textit{(2.6.1.4.9)},   
 
 \begin{align*}
          C^1_{\prime\prime,\prime\prime\prime}  C^2_{\prime\prime,\prime\prime\prime}                 \mathfrak{C}^{\prime}_{k, l^{\prime\prime}}    \text{ } \text{ . } \text{ }  
 \end{align*}
  
  \noindent from which we conclude the argument after having set $\big(  \mathscr{C}_4  \big)_i \equiv (  C^{1,2,\prime} )^{\prime}$, for $C^1_{\prime\prime,\prime\prime\prime} C^2_{\prime\prime,\prime\prime\prime}                   \mathfrak{C}^{\prime}_{k, l^{\prime\prime}} \geq ( C^{1,2} )^{\prime}$. \boxed{}

   \bigskip
   
   \noindent We return to the lower bound for the vertical crossing in the first case of \textbf{Lemma} \textit{2.6}.
   
   \noindent

   \bigskip
   
   \noindent \textbf{Lemma} $\textit{2.6.1.5}$ (\textit{transitive lower bound for vertical crossings across the strip}). Under the assumptions of \textbf{Lemma} $\textit{2.6.1.4}$, there exists a strictly positive constant satisfying,

   \begin{align*}
         \textbf{P}^{\xi^{\mathrm{Sloped}}}_{[-m,m] \times [0,n^{\prime} N]} \big[     \mathscr{F}\mathscr{C}_i \overset{h \geq ck}{\not\longleftrightarrow} \mathscr{F}\mathscr{C}_{N}            \big] \geq  \mathscr{C}_5   \text{ } \text{ , } \text{ }    \\
   \end{align*}

   \noindent for suitable $\mathscr{C}_5$.

   \bigskip
   
   \noindent \textit{Proof of Lemma 2.6.1.5}. Write,
   
   \begin{align*}
   \textbf{P}^{\xi^{\mathrm{Sloped}}}_{[-m,m] \times [0,n^{\prime}N]} \big[      \mathscr{F}\mathscr{C}_i \overset{h \geq ck}{\not\longleftrightarrow} \mathscr{F}\mathscr{C}_{N}     \big]  \geq     \textbf{P}^{\xi^{\mathrm{Sloped}}}_{[-m,m] \times [0,n^{\prime}N]} \big[     \bigcap_{i^{\prime} \in \mathcal{I}}        \big\{   \mathscr{F}\mathscr{C}_{i^{\prime}}  \overset{h \geq ck}{\not\longleftrightarrow}    \mathscr{F}\mathscr{C}_{i^{\prime}+1}     \big\}     \big] \text{ } \end{align*}
   
   \begin{align*}
   \overset{\mathrm{(FKG)}}{\geq} \prod_{i^{\prime}\in\mathcal{I}} \text{ }      \textbf{P}^{\xi^{\mathrm{Sloped}}}_{[-m,m] \times [0,n^{\prime}N]} \big[ \mathscr{F}\mathscr{C}_{i^{\prime}}  \overset{h \geq ck}{\not\longleftrightarrow}    \mathscr{F}\mathscr{C}_{i^{\prime}+1}     \big]  \text{ , } \text{ } \\ \overset{(\textbf{Lemma} \text{ } \textit{2.6.1.4})}{\geq}   \prod_{i^{\prime}\in\mathcal{I}}   \big( \mathscr{C}_4 \big)_{i^{\prime}} 
   \end{align*}

   \noindent yielding the desired estimate, upon setting $\mathscr{C}_5$ equal to,
   
      \begin{align*}
    \text{ }    \mathscr{C}_{4 , \mathrm{min}}^{|\mathcal{I}|}  \text{ } \text{ , } \text{ } \\
   \end{align*}
   
   \noindent for,
   
   \begin{align*}
  \text{ }   \mathscr{C}_{4,\mathrm{min}} =         \underset{i^{\prime} \in \mathcal{I}}{\mathrm{inf}} \text{ } \big(  \mathcal{C}_4  \big)_{i^{\prime}}        \text{ }    \text{ . } \\
   \end{align*}
   
   \noindent Hence we conclude the argument. \boxed{}

     \bigskip
   
  \noindent From previous results, we conclude \textbf{Case one} of the proof by  make use of the event containment, namely that the number of faces required for the intersection of crossing events,  
   
   \begin{align*}
        \bigg\{     \big\{      \mathscr{L}_1                    \overset{h \leq ck}{\longleftrightarrow}    F_{\mathrm{crit}}     \big\}    \cap     \big\{     F_{\mathrm{supcrit}}                        \overset{h \leq ck }{\longleftrightarrow}        \mathscr{L}_2       \big\} \cap            \big\{ \mathscr{F}\mathscr{C}_{i} \overset{h \geq ck}{\not\leftrightarrow}   \mathscr{F}\mathscr{C}_{i+1}           \big\}        \bigg\}   \text{ }         \text{ } \text{ , }  \\
   \end{align*}
   
   \noindent is contained within the number of faces required for the intersection of crossing events,

   \begin{align*} 
 \bigg\{   \big\{     F_1    \underset{\mathscr{S}}{\overset{h \leq ck }{\longleftrightarrow}}  \text{ }    F_n     \big\}       \cap         \big\{ \mathscr{F}\mathscr{C}_{i} \overset{h \geq ck}{\not\leftrightarrow}   \mathscr{F}\mathscr{C}_{i+1}                \big\}       \bigg\}        \text{ }  \text{ . } \text{ } \\
   \end{align*}
   
   \noindent in which the \textit{critical} face $F_{\mathrm{crit}}$, and \textit{supcritical} face $F_{\mathrm{supcrit}}$, are faces bounded within the $\textbf{Z}^2$ strip for which $\{ \text{ }    F_{\mathrm{Crit}+1}       \text{ }   \overset{h \geq ck}{\not\longleftrightarrow}  \text{ }   F_{\mathrm{supcrit}-1}   \text{ }  \}$. Hence it suffices to show that,

   \begin{align*}
   \textbf{P}^{\xi^{\mathrm{Sloped}}}_{[-m,m] \times [0,n^{\prime}N]}  \bigg[  \big\{    F_1    \underset{\mathscr{S}}{\overset{h \leq ck }{\longleftrightarrow}}    F_n     \big\}      \cap         \big\{  \mathscr{F}\mathscr{C}_{i} \overset{h \geq ck}{\not\leftrightarrow}   \mathscr{F}\mathscr{C}_{i+1}              \big\}  \bigg] \text{ }            \text{ } \text{ , } \text{ } \\
   \end{align*}

   \noindent is bound below by the same constant as for the intersection of crossing probabilities, 
   
    \begin{align*}
        \textbf{P}^{\xi^{\mathrm{Sloped}}}_{[-m,m] \times [0,n^{\prime}N]}  \bigg[                       \big\{        \mathscr{L}_1                    \overset{h \leq ck}{\longleftrightarrow}    F_{\mathrm{crit}}     \big\}      \cap    \big\{       F_{\mathrm{supcrit}}                        \overset{h \leq ck }{\longleftrightarrow}            \mathscr{L}_2       \big\}  \cap           \big\{ \mathscr{F}\mathscr{C}_{i} \overset{h \geq ck}{\not\leftrightarrow}   \mathscr{F}\mathscr{C}_{i+1}           \big\} \bigg]             \text{ } \text{ . } \text{ } \\
   \end{align*}

   \noindent To this end, we incorporate the product estimate from a previous expression for one crossing probability given in (\textit{2.6.1}),
   
   \begin{figure}
\begin{align*}
\includegraphics[width=1.08\columnwidth]{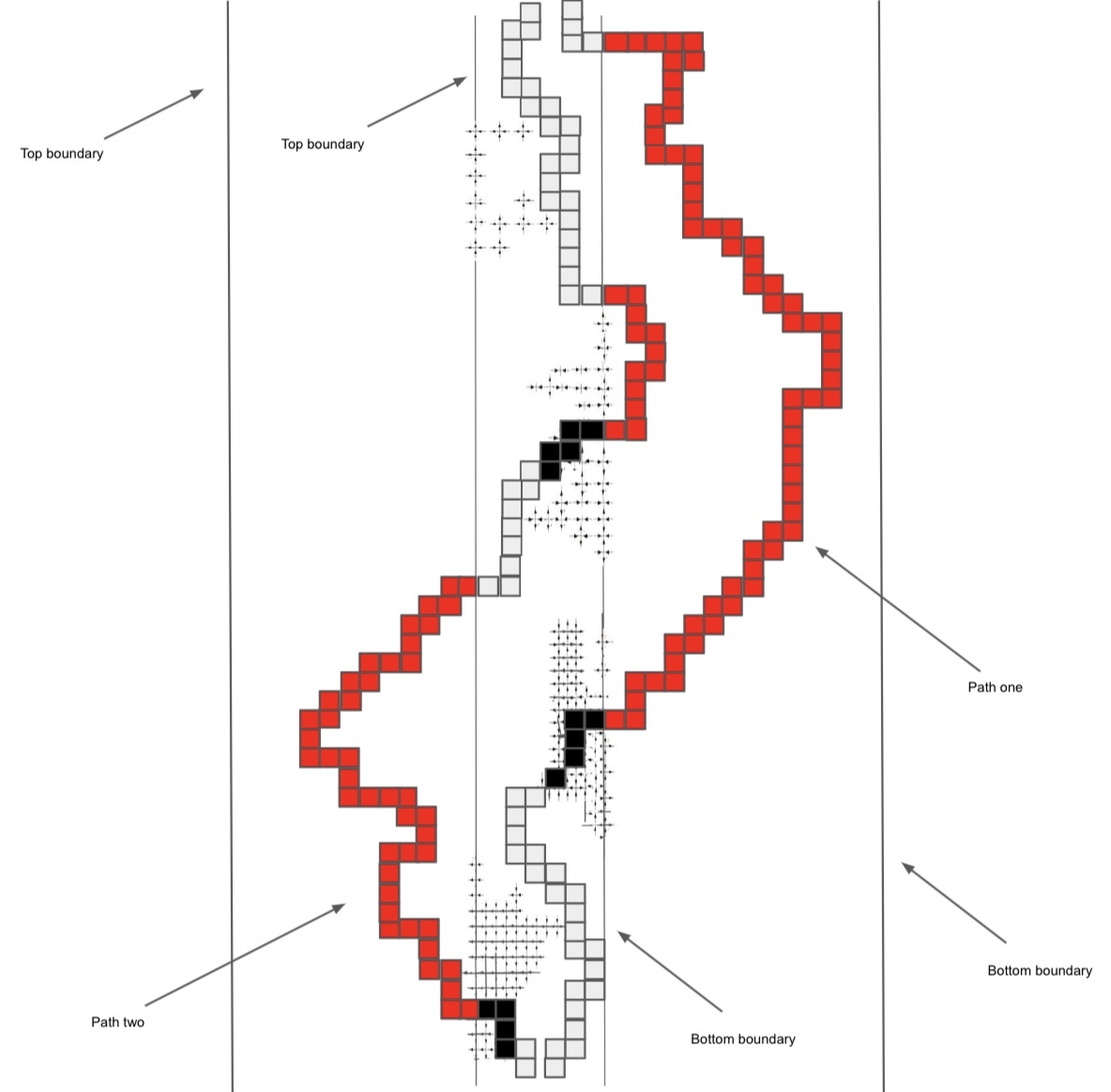}\\
\end{align*}
\caption{\textit{Removing the requirement that paths within strips of the square lattice not intersect freezing clusters can result in paths failing to remain within the top boundary and bottom boundary}. \underline{Path one}, as given from grey faces, begins by intersecting the second \textit{freezing cluster}, as given from black faces, afterwards intersecting the top boundary of the smaller square lattice strip, as given from red faces, and finally, intersects the top boundary of the smaller square lattice strip, as given from grey faces. \underline{Path two}, as given from grey faces, begins by intersecting the first \textit{freezing cluster}, as given from black faces, afterwards intersecting the bottom boundary of the square lattice strip, as given from red faces, afterwards intersecting the bottom boundary of the square lattice strip again, as given from grey faces, and finally, exits and re-enters the smaller square lattice strip, respectively given from black, red and grey faces.}
\end{figure}

   \begin{align*} 
\prod_{2 \leq i \leq N-2}\textbf{P}^{\xi^{\mathrm{Sloped}}}_{[-m,m] \times [0,n^{\prime}N]} [   \mathscr{F}\mathscr{C}_i \overset{h \geq ck}{\not\longleftrightarrow} \mathscr{F}\mathscr{C}_{i+1} ]     \geq   (c_{\mathscr{F}\mathscr{C}}   )^{N-3}         \text{ } \text{ , } \text{ } \\
   \end{align*}
   
   \noindent each of which exist almost surely by finite energy, for suitable $c_{\mathscr{F}\mathscr{C}}$. As an intersection of crossing events occurring simultaneously, we exhibit the lower bound for the vertical domain crossings, as,
   
   \begin{align*}
        \textbf{P}^{\xi^{\mathrm{Sloped}}}_{[-m,m] \times [0,n^{\prime}N]} \bigg[   \big\{   \gamma_L  \overset{h\geq ck}{\leftrightarrow} \mathscr{F}\mathscr{C}_1   \big\}  \cap  \big\{    \mathscr{F}\mathscr{C}_i \overset{h \geq ck}{\not\longleftrightarrow} \mathscr{F}\mathscr{C}_{i+1}    \big\}   \cap   \big\{   \mathscr{F}\mathscr{C}_{N}               \overset{h\geq ck}{\leftrightarrow}    \gamma_R   \big\}           \bigg]  \text{ }
        \end{align*}
        
        \noindent by $\mathrm{(FKG)}$, yields the lower bound,

        \begin{align*}
     \text{ }   \textbf{P}^{\xi^{\mathrm{Sloped}}}_{[-m,m] \times [0,n^{\prime}N]} \big[   \gamma_L  \overset{h\geq ck}{\leftrightarrow} \mathscr{F}\mathscr{C}_1          \big]   \textbf{P}^{\xi^{\mathrm{Sloped}}}_{[-m,m] \times [0,n^{\prime}N]} \big[   \mathscr{F}\mathscr{C}_i \overset{h \geq ck}{\not\longleftrightarrow} \mathscr{F}\mathscr{C}_{i+1}   \big] 
      \textbf{P}^{\xi^{\mathrm{Sloped}}}_{[-m,m] \times [0,n^{\prime}N]} \big[     \mathscr{F}\mathscr{C}_{N}               \overset{h\geq ck}{\leftrightarrow}    \gamma_R            \big]   \text{ } \text{ , } \text{ } 
   \end{align*}

   \noindent after which incorporating previously obtained lower bounds for each probability yields,

   \begin{align*}
   \text{ }  (c_{\mathscr{F}\mathscr{C}}   )^{N-3}  \text{ }       \textbf{P}^{\xi_{\mathrm{Sloped}}}_{[-m,m]\times[0,n^{\prime}N]} \big[     \gamma_L \overset{h\geq ck}{\leftrightarrow} \mathscr{F}\mathscr{C}_1  \big]  \text{ }              \text{ }          \textbf{P}^{\xi_{\mathrm{Sloped}}}_{[-m,m]\times[0,n^{\prime}N]} \big[  \mathscr{F}\mathscr{C}_{N}               \overset{h\geq ck}{\leftrightarrow}    \gamma_R    \big]   \\ \text{ }  \overset{\mathrm{(\textbf{Lemma } \textit{2.2})\text{ } }}{\geq} \text{ }   (c_{\mathscr{F}\mathscr{C}}   )^{N-3}  \bigg|  F (    \gamma_L     \cap     \mathcal{L}^1_{\mathscr{F}\mathscr{C}_1}   ) \bigg| \bigg(   \big( \delta^1_L \big)^{-1} \delta^2_L    \bigg)          \textbf{P}^{\xi_{\mathrm{Sloped}}}_{[-m,m]\times[0,n^{\prime}N]} \big[  \mathscr{F}\mathscr{C}_{N}               \overset{h\geq ck}{\leftrightarrow}    \gamma_R    \big]     \\ \overset{\mathrm{\text{ } (\textbf{Lemma } \textit{2.3})}}{\geq}    (c_{\mathscr{F}\mathscr{C}}   )^{N-3}          {\bigg|  F (   \gamma_L     \cap    \mathcal{L}^1_{\mathscr{F}\mathscr{C}_1}   )  F (    \gamma_R    \cap  \mathcal{L}^2_{\mathscr{F}\mathscr{C}_N} )   \bigg| \bigg[  \big( \delta^1_L \big)^{-1} \delta^2_L \bigg]  }    {   \bigg[ \big( \delta^1_R \big)^{-1} \text{ }  \delta^2_R }     \bigg] \text{ } \\   \equiv     (c_{\mathscr{F}\mathscr{C}}   )^{N-3}       \bigg| \text{ } F (    \gamma_L     \cap    \mathcal{L}^1_{\mathscr{F}\mathscr{C}_1}  ) \bigg|            \bigg| F (    \gamma_R    \cap   \mathcal{L}^2_{\mathscr{F}\mathscr{C}_N}   )  \bigg|            \bigg[   \frac{\delta^2_L \delta^2_R  }{\delta^1_L \delta^1_R }     \bigg]         \text{ }  \\   \text{ } \geq \text{ }  C_{\mathscr{F}\mathscr{C}} \text{ }    C_L         \text{ }    C_R    \text{ }  \\  \geq   C_{\mathscr{F}\mathscr{C}}   \mathscr{C}^2_4   \text{ }  \text{ , } \text{ } 
    \end{align*}

    \noindent where the estimate in the last inequality is provided from a combination of results from \textbf{Lemma} \textit{2.2} and \textbf{Lemma} \textit{2.3}, while in the ultimate lower bound, each respective estimate can be bound below with,
    
    \begin{align*}
        \mathscr{C}_4 \text{ } = \text{ } \underset{\mathrm{sides } \text{ } L,R}{\mathrm{inf}} \text{ }  \big\{        C_L     ,      C_R       \big\}           \text{ } \text{ , } \text{ } \\
    \end{align*}

    \noindent for $C_L,C_R >0$, and $C_{\mathscr{F}\mathscr{C}}$ for which,
    
    \begin{align*}
   c_{\mathscr{F}\mathscr{C}} \text{ } \geq \text{ }     \sqrt[N-3]{c_{\mathscr{F}\mathscr{C}}}    \text{ } \text{ . } \text{ } \\
    \end{align*}

     \noindent Finally, for $FV \equiv \mathscr{D}$, and $\mathcal{V}^{h \leq (1-c)k} (\mathscr{D}) \equiv \mathcal{V}^{h \leq (1-c)k}$ as given in the statement of \textbf{Lemma} $\textit{2.6}$,

    \begin{align*}
       \text{ }  \textbf{P}^{\chi^{\prime}_{\mathrm{Sloped}}|_{\partial (  [-m,m] \times [0,n^{\prime}_2 N]   ) }}_{\mathscr{D}} \big[  \mathcal{V}^{h \leq (1-c)k}   \big]   \geq     \textbf{P}^{\chi^{\prime}_{\mathrm{Sloped}}|_{\partial (  [-m,m] \times [0,n^{\prime}_2 N]   ) }}_{\mathscr{D}} \big[ \mathcal{V}^{h \leq (1-c)k}  \big]    \text{ } \\  {\geq}  \textbf{P}^{\chi^{\prime}_{\mathrm{Sloped}}|_{\partial (  [-m,m] \times [0,n^{\prime}_2 N]   ) }}_{[-m,m]\times[0,n_2^{\prime}N]} \big[  \mathcal{V}^{h \leq (1-c)k}  \big]      \text{ } \text{ , } \text{ }  \\
    \end{align*}

             \begin{figure}
\begin{align*}
\includegraphics[width=0.88\columnwidth]{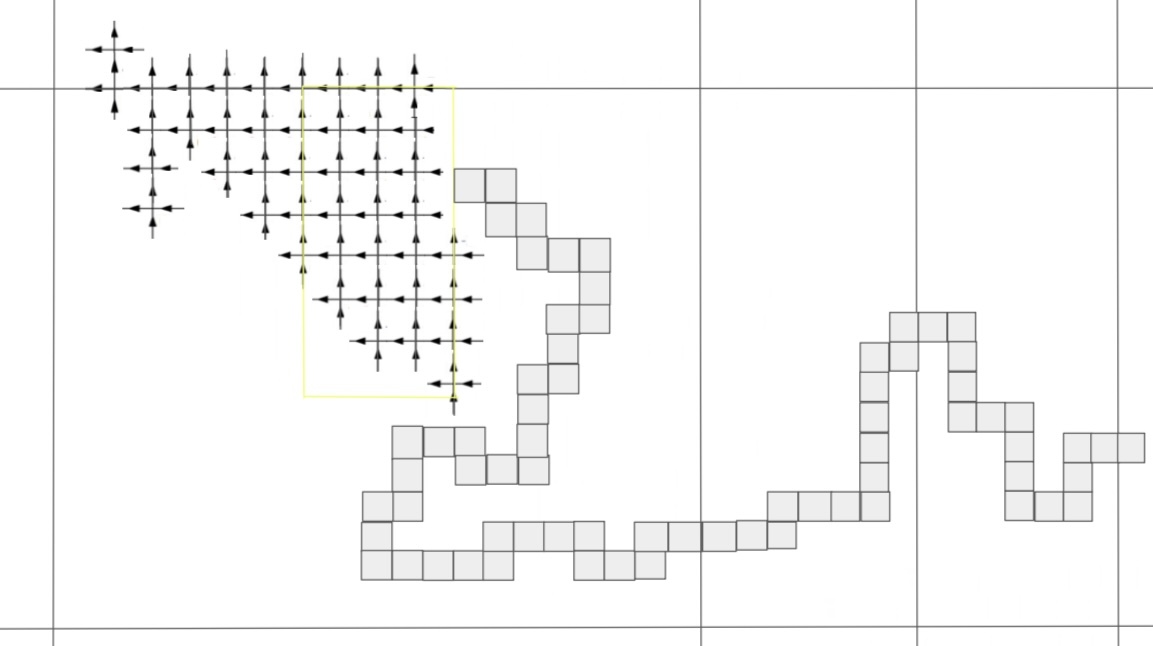}\\
\end{align*}
\caption{\textit{A finite volume configuration with yellow finite sub-volume surrounding a boundary of the freezing cluster in the strip}. Another possible configuration, exhibited above, consists of a crossing intersecting three lines to the left of the \textit{freezing cluster}.}
\end{figure}

    \noindent where $\mathcal{H}^{h \leq (1-c)k}(\mathscr{D}) = \{ \mathscr{F}_L \underset{\mathscr{D}}{\overset{h\geq ck}{\leftrightarrow}}    \mathscr{F}_R    \} $, where $\mathscr{F}_L, \mathscr{F}_R \in F( \text{ } \textbf{Z}^2 \text{ } )$ denote the leftmost and rightmost faces of the crossing enclosed within the interior of $\mathscr{D}$. Because the number of faces required for the crossing for the probability measure supported over $\mathscr{D}$, in comparison to the probability measure in the lower bound, requires less faces for the height-function crossing at, or above, threshold $ck$ to occur, one has that the collection of faces enclosed within $\mathscr{D}$ for the crossing event below to occur,
    
    \begin{align*}
    \text{ }    \bigg| \big\{   F ,  G \in F(\mathscr{D})  : \{  F \overset{h \geq ck}{\leftrightarrow}  G   \} \cap  \{    F_L      \overset{h \geq ck}{\leftrightarrow}    F_R       \} \neq \emptyset  \big\}   \bigg|      \text{ } \text{ , } 
    \end{align*}

  \noindent is strictly dominated by the number of faces within the strip for the same crossing event to occur,

    \begin{align*}
  \bigg|             \big\{    F^{\prime} ,  G^{\prime} \in F([-m,m] \times [0,n^{\prime}_2 N])  :    \{  F^{\prime}  \overset{h \geq ck}{\leftrightarrow}  G^{\prime}    \}  \cap   \{          F_L   \overset{h \geq ck}{\leftrightarrow}      F_R        \} \neq \emptyset   \big\}    \bigg|  \text{ } \text{ , }  
    \end{align*}
    
    \noindent where $F_L$ and $F_R$ are respectively denote the leftmost and rightmost faces belonging to $\mathcal{H}$. Concluding, this implies,
    
    \begin{align*}
 \textbf{P}^{\chi^{\prime}_{\mathrm{Sloped}}|_{\partial (  [-m,m] \times [0,n^{\prime}_2 N]   ) }}_{\mathscr{D}} \big[   \mathcal{I}_0      \overset{h \geq (1-c)k}{\longleftrightarrow}       \widetilde{\mathcal{I}_0}   \big]  \geq\textbf{P}^{\chi^{\prime}_{\mathrm{Sloped}}|_{\partial (  [-m,m] \times [0,n^{\prime}_2 N]   ) }}_{\mathscr{D}} \big[     \mathcal{V}^{h \leq (1-c)k}   \big]   \text{ }  \text{ , } \text{ } \\
    \end{align*}

    \noindent in turn giving,
    
    \begin{align*}
          \textbf{P}^{\chi^{\prime}_{\mathrm{Sloped}}|_{\partial (  [-m,m] \times [0,n^{\prime}_2 N]   ) }}_{\mathscr{D}} \big[     \mathcal{V}^{h \leq (1-c)k}  \big]    \equiv \textbf{P}^{\xi^{\mathrm{Sloped}}}_{\mathscr{D}} \big[      \mathcal{V}^{h \leq (1-c)k}      \big]  \leq     C^{\emptyset}_{\mathcal{V}}       \text{ } \text{ , } \text{ } \\
    \end{align*}
    
    \noindent which completes the argument for \textbf{Case one} of \textbf{Lemma} $\textit{2.6}$ with $\gamma_L \cap \gamma_R = \emptyset$.

 \bigskip

    \noindent  \textbf{Case two} (\textit{intersecting left and right boundaries of the strip domain}). In comparison to \textbf{Case one} above, we decompose the vertical crossing event of $h$ across the complement of all \textit{freezing clusters} which have nonempty intersection with each $\mathscr{D}_i$, as,
    
    \begin{align*}
        \textbf{P}^{\xi_{\mathrm{Sloped}}}_{[-m,m]\times[0,n^{\prime}N]} \big[      \mathcal{V}^{h \leq (1-c)k}_{(  \mathscr{D}_i  \cap  \mathscr{F} \mathscr{C}_i )^c}      \big]     =             \textbf{P}^{\xi_{\mathrm{Sloped}}}_{[-m,m] \times [0,n^{\prime} N]} \bigg[ \bigcap_{i \in \mathcal{I}^{\prime}}   \bigg\{  \{  \mathcal{L}^1_{\mathscr{F}\mathscr{C}_i} \underset{\mathscr{D}_i}{\overset{h\geq ck}{\longleftrightarrow}}  \widetilde{\gamma_L}  \}     \cap    \big\{ \widetilde{\gamma_L} \underset{(  S    \cap  \mathscr{D}_i   )^c }{\overset{h \geq ck}{\longleftrightarrow}} \mathscr{F}\mathscr{C}_{i+1}     \big\} \\ \cap    \big\{    \mathscr{F}\mathscr{C}_N \underset{\mathscr{D}_N \text{ }  \cap     \mathcal{L}^N_{\mathscr{F}\mathscr{C}_N}       }{\overset{h \geq ck}{\longleftrightarrow}} \mathcal{L}^N_{\mathscr{F}\mathscr{C}_N}   \big\}     \bigg\}   \bigg] \text{ , } 
      \end{align*}

        \noindent the first part of which can be respectively decomposed as, from one application of $\mathrm{(FKG)}$,

        \begin{align*}
 \text{ }    \prod_{i \in \mathcal{I}^{\prime}} \text{ } \textbf{P}^{\xi_{\mathrm{Sloped}}}_{[-m,m] \times [0,n^{\prime} N]} \bigg[ \big\{  \mathcal{L}^1_{\mathscr{F}\mathscr{C}_1} \underset{\mathscr{D}_i}{\overset{h\geq ck}{\leftrightarrow}}  \widetilde{\gamma_L}  \big\}        \cap    \big\{ \widetilde{\gamma_L} \underset{(  S    \cap  \mathscr{D}_i   )^c }{\overset{h \geq ck}{\longleftrightarrow}} \mathscr{F}\mathscr{C}_{i+1}  \big\}   \bigg]  \text{ } \\ 
    \end{align*}
    
    \noindent which can be lower bounded with,
    
    \begin{align*}
 \text{ }          \prod_{i \in \mathcal{I}^{\prime}} \text{ } \textbf{P}^{\xi_{\mathrm{Sloped}}}_{[-m,m] \times [0,n^{\prime} N]} \text{ } \big[          \mathcal{L}^1_{\mathscr{F}\mathscr{C}_1} \underset{\mathscr{D}_1}{\overset{h\geq ck}{\leftrightarrow}}  \widetilde{\gamma_L}  \big] \textbf{P}^{\xi_{\mathrm{Sloped}}}_{[-m,m] \times [0,n^{\prime} N]}  \big[      \widetilde{\gamma_L} \underset{( S  \cap  \mathscr{D}_i  )^c }{\overset{h \geq ck}{\leftrightarrow}} \mathscr{F}\mathscr{C}_{i+1}           \big]      \textbf{P}^{\xi_{\mathrm{Sloped}}}_{[-m,m] \times [0,n^{\prime} N]} \big[     \mathscr{F}\mathscr{C}_N \underset{\mathscr{D}_N   \cap       \mathcal{L}^N_{\mathscr{F}\mathscr{C}_N}       }{\overset{h \geq ck}{\leftrightarrow}}              \mathcal{L}^N_{\mathscr{F}\mathscr{C}_N}        \big]     \text{ } \text{ , } \text{ } \\
\tag{\textit{2.6.2}}    \end{align*}

    \noindent following a second application of $\mathrm{(FKG)}$, where in the first and third connectivity events, we respectively denote lines oriented with respect to \textit{freezing clusters} $\mathscr{F}\mathscr{C}_1$ and $\mathscr{F}\mathscr{C}_N$ with $\mathcal{L}^1_{\mathscr{F}\mathscr{C}_1}$ and $\mathcal{L}^N_{\mathscr{F}\mathscr{C}_N}$; the product in (\textit{2.6.2}) is taken over a countable index set $\mathcal{I}^{\prime}$, where $\widetilde{\gamma_L}$ is one part of the boundary of the strip from crossings of $h \geq ck$.
    
    \bigskip
    
    \noindent To lower bound $\textit{(2.6.2)}$, we implement the following series of results for each individual term, similarly as done in \textbf{Case} one for $\textit{(2.6.1)}$.
    
    \bigskip

   \noindent \textbf{Lemma} $\textit{2.6.2.1}$ (\textit{a lower bound for the first crossing probability of \textit{(2.6.2)}}). For the line $\mathscr{L}_1$ appearing before the first \textit{freezing cluster} in the finite-volume strip, one has the strictly positive lower bound,
   
   \begin{align*}
          \textbf{P}^{\xi^{\mathrm{Sloped}}}_{[-m,m]\times[0,n^{\prime}N]} \big[       \mathcal{L}^1_{\mathscr{F}\mathscr{C}_1} \underset{\mathscr{D}_1}{\overset{h\geq ck}{\longleftrightarrow}}  \widetilde{\gamma_L}         \big]   \geq        \mathscr{C}^{\prime}_1   \text{ } \text{ , } \text{ } \\
   \end{align*}

   \bigskip
   
   \noindent \textit{Proof of Lemma 2.6.2.1}. The probability of the connectivity event between the first line and the critical \textit{freezing cluster} occurring will be analyzed through the following inequality, 
   
   \begin{align*}
     \textbf{P}^{\xi^{\mathrm{Sloped}}}_{[-m,m]\times[0,n^{\prime}N]} \big[     \mathcal{L}^1_{\mathscr{F}\mathscr{C}_1} \underset{\mathscr{D}_1}{\overset{h\geq ck}{\longleftrightarrow}}  \widetilde{\gamma_L}       \big]  \text{ } \geq \text{ } \textbf{P}^{\xi^{\mathrm{Sloped}}}_{[-m,m]\times[0,n^{\prime}N]} \big[      \bigcap_{I \in \mathcal{I}_{\mathscr{D}}} \big\{       \mathcal{L}^1_{\mathscr{F}\mathscr{C}_1} \underset{\mathscr{D}_1}{\overset{h\geq ck}{\longleftrightarrow}}   \widetilde{\gamma_I}         \big\}      \big] \text{ }
     \end{align*}
     
     \noindent which admits the lower bound, by $\mathrm{(FKG)}$,
     
     \begin{align*}
     \prod_{I \in \mathcal{I}_{\mathscr{D}}} \text{ }   \textbf{P}^{\xi^{\mathrm{Sloped}}}_{[-m,m]\times[0,n^{\prime}N]} \big[     \mathcal{L}^1_{\mathscr{F}\mathscr{C}_1} \underset{\mathscr{D}_1}{\overset{h\geq ck}{\longleftrightarrow}}   \widetilde{\gamma_I}           \big]     \text{ } \text{ , } \text{ } \\
   \end{align*}
   
   \noindent where connectivity is quantified through the number of faces acquiring the required height for the crossing up to the reflection of the left boundary $\widetilde{\gamma}_{I_{\mathscr{D}}} \equiv \widetilde{\gamma}_I$ with $\mathcal{L}^1_{\mathscr{F}\mathscr{C}_1}$. The index set over which the intersection of crossing events is,
   
   \begin{align*}
         \mathcal{I}_{\mathscr{D}} \text{ } \equiv  \text{ }     \big|            \big\{   u > 0  :  \mathscr{D}_u  \subset \mathscr{D}_1   ,  d  \big(  \mathcal{L}^1_{\mathscr{F}\mathscr{C}_1}  ,  \mathscr{D}_u  \big)  <    d   \big(  \mathcal{L}^1_{\mathscr{F}\mathscr{C}_1} ,  \mathscr{D}_1   \big) \big\}    \big| \text{ } \text{ , } \text{ }  \\
   \end{align*}
   
   \noindent which readily yields the lower bound, 
   
   \begin{align*}
           \prod_{I \in \mathcal{I}_{\mathscr{D}}} \text{ }   \mathfrak{C}_{\mathscr{D}} \geq   \mathfrak{C}^{\prime}_{\mathscr{D}} \big(     | \mathcal{I}_{\mathscr{D}}  |     \big)    \text{ } \text{ , }  \tag{\textit{2.6.2.1.1}}\\
   \end{align*}
   
   \noindent upon setting $\mathscr{C}^{\prime}_1 \equiv \mathfrak{C}^{\prime}_{\mathscr{D}} \big( \text{ }    | \text{ } \mathcal{I}_{\mathscr{D}} \text{ } |    \text{ } \big)$, where each probability occurs with finite energy, as

   \begin{align*}
     \textbf{P}^{\xi^{\mathrm{Sloped}}}_{[-m,m]\times[0,n^{\prime}N]} \big[    \mathcal{L}^1_{\mathscr{F}\mathscr{C}_1} \underset{\mathscr{D}_1}{\overset{h\geq ck}{\longleftrightarrow}}   \widetilde{\gamma_I}              \big] \geq      \mathfrak{C}_{\mathscr{D}}   \text{ } \text{ , } \text{ } 
   \end{align*}
   
   \noindent hence concluding the argument. \boxed{}
   
   \bigskip
   
      \noindent \textbf{Lemma} $\textit{2.6.2.2}$ (\textit{a lower bound for the second crossing probability of \textit{(2.6.2)}}). For the line $\mathscr{L}_1$ appearing before the first \textit{freezing cluster} in the finite-volume strip, one has the lower bound,
   
   \begin{align*}
        \textbf{P}^{\xi_{\mathrm{Sloped}}}_{[-m,m] \times [0,n^{\prime} N]}   \big[         \widetilde{\gamma_L} \underset{(  S    \cap \mathscr{D}_i  )^c}{\overset{h \geq ck}{\longleftrightarrow}} \mathscr{F}\mathscr{C}_{i+1}                                            \big]  \geq     \mathscr{C}^{\prime}_2              \text{ } \text{ , } \text{ }
   \end{align*}
   
   \bigskip
   
   \noindent \textit{Proof of Lemma 2.6.2.2}. As a reflection of \textit{(2.6.1.2 C)} which provided a lower bound for crossing probabilities between supercritical lines in the strip and supercritical \textit{freezing clusters} used in arguments for \textbf{Lemma} $\textit{2.6.1.2}$, we write,
   
   \begin{align*}
        \textbf{P}^{\xi_{\mathrm{Sloped}}}_{[-m,m] \times [0,n^{\prime} N]}   \big[         \widetilde{\gamma_L} \underset{( S   \cap  \mathscr{D}_i   )^c}{\overset{h \geq ck}{\longleftrightarrow}} \mathscr{F}\mathscr{C}_{i+1}                                                  \big]  \text{ } \geq \text{ }    \textbf{P}^{\xi_{\mathrm{Sloped}}}_{[-m,m] \times [0,n^{\prime} N]}   \big[    \bigcap_{k^{\prime} \in \mathcal{I}}       \big\{       \widetilde{\gamma_L}         \underset{(  S    \cap  \mathscr{D}_i  )^c}{ \overset{h \geq ck}{\longleftrightarrow} 
        }\mathscr{L}_{k^{\prime}+1}           \big\}
      \big]    \text{ } \\ \overset{\mathrm{(FKG)}}{\geq} \text{ } \prod_{k^{\prime} \in \mathcal{I}} \text{ }  \textbf{P}^{\xi_{\mathrm{Sloped}}}_{[-m,m] \times [0,n^{\prime} N]}   \big[      \widetilde{\gamma_L}    \underset{(  S    \cap  \mathscr{D}_i  )^c}{      \overset{h \geq ck}{\longleftrightarrow} } \mathscr{L}_{k^{\prime}+1}          \big]   \text{ , } \text{ } 
   \end{align*}
   
   \noindent where $\mathscr{L}_{k^{\prime}+1}$ indicate a series of lines in the strip preceding $\mathscr{F}\mathscr{C}_{i+1}$, with $k^{\prime} > 0$. Finally, applying similar estimates by finite energy, in addition to the number of faces for which the connectivity event with the left reflected boundary of the domain and $\mathscr{L}_{k^{\prime}+1}$ occurs with $\mathfrak{C}^{\prime}$, we conclude the argument upon setting,

   \begin{align*}
         \mathfrak{C}^{\prime}     \equiv \sqrt[|\mathcal{I}|]{\mathfrak{C}^{\prime\prime}}    \text{ } \text{ , } \text{ } \tag{\textit{2.6.2.2.1}} \\
   \end{align*}
   
   \noindent equal to $\mathscr{C}^{\prime}_2$. \boxed{}
   
   \bigskip
   
      \noindent \textbf{Lemma} $\textit{2.6.2.3}$ (\textit{a lower bound for the third crossing probability of \textit{(2.6.2)}}). For the line $\mathscr{L}_1$ appearing before the first \textit{freezing cluster} in the finite-volume strip, one has the strictly positive lower bound,
   
   \begin{align*}
        \textbf{P}^{\xi_{\mathrm{Sloped}}}_{[-m,m] \times [0,n^{\prime} N]}   \big[       \mathscr{F}\mathscr{C}_N \underset{\mathscr{D}_N \cap      \mathcal{L}^N_{\mathscr{F}\mathscr{C}_N}}{\overset{h \geq ck}{\longleftrightarrow}}              \mathcal{L}^N_{\mathscr{F}\mathscr{C}_N}                          \big]     \geq        \mathscr{C}^{\prime}_3  \text{ }  \text{ , } \text{ } \\
   \end{align*}
   
   \bigskip
   
   \noindent \textit{Proof of Lemma 2.6.2.3}. As a reflection of arguments for lower bounding the first probability,
   
   \begin{align*}
       \textbf{P}^{\xi^{\mathrm{Sloped}}}_{[-m,m] \times [0,n^{\prime}N]} \big[    \mathscr{F}\mathscr{C}_i    \underset{\mathscr{D}_N   \cap      \mathcal{L}^N_{\mathscr{F}\mathscr{C}_N}}{   \overset{h \geq ck}{\longleftrightarrow}    }  \mathscr{L}^{\prime}        \big]       \text{ }     \text{ , } \text{ } \\
   \end{align*}
   
  \noindent of \textbf{Lemma} \textit{2.6.1.4} which is a quantification of connectivity between the \textit{freezing cluster} and a line in the complementary region of the strip to the \textit{freezing cluster}, we write, for $\mathcal{L}^N_{\mathscr{F}\mathscr{C}_N}  \equiv\mathcal{L}$, and,

        \begin{align*}
           \mathscr{F}\mathscr{C}_N  \cap    \mathcal{L}  \cap S  \equiv  \big\{ \forall S \subsetneq \textbf{Z}^2 ,   \exists     \mathscr{F} \in F(  S  )  :      \mathscr{F}  \cap  \big(   F_L( \mathscr{F}\mathscr{C}_N  )    \cap         F_R(   \mathcal{L}    ) \big)  \neq \emptyset           \big\}             \text{ } \text{ , } \text{ }   \\
    \end{align*}

    \noindent for the number of faces contained within the finite volume strip appearing within $\mathscr{F}\mathscr{C}_N$ and $\mathcal{L}$, given a partition of the faces in $S$ appearing to the right or left of a \textit{freezing cluster} or line, each of which are respectively denoted with  $F_L(\text{ } \mathscr{F}\mathscr{C}_N \text{ } )$, and with $F_R(\text{ }   \mathcal{L}    \text{ } )$. Across a countable intersection of crossing events,
  
  \begin{align*}
       \textbf{P}^{\xi_{\mathrm{Sloped}}}_{[-m,m] \times [0,n^{\prime} N]}   \big[       \mathscr{F}\mathscr{C}_N \underset{\mathscr{D}_N}{\overset{h \geq ck}{\leftrightarrow}}              \mathcal{L}                   \big]    \geq \textbf{P}^{\xi_{\mathrm{Sloped}}}_{[-m,m] \times [0,n^{\prime} N]}   \big[   \bigcap_{k^{\prime\prime} , j^{\prime\prime} \in \mathcal{I}}  \big\{             \mathscr{L}_{k^{\prime\prime}}    \underset{\mathscr{F}\mathscr{C}_N \text{ } \cap \text{ }   \mathcal{L} \text{ } \cap S }{\overset{h \geq ck}{\longleftrightarrow}}   \mathscr{L}_{j^{\prime\prime}}       \big\}      \big]  \text{ } \text{ , } \text{ } \\
  \end{align*}
  
  \noindent where beginning in the intersection of crossing events preceding the second inequality, we compare the crossing probabilities in the uppermost bound, of height $h \geq ck$ and occurring within the $n$th domain $\mathscr{D}_N$, with crossing probabilities in the second most upper bound, of height $h \geq ck$ and occurring within $\mathscr{F}\mathscr{C}_N \text{ } \cap \text{ }   \mathcal{L} \text{ } \cap S$, where,
  
  \begin{align*}
 S  \equiv    [-m,m] \times [0,n^{\prime} N]  \text{ } \text{ , } \text{ } \\
  \end{align*}
   
   \noindent and the series of lines $\mathscr{L}_{k^{\prime\prime}}$ are configured within the strip so as to appear on the left of $\mathscr{F}\mathscr{C}_N$, in which there can exist some index threshold $k^{\prime\prime}_{\mathrm{crit}} \in \mathcal{I}$ for which $\mathscr{L}_{k^{\prime\prime}_{\mathrm{crit}}} \cap \mathscr{F} \mathscr{C}_N \neq \emptyset$, in addition to another series of lines appearing to the right of $\mathcal{L}$ in the strip in which there can exist some other threshold, $j^{\prime\prime}_{\mathrm{crit}}$, for which $\mathscr{L}_{j^{\prime\prime}_{\mathrm{crit}}} \cap \mathcal{L} \neq \emptyset$. Furthermore, along the lines of \textit{(2.6.1.3 A)} for arguments in \textbf{Lemma} $\textit{2.6.1.3}$, as well as the application of $\textbf{Lemma}$ $\textit{2.6.1.4}$ for lower bounding a product of crossing probabilities by $\mathrm{(FKG)}$ for arguments in $\textbf{Lemma}$ \textit{2.6.1.5},  
   
   \begin{align*}
     \prod_{k^{\prime\prime},j^{\prime\prime} \in\mathcal{I}}  \textbf{P}^{\xi^{\mathrm{Sloped}}}_{[-m,m] \times [0,n^{\prime}N]} \big[   \mathscr{L}_{k^{\prime\prime}}    \underset{\mathscr{F}\mathscr{C}_N  \cap   \mathcal{L}  \cap S }{\overset{h \geq ck}{\longleftrightarrow}}   \mathscr{L}_{j^{\prime\prime}}                               \big]  \geq      \prod_{k^{\prime\prime},j^{\prime\prime} \in \mathcal{I}} \mathfrak{C}^{\prime\prime}  \geq  \big(   \mathfrak{C}^{\prime\prime}   \big)^{|\mathcal{I}|} \equiv \mathfrak{C}^{\prime\prime\prime}  \text{ }  \text{ , } \text{ }  \tag{\textit{2.6.2.3.1}} \\
   \end{align*}

    \noindent also by $\mathrm{(FKG)}$, for strictly positive $\mathfrak{C}^{\prime\prime\prime}$, where the power to which the constant in the lower bound is raised satisfy,
    
    \begin{align*}
          \mathfrak{C}^{\prime\prime}  \leq  {\underset{k^{\prime\prime} , j^{\prime\prime}: k^{\prime\prime} , j^{\prime\prime} \in \textbf{Z}}{\mathrm{inf}}} \text{ } \big\{      \textbf{P}^{\xi^{\mathrm{Sloped}}}_{[-m,m] \times [0,n^{\prime}N]} \big[   \mathscr{L}_{k^{\prime\prime}}    \underset{\mathscr{F}\mathscr{C}_N  \cap   \mathcal{L}  \cap S }{\overset{h \geq ck}{\longleftrightarrow}}   \mathscr{L}_{j^{\prime\prime}}        \big]  
         \big\}                 \text{ } \text{ , } \text{ }   \\
    \end{align*}

   \noindent  from which we conclude the argument upon setting $\mathscr{C}^{\prime}_3 \equiv  \mathfrak{C}^{\prime\prime\prime}$. \boxed{}

    \bigskip

 \noindent Given the existence of finitely many \textit{freezing clusters} for which $\mathscr{F}\mathscr{C}_i  \cap  \mathscr{D}_i  \neq \emptyset$ for each $i$, as well as the boundary of the $i$th \textit{freezing cluster} being defined as the collection of incident faces to the complement of the \textit{freezing cluster} in the strip,

    \begin{align*}
      \partial \mathscr{F} \mathscr{C}_i  \equiv \big\{     \mathscr{F}  \in \mathscr{F} (     \mathscr{D}_i           )  :            ( \widetilde{     \mathscr{F} \mathscr{C}_i    }     ) \cap \big(      \mathscr{F} \mathscr{C}_i    \cap S  \big)^c \text{ } \neq \text{ } \emptyset        \big\}           \text{ }   \text{ , }    \text{ } \tag{\textit{INT}} \\
    \end{align*}
    
    \noindent where $  ( \widetilde{     \mathscr{F} \mathscr{C}_i   }      ) $ denotes the interior of faces within the $i$th \textit{freezing cluster}, in addition to the complement of the intersection of the strip with $\mathscr{F}\mathscr{C}_i$.  Also, given $\emptyset \neq \big(\mathscr{D}_i \cap \mathscr{F}\mathscr{C}_i\big) \subset \mathscr{D}_i$ for each $i$, the probability corresponding to the intersection of vertical crossing events over $\mathcal{I}$, 
    
    \begin{align*}
\bigcap_{i \in \mathcal{I}} \text{ }  \mathcal{V}^{h \leq (1-c)k}_{(  \mathscr{D}_i  \cap  \mathscr{F} \mathscr{C}_i )^c }  =   \bigcap_{i \in \mathcal{I}}   
\text{ } \bigg\{  \{  \mathcal{L}^1_{\mathscr{F}\mathscr{C}_i} \underset{\mathscr{D}_i}{\overset{h\geq ck}{\longleftrightarrow}}  \widetilde{\gamma_L}  \}     \cap    \big\{ \widetilde{\gamma_L} \underset{( \text{ } S  \text{ }  \cap \text{ } \mathscr{D}_i  \text{ } )^c }{\overset{h \geq ck}{\longleftrightarrow}} \mathscr{F}\mathscr{C}_{i+1}   \big\}  \cap  \big\{    \mathscr{F}\mathscr{C}_N \underset{\mathscr{D}_N \text{ }  \cap \text{ }      \mathcal{L}^N_{\mathscr{F}\mathscr{C}_N} \text{ }       }{\overset{h \geq ck}{\longleftrightarrow}} \mathcal{L}^N_{\mathscr{F}\mathscr{C}_N}   \big\}    \bigg\}  \text{ } \text{ , } \text{ } 
    \end{align*}

         \begin{figure}
\begin{align*}
\includegraphics[width=0.88\columnwidth]{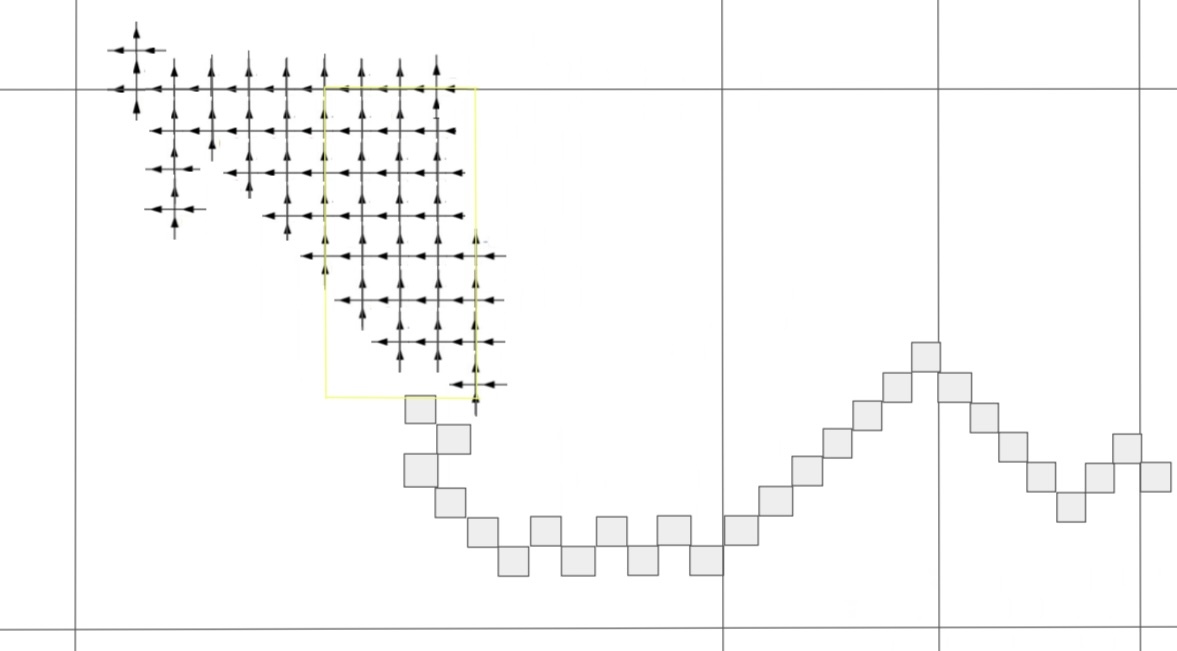}\\
\end{align*}
\caption{\textit{A finite volume configuration with yellow finite subvolume surrounding a boundary of the freezing cluster in the strip}. An horizontal $\mathrm{x}$-crossing between the yellow finite volume boundary and three lines to the left of the \textit{freezing cluster} is depicted, which is a horizontal $\mathrm{x}$-crossing counterpart to a previously shown horizontal crossing between lines to the left of the \textit{freezing cluster}.}
\end{figure}

    \noindent has a lower bound satisfying,
    
    \begin{align*}
     \textbf{P}^{\xi^{\mathrm{Sloped}}}_{[-m,m] \times [0,n^{\prime} N]} \big[     \bigcap_{i \in \mathcal{I}} \text{ }  \mathcal{V}^{h \leq (1-c)k}_{(  \mathscr{D}_i  \cap  \mathscr{F} \mathscr{C}_i )^c }                          \big]                    \overset{\mathrm{(FKG)}}{\geq}  \text{ }        \prod_{i=1}^3    \mathscr{C}^{\prime}_i  \text{ } \equiv     \mathfrak{C}^{\prime}_{\mathscr{D}}                   \mathfrak{C}^{\prime\prime}     \mathfrak{C}^{\prime\prime\prime}          \text{ }  \text{ , } \text{ } 
    \end{align*}
    
    \noindent because each one of the three terms in the product,
    
    \begin{align*}
{\underset{k^{\prime} , k^{\prime\prime},j^{\prime\prime} \in\mathcal{I}}{\prod_{I \in \mathcal{I}_{\mathscr{D}} }}}  \text{ }   \underset{(\textbf{Lemma} \text{ }  \textit{2.6.2.1})}{\underbrace{\textbf{P}^{\xi^{\mathrm{Sloped}}}_{[-m,m]\times[0,n^{\prime}N]} \bigg[         \mathcal{L}^1_{\mathscr{F}\mathscr{C}_1} \underset{\mathscr{D}_1}{\overset{h\geq ck}{\longleftrightarrow}}  \widetilde{\gamma_L}           \bigg] }  }   \text{ }      \underset{(\textbf{Lemma} \text{ }  \textit{2.6.2.2})}{\underbrace{\textbf{P}^{\xi^{\mathrm{Sloped}}}_{[-m,m] \times [0,n^{\prime} N]} \big[                    \widetilde{\gamma_L} \overset{h \geq ck}{\longleftrightarrow} \mathscr{F}\mathscr{C}_{i+1}                                            \big]  }}   \text{ }                           \underset{(\textbf{Lemma} \text{ }  \textit{2.6.2.3})}{\underbrace{\textbf{P}^{\xi^{\mathrm{Sloped}}}_{[-m,m] \times [0,n^{\prime}N]} \big[  \mathscr{L}_{k^{\prime\prime}}    \underset{\mathscr{F}\mathscr{C}_N  \cap  \mathcal{L}  \cap S }{\overset{h \geq ck}{\longleftrightarrow}}   \mathscr{L}_{j^{\prime\prime}}                              \big] }}     \text{ , } 
 \end{align*}

 \noindent can be bound below by gathering constants given in \textbf{Lemma} $\textit{2.6.2.1}$, \textbf{Lemma} $\textit{2.6.2.2}$, and \textbf{Lemma} $\textit{2.6.2.3}$, giving the product lower bound,

 \begin{align*}
 \mathfrak{C}^{\prime}_{\mathscr{D}}                  \mathfrak{C}^{\prime\prime}     \mathfrak{C}^{\prime\prime\prime}          \text{ } \text{ , } \text{ } 
    \end{align*}

    \noindent upon collecting estimates, including $\textit{(2.6.2.1.1)}$ from $\textbf{Lemma}$ \textit{2.6.2.1}, $\textit{(2.6.2.2.1)}$ from $\textbf{Lemma}$ \textit{2.6.2.2}, and $\textit{(2.6.3.2.1)}$ from $\textbf{Lemma}$ \textit{2.6.2.3}, respectively. Concluding, there exists strictly positive $\mathfrak{C}_{\mathscr{D}}^{\prime,\prime\prime,\prime\prime\prime}$ for which,

        \begin{align*}
      \mathfrak{C}^{\prime}_{\mathscr{D}}                   \mathfrak{C}^{\prime\prime}     \mathfrak{C}^{\prime\prime\prime}     \geq       \mathfrak{C}_{\mathscr{D}}^{\prime,\prime\prime,\prime\prime\prime}   \text{ }  \text{ . } 
        \end{align*}

    \noindent Finally, as demonstrated in \textbf{Case one}, the vertical crossing probability for intersecting left and right boundaries of the domain is upper bounded by the segment connectivity event between $\mathcal{I}_0$ and $\widetilde{\mathcal{I}_0}$,
    
    \begin{align*}
        \textbf{P}^{\xi^{\mathrm{Sloped}}}_{[-m,m] \times [0,n^{\prime} N ]} \big[    \mathcal{I}_0      \overset{h \geq (1-c)k}{\longleftrightarrow}       \widetilde{\mathcal{I}_0}                    \big] \text{ }           \text{ } \text{ , } \text{ } 
    \end{align*}
    
    \noindent from which one obtains a lower-bound dependent upon the vertical crossing across the domain, for strictly positive $\mathscr{C}^{\neq \emptyset}$,

    \begin{align*}
     \mathscr{C}^{\neq \emptyset}     \text{ } \textbf{P}^{\chi^{\prime}_{\mathrm{Sloped}}|_{\partial (  [-m,m] \times [0,n^{\prime}_2 N]   ) }}_{\mathscr{D}} \big[    \mathcal{V}^{h \leq (1-c)k}  \big]  \equiv \mathscr{C}^{\neq \emptyset} \textbf{P}^{\xi^{\mathrm{Sloped}}}_{\mathscr{D}} \big[     \mathcal{V}^{h \leq (1-c)k}         \big]  \leq     C_{\mathcal{V}}       \text{ } \text{ . } \text{ } \\
    \end{align*}
    
 \noindent To establish that the second case, \textbf{Case two}, of the result holds for intersecting boundaries of the \textit{symmetric strip domain}, consider the following. If an intersection occurs with the left boundary of the \textit{strip symmetric domain} and some line situated within the strip of the square lattice across which the probability of obtaining a long crossing is quantified, then another symmetric domain, $\mathscr{D}^{\prime}$, can be constructed by performing reflections, which preserve parity of the six-vertex configuration about the line of intersection, so that $\partial \mathscr{D}^{\prime} \cap \mathscr{D} \neq \emptyset$ have the same boundary conditions, in addition to $\partial \mathscr{D}^{\prime\prime} \cap \partial \mathscr{D}^{\prime} \neq \emptyset$ having the same boundary conditions. To this end, for a subdomain $\mathscr{D}^{\prime\prime} \subsetneq \mathscr{D}^{\prime} \subsetneq \mathscr{D}$ satisfying $\mathscr{D} \supsetneq \mathscr{D}^{\prime\prime} \cap \mathscr{D}^{\prime} \neq \emptyset$, introduce the four demarcations,

 \[
F \big( \textbf{Z}^2 \big) \supsetneq \mathscr{D}^{\prime\prime} \equiv \text{ } 
\left\{\!\begin{array}{ll@{}>{{}}l}       \gamma_L     \text{ } \text{ , }   \\  \gamma^{\prime}_L
          \text{ } \text{ , } \\        \mathscr{B} \text{ } \text{ , } \\  \mathscr{B}^{\prime} \text{ } \text{ , } 
\end{array}\right.
\]

\noindent which is taken to the union of faces bound within the union of boundaries $\gamma_L$, $\gamma^{\prime}_L$, $\mathscr{B}$, and $\mathscr{B}^{\prime}$,

\begin{align*}
 \textbf{P}^{\xi^{\mathrm{Sloped}}}_{\mathscr{D}^{\prime\prime}} \big[     \mathcal{I}_0      \overset{h \geq (1-c)k+j}{\longleftrightarrow}       \widetilde{\mathcal{I}_0}             \big] \overset{(\mathrm{CBC})}{\leq}  \textbf{P}^{(\xi^{\mathrm{Sloped})^{\prime}}}_{\mathscr{D}^{\prime\prime}} \big[   \mathcal{I}_0      \overset{h \geq (1-c)k+j}{\longleftrightarrow}       \widetilde{\mathcal{I}_0}      \big]  \leq       1 -       \textbf{P}^{(\xi^{\mathrm{Sloped})^{\prime}}}_{\mathscr{D}^{\prime\prime}} \big[  \gamma_L \overset{ h \geq (1-c)k + j -1}{\longleftrightarrow}  \mathscr{B}^{\prime}  \big] \\ < 1 - \mathfrak{C}^{\prime\prime}   \text{ } \text{ , } 
\end{align*}

\noindent because, for the first term appearing in the lower bound above, 

\begin{align*}
 \textbf{P}^{\xi^{\mathrm{Sloped}}}_{\mathscr{D}^{\prime\prime}} \big[    \mathcal{I}_0      \overset{h \geq (1-c)k+j}{\longleftrightarrow}       \widetilde{\mathcal{I}_0}                    \big]      \overset{(\textbf{Corollary}    \textit{1.2}                 )}{\leq} \text{ } \mathscr{C}^{\neq \emptyset}   \text{ }  \textbf{P}^{\xi^{\mathrm{Sloped}}}_{\mathscr{D}}  \big[               \mathcal{I}_0      \overset{h \geq (1-c)k }{\longleftrightarrow}       \widetilde{\mathcal{I}_0}   \big]  \text{ } \text{ , } \text{ } 
\end{align*}

\noindent after having applied \textbf{Corollary}  \textit{1.2} with $\Lambda \equiv \mathscr{D}$ and $\Lambda^{\prime} \equiv \mathscr{D}^{\prime\prime}$, in addition to the fact that there exists sloped boundary conditions for which,

\begin{align*}
  \textbf{P}^{(\xi^{\mathrm{Sloped})^{\prime}}_1}_{\mathscr{D}^{\prime\prime}}  \big[  \gamma_L \overset{ h \geq (1-c)k + j -1}{\longleftrightarrow}  \mathscr{B}^{\prime}   \big]  \overset{(\mathrm{CBC})}{\geq}   \textbf{P}^{(\xi^{\mathrm{Sloped})^{\prime}}_2}_{\mathscr{D}^{\prime\prime}}  \big[  \gamma_L \overset{ h \geq (1-c)k + j -1}{\longleftrightarrow}  \mathscr{B}^{\prime}   \big]  \\ \overset{(\mathrm{SMP})}{\equiv}      \textbf{P}^{(\xi^{\mathrm{Sloped}}_3)^{\prime}}_{\mathscr{D}^{\prime}}   \big[  \gamma_L \overset{ h \geq (1-c)k + j -1}{\longleftrightarrow}  \mathscr{B}^{\prime} \big| h \equiv ( \xi^{\mathrm{Sloped}}_2 )^{\prime} \text{ over }  \big( \mathscr{D}^{\prime\prime} \big)^c \cup  \partial \mathscr{D}^{\prime\prime} \big] \\ \geq    \textbf{P}^{(\xi^{\mathrm{Sloped})^{\prime}}_1}_{\mathscr{D}^{\prime\prime}}  \big[  \gamma_L \overset{ h \geq (1-c)k + j -1}{\longleftrightarrow}  \mathscr{B}^{\prime}   \big]       \text{ } \text{ , } 
\end{align*}

\noindent where $\mathscr{C}^{\neq \emptyset}$ denotes a constant from \textbf{Corollary} \textit{1.2} for nonintersecting left and right boundaries of the domain, in addition to the domain, which takes the form,

\begin{align*}
 \textbf{P}  \big[  \gamma_L \overset{ h \geq (1-c)k + j -1}{\longleftrightarrow}  \mathscr{B}^{\prime}   \big]  \geq \mathfrak{C}_{\mathscr{D}}^{\prime,\prime\prime,\prime\prime\prime}  > 0     \text{ } \text{ , } 
\end{align*}

\noindent which completes the argument for \textbf{Case two} of \textbf{Lemma} $\textit{2.6}$ with $\gamma_L \cap \gamma_R \neq \emptyset$, in which the upper bound takes the form,

\begin{align*}
    \textbf{P}^{\xi^{\mathrm{Sloped}}}_{\mathscr{D}^{\prime\prime}} \big[      \mathcal{I}_0      \overset{h \geq (1-c)k}{\longleftrightarrow}       \widetilde{\mathcal{I}_0}                  \big]  \leq   C_{\mathcal{V}}       \text{ } \text{ , } \text{ } 
\end{align*}

\noindent from which we conclude the argument. \boxed{}

\subsection{Bridging events in the strip}

\noindent We introduce a modification to bridging events, defined in {\color{blue}[11]} from the quantity, 

\begin{align*}
    \mathcal{B}_{h \geq l}(j)  \equiv  \big\{ \mathcal{I}_{j-1} \overset{h \geq k \text{ in } \textbf{Z} \times [0,n]}{\longleftrightarrow} \mathcal{I}_{j+1}  \big\} \text{ , } \text{ } 
\end{align*}

\noindent which also are defined for the absolute value $|h|$ of the height function instead of only for $h$, with,

\begin{align*}
       \mathcal{B}^{\prime}_{h \geq l}(j)    \equiv  \big\{      \mathcal{I}_{j-1}                                       \overset{h \geq k \text{ in } \textbf{Z} \times [0,n^{\prime} N ]}{\longleftrightarrow}  CC_k (  \mathscr{F}\mathscr{C}_i    )                   \big\}   \cap  \big\{    CC_l (  \mathscr{F} \mathscr{C}_i  )               \overset{h \geq k \text{ in } \textbf{Z} \times [0,n^{\prime} N ]}{\longleftrightarrow}   \mathcal{I}_{j+1}     \big\}   \text{ } \text{ , } \text{ } 
\end{align*}

\noindent in which, with positive probability, there exists blocking connected components, given for $k,l >0$, with $CC_k ( \mathscr{F} \mathscr{C}_i )$ and $CC_l (  \mathscr{F} \mathscr{C}_i )$, respectively, of the \textit{freezing cluster} $\mathscr{F}\mathscr{C}_i$, with any of the faces required for connectivity between $\mathcal{I}_{j-1}$ and $\mathcal{I}_{j+1}$, with the following items. Readily, if there are a countable number $N$ of connected components belonging to the subset of the strip over which the weights do not form infinitely many disjoint circuits, such connected components obstructing $\{   \mathcal{I}_{j-1}        \overset{h \geq ck \text{ in } \textbf{Z} \times [0,n^{\prime}N]}{\longleftrightarrow}       \mathcal{I}_{j+1}      \}$ from occurring imply that the modified bridging event would instead take the form,

\begin{align*}
      \mathcal{B}^{\prime}_{h \geq l}(j)    \equiv \big\{                       \mathcal{I}_{j-1}                        \overset{h \geq k }{\underset{ \mathscr{F} \mathscr{C}^c \cap \Lambda }{\longleftrightarrow}}          \mathcal{I}_{j+1}                       \big\}        \text{ , } \text{ } 
\end{align*}

\noindent which, dependent upon the number of frozen faces in finite volume, can be strictly contained within the intersection of events,

\begin{align*}
 \big\{   \mathcal{I}_{j-1}            \overset{h \geq k \text{ in } \textbf{Z} \times [0,n^{\prime} N ]}{\longleftrightarrow}    CC_1 ( \mathscr{F}\mathscr{C}_i   )                      \big\}    \cap \underset{\text{Intersection with } CC_{i+1} , \cdots , CC_{N-1}}{\underbrace{\cdots}} \text{ } \cap    \text{ } \big\{   CC_N ( \mathscr{F}\mathscr{C}_i   )     \overset{h \geq k \text{ in } \textbf{Z} \times [0,n^{\prime} N ]}{\longleftrightarrow}   \mathcal{I}_{j+1}      \big\} \\  \text{ } \text{ . } 
\end{align*}

\noindent \textbf{Definition} \textit{8} (\textit{neighborhoods of faces surrounding connected components of the freezing cluster}). Denote,

\begin{align*}
           \mathfrak{N}_{k^{\prime},k,i}   \equiv  \mathfrak{N}_{k^{\prime}} \big(   CC_k  (   \mathscr{F}\mathscr{C}_i          )    \big)  \equiv      \big\{     \mathscr{F}  \in          F (     S  \cap \mathscr{D}_i     )                :             \big|              \mathscr{F}                \overset{h \geq ck}{\longleftrightarrow}                       CC_k ( \mathscr{F} \mathscr{C}_i)                      \big|    <   k                      \big\}         \text{ } \text{ , } \text{ } 
\end{align*}

\noindent as the \textit{neighborhood} of strictly positive number of faces surrounding $CC_k( \text{ } \mathscr{F} \mathscr{C}_i \text{ } )$, that lie within a strictly $k^{\prime}$ distance of the \textit{freezing cluster} boundary, as given in \textit{(INT)}, following \textit{(2.6.2.3.1)}, which is indexed by $i$. As expected, the boundary of such a neighborhood is given by,

\begin{align*}
      \partial    \mathfrak{N}_{k^{\prime},k,i}        \equiv     \partial \big(     CC_k \big( \mathscr{F}\mathscr{C}_i              \big) \big)  \equiv \big\{         \mathscr{F}  \in         F (    S \cap  \mathscr{D}_i      )            :             \big|          \mathscr{F}                  \overset{h \geq ck}{\longleftrightarrow}                       CC_k ( \mathscr{F} \mathscr{C}_i  )                   \big|    =  k             \big\}               \text{ } \text{ . } \text{ } 
\end{align*}

\noindent Finally, the interior is denoted with $\widetilde{\mathfrak{N}_{k^{\prime},k,i}} \equiv \big(     \partial  \mathfrak{N}_{k^{\prime},k,i}  \big)^c \cap \mathfrak{N}_{k^{\prime},k,i}$.

\bigskip

 \begin{figure}
\begin{align*}
\includegraphics[width=0.85\columnwidth]{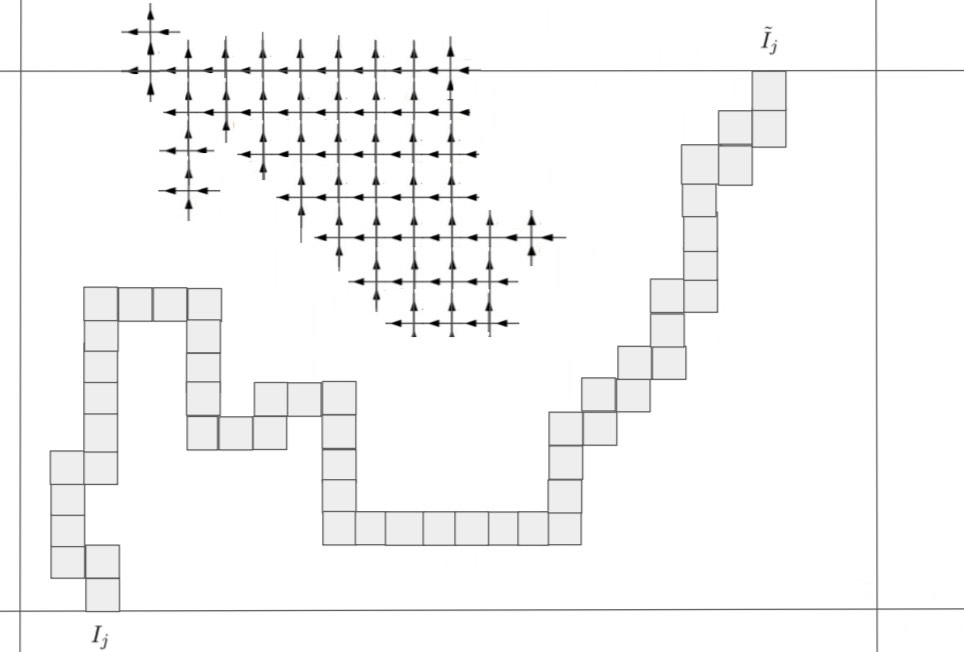}\\\includegraphics[width=0.85\columnwidth]{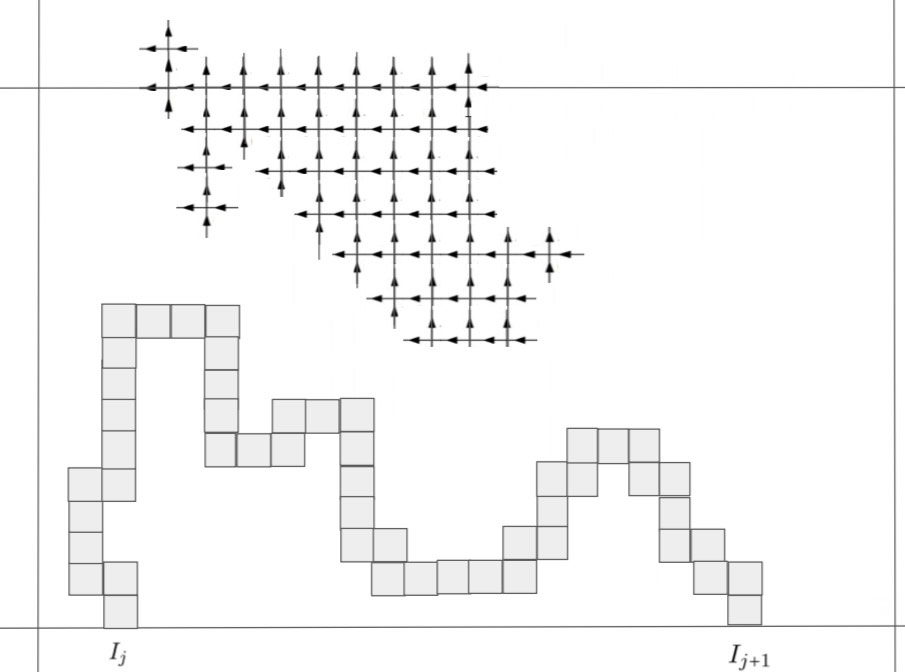}\\
\end{align*}
\caption{\textit{(TOP PANEL). A depiction of connected components of faces between $\mathcal{I}_j$ and $\widetilde{\mathcal{I}_j}$ avoiding faces belonging to the \textit{cluster} of frozen faces}. \textit{(BOTTOM PANEL). Connected component of face connectivity event between $\mathcal{I}_j$ and $\mathcal{I}_{j+1}$, where the latter interval is displaced some sufficiently small distance away from $\mathcal{I}_j$.}}
\end{figure}

\noindent \textbf{Lemma} \textit{2.7.4} (\textit{simultaneously incorporating estimates from the previous three lemmas}). The intersection of connectivity events within the face neighborhood and the complement from $\textbf{Lemma}$ \textit{2.7.1}, and $\textbf{Lemma}$ \textit{2.7.2}, appears in a lower bound over the portion of the strip for which the modified bridging event $\mathcal{B}^{\prime}$ in $\textbf{Definition}$ \textit{7} occurs satisfies,

\begin{align*}
   \textbf{P}^{\xi^{\mathrm{Sloped}}}_{\textbf{Z} \times [0,n^{\prime} N]} \big[                         \mathcal{B}^{\prime}_{h \geq l}(j)       \big]  \geq     \mathcal{C}_{\mathcal{B}} 
   \big(   \textbf{P}^{\xi^{\mathrm{Sloped}}}_{\textbf{Z} \times [0,n^{\prime}N]}     \big[                             \mathcal{I}_0 \overset{ h \geq ck}{\longleftrightarrow} \widetilde{\mathcal{I}_0}                     \big]   \big)^2                  \text{ } \text{ , } \text{ } 
\end{align*}

\noindent for a suitably chosen, strictly positive, constant.

\bigskip


\noindent \textit{Proof of Lemma 2.7.4}. We begin with comparing the height crossing, irrespective of the absolute value of the height function,

\begin{align*}
      \textbf{P}^{\xi^{\mathrm{Sloped}}}_{\textbf{Z} \times [0,n^{\prime} N]} \big[                       \mathcal{B}^{\prime}_{h \geq l}(j)    \big] \text{ } 
     \end{align*}
     
   \noindent with,
     
     \begin{align*}
      \textbf{P}^{\xi^{\mathrm{Sloped}}}_{\textbf{Z} \times [0,n^{\prime} N]} \big[                       \mathcal{B}^{\prime}_{\text{ } | h  | \geq l}(j)   \big]   \text{ } \text{ , } 
     \end{align*}
     
     \noindent corresponding to the crossing event dependent upon $|h|$. For a lower bound of the crossing event conditionally upon the absolute value of the height function, the event measurable from $|h|$ over all faces contained within finite volume takes the form,
     
     \begin{align*}
         \textbf{P}^{\xi^{\mathrm{Sloped}}}_{\textbf{Z} \times [0,n^{\prime}N]} \big[             \mathcal{B}^{\prime}_{|h|\geq l }        \big|           \mathscr{C}         \big]         \text{ } \text{ , } 
     \end{align*}

\noindent where the conditional event $\mathscr{C}$, from previous descriptions of domains for sloped boundary conditions in the strip, is dependent upon the number of faces along the left and right boundaries of the domain for which $h \geq ck$. Conditionally upon the existence of leftmost, and rightmost, top to bottom crossings in the strip, the leftmost crossing must satisfy that none of the faces do not intersect any faces contained within the \textit{freezing cluster}, in which,

\begin{align*}
    \forall \text{ } \mathscr{F} \subset \underset{\text{ Leftmost crossing for inducing } \gamma_L \text{ of the domain}}{\underbrace{\big\{   \mathcal{I}_{-j}          \overset{h \geq ck}{\longleftrightarrow}      \widetilde{ \mathcal{I}_j } \big\}}}      : \mathscr{F}  \cap  \mathscr{F} \mathscr{C}_i \equiv  \emptyset        \text{ , } \text{ } 
\end{align*}

\noindent for some natural $j$, while similarly, the rightmost crossing must also satisfy that none of the faces contained within the left boundary do not intersect any of the faces contained in another \textit{freezing cluster},

\begin{align*}
       \forall \text{ } \mathscr{F} \subset \underset{\text{ Rightmost crossing for inducing } \gamma_R \text{ of the domain}}{\underbrace{\big\{   \mathcal{I}_{-k}          \overset{h \geq ck}{\longleftrightarrow}     \text{ }  \widetilde{ \mathcal{I}_k } \big\}}} \text{ }     : \text{ } \mathscr{F} \text{ } \cap \text{ } \mathscr{F} \mathscr{C}_i \text{ } \equiv \text{ } \emptyset        \text{ , } \text{ }         \text{ } \text{ . } \text{ }
\end{align*}

\noindent for some natural $k > j$. Explicitly,

\begin{align*}
 \mathscr{C} \equiv      \big\{  \mathcal{I}_{-j}          \overset{h \geq ck}{\longleftrightarrow}     \widetilde{ \mathcal{I}_j }   \big\}        \cap  \big\{    \mathcal{I}_{-k}          \overset{h \geq ck}{\longleftrightarrow}     \widetilde{ \mathcal{I}_k }    \big\}           \text{ } \text{ . } 
\end{align*}

\noindent Consequently, $\textbf{P}^{\xi^{\mathrm{Sloped}}}_{\textbf{Z} \times [0,n^{\prime} N]} [     \mathscr{C}     ]$ admits the following probability as a lower bound,

\begin{align*}
  \textbf{P}^{\xi^{\mathrm{Sloped}}}_{\textbf{Z} \times [0,n^{\prime}N]} \big[  \big\{    \mathcal{I}_{-j}          \overset{h \geq ck}{\longleftrightarrow}   \widetilde{ \mathcal{I}_j }  \big\}\cap   \big\{   \mathcal{I}_{-k}          \overset{h \geq ck}{\longleftrightarrow}      \widetilde{ \mathcal{I}_k }  \big\}      \big]  \overset{\mathrm{(FKG)}}{\geq}       \textbf{P}^{\xi^{\mathrm{Sloped}}}_{\textbf{Z} \times [0,n^{\prime}N]} \big[    \mathcal{I}_{-j}          \overset{h \geq ck}{\longleftrightarrow}     \widetilde{ \mathcal{I}_j } \big]   \\ \times       \textbf{P}^{\xi^{\mathrm{Sloped}}}_{\textbf{Z} \times [0,n^{\prime}N]} \big[      \mathcal{I}_{-k}          \overset{h \geq ck}{\longleftrightarrow}       \widetilde{ \mathcal{I}_k }  \big]   \text{ } \text{ . } \end{align*}

\noindent The following sequence of inequalities above is equivalent to the fact that the following conditional bridging probability admits the lower bound $\mathcal{C}_{\mathcal{B}}$, in which, 

  \begin{align*}
      \textbf{P}^{\xi^{\mathrm{Sloped}}}_{\textbf{Z} \times [0,n^{\prime}N]} \big[         \mathcal{B}^{\prime}_{\text{ } |h|\text{} \geq l}(j)                     \big|   
      \big\{           \mathcal{I}_{-j}          \overset{h \geq ck}{\longleftrightarrow}      \widetilde{ \mathcal{I}_j }            \big\}      \cap    \big\{  \mathcal{I}_{-k}          \overset{h \geq ck}{\longleftrightarrow}      \widetilde{ \mathcal{I}_k}   \big\}  \big]  \geq  \mathcal{C}_{\mathcal{B}}        \text{ } \text{ , } \text{ } 
\end{align*}

\noindent which implies that the probability in the upper bound form the last inequality can be analyzed by decomposing the intersection of three crossing events appearing to lower bound the previous probability, as from the sequence of inequalities below,

\begin{align*}
  \textbf{P}^{\xi^{\mathrm{Sloped}}}_{\textbf{Z} \times [0,n^{\prime} N]} \big[   \mathcal{B}^{\prime}_{ | h  | \geq c_0 k }(j)       | |h|_{\Lambda^{\mathrm{Out}}}        \big]  \geq \textbf{P}^{\xi^{\mathrm{Sloped}}}_{\textbf{Z} \times [0,n^{\prime} N]} \big[               \mathcal{H}_{h \geq c_0 k} (  \mathscr{D}_i  )   |       |h|_{\Lambda^{\mathrm{Out}}}     \big]  \tag{\textit{2.7.4.1}} \\ \geq         \textbf{P}^{\xi^{\mathrm{Sloped}}}_{\textbf{Z} \times [0,n^{\prime} N]} \big[   \mathcal{H}_{ |h|   \geq c_0 k} (  \mathscr{D}_i  )  |        |h|_{\Lambda^{\mathrm{Out}}}  \equiv  \xi        \big] \tag{\textit{2.7.4.2}} \\ \geq  \textbf{P}^{\xi^{\mathrm{Sloped}}}_{\textbf{Z} \times [0,n^{\prime} N]}  \big[     \mathcal{H}_{ h   \geq c_0 k} ( \mathscr{D}_i  )  |      |h|_{\Lambda^{\mathrm{Out}}}  \equiv  \xi       \big]     \text{ } \text{ , } \text{ } \tag{\textit{2.7.4.3}}
\end{align*}


\noindent where in the lower bound given by $\textit{(2.7.4.1)}$, the probability of obtaining a conditional horizontal crossing event across $\mathscr{D}$ is smaller than the probability of a conditional crossing with $\mathcal{B}^{\prime}_{ | h |  \geq c_0 k}$ occurring, while the two subsequent lower bounds in \textit{(2.7.4.2)}, and in \textit{(2.7.4.3)}, are obtained from conditioning that the height function satisfy $\xi \sim \textbf{BC}^{\mathrm{Sloped}}$, in which, 

\begin{align*}
  \xi  \equiv   \big\{ \forall   i  \in \textbf{Z}  , \exists \mathscr{F}_i  \in   F(         \textbf{Z}^2  \cap  \Lambda       \cap    \mathscr{D}_i     )       :  \mathscr{F}_i  \equiv \ h_i  \big\}          \text{ } \text{ , } \\
\end{align*}

\noindent and that $\mathcal{H}_{h \geq c_0 k } \subset \mathcal{H}_{ | h | \geq c_0 k}$, respectively. To remove the conditioning upon the absolute value of the height function, observe,

\begin{align*}
    \textbf{P}^{\xi^{\mathrm{Sloped}}}_{\textbf{Z} \times [0,n^{\prime} N]}  \big[     \mathcal{H}_{ h   \geq c_0 k} ( \mathscr{D}_i )   |        |h|_{\Lambda^{\mathrm{Out}}}  \equiv  \xi         \big]    \geq  \textbf{P}^{\xi^{\mathrm{Sloped}}}_{\textbf{Z} \times [0,n^{\prime} N]}  \big[ \big\{      \mathcal{H}_{ h   \geq c_0 k} ( \mathscr{D}_i  )   \big\}  \cap  \big\{       |h|_{\Lambda^{\mathrm{Out}}}  \equiv \xi     \big\}    \big]  \text{ } \\  \geq \text{ }    \textbf{P}^{\xi^{\mathrm{Sloped}}}_{\textbf{Z} \times [0,n^{\prime} N]}  \big[   \mathcal{H}_{ h  \geq c_0 k} (\mathscr{D}_i )  \big]  \textbf{P}^{\xi^{\mathrm{Sloped}}}_{\textbf{Z} \times [0,n^{\prime} N]}  \big[  |h|_{\Lambda^{\mathrm{Out}}}  \equiv  \xi   \big]  \\ \geq \text{ }                          \text{ }  \textbf{P}^{\xi^{\mathrm{Sloped}}}_{\textbf{Z} \times [0,n^{\prime} N]}  \big[    \mathcal{H}_{ h  \geq c_0 k} (  \mathscr{D}_i  )  \big]        \text{ } \text{ , } \tag{\textit{2.7.4.4}}
\end{align*}

\noindent where in \textit{(2.7.4.4)}, as for flat boundary conditions, the absolute value of the height function for the six-vertex model satisfies a ferromagnetic Ising model, hence admitting the lower bound, 

\begin{align*}
        \textbf{P}^{\xi^{\mathrm{Sloped}}}_{\textbf{Z} \times [0,n^{\prime} N]}  \big[  |h|_{\Lambda^{\mathrm{Out}}} \text{ } \equiv \text{ } \xi \big] \text{ } \geq \text{ }            4^{-1}       \text{ } \text{ . } 
\end{align*}

\noindent To lower bound the horizontal crossing probability across $\mathscr{D}_1$, we make use of the following series of results.

\bigskip

\noindent The aim of such estimates is to put together the three following estimates to ensure, across the regions of the finite strip, that crossing events between $\mathcal{I}_0$ and a sufficiently large neighborhood of faces surrounding each \textit{freezing cluster} occur with positive probability, and also that crossing events between $\widetilde{\mathcal{I}_0}$ and another sufficiently large neighborhood of another \textit{freezing cluster} occur with positive probability.

       \begin{figure}
\begin{align*}
\includegraphics[width=1.05\columnwidth]{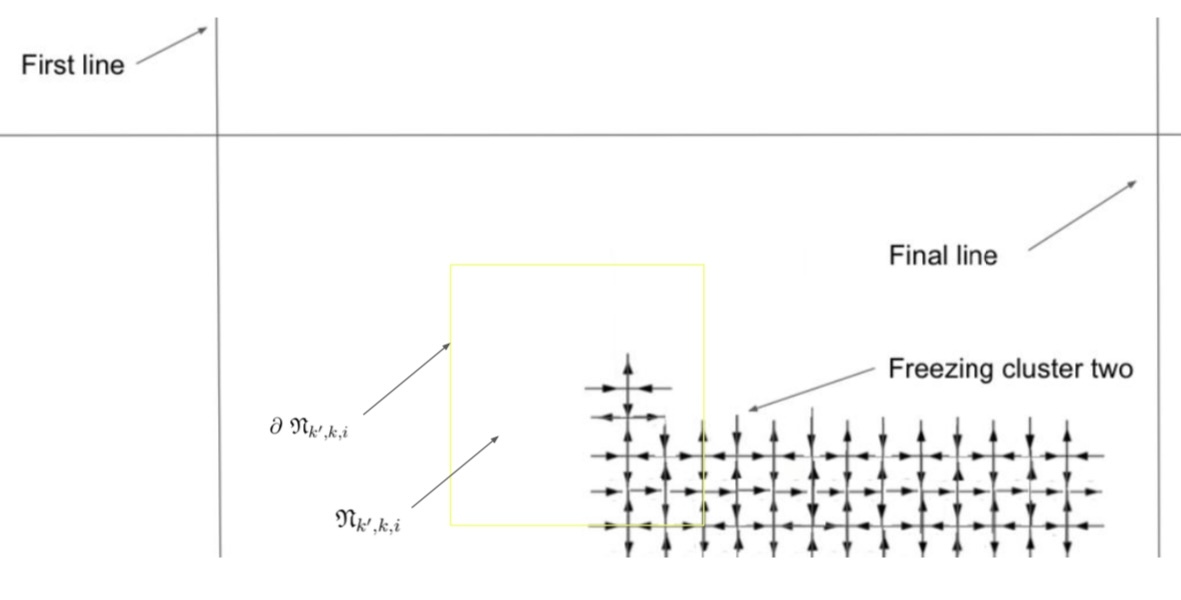}\\
\end{align*}
\caption{\textit{A depiction of the collection of faces bound within $   \partial \mathfrak{N}_{k^{\prime},k,i}$, and within $\mathfrak{N}_{k^{\prime},k,i}$ embedded within the strip.} With a portion of the \textit{freezing cluster} from the previous \textit{figure}, the boundary of the yellow finite volume in the strip indicates the boundary of the face neighborhood, while the interior of the finite volume is also labeled. Connectivity is quantified through the probability of $\big\{  \partial \mathfrak{N}_{k^{\prime},k,i}  \overset{h \geq ck}{\longleftrightarrow}  \mathfrak{N}_{k^{\prime},k,i}  \big\} \in \Omega^{\mathrm{Sloped}}$ occurring, within the interior of the \textit{face neighborhood} intersected $\big(  \Lambda \cap \mathscr{F}\mathscr{C}_i \big)\subset \textbf{Z}^2$.}
\end{figure}

\bigskip

\noindent \textbf{Proposition} \textit{2.7.4.1} (\textit{horizontal crossing probability across strip domains}). From the probability

\noindent $\textbf{P}^{\xi^{\mathrm{Sloped}}}_{[-m,m]\times[0,n^{\prime}N]} [   \mathcal{H}(\mathscr{D}_i) ]$ in \textit{(2.7.4.4)}, the crossing event can be decomposed into the intersection of crossings,

\begin{align*}
 \mathscr{I}   \equiv   \big\{     \mathcal{I}_0 \underset{\Lambda \cap  ( \text{ }   \mathfrak{N}_{k^{\prime},k,i}          \text{ } )^c   \cap \mathscr{F}\mathscr{C}_i}{\overset{h \geq ck}{\longleftrightarrow}}        \partial  \text{ }  \mathfrak{N}_{k^{\prime},k,i}                                          \big\}     \cap \big\{        \partial    \mathfrak{N}_{k^{\prime},k,i}                           \underset{\widetilde{\mathfrak{N}_{k^{\prime},k,i}}}{\overset{h \geq ck}{\longleftrightarrow}}       \mathfrak{N}_{k^{\prime},k,i}            \big\}   \cap\big\{                     \partial  \mathfrak{N}_{k^{\prime},k,i}       \underset{\mathscr{D}_i \text{ } \cap  (  \mathscr{F} \mathscr{C}_i )^c }{\overset{h \geq ck}{\longleftrightarrow}}      \widetilde{\mathcal{I}_0}        \big\}       \text{ }    \text{ , } \text{ } 
\end{align*}

\noindent with the corresponding lower bound,

\begin{align*}
  \textbf{P}^{\xi^{\mathrm{Sloped}}}_{\textbf{Z} \times [0,n^{\prime}N]} \big[  \mathscr{I}  \big]  \geq     \underset{o  \in  \{ \mathfrak{n}, \mathfrak{n}^{\prime} , \mathfrak{n}^{\prime\prime}  \} }{\prod}  \mathcal{C}_{o}      \equiv   \mathcal{C}_{\mathfrak{n}} \mathcal{C}_{\mathfrak{n}^{\prime}}  \mathcal{C}_{\mathfrak{n}^{\prime\prime\prime}}  \text{ . } 
\end{align*}

\bigskip

\noindent \textit{Proof of Proposition 2.7.4.1}. We isolate each crossing event through the following series of results.

\bigskip

\noindent \textbf{Lemma} \textit{2.7.1} (\textit{strict positivity of connectivity between the boundary of the face neighborhood and the complement of the face neighborhood}). The probability of connectivity occurring between $\partial  \mathfrak{N}_{k^{\prime},k,i}$ and any face in $\Lambda \cap \mathfrak{N}_{k,k^{\prime},i}$, within $\Lambda \cap  ( \mathfrak{N}_{k^{\prime},k,i}       )^c   \cap \mathscr{F}\mathscr{C}_i$ satisfies,

\begin{align*}
 \textbf{P}^{\xi^{\mathrm{Sloped}}}_{[-m,m] \times [0,n^{\prime} N]} \big[    \mathcal{I}_0 \underset{\Lambda \cap  (   \mathfrak{N}_{k^{\prime},k,i}        )^c   \cap \mathscr{F}\mathscr{C}_i}{\overset{h \geq ck}{\longleftrightarrow}}        \partial   \mathfrak{N}_{k^{\prime},k,i}                                                     \big]    \geq   \mathcal{C}_{\mathfrak{N}}  \text{ } \text{ , } \text{ } 
\end{align*}

\noindent for a suitably chosen, strictly positive, constant.

\bigskip

\noindent \textit{Proof of Lemma 2.7.1}. From typical argumentation,

\begin{align*}
   \textbf{P}^{\xi^{\mathrm{Sloped}}}_{[-m,m] \times [0,n^{\prime} N]} \big[     \mathcal{I}_0 \underset{\Lambda \cap  ( \mathfrak{N}_{k^{\prime},k,i}          )^c   \cap \mathscr{F}\mathscr{C}_i}{\overset{h \geq ck}{\longleftrightarrow}}        \partial  \mathfrak{N}_{k^{\prime},k,i}                                                  \big]   \geq   \textbf{P}^{\xi^{\mathrm{Sloped}}}_{[-m,m] \times [0,n^{\prime} N]} \big[  \bigcap_{\mathfrak{n} \in \mathcal{I}}  \big\{     \partial  \mathfrak{N}_{\mathfrak{n} , k , i}     \underset{\Lambda \cap  (    \mathfrak{N}_{k^{\prime},k,i}           )^c   \cap \mathscr{F}\mathscr{C}_i}{\overset{h \geq ck}{\longleftrightarrow}}            \partial  \mathfrak{N}_{\mathfrak{n}+1,k,i}                         \big\}  \big] \text{ } \\ \overset{\mathrm{(FKG)}}{\geq}  \prod_{\mathfrak{n} \in \mathcal{I}} \textbf{P}^{\xi^{\mathrm{Sloped}}}_{[-m,m] \times [0,n^{\prime} N]} \big[               \partial \mathfrak{N}_{\mathfrak{n}^{\prime} , k , i}     \underset{\Lambda \cap  (  \mathfrak{N}_{k^{\prime},k,i}          )^c   \cap \mathscr{F}\mathscr{C}_i}{\overset{h \geq ck}{\longleftrightarrow}}            \partial \mathfrak{N}_{\mathfrak{n}^{\prime}+1,k,i}     \big] \text{ } \\  \geq   \prod_{\mathfrak{n} \in \mathcal{I}}          \mathcal{C}_{\mathfrak{n}} \equiv \mathcal{C}_{\mathfrak{N}}               \text{ } \text{ , } \text{ } 
\end{align*}

\noindent from which we conclude the argument, for the smallest, and largest, respectively, indices $\mathfrak{n}^{0}$ and $\mathfrak{n}^{| \text{ } \mathcal{I}_{\mathfrak{n}} \text{ } |}$ from $\mathcal{I}$, given $\partial \text{ } \mathfrak{N}_{\mathfrak{n}^{0},k,i} \subset S$ such that $\partial \text{ } \mathfrak{N}_{\mathfrak{n}^{0},k,i} \text{ } \cap \text{ }  \mathcal{I}_0 \neq \emptyset$, and $\partial \text{ } \mathfrak{N}_{\mathfrak{n}^{| \text{ } \mathcal{I}_{\mathfrak{n}} \text{ } |},k,i} \subset S$ such that $\partial \text{ } \mathfrak{N}_{\mathfrak{n}^{| \text{ } \mathcal{I}_{\mathfrak{n}} \text{ } |},k,i} \text{ } \cap \text{ } \partial \text{ } \mathfrak{N}_{k^{\prime},k,i} \neq \emptyset$, with $\partial \text{ } \mathfrak{N}_{\mathfrak{n}^{\prime}, k , i} \text{ } \cap \text{ } \partial \text{ } \mathfrak{N}_{\mathfrak{n}^{\prime}+1,k,i} \neq \emptyset$, for every $\mathfrak{n}^{\prime} > 0$. \boxed{}

\bigskip

\noindent \textbf{Lemma} \textit{2.7.1.1} (\textit{strict positivity of connectivity between the interior and boundary of the face neighborhood}). The probability of connectivity occurring between $\mathfrak{N}_{k^{\prime},k,i}$, and $\partial    \text{ }  \mathfrak{N}_{k^{\prime},k,i}$, satisfies,

\begin{align*}
     \textbf{P}^{\xi^{\mathrm{Sloped}}}_{[-m,m] \times [0,n^{\prime} N]} \big[                        \partial      \mathfrak{N}_{k^{\prime},k,i}                           \underset{\widetilde{\mathfrak{N}_{k^{\prime},k,i}}}{\overset{h \geq ck}{\longleftrightarrow}}       \mathfrak{N}_{k^{\prime},k,i}                                           \big]     \geq               \mathcal{C}_{\mathfrak{N}^{\prime}}               \text{ } \text{ , } \text{ } \\
\end{align*}

\noindent for a suitably chosen, strictly positive, constant.

\bigskip

\noindent \textit{Proof of Lemma 2.7.1.1}. From typical argumentation,

\begin{align*}
      \textbf{P}^{\xi^{\mathrm{Sloped}}}_{[-m,m] \times [0,n^{\prime} N]} \big[              \partial    \mathfrak{N}_{k^{\prime},k,i}         \underset{\widetilde{\mathfrak{N}_{k^{\prime},k,i}}}{\overset{h \geq ck}{\longleftrightarrow}}                  \mathfrak{N}_{k^{\prime},k,i}                                                            \big]  \geq     \textbf{P}^{\xi^{\mathrm{Sloped}}}_{[-m,m] \times [0,n^{\prime} N]} \big[    \bigcap_{\mathfrak{n}^{\prime} \in \mathcal{I}} \big\{          \partial    \mathfrak{N}_{n^{\prime}+1,k,i}                  \underset{\widetilde{\mathfrak{N}_{n^{\prime},k,i}}}{\overset{h \geq ck}{\longleftrightarrow}}            \mathfrak{N}_{k^{\prime},k,i}     \big\}                                                        \big]   \text{ } \\ \overset{\mathrm{(FKG)}}{\geq}     \prod_{\mathfrak{n}^{\prime} \in \mathcal{I}}       \textbf{P}^{\xi^{\mathrm{Sloped}}}_{[-m,m] \times [0,n^{\prime} N]} \big[     \partial      \mathfrak{N}_{n^{\prime}+1,k,i}                    \underset{\widetilde{\mathfrak{N}_{k^{\prime},k,i}}}{\overset{h \geq ck}{\longleftrightarrow}}        \mathfrak{N}_{n^{\prime},k,i}                                         \big] \geq      \prod_{\mathfrak{n}^{\prime} \in \mathcal{I}}     \mathcal{C}_{\mathfrak{n}^{\prime}}        \equiv \mathcal{C}_{\mathfrak{N}^{\prime}}     \text{ } \text{ , } \text{ } \\
\end{align*}

\noindent from which we conclude, in which from the intersection of crossing events given above, for the smallest, and largest, respectively, indices $(\mathfrak{n}^{\prime})^0$ and $(\mathfrak{n}^{\prime})^{|  \mathcal{I}_{\mathfrak{n}^{\prime}} |} $ from $\mathcal{I}$, given $\partial \mathfrak{N}_{(\mathfrak{n}^{\prime})^0,k,i} \subset S$ with $\partial \mathfrak{N}_{(\mathfrak{n}^{\prime})^0,k,i} \text{ } \cap \partial   \mathfrak{N}_{k^{\prime},k,i} \neq \emptyset$, and $\mathfrak{N}_{n^{\prime},k,i} \subset S$, for every $\mathfrak{n}^{\prime\prime} > 0$. \boxed{}

\bigskip

\noindent \textbf{Lemma} \textit{2.7.1.2} (\textit{strict positivity of connectivity between the boundary of the face neighborhood and the final interval $\widetilde{\mathcal{I}_0}$ in the segment connectivity event}). The probability of connectivity occurring between $\partial \text{ } \mathfrak{N}_{k^{\prime},k,i}$ and $\widetilde{\mathcal{I}_0}$ satisfies,

\begin{align*}
     \textbf{P}^{\xi^{\mathrm{Sloped}}}_{[-m,m] \times [0,n^{\prime} N]} \big[    \partial  \mathfrak{N}_{k^{\prime},k,i}       \underset{\mathscr{D}_i \cap  (  \mathscr{F} \mathscr{C}_i  )^c }{\overset{h \geq ck}{\longleftrightarrow}}      \widetilde{\mathcal{I}_0}  \big]  \geq   \mathcal{C}_{\mathfrak{N}^{\prime\prime}}    \text{ } \text{ , } \text{ } 
\end{align*}

\noindent for a suitably chosen, strictly positive, constant.

\bigskip

\noindent \textit{Proof of Lemma 2.7.1.2}. From typical argumentation,

\begin{align*}
     \textbf{P}^{\xi^{\mathrm{Sloped}}}_{[-m,m] \times [0,n^{\prime} N]} \big[    \partial \mathfrak{N}_{k^{\prime},k,i}       \underset{\mathscr{D}_i  \cap  (  \mathscr{F} \mathscr{C}_i \text{ } )^c }{\overset{h \geq ck}{\longleftrightarrow}}      \widetilde{\mathcal{I}_0}   \big]  \geq     \textbf{P}^{\xi^{\mathrm{Sloped}}}_{[-m,m] \times [0,n^{\prime} N]} \big[     \bigcap_{\mathfrak{n}^{\prime\prime} \in \mathcal{I}} \big\{       \mathfrak{N}_{n^{\prime\prime},k,i}     \underset{\mathscr{D}_i  \cap (  \mathscr{F} \mathscr{C}_i  )^c }{\overset{h \geq ck}{\longleftrightarrow}}   \mathfrak{N}_{\mathfrak{n}^{\prime\prime}+1,k,i}     \big\}                \big] \text{ } \\  \overset{\mathrm{(FKG)}}{\geq} \prod_{\mathfrak{n}^{\prime\prime} \in \mathcal{I}}   \textbf{P}^{\xi^{\mathrm{Sloped}}}_{[-m,m] \times [0,n^{\prime} N]} \big[  \mathfrak{N}_{\mathfrak{n}^{\prime\prime},k,i}     \underset{\mathscr{D}_i  \cap  (  \mathscr{F} \mathscr{C}_i  )^c }{\overset{h \geq ck}{\longleftrightarrow}}   \mathfrak{N}_{\mathfrak{n}^{\prime\prime}+1,k,i}    \big] \geq  \prod_{\mathfrak{n}^{\prime\prime} \in \mathcal{I}}  \mathcal{C}_{\mathfrak{n}^{\prime\prime}} \equiv  \mathcal{C}_{\mathfrak{N}^{\prime\prime}} \text{ }  \text{ , } \text{ } 
\end{align*}

\noindent from which we conclude the argument. \boxed{}

\bigskip

\noindent To conclude the proof of \textbf{Proposition} \textit{2.7.4.1}, combining the three previous estimates yields,

\begin{align*}
       \textbf{P}^{\xi^{\mathrm{Sloped}}}_{[-m,m] \times [0,n^{\prime}N]} \big[  \mathscr{I}    \big]  \overset{\mathrm{(FKG)}}{\geq}   \mathcal{C}_{\mathfrak{n}}     \mathcal{C}_{\mathfrak{n}^{\prime}}  \mathcal{C}_{\mathfrak{n}^{\prime\prime\prime}}   \text{ . } \tag{\textit{2.7.1}}
\end{align*}

\noindent while, on the other hand, besides the horizontal crossing estimate in $\textit{(2.7.4.4)}$, for vertical crossings,

\begin{align*}
   \textbf{P}^{\xi^{\mathrm{Sloped}}}_{\mathscr{D}_i} [    \mathcal{V} (  \mathscr{D}_i  )         ]  \leq      \textbf{P}^{(\xi^{\mathrm{Sloped}})_2}_{\textbf{Z} \times [0,n^{\prime}]} [           \mathcal{V} ( \mathscr{D}_i  )          ]     \tag{\textit{2.7.4.5}} \text{ } \text{ , }  \end{align*}
   
   \noindent which can be further upper bounded with the probability,
   
   \begin{align*}
   \textbf{P}^{(\xi^{\mathrm{Sloped}})_2}_{\textbf{Z} \times [0,n^{\prime}]} \big[  [ - \delta n , 2 \delta n ]     \times \{ 0 \}  {\overset{h \geq ck}{\longleftrightarrow}}     [ - \delta n , 2 \delta n ]             \times \{ n \}     \big]           \text{ } \text{ , }      \tag{\textit{2.7.4.6}} 
\end{align*}

\noindent in which, for \textit{(2.7.4.5)}, the probability of a vertical crossing occurring is less than the probability of the same vertical crossing event occurring, in which, 

\begin{align*}
     \textbf{P}^{\xi^{\mathrm{Sloped}}}_{[-m,m] \times [0,n^{\prime}N]} \big[         \mathcal{B}^{\prime}_{ |h| \geq l}(j)                      \text{ } \big| \text{ }  
      \big\{             \mathcal{I}_{-j}          \overset{h \geq ck}{\longleftrightarrow}      \widetilde{ \mathcal{I}_j }             \big\}   \cap   \big\{   \mathcal{I}_{-k}          \overset{h \geq ck}{\longleftrightarrow}     \widetilde{ \mathcal{I}_k}   \big\}  \big] \text{ } \\  \geq \text{ }    \textbf{P}^{\xi^{\mathrm{Sloped}}}_{[-m,m] \times [0,n^{\prime}N]} \big[  \mathcal{B}^{\prime}_{ |h|\text{} \geq l}(j)      \cap   \big\{               \mathcal{I}_{-j}          \overset{h \geq ck}{\longleftrightarrow}    \widetilde{ \mathcal{I}_j }             \big\}   \cap   \big\{   \mathcal{I}_{-k}          \overset{h \geq ck}{\longleftrightarrow}       \widetilde{ \mathcal{I}_k}  \big\}     \big]  \\ \text{ } \overset{\mathrm{(FKG)}}{\geq} \text{ }     \textbf{P}^{\xi^{\mathrm{Sloped}}}_{[-m,m] \times [0,n^{\prime}N]} \big[  \mathcal{B}^{\prime}_{\text{ } |h| \geq l}(j)    \big]         \textbf{P}^{\xi^{\mathrm{Sloped}}}_{[-m,m] \times [0,n^{\prime}N]} \big[            \mathcal{I}_{-j}          \overset{h \geq ck}{\longleftrightarrow}     \widetilde{ \mathcal{I}_j }          \big]         \textbf{P}^{\xi^{\mathrm{Sloped}}}_{[-m,m] \times [0,n^{\prime}N]}        \big[   \mathcal{I}_{-k}          \overset{h \geq ck}{\longleftrightarrow}   \widetilde{ \mathcal{I}_k}   \big]  \text{ . } 
\end{align*}

\noindent To bound the product of probabilities above obtained from $\mathrm{(FKG)}$, it suffices to provide estimates for each of the following two quantities, as,

   \begin{figure}
\begin{align*}
\includegraphics[width=0.52\columnwidth]{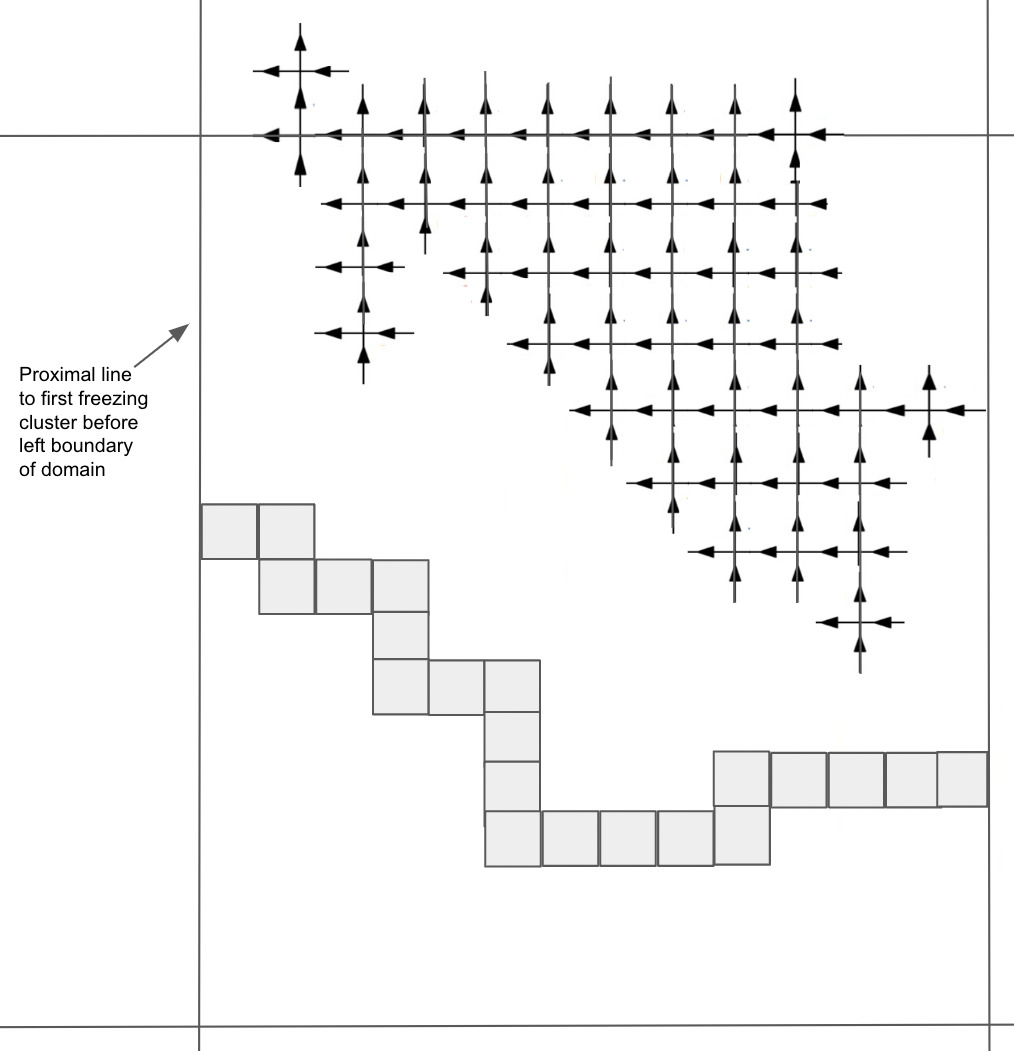}\\
\end{align*}
\caption{\textit{Configuration of the freezing cluster with the strip in finite volume}. Outside of the faces that are bound within the macroscropic \textit{freezing cluster} which intersects the boundary of the strip, a possible crossing outside of the \textit{freezing cluster}, within $\Lambda \cap  \big( \text{ } \mathscr{F}\mathscr{C}_i\text{ } \big)^c$, can be obtained for faces that are highlighted in gray above.}
\end{figure}

\begin{align*}
 \textbf{P}^{\xi^{\mathrm{Sloped}}}_{[-m,m] \times [0,n^{\prime}N]} \big[          \mathcal{B}^{\prime}_{|h| \geq l}(j)   \big] \text{ }     \text{ , } 
\end{align*}

\noindent will be shown to admit the lower bound as provided in the statement,

\begin{align*}
   \text{ }      \mathcal{C}_{\mathcal{B}}  \big(                           \textbf{P}^{\xi^{\mathrm{Sloped}}}_{[-m,m] \times [0,n^{\prime} N]} [        \mathcal{I}_0 \overset{h \geq ck}{\longleftrightarrow} \widetilde{\mathcal{I}_0} ]         \big)^2                                                                                     \text{ } \text{ . } 
\end{align*}

\noindent To achieve such an estimate, observe, first,

\begin{align*}
   \textbf{P}^{\xi^{\mathrm{Sloped}}}_{[-m,m] \times [0,n^{\prime}N]} \big[            \mathcal{B}^{\prime}_{ |h| \geq l}(j)   \big]  \geq   \textbf{P}^{\xi^{\mathrm{Sloped}}}_{[-m,m] \times [0,n^{\prime}N]}    \bigg[ \big\{  \mathcal{I}_0 \underset{\Lambda \cap  (   \mathfrak{N}_{k^{\prime},k,i}       )^c   \cap \mathscr{F}\mathscr{C}_i}{\overset{h \geq ck}{\longleftrightarrow}}        \partial  \mathfrak{N}_{k^{\prime},k,i}  \big\}  \cap  \big\{             \partial     \mathfrak{N}_{k^{\prime},k,i}                           \underset{\widetilde{\mathfrak{N}_{k^{\prime},k,i}}}{\overset{h \geq ck}{\longleftrightarrow}}       \mathfrak{N}_{k^{\prime},k,i}        \big\}  \\ \cap \big\{       \partial  \mathfrak{N}_{k^{\prime},k,i}       \underset{\mathscr{D}_i  \cap ( \mathscr{F} \mathscr{C}_i  )^c }{\overset{h \geq ck}{\longleftrightarrow}}      \widetilde{\mathcal{I}_0}    \big\}     \bigg]     \\ \overset{\mathrm{(FKG)}}{\geq} \textbf{P}^{\xi^{\mathrm{Sloped}}}_{[-m,m] \times [0,n^{\prime}N]}    \big[  \mathcal{I}_0 \underset{\Lambda \cap  ( \text{ }   \mathfrak{N}_{k^{\prime},k,i}          \text{ } )^c   \cap \mathscr{F}\mathscr{C}_i}{\overset{h \geq ck}{\longleftrightarrow}}        \partial  \text{ }  \mathfrak{N}_{k^{\prime},k,i} \big] \text{ }  \text{ } \textbf{P}^{\xi^{\mathrm{Sloped}}}_{[-m,m] \times [0,n^{\prime}N]} \big[            \partial    \mathfrak{N}_{k^{\prime},k,i}                           \underset{\widetilde{\mathfrak{N}_{k^{\prime},k,i}}}{\overset{h \geq ck}{\longleftrightarrow}}       \mathfrak{N}_{k^{\prime},k,i}      \big] \\ \times   \textbf{P}^{\xi^{\mathrm{Sloped}}}_{[-m,m] \times [0,n^{\prime}N]} \big[        \partial \mathfrak{N}_{k^{\prime},k,i}       \underset{\mathscr{D}_i \cap  (  \mathscr{F} \mathscr{C}_i \text{ } )^c }{\overset{h \geq ck}{\longleftrightarrow}}      \widetilde{\mathcal{I}_0}        \big]    \tag{\textit{2.7.1} I}     \text{ , } 
\end{align*}

\noindent will be shown to admit a lower bound dependent upon the probability of height $\geq ck$,

\begin{align*}
      \textbf{P}^{\xi^{\mathrm{Sloped}}}_{[-m,m] \times [0,n^{\prime} N]} [        \mathcal{I}_0 \overset{h \geq ck}{\longleftrightarrow} \widetilde{\mathcal{I}_0} ]       \text{ }  \text{ . } 
\end{align*}

\noindent Proceeding from (\textit{2.7.1} I), 

\begin{align*}
 (\textit{2.7.1} \text{ } \mathrm{I}) \equiv   \textbf{P}^{\xi^{\mathrm{Sloped}}}_{[-m,m] \times [0,n^{\prime}N]}    \big[  \mathcal{I}_0 \underset{\Lambda \cap  (    \mathfrak{N}_{k^{\prime},k,i}          )^c   \cap \mathscr{F}\mathscr{C}_i}{\overset{h \geq ck}{\longleftrightarrow}}        \partial    \mathfrak{N}_{k^{\prime},k,i} \big]  \textbf{P}^{\xi^{\mathrm{Sloped}}}_{[-m,m] \times [0,n^{\prime}N]} \big[            \partial    \mathfrak{N}_{k^{\prime},k,i}                           \underset{\widetilde{\mathfrak{N}_{k^{\prime},k,i}}}{\overset{h \geq ck}{\longleftrightarrow}}       \mathfrak{N}_{k^{\prime},k,i}      \big] \\ \times    \textbf{P}^{\xi^{\mathrm{Sloped}}}_{[-m,m] \times [0,n^{\prime}N]} \big[        \partial  \mathfrak{N}_{k^{\prime},k,i}       \underset{\mathscr{D}_i \cap  ( \mathscr{F} \mathscr{C}_i )^c }{\overset{h \geq ck}{\longleftrightarrow}}      \widetilde{\mathcal{I}_0}        \big]   \\ \overset{(\textbf{Lemma} \textit{ 2.7.1})}{\geq} \text{ }     \mathcal{C}_{\mathfrak{N}}                   \textbf{P}^{\xi^{\mathrm{Sloped}}}_{[-m,m] \times [0,n^{\prime}N]}    \big[  \mathcal{I}_0 \underset{\Lambda \cap  (    \mathfrak{N}_{k^{\prime},k,i}          )^c   \cap \mathscr{F}\mathscr{C}_i}{\overset{h \geq ck}{\longleftrightarrow}}        \partial  \text{ }  \mathfrak{N}_{k^{\prime},k,i} \big]      \textbf{P}^{\xi^{\mathrm{Sloped}}}_{[-m,m] \times [0,n^{\prime}N]} \big[        \partial  \mathfrak{N}_{k^{\prime},k,i}       \underset{\mathscr{D}_i  \cap  (  \mathscr{F} \mathscr{C}_i  )^c }{\overset{h \geq ck}{\longleftrightarrow}}      \widetilde{\mathcal{I}_0}        \big] \text{ , }  \end{align*}

  \noindent which can be further analyzed by observing, for the probability with the connectivity event involving $\mathcal{I}_0$,

  \begin{align*}
   \textbf{P}^{\xi^{\mathrm{Sloped}}}_{[-m,m] \times [0,n^{\prime}N]}    \big[  \mathcal{I}_0 \underset{\Lambda \cap  (  \mathfrak{N}_{k^{\prime},k,i}       )^c   \cap \mathscr{F}\mathscr{C}_i}{\overset{h \geq ck}{\longleftrightarrow}}        \partial    \mathfrak{N}_{k^{\prime},k,i} \big]    \geq  \textbf{P}^{\xi^{\mathrm{Sloped}}}_{[-m,m] \times [0,n^{\prime}N]}    \bigg[   \big\{ \mathcal{I}_0 \underset{\Lambda \cap  (   \mathfrak{N}_{k^{\prime},k,i}           )^c   \cap \mathscr{F}\mathscr{C}_i}{\overset{h \geq ck}{\longleftrightarrow}}        \partial   \mathfrak{N}_{k^{\prime},k,i} \big\}  \\ \cap     \big\{        \widetilde{\mathcal{I}_0} \underset{\Lambda  \cap  (            \mathscr{D}_i \cap  (  \mathscr{F} \mathscr{C}_i  )^c )}{\overset{h \geq ck}{\longleftrightarrow }    }                 \partial \mathfrak{N}_{k^{\prime},k,i}      \big\}    \bigg] \\ \overset{\mathrm{(FKG)}}{\geq} \textbf{P}^{\xi^{\mathrm{Sloped}}}_{[-m,m] \times [0,n^{\prime}N]}    \big[     \mathcal{I}_0 \underset{\Lambda   \cap   (    \mathfrak{N}_{k^{\prime},k,i}           )^c   \cap \mathscr{F}\mathscr{C}_i}{\overset{h \geq ck}{\longleftrightarrow}}        \partial    \mathfrak{N}_{k^{\prime},k,i}       \big]   \textbf{P}^{\xi^{\mathrm{Sloped}}}_{[-m,m] \times [0,n^{\prime}N]}  \big[               \widetilde{\mathcal{I}_0} \underset{\Lambda  \cap  \big(            \mathscr{D}_i \cap (  \mathscr{F} \mathscr{C}_i )^c \big)}{\overset{h \geq ck}{\longleftrightarrow }    }                 \partial \mathfrak{N}_{k^{\prime},k,i}                                        \big]                                      \\ \overset{(\textbf{Lemma}\text{ } \textit{2.7.1.1})}{\geq}   \mathcal{C}_{\mathfrak{N}^{\prime}}      \textbf{P}^{\xi^{\mathrm{Sloped}}}_{[-m,m] \times [0,n^{\prime}N]}  \big[                  \mathcal{I}_0\underset{\Lambda  \cap  (  \mathfrak{N}_{k^{\prime},k,i}          )^c   \cap \mathscr{F}\mathscr{C}_i}{\overset{h \geq ck}{\longleftrightarrow}}  \widetilde{\mathcal{I}_0}                                        \big]      \text{ } 
   \\ \geq  \mathcal{C}_{\mathfrak{N}^{\prime}}   \textbf{P}^{\xi^{\mathrm{Sloped}}}_{[-m,m] \times [0,n^{\prime}N]}     \big[           \mathcal{I}_0 \underset{\Lambda }{\overset{h \geq ck}{\longleftrightarrow}}  \widetilde{\mathcal{I}_0}   \big]      \\ \equiv \text{ }  \mathcal{C}_{\mathfrak{N}^{\prime}}  \textbf{P}^{\xi^{\mathrm{Sloped}}}_{[-m,m] \times [0,n^{\prime}N]}     \big[         \mathcal{I}_0 {\overset{h \geq ck}{\longleftrightarrow}}  \widetilde{\mathcal{I}_0}             \big]      \text{ }     \text{ , } \text{ } \tag{\textit{2.7.1} II}
      \end{align*}
  
  \noindent while, for the probability with the connectivity event involving $\widetilde{\mathcal{I}_0}$, similarly, observe,
  
   \begin{align*}
    \textbf{P}^{\xi^{\mathrm{Sloped}}}_{[-m,m] \times [0,n^{\prime}N]} \big[        \partial  \mathfrak{N}_{k^{\prime},k,i}       \underset{\mathscr{D}_i  \cap  (  \mathscr{F} \mathscr{C}_i )^c }{\overset{h \geq ck}{\longleftrightarrow}}      \widetilde{\mathcal{I}_0}        \big]  \geq     \textbf{P}^{\xi^{\mathrm{Sloped}}}_{[-m,m] \times [0,n^{\prime}N]} \bigg[        \big\{ \partial  \mathfrak{N}_{k^{\prime},k,i}       \underset{\mathscr{D}_i  \cap  (  \mathscr{F} \mathscr{C}_i  )^c }{\overset{h \geq ck}{\longleftrightarrow}}      \widetilde{\mathcal{I}_0} \big\}\\ \cap \big\{        \mathcal{I}_0 \overset{h \geq ck}{\underset{\Lambda   \cap \big(            \mathscr{D}_i  \cap  (  \mathscr{F} \mathscr{C}_i )^c \big) }{\longleftrightarrow} }               \partial \mathfrak{N}_{k^{\prime},k,i}        \big\}      \bigg] \text{ } \\ \overset{\mathrm{(FKG)}}{\geq}        \textbf{P}^{\xi^{\mathrm{Sloped}}}_{[-m,m] \times [0,n^{\prime}N]}     \big[      \partial  \mathfrak{N}_{k^{\prime},k,i}       \underset{\mathscr{D}_i  \cap ( \mathscr{F} \mathscr{C}_i  )^c }{\overset{h \geq ck}{\longleftrightarrow}}      \widetilde{\mathcal{I}_0}                      \big]    \textbf{P}^{\xi^{\mathrm{Sloped}}}_{[-m,m] \times [0,n^{\prime}N]}     \big[              \mathcal{I}_0 \overset{h \geq ck}{\underset{\Lambda   \cap \big(            \mathscr{D}_i  \cap (  \mathscr{F} \mathscr{C}_i  )^c \big) }{\longleftrightarrow}}                \partial \mathfrak{N}_{k^{\prime},k,i}                 \big]  \text{ } \\ \overset{(\textbf{Lemma}\text{ } \textit{2.7.1.2})}{\geq}      \mathcal{C}_{\mathfrak{N}^{\prime\prime\prime}}      \textbf{P}^{\xi^{\mathrm{Sloped}}}_{[-m,m] \times [0,n^{\prime}N]}     \big[      \mathcal{I}_0 \overset{h \geq ck}{\underset{\Lambda \cap \big(            \mathscr{D}_i \cap  (  \mathscr{F} \mathscr{C}_i  )^c \big) }{\longleftrightarrow}}  \widetilde{\mathcal{I}_0}                       \big]       \\\geq     \mathcal{C}_{\mathfrak{N}^{\prime\prime\prime}}   \textbf{P}^{\xi^{\mathrm{Sloped}}}_{[-m,m] \times [0,n^{\prime}N]}     \big[           \mathcal{I}_0 \overset{h \geq ck}{\underset{\Lambda }{\longleftrightarrow} } \widetilde{\mathcal{I}_0}            \big]         \text{ }      \\ \equiv       \mathcal{C}_{\mathfrak{N}^{\prime\prime\prime}}   \textbf{P}^{\xi^{\mathrm{Sloped}}}_{[-m,m] \times [0,n^{\prime}N]}     \big[       \mathcal{I}_0 \overset{h \geq ck}{\longleftrightarrow}  \widetilde{\mathcal{I}_0}           \big]      \text{ . } \text{ } \tag{\textit{2.7.1} III}
  \end{align*}
  
  \noindent Hence,

\begin{figure}
    \centering
    \includegraphics[width=0.84\columnwidth]{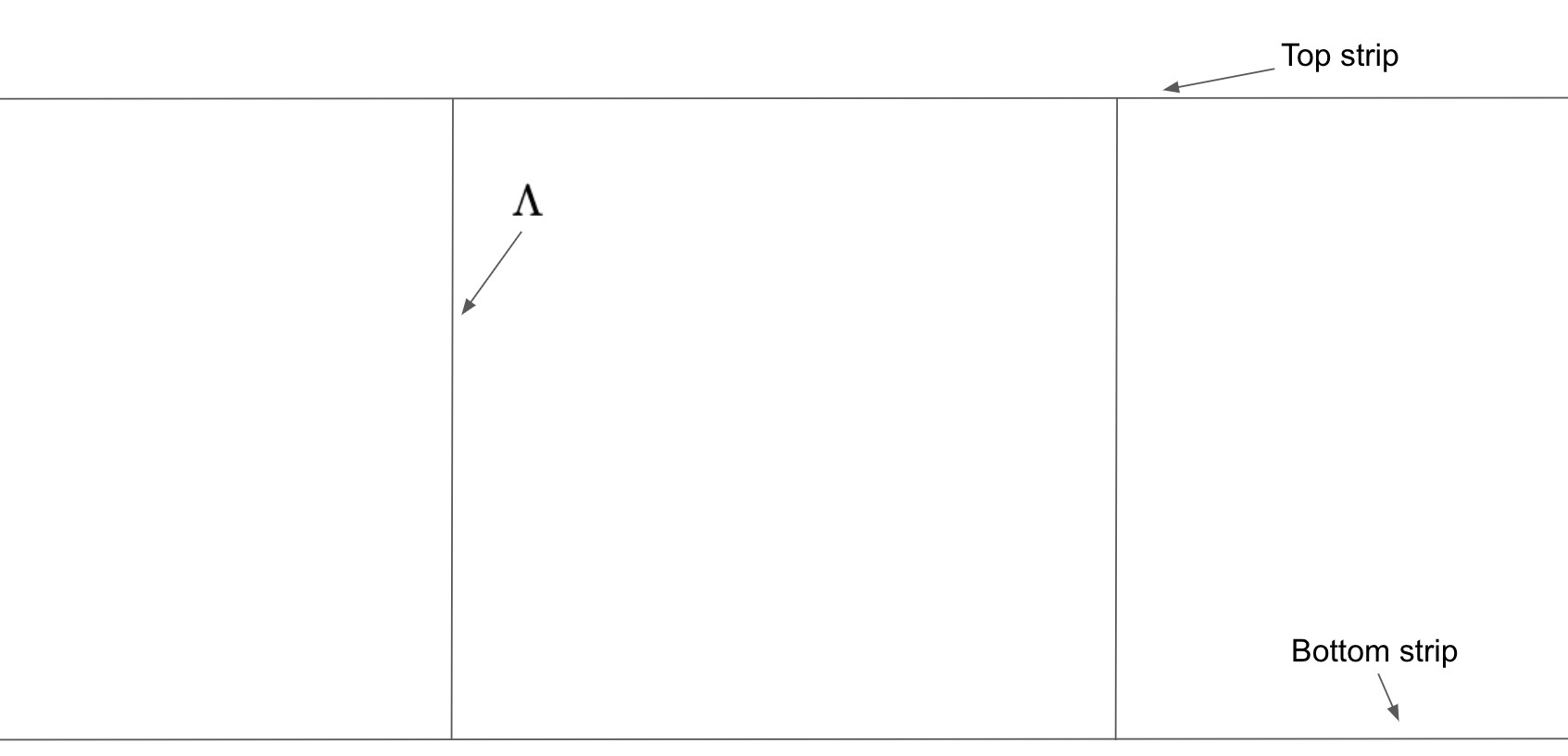}
    \caption{\textit{Configuration of finite volumes for arguments involving} $\Lambda$ \textit{and complements of freezing clusters}. From a finite-volume restriction of the strip with $\Lambda$, depending upon the slope of boundary conditions on incident faces to the strip, with positive probability there will exist a frozen distribution of six-vertex weights within the interior of $\Lambda$. Also within the interior of $\Lambda$, fix $\mathscr{D}_i$, $\mathscr{F}\mathscr{C}_i$, $\mathfrak{N}_{k^{\prime},k,i}$, and $\widetilde{\mathfrak{N}_{k^{\prime},k,i}}$.}
\end{figure}

  \begin{align*}
                   \text{ }    \textbf{P}^{\xi^{\mathrm{Sloped}}}_{[-m,m] \times [0,n^{\prime} N]} \big[                         \mathcal{B}^{\prime}_{h \geq l}(j)    \big]       \text{ } \overset{m \rightarrow+ \infty}{\geq} \text{ } \textbf{P}^{\xi^{\mathrm{Sloped}}}_{\textbf{Z} \times [0,n^{\prime} N]} \big[                         \mathcal{B}^{\prime}_{h \geq l}(j)    \big]    \text{ }   \overset{(\textit{2.7.1}\text{ } \mathrm{I})}\geq \big(  \mathcal{C}_{\mathfrak{n}^{\prime}} \text{ } \mathcal{C}_{\mathfrak{n}^{\prime\prime\prime}}      \big)^{-1} \text{ }  \textbf{P}^{\xi^{\mathrm{Sloped}}}_{\textbf{Z} \times [0,n^{\prime} N]} \big[  \mathscr{I}   \cap \text{ } \widetilde{\mathscr{I}}  \big]       \\  \overset{(\textit{2.7.1})\text{ },\text{ } (\textit{2.7.1}\text{ } \mathrm{II})\text{ } , \text{ } (\textit{2.7.1}\text{ } \mathrm{III})\text{ } , \text{ } \mathrm{(FKG)}}{\geq} \text{ }     \mathcal{C}_{\mathfrak{N}} \text{ } \mathcal{C}_{\mathfrak{N}^{\prime}} \text{ } \mathcal{C}_{\mathfrak{N}^{\prime\prime\prime}}  \text{ }     \big(      \textbf{P}^{\xi^{\mathrm{Sloped}}}_{\textbf{Z} \times [0,n^{\prime}N]}     \big[                         \mathcal{I}_0 \underset{ h \geq ck}{\longleftrightarrow} \widetilde{\mathcal{I}_0}                      \big]   \big)^2      \text{ }  \\  \geq \mathcal{C}_B \big(   \textbf{P}^{\xi^{\mathrm{Sloped}}}_{\textbf{Z} \times [0,n^{\prime}N]}     \big[                           \mathcal{I}_0 \underset{ h \geq ck}{\longleftrightarrow} \widetilde{\mathcal{I}_0}                    \big] \big)^2  \text{ , } 
\end{align*}

\noindent for a strictly positive constant with upper bound given by the product of three constants, respectively given in \textbf{Lemma} \textit{2.7.1}, \textbf{Lemma} \textit{2.7.1.1}, and \textbf{Lemma} \textit{2.7.1.2}, taken over $\mathfrak{N}$, $\mathfrak{N}^{\prime}$, and $\mathfrak{N}^{\prime\prime}$,

\begin{align*}
      \mathcal{C}_{\mathfrak{N}} \mathcal{C}_{\mathfrak{N}^{\prime}}  \mathcal{C}_{\mathfrak{N}^{\prime\prime\prime}}   \geq  C_B      \text{ , } 
\end{align*}

\noindent and,

\begin{align*}
  \widetilde{\mathscr{I}} \equiv \big\{       \mathcal{I}_0 \underset{\Lambda   ( \text{ }   \mathfrak{N}_{k^{\prime},k,i}          \text{ } )^c   \cap \cap   \mathscr{F}\mathscr{C}_i}{\overset{h \geq ck}{\longleftrightarrow}}        \partial    \mathfrak{N}_{k^{\prime},k,i}             \big\} \cap  \big\{             \partial \mathfrak{N}_{k^{\prime},k,i}       \underset{\mathscr{D}_i \cap  (  \mathscr{F} \mathscr{C}_i  )^c }{\overset{h \geq ck}{\longleftrightarrow}}      \widetilde{\mathcal{I}_0}               \big\}   \text{ , } 
\end{align*}

\noindent completing the argument. \boxed{}
 
 \bigskip
 
 \noindent Next we provide arguments for the following item.

 \subsection{RSW type estimate before transferring crossing probability estimates to an annulus of the square lattice from the strip}

\noindent \textbf{Proposition} \textit{2.7.5} (\textit{crossing probability estimate between two segment connectivity events}). For some fixed $c$ as in previous results, and $n^{\prime} > n > 0$, the probability, under sloped boundary conditions, of a segment connectivity event occurring over scale $[-n^{\prime} , 2n^{\prime}]$, satisfies,

\begin{align*}
  \textbf{P}^{\xi^{\mathrm{Sloped}}}_{\textbf{Z} \times [- n^{\prime} , 2n^{\prime}]} \bigg[            [-n^{\prime} , 2n^{\prime}]         \times \{ 0\}   \underset{\textbf{Z} \times [ 0 ,n^{\prime} N]}{\overset{h \geq ck}{\longleftrightarrow}}              [-n^{\prime} , 2n^{\prime}]             \times \{ \rho  \lfloor \delta^{\prime} n \rfloor \}  \bigg] \geq \cdots \cdots \cdots  \\  \cdots \cdots   \bigg( \sqrt[\mathscr{N}]{\mathscr{O}}  \textbf{P}^{\xi^{\mathrm{Sloped}}}_{\textbf{Z} \times [- n^{\prime} , 2n^{\prime}]}  \bigg[    [0 , \lfloor \delta^{\prime} n \rfloor ] \times \{0 \}     \underset{\textbf{Z} \times [ 0 ,n^{\prime} N]}{\overset{h \geq ck}{\longleftrightarrow}}    \textbf{Z}   \times \{ n\}         \bigg] \bigg)^{\text{ } \mathscr{N}}  \text{ , } 
\end{align*}

\noindent for $\rho >0$, and $\delta^{\prime}$ as previously mentioned, $\mathscr{O} < 1$, and strictly positive $\mathscr{N}$.

\bigskip

\noindent \textit{Proof of Proposition 2.7.5}. Under the existence of modified bridging events in $[-m,m] \times [0,n^{\prime} N]$, for $\rho > 0$,

\begin{align*}
    \textbf{P}^{\xi^{\mathrm{Sloped}}}_{[-m,m] \times [0,n^{\prime} N]} \big[   \{ 0 \} \times [n^{\prime} , 2 n^{\prime} ] \overset{ h \geq ck}{\longleftrightarrow } \{ \rho \delta n \}   \times [n^{\prime} , 2 n^{\prime} ] \big]  \geq \textbf{P}^{\xi^{\mathrm{Sloped}}}_{[-m,m] \times [0,n^{\prime} N]} \big[    \{ 0 \} \times [n^{\prime} , 2 n^{\prime} ] \overset{ h \geq ck}{\longleftrightarrow } \{ \rho \delta n \}   \times [n^{\prime} , 2 n^{\prime}  ]    \big] \text{ } \\  \overset{m \rightarrow + \infty}{\geq } \textbf{P}^{\xi^{\mathrm{Sloped}}}_{\textbf{Z} \times [0,n^{\prime} N]} \big[        \{ 0 \} \times [n^{\prime} , 2 n^{\prime} ] \overset{ h \geq ck}{\longleftrightarrow } \{ \rho \delta n \}   \times [n^{\prime} , 2 n^{\prime}  ]      \big] \text{ , } 
\end{align*}

\noindent further implying that the ultimate lower bound from above satisfies,

\begin{align*}
     \textbf{P}^{\xi^{\mathrm{Sloped}}}_{\textbf{Z} \times [0,n^{\prime} N]} \big[         \{ 0 \} \times [n^{\prime} , 2 n^{\prime} ] \overset{ h \geq ck}{\longleftrightarrow } \{ \rho \delta n \}   \times [n^{\prime} , 2 n^{\prime}  ]     \big] \geq          \big(     \textbf{P}^{\xi^{\mathrm{Sloped}}}_{\textbf{Z} \times [0,n^{\prime} N]} \big[    \mathcal{B}^{\prime}_{h \geq l}(j)         \big]        \big)^{\rho} \text{ , } 
\end{align*}

\noindent because,

\begin{align*}
    \big(   \mathcal{B}^{\prime}_{h \geq l}(j)   \big)^{\rho} \text{ } \subsetneq \text{ }    \underset{ i \geq 0}{\bigcap_{\rho_i < \rho }} \big\{      \{ 0 \} \times [n^{\prime} , 2 n^{\prime} ] \overset{ h \geq ck}{\longleftrightarrow } \{ \rho_i \delta n \}   \times [n^{\prime} , 2 n^{\prime}  ]      \big\} \text{ }                \text{ }   \text{ . } 
\end{align*}

\noindent As a result, to show that a lower bound of the form,

\begin{align*}
          \textbf{P}^{\xi^{\mathrm{Sloped}}}_{\textbf{Z} \times [-n^{\prime}, 2n^{\prime}]} \big[  \{ 0 \}   \times [ n^{\prime} , 2 n^{\prime}  ] \overset{ h \geq ck}{\longleftrightarrow} \{n \}  \times   [ n^{\prime} , 2 n^{\prime}  ]    \big]      \geq        \mathscr{O} \big(  \textbf{P}^{\xi^{\mathrm{Sloped}}}_{[-n^{\prime}, 2n^{\prime}]} \big[    \mathcal{I}_0  \overset{h \geq ck}{ \longleftrightarrow } \widetilde{\mathcal{I}_0} \big]         \big)^{\mathscr{N}^{\prime}}      \text{ , } \tag{$\mathscr{O}$\text{ }-\text{ }$\mathscr{N}$ \textit{bound}}
\end{align*}

\noindent holds for some strictly positive prefactor, and power to which the crossing probability in the lower bound is raised, consider the following intersection, with corresponding estimates,

\begin{align*}
      \textbf{P}^{\xi^{\mathrm{Sloped}}}_{\textbf{Z} \times [-n^{\prime}, 2n^{\prime}]} \big[    \mathcal{I}_0  \overset{h \geq ck}{ \longleftrightarrow } \widetilde{\mathcal{I}_0} \big]      \text{ } \geq \text{ } \textbf{P}^{\xi^{\mathrm{Sloped}}}_{\textbf{Z} \times [-n^{\prime}, 2n^{\prime}]} \big[ \big\{    \mathcal{I}_0  \overset{h \geq ck}{ \longleftrightarrow } \widetilde{\mathcal{I}_0} \big\} \text{ } \cap \text{ }   \big\{        \mathcal{I}_0            \overset{h \geq ck}{ \longleftrightarrow }     ( \rho \delta n , 0 \text{ } ] \times   \{ n \}         \big\}              \big]   \text{ }  \tag{$\delta$ \textit{LIM}} \\ \overset{ \delta \rightarrow - \infty}{\geq} \text{ }    \textbf{P}^{\xi^{\mathrm{Sloped}}}_{\textbf{Z} \times [-n^{\prime}, 2n^{\prime}]} \big[ \big\{    \mathcal{I}_0  \overset{h \geq ck}{ \longleftrightarrow } \widetilde{\mathcal{I}_0} \big\} \text{ } \cap \text{ }   \big\{        \mathcal{I}_0            \overset{h \geq ck}{ \longleftrightarrow }     ( - \infty , 0 \text{ } ] \times   \{ n \}         \big\}              \big]   \\ \overset{\mathrm{(FKG)}}{\geq} \text{ }            \textbf{P}^{\xi^{\mathrm{Sloped}}}_{\textbf{Z} \times [-n^{\prime}, 2n^{\prime}]} [      \mathcal{I}_0  \overset{h \geq ck}{ \longleftrightarrow } \widetilde{\mathcal{I}_0}               ] \text{ }    \textbf{P}^{\xi^{\mathrm{Sloped}}}_{\textbf{Z} \times [-n^{\prime}, 2n^{\prime}]} [               \mathcal{I}_0            \overset{h \geq ck}{ \longleftrightarrow }     ( - \infty , 0 \text{ } ] \times   \{ n \}        ]  \\ \text{ }    \geq \text{ } \big(          \textbf{P}^{\xi^{\mathrm{Sloped}}}_{\textbf{Z} \times [-n^{\prime}, 2n^{\prime}]} [               \mathcal{I}_0            \overset{h \geq ck}{ \longleftrightarrow }     ( - \infty , 0 \text{ } ] \times   \{ n \}        ]  \big)^2          \text{ , } 
\end{align*}

\noindent for the configuration over the square lattice obtained from the intersection $\big\{    \mathcal{I}_0  \overset{h \geq ck}{ \longleftrightarrow } \widetilde{\mathcal{I}_0} \big\} \text{ } \cap \text{ }   \big\{        \mathcal{I}_0            \overset{h \geq ck}{ \longleftrightarrow }     ( \rho \delta n , 0 \text{ } ] \times   \{ n \}         \big\}$ occurring, while the probability  $  \textbf{P}^{\xi^{\mathrm{Sloped}}}_{[-n^{\prime}, 2n^{\prime}]} \big[    \mathcal{I}_0  \overset{h \geq ck}{ \longleftrightarrow } \widetilde{\mathcal{I}_0} \big]$ can also be bound below with,

\begin{align*}
      \textbf{P}^{\xi^{\mathrm{Sloped}}}_{\textbf{Z} \times [-n^{\prime}, 2n^{\prime}]} \big[    \mathcal{I}_0  \overset{h \geq ck}{ \longleftrightarrow } \widetilde{\mathcal{I}_0} \big]  \text{ } \geq \text{ }   \textbf{P}^{\xi^{\mathrm{Sloped}}}_{\textbf{Z} \times [-n^{\prime}, 2n^{\prime}]} \big[ \big\{    \mathcal{I}_0  \overset{h \geq ck}{ \longleftrightarrow } \widetilde{\mathcal{I}_0} \big\} \text{ } \cap \text{ }   \big\{        \mathcal{I}_0            \overset{h \geq ck}{ \longleftrightarrow }     ( \rho \delta n , 0 \text{ } ] \times   \{ n \}         \big\}  \big] \\ \overset{\mathrm{(FKG)}}{\geq} \text{ }   \textbf{P}^{\xi^{\mathrm{Sloped}}}_{\textbf{Z} \times [-n^{\prime}, 2n^{\prime}]} \big[     \mathcal{I}_0  \overset{h \geq ck}{ \longleftrightarrow } \widetilde{\mathcal{I}_0} \big] \text{ }   \textbf{P}^{\xi^{\mathrm{Sloped}}}_{\textbf{Z} \times [-n^{\prime}, 2n^{\prime}]} \big[         \mathcal{I}_0            \overset{h \geq ck}{ \longleftrightarrow }     ( \rho \delta n , 0 \text{ } ] \times   \{ n \}      \big]  \text{ } \\ \geq  \text{ }       \big( \textbf{P}^{\xi^{\mathrm{Sloped}}}_{\textbf{Z} \times [-n^{\prime}, 2n^{\prime}]} \big[         \mathcal{I}_0            \overset{h \geq ck}{ \longleftrightarrow }     ( \rho \delta n , 0 \text{ } ] \times   \{ n \}      \big]  \big)^2      \text{ } \text{ , } 
\end{align*}

\noindent corresponding to the case that $\delta \not\longrightarrow -\infty$, as given in the weak infinite volume measure below $(\delta \text{ } \textit{LIM})$. In each case, for the final lower bound dependent upon $\big(          \textbf{P}^{\xi^{\mathrm{Sloped}}}_{\textbf{Z} \times [-n^{\prime}, 2n^{\prime}]} [               \mathcal{I}_0            \overset{h \geq ck}{ \longleftrightarrow }     ( - \infty , 0 \text{ } ] \times   \{ n \}        ]  \big)^2$, in the case that $\delta \longrightarrow -\infty$, or that $\delta \not\longrightarrow -\infty$,

\begin{align*}
\delta \not\longrightarrow -\infty \text{ } \Rightarrow \text{ } \big\{    \mathcal{I}_0            \overset{h \geq ck}{ \longleftrightarrow }     ( \rho \delta n , 0 \text{ } ] \times   \{ n \}       \big\} \text{ }  \subset \text{ } \big\{ \mathcal{I}_0 \overset{ h \geq ck}{\longleftrightarrow}        \widetilde{\mathcal{I}_0} \big\} \\ \delta \longrightarrow -\infty \text{ } \Rightarrow  \text{ } \big\{    \mathcal{I}_0            \overset{h \geq ck}{ \longleftrightarrow }     ( - \infty , 0 \text{ } ] \times   \{ n \}       \big\} \text{ }  \subset \text{ } \big\{ \mathcal{I}_0 \overset{ h \geq ck}{\longleftrightarrow}        \widetilde{\mathcal{I}_0} \big\}
\end{align*}

\noindent Finally, if $\{ \mathcal{I}_0 \longleftrightarrow \widetilde{\mathcal{I}_0} \}$ occurs, then there is no lower bound to show. Otherwise, for the remaining case,

\begin{align*}
    \textbf{P}^{\xi^{\mathrm{Sloped}}}_{\textbf{Z} \times [-n^{\prime} , 2n^{\prime}} [      \mathcal{I}_0 \overset{h \geq ck}{\longleftrightarrow} \widetilde{\mathcal{I}_0}         ]  \geq  \textbf{P}^{\xi^{\mathrm{Sloped}}}_{\textbf{Z} \times [-n^{\prime} , 2n^{\prime}} \bigg[   \big\{   \mathcal{I}_0 \overset{h \geq ck}{\longleftrightarrow}    ( - \infty , 0 \text{ } ]         \times  \{n \}     \big\} \cap \big\{    \widetilde{\mathcal{I}_0} \overset{h \geq ck}{\longleftrightarrow}     [ \delta n , + \infty )             \times  \{n \}      \big\}           \bigg] \\ \overset{\mathrm{(FKG)}}{\geq}\textbf{P}^{\xi^{\mathrm{Sloped}}}_{\textbf{Z} \times [-n^{\prime} , 2n^{\prime}} \big[    \mathcal{I}_0 \overset{h \geq ck}{\longleftrightarrow}    ( - \infty , 0 \text{ } ]         \times  \{n \}       \big]  \text{ }   \textbf{P}^{\xi^{\mathrm{Sloped}}}_{\textbf{Z} \times [-n^{\prime} , 2n^{\prime}]} \big[      \widetilde{\mathcal{I}_0} \overset{h \geq ck}{\longleftrightarrow}     [ \delta n , + \infty )             \times  \{n \}      \big] \\  \geq   \bigg[  \underset{\delta \in \textbf{R}}{\mathrm{inf}} \bigg\{     \textbf{P}^{\xi^{\mathrm{Sloped}}}_{\textbf{Z} \times [-n^{\prime} , 2n^{\prime}]} \big[    \mathcal{I}_0 \overset{h \geq ck}{\longleftrightarrow}    ( - \infty , 0 \text{ } ]         \times  \{n \}       \big]       ,  \textbf{P}^{\xi^{\mathrm{Sloped}}}_{\textbf{Z} \times [-n^{\prime} , 2n^{\prime}} \big[      \widetilde{\mathcal{I}_0} \overset{h \geq ck}{\longleftrightarrow}     [ \delta n , + \infty )             \times  \{n \}      \big]  \bigg\}  \bigg]^2      \text{ . }          
\end{align*}

\noindent From the infimum given in the lower bound above, if $  \textbf{P}^{\xi^{\mathrm{Sloped}}}_{\textbf{Z} \times [-n^{\prime} , 2n^{\prime}} \big[    \mathcal{I}_0 \overset{h \geq ck}{\longleftrightarrow}    ( - \infty , 0 \text{ } ]         \times  \{n \}       \big] \geq \textbf{P}^{\xi^{\mathrm{Sloped}}}_{\textbf{Z} \times [-n^{\prime} , 2n^{\prime}} \big[      \widetilde{\mathcal{I}_0} \overset{h \geq ck}{\longleftrightarrow}     [ \delta n , + \infty )             \times  \{n \}      \big] $,

\begin{align*}
  \textbf{P}^{\xi^{\mathrm{Sloped}}}_{\textbf{Z} \times [-n^{\prime} , 2n^{\prime}} \big[    \mathcal{I}_0 \overset{h \geq ck}{\longleftrightarrow}    ( - \infty , 0 \text{ } ]         \times  \{n \}       \big]      \textbf{P}^{\xi^{\mathrm{Sloped}}}_{\textbf{Z} \times [-n^{\prime} , 2n^{\prime}} \big[      \widetilde{\mathcal{I}_0} \overset{h \geq ck}{\longleftrightarrow}     [ \delta n , + \infty )             \times  \{n \}      \big]         \geq    \big(   \textbf{P}^{\xi^{\mathrm{Sloped}}}_{\textbf{Z} \times [-n^{\prime} , 2n^{\prime}} \big[      \widetilde{\mathcal{I}_0} \overset{h \geq ck}{\longleftrightarrow}     [ \delta n , + \infty )             \times  \{n \}      \big]   \big)^2 \text{ , } 
\end{align*}

\noindent while, if $  \textbf{P}^{\xi^{\mathrm{Sloped}}}_{\textbf{Z} \times [-n^{\prime} , 2n^{\prime}} \big[    \mathcal{I}_0 \overset{h \geq ck}{\longleftrightarrow}    ( - \infty , 0 \text{ } ]         \times  \{n \}       \big] \leq \textbf{P}^{\xi^{\mathrm{Sloped}}}_{\textbf{Z} \times [-n^{\prime} , 2n^{\prime}} \big[      \widetilde{\mathcal{I}_0} \overset{h \geq ck}{\longleftrightarrow}     [ \delta n , + \infty )             \times  \{n \}      \big] $, 

\begin{align*}
  \textbf{P}^{\xi^{\mathrm{Sloped}}}_{\textbf{Z} \times [-n^{\prime} , 2n^{\prime}} \big[    \mathcal{I}_0 \overset{h \geq ck}{\longleftrightarrow}    ( - \infty , 0 \text{ } ]         \times  \{n \}       \big]   \text{ }      \textbf{P}^{\xi^{\mathrm{Sloped}}}_{\textbf{Z} \times [-n^{\prime} , 2n^{\prime}} \big[      \widetilde{\mathcal{I}_0} \overset{h \geq ck}{\longleftrightarrow}     [ \delta n , + \infty )             \times  \{n \}      \big]      \geq    \big(   \textbf{P}^{\xi^{\mathrm{Sloped}}}_{\textbf{Z} \times [-n^{\prime} , 2n^{\prime}} \big[    \mathcal{I}_0 \overset{h \geq ck}{\longleftrightarrow}    ( - \infty , 0 \text{ } ]         \times  \{n \}       \big]    \big)^2 \text{ . }   
\end{align*}

\noindent With $\delta \longrightarrow - \infty$, from conditions for the first lower bound in which $\textbf{P}^{\xi^{\mathrm{Sloped}}}_{\textbf{Z} \times [-n^{\prime} , 2n^{\prime}} \big[    \mathcal{I}_0 \overset{h \geq ck}{\longleftrightarrow}    ( - \infty , 0 \text{ } ]         \times  \{n \}       \big] > \textbf{P}^{\xi^{\mathrm{Sloped}}}_{\textbf{Z} \times [-n^{\prime} , 2n^{\prime}} \big[      \widetilde{\mathcal{I}_0} \overset{h \geq ck}{\longleftrightarrow}     [ \delta n , + \infty )             \times  \{n \}      \big]$,

\begin{align*}
         \textbf{P}^{\xi^{\mathrm{Sloped}}}_{\textbf{Z} \times [-n^{\prime} , 2n^{\prime}} \big[    \mathcal{I}_0 \overset{h \geq ck}{\longleftrightarrow}    ( - \infty , 0 \text{ } ]         \times  \{n \}       \big]      \textbf{P}^{\xi^{\mathrm{Sloped}}}_{\textbf{Z} \times [-n^{\prime} , 2n^{\prime}} \big[      \widetilde{\mathcal{I}_0} \overset{h \geq ck}{\longleftrightarrow}     [ \delta n , + \infty )             \times  \{n \}      \big]         \geq    \big(   \textbf{P}^{\xi^{\mathrm{Sloped}}}_{\textbf{Z} \times [-n^{\prime} , 2n^{\prime}} \big[      \widetilde{\mathcal{I}_0} \overset{h \geq ck}{\longleftrightarrow}   \textbf{Z}             \times  \{n \}      \big]   \big)^2                    \text{ }   \text{ , } 
\end{align*}

\noindent From the inequality given in ($\mathscr{O}$\text{ }-\text{ }$\mathscr{N}$ \textit{bound}), we read the strictly positive prefactor $\mathscr{O}$, and the power to which the probability in the lower bound is raised, below,

\begin{align*}
          \text{ } \mathscr{O} \text{ }      \bigg[   \textbf{P}^{\xi^{\mathrm{Sloped}}}_{\textbf{Z} \times [-n^{\prime} , 2n^{\prime}} \big[      \widetilde{\mathcal{I}_0} \overset{h \geq ck}{\longleftrightarrow}   \textbf{Z}             \times  \{n \}      \big]   \bigg]^{\mathscr{N}}          \text{ }                    \text{ }   \text{ , } \text{ } 
\end{align*}

\noindent which in turn implies,

\begin{align*}
        \textbf{P}^{\xi^{\mathrm{Sloped}}}_{\textbf{Z} \times [-n^{\prime} , 2n^{\prime} ]} \big[                    \mathcal{I}_0   \underset{\textbf{Z} \times [0,n^{\prime}N]}{\overset{h \geq k}{\longleftrightarrow}}                \widetilde{\mathcal{I}_0}    \big]  \geq  \mathscr{O}    \bigg[   \textbf{P}^{\xi^{\mathrm{Sloped}}}_{[-m,m] \times [0,n^{\prime}N]} \big[   \mathcal{I}_0   \underset{\textbf{Z} \times [0,n^{\prime}N]}{\overset{h \geq k}{\longleftrightarrow}}                \widetilde{\mathcal{I}_0}           \big]                \bigg]^{\mathscr{N}} \\  \overset{m \longrightarrow + \infty}{\geq} \text{ }   \mathscr{O}     \bigg[   \textbf{P}^{\xi^{\mathrm{Sloped}}}_{\textbf{Z} \times [-n^{\prime} , 2n^{\prime}]} \big[      \widetilde{\mathcal{I}_0} \underset{\textbf{Z} \times [0,n^{\prime} N]}{\overset{h \geq ck}{\longleftrightarrow}}   \textbf{Z}             \times  \{n \}      \big]   \bigg]^{ \mathscr{N}} \text{ }  \\ \equiv          \bigg[  \sqrt[\mathscr{N}]{\mathscr{O}}   \textbf{P}^{\xi^{\mathrm{Sloped}}}_{\textbf{Z} \times [-n^{\prime} , 2n^{\prime}]} \big[      \widetilde{\mathcal{I}_0} \underset{\textbf{Z} \times [0,n^{\prime} N]}{\overset{h \geq ck}{\longleftrightarrow}}   \textbf{Z}             \times  \{n \}      \big]   \bigg]^{\mathscr{N}}             \text{ } \text{ , } 
\end{align*}

\noindent hence concluding the argument, for $\mathscr{O} < 1$ and strictly positive $\mathscr{N}$. \boxed{}

\section{From strip crossing probabilities to annulus crossing probabilities}

     \begin{figure}
\begin{align*}
\includegraphics[width=0.61\columnwidth]{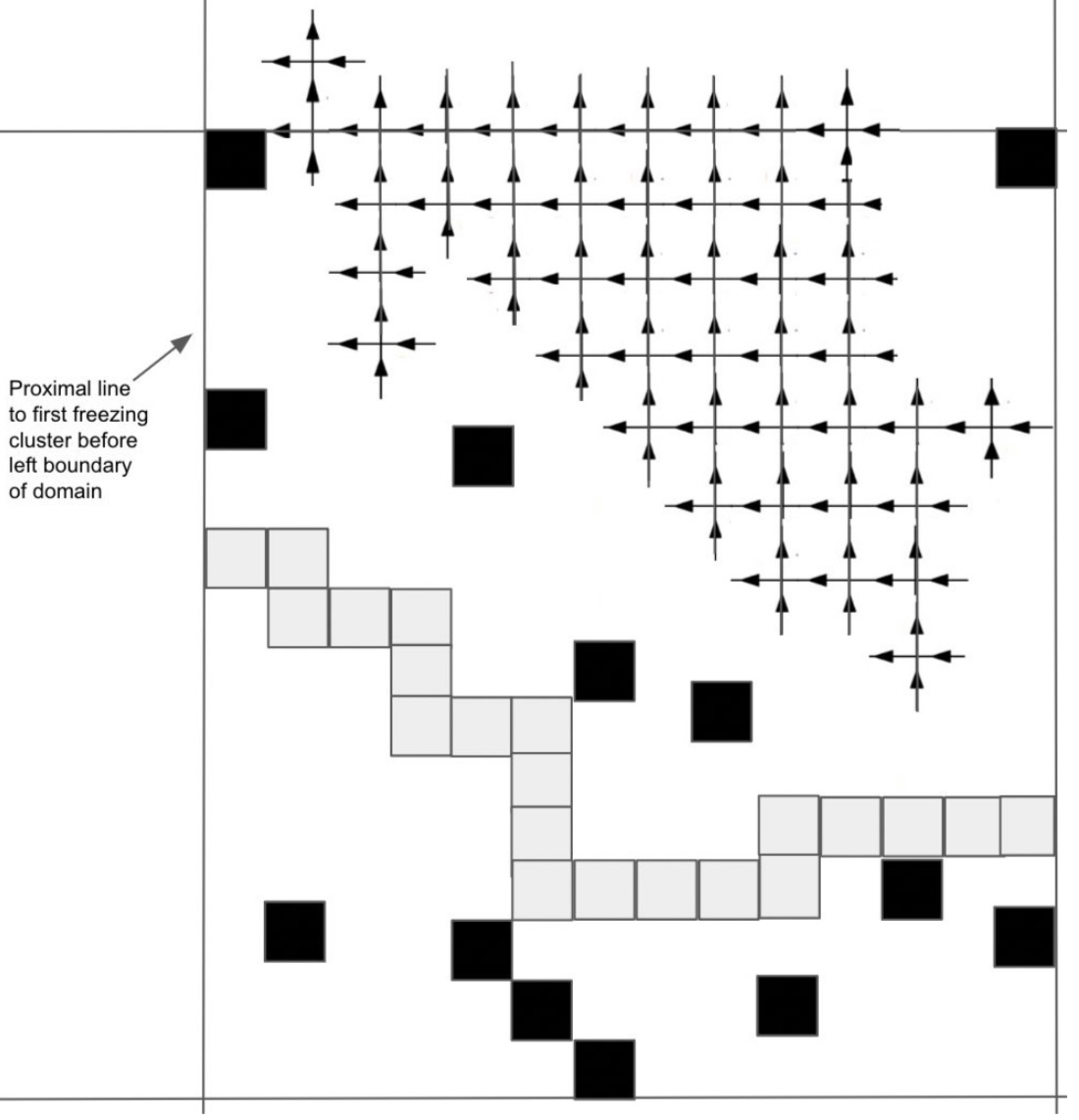}\\
\includegraphics[width=0.61\columnwidth]{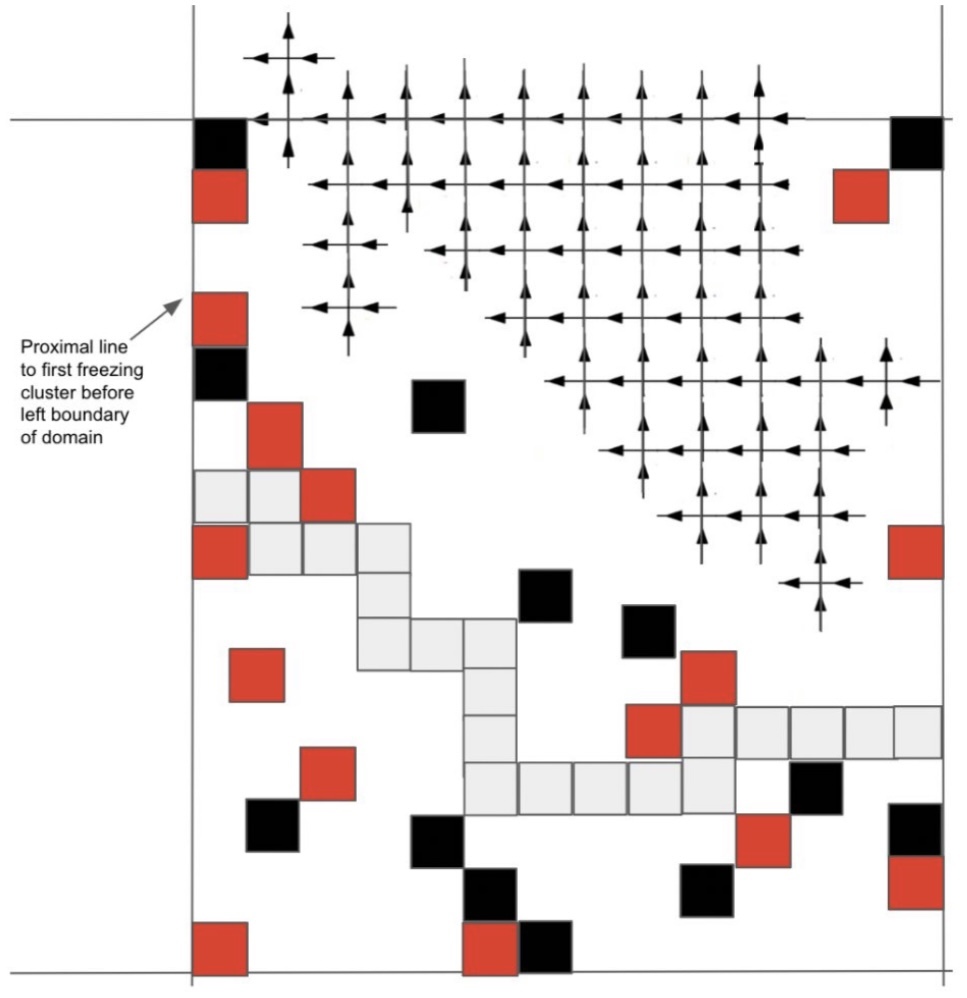}\\
\end{align*}
\caption{\textit{Configuration of the freezing cluster with the strip in finite volume}. As the constant $c$ required in the threshold $ck$ of the height function is decreased to $0$, additional faces, rather than those in the highlighted grey connected component, highlighted in black, can also contribute to the probability of obtaining a vertical crossing in the strip. As the crossing threshold is lowered further, additional faces bound within the finite volume strip, highlighted in red, will also contribute to a vertical crossing. Such geometric observations are }
\end{figure}

     \begin{figure}
\begin{align*}
\includegraphics[width=0.58\columnwidth]{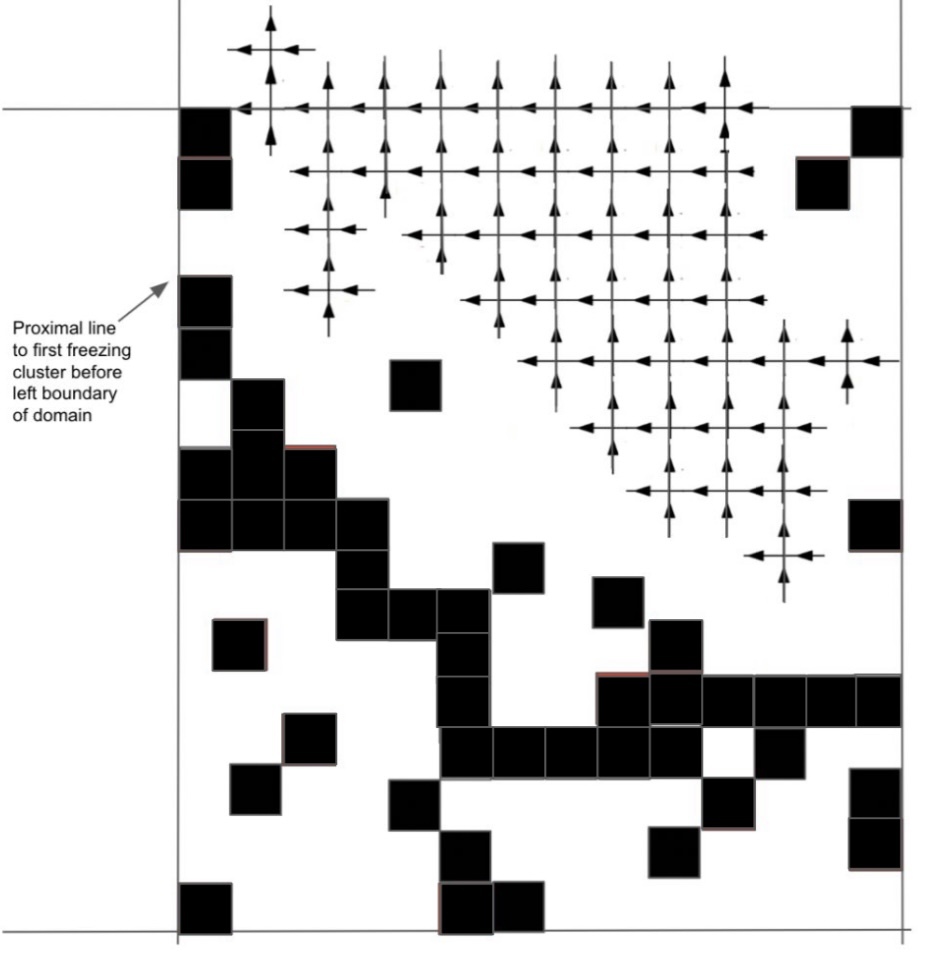}\\
\end{align*}
\caption{\textit{Configuration of the freezing cluster with the strip in finite volume}. The union of all connected components from images of the homomorphism from the faces highlighted in black, and in red, from the previous figure are together displayed as a unified black connected component of faces.}
\end{figure}

\noindent With results from the previous section, we proceed to analyze the impact of sloped boundary conditions with additional results  starting below.

\bigskip

\noindent The main result to be proved, in \textit{3.1}, is the following,

\bigskip

\noindent \textbf{Theorem} $6V$ \textit{1} (\textit{annulus-strip Russo-Seymour-Welsh}). Under $\xi^{\mathrm{Sloped}} \sim \textbf{B}\textbf{C}^{\mathrm{Sloped}}$, for $\delta > \delta^{\prime} >0$, $C^{\prime}>0$ and $n_R < 12$,

\begin{align*}
    \textbf{P}^{\xi^{\mathrm{Sloped}}}_{\textbf{Z} \times [0,n^{\prime} N]} \big[       \mathcal{O}_{h\geq ck} ( 6n , n_R n)            \big]  \geq c \bigg[    \textbf{P}^{\xi^{\mathrm{Sloped}}}_{\textbf{Z} \times [-n^{\prime} , 2 n^{\prime} ] }   \bigg[     [ 0 , \lfloor \delta^{\prime} n \rfloor ]           \times \{ 0 \}       \underset{ \textbf{Z} \times [0,n^{\prime} N]}{\overset{h\geq ck}{\longleftrightarrow}}          \textbf{Z}         \times \{ n \}    \bigg]    \bigg]^{C^{\prime}}   \text{ } \text{ , } 
\end{align*}

\noindent in which the probability of obtaining a crossing across the annulus, is lower bounded by the probability of obtaining a crossing between $    [ 0 , \lfloor \delta^{\prime} n \rfloor ]           \times \{ 0 \}$, and between $  \textbf{Z}         \times \{ n \}$, raised to some $C^{\prime}$.

\bigskip
 
\noindent In comparison to boundary conditions on the top and bottom of the strip from {\color{blue}[11]} over which the height function attains a value of at most $k$ and descends to $0,1$ boundary conditions for other faces incident to the boundary, arguments for upper bounding crossing probabilities across the strip, under sloped boundary conditions, with crossing probabilities across the annulus are still required for estimates of crossing probabilities under sloped boundary conditions on the cylinder. To this end, we adapt arguments for \textit{Theorem 3.1} in {\color{blue}[11]} with the result below.

\bigskip

\noindent \textit{Proof of Theorem 6V 1}. We make use of previously characterized domains in the strip, in which

\begin{align*}
     \textbf{P}^{\xi^{\mathrm{Sloped}}}_{[-m,m] \times [0,n^{\prime} N]} \big[     \mathcal{H}_k     \big]          \text{ } \text{ , } \text{ } 
\end{align*}

\noindent for any horizontal crossing event $\mathcal{H}_k$ can be bound below with the probability,

\begin{align*}
    \frac{1}{\mathscr{C}^{\prime}}  \text{ }   \textbf{P}^{(\xi^{\mathrm{Sloped}})_2}_{\Lambda} \big[     \mathcal{H}_{k-j}   \big]             \text{ } \text{ , } \text{ } 
\end{align*}

\noindent from the conditions of \textbf{Corollary} \textit{1.2}, with $\xi^{\mathrm{Sloped}} \geq  (\xi^{\mathrm{Sloped}})_2 \sim \textbf{B}\textbf{C}^{\mathrm{Sloped}}$ and $\Lambda \subset [-m,m] \times [0,n^{\prime}N]$. Within a strip domain, we introduce $\Gamma^{\prime}_L$, $\Gamma^{\prime}_R$, and $\Gamma^{\prime}_C$, with $\Gamma^{\prime}_L , \Gamma^{\prime}_R, \Gamma^{\prime}_C \subset \mathscr{D}_i$. In order to lower bound the conditional intersection of vertical crossing probabilities, we begin with,

\begin{align*}
 \textbf{P}^{\xi^{\mathrm{Sloped}}}_{[-m,m ]\times [0,n^{\prime}N]} \big[  \mathcal{V}(\Gamma^{\prime}_L)  \cap  \mathcal{V}(\Gamma^{\prime}_R)      \big|  \Gamma = \gamma   \big] \text{ }   \text{ , } \text{ } \\
\end{align*}

\noindent where left and right boundaries of the domain appear in the conditioning for the event defined above, for a realization of the domain $\gamma$ from $\Gamma$. Our strategy is to randomly reveal the faces within $\Gamma^{\prime}_L$, $\Gamma^{\prime}_R$, and $\Gamma^{\prime}_C$, in which the crossing will occur with high probability given a lower bound on horizontal crossing events across the domain. We observe that the above crossing probability strictly dominates the product,

\begin{align*}
     \textbf{P}^{\xi^{\mathrm{Sloped}}}_{\Gamma^{\prime}_L} \big[     \mathcal{V}(\Gamma^{\prime}_L)   \big]    \text{ }   \textbf{P}^{\xi^{\mathrm{Sloped}}}_{\Gamma^{\prime}_R} \big[     \mathcal{V}(\Gamma^{\prime}_R)      \big]    \text{ } \text{ , } \text{ } \\
\end{align*}

\noindent of crossing probabilities across left and right subdomains of $\mathscr{D}_i$, leading to the estimate obtained from rearrangement,

\begin{align*}
\mathscr{P} \text{ }  \equiv \text{ }       \textbf{P}^{\xi^{\mathrm{Sloped}}}_{[-m,m] \times [0,n^{\prime}N]} \big[   \mathcal{H}(\Gamma^{\prime}_C)         \text{ }  \big| \text{ }      \mathcal{V}(\Gamma^{\prime}_L) \cap \mathcal{V}(\Gamma^{\prime}_R)    \big]  \text{ }  \geq \text{ }    \textbf{P}^{\xi^{\mathrm{Sloped}}}_{[-m,m] \times [0,n^{\prime}N]} \big[     \mathcal{H}(\Gamma^{\prime}_C)  \cap     \mathcal{V}(\Gamma^{\prime}_L)  \cap   \mathcal{V}(\Gamma^{\prime}_R)      \big]         \text{ } \text{ , } \text{ }  \tag{\textit{*}}\\
\end{align*}

\noindent which, finally, itself strictly dominates the product of probabilities below,

\begin{align*}
    \textbf{P}^{\xi^{\mathrm{Sloped}}}_{[-m,m] \times [0,n^{\prime} N]} \big[       \mathcal{H}(\Gamma^{\prime}_C)      \big] \text{ }   \textbf{P}^{\xi^{\mathrm{Sloped}}}_{[-m,m] \times [0,n^{\prime} N]} \big[  \mathcal{V}(\Gamma^{\prime}_L)  \big] \text{ } \textbf{P}^{\xi^{\mathrm{Sloped}}}_{[-m,m] \times [0,n^{\prime} N]} \big[  \mathcal{V}(\Gamma^{\prime}_R) \big]  \text{ } \text{ , } \text{ } \tag{\textit{**}}\\
\end{align*}

\noindent by (FKG). To complete the lower bound, we make use of the fact that the conditional measured can be expressed as the convex combination of crossing probabilities, in which each probability measure is supported over $\mathscr{D}_i$,

\begin{align*}
    \mathscr{C}_1 \textbf{P}^{\xi^{\mathrm{Sloped}}}_{\mathscr{D}_i}  \big[           \mathcal{H}(\Gamma^{\prime}_C)  \big]  +  \mathscr{C}_2 \textbf{P}^{\xi^{\mathrm{Sloped}}}_{\mathscr{D}_i}  \big[        (\mathcal{H}(\Gamma^{\prime}_C))^c         \big]   \text{ } \text{ , } \text{ }  \tag{CONV}
\end{align*}

\noindent for suitably chosen constants $\mathscr{C}_1$ and $\mathscr{C}_2$ with $\mathscr{C}_1+\mathscr{C}_2 =1$. To lighten notation for the strip connectivity event below, we suppress conditions relating to the domain over which connectivity occurs, in which,

\begin{align*}
 \big\{ [0,\lfloor{\delta n}\rfloor{}  ]  \times \{ 0 \} \underset{{\textbf{Z} \times[0,n^{\prime}N]}}{\overset{h\geq ck}{\longleftrightarrow}} [-m,m] \times \{ n \}      \big\}  \equiv      \big\{[0,\lfloor{\delta n}\rfloor{}  ]  \times  \{ 0 \}    {\overset{h\geq ck}{\longleftrightarrow}}   [-m,m] \times \{ n \}     \big\}         \text{ } \text{ . } \text{ } 
\end{align*}

\noindent These rearrangements allow us to conclude, by (CBC) that there exists a suitable domain across which the horizontal crossing probability can be lower bounded, with the estimate beginning from $\textit{(*)}$,

\begin{align*}
        \mathscr{P}   \overset{\textit{(**)}}{\geq}   \textbf{P}^{\xi^{\mathrm{Sloped}}}_{[-m,m] \times [0,n^{\prime} N]} \big[    \mathcal{H}(\Gamma^{\prime}_C)     \big]   \tag{\textit{***}}  \end{align*}
        
        \noindent which can be further lower bounded with the probability,
        
        \begin{align*}
        \mathcal{C}_{L,R}     \big( \textbf{P}^{\xi^{\mathrm{Sloped}}}_{[-m,m] \times [0,n^{\prime} N]} \big[  [0,\lfloor{\delta n}\rfloor{}  ] \text{ } \times \text{ } \{ 0 \}\text{ }     {\overset{h\geq ck}{\longleftrightarrow}} \text{ }   [-m,m] \text{ } \times \{ n \}    \big]      \big)^{\mathscr{N}}       \text{ } \text{ , } \text{ } \tag{\textit{****}}
\end{align*}

\noindent where, first in $\textit{(***)}$, there exists strictly positive constants for which,

\begin{align*}
  \textbf{P}^{\xi^{\mathrm{Sloped}}}_{[-m,m] \times [0,n^{\prime} N]} \big[ \mathcal{V}(\Gamma^{\prime}_L)  \big]  \geq \mathcal{C}_L  \text{ } \text{ , } 
\end{align*}

\noindent and also for which,

\begin{align*}
   \textbf{P}^{\xi^{\mathrm{Sloped}}}_{[-m,m] \times [0,n^{\prime} N]} \big[  \mathcal{V}(\Gamma^{\prime}_R) \big]  \geq \mathcal{C}_R \text{ } \text{ , } 
\end{align*}

\noindent each of which respectively exist almost surely by finite energy. Second, in \textit{(****)},

\begin{align*}
    \textbf{P}^{\xi^{\mathrm{Sloped}}}_{[-m,m] \times [0,n^{\prime}N]} \big[    \mathcal{H}(\Gamma^{\prime}_C)        ]   \overset{\mathrm{(FKG)}}{\geq} \prod_{\mathscr{N} \text{ }  \mathrm{ Subdomains \text{ } indexed\text{ }  in } \text{ } i}  \textbf{P}^{\xi^{\mathrm{Sloped}}}_{[-m,m] \times [0,n^{\prime}N]} \big[   \mathcal{H}(\Gamma^{\prime}_i)      ]  \geq   \bigg[     \textbf{P}^{\xi^{\mathrm{Sloped}}}_{[-m,m] \times [0,n^{\prime}N]}  \big[       \mathcal{H}(\Gamma^{\prime}_i )        \big]  \bigg]^{\mathscr{N}} \text{ } \\  \geq  \bigg[    \textbf{P}^{\xi^{\mathrm{Sloped}}}_{[-m,m] \times [0,n^{\prime}N]}    \big[     [0,\lfloor{\delta n}\rfloor{}  ]  \times  \{ 0 \}    {\overset{h\geq ck}{\longleftrightarrow}}   [-m,m]  \times \{ n \}     \big]     \bigg]^{\mathscr{N}}     \text{ } \text{ , } \text{ } 
\end{align*}

\noindent where $\mathscr{N}$ denotes the number of subdomains indexed in $i$ over which the produce of horizontal crossing events is taken, hence yielding the desired estimate. The final constant lower bound the product of both of the probabilities above satisfies, 

\begin{align*}
  \mathcal{C}_L  \mathcal{C}_R \geq  \mathcal{C}_{L,R}  \text{ } \text{ , } \text{ } 
\end{align*}

\noindent where,

\begin{align*}
  \mathcal{C}_{L,R} \equiv  \underset{\mathrm{sides }\text{ }  L , R}{\mathrm{inf}} \big\{ \mathcal{C}_L  , \mathcal{C}_R  \big\}    \text{ } \text{ , } \text{ } 
\end{align*}

\noindent as a result implying, with $m \rightarrow +\infty$, that,

\begin{align*}
     \mathscr{P}        \geq  \mathcal{C}_{L,R}   \bigg[   \textbf{P}^{\xi^{\mathrm{Sloped}}}_{\textbf{Z} \times [0,n^{\prime} N]} \big[  [0,\lfloor{\delta n}\rfloor{}  ]  \times  \{ 0 \}    {\overset{h\geq ck}{\longleftrightarrow}}  [-m,m] \times \{ n \}    \big]  \bigg]^{\mathscr{N}} \text{ } \text{ . } 
\end{align*}

\noindent After obtaining this estimate, observe, from $\mathrm{(CONV)}$,

\begin{align*}
  \mathscr{P} \underset{\mathscr{D}_i \subset [-m,m] \times [0,n^{\prime} N]}{\overset{\mathrm{(SMP)}}{\equiv}} \text{ }    \textbf{P}^{\xi^{\mathrm{Sloped}}}_{\mathscr{D}_i} \big[   \mathcal{H}(\Gamma^{\prime}_C)         \text{ }  \big| \text{ }      \mathcal{V}(\Gamma^{\prime}_L) \text{ } \cap \text{ } \mathcal{V}(\Gamma^{\prime}_R) ,  h = \xi \text{ on } \partial \mathscr{D}_i  \cup      \mathscr{D}_i^c    \big]        \\ \equiv  \textbf{P}^{\xi^{(\mathrm{Sloped})_{\mathrm{CONV}}}}_{\mathscr{D}_i} \big[   \mathcal{H}(\Gamma^{\prime}_C)         \text{ }  \big| \text{ }      \mathcal{V}(\Gamma^{\prime}_L) \cap  \mathcal{V}(\Gamma^{\prime}_R)   \big] \\ \overset{\textit{(*)}}{\geq}       \textbf{P}^{\xi^{(\mathrm{Sloped})_{\mathrm{CONV}}}}_{\mathscr{D}_i} \big[      \mathcal{H}(\Gamma^{\prime}_C)  \cap      \mathcal{V}(\Gamma^{\prime}_L) \cap     \mathcal{V}(\Gamma^{\prime}_R)    \big]                                               \\
  \overset{\textit{(**)}}{\geq} \textbf{P}^{\xi^{(\mathrm{Sloped})_{\mathrm{CONV}}}}_{\mathscr{D}_i} \big[    \mathcal{H}(\Gamma^{\prime}_C)       \big]   \textbf{P}^{\xi^{(\mathrm{Sloped})_{\mathrm{CONV}}}}_{[-m,m] \times [0,n^{\prime} N]} \big[  \mathcal{V}(\Gamma^{\prime}_L)  \big]  \textbf{P}^{\xi^{(\mathrm{Sloped})_{\mathrm{CONV}}}}_{[-m,m] \times [0,n^{\prime} N]} \big[  \mathcal{V}(\Gamma^{\prime}_R) \big]  \\ \overset{\mathrm{(FKG)}}{\geq}  \underset{k \in \mathcal{I}_3}{\underset{j \in \mathcal{I}_2}{\prod_{i_1 \in \mathcal{I}_1} }}   \textbf{P}^{\xi^{(\mathrm{Sloped})_{\mathrm{CONV}}}}_{\mathscr{D}_i} \big[      (\mathcal{H}(\Gamma^{\prime}_C))_{i_1}     \big]   \textbf{P}^{\xi^{(\mathrm{Sloped})_{\mathrm{CONV}}}}_{\mathscr{D}_i} \big[ (\mathcal{V}(\Gamma^{\prime}_L))_j  \big] \text{ } \textbf{P}^{\xi^{(\mathrm{Sloped})_{\mathrm{CONV}}}}_{\mathscr{D}_i} \big[  (\mathcal{V}(\Gamma^{\prime}_R))_k  \big]                      \tag{\textit{i-j-k FKG}} \\ \geq           \underset{k \in \mathcal{I}_3}{\underset{j \in \mathcal{I}_2}{\prod_{i_1 \in \mathcal{I}_1} }} \text{ }   \bigg[      \mathscr{C}_1 \textbf{P}^{\xi^{(\mathrm{Sloped})_{\mathrm{CONV}}}}_{\mathscr{D}_i} \big[         (  \mathcal{H}(\Gamma^{\prime}_C) )_{i_1}                                 \big]  +  \mathscr{C}_2    \textbf{P}^{\xi^{(\mathrm{Sloped})_{\mathrm{CONV}}}}_{\mathscr{D}_i} \big[                                               ( (\mathcal{H}(\Gamma^{\prime}_C))^c )_{i_1}     \big]   \bigg]   \textbf{P}^{\xi^{(\mathrm{Sloped})_{\mathrm{CONV}}}}_{\mathscr{D}_i} \big[  (\mathcal{V}(\Gamma^{\prime}_L))_j  \big] \text{ }   \\ \times         \textbf{P}^{\xi^{(\mathrm{Sloped})_{\mathrm{CONV}}}}_{\mathscr{D}_i} \big[ (\mathcal{V}(\Gamma^{\prime}_R))_k                     \big] \text{ }  \text{ , }  \end{align*}

 \noindent in which from $\textit{(i-j-k-FKG)}$, product indices $i_1$, $j$ and $k$ are respectively introduced for crossings $\mathcal{H}(\Gamma^{\prime}_C)$, $\mathcal{V}(\Gamma^{\prime}_L)$ and $\mathcal{V}(\Gamma^{\prime}_R)$. Proceeding, bounding each vertical crossing probability individually with a strictly positive constant yields, 
  
  \begin{align*}
\cdots \overset{\mathrm{(FKG)}}{\geq}  \underset{k \in \mathcal{I}_3}{\underset{j \in \mathcal{I}_2}{\prod_{i_1 \in \mathcal{I}_1} }}      \bigg[      \mathscr{C}_1  \textbf{P}^{\xi^{(\mathrm{Sloped})_{\mathrm{CONV}}}}_{\mathscr{D}_i} \big[      (   \mathcal{H}(\Gamma^{\prime}_C) )_{i_1}                              \big] + \mathscr{C}_2  \textbf{P}^{\xi^{(\mathrm{Sloped})_{\mathrm{CONV}}}}_{\mathscr{D}_i} \big[                                        ((\mathcal{H}(\Gamma^{\prime}_C))^c  )_{i_1}    \big]        \bigg]     \underset{\geq\text{ } \mathcal{C} (\mathcal{V}(\Gamma^{\prime}_L) )     }{\underbrace{\textbf{P}^{\xi^{(\mathrm{Sloped})_{\mathrm{CONV}}}}_{\mathscr{D}_i} \big[  (\mathcal{V}(\Gamma^{\prime}_L))_j  \big]  }}   \\  \times \underset{\geq \mathcal{C} (\mathcal{V}(\Gamma^{\prime}_R) )         }{\underbrace{ \textbf{P}^{\xi^{(\mathrm{Sloped})_{\mathrm{CONV}}}}_{\mathscr{D}_i} \big[  (\mathcal{V}(\Gamma^{\prime}_R))_k        \big]}} \end{align*}

\noindent which result as multiplicative constants to the following product over $i_1$,

\begin{align*}
\cdots \geq  \mathcal{C}(\mathcal{V}(\Gamma^{\prime}_L)) \mathcal{C}(\mathcal{V}(\Gamma^{\prime}_R))  \prod_{i_1 \in \mathcal{I}_1}  \bigg[       \mathscr{C}_1 \textbf{P}^{\xi^{(\mathrm{Sloped})_{\mathrm{CONV}}}}_{\mathscr{D}_i} \big[      (   \mathcal{H}(\Gamma^{\prime}_C) )_{i_1}                               \big] + \mathscr{C}_2 \textbf{P}^{\xi^{(\mathrm{Sloped})_{\mathrm{CONV}}}}_{\mathscr{D}_i} \big[                                                ((\mathcal{H}(\Gamma^{\prime}_C))^c  )_{i_1}    \big]      \bigg]    \text{ , } 
       \end{align*}
      
      \noindent providing additional constants for the lower bound, given below as,
      
      \begin{align*} \cdots  \geq   \mathcal{C}(\mathcal{V}(\Gamma^{\prime}_L))  \mathcal{C}(\mathcal{V}(\Gamma^{\prime}_R)) \mathscr{C}_{1,2}  \prod_{i_1 \in \mathcal{I}_1 }     \bigg[   \textbf{P}^{\xi^{(\mathrm{Sloped})_{\mathrm{CONV}}}}_{\mathscr{D}_i} \big[        (   \mathcal{H}(\Gamma^{\prime}_C) )_{i_1}                              \big]  + \textbf{P}^{\xi^{(\mathrm{Sloped})_{\mathrm{CONV}}}}_{\mathscr{D}_i} \big[                                               ((\mathcal{H}(\Gamma^{\prime}_C))^c  )_{i_1}     \big]   \bigg]   \text{ } \tag{\textit{First constant}} \text{ . } 
      \end{align*}

   \noindent The product of probabilities above admits a lower bound,
   
   \begin{align*}
     \prod_{i_1 \in \mathcal{I}_1 }   \bigg[     \textbf{P}^{\xi^{(\mathrm{Sloped})_{\mathrm{CONV}}}}_{\mathscr{D}_i} \big[      (   \mathcal{H}(\Gamma^{\prime}_C) )_{i_1}                            \big]  +  \textbf{P}^{\xi^{(\mathrm{Sloped})_{\mathrm{CONV}}}}_{\mathscr{D}_i} \big[                                                ((\mathcal{H}(\Gamma^{\prime}_C))^c  )_{i_1}   \big]   \bigg] \geq       \prod_{i_1 \in \mathcal{I}_1 }  \textbf{P}^{\xi^{\mathrm{Sloped}}_{\mathrm{CONV}}}_{\mathscr{D}_i} \big[   (   \mathcal{H}(\Gamma^{\prime}_C) )_{i_1}      \big] \text{ }   \text{ , } 
   \end{align*}
      
      \noindent in addition to the pushforward of the collection of events $\{ \text{ } \mathcal{H}(\Gamma^{\prime}_C)_{i_1}\text{ } \}$, under $\textbf{P}^{\xi^{\mathrm{Sloped}}} [ \text{ } \cdot \text{ } ]$, satisfying,

      \begin{align*}
      \textbf{P}^{\xi^{\mathrm{Sloped}}_{\mathrm{CONV}}}_{\mathscr{D}_i} \big[  (   \mathcal{H}(\Gamma^{\prime}_C) )_{i_1}    \big]  \equiv    \textbf{P}^{\xi^{\mathrm{S}}_{\mathrm{C}}}_{\mathscr{D}_i} \big[    (   \mathcal{H}(\Gamma^{\prime}_C) )_{i_1}   \big]   \overset{\mathrm{(SMP)}}{\equiv}   \textbf{P}^{\xi^{\mathrm{S}}_{\mathrm{C}}|_{\partial (\mathscr{D}_0)_{i_1}}}_{(\mathscr{D}_0)_{i_1}} \big[      (   \mathcal{H}(\Gamma^{\prime}_C) )_{i_1}     \big| \text{ }  h =  \xi^{\mathrm{Sloped}}_{\mathrm{CONV}} \text{ on }      ((\mathscr{D}_0)_{i_1})^c  \cup  \partial ( (\mathscr{D}_0)_{i_1})    \big] \text{ }    \\ \overset{\xi^{\mathrm{S}}_{\mathrm{C}}|_{\partial (\mathscr{D}_0)_{i_1}}\text{ }  \leq \text{ } (\xi^{\mathrm{S}}_{\mathrm{C}})_2 \text{ } \sim \textbf{B}\textbf{C}^{\mathrm{Sloped}}}{\equiv} \text{ }   \textbf{P}^{(\xi^{\mathrm{S}}_{\mathrm{C}})_2}_{(\mathscr{D}_0)_{i_1}} \big[ (   \mathcal{H}(\Gamma^{\prime}_C) )_{i_1}              \big] \text{ } 
      \end{align*}

      \noindent finally, from which there exists another pair of sloped boundary conditions, with $(\xi^{\mathrm{S}}_{\mathrm{C}})_2 \leq (\xi^{\mathrm{Sloped}}_{\mathrm{CONV}})_3 \sim \textbf{B}\textbf{C}^{\mathrm{Sloped}}$, for which, 
      
      \begin{align*}
      \text{ }    \textbf{P}^{(\xi^{\mathrm{S}}_{\mathrm{C}})_2}_{(\mathscr{D}_0)_{i_1}}  \big[   (   \mathcal{H}(\Gamma^{\prime}_C) )_{i_1}                  \big]    \leq     \textbf{P}^{(\xi^{\mathrm{Sloped}}_{\mathrm{CONV}})_3}_{(\mathscr{D}_0)_{i_1}}     \big[       (   \mathcal{H}(\Gamma^{\prime}_C) )_{i_1}             \big]      \text{ } \text{ , } 
      \end{align*}
      
      \noindent for $(\mathscr{D}_0)_{i_1} \subset \Gamma^{\prime}_C$. Altogether, 
      
      \begin{align*}
\mathcal{C}(\mathcal{V}(\Gamma^{\prime}_L)) \text{ } \mathcal{C}(\mathcal{V}(\Gamma^{\prime}_R)) \text{ }  \mathscr{C}_{1,2} \text{ }     \prod_{i_1 \in \mathcal{I}_1} \text{ }    \textbf{P}^{(\xi^{\mathrm{Sloped}}_{\mathrm{CONV}})_3}_{(\mathscr{D}_0)_{i_1}}     \big[       (   \mathcal{H}(\Gamma^{\prime}_C) )_{i_1}              \big]        \geq  \mathcal{C}(\mathcal{V}(\Gamma^{\prime}_L)) \text{ } \mathcal{C}(\mathcal{V}(\Gamma^{\prime}_R))  \mathscr{C}_{1,2} \text{ } \\ \times  \bigg[  \textbf{P}^{\xi^{(\mathrm{Sloped})_{\mathrm{CONV}}}}_{[-m,m] \times [0,n^{\prime} N]}  \big[  [0,\lfloor{\delta_i n}\rfloor{}  ]  \times  \{ 0 \}  {\overset{h\geq ck}{\longleftrightarrow}}   [-m,m] \times  \{ n \}     \big] \bigg]^{\mathscr{N}}  \\             \geq   \mathscr{C}^{\prime}_{\mathrm{min}}                  \bigg[   \textbf{P}^{\xi^{(\mathrm{Sloped})_{\mathrm{CONV}}}}_{[-m,m] \times [0,n^{\prime} N]}  \big[ [0,\lfloor{\delta_i n}\rfloor{}  ]  \times  \{ 0 \}  {\overset{h\geq ck}{\longleftrightarrow}}  [-m,m] \times  \{ n \}  \big] \bigg]^{\mathscr{N}} \text{ , } \tag{\textit{Second constant}}
\end{align*}    

\noindent where the strictly positive constant $\mathscr{C}_{1,2}$, introduced in $\textit{(First constant)}$, satisfies,

\begin{align*}
    \mathscr{C}_{1,2}  \equiv   \mathrm{inf} \big\{  \mathscr{C}_1 , \mathscr{C}_2 \big\}   \text{ , } 
\end{align*}

\noindent the strictly positive constant $\mathscr{C}^{\prime}_{\mathrm{min}}$, introduced in $\textit{(Second constant)}$, satisfies,

\begin{align*}
    \mathscr{C}^{\prime}_{\mathrm{min}}  \equiv      \underset{1,2}{\underset{L,R}{\mathrm{inf}}} \big\{  \mathcal{C}(\mathcal{V}(\Gamma^{\prime}_L))  \mathcal{C}(\mathcal{V}(\Gamma^{\prime}_R))     ,  \mathscr{C}_{1,2}    \big\}   \Leftrightarrow    \underset{1,2}{\underset{L,R}{\mathrm{inf}}}   \big\{   \mathcal{C}(\mathcal{V}(\Gamma^{\prime}_L)) \mathcal{C}(\mathcal{V}(\Gamma^{\prime}_R))     ,   \underset{1,2}{\mathrm{inf}} \big\{ \mathscr{C}_1 , \mathscr{C}_2 \big\}     \big\}              \text{ . } 
\end{align*}

\noindent Finally, the vertical and horizontal crossing events introduced in \textit{(i-j-k FKG)}, $(\mathcal{H}(\Gamma^{\prime}_C))_i$, $(\mathcal{V}(\Gamma^{\prime}_L))_j$, and $(\mathcal{V}(\Gamma^{\prime}_R))_k$, respectively denote partitions of the crossing events $\mathcal{H}(\Gamma^{\prime}_C)$, $\mathcal{V}(\Gamma^{\prime}_L)$, and $\mathcal{V}(\Gamma^{\prime}_R)$ with $i \in \mathcal{I}_1,j \in \mathcal{I}_2$, and $k \in \mathcal{I}_3$, for index sets $\mathcal{I}_1, \mathcal{I}_2$, and $\mathcal{I}_3$. To obtain $\textit{(First constant)}$, we express the crossing event as a product over $i_1$, after having expressed the conditional crossing event from lower bounds provided in \textit{(*)} and \textit{(**)}. 

\bigskip

\noindent Next, for horizontal crossings across $\Gamma^{\prime}_C$, to make use of the estimate,

\begin{align*}
    \textbf{P}^{\xi^{\mathrm{Sloped}}}_{\mathscr{D}_0} \big[     \mathcal{H}(\Gamma^{\prime}_C)  \big]  \geq      \mathscr{C}^{\prime}_{\mathrm{min}}   \bigg[   \textbf{P}^{\xi^{(\mathrm{Sloped})_{\mathrm{CONV}}}}_{[-m,m] \times [0,n^{\prime}N]}  \big[  [0,\lfloor{\delta_i n}\rfloor{}  ]  \times \{ 0 \}  {\overset{h\geq ck}{\longleftrightarrow}}  [-m,m] \times  \{ n \}    \big]  \bigg]^{\mathscr{N}} \text{ , } \text{ } 
\end{align*}

\noindent consider $\frac{\pi}{2}$, $\pi$, and $\frac{3 \pi}{2}$ rotations of the centermost region $\Gamma^{\prime}_C$, respectively denoted with $(\Gamma^{\prime}_C)_{\frac{\pi}{2}}$, $(\Gamma^{\prime}_C)_{\pi}$, and $(\Gamma^{\prime}_C)_{\frac{3 \pi}{2}}$. Across the annulus, further rearrangements yield the estimate,

\begin{align*}
      \textbf{P}^{\xi^{\mathrm{Sloped}}}_{\Lambda} \big[      \mathcal{H}( (\Gamma^{\prime}_C)_{\frac{\pi}{2}} )     \big]  \text{ } \overset{\mathrm{(SMP)}}{\equiv} \text{ }           \text{ }     \textbf{P}^{\xi^{\mathrm{Sloped}}|{\partial (\mathscr{D}_0)_{i_1}}}_{(\mathscr{D}_0)_{i_1}} \big[    \mathcal{H}( (\Gamma^{\prime}_C)_{\frac{\pi}{2}} )      \big| \text{ } h =   \text{ }  \xi^{\mathrm{Sloped}} \text{ } \text{ on } \text{ } (\mathscr{D}_0)_{i_1}^c \text{ } \cup \text{ }  \partial (\mathscr{D}_0)_{i_1}  \big] \text{ }        \\ \text{ } \equiv    \textbf{P}^{(\xi^{\mathrm{Sloped}})_3}_{(\mathscr{D}_0)_{i_1}} \big[    \mathcal{H}( (\Gamma^{\prime}_C)_{\frac{\pi}{2}} )              \big] \text{ } \\ \text{ } \geq      \textbf{P}^{(\xi^{\mathrm{Sloped}})_3}_{[-m,m] \times [0,n^{\prime}N]} \big[         \mathcal{H}( (\Gamma^{\prime}_C)_{\frac{\pi}{2}} )     \big]         \text{ } \\ \text{ } \overset{m \rightarrow +\infty}{\geq}   \textbf{P}^{(\xi^{\mathrm{Sloped}})_3}_{\textbf{Z} \times [0,n^{\prime}N]} \big[        \mathcal{H}( (\Gamma^{\prime}_C)_{\frac{\pi}{2}} )     \big]       \text{ }  \text{ , } \text{ } 
\end{align*}

\noindent corresponding to the probability of crossing $(\Gamma^{\prime}_C)_{\frac{\pi}{2}}$ in the longer direction, while along similar lines,

\begin{align*}
      \textbf{P}^{\xi^{\mathrm{Sloped}}}_{\Lambda} \big[       \mathcal{H}( (\Gamma^{\prime}_C)_{\pi} )         \big]   \overset{\mathrm{(SMP)}}{\equiv}    \textbf{P}^{\xi^{\mathrm{Sloped}}}_{(\mathscr{D}_0)_{i_1}} \big[   \mathcal{H}( (\Gamma^{\prime}_C)_{\pi} )             \big| \text{ } h =    \text{ } \xi^{\mathrm{Sloped}} \text{ on }   (\mathscr{D}_0)_{i_1}^c  \cup   \partial (\mathscr{D}_0)_{i_1}  \big] \text{ } \\  \text{ }  \equiv   \textbf{P}^{((\xi^{\mathrm{Sloped}})_3)^{\prime}}_{(\mathscr{D}_0)_{i_1}} \big[    \mathcal{H}( (\Gamma^{\prime}_C)_{\pi} )      \big] \\ \text{ }     \geq   \textbf{P}^{((\xi^{\mathrm{Sloped}})_3)^{\prime}}_{[-m,m] \times [0,n^{\prime}N]} \big[      \mathcal{H}( (\Gamma^{\prime}_C)_{\pi} )      \big]   \text{ }    \\ \overset{m \rightarrow +\infty}{\geq} \textbf{P}^{((\xi^{\mathrm{Sloped}})_3)^{\prime}}_{\textbf{Z} \times [0,n^{\prime}N]} \big[       \mathcal{H}( (\Gamma^{\prime}_C)_{\pi} )        \big]    \text{ }  \text{ , } \text{ } 
\end{align*}

\noindent yield an estimate corresponding to the probability of crossing $(\Gamma^{\prime}_C)_{\pi}$ in the longer direction, and finally,

\begin{align*}
      \textbf{P}^{\xi^{\mathrm{Sloped}}}_{\Lambda} \big[        \mathcal{H}( (\Gamma^{\prime}_C)_{\frac{3\pi}{2}} )       \big]    \overset{\mathrm{(SMP)}}{\equiv}  \textbf{P}^{\xi^{\mathrm{Sloped}}}_{(\mathscr{D}_0)_{i_1}} \big[    \mathcal{H}( (\Gamma^{\prime}_C)_{\frac{3\pi}{2}} )          \big|  h =    \xi^{\mathrm{Sloped}}  \text{ } \text{ on }   (\mathscr{D}_0)_{i_1}^c \text{ } \cup   \partial (\mathscr{D}_0)_{i_1}  \big] \text{ }      \\ \equiv        \textbf{P}^{((\xi^{\mathrm{Sloped}})_3)^{\prime\prime}}_{(\mathscr{D}_0)_{i_1}} \big[    \mathcal{H}( (\Gamma^{\prime}_C)_{\frac{3\pi}{2}} )  \text{ } \big] \\ \text{ } \geq  \textbf{P}^{((\xi^{\mathrm{Sloped}})_3)^{\prime\prime}}_{[-m,m] \times [0,n^{\prime}N]} \big[    \mathcal{H}( (\Gamma^{\prime}_C)_{\frac{3\pi}{2}} )  \big] \text{ } \\ \text{ }  \overset{m \rightarrow +\infty}{\geq}   \textbf{P}^{((\xi^{\mathrm{Sloped}})_3)^{\prime}}_{\textbf{Z} \times [0,n^{\prime}N]} \big[     \mathcal{H}( (\Gamma^{\prime}_C)_{\pi} )         \big]       \text{ }    \text{ , } \text{ }   \tag{\textit{Horizontal crossing one}}
\end{align*}

\noindent yield an estimate corresponding to the probability of crossing $(\Gamma^{\prime}_C)_{\frac{3 \pi}{2}}$ in the longer direction, for $(\xi^{\mathrm{Sloped}})_3, ((\xi^{\mathrm{Sloped}})_3)^{\prime}$, $((\xi^{\mathrm{Sloped}})_3)^{\prime\prime} \sim \textbf{B}\textbf{C}^{\mathrm{Sloped}}$. Under the circumstance that all of the three crossings given above occur simultaneously,

\begin{align*}
   \text{ }   \mathscr{C}^{\prime}_{\mathrm{min}} \text{ }    \textbf{P}^{\xi^{\mathrm{Sloped}}}_{\Lambda} \big[      \mathcal{H}((\Gamma^{\prime}_C)_{\frac{\pi}{2}})   \cap\mathcal{H}((\Gamma^{\prime}_C)_{\pi}) \cap    \mathcal{H}((\Gamma^{\prime}_C)_{\frac{3\pi}{2}})    \big]   \overset{\mathrm{(FKG)}}{\geq} \mathscr{C}^{\prime}_{\mathrm{min}}  \textbf{P}^{\xi^{\mathrm{Sloped}}}_{\Lambda} \big[ \text{ }      \mathcal{H}((\Gamma^{\prime}_C)_{\frac{\pi}{2}}) \text{ } \big] \text{ }  \textbf{P}^{\mathrm{Sloped}}_{\Lambda} \big[ \text{ }    \mathcal{H}((\Gamma^{\prime}_C)_{\pi})      \text{ } \big] \\ \times  \text{ }     \text{ }  \textbf{P}^{\xi^{\mathrm{Sloped}}}_{\Lambda} \big[      \mathcal{H}((\Gamma^{\prime}_C)_{\frac{3\pi}{2}})    \big] \\  \geq \text{ } \mathscr{C}^{\prime}_{\mathrm{min}}   \text{ }    \prod_{\mathcal{D} \text{ } \in\text{ } \{    (\Gamma^{\prime}_C)_0  ,    (\Gamma^{\prime}_C)_{\frac{\pi}{2}}    ,   (\Gamma^{\prime}_C)_{\pi} ,  (\Gamma^{\prime}_C)_{\frac{3\pi}{2}}  \} } \text{ }          \textbf{P}^{\xi^{\mathrm{Sloped}}}_{\Lambda} \big[   \mathcal{H}(\mathcal{D})   ]       \text{ }  \\ \geq \text{ } \mathscr{C}^{\prime}_{\mathrm{min}}  \text{ }           \prod_{\mathcal{D} \text{ } \in\text{ }  \{    (\Gamma^{\prime}_C)_0  ,    (\Gamma^{\prime}_C)_{\frac{\pi}{2}}    ,  (\Gamma^{\prime}_C)_{\pi}  , (\Gamma^{\prime}_C)_{\frac{3\pi}{2}}  \} } \text{ }  \textbf{P}^{\xi^{\mathrm{Sloped}}}_{[-m,m] \times [0,n^{\prime}N]}   \big[     \mathcal{H}(\mathcal{D})             \big]         \\ \overset{ m \rightarrow + \infty}{\geq} \text{ }    \mathscr{C}^{\prime}_{\mathrm{min}}  \text{ }          \prod_{\mathcal{D} \text{ } \in\text{ }  \{     (\Gamma^{\prime}_C)_0  ,    (\Gamma^{\prime}_C)_{\frac{\pi}{2}}  ,  (\Gamma^{\prime}_C)_{\pi}  , (\Gamma^{\prime}_C)_{\frac{3\pi}{2}}  \} } \text{ }  \textbf{P}^{\xi^{\mathrm{Sloped}}}_{\textbf{Z} \times [0,n^{\prime}N]}   \big[    \mathcal{H}(\mathcal{D})               \big]        \text{ } \text{ . } \tag{\textit{Horizontal crossing two}}
      \end{align*}

      \noindent From the domains over which the final product above is taken, observe the lower bound,

     \begin{align*}
\text{ }   \mathscr{C}^{\prime}_{\mathrm{min}}  \text{ }   \prod_{\mathcal{D} \text{ } \in\text{ }  \{     (\Gamma^{\prime}_C)_0  ,    (\Gamma^{\prime}_C)_{\frac{\pi}{2}}   ,   (\Gamma^{\prime}_C)_{\pi}  ,  (\Gamma^{\prime}_C)_{\frac{3\pi}{2}}  \}}  \textbf{P}_{\textbf{Z} \times [0,n^{\prime}N]}^{\xi^{\mathrm{Sloped}}} \big[  [0,\lfloor{\delta_i n}\rfloor{}  ]  \times  \{ 0 \}  {\overset{h\geq ck}{\longleftrightarrow}}   \textbf{Z} \times  \{ n \}     \big] \text{ }  \\ \equiv \mathscr{C}^{\prime}_{\mathrm{min}}  \bigg[  \textbf{P}_{\textbf{Z} \times [0,n^{\prime}N]}^{\xi^{\mathrm{Sloped}}} \big[  [0,\lfloor{\delta_i n}\rfloor{}  ]  \times \{ 0 \}  {\overset{h\geq ck}{\longleftrightarrow}}   \textbf{Z} \times \{ n \}    \text{ } \big]   \bigg]^{ | \mathcal{D} | } \\ 
\geq \text{ }     \mathscr{C}^{\prime}_{\mathrm{min}}                  \bigg[  \textbf{P}^{\xi^{\mathrm{Sloped}}}_{\textbf{Z} \times [0,n^{\prime}N]}  \big[  [0,\lfloor{\delta_i n}\rfloor{}  ]  \times  \{ 0 \}  {\overset{h\geq ck}{\longleftrightarrow}}   [-m,m] \times  \{ n \}    \big]  \bigg]^{4  \mathscr{N}}    \text{ , } \text{ } \tag{\textit{Third constant}}
\end{align*}

\noindent where to obtain the expression given in the $\textit{(Third constant)}$, we bound the horizontal crossing across each $\mathcal{D}$ with,

\begin{align*}
  \textbf{P}^{\xi^{\mathrm{Sloped}}}_{\Lambda} \big[  \mathcal{H}(\mathcal{D}) \big]  \geq         \textbf{P}_{\textbf{Z} \times [0,n^{\prime}N]}^{\xi^{\mathrm{Sloped}}} \big[   [0,\lfloor{\delta_i n}\rfloor{}  ]  \times \{ 0 \}  {\overset{h\geq ck}{\longleftrightarrow}}  \text{ } [-m,m] \times  \{ n \}     \big]          \text{ } \text{ , } 
\end{align*}

\noindent in which the annulus crossing event admits the following lower bound, from (\textit{Third constant}),

\begin{align*}
\textbf{P}^{\xi^{\mathrm{Sloped}}}_{\Lambda_{12n}} \big[ \mathcal{O}_{h \geq ck - l^{\prime} - n}(6n, n_R n )  \big]     \text{ }  \overset{m \longrightarrow +\infty}{\geq}    \mathscr{C}^{\prime}_{\mathrm{min}}                  \bigg[    \textbf{P}^{\xi^{\mathrm{Sloped}}}_{\textbf{Z} \times [0,n^{\prime}N]}  \big[ [0,\lfloor{\delta^{\prime} n}\rfloor{}  ]  \times \{ 0 \}  {\overset{h\geq ck}{\longleftrightarrow}}  \textbf{Z} \times  \{ n \}    \big] \bigg]^{4  \mathscr{N}}      \text{ , } 
\end{align*}

\noindent with the above height required by the annulus crossing $\mathcal{O}$, where,

\begin{align*}
    l^{\prime} \equiv ((\xi^{\mathrm{Sloped}})_3)^{\prime\prime} - \xi^{\mathrm{Sloped}} \text{ } \text{ , } 
\end{align*}

\noindent is the change in boundary conditions, between the measures $\textbf{P}^{((\xi^{\mathrm{Sloped}})_3)^{\prime\prime}}[ \text{ } \cdot \text{ }]$, and $\textbf{P}^{\xi^{\mathrm{Sloped}}_3}[ \text{ } \cdot \text{ }]$, used to previously obtain (\textit{Horizontal crossing one}), and (\textit{Horizontal crossing two}), with $n$ satisfying,

\begin{align*}
     n  < h - ck + l^{\prime}       \text{ } \text{ , } 
\end{align*}

\noindent from which we conclude the argument, upon potential readjustment of $c$, $l^{\prime}$, and $n$ defined above. \boxed{}

\bigskip

   \begin{figure}
\begin{align*}
\includegraphics[width=0.78\columnwidth]{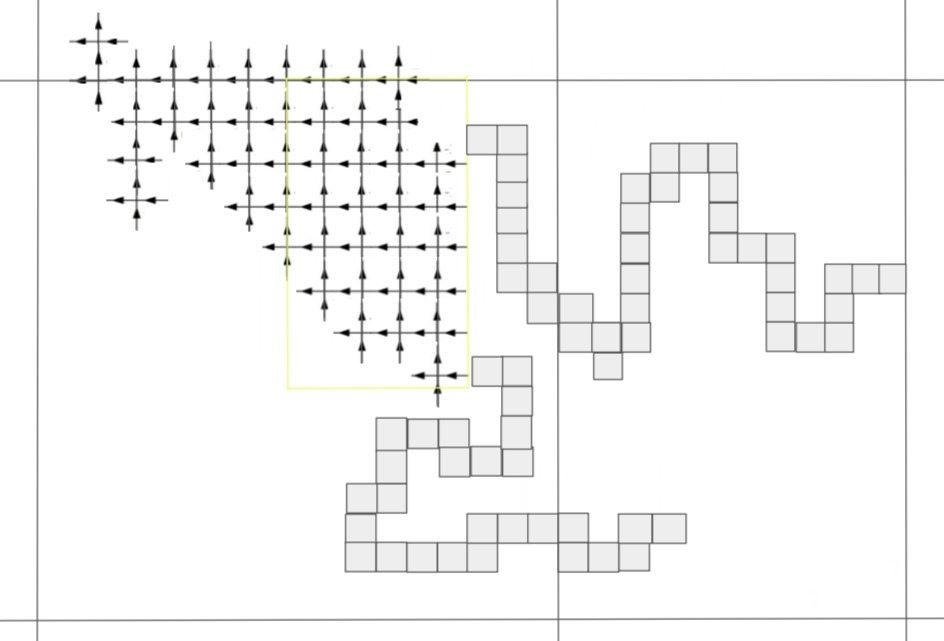}\\
\end{align*}
\caption{\textit{Crossing events that are reliant upon more faces, as shown with the path above that intersects a second line to the left of the \textit{freezing cluster}, occur with higher probability than crossing events that are reliant upon fewer faces, as shown with the path that intersects only the first line to the left of the \textit{freezing cluster}.}}
\end{figure}

  \begin{figure}
\begin{align*}
\includegraphics[width=0.98\columnwidth]{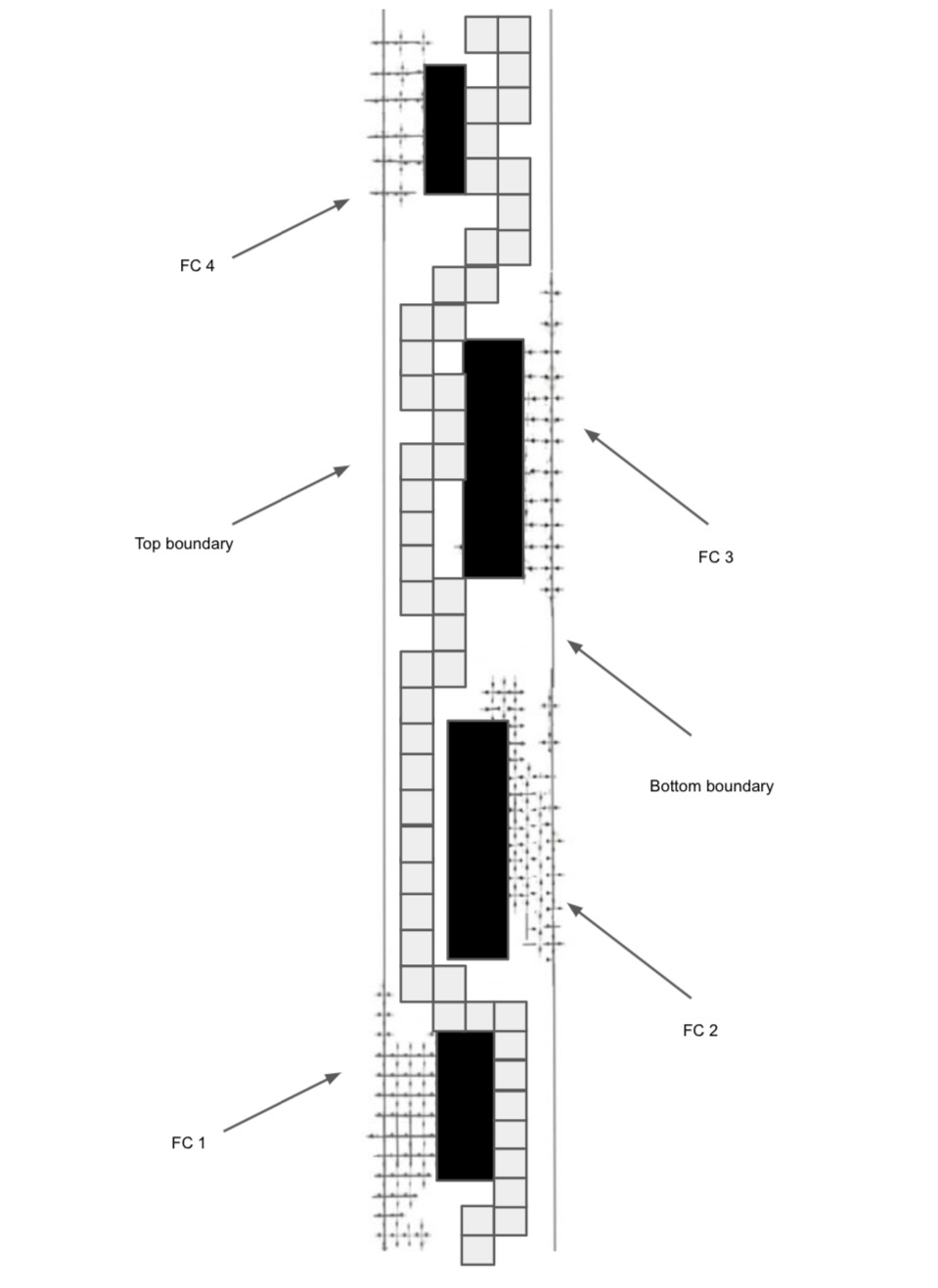}\\
\end{align*}
\caption{\textit{A depiction of four freezing clusters within the strip from Figure 5}. Above, four black finite subvolumes in the strip intersect a proportion of the frozen faces within each \textit{freezing cluster}. A macroscropically crossing face path, with individual faces constituting the path, is depicted in gray.}
\end{figure}

\section{Cylindrical crossing events of the sloped free-energy landscape}

\subsection{Overview}

\noindent From comparisons between the crossing probability across the annulus with the crossing probability across the strip, we now incorporate estimates for sloped Gibbs states in the cylinder. In doing so, we provide a quantification on the exponent to which the free energy is raised, and establish comparisons between the free energy exponent provided in {\color{blue}[11]} for flat boundary conditions, in which a map $\textbf{T}$ is constructed for reversing the orientation of loops, and hence, of six-vertex, configurations in the cylinder. Conditionally upon the existence of two $\mathrm{x}$-crossings which are separated by some strictly positive distance $L$, in contrast to the map $\textbf{T}$ that is introduced for flat boundary conditions, for sloped boundary conditions another map $\textbf{T}^{\prime}$ satisfies a different set of properties, namely that the map still reverse the orientation of a portion of the directed loops in the cylinder, however, with the additional stipulation that the number of possible preimages of the directed loop configurations satisfy a different condition. That is, instead of the number of preimages being dependent upon $N^2 2^{2M/\alpha}$, which is a factor obtained from $N^2$ possible ways of choosing the first edges, and from $2M/\alpha$ ways of choosing edges between $\gamma_{i^{*}}$ and $\gamma_{L+i^{*}}$, the total number of possible preimages, sampled from the event intersection $\Omega^{{\mathrm{S}},0} \cap \mathcal{B}(L)$ given $\omega \in \text{ } \Omega^{{\mathrm{S}},0} \text{ } \equiv \text{ } \Omega^{{\mathrm{S}},(0)}_{\mathrm{6V}}$, and hence satisfies a similar upper bound dependent upon $N^2$. As the slope of the boundary conditions increases, the upper bound on the total number of possible preimages which is given by the pullback $(\textbf{T}^{\prime})^{-1}(\omega)$. Formally, we define $\textbf{T}^{\prime}$ below, and afterwards, characterize differences in the orientation-reversing procedure for directed loops. Furthermore, in what follows, denote the probability measure supported over the cylinder with $\textbf{P}_{\mathscr{O}_{N,M}} [ \text{ } \cdot \text{ } ] \equiv \textbf{P}_{\mathscr{O}}[ \text{ } \cdot \text{ } ]\text{ } $, for $N$ even, where the finite cylindrical volume is bound by the top and bottom boundaries of the cylinder, in addition to the intermediate cylindrical volume between the bottom and top faces, consisting of $NM$ faces, namely $\mathscr{O} \equiv \mathscr{O}_{N,M}\text{ } \equiv\text{ }  \mathcal{T}\big( \text{ }    \mathscr{O}_{N,M}     \text{ } \big) \text{ }  \cup \text{ }   \mathcal{I}\big( \text{ }    \mathscr{O}_{N,M}     \text{ } \big)  \text{ } \cup \text{ } \partial \big( \text{ }   \mathcal{I}\big( \text{ }    \mathscr{O}_{N,M}     \text{ } \big)        \text{ } \big) \text{ }  \cup \text{ }   \mathcal{B}\big( \text{ }    \mathscr{O}_{N,M}     \text{ } \big)\text{ } $, with $|\mathcal{B}\big( \text{ }    \mathscr{O}_{N,M}     \text{ } \big) | = 2n$. Finally, for arguments in \textbf{Lemma} \textit{4.5} and \textbf{Lemma} \textit{4.6} closer to the end of the section, denote the union of vertical $\mathrm{x}$-crossing events with,

\begin{align*}
 \mathscr{U}\mathscr{V} \text{ } \equiv \text{ } \underset{j \in \textbf{N} :  \mathscr{F}_1 \in \mathcal{B} (\mathscr{O} )  ,  \mathscr{F}_2 \in \mathcal{T} (  \mathscr{O} )}{\bigcup}     \big\{   \mathscr{F}_1 \underset{\mathcal{I}_j ( \mathscr{O} )}{\overset{h \geq ck}{\longleftrightarrow}}      \mathscr{F}_2  \big\}          \text{ } \text{ , } 
\end{align*}

\noindent in which, from the decomposition of $\mathscr{O}$ into top, interior, interior boundary, and bottom domains, $\mathcal{I}_j \big( \text{ } \mathscr{O} \text{ } \big) \subset \mathscr{O}$ for all $j$. Finally, $\mathcal{I} \big( \text{ } \mathscr{O} \text{ } \big)$, as provided above in the decomposition of $\mathscr{O}$, itself can be partitioned into countably many regions, each of which are contained within $\mathscr{O}$, so that only the boundaries of each $\mathcal{I}_j \big( \text{ } \mathscr{O} \text{ } \big)$ intersect,

\begin{align*}
     \mathcal{I} \big(  \mathscr{O}  \big)     \text{ } \equiv \text{ }  {\underset{\text{countably many }j : j \in \textbf{N}}{\bigcup}} \text{ } \mathcal{I}_j \big( \mathscr{O}  \big)    \text{ } \text{ . } 
\end{align*}

\subsection{Defining the map}

\noindent \textbf{Definition} \textit{9} (\textit{orientation-reversing map for directed loops under sloped boundary conditions in the limit of full packing}). We define the map, for $L>0$ and from the event $\mathcal{B}(L)$ below, with the associated path $\gamma_{i^{*}}$ located on $\mathcal{B}\big( \text{ } \mathscr{O}_{N,M}  \text{ } \big)$ at position $i^{*}$,

\begin{align*}
      \mathcal{B}(L) \equiv  \bigg\{  \forall \text{ } 0 \leq i^{*} \leq 2 n  ,  L <  2n  , \text{ }  \exists \text{ } \gamma_{i^{*}}         ,     \gamma_{L+i^{*}}  :  \gamma_{i^{*}} \equiv \big\{  \mathscr{F}_{\mathrm{Root}_1} \overset{h \geq ck}{\longleftrightarrow}  \mathscr{F}^{\mathcal{T}}   \big\}    ,   \gamma_{L+i^{*}} \equiv     \big\{ \mathscr{F}_{\mathrm{Root}_2}    \overset{h \geq ck}{\longleftrightarrow}    \mathscr{F}^{\mathcal{T}^{\prime}}      \big\}    \bigg\}          \text{ , } 
\end{align*}

\begin{figure}
\begin{align*}
\includegraphics[width=0.75\columnwidth]{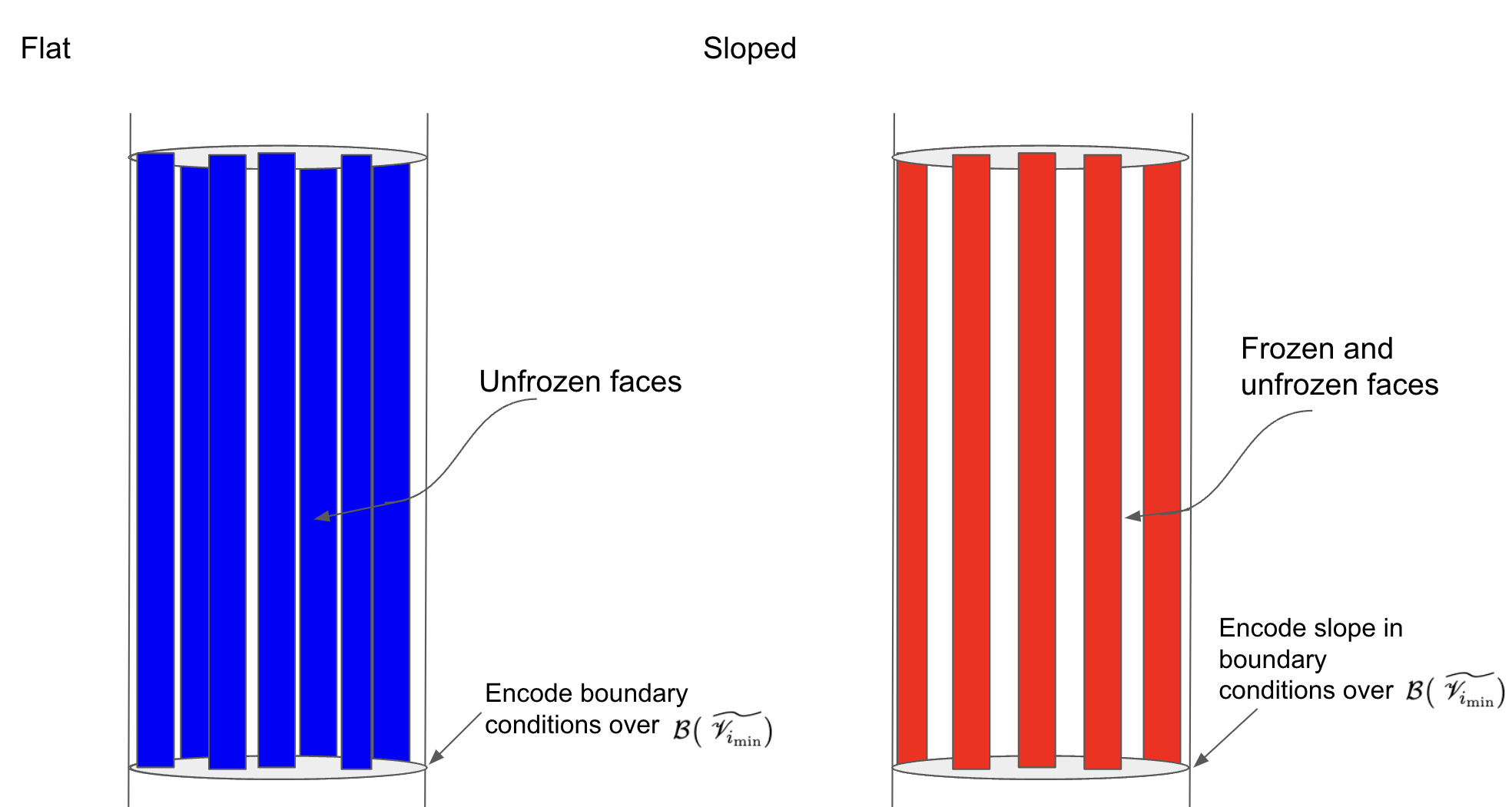}\\
\end{align*}
\caption{\textit{A heuristic of the expected dependency of the distribution of symmetric domains in the cylinder with flat, and with sloped, boundary conditions}. From previous results for flat boundary conditions, the interior of the cylinder above the bottom of the lowest occurring \textit{symmetric} domain intersecting $\mathrm{Slice}_{y-1,y,y+1}$ consists entirely of unfrozen faces, which are contained in some of the paths going from the bottom to top of the cylinder in blue. In comparison, under sloped boundary conditions, through boundary conditions of the height function over $\mathcal{B} \big( \text{ } \widetilde{\mathscr{V_{i_{\mathrm{min}}}}} \text{ } \big)$, it is expected that the distribution of \textit{symmetric} domains in the interior of the cylinder above $\mathcal{B} \big( \text{ } \widetilde{\mathscr{V}_{i_{\mathrm{min}}}} \text{ } \big)$, and intersecting $\mathrm{Slice}_{y-1,y,y+1}$, consist of frozen and unfrozen faces. A priori, one could expect that the distribution of blue regions, corresponding to points which are deemed to satisfy the \textit{good} property about which domains are placed, would contain more unfrozen faces in the cylinder under flat boundary conditions than in the cylinder under sloped boundary conditions.}
\end{figure}

\noindent where $\mathscr{F}_{\mathrm{Root}_1}$ and $\mathscr{F}_{\mathrm{Root}_2}$ are two faces belonging to the bottom of the cylinder that are separated by distance $L$, and,  $\mathscr{F}^{\mathcal{T}}$ and $\mathscr{F}^{\mathcal{T}^{\prime}}$. In the arguments throughout the section, fix $\mathscr{F}_{\mathrm{Root}_1} \equiv \mathscr{F}_1$, and $\mathscr{F}^{\mathcal{T}} \equiv \mathscr{F}_2$, which denote two faces on the top of the cylinder. The map,

\begin{align*}
    \mathcal{T}^{\prime}  :                       \Omega^{{\mathrm{S}},2L}             \longrightarrow           \Omega^{0,2L}          \cap     \mathcal{B}(L) \text{ }     \text{ , } 
\end{align*}

\noindent in which the image of the configuration with reversed loop orientation satisfies the inequality,

\begin{align*}
  w_{6\mathrm{V}}\big(\mathcal{T}^{\prime}(\omega)\big)   \geq         c^{-2M/\alpha-i }             w_{6\mathrm{V}} \big( \omega \big)    \text{ } \text{, for }     1 < i <   -2M/\alpha       \text{ }   \text{ , } \tag{\textit{Property One}} 
\end{align*}

\noindent which is denoted as (\textit{Property One}), where the change in six-vertex weighs in the lower bound for the weight of the reverse directed loop configuration, $w_{6\mathrm{V}}\big(\mathcal{T}^{\prime}(\omega)\big)$, is expressed in terms of the change in weights within the frozen, and unfrozen, regions, with respective number of edges $MN$, $M^{\prime} N^{\prime}$

\noindent The second requirement states,

\begin{align*}
      \big|  (\mathcal{T}^{\prime})^{-1} \big(   \omega       \big)       \big|    \leq     N^2            2^{-2M/\alpha-i }      \text{ } \text{, for }     1 < i <   -2M/\alpha     \text{ }   \text{ , } \tag{\textit{Property Two}} 
\end{align*}

\noindent $\forall \text{ } \omega \in \Omega^{\mathrm{S},2L}$, where the number of edges belonging to the frozen portion of the finite volume is denoted with the proportion of edges belonging to frozen faces over the edge set $E\big(\mathscr{O}     \big)   $,

\begin{align*}
       \mathscr{E}_{\mathrm{frozen}}  \equiv              \big\{     \mathscr{P}_1  ,  \mathscr{P}_2  \in F\big( \mathscr{O}      \big)        :       \mathscr{P}_1 \text{ and } \mathscr{P}_2 \text{ have identically oriented arrows} \big\}        \text{ }  \text{ , } 
\end{align*}

\noindent where a \textit{plaquette} $\mathscr{P}$ is denoted as the union of four faces incident to every vertex satisfying the ice rule, such that for any collection of \textit{four} edges in each such $\mathscr{P}$, denoted by $(e_{\mathscr{P}})_1 , \cdots , (e_{\mathscr{P}})_4$,

\begin{align*}
  \bigcap_{1 \leq i \leq 4}  \text{ } e_{\mathscr{P}_i }        \text{ } \equiv\text{ }     v_{\mathscr{P}}   \text{ }  \text{ }   \text{ , } 
\end{align*}

\noindent where $v_{\mathscr{P}} \in V (\text{ } \mathscr{O} \text{ } )$.

\subsection{Section results}

\noindent First, we introduce the \textit{sloped} free energy function.

\bigskip

\noindent \textbf{Definition} \textit{10} (\textit{sloped free energy function}). Define the function $g_c : ( -1 , 1 ) \longrightarrow \textbf{R}^{+}$, which, for a strictly positive arrow unbalance parameter $\alpha$, is defined as,

\begin{align*}
g_c (\alpha)  \text{ } \equiv \text{ } \underset{N \longrightarrow + \infty}{\underset{N \text{ } \text{even}}{\mathrm{lim}}}  \underset{M \longrightarrow + \infty}{\mathrm{lim}}  \text{ }  \big( \text{ } NM \text{ } \big)^{-1} \text{ } \mathrm{log} \big(   Z^{\mathrm{Sloped},\alpha}_{M,N}  \big)  \text{ } \equiv \text{ } \underset{N \longrightarrow + \infty}{\underset{N \text{ } \text{even}}{\mathrm{lim}}} \text{ } \underset{M \longrightarrow + \infty}{\mathrm{lim}}  \big( \text{ } NM \text{ } \big)^{-1} \text{ } \mathrm{log} \big(    Z^{\mathrm{S},\alpha}_{M,N} \big)  \text{ } \text{ , } 
\end{align*}

\noindent for the partition function $Z^{\mathrm{S},\alpha}_{M,N}$ taken under sloped boundary conditions.

\bigskip

\noindent \textit{Remark}. In comparison to $f_c$ of {\color{blue}[11]}, $g_c$ is defined so that it coincides with $f_c$ over $(-\frac{1}{2} , \frac{1}{2} ) \subset (-1,1)$ in the codomain.

\bigskip

\noindent Next, introduce how $g_c$ can serve as a lower bound to crossing probabilities in $\mathscr{O}$.

\bigskip

\noindent \textbf{Lemma} \textit{4.1} (\textit{estimation of the six-vertex free energy from cylindrical crossings}). From $\mathcal{T}^{\prime}$, denote,

\begin{align*}
  \mathcal{S}_{\mathscr{F}} \equiv   \big\{      \mathscr{F} \in F(  \textbf{Z}^2  )  :  \mathscr{F} \cap    \mathcal{B} \big(  \mathscr{O} \big) \neq \emptyset   \big\}       \text{ }   \text{ , } 
\end{align*}

\noindent and the corresponding cylindrical crossing event from some face on $\mathcal{T}(\mathscr{O}_{N,M})$,

\begin{align*}
      \mathscr{C}_{\mathscr{O}} \equiv                \big\{  \mathscr{F}_1  \in   \mathcal{B} \big(   \mathscr{O}    \big)  ,  \mathscr{F}_2  \in  \mathcal{T}\big(   \mathscr{O}     \big)  :  \big\{  \mathscr{F}_1 \underset{\mathcal{I} (  \mathscr{O} )}{\overset{ h \geq ck}{\longleftrightarrow} }     \text{ } \mathscr{F}_2 \big\}    \big\}    \text{ , } 
\end{align*}

\noindent from which the limit infimum of the natural logarithm, scaled in the cardinality of the cylindrical finite volume in the weak limit, of the supremum of the crossing probability supported over $\mathscr{O}$,

\begin{align*}
\text{ }     \underset{ N \longrightarrow + \infty}{\underset{N \text{ even}}{\mathrm{lim\text{ }  inf}}}  \text{ }    \underset{M \longrightarrow \text{ } +\infty}{\mathrm{lim \text{ } inf}}     \text{ }   \big( \text{ }   N M     \text{ } \big)^{-1}  \mathrm{log} \text{ } \bigg[   \underset{\mathscr{F}_2 \in \mathcal{T}(  \mathscr{O} )}{\underset{\mathscr{F}_1 \in \mathcal{S}_{\mathscr{F}}}{\mathrm{sup}}} \textbf{P}^{\xi^{\mathrm{Sloped}}}_{\mathscr{O}}  \big[    \mathscr{F}_1 \underset{\mathcal{I}( \mathscr{O})}{\overset{ h \geq ck}{\longleftrightarrow} }   \mathscr{F}_2 \big]    \bigg]       \text{ }         \text{ , } 
\end{align*}

\noindent admits the free-energy lower bound, upon first sending $M \longrightarrow +\infty$ and then $N \longrightarrow +\infty$, as,

\begin{align*}
      g_c (     \beta                  )    -    g_c ( 0 )   \text{ , } 
\end{align*}

\noindent where $f$ denotes the \textit{cylindrical free-energy} of the six-vertex model provided in the result of \textbf{Theorem} $\textit{1.4}$, and,

\begin{align*}
 \textbf{P}^{\xi^{\mathrm{Sloped}}}_{\mathscr{O}}  \big[   \mathscr{F}_1 \underset{\mathcal{I}( \mathscr{O} )}{\overset{ h \geq ck}{\longleftrightarrow} }     \text{ } \mathscr{F}_2  \big]  \equiv   \textbf{P}^{\xi^{\mathrm{Sloped}}}_{\mathscr{O}} \big[   {\underset{\text{countably many }j}{\bigcap}} \big\{ \text{ }       \mathscr{F}_1 \underset{\mathcal{I}_j (  \mathscr{O} )}{\overset{h \geq ck}{\longleftrightarrow}}      \mathscr{F}_2         \big\}    \big]     \equiv      \textbf{P}^{\xi^{\mathrm{Sloped}}}_{\mathscr{O}} \big[\mathscr{U} \mathscr{V}   \big]      \text{ } \text{ , }
\end{align*}

\noindent under the same choice of $\mathscr{F}_1$ and $\mathscr{F}_2$ introduced in \textit{4.1}, namely for $\mathscr{F}_1 \in \mathcal{B} \big(  \mathscr{O}  \big)$, and $\mathscr{F}_2 \in \mathcal{T} \big(  \mathscr{O}  \big)$. 

\bigskip

\noindent Instead of immediately requiring that the natural logarithm of the crossing probability between $\mathscr{F}_1$ and $\mathscr{F}_2$, as defined above, below we provide a free-energy lower bound dependent upon $\mathcal{B}(L)$.

\bigskip

\noindent \textbf{Lemma} \textit{4.2} (\textit{estimation of the six-vertex free energy from restricted cylindrical crossings}). If $\mathcal{B}(L)$ occurs, the natural logarithm of the probability given below,

\begin{align*}
        \underset{ N \longrightarrow + \infty}{\underset{N \text{ even}}{\mathrm{lim\text{ }  inf}}} \text{ }    \underset{M \longrightarrow \text{ } +\infty}{\mathrm{lim \text{ } inf}}     \text{ }   \big( \text{ }   N M     \text{ } \big)^{-1}  \mathrm{log} \big\{  \textbf{P}^{\xi^{\mathrm{Sloped}}}_{\mathscr{O}}  \big[        \mathcal{B}(L)
          \big]        \big\}       \text{ , } 
\end{align*}

\noindent admits an identical free-energy lower bound provided in $\textbf{Lemma}$ \textit{4.1}, upon first sending $M \longrightarrow +\infty$ and then $N \longrightarrow +\infty$, as,

\begin{align*}
       g_c (     \beta                  )       \text{ } - \text{ }    g_c ( 0 )   \text{ , } 
\end{align*}

\noindent where $f^{\prime}_c$ denotes the sloped free energy of the six-vertex model, as a counterpart to the free energy $f_c$ for flat boundary conditions provided in $\textbf{Theorem}$ \textit{1.4}.

\bigskip

\noindent With the two previous results, we demonstrate that the two properties of $\mathcal{T}^{\prime}$ hold with the following \textit{Proof}.

\bigskip

\noindent \textit{Proof of Lemma 4.2}. From ($\textit{Property One}$) and ($\textit{Property Two}$) of $\mathcal{T}^{\prime}$ provided in $\textbf{Definition}$ \textit{8}, we argue that each of the three items given below holds.  Fix some $\omega \in \Omega^{S, 2L}$:

\begin{figure}
\begin{center}
\begin{tikzpicture}
\filldraw[color=blue!80, fill=red!0, very thick](-1,0) circle (6);\node[inner sep=0pt] (test) at (0,10.5)
    {\includegraphics[width=.07\textwidth]{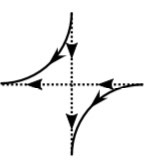}};
\node[inner sep=0pt] (test) at (0,0)
    {\includegraphics[width=.08\textwidth]{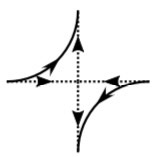}};\node[inner sep=0pt] (test) at (-5,0)
    {\includegraphics[width=.08\textwidth]{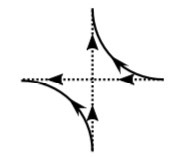}};\node[inner sep=0pt] (test) at (-10,-0.5)
    {\includegraphics[width=.08\textwidth]{IMG_0122.jpeg}};\node[inner sep=0pt] (test) at (-3,-2.5)
    {\includegraphics[width=.08\textwidth]{IMG_0122.jpeg}};\node[inner sep=0pt] (test) at (-12,-5.5)
    {\includegraphics[width=.08\textwidth]{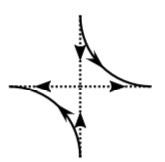}};\node[inner sep=0pt] (test) at (-10,-2.5)
    {\includegraphics[width=.08\textwidth]{IMG_0122.jpeg}};\node[inner sep=0pt] (test) at (-5,-2.5)
    {\includegraphics[width=.07\textwidth]{IMG_0123.jpeg}};\node[inner sep=0pt] (test) at (-5,6.8)
    {\includegraphics[width=.08\textwidth]{IMG_0123.jpeg}};\node[inner sep=0pt] (test) at (-2,7.8)
    {\includegraphics[width=.08\textwidth]{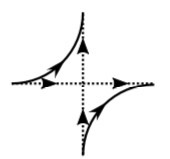}};\node[inner sep=0pt] (test) at (-7.8,5.8)
    {\includegraphics[width=.07\textwidth]{IMG_0123.jpeg}};\node[inner sep=0pt] (test) at (-1,5.8)
    {\includegraphics[width=.08\textwidth]{IMG_0124.jpeg}};\node[inner sep=0pt] (test) at (-10,5.8)
    {\includegraphics[width=.08\textwidth]{IMG_0124.jpeg}};\node[inner sep=0pt] (test) at (8,-3.5)
    {\includegraphics[width=.08\textwidth]{IMG_0124.jpeg}};\node[inner sep=0pt] (test) at (-2,7.8)
    {\includegraphics[width=.08\textwidth]{IMG_0124.jpeg}};\node[inner sep=0pt] (test) at (-1,-7.8)
    {\includegraphics[width=.08\textwidth]{IMG_0122.jpeg}};\node[inner sep=0pt] (test) at (-3,-5.8)
    {\includegraphics[width=.08\textwidth]{IMG_0125.jpeg}};\node[inner sep=0pt] (test) at (0,-8.8)
    {\includegraphics[width=.08\textwidth]{IMG_0124.jpeg}};\node[inner sep=0pt] (test) at (16.5,3.5)
    {\includegraphics[width=.08\textwidth]{IMG_0125.jpeg}};
\end{tikzpicture}
\end{center}
\caption{\textit{A visual schematic for one realization of a fully packed loop configuration over the topmost boundary of the cylinder for free energy estimates of the six-vertex model}. With an above view of the topmost boundary of the cylinder which includes the final face of each path beginning from the bottommost face of the cylinder, taking the weak volume limit yields a series of cylinders with increasing diameter. Within the cylindrical finite volume that is highlighted in blue, a portion of the fully packed loop configurations corresponding to six-vertex weights are depicted, in which directionality of full packed loop configurations is established in correspondence with directionality of arrows about each vertex. In the complementary region of the square lattice outside of the higlighted blue finite volume, in  the weak infinite volume limit additional statistical contributions to the probability of a crossing event occurring will also be dependent upon the values of the height function for all faces below the topmost faces on the topmost boundary of the cylinder. In comparison to flat Gibbs states within the cylinder, sloped Gibbs states in the cylinder exhibit differences in the value of the height function that is influenced by boundary faces, two of which are depicted, which have nonempty intersection with the topmost boundary.}
\end{figure}

\begin{itemize}
\item[$\bullet$] \underline{\textit{Preimage and Image Spaces of the loop-reversing map} $\mathcal{T}^{\prime}$}: To verify that configurations in $\Omega^{S,2L}$ are pushed forwards to configurations $\Omega^{S,2L} \cap \mathcal{B}(L)$ under $\mathcal{T}^{\prime}$, observe that the distribution of loops with reversed orientation, under $w_{\mathrm{6V}} \big( \mathcal{T}^{\prime}(\omega) \big)$, is determined by the unbalance between upwards and downwards facing arrows in the cylinder. From the lower bound provided on the weight of the orientation reversed configuration, $w_{\mathrm{6V}} \big( \mathcal{T}^{\prime}(\omega) \big)$, the change in statistical weight when reversing the loop orientation, namely performing the transformation $\omega \mapsto \mathcal{T}^{\prime}(\omega)$, is expressed in terms of contributions from the change of six-vertex weight from reversing the orientation of all faces belonging to a \textit{freezing cluster}, as introduced in \textbf{Definition} \textit{3}. The preimage $\big( \mathcal{T}^{\prime}\big)^{-1}(\omega)$ of the orientation-reversed configuration under the inverse map $\big(\mathcal{T}^{\prime}\big)^{-1}  :  \Omega^{S,2L} \cap \mathcal{B}(L) \longrightarrow  \Omega^{S,2L}$ hence satisfies the property that the configuration in $\Omega^{S,2L} \cap \mathcal{B}(L)$ has as many upward as downards paths.

    \item[$\bullet$] \underline{\textit{Property One}}. From such a configuration $\omega$, the reversal of the orientation of directed loops obtained from the existence of a leftmost, and rightmost, top to bottom crossings as stipulated by the two crossing paths of $\mathcal{B}(L)$, implies that the lower bound in the weight of the configuration obtained from $\mathcal{T}^{\prime}$, $w_{\mathrm{6V}}\big(  \mathcal{T}^{\prime}( \omega)            \big)$, can at most be bound below with $c^{-2M/\alpha-i }             \text{ }    w_{6\mathrm{V}} \big( \omega \big)    \text{ } \text{, for }     1 < i <   -2M/\alpha$, as each six-vertex weight changes with a factor of at most $c$. 
    
      \item[$\bullet$] \underline{\textit{Property Two}}. From such a configuration $\omega$, to justify that the upper bound on the number of possible preimages of $\omega$, under $(\mathcal{T}^{\prime})^{-1}$, is of the given form, observe the following. The finite volume factors present in the upper bound, over a region with area $NM$ contained within $\mathscr{O}$ represent the total area spanned by the cylindrical finite volume. For the remaining factor in the upper bound, observe that this term represents the total number of ways to pick edges from $2^{-2M/\alpha-i }$ total possible number of edges, for $1 < i <   -2M/\alpha$. Hence, $N^2 2^{-2M/\alpha-i }$ represents all of the manners in which such configurations from the preimage can be determined.
\end{itemize}

\noindent From \textit{Property One} and \textit{Property Two}, to conclude, the free energy lower bound, for $N,M < + \infty$ with $N$ even,

\begin{align*}
       \big( \text{ }   N M     \text{ } \big)^{-1}      \mathrm{log}  \bigg[     \underset{\mathscr{F}_2 \in \mathcal{T} (  \mathscr{O} )}{\underset{\mathscr{F}_1 \in \mathcal{S}_{\mathscr{F}}}{\mathrm{sup}}}\textbf{P}^{\xi^{\mathrm{Sloped}}}_{\mathscr{O}}  \big[    \mathscr{F}_1 \underset{\mathcal{I}( \mathscr{O})}{\overset{ h \geq ck}{\longleftrightarrow} }    \mathscr{F}_2 \big]    \bigg]   \text{ } \text{ , } \end{align*}
       
       \noindent which is equivalent to the following probability below, also normalized in $ \big( \text{ }   N M     \text{ } \big)^{-1}$, which admits the lower bound from an application of ($\mathrm{FKG}$),

       \begin{figure}
    \centering
    \includegraphics[width=0.86\columnwidth]{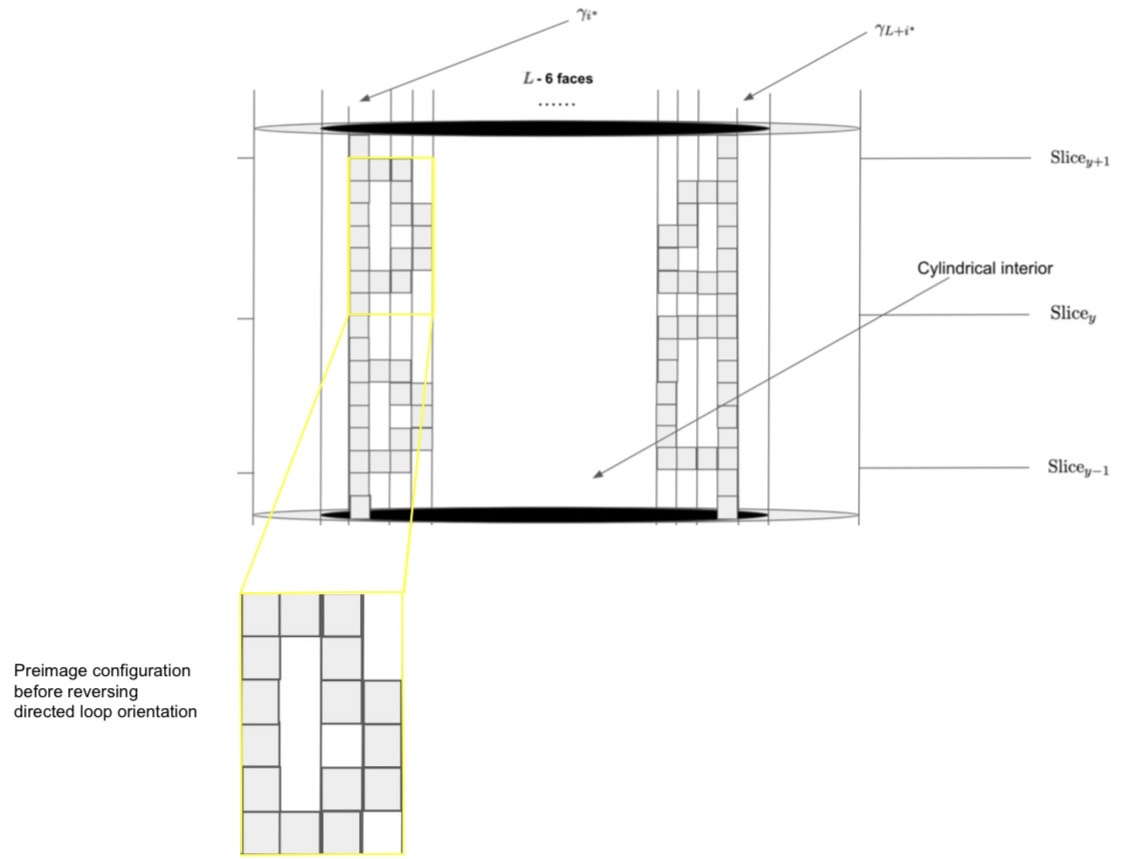}
    \caption{\textit{A depiction of the faces near the left and right boundaries of the cylinder induced by crossings $\gamma_L$, and by $\gamma_{L+i^{*}}$, respectively, each of which is highlighted in grey.}}
    \end{figure}

       \begin{figure}
    \centering
    \includegraphics[width=0.84\columnwidth]{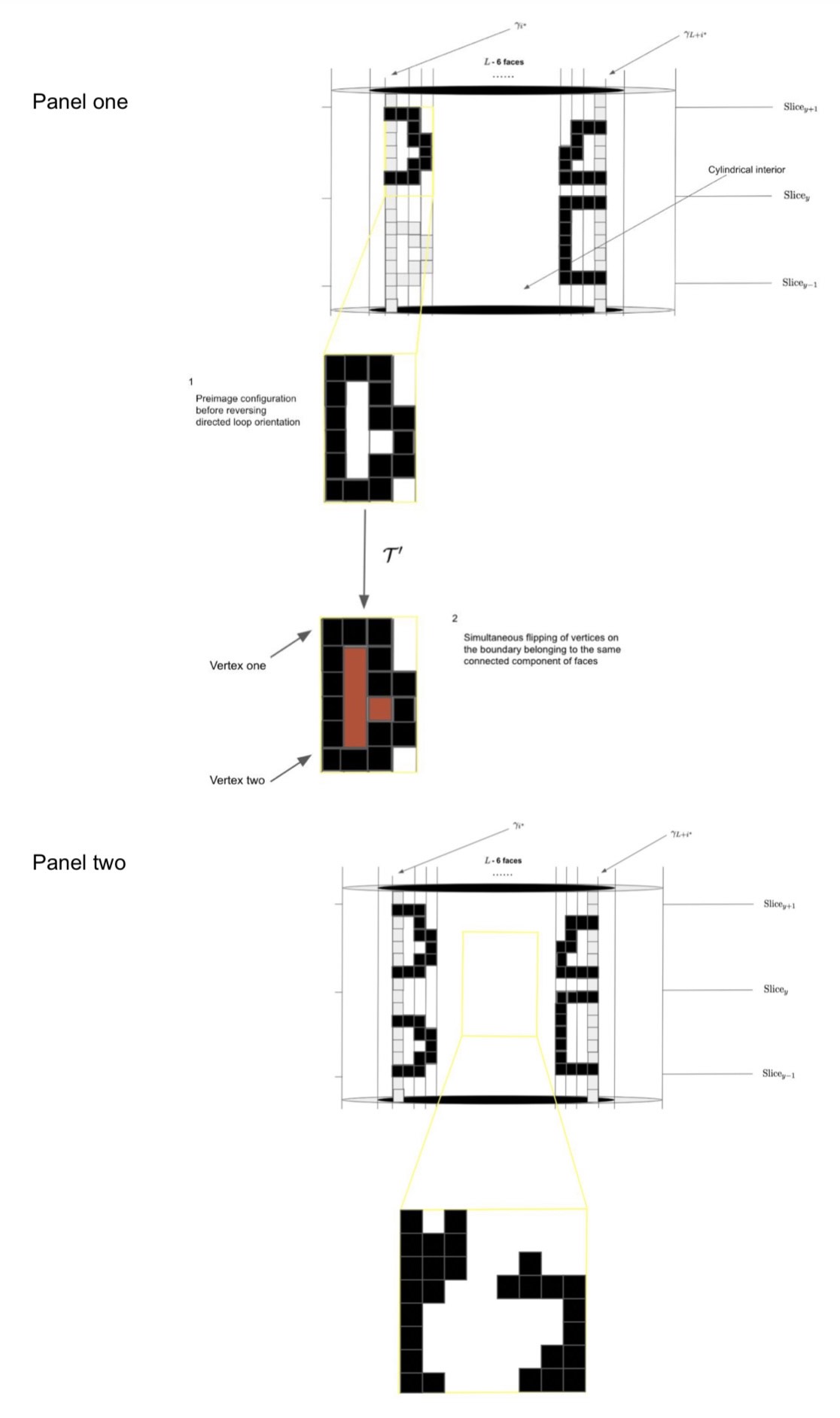}
    \caption{\textit{A further adaptation of the depiction of the faces provided in the previous figure.} In Panel one, part of the boundary faces to the left and right of the cylinder are highlighted in black, whose directed orientation in the full-packing limit is reversed when pushed forwards under $\mathcal{T}^{\prime}$. In Panel two, the distribution of frozen faces within the interior of the cylinder spanned by the yellow finite volume, is shown. }
    \end{figure}

       \begin{align*}
       \big( \text{ }   N M     \text{ } \big)^{-1}             \mathrm{log}  \bigg[     \underset{\mathscr{F}_2 \in \mathcal{T} ( \mathscr{O})}{\underset{\mathscr{F}_1 \in \mathcal{S}_{\mathscr{F}}}{\mathrm{sup}}} \textbf{P}^{\xi^{\mathrm{Sloped}}}_{\mathscr{O}}  \big[     \underset{\text{countably many }j: j \in \textbf{N}}{\bigcap}\text{ }  \big\{ \text{ }     \mathscr{F}_1 \underset{\mathcal{I}_j (  \mathscr{O} )}{{\overset{h \geq ck}{\longleftrightarrow }}}    \mathscr{F}_2   \big\}    \big]  \bigg] \\ \text{ } \overset{\mathrm{(FKG)}}{\geq} \text{ }      \big( \text{ }   N M     \text{ } \big)^{-1}           \text{ }     \mathrm{log}  \bigg[  \underset{\mathscr{F}_2 \in \mathcal{T}(  \mathscr{O} )}{\underset{\mathscr{F}_1 \in \mathcal{S}_{\mathscr{F}}}{\mathrm{sup}}}     \bigg[              \underset{\text{countably many }j: j \in \textbf{N}}{\prod}     \textbf{P}^{\xi^{\mathrm{Sloped}}}_{\mathscr{O}}  \big[    \mathscr{F}_1 \underset{\mathcal{I}_j (  \mathscr{O} )}{{\overset{h \geq ck}{\longleftrightarrow }}}    \mathscr{F}_2                  \big]  \bigg]   \bigg] \text{ } \\  \geq  g_c (     \beta                  )       -    g_c ( 0 )    
\end{align*}

\noindent hence concluding the argument. \boxed{}

\bigskip

\noindent We turn to the \textit{Proof} for the first $\textbf{Lemma}$, $\textbf{Lemma}$ \textit{4.1}, provided in the section, which is primarily reliant upon exhibiting that the supremum of the crossing probability of $\mathcal{B}(L)$ occurring, under $\textbf{P}_{\mathscr{O}} \big[ \cdot \big]$, upper bounds the crossing probability of the faces $\mathscr{F}_1$ and $\mathscr{F}_2$ being connected within the cylindrical interior $\mathcal{I} \big(   \mathscr{O}_{N,M}   \big)$.

\bigskip

\noindent \textit{Proof of Lemma 4.1}. Fix $j \geq 0$ and $n,k \geq 1$. We begin by defining crossing events from root faces on $\mathcal{B}(\mathscr{O})$. Specifically, we require that there exist a strictly positive number, $2j+2$, of $\mathrm{x}$-crossings in $\mathcal{I}(\mathscr{O})$, denoted with $\gamma_i$. Along each such $\gamma_i$, the value of the height function is dependent upon whether the path traverses an odd or even face; if the path begins on an even face on $\mathcal{B}\big(  \mathscr{O} \big)$, then $\gamma_i$, the height function $h \equiv \mathrm{Even}$, ie the image of the height function on every even face is an even integer, and otherwise, $h \equiv \mathrm{Odd}$ if the path traverses an odd face. These conventions for the image of difference faces of the height function over the cylinder is adopted from {\color{blue}[11]}.

\bigskip

\noindent Furthermore, to distinguish between crossings across any of the $2j+2$ paths, introduce the intersection of the interior of the cylinder with the line centered for some $j>0$ such that $\mathcal{L}_j \cap \mathcal{B}\big(  \mathscr{O}  \big) \neq \emptyset$, namely,

\begin{align*}
    \mathcal{I}_j \big(   \mathscr{O}  \big) \equiv  \mathcal{I} \big(   \mathscr{O}     \big)    \cap      \mathcal{L}_j       \text{ }   \text{ . } 
\end{align*}

\bigskip

\noindent From the result of $\textbf{Lemma}$ \textit{4.2}, we exhibit that the probability of $\mathscr{F}_1$ and $\mathscr{F}_2$ being connected with height at least $\geq (n-j)k$, within $\mathcal{I}_j\big(     \mathscr{O}   \big)$, or within $\mathcal{I}_{j+1}\big(    \mathscr{O}   \big)$, for some $j < j+1$, satisfies,

\begin{align*}
     \textbf{P}^{\xi^{\mathrm{Sloped}}}_{\mathscr{O}} \big[         \mathcal{V}^{\text{ } \mathrm{x}}_{h \geq (n-j)k} \text{ } \big(  \mathscr{O}^{\prime}      \big) \big]  \leq      \textbf{P}^{\xi^{\mathrm{Sloped}}}_{\mathscr{O}}\big[          \mathcal{V}^{\text{ } \mathrm{x}}_{h \leq - (  n - (j+1) +1       ) k    } \big(  \mathscr{O}^{\prime} \big)     \big]       \text{ }  \text{ , } 
\end{align*}

\noindent in which the crossing probability $\mathcal{B}(L)$ dominates the probability of the crossing event $\big\{  \mathscr{F}_1 \underset{\mathcal{I}_j(\mathscr{O}_{N,M})}{\overset{ h \geq ck}{\longleftrightarrow} }     \mathscr{F}_2 \big\}$ occurring. To achieve such an upper bound as provided above, introduce the decomposition for the connectivity event between faces beginning on $\mathcal{B}(\mathscr{O})$, in terms of a summation of vertical $\mathrm{x}$-crossing events,

\begin{align*}
        \textbf{P}^{\xi^{\mathrm{Sloped}}}_{\mathscr{O}} \big[   \mathscr{F}_1 \underset{\mathcal{I}_j ( \mathscr{O}_{N,M})}{\overset{ h \geq ck}{\longleftrightarrow} }   \mathscr{F}_2       \big]  \equiv     \underset{\mathscr{D}_1 \sim \mathscr{D}_{\mathscr{O}} }{\underset{\xi^{\mathscr{D}_1}\text{ } \sim \text{ }  \textbf{B}\textbf{C}^{\mathrm{Sloped}}}{\underset{\big( \text{ }  \mathscr{D}_1        \text{ }  , \text{ }   \xi^{\mathscr{D}_1}      \text{ } \big)}{\sum}} }  \textbf{P}^{\xi^{\mathrm{Sloped}}}_{\mathscr{O}} \bigg[        \mathcal{V}^{\text{ } \mathrm{x}}_{h \geq (n-j)k} \text{ } \big( \mathscr{O}^{\prime}        \big)    \big|       h|_{\mathscr{D}_1}  =  \xi^{\mathscr{D}_1}  \bigg]     \textbf{P}^{\xi^{\mathrm{Sloped}}}_{\mathscr{O}} \big[     h|_{\mathscr{D}_1} = \xi^{\mathscr{D}_1}      \big]    \text{ }   \text{ , } 
\end{align*}

\noindent of the crossing probability across cylindrical domains where $\mathscr{F}_1$ is an even face, from the set of all such admissible domains $\mathscr{A}\mathscr{D}$, and $\mathscr{O}^{\prime} \subset \mathscr{O}$, while for the closely related crossing probability, introduce the decomposition,

\begin{align*}
        \textbf{P}^{\xi^{\mathrm{Sloped}}}_{\mathscr{O}} \big[  \mathscr{F}_1 \underset{\mathcal{I}_{j+1} ( \mathscr{O}_{N,M})}{\overset{ h \geq ck}{\longleftrightarrow} }      \mathscr{F}_2     \big] \equiv     \underset{\mathscr{D}_1 \sim \mathscr{D}_{\mathscr{O}} }{\underset{\xi^{\mathscr{D}_1}\text{ } \sim \text{ }  \textbf{B}\textbf{C}^{\mathrm{Sloped}}}{\underset{\big( \text{ }  \mathscr{D}_1        \text{ }  , \text{ }   \xi^{\mathscr{D}_1}      \text{ } \big)}{\sum}} }  \textbf{P}^{\xi^{\mathrm{Sloped}}}_{\mathscr{O}} \bigg[        \mathcal{V}^{\text{ } \mathrm{x}}_{h \leq - (  n - (j+1) +1       ) k    } \text{ } \big(   \mathscr{O}^{\prime}     \big)   \big|       h|_{\mathscr{D}_1}  = \xi^{\mathscr{D}_1}       \bigg]    \textbf{P}^{\xi^{\mathrm{Sloped}}}_{\mathscr{O}} \big[    h|_{\mathscr{D}_1}  =  \xi^{\mathscr{D}_1}      \big]                 \text{ }    \text{ , } 
\end{align*}

\noindent for connectivity beginning on an odd face $\mathscr{F}_1$, also for $\mathscr{O}^{\prime} \subset \mathscr{O}$. We proceed by observing that, the conditional probability on the value of the height function achieved over $\mathscr{D}_1$, from previous applications of $\mathrm{(SMP)}$, yields,

\begin{align*}
  \textbf{P}^{\xi^{\mathrm{Sloped}}}_{\mathscr{O}} \bigg[         \mathcal{V}^{\text{ } \mathrm{x}}_{h \geq (n-j)k} \text{ } \big(     \mathscr{O}^{\prime}      \big)  \big|     h|_{\mathscr{D}_1}  \text{ } = \text{ } \xi^{\mathrm{Sloped}}         \bigg]  \text{ } \overset{\mathrm{(SMP)}}{\equiv} \text{ }      \textbf{P}^{\xi^{\mathscr{D}_1}}_{\mathscr{O}} \big[     \mathcal{V}^{\text{ } \mathrm{x}}_{h \geq (n-j)k} \text{ } \big(  \mathscr{O}^{\prime}      \big)        \big]  \text{ }   \\ \text{ } \underset{(\xi^{\mathscr{D}_1})^{\prime} \geq \xi^{\mathscr{D}_1}}{\overset{\mathrm{(CBC)}}{\leq}}      \textbf{P}^{(\xi^{\mathscr{D}_1})^{\prime}}_{\mathscr{O}} \big[      \mathcal{V}^{\text{ } \mathrm{x}}_{h \geq (n-j)k} \text{ } \big(    \mathscr{O}^{\prime}        \big)             \big]   \text{ } \\ \text{ }  \equiv         \textbf{P}^{(\xi^{\mathscr{D}_1})^{\prime}}_{\mathscr{O}} \big[     \mathcal{V}^{\text{ } \mathrm{x}}_{h \leq - (  n - (j+1) +1       ) k    }  \big]    \\ \text{ } \overset{\mathrm{(SMP)}}{\equiv}     \textbf{P}^{\xi^{\mathrm{Sloped}}}_{\mathscr{O}}\bigg[           \mathcal{V}^{\text{ } \mathrm{x}}_{h \leq - (  n - (j+1) +1       ) k    }    \text{ }  \big| \text{ }   h|_{\mathscr{D}_1} = \xi^{\mathscr{D}_1}      \bigg]  \text{ , } 
\end{align*}

\noindent where in the final expression in the sequence of inequalities above, we incorporate conditioning on the value of the height function, so that $h$ has boundary conditions $\xi^{\mathscr{D}_1}$ when restricted to $\mathscr{D}_1$, while the base probability measure $\textbf{P}^{\xi}_{\mathscr{O}}[ \text{ } \cdot \text{ } ]$ has boundary conditions $\xi \equiv \xi^{\mathrm{Sloped}} \sim \textbf{B}\textbf{C}^{\mathrm{Sloped}}$, hence implying that the conditionally defined vertical crossing probability across $\mathscr{O}^{\prime}$ can be upper bounded with,

\begin{align*}
      \textbf{P}^{\xi^{\mathrm{Sloped}}}_{\mathscr{O}} \bigg[         \mathcal{V}^{\text{ } \mathrm{x}}_{h \geq (n-j)k} \text{ } \big(   \mathscr{O}^{\prime}     \big) \text{ }   \big|  \text{ }    h|_{\mathscr{D}_1}  \text{ } = \text{ } \xi^{\mathrm{Sloped}}         \bigg]   \leq    \textbf{P}^{\xi^{\mathrm{Sloped}}}_{\mathscr{O}}\bigg[         \mathcal{V}^{\text{ } \mathrm{x}}_{h \leq - (  n - (j+1) +1       ) k    }    \text{ }  \big| \text{ }   h|_{\mathscr{D}_1} = \xi^{\mathscr{D}_1}      \bigg]        \text{ }     \text{ , } 
\end{align*}

\noindent from which further rearrangements also imply that a similar bound holds, independently of conditioning on the value of the height function appearing in the upper and lower bounds of the previous inequality, by first multiplying the inequality with the unconditioned probability on the value of the height function,

\begin{align*}
  \text{ }    \textbf{P}^{\xi^{\mathrm{Sloped}}}_{\mathscr{O}} \bigg[  \mathcal{V}^{\text{ } \mathrm{x}}_{h \geq (n-j)k} \text{ } \big(    \mathscr{O}^{\prime}     \big)  \text{ }  \big| \text{ }      h|_{\mathscr{D}_1}  \text{ } = \text{ } \xi^{\mathrm{Sloped}}           \bigg]     \textbf{P}_{\mathscr{O}}   \big[ h|_{\mathscr{D}_1} = \xi^{\mathscr{D}_1}  \big]    \leq      \textbf{P}^{\xi^{\mathrm{Sloped}}}_{\mathscr{O}} \bigg[         \mathcal{V}^{\text{ } \mathrm{x}}_{h \leq - (  n - (j+1) +1       ) k    }    \text{ }  \big| \text{ }   h|_{\mathscr{D}_1} = \xi^{\mathscr{D}_1}        \bigg]     \\ \times       \textbf{P}_{\mathscr{O}} \text{  } \big[    h|_{\mathscr{D}_1} = \xi^{\mathscr{D}_1}     \big]    \text{ }   \text{ , } 
\end{align*}

\noindent and from this expression, taking the summation over all realizations of the domains and boundary conditions over which the crossing probability occurs,

\begin{align*}
 \underset{\mathscr{D}_1 \sim \mathscr{D}_{\mathscr{O}} }{\underset{\xi^{\mathscr{D}_1}\text{ } \sim \text{ }  \textbf{B}\textbf{C}^{\mathrm{Sloped}}}{\underset{\big( \text{ }  \mathscr{D}_1        \text{ }  , \text{ }   \xi^{\mathscr{D}_1}      \text{ } \big)}{\sum}} }      \text{ }    \textbf{P}^{\xi^{\mathrm{Sloped}}}_{\mathscr{O}} \bigg[    \mathcal{V}^{\text{ } \mathrm{x}}_{h \geq (n-j)k} \text{ } \big(      \mathscr{O}^{\prime}        \big)  \text{ }  \big| \text{ }      h|_{\mathscr{D}_1}   =  \xi^{\mathrm{Sloped}}            \bigg]         \textbf{P}_{\mathscr{O}}    \big[   h|_{\mathscr{D}_1} = \xi^{\mathscr{D}_1}      \big]  \text{ }   \\     \leq   \underset{\mathscr{D}_1 \sim \mathscr{D}_{\mathscr{O}} }{\underset{\xi^{\mathscr{D}_1}\text{ } \sim \text{ }  \textbf{B}\textbf{C}^{\mathrm{Sloped}}}{\underset{\big( \text{ }  \mathscr{D}_1        \text{ }  , \text{ }   \xi^{\mathscr{D}_1}      \text{ } \big)}{\sum}} }  \textbf{P}^{\xi^{\mathrm{Sloped}}}_{\mathscr{O}} \bigg[          \mathcal{V}^{\text{ } \mathrm{x}}_{h \leq - (  n - (j+1) +1       ) k    }    \text{ }  \big| \text{ }   h|_{\mathscr{D}_1} = \xi^{\mathscr{D}_1}    \bigg]             \textbf{P}_{\mathscr{O}}  \big[     h|_{\mathscr{D}_1} = \xi^{\mathscr{D}_1}      \big]             \text{ , } 
\end{align*}

\noindent from which the final inequality above is equivalent to,

\begin{align*}
     \textbf{P}^{\xi^{\mathrm{Sloped}}}_{\mathscr{O}} \big[    \mathscr{F}_1 \underset{\mathcal{I}_{j}( \mathscr{O}_{N,M})}{\overset{ h \geq ck}{\longleftrightarrow} }      \mathscr{F}_2        \big]  \text{ } \leq \text{ }  \textbf{P}^{\xi^{\mathrm{Sloped}}}_{\mathscr{O}} \big[    \mathscr{F}_1 \underset{\mathcal{I}_{j+1}(\mathscr{O}_{N,M} )}{\overset{ h \geq ck}{\longleftrightarrow} }      \mathscr{F}_2      \big]   \text{ . } 
\end{align*}

\noindent To proceed, fix some root face on $\mathcal{B}\big(  \mathscr{O}  \big)$, from which the remainder of faces along each $\gamma_i$, either within $\mathcal{I}_j\big( \mathscr{O}  \big)$, or within $\mathcal{I}_{j+1}\big(  \mathscr{O}  \big)$, is quantified with,

\begin{align*}
        \textbf{P}^{\xi^{\mathrm{Sloped}}}_{\mathscr{O}} \big[     \mathcal{B} \big(   \lceil         L N     \rceil    \big)   \big]      \text{ }      \text{ , } 
\end{align*}

\noindent for $j \equiv 0$, $n\text{ } \equiv \text{ } \lceil\frac{\lfloor L N \rfloor}{k}\rceil$, such that $  \mathcal{B} ( \text{ }   \lceil \text{ }        L \text{ } N    \text{ }  \rceil   \text{ }  )  \subset \mathcal{B} ( \text{ }     n k       \text{ } )$. To iteratively applying the same argument for any $j >0$ begin with the following summation,

\begin{align*}
\text{ } N^{-1} \text{ }  \sum_{j \geq 0}    \text{ }   \textbf{P}^{\xi^{\mathrm{Sloped}}}_{\mathscr{O}} \big[     \mathscr{F}_1 \underset{\mathcal{I}_j ( \mathscr{O}_{N,M} )}{\overset{ h \geq ck}{\longleftrightarrow} }   \mathscr{F}_2     \big] \text{ }                   \text{ , } \tag{\textit{Prob Sum}}
\end{align*}

\noindent normalized in $N^{-1}$, can be upper bounded as follows,

\begin{align*}
 (\textit{Prob Sum})\leq    \textbf{P}^{\xi^{\mathrm{Sloped}}}_{\mathscr{O}} \big[      \mathscr{F}_1 \underset{\mathcal{I}_0( \mathscr{O}_{N,M})}{\overset{ h \geq ck}{\longleftrightarrow} }   \mathscr{F}_2  \big] \text{ }  \\ \overset{(*)}{\leq}  \textbf{P}^{\xi^{\mathrm{Sloped}}}_{\mathscr{O}} \big[    \mathscr{F}_1 \underset{\mathcal{I}_{j-1}( \mathscr{O}_{N,M})}{\overset{ h \geq ck}{\longleftrightarrow} }   \mathscr{F}_2     \big]  \leq   2^N  \underset{\mathscr{F}_2 \in \mathcal{T}( \mathscr{O} ) }{\underset{\mathscr{F}_1 \in \mathcal{B}( \mathscr{O})}{\mathrm{sup}}} \text{ } \textbf{P}^{\xi^{\mathrm{Sloped}}}_{\mathscr{O}} \big[      \mathscr{F}_1 \underset{\mathcal{I}( \mathscr{O}_{N,M})}{\overset{ h \geq ck}{\longleftrightarrow} }   \mathscr{F}_2      \big]  \text{ } \text{ , } \end{align*}
 
 \noindent where the final expression obtained above is equivalent to the probability, of the intersection over $j$,
 
 \begin{align*}
 2^N  \underset{\mathscr{F}_2 \in \mathcal{T}( \mathscr{O} ) }{\underset{\mathscr{F}_1 \in \mathcal{B}( \mathscr{O})}{\mathrm{sup}}}      \textbf{P}^{\xi^{\mathrm{Sloped}}}_{\mathscr{O}} \big[     \underset{\mathrm{countably\text{ }  many} \text{ } j}{\bigcup}  \big\{  \mathscr{F}_1 \underset{\mathcal{I}_j ( \mathscr{O}_{N,M})}{\overset{ h \geq ck}{\longleftrightarrow} }   \mathscr{F}_2     \big\}    \big]      \text{ }     \text{ , } 
\end{align*}

\noindent where in the sequence of inequalities above, $(*)$ holds because,

\begin{align*}
 \big\{    \mathscr{F}_1 \underset{\mathcal{I}_0( \mathscr{O}_{N,M})}{\overset{ h \geq ck}{\longleftrightarrow} }   \mathscr{F}_2     \big\}  \subsetneq  \big\{        \mathscr{F}_1 \underset{\mathcal{I}_{j-1}( \mathscr{O}_{N,M})}{\overset{ h \geq ck}{\longleftrightarrow} }   \mathscr{F}_2     \big\} \text{ } \text{ , } 
\end{align*}

\noindent In the final term of the sequence of inequalities above, observe the containment of vertical $\mathrm{x}$-crossings across the cylinder, with the last face of the path, $\mathscr{F}_2$, along $\mathcal{T} \big(  \mathscr{O} \big)$ held fixed,

\begin{align*}
    \big\{  \mathscr{F}_1   \underset{\mathcal{I} (  \mathscr{O} )}{\overset{h \geq ck}{\longleftrightarrow}}   \mathscr{F}_2   \big\}     \subsetneq    \underset{\mathscr{F}_i \in \mathcal{B}(  \mathscr{O} )}{\bigcup}   \big\{   \mathscr{F}_i  \underset{\mathcal{I} (  \mathscr{O} )}{\overset{h \geq ck}{\longleftrightarrow}}   \mathscr{F}_2   \big\}      \text{ }   \text{ , } 
\end{align*}

\noindent and also the fact that $\mathcal{I}_j \big(  \mathscr{O}  \big) \text{ } \subset \text{ } \mathcal{I} \big(  \mathscr{O} \big)$. Finally, the prefactor to the probability in the upper bound also results from the same reasoning in the arguments for the proof of \textbf{Proposition} \textit{4.1} in {\color{blue}[11]}, namely that,

\begin{align*}
    { N \choose 2n}    \text{ }  \leq \text{ } 2^N \text{ } \text{ , }
\end{align*}

\noindent resulting from the total number of ways of choosing $2n$ faces from a total of $N$ faces on $\mathcal{B} \big(  \mathscr{O}  \big)$. Finally, from the previous result, the free-energy lower bound hence transfers to the desired probability, from the fact that,

\begin{align*}
        \big( \text{ }   N M     \text{ } \big)^{-1} \text{ }  \mathrm{log} \text{ } \bigg[ \text{ }    \underset{\mathscr{F}_2 \in \mathcal{T} (  \mathscr{O}_{N,M})}{\underset{\mathscr{F}_1 \in \mathcal{S}_{\mathscr{F}}}{\mathrm{sup}}} \text{ } \textbf{P}^{\xi^{\mathrm{Sloped}}}_{\mathscr{O}}  \big[    \mathscr{F}_1 \underset{\mathcal{I}(\mathscr{O}_{N,M} )}{\overset{ h \geq ck}{\longleftrightarrow} }     \text{ } \mathscr{F}_2  \big]       \bigg]   \geq     \big( \text{ }   N M     \text{ } \big)^{-1}   \mathrm{log}  \big\{ \textbf{P}^{\xi^{\mathrm{Sloped}}}_{\mathscr{O}}  \big[        \mathcal{B}(L)
          \big]       \big\}    \text{ }   \text{ , } 
\end{align*}

\noindent as the limit infimum for $N$ even, and then $M$, each of which are taken $\longrightarrow + \infty$, respectively, yielding the desired lower bound dependent upon the free energy,

\begin{align*}
  \text{ }    \underset{ N \longrightarrow + \infty}{\underset{N\text{ }  \mathrm{even}}{\mathrm{lim\text{ }  inf}}} \text{ }   \underset{ M \longrightarrow + \infty}{\mathrm{lim\text{ } inf}}  \big( \text{ } N M  \text{ } \big)^{-1} \text{ }      \mathrm{log} \bigg[     \underset{\mathscr{F}_2 \in \mathcal{T}(   \mathscr{O}_{N,M})}{\underset{\mathscr{F}_1 \in \mathcal{S}_{\mathscr{F}}}{\mathrm{sup}}}  \textbf{P}^{\xi^{\mathrm{Sloped}}}_{\mathscr{O}}  \big[  \mathscr{F}_1 \underset{\mathcal{I}( \mathscr{O}_{N,M})}{\overset{ h \geq ck}{\longleftrightarrow} }   \mathscr{F}_2 \big]        \bigg] \text{ } \\ =               \underset{ N \longrightarrow + \infty}{\underset{N\text{ }  \mathrm{even}}{\mathrm{lim\text{ }  inf}}} \text{ }   \underset{ M \longrightarrow + \infty}{\mathrm{lim\text{ } inf}}  \big( \text{ } N M  \text{ } \big)^{-1}      \mathrm{log} \bigg[    \underset{\mathscr{F}_2 \in \mathcal{T}(   \mathscr{O}_{N,M})}{\underset{\mathscr{F}_1 \in \mathcal{S}_{\mathscr{F}}}{\mathrm{sup}}}  \textbf{P}^{\xi^{\mathrm{Sloped}}}_{\mathscr{O}}  \big[    \underset{\mathrm{countably\text{ }  many} \text{ } j: j \in \textbf{N}}{\bigcup}  \big\{   \mathscr{F}_1 \underset{\mathcal{I}_j( \mathscr{O}_{N,M})}{\overset{ h \geq ck}{\longleftrightarrow} }      \mathscr{F}_2  \big\} \big]      \bigg]            \text{ }  \text{ , } 
\end{align*}

\noindent which, after applying the previous result which gave a free energy lower bound, yields the identical lower bound from the sloped free energy function,

\begin{align*}
               g_c (     \beta                  )   -     g_c ( 0 )       \text{ } \text{ , } 
\end{align*}

\noindent hence concluding the proof. \boxed{}

       \begin{figure}
    \centering
    \includegraphics[width=1.04\columnwidth]{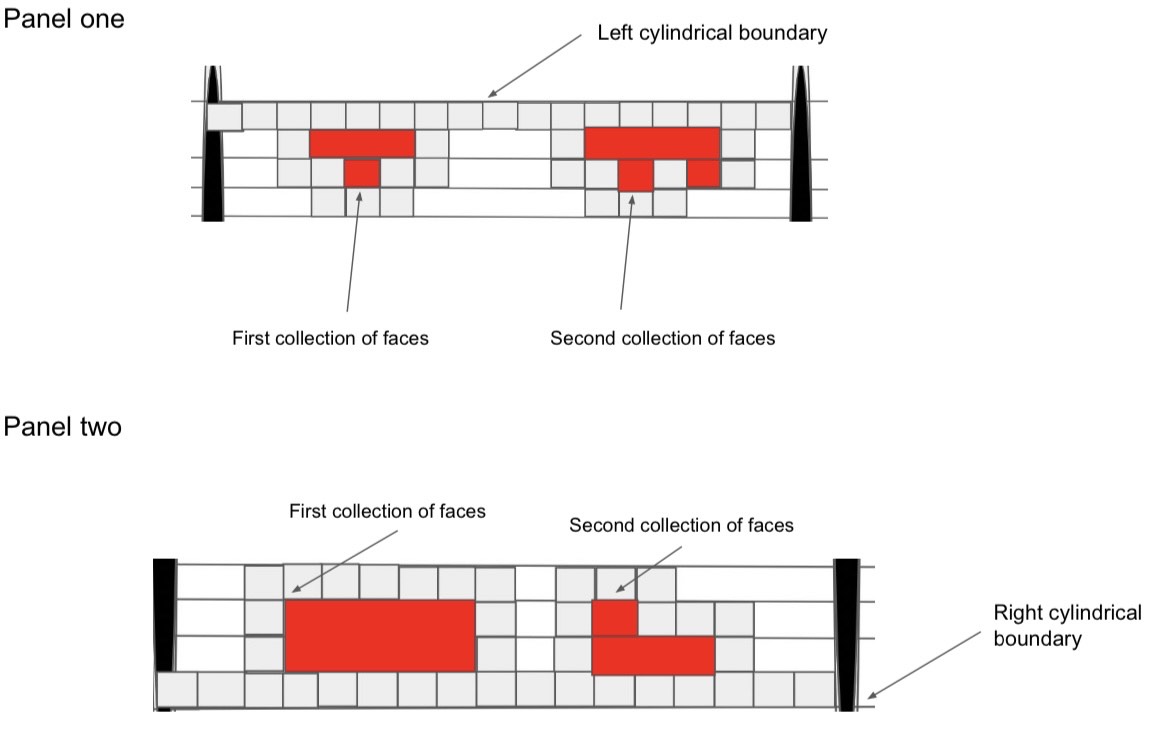}
    \caption{\textit{A depiction of the faces incident to the left and right boundaries of the cylinder, after a } $\frac{\pi}{2}$ \textit{ rotation, from Panel one, and Panel two, respectively.}}
    \end{figure}

\bigskip

\noindent In the next section, criteria are imposed on constants appearing along with the free energy.

\subsection{Towards obtaining the power of the exponential lower bound for the annulus crossing probability}

\noindent Below, we present the theorem for the free energy exponent, which is an analogue to \textbf{Theorem} \textit{1.5}, of {\color{blue}[11]} that was introduced in \textit{1.3}. The following lower bound estimate is dependent upon the free energy function $f^{\prime}_c$ for sloped boundary conditions, from \textbf{Definition} \textit{10}.

\bigskip

\noindent \textbf{Theorem} $6V$ \textit{2} (\textit{sloped free energy perturbation away from }$0$). Given the sloped free energy function provided in \textbf{Definition} \textit{10}, and parameters in \textbf{Theorem} \textit{6V} \textit{1}, the crossing probability across the annulus satisfies,

\begin{align*}
             \textbf{P}^{\xi^{\mathrm{Sloped}}} \big[        \mathcal{O}_{h \geq ck}  ( 6n , n_R n)  \big]     \geq      c \text{ }  \mathrm{exp} \big[       C^{\prime}  \big(  r^{\prime} \big)^2    \big(           g_c (     \beta                  )   -  g_c ( 0 )    \big)               \big]                                \text{ } \text{ , } 
\end{align*}

\noindent in which the lower bound depends upon strictly positive constants, where,

\begin{align*}
  C^{\prime}  \equiv      M N     \text{ } \text{ , } 
\end{align*}

\noindent is a prefactor from the number of preimages of a directed, oriented loop configuration, given from the pullback under $\big( \text{ } \mathcal{T}^{\prime}\text{ } \big)^{-1}$, from (\textit{Property Two}) of $\big( \text{ } \mathcal{T}^{\prime}\text{ } \big)^{-1}$,

\begin{align*}
      g_c (     \beta                  )     -   g_c ( 0 ) \text{ } \text{ , } 
\end{align*}

\noindent is the perturbation of the sloped free energy function away from $0$, and finally,

\begin{align*}
       \beta \equiv     \frac{k}{\eta     r^{\prime}    }        \text{ } \text{ , } 
\end{align*}

\noindent is a strictly positive parameter at which the sloped free energy function, $g_c$, is evaluated away from $0$.

\bigskip

\noindent \textit{Proof of Theorem }. We incorporate previous items used in the proofs of $\textbf{Lemma}$ \textit{4.1} and $\textbf{Lemma}$ \textit{4.2}, namely, the decomposition of the crossing event across $\mathcal{I}\big(  \mathscr{O}  \big)$,

\begin{align*}
        \textbf{P}^{\xi^{\mathrm{Sloped}}}_{\mathscr{O}} \big[   \mathscr{F}_1 \underset{\mathcal{I}( \mathscr{O}_{N,M} )}{\overset{ h \geq ck}{\longleftrightarrow} }   \mathscr{F}_2     \big]  \equiv     \underset{\mathscr{D}_1 \sim \mathscr{D}_{\mathscr{O}} }{\underset{\xi^{\mathscr{D}_1}\text{ } \sim \text{ }  \textbf{B}\textbf{C}^{\mathrm{Sloped}}}{\underset{\big( \text{ }  \mathscr{D}_1        \text{ }  , \text{ }   \xi^{\mathscr{D}_1}      \text{ } \big)}{\sum}} }  \textbf{P}^{\xi^{\mathrm{Sloped}}}_{\mathscr{O}} \big[       \mathcal{V}^{\text{ } \mathrm{x}}_{h \geq (n-j)k} \text{ } \big(     \mathscr{O}^{\prime}      \big)  \big|       h|_{\mathscr{D}_1}   =  \xi^{\mathscr{D}_1}      \big] \text{ }    \textbf{P}^{\xi^{\mathrm{Sloped}}}_{\mathscr{O}} \big[     h|_{\mathscr{D}_1} \text{ } = \text{ } \xi^{\mathscr{D}_1}      \big]    \text{ }           \text{ , } 
\end{align*}

\noindent and a suitable realization of a domain within the cylinder, with $\mathscr{D}_1$. Next, from such a $\mathscr{D}_1$, to show that the conditional probability, supported over $\mathscr{O}$, satisfies,

\begin{align*}
  \textbf{P}^{\xi^{\mathrm{Sloped}}}_{\mathscr{O}} \big[        \mathscr{F}_1 \underset{\mathcal{I}( \mathscr{O}_{N,M} )}{\overset{ h \geq ck}{\longleftrightarrow} }   \mathscr{F}_2            \big|       \big(  \mathscr{D}_1 \big)^{\prime} \equiv \mathscr{D}  \big]  \geq  \textbf{P}^{\xi^{\mathrm{Sloped}}}_{\mathscr{O}} \big[        \mathscr{F}_1 \underset{\mathcal{I}(\mathscr{O}_{N,M})}{\overset{ h \geq ck}{\longleftrightarrow} }   \mathscr{F}_2        \big]   \text{ . } 
\end{align*}

\noindent To this end, for such a realization $\big( \mathscr{D}_1 \big)^{\prime}$, it suffices to argue that, for the odd faces along $\mathcal{B} \big(  \mathscr{O} \big)$, the collection of events,

\begin{align*}
      \text{ }        \big\{      \textbf{P}^{\xi^{\mathrm{Sloped}}}_{\mathscr{O}} \big[      \mathscr{F}_1 \underset{\mathcal{I}_{j} ( \mathscr{O}_{N,M} ) \cap  (  \mathscr{D}_1)^{\prime}}{\overset{ h \geq ck}{\longleftrightarrow} }   \mathscr{F}_2       \big]       \big\}_{j \text{ } \mathrm{odd}, \text{ } 1 \leq j \leq 2n-2 }                          \text{ }       \text{ , } \tag{\textit{Odd Face Connectivity}}
\end{align*}

\noindent connecting $\mathscr{F}_1$ and $\mathscr{F}_2$, of height $\geq ck$, within $\mathcal{I}_j \big( \mathscr{O}  \big)  \cap \big( \mathscr{D}_1 \big)^{\prime}$, establish a partition of $\mathscr{O}$ into $n$ volumes $\widetilde{\mathscr{V}_1} , \cdots , \widetilde{\mathscr{V}_n}$. Fix some index $i^{\prime}$ for each $\widetilde{\mathscr{V}_{i^{\prime}}}$. From the left and right crossings required for each $\widetilde{\mathscr{V}_{i^{\prime}}}$ to exist, the four boundaries enclosing each finite volume yield the decomposition $\widetilde{\mathscr{V}_{i^{\prime}}}  \equiv (\widetilde{\gamma_{i^{\prime}}})_L   \cup  (\widetilde{\gamma_{i^{\prime}}})_R    \cup   \big(  \mathcal{B}\big(  \mathscr{O}  \big) \big)^{\prime}   \cup        \big(  \mathcal{T}\big( \mathscr{O} \big)   \big)^{\prime}      \text{ } $, where $  \big(  \mathcal{B}\big(  \mathscr{O}  \big)   \big)^{\prime} \subset \mathcal{B} \big( \text{ } \mathscr{O} \text{ } \big)$, and $  \big(  \mathcal{T}\big(  \mathscr{O}  \big)   \big)^{\prime}  \subset \mathcal{T}\big( \mathscr{O}  \big)$. Besides the left and right boundaries of the cylindrical domain indexed in $i$, $ \big(  \mathcal{B}\big( \mathscr{O}  \big) \big)^{\prime} $, and $ \big(  \mathcal{T}\big(  \mathscr{O}  \big)   \big)^{\prime}$, are respectively defined with, 

\begin{align*}
  \big(  \mathcal{B}\big(  \mathscr{O}  \big)  \big)^{\prime}  \text{ } \equiv   \text{ }    \mathcal{B} \big(  \mathscr{O}  \big)  \text{ }   \cap \text{ }       \big( \mathcal{L}_j \cup \mathcal{L}_{j+1}  \big)               \text{ , }     \\  \big(  \mathcal{T}\big( \mathscr{O} \big)   \big)^{\prime} \text{ } \equiv   \text{ }  \mathcal{T}\big(  \mathscr{O}  \big) \text{ } \cap \text{ }        \big(  \mathcal{L}_j \cup \mathcal{L}_{j+1}  \big)               \text{ , } 
\end{align*}

\noindent where $\mathcal{L}_j \cup \mathcal{L}_{j+1}$ denote the finite subvolume spanned within the cylinder by the lines $\mathcal{L}_j$ and by $\mathcal{L}_{j+1}$.

\bigskip

\noindent Within each $\widetilde{\mathscr{V}_{i^{\prime}}}$, with positive probability the crossing restricted across each such domain indexed in $i^{\prime}$,

\begin{align*}
  \textbf{P}^{\xi^{\mathrm{Sloped}}}_{\mathscr{O}} \big[    \mathscr{F}_1    \underset{\widetilde{\mathscr{V}_{i^{\prime}}}}{\overset{h \geq ck}{\longleftrightarrow}}          \mathscr{F}_2 \big] \text{ }  \text{ , } 
\end{align*}

\noindent occurs,\footnote{With some abuse of notation, instead of requiring that $\mathscr{F}_1 \in \mathcal{B} \big( \text{ } \mathscr{O} \text{ } \big)$, and that $\mathscr{F}_2 \in \mathcal{T} \big( \text{ } \mathscr{O} \text{ } \big)$, for the collection of connectivity events for $j$ odd, $\mathscr{F}_1 \in \mathcal{B} \big( \text{ } \mathscr{O} \text{ } \big) \text{ } \cap \text{ } \big( \text{ } \mathcal{L}_j  \cup   \mathcal{L}_{j+1} \text{ } \big) $, and $\mathscr{F}_2 \in \mathcal{T} \big( \text{ }  \mathscr{O}    \text{ } \big) \text{ } \cap \text{ } \big( \text{ }  \mathcal{L}_j  \cup   \mathcal{L}_{j+1}  \text{ } \big)$, for $j$ odd. That is, the intersection of the bottom of the cylinder $\mathscr{O}$ with the union of two adjacent lines which have nonempty intersection with the top and bottom.} as the connectivity event within faces formed by the collection of events along odd faces, from ($\textit{Odd Face Connectivity}$), establishes left and right boundaries for each $\widetilde{\mathscr{V}_{i^{\prime}}}$, and hence also for each $\mathscr{V}_{i^{\prime}}$. In particular,

\begin{align*}
   \text{ }      \textbf{P}^{\xi^{\mathrm{Sloped}}}_{\mathscr{O}} \big[          \mathscr{F}_1    \underset{\widetilde{\mathscr{V}_{i^{\prime}}}}{\overset{h \geq ck}{\longleftrightarrow}}          \mathscr{F}_2   \big]   \geq    \textbf{P}^{\xi^{\mathrm{Sloped}}}_{\mathscr{O}} \big[  \mathscr{F}_1    \underset{\mathscr{V}_{i^{\prime}}}{\overset{h \geq ck}{\longleftrightarrow}}          \mathscr{F}_2  \big]    \text{ }   \text{ . } 
\end{align*}

\noindent because there can exist countably many $\mathscr{V}_{i^{\prime}}$ that are contained within each $\widetilde{\mathscr{V}_{i^{\prime}}}$.

            \begin{figure}
\begin{align*}
\includegraphics[width=0.94\columnwidth]{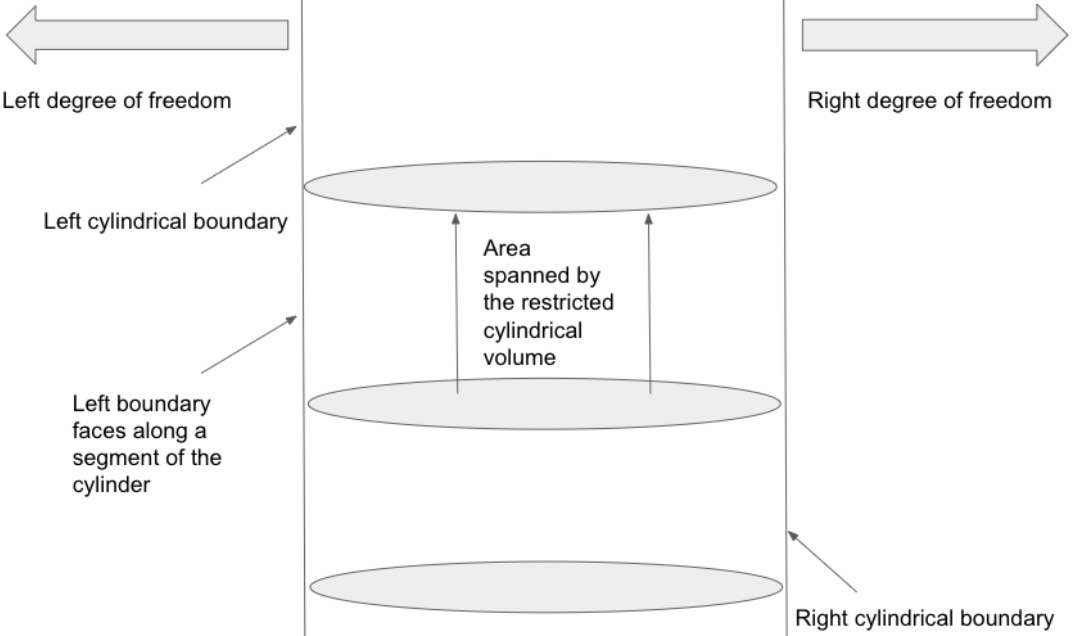}\\
\end{align*}
\caption{\textit{Quantification of the six-vertex free energy on the cylinder}. From strip estimates for crossing probabilities and the corresponding domains that are constructed, sloped boundary conditions are carried over into the cylinder by enforcing values of the height function on the left and right boundaries of the cylinder. Within each isolated segment of the cylinder depending upon the cardinality of the volume in the weak volume limit, crossing estimates are derived in restrictions to smaller cylindrical volumes intersecting each segment.}
\end{figure}

\begin{figure}
    \centering
    \includegraphics[width=0.97\columnwidth]{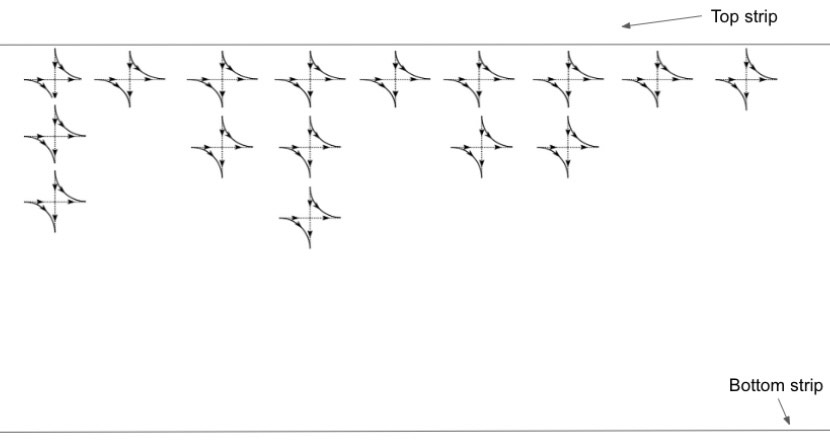}
    \caption{\textit{A fully packed loop configuration in the strip corresponding to previous configurations used in arguments throughout Section 3 for strip crossing probabilities before estimations of the cylindrical free-energy depending upon the arrow imbalance.} Surrounding the top boundary of the strip, selected faces of the six-vertex height function are depicted in the full packing limit, in which one obtains directed loops oriented towards the bottom boundary of the strip.}
\end{figure}

\bigskip

\noindent In comparison to the distribution of domains $\widetilde{\mathscr{V}_{i^{\prime}}}$ for sufficiently flat boundary conditions, from arguments that have been presented for sloped boundary conditions, with positive probability there can exist \textit{cylindrical freezing clusters}, which is the cylindrical counterpart of the \textit{freezing cluster} whose connectivity, and disconnectivity, properties were quantified throughout the previous section in the setting of the strip. 

\bigskip

\noindent For such cylindrical domains within each $\mathscr{V}_{i^{\prime}}$, with positive probability there exists two paths, one to the left and another to the right, such that, for $\widetilde{\mathscr{V}_{i^{\prime}}} \subsetneq \mathscr{V}_{i^{\prime}}$, $\widetilde{\mathscr{V}_{i^{\prime}}} \text{ } \cap \text{ } \mathscr{O} \text{ } \equiv \text{ } \widetilde{\mathscr{V}_{i^{\prime}}}$, and similarly, also that $\mathscr{V}_{i^{\prime}} \text{ } \cap \text{ } \mathscr{O} \text{ } \equiv \text{ } \mathscr{V}_{i^{\prime}}$.

\subsection{Introducing properties I-XIV that are satisfied by cylindrical symmetric domains}

\subsubsection{Finite volume objects}

\noindent To characterize the distribution of \textit{cylindrical} domains within a sequence of finite volumes $\mathscr{O}$, and to also obtain the power to which the exponent is raised in $Z^s$, consider the following. As $\mathscr{O} \longrightarrow \mathscr{O}_{+ \infty}$, partition each $\mathscr{O}$ into slices incremented in the $y$ position, namely $\mathrm{Slice}_y$ as introduced in {\color{blue}{[11]}} in arguments for the proof of the free energy exponent. Each $\mathrm{Slice}_y$ admits the decomposition $\mathrm{Slice}_y \text{ } \equiv \mathcal{T} \big( \text{ } \mathrm{Slice}_y \text{ } \big)  \text{ } \cup  \text{ } \mathcal{I} \big( \text{ } \mathrm{Slice}_y \text{ } \big) \text{ } \cup \text{ } \mathcal{B} \big( \text{ } \mathrm{Slice}_y \text{ } \big)$ into the top, interior, and bottom parts within $\mathscr{O}$. Additionally, introduce the notation for the union of Slices, in which over countably many $y$, set $\mathrm{Slice}_{y_0 \text{ } , \text{ } \cdots \text{ } ,\text{ } y_f} \text{ } \equiv \text{ } \underset{0 \leq i\leq f}{\underset{{\mathrm{countably \text{ } many }\text{ } i \text{ } : \text{ } \lfloor   N / i      \rfloor \in \textbf{Z}_{N-1}}}{\bigcup}} \mathrm{Slice}_{y_i}$. Necessarily, $\mathrm{Slice}_{y_j , y_{j+1} , y_{j+2}} \text{ } \cap \text{ }        \mathscr{O}              \text{ } \equiv \text{ } \mathscr{V}_{i^{\prime}-1} \text{ } \cup \text{ }    \mathscr{V}_{i^{\prime}}        \text{ } \cup   \text{ } \mathscr{V}_{i^{\prime}+1} \text{ }  \subset \mathscr{O} \text{ }$, for suitable $1 \leq j \leq j+1 \leq j+2 \leq N$. Similarly, for the same choice of $j,j+1,j+2$, there exists $\mathscr{O}^{\prime} \subset \mathscr{O}$ such that $\mathrm{Slice}_{y_j , y_{j+1} , y_{j+2}} \text{ } \cap \text{ }        \mathscr{O}^{\prime}              \text{ } \equiv \text{ } \widetilde{\mathscr{V}}_{i^{\prime}-1} \text{ } \cup \text{ }    \widetilde{\mathscr{V}}_{i^{\prime}}        \text{ } \cup   \text{ } \widetilde{\mathscr{V}}_{i^{\prime}+1} \text{ }  \subset \mathscr{O}^{\prime} \text{ }$.

\bigskip

With the following sequence of statements, we characterize the probability with which crossing events occur over the cylinder under sloped boundary conditions. In comparison to crossing probability estimates from the bottom to top of the cylinder which have previously been shown to occur with sufficiently good probability for sufficiently flat boundary conditions, {\color{blue}[11]}, we specify additional properties that cylindrical symmetric domains satisfy, ranging from: studying configurations for which bottom to top crossings occur with sufficiently good probability in spite of the fact that boundary conditions along the circumference of the bottom cylindrical face are frozen from boundary conditions of the height function; further characterizing configurations over the cylinder for which additional properties, with respect to tightness, and good, properties of such symmetric cylindrical domains are expected to hold; a violation of the modulo $6$ property of cylindrical symmetric domains which was originally formulated for cylindrical symmetric domains under sufficiently flat boundary conditions as provided in {\color{blue}[11]}; implications for encoding boundary conditions in several other models of Statistical Mechanics, ranging from the Ashkin-Teller, generalized random-cluster, and $\big(q_{\sigma} , q_{\tau} \big)$-spin models, which are considered after transfering results from crossing probabilities over the cylinder back to events over $\textbf{T}$.

\bigskip

More specifically, the environment of the random geometry can be further characterized with the following sequence of statements,

\begin{align*}
         \textbf{P}^{\xi^{\mathrm{Sloped}}}_{\mathscr{O}} \bigg[ \text{ }        \big| \mathscr{V}_{i^{\prime}}    \text{ }  \cap \text{ }  \mathrm{Slice}_y \big| \neq \emptyset                         \text{ }     \bigg]  \text{ }  \geq \text{ }  0        \text{ }   \text{ , } \text{ }\tag{\textit{I: Positive probability of intersection with} $\mathrm{Slice}_y$}
\end{align*}

\noindent which states that there is positive probability for one of the lines $\mathrm{Slice}_y$, of which the union $\underset{y}{\bigcup} \text{ } \mathrm{Slice}_y \supset \mathrm{Slice}_y$ establishes a partition of $\mathscr{O}$, have nonempty intersection with $\mathscr{V}_{i^{\prime}}$,

\begin{align*}
    \mathcal{I}_{\mathscr{V}_{i^{\prime}},y} \text{ } \equiv \text{ } \mathcal{I}_{i^{\prime},y} \text{ } \equiv \text{ }    \big\{     \forall \text{ } i^{\prime} \in \textbf{N}  \text{ } , \text{ } y \in \textbf{R} \text{ } ,\text{ }  \exists     \text{ }         \mathscr{F} \in F \big(  \mathscr{O} \big) \text{ } : \text{ } \mathscr{F} \text{ } \cap \text{ }         F \big( \mathscr{V}_{i^{\prime}}  \big)     \text{ } \neq \text{ } \emptyset     \big\}               \text{ }  \text{ , } \tag{\textit{II: Index set for domains in} $i^{\prime}$}
\end{align*}

\noindent which states that there is a natural ordering of the $\widetilde{\mathscr{V}_{i^{\prime}}}$ depending upon whether the connected components of each $\widetilde{\mathscr{V}_{i^{\prime}}}$ intersect $\mathrm{Slice}_y$ with the index set $\mathcal{I}_{i^{\prime},y}$. More explicitly, in terms of probabilities,

\begin{align*}
     1 \text{ } > \text{ }     \textbf{P}^{\xi^{\mathrm{Sloped}}}_{\mathscr{O}} \big[ \text{ }  \big| \mathcal{I}_{i^{\prime},y} \big| \text{ } \neq \emptyset   \text{ } \big] \text{ } > \text{ }    \textbf{P}^{\xi^{\mathrm{Sloped}}}_{\mathscr{O}} \big[ \text{ }    \big| \mathcal{I}_{i^{\prime\prime},y} \big| \equiv \emptyset  \text{ } \big] \text{ } \equiv \text{ } 0      \text{ } \text{ , }  \tag{\textit{II: Index set for domains in} $i^{\prime}$}
\end{align*}

 \noindent for $i^{\prime\prime} > i^{\prime}$, in which the connected components of the domain $\mathscr{V}_{i^{\prime}}$ are below the connected components of the domain $\mathscr{V}_{i^{\prime\prime}}$ in $\mathscr{O}$,  

\begin{align*}
           \textbf{P}^{\xi^{\mathrm{Sloped}}}_{\mathscr{O}} \big[     y \text{ }\equiv 2 \text{ } \mathrm{mod}\text{ }  3  :      \text{ }    \big| \text{ }  \widetilde{\mathscr{V}_{i^{\prime}}} \text{ }  \cap \text{ } \mathrm{Slice}_{y-1,y,y+1}         \text{ } \big| \neq \emptyset             \text{ }     \big]   \text{ } \geq \text{ }  0              \text{ }   \text{ , } \text{ } \tag{\textit{III: Positive probability of intersection with } $\mathrm{Slice}_{y-1,y,y+1}$}
\end{align*}

\noindent which states that a previously shown property, (\textit{Positive probability of intersection with } $\mathrm{Slice}_y$), also holds for the intersection between the subpartition $ \mathrm{Slice}_{y-1,y,y+1}      \text{ } \equiv \text{ } \mathrm{Slice}_{y-1} \cup \mathrm{Slice}_y \cup \mathrm{Slice}_{y+1} \text{ }  \subset \text{ } \underset{y}{\bigcup} \text{ }  \mathrm{Slice}_y \text{ } $ and $\widetilde{\mathscr{V}_{i^{\prime}}}$,

\begin{align*}
           \textbf{P}^{\xi^{\mathrm{Sloped}}}_{\mathscr{O}} \bigg[ \forall \text{ }  i^{\prime} \in \textbf{N} \text{ } , \text{ }          \exists \text{ }  (\widetilde{\gamma_{i^{\prime}}})_L , (\widetilde{\gamma_{i^{\prime}}})_R \in F \big( \mathscr{O} \big)   : \mathscr{O} \supset \bigg[  \text{ }  (\widetilde{\gamma_{i^{\prime}}})_L \text{ } \cap \text{ }   \mathcal{T} \big(     {\mathscr{V}_{i^{\prime}}}   \big) \text{ } \cap \text{ } \mathcal{B} \big(   {\mathscr{V}_{i^{\prime}}}   \big)    \text{ } \cap \text{ }      (\widetilde{\gamma_{i^{\prime}}})_R   \text{ } \bigg]  \supset       \text{ }  \widetilde{\mathscr{V}_{i^{\prime}}}   \bigg]   \text{ } \geq \text{ } 0               \text{ }   \text{ , } \text{ } \tag{\textit{IV: Existence of symmetric domains in the cylinder that are bound by left and right boundaries}}
\end{align*}

\noindent which states that left and right boundaries, respectively given by $(\widetilde{\gamma_{i^{\prime}}})_L$ and $(\widetilde{\gamma_{i^{\prime}}})_R$, exist as a subset of faces of $\mathscr{O}$, and, $\forall \text{ } i_{\mathrm{min}} \equiv \mathrm{inf} \text{ } i^{\prime} \text{ }$,

\begin{align*}
     \text{ }  \textbf{P}^{\xi^{\mathrm{Sloped}}}_{\mathscr{O}} \big[  i_{\mathrm{min}} > i \text{ }     :      \big| \mathcal{I}_{i,y}        \big| \equiv \emptyset                \big] \text{ } \geq \text{ } 0     \text{ }  \text{ , } \tag{\textit{V: Existence of the lowest symmetric domain intersecting } $\mathrm{Slice}_{y-1,y,y+1}$ }
\end{align*}

\noindent which states that there is a minimum value, $i_{\mathrm{min}}$, below which there is zero probability of obtaining any faces $\mathscr{F} \in F \big( \text{ } \mathscr{O} \text{ } \big)$, in which $\mathscr{F} \cap F \big( \text{ } \widetilde{\mathscr{V}_{i^{\prime}}} \text{ } \big) \equiv \emptyset$. Denote the left and right boundaries of the lowest of such domains $\widetilde{\mathscr{V}_{i^{\prime}}}$ with $(\widetilde{\gamma_{i_{\mathrm{min}}}})_L$, and with $(\widetilde{\gamma_{i_{\mathrm{min}}}})_R$, respectively, which exist with positive probability given one previous item in the sequence of inequalities above, (\textit{Existence of symmetric domains in the cylinder that are bound by left and right boundaries}). From the construction of $\gamma_i$, the left and right boundaries for the lowest such domain are given by,

\begin{align*}
    (\widetilde{\gamma_{i_{\mathrm{min}}}})_L  \equiv \underset{\emptyset \neq \mathscr{F}}{\bigcup} \text{ }     \big\{   \forall L > 0  ,  \mathscr{F} \in F \big( \mathscr{O}  \big) :   \mathscr{F} \cap   F \big(  \gamma_{L+i^{*}} \cap \gamma_{L+ 3\frac{i^{*}}{2}}  \big) \neq \emptyset \big\}          \text{ } \text{ , } 
\end{align*}

\noindent and also by,

\begin{align*}
  (\widetilde{\gamma_{i_{\mathrm{min}}}})_R   \equiv     \underset{\emptyset \neq \mathscr{F}}{\bigcup} \text{ }     \big\{    \forall L > 0 ,     \mathscr{F} \in F \big(  \mathscr{O} \big) :   \mathscr{F} \cap    F \big(    \gamma_{L + 3\frac{i^{*}}{2}} \cap \gamma_{L+2i^{*}}  \big) \neq \emptyset      \big\}               \text{ } \text{ , } 
\end{align*}

\noindent where $F \big(   \gamma_{L+i^{*}} \cap \gamma_{L+ 3\frac{i^{*}}{2}}    \big)$ and $F \big(    \gamma_{L + 3\frac{i^{*}}{2}} \cap \gamma_{L+2i^{*}}    \big)$, where $\big\{  \gamma_{L+i^{*}} \cap \gamma_{L+ 3\frac{i^{*}}{2}}\big\} , \big\{   \gamma_{L + 3\frac{i^{*}}{2}} \cap \gamma_{L+2i^{*}}  \big\}  \subset \mathscr{O}_{+ \infty}$, denote the subset of faces from $F \big(\mathscr{O}  \big)$ satisfying,

\begin{align*}
   F \big(  \gamma_{L+i^{*}} \cap \gamma_{L+ 3\frac{i^{*}}{2}}           \big) \equiv   \big\{     F \in F \big(  \mathscr{O} \big) :  F \in \big(  F \big(  \gamma_{L+i^{*}}  \big) \cap F \big(  \gamma_{L+\frac{i^{*}}{2}}  \big)   \big)         \big\}    \text{ }   \text{ , } \text{ } 
\end{align*}

\noindent and similarly that,

\begin{align*}
       F \big(      \gamma_{L + 3\frac{i^{*}}{2}} \cap \gamma_{L+2i^{*}}     \big)  \equiv     \big\{     F \in F \big(  \mathscr{O}  \big) :  F \in \big(  F \big( \gamma_{L+3\frac{i^{*}}{2}} \big) \cap F \big(  \gamma_{L+2i^{*}}\big) \big)     \big\}      \text{ }   \text{ . } \text{ } 
\end{align*}

\bigskip

  \begin{figure}
\begin{align*}
\includegraphics[width=0.88\columnwidth]{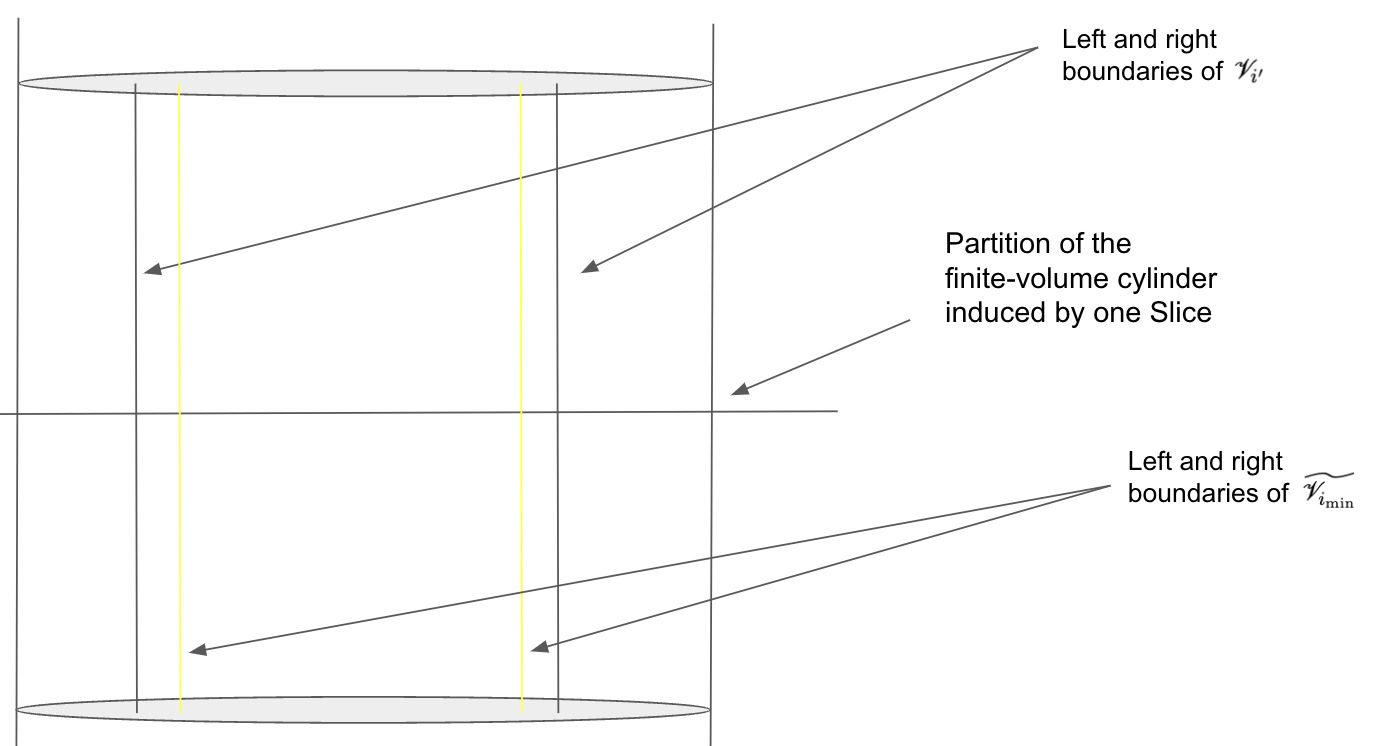}\\
\end{align*}
\caption{\textit{An illustration of the left and right boundaries induced by the symmetric domains} $\mathscr{V}_i$ and $\widetilde{\mathscr{V}_{i_{\mathrm{min}}}}$. Given a partition of $\mathscr{O}$ induced by some intersection of $\mathrm{Slice}$ as indicated above, the left and right boundaries of $\mathscr{V}_i$, in addition to the left and right boundaries of $\widetilde{\mathscr{V}_{i_{\mathrm{min}}}}$, which can be further displaced within the interior of $\mathscr{O}$ as indicated by the two yellow lines, induce a smaller domain $\bar{\mathscr{V}_i} \equiv \mathscr{V}_i \cap \big(  (\widetilde{\gamma_{i_{\mathrm{min}}}})_L \cap  \mathcal{T} \big( \widetilde{\mathscr{V}_{i_{\mathrm{min}}}}\big) \cap (\widetilde{\gamma_{i_{\mathrm{min}}}})_R \cap  \mathcal{B} \big( \widetilde{\mathscr{V}_{i_{\mathrm{min}}}}\big)  \big)$ spanned by the left and right boundaries of $\widetilde{\mathscr{V}_{i_{\mathrm{min}}}}$. $\mathcal{T} \big( \widetilde{\mathscr{V}_{i_{\mathrm{min}}}}\big)$ and $\mathcal{B} \big( \widetilde{\mathscr{V}_{i_{\mathrm{min}}}}\big)$ denote subsets of $T \big(  \mathscr{O} \big)$, and of $\mathcal{B} \big(  \mathscr{O}  \big)$, introduced in (\textit{Property VI}) below.}
\end{figure}

\noindent From $(\widetilde{\gamma_{i_{\mathrm{min}}}})_L$ and $(\widetilde{\gamma_{i_{\mathrm{min}}}})_R$ introduced in the previous paragraph, 

\begin{align*}
  \text{ }      \textbf{P}^{\xi^{\mathrm{Sloped}}}_{\mathscr{O}} \bigg[   \bigg[  (\widetilde{\gamma_{i^{\prime}}})_L \text{ } \cap \text{ }    \mathcal{T} \big( \mathscr{V}_{i^{\prime}} \text{ } \big) \text{ } \cap \text{ }      \mathcal{B} \big( \mathscr{V}_{i^{\prime}} \text{ } \big)  \text{ }  \cap \text{ } (\widetilde{\gamma_{i^{\prime}}})_R \bigg] \supset  \bigg[      (\widetilde{\gamma_{i_{\mathrm{min}}}})_L    \text{ } \cap \text{ } \mathcal{T} \big(  \widetilde{\mathscr{V}_{i_{\mathrm{min}}}} \big)  \text{ } \cap \text{ }        \mathcal{B} \big(\widetilde{\mathscr{V}_{i_{\mathrm{min}}}} \big)    \text{ } \cap    (\widetilde{\gamma_{i^{\mathrm{min}}}})_R       \bigg]  \bigg] \text{ } \geq \text{ } 0    \text{ , } \tag{\textit{VI: Containment of left and right boundaries for the lowest symmetric domain intersecting} $\mathrm{Slice}_y$}
\end{align*}

\noindent states that the left and right boundaries of the lowest domain, $\widetilde{\mathscr{V}_{i_{\mathrm{min}}}}$, with positive probability are contained within the subregion of the cylinder given by the union of the top and bottom parts of $\mathrm{Slice}_{y-1,y,y+1}$, $(\widetilde{\gamma_{i^{\prime}}})_L$, and $(\widetilde{\gamma_{i^{\prime}}})_R$, where,

\begin{align*}
         \mathcal{T} \big( \widetilde{\mathscr{V}_{i_{\mathrm{min}}}}\big) \subset \mathcal{T} \big( \mathscr{V}_{i^{\prime}}          \big) \text{ , } \\ \text{ }  \mathcal{B} \big( \widetilde{\mathscr{V}_{i_{\mathrm{min}}}} \big) \subset \mathcal{B} \big( \mathscr{V}_{i^{\prime}}   \big)    \text{ , }    \text{ } 
\end{align*}    

\noindent see the figure on the top of this page for a depiction of the containment described above between boundaries of \textit{cylindrical domains}. Finally, fix $ (\widetilde{\gamma_{i_{\mathrm{min}}}})_L , (\widetilde{\gamma_{i_{\mathrm{min}}}})_R  \subset F ( \text{ } \mathscr{O} \text{ } )$, from which the definition introduced for $\bar{\mathscr{V}_i}$ in the figure above, 

\begin{align*}
   \textbf{P}^{\xi^{\mathrm{Sloped}}}_{\mathscr{O}} \bigg[ \forall   (\widetilde{\gamma_{i_{\mathrm{min}}}})_L , (\widetilde{\gamma_{i_{\mathrm{min}}}})_R \text{ } , \text{ }  \exists \bar{\mathscr{V}_i} :  \bigg[  \text{ } \mathrm{Slice}_{y-1,y,y+1} \cap \mathscr{O} \cap \big(      (\widetilde{\gamma_{i_{\mathrm{min}}}})_L     \cup  (\widetilde{\gamma_{i_{\mathrm{min}}}})_R  \big)   \bigg]  \supset  \text{ } \bar{\mathscr{V}_i}   \bigg] \geq 0  \text{ } \text{ , } \tag{\textit{VII: The restriction of the symmetric domain} $\mathscr{V}_i$ \textit{to} $\bar{\mathscr{V}_i}$ \textit{exists with positive probability}}
\end{align*}

\noindent states that, with positive probability, there exists $\bar{\mathscr{V}_i} \subset \mathscr{V}_i$, which is bound by the intersection of $ \mathcal{T} \big( \widetilde{\mathscr{V}_{i_{\mathrm{min}}}}\big) \cup  \mathcal{B} \big( \widetilde{\mathscr{V}_{i_{\mathrm{min}}}}\big)$ with the left and right boundaries of $\widetilde{\mathscr{V}_{i_{\mathrm{min}}}}$.

\bigskip

\noindent Otherwise, for $y \neq 2 \text{ } \mathrm{mod} \text{ } 3$, the following sequence of additional requirements on each $\widetilde{\mathscr{V}_{i^{\prime}}} \subset \mathscr{V}_{i^{\prime}}$ hold, in which,

\begin{align*}
     \text{ }  \textbf{P}^{\xi^{\mathrm{Sloped}}}_{\mathscr{O}} \bigg[    y  \text{ }   \equiv   0  \text{ }  \mathrm{mod} \text{ }   3    :          \big| \text{ }      \mathrm{Slice}_{y-1,y,y+1}             \text{ }  \cap   \text{ }    \mathcal{T} \big(  \widetilde{\mathscr{V}_{i_{\mathrm{min}}}}  \big)    \text{ } \big| \equiv 1                             
    \text{ } \bigg] \text{ } \geq 0        \tag{\textit{VIII: Positive probability of connectivity between when} $y \equiv  0 \text{ } \mathrm{mod} \text{ } 3$} \text{ , }
\end{align*}

\noindent states that for $ y  \text{ }   \equiv \text{ }   0 \text{ } \mathrm{mod} \text{ } 3$, the top part of the lowest domain intersecting $\mathrm{Slice}_{y-1,y,y+1}$ only has one point of intersection, while,

\begin{align*}
\textbf{P}^{\xi^{\mathrm{Sloped}}}_{\mathscr{O}} \bigg[      y  \text{ }   \equiv \text{ }   1 \text{ } \mathrm{mod} \text{ } 3   :          \big| \text{ }      \mathrm{Slice}_{y-1,y,y+1} \text{ } \cap \text{ }          \mathcal{B} \big( \widetilde{\mathscr{V}_{i_{\mathrm{min}}}} \big)       \text{ } \big|    \equiv 1          \text{ } \bigg] \text{ }  \geq 0      \text{ , }     \tag{\textit{VIV: Positive probability of connectivity between} $\mathrm{Slice}$ \textit{ and the bottom of the lowest domain when} $y \equiv  1 \text{ } \mathrm{mod} \text{ } 3$}
\end{align*}

\noindent states that for $   y  \text{ }   \equiv \text{ }   1 \text{ } \mathrm{mod} \text{ } 3 \text{ }$,  the bottom part of the lowest domain intersecting $\mathrm{Slice}_{y-1,y,y+1}$ only has one point of intersection,

\begin{align*}
       \textbf{P}^{\xi^{\mathrm{Sloped}}}_{\mathscr{O}} \bigg[  y \equiv 0 \text{ } \mathrm{mod} \text{ } 3 ,   1 \text{ } \mathrm{mod} \text{ } 3 :     \big| \text{ }     \mathrm{Slice}_{y-1,y,y+1}   \text{ }   \cap   \text{ }    \mathcal{T} \big( \text{ } \widetilde{\mathscr{V}_{i_{\mathrm{min}}}} \text{ } \big)  \text{ } \big|  \equiv 1    \text{ } \bigg] \geq 0   \tag{\textit{X: Positive probability of connectivity between} $\mathrm{Slice}$ \textit{ and the top of the lowest domain when} $y \equiv  0 \text{ } \mathrm{mod} \text{ } 3$} 
\end{align*}

\noindent states that for $y \equiv  0 \text{ } \mathrm{mod} \text{ } 3$, the bottom part of the lowest domain intersecting $\mathrm{Slice}_{y-1,y,y+1}$ only has one point of intersection,

\begin{align*}
        \textbf{P}^{\xi^{\mathrm{Sloped}}}_{\mathscr{O}} \bigg[    \forall \text{ }   i^{\prime} \in \textbf{N}  ,          \exists \text{ }  \bar{\mathscr{V}_{i^{\prime}}} \subset F \big(  \mathscr{O} \cap \mathrm{Slice}_{y-1,y,y+1}  \big) : \big| \text{ }  \bar{\mathscr{V}_{i^{\prime}}}\text{ } \cap\text{ }                 \mathrm{Slice}_{y-1,y,y+1}   \text{ } \big|  \neq \emptyset                    \text{ }     \bigg] \geq 0      \text{ } \text{ , } \tag{\textit{XI: Existence of the union of symmetric domains} $\bar{\mathscr{V}_{i^{\prime}}}$ \textit{intersecting} $\mathrm{Slice}_{y-1,y,y+1}$}
\end{align*}

\noindent states that the union of all \textit{cylindrical} domains intersecting $\mathrm{Slice}_{y-1,y,y+1}$,

\begin{align*}
  \emptyset \neq \bar{\mathscr{V}_{i^{\prime}}} \text{ } \equiv \text{ }   \bar{\mathscr{V}_{i^{\prime},y-1,y,y+1}} \equiv \bigg[  \underset{i^{\prime} \in \textbf{N}}{\bigcup} \text{ }          \widetilde{\mathscr{V}_{i^{\prime}}}           \bigg]      \cap  \mathrm{Slice}_{y-1,y,y+1} \text{ } \text{ , } 
\end{align*}
 
 \noindent so that $\bar{\mathscr{V}_{i^{\prime},y-1,y,y+1}} \supset \widetilde{\mathscr{V}_{i^{\prime}}} \text{ }  \forall  \text{ }  i^{\prime}$. Before defining each component of $\bar{\mathscr{V}_{i^{\prime},y-1,y,y+1}}$, introduce, the index set,
 
 \begin{align*}
        \mathcal{I}_{\mathscr{O}} \equiv \text{ } \big\{   i^{\prime} : \text{ } \bar{\mathscr{V}_{i^{\prime}}} \cap \mathrm{Slice}_{y-1,y,y+1} \neq \emptyset       \big\}        \text{ } \text{ . } 
 \end{align*}

 \noindent For the following arguments, denote the top, left, right, and bottom, respectively, boundaries of the union $\bar{\mathscr{V}_{i^{\prime},y-1,y,y+1}}$ of domains $\widetilde{V_{i^{\prime}}}$ with,
 
 \begin{align*}
  \mathrm{Top}_{x,y}  \equiv \text{ } \mathcal{T} \big(\bar{\mathscr{V}_{i^{\prime},y}}  \big) \text{ } \equiv \text{ }  \mathcal{T}\big(  \bar{\mathscr{V}_{i^{\prime}}}     \big) \text{ } \equiv \text{ }   \underset{y}{\underset{i^{\prime} \in \mathcal{I}_{\mathscr{O}}}{\bigcup}} \text{ } \mathcal{T} \big(  \widetilde{\mathscr{V}_{i^{\prime},y}}\big)       \text{ } \equiv \text{ }   \underset{i^{\prime} \in \mathcal{I}_{\mathscr{O}}}{\bigcup} \text{ } \mathcal{T} \big( \widetilde{\mathscr{V}_{i^{\prime}}}\big) \text{ }   
  \\ \bar{(\gamma_{i^{\prime}})_L} \text{ } \equiv \text{ }     \underset{y}{\underset{i^{\prime}\in \mathcal{I}_{\mathscr{O}}}{\bigcup }}    \text{ }  \widetilde{(\gamma_{i^{\prime},y})}_L        \text{ }   \equiv \text{ }  \underset{i^{\prime}\in \mathcal{I}_{\mathscr{O}}}{\bigcup }    \text{ }  \widetilde{(\gamma_{i^{\prime}})}_L \text{ } \text{ , } \\ 
  \mathrm{Bottom}_{x,y} \text{ } \equiv \text{ }  \mathcal{B} \big(  \bar{\mathscr{V}_{i^{\prime},y}}  \big)           \text{ }   \equiv \text{ }     \mathcal{B} \big( \bar{\mathscr{V}_{i^{\prime}}}  \big)         \text{ } \equiv  \text{ }    \underset{y}{\underset{i^{\prime}\in \mathcal{I}_{\mathscr{O}}}{\bigcup}} \text{ } \mathcal{B} \big(  \widetilde{\mathscr{V}_{i^{\prime},y}}    \big)          \text{ }  \equiv \text{ } \underset{i^{\prime}\in \mathcal{I}_{\mathscr{O}}}{\bigcup} \text{ } \mathcal{B} \big(  \widetilde{\mathscr{V}_{i^{\prime}}}     \big)  \\  \bar{(\gamma_{i^{\prime}})_R}   \equiv \text{ }      \underset{y}{\underset{i^{\prime}\in \mathcal{I}_{\mathscr{O}}}{\bigcup }}  \text{ }     \widetilde{(\gamma_{i^{\prime},y})}_R        \text{ } \equiv \text{ }  \underset{i^{\prime}\in \mathcal{I}_{\mathscr{O}}}{\bigcup }  \text{ }     \widetilde{(\gamma_{i^{\prime}})}_R     \text{ } \text{ , } 
 \end{align*}

 \noindent so that $\bar{\mathscr{V}_{i^{\prime},y-1,y,y+1}}$ is composed of the boundaries contained within $\mathcal{I} \big( \text{ } \mathscr{O} \text{ } \big)$,

 \begin{align*}
 \bar{\mathscr{V}_{i^{\prime},y-1,y,y+1}} \equiv \text{ }   \mathrm{Top}_{x,y} \text{ } \cup \text{ }    \bar{(\gamma_{i^{\prime}})_L} \text{ } \cup \text{ }    \mathrm{Bottom}_{x,y} \text{ } \cup \text{ } \bar{(\gamma_{i^{\prime}})_R}     \text{ }  \text{ , } 
 \end{align*}

  \noindent with $\mathrm{Top}_{x,y} , \bar{(\gamma_{i^{\prime}})_L} , \mathrm{Bottom}_{x,y} , \bar{(\gamma_{i^{\prime}})_R} \subset \big( \mathscr{O} \cap   \mathrm{Slice}_{y-1,y,y+1}  \big)$, and,

  \begin{align*}
\big( \mathcal{T} \big( \text{ } \widetilde{\mathscr{V}_{i^{\prime}}} \text{ } \big) \text{ }  \cap \text{ } \mathcal{B} \big( \text{ } \widetilde{\mathscr{V}_{i^{\prime}}} \text{ } \big) \text{ } \cap \text{ }      \bar{(\gamma_{i^{\prime}})_L}   \text{ } \cap \text{ }  \bar{(\gamma_{i^{\prime}})_R}    \big)            \supset \big( \mathcal{T} \big( \text{ } \widetilde{\mathscr{V}_{i^{\prime}}} \text{ } \big) \text{ }  \cap \text{ }  \mathcal{B} \big( \text{ } \widetilde{\mathscr{V}_{i^{\prime}}} \text{ } \big) \text{ }  \cap \text{ }  \widetilde{(\gamma_{i^{\prime}})}_L \text{ } \cap \text{ }    \widetilde{(\gamma_{i^{\prime}})}_R   \big) \text{ } \text{ , } \text{ }  \end{align*}

\begin{align*}
\big( \mathcal{T} \big( \text{ } \bar{\mathscr{V}_{i^{\prime}}} \text{ } \big) \text{ }  \cap \text{ } \mathcal{B} \big( \text{ } \bar{\mathscr{V}_{i^{\prime}}} \text{ } \big) \text{ } \cap \text{ }      \bar{(\gamma_{i^{\prime}})_L}   \text{ } \cap \text{ }  \bar{(\gamma_{i^{\prime}})_R}    \big)            \supset \big( \mathcal{T} \big( \text{ } \widetilde{\mathscr{V}_{i^{\prime}}} \text{ } \big) \text{ }  \cap \text{ }  \mathcal{B} \big( \text{ } \widetilde{\mathscr{V}_{i^{\prime}}} \text{ } \big) \text{ }  \cap \text{ }  \widetilde{(\gamma_{i^{\prime}})}_L \text{ } \cap \text{ }    \widetilde{(\gamma_{i^{\prime}})}_R   \big) \text{ } \text{ , } \text{ } \\ \big( \mathcal{T} \big( \text{ } \bar{\mathscr{V}_{i^{\prime}}} \text{ } \big) \text{ }  \cap \text{ } \mathcal{B} \big( \text{ } \bar{{\mathscr{V}_{i^{\prime}}}} \text{ } \big) \text{ } \cap \text{ }      \widetilde{(\gamma_{i^{\prime}})}_L   \text{ } \cap \text{ }  \bar{(\gamma_{i^{\prime}})_R}    \big)            \supset \big( \mathcal{T} \big( \text{ } \widetilde{\mathscr{V}_{i^{\prime}}} \text{ } \big) \text{ }  \cap \text{ }  \mathcal{B} \big( \text{ } \widetilde{\mathscr{V}_{i^{\prime}}} \text{ } \big) \text{ }  \cap \text{ }  \widetilde{(\gamma_{i^{\prime}})}_L \text{ } \cap \text{ }    \bar{(\gamma_{i^{\prime}})}_R   \big) \text{ } \text{ , } \text{ }  \\     \big( \mathcal{T} \big( \text{ } \bar{\mathscr{V}_{i^{\prime}}} \text{ } \big) \text{ }  \cap \text{ } \mathcal{B} \big( \text{ } \bar{\mathscr{V}_{i^{\prime}}} \text{ } \big) \text{ } \cap \text{ }      \bar{(\gamma_{i^{\prime}})_L}   \text{ } \cap \text{ }  \widetilde{(\gamma_{i^{\prime}})}_R    \big)            \supset \big( \mathcal{T} \big( \text{ } \widetilde{\mathscr{V}_{i^{\prime}}} \text{ } \big) \text{ }  \cap \text{ }  \mathcal{B} \big( \text{ } \widetilde{\mathscr{V}_{i^{\prime}}} \text{ } \big) \text{ }  \cap \text{ }  \bar{(\gamma_{i^{\prime}})}_L \text{ } \cap \text{ }    \widetilde{(\gamma_{i^{\prime}})}_R   \big) \text{ } \text{ , }  \\ \big( \mathcal{T} \big( \text{ } \bar{\mathscr{V}_{i^{\prime}}} \text{ } \big) \text{ }  \cap \text{ } \mathcal{B} \big( \text{ } \bar{\mathscr{V}_{i^{\prime}}} \text{ } \big) \text{ } \cap \text{ }      \bar{(\gamma_{i^{\prime}})_L}   \text{ } \cap \text{ }  \bar{(\gamma_{i^{\prime}})_R}    \big)            \supset \big( \mathcal{T} \big( \text{ } \widetilde{\mathscr{V}_{i^{\prime}}} \text{ } \big) \text{ }  \cap \text{ }  \mathcal{B} \big( \text{ } \widetilde{\mathscr{V}_{i^{\prime}}} \text{ } \big) \text{ }  \cap \text{ }  \widetilde{(\gamma_{i^{\prime}})}_L \text{ } \cap \text{ }    \bar{(\gamma_{i^{\prime}})}_R   \big) \text{ } \text{ , } \text{ }   \\     \big( \mathcal{T} \big( \text{ } \bar{\mathscr{V}_{i^{\prime}}} \text{ } \big) \text{ }  \cap \text{ } \mathcal{B} \big( \text{ } \bar{\mathscr{V}_{i^{\prime}}} \text{ } \big) \text{ } \cap \text{ }      \bar{(\gamma_{i^{\prime}})_L}   \text{ } \cap \text{ }  \bar{(\gamma_{i^{\prime}})_R}    \big)            \supset \big( \mathcal{T} \big( \text{ } \widetilde{\mathscr{V}_{i^{\prime}}} \text{ } \big) \text{ }  \cap \text{ }  \mathcal{B} \big( \text{ } \widetilde{\mathscr{V}_{i^{\prime}}} \text{ } \big) \text{ }  \cap \text{ }  \bar{(\gamma_{i^{\prime}})}_L \text{ } \cap \text{ }    \widetilde{(\gamma_{i^{\prime}})}_R   \big) \text{ } \text{ , }    \\   \big( \mathcal{T} \big( \text{ } \bar{\mathscr{V}_{i^{\prime}}} \text{ } \big) \text{ }  \cap \text{ } \mathcal{B} \big( \text{ } \bar{\mathscr{V}_{i^{\prime}}} \text{ } \big) \text{ } \cap \text{ }      \widetilde{(\gamma_{i^{\prime}})}_L   \text{ } \cap \text{ }  \widetilde{(\gamma_{i^{\prime}})}_R   \big)            \supset \big( \mathcal{T} \big( \text{ } \widetilde{\mathscr{V}_{i^{\prime}}} \text{ } \big) \text{ }  \cap \text{ }  \mathcal{B} \big( \text{ } \widetilde{\mathscr{V}_{i^{\prime}}} \text{ } \big) \text{ }  \cap \text{ }  \widetilde{(\gamma_{i^{\prime}})}_L \text{ } \cap \text{ }    \widetilde{(\gamma_{i^{\prime}})}_R   \big) \text{ } \text{ , } \text{ }  \forall \text{ } i^{\prime}       \text{ } \text{ , } 
  \end{align*}

  \noindent which are a series of several containment between either $\widetilde{\mathscr{V}_{i^{\prime}}}$, or $\bar{\mathscr{V}_{i^{\prime}}}$, from the faces of the cylinder intersecting $\mathrm{Slice}_{y-1,y,y+1}$, $F \big( \text{ } \mathscr{O}\text{ }  \cap \text{ }  \mathrm{Slice}_{y-1,y,y+1} \text{ } \big) \subset F \big( \text{ } \mathscr{O} \text{ } \big)$.
  
  \bigskip
  
  \noindent The corresponding connectivity event in the union of \textit{symmetric} domains takes the form,
 
 \begin{align*}
     \bigg\{ \text{ } \mathrm{Top}_{x,y} \underset{\bar{\mathscr{V}_{i^{\prime}}}}{\overset{h\geq ck}{\longleftrightarrow}} \mathrm{Bottom}_{x,y}   \text{ } \bigg\}        \text{ }  \text{ , } 
 \end{align*}

 \noindent which is pushed forwards under $\textbf{P}^{\xi^{\mathrm{Sloped}}}_{\mathscr{O}} \big[  \cdot \big]$ as,

 \begin{align*}
  \textbf{P}^{\xi^{\mathrm{Sloped}}}_{\mathscr{O}} \bigg[ \text{ } \mathrm{Top}_{x,y} \underset{\bar{\mathscr{V}_{i^{\prime}}}}{\overset{h\geq ck}{\longleftrightarrow}} \mathrm{Bottom}_{x,y}  \text{ } \bigg]    \text{ } \text{ . } 
 \end{align*}

 \noindent Moreover, by construction, recall,

 \begin{align*}
         \mathrm{Top}_{x,y}   \supseteq \mathcal{T} \big(  \bar{\mathscr{V}_{i^{\prime}}}  \big) \equiv  \underset{i^{\prime} \in \mathcal{I}_{\mathscr{O}}}{\bigcup }      \mathcal{T} \big( \widetilde{\mathscr{V}_{i^{\prime}}}  \big)          \supset   \underset{i^{\prime}\in \mathcal{I}_{\mathscr{O}} }{\bigcup}  \underset{j^{\prime}\in \mathcal{I}_{\mathscr{O}}}{\bigcap}  \bigg[    \mathcal{T} \big( \widetilde{\mathscr{V}_{i^{\prime}}}  \big)  \cap  \mathcal{T} \big( \widetilde{\mathscr{V}_{j^{\prime}}}  \big)  \bigg]  \supset \mathcal{T} \big( \widetilde{\mathscr{V}_{i^{\prime}}}  \big)  \supset  \mathcal{T} \big( \widetilde{\mathscr{V}_{j^{\prime}}}  \big) \text{ , } 
 \end{align*}
 
 \noindent and also that,
 
 \begin{align*}
       \mathrm{Bottom}_{x,y}  \supseteq  \mathcal{B} \big(  \bar{\mathscr{V}_{i^{\prime}}} \big) \equiv  \underset{i^{\prime} \in \mathcal{I}_{\mathscr{O}}}{\bigcup }      \mathcal{B} \big(  \widetilde{\mathscr{V}_{i^{\prime}}}  \big)  \supset \underset{i^{\prime}\in \mathcal{I}_{\mathscr{O}}}{\bigcup}  \underset{j^{\prime} \in \mathcal{I}_{\mathscr{O}}}{\bigcap}  \bigg[    \mathcal{B} \big(  \widetilde{\mathscr{V}_{i^{\prime}}}  \big) \cap \mathcal{B} \big( \widetilde{\mathscr{V}_{j^{\prime}}}  \big)  \bigg]   \supset  \mathcal{B} \big(  \widetilde{\mathscr{V}_{i^{\prime}}}  \big)   \supset        \mathcal{B} \big(  \widetilde{\mathscr{V}_{j^{\prime}}}  \big)   \text{ } \text{ , }
 \end{align*}
 
\noindent each of which hold $\forall \text{ } i^{\prime} > j^{\prime}$.

\bigskip

\noindent Next, from properties of each $\bar{\mathscr{V}_{i^{\prime}}}$, fix a suitable \textit{absolute} constant $\delta >\delta^{\prime} >0$, with the same choice of $j$ as provided throughout results for the strip, so that,

\begin{align*}
  \textbf{P}^{\xi^{\mathrm{Sloped}}}_{\textbf{Z} \times [0,n^{\prime}N]} \big[ [ 0 , \lfloor \delta^{\prime} n \rfloor ]      \times \{ 0 \} \overset{h\geq ck}{\longleftrightarrow}           [ j , \big( \text{ } j + 1 \text{ } \big) \text{ } \lfloor  \delta^{\prime} n     \rfloor     ]  \times \{ n \}        \big] \text{ }       \text{ } \text{ , } 
\end{align*}

\noindent occurs with positive probability, which, recall, was related to several previous estimates and arguments for the connectivity event between $\mathcal{I}_j$ and $\widetilde{\mathcal{I}_j}$,

\begin{align*}
    \textbf{P}^{\xi^{\mathrm{Sloped}}}_{\textbf{Z} \times [0,n^{\prime}N]} \big[  \mathcal{I}_j \overset{h \geq ck}{\longleftrightarrow}  \widetilde{\mathcal{I}_j} \big]           \text{ } \text{ . } 
\end{align*}

\noindent The \textit{tightness} property,

\begin{align*}
\textbf{P}^{\xi^{\mathrm{Sloped}}}_{\mathscr{O}} \bigg[                 \forall \text{ }     \delta > \delta^{\prime\prime} > 0 \text{ } ,   \exists \text{ }    r >      r^{\prime} >     \frac{2k}{\eta}            \text{ }   :      \big| \text{ }        \bar{\mathscr{V}_{i^{\prime},y}}      \cap           \mathcal{T} \big(  \bar{\mathscr{V}_{i^{\prime},y}}      \big)       \text{ } \big|  \leq \lfloor  \delta^{\prime\prime} r^{\prime}                  \rfloor     \bigg] \text{ } \geq 0              \text{ } \text{ , }          \tag{\textit{XII: Tightness of cylindrical domains intersecting} $\mathrm{Slice}_{y-1,y,y+1}$ \text{ from the top of the domain}}      
\end{align*}

\noindent states that the number of faces contained in the intersection $\big| \text{ }     \bar{\mathscr{V}_{i^{\prime},y}}      \cap           \mathcal{T} \big( \text{ }  \bar{\mathscr{V}_{i^{\prime},y}}     \text{ } \big)         \text{ } \big|$ does not exceed the floor of the smallest integer $x$ so that $x < \delta^{\prime\prime}r$, while similarly,

 \begin{figure}
\begin{align*}
\includegraphics[width=1.04\columnwidth]{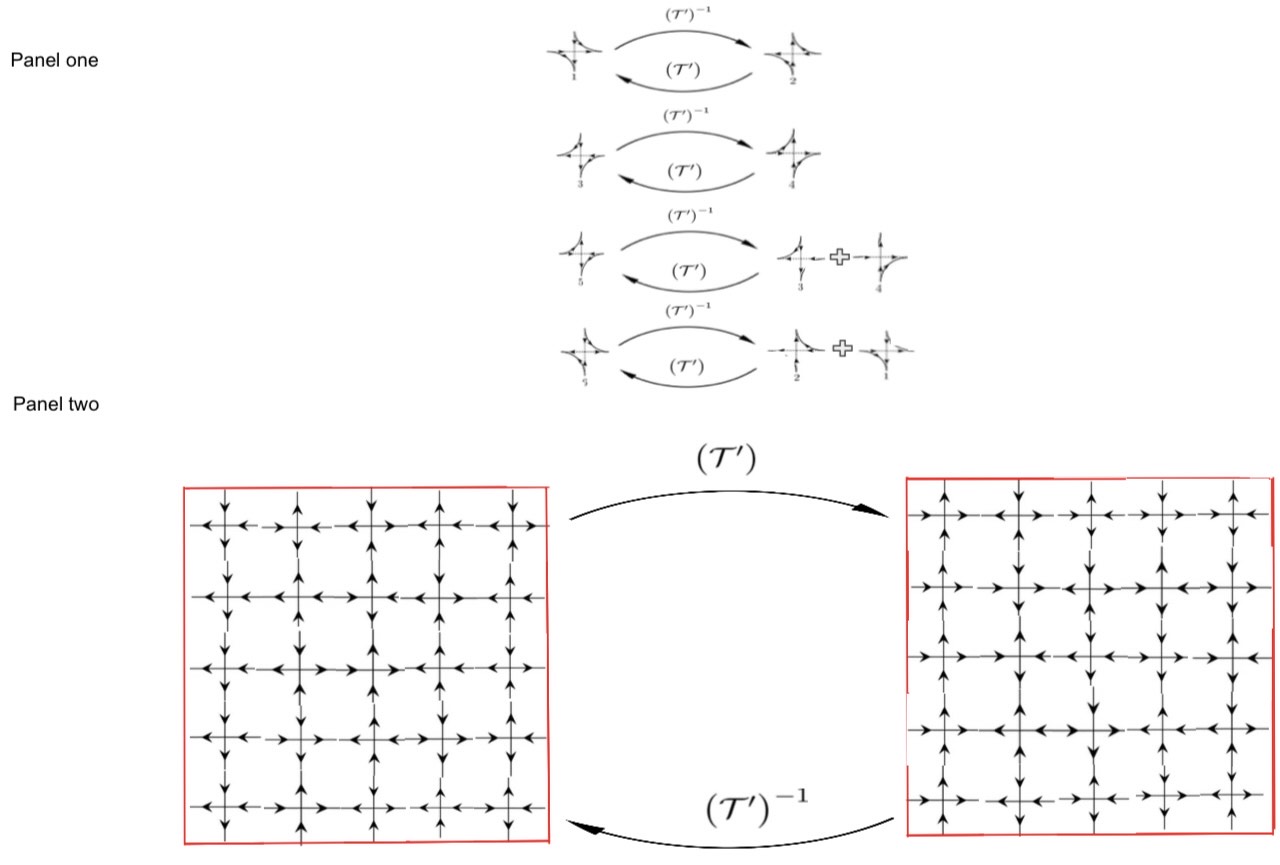}\\
\end{align*}
\caption{\textit{Panel one, A depiction of the action of the loop reversing map on each six-vertex configuration}. For the vertices of the first two types, the configurations under the pullback $(\mathcal{T}\big)^{-1}$ is given by the remaining vertex within the same type. For vertices of types $5$ and $6$, from the left-turning splitting convention originally introduced in {\color{blue}[11]}, the reversed loop orientation corresponding to these vertices of the remaining two types are given by combinations of directed loops from vertices of types $3$ and $4$, in addition to combinations of directed loops from vertices of types $2$ and $1$, respectively. \textit{Panel two, A depiction of the reversal of loops for every other vertex in a sample finite volume contained between} $\gamma_L$ and $\gamma_{L+i^{*}}$. One possible reconstruction of the finite volume provided in \textit{Property two} of $\mathcal{T}^{\prime}$ is provided, with faces outside of $\gamma_L$ and $\gamma_{L+i^{*}}$ excluded from the illustration.}
\end{figure}

\begin{align*}
    \textbf{P}^{\xi^{\mathrm{Sloped}}}_{\mathscr{O}} \bigg[                          \forall \text{ }     \delta > \delta^{\prime\prime} > 0      \text{ } ,  \exists \text{ }  r >    r^{\prime} >       \frac{2k}{\eta}        \text{ }   :                  \text{ }           \big| \text{ }        \bar{\mathscr{V}_{i^{\prime},y}}      \cap        \mathcal{B} \big(   \bar{\mathscr{V}_{i^{\prime},y}}   \big)      \text{ } \big|   \leq \lfloor        \delta^{\prime\prime} r^{\prime}                 \rfloor                 \bigg] \text{ } \geq 0              \text{ } \text{ , }          \tag{\textit{XIII: Tightness of cylindrical domains intersecting} $\mathrm{Slice}_{y-1,y,y+1}$ \text{ from the bottom of the domain}}    
\end{align*}

\noindent states that the number of faces contained in the intersection $  \big| \text{ }        \bar{\mathscr{V}_{i^{\prime},y}}      \cap        \mathcal{B} \big( \bar{\mathscr{V}_{i^{\prime},y}}     \big)      \text{ } \big|$ does not exceed the floor of the smallest integer $y$ so that $y < \delta^{\prime\prime}r$. To conclude, the \textit{goodness} property,

\begin{align*}
       \textbf{P}^{\xi^{\mathrm{Sloped}}}_{\mathscr{O}} \bigg[  \forall \text{ } i^{\prime} \text{ } ,  \exists  \text{ }     \bar{\mathscr{V}_{i^{\prime},y-1}} ,  \bar{\mathscr{V}_{i^{\prime},y}} , \bar{\mathscr{V}_{i^{\prime},y+1}} \subset F \big( \mathscr{O} \cap \mathrm{Slice}_{y-1,y,y+1}  \big)  : \bar{\mathscr{V}_{i^{\prime},y-1}} ,  \bar{\mathscr{V}_{i^{\prime},y}} , \bar{\mathscr{V}_{i^{\prime},y+1}}  \text{ satisfy \textit{Property XIII}}        \text{ }  \bigg]  \geq 0  \text{ }  \text{ , } \tag{\textit{XIV: Good property from the tightness property}}
\end{align*}

\noindent states that $\bar{\mathscr{V}_{i^{\prime},y-1}}$, $\bar{\mathscr{V}_{i^{\prime},y}}$, and $\bar{\mathscr{V}_{i^{\prime},y+1}}$ are simultaneously \textit{tight}, where each domain satisfies, by construction,

\begin{align*}
   \bar{\mathscr{V}_{i^{\prime},y-1}} \cap \mathrm{Slice}_{y-1}      \neq \emptyset  \text{ } \text{ , }  \\    \bar{\mathscr{V}_{i^{\prime},y}} \cap \mathrm{Slice}_{y}       \neq \emptyset  \text{ } \text{ , }  \\     \bar{\mathscr{V}_{i^{\prime},y+1}} \cap \mathrm{Slice}_{y+1}    \neq \emptyset      \text{ } \text{ . } 
\end{align*}

\bigskip

\noindent Following all properties listed, the result below captures the fraction of domains contained in $\bar{\mathscr{V}_{i^{\prime}}}$ which are deemed to satisfy the \textit{good} property, which was introduced in {\color{blue}[11]} as a condition that is satisfied when unions of \textit{symmetric} domains intersecting $\mathrm{Slice}_{y-1}$, $\mathrm{Slice}_y$, and $\mathrm{Slice}_{y+1}$, are each \textit{tight}. The analogue of this condition for different boundary conditions is satisfied with $\bar{\mathscr{V}_{i^{\prime},y-1}}$, $\bar{\mathscr{V}_{i^{\prime},y}}$, and $\bar{\mathscr{V}_{i^{\prime},y+1}}$.

\subsection{Quantification of the fraction of $(x,y)$ in $\mathscr{O}$ which satisfy the \textit{good} property}

The failure of symmetric cylindrical domains with strictly less than $\frac{mn}{6}$ \textit{good} points

\bigskip

\noindent \textbf{Lemma} \textit{4.4} ({\textit{strictly less than half of the admissible points about which symmetric domains can be situated in the finite volume cylinder are good}}). Fix,

\begin{align*}
n^{\prime} \equiv \frac{N}{\eta r^{\prime}} \text{ } \text{ , } 
\end{align*}

\noindent and,

\begin{align*}
m^{\prime} \equiv \frac{M}{r^{\prime}} \text{ } \text{ . } 
\end{align*}

\noindent For $1 \leq x \leq n$ and $1 \leq y \leq m$, with $y \equiv 2 \text{ } \mathrm{mod} \text{ } 3$, strictly less than half of the pairs $(x,y)$ are \textit{good}.

\bigskip

\noindent \textit{Proof of Lemma 4.4}. To demonstrate that the \textit{good} property holds for pairs $(x,y)$ satisfying the conditions provided above, first, it suffices to show that the collection of $\bar{\mathscr{V}_{i^{\prime}}}$ satisfy,

\begin{align*}
    \underset{i^{\prime}}{\bigcap} \text{ }   \bar{\mathscr{V}_{i^{\prime}}}   \equiv \emptyset  \text{ } \text{ , } 
\end{align*}

\noindent where the intersection taken over $i^{\prime}$, in which, for $0 <i^{\prime\prime\prime} < i^{\prime\prime}< i^{\prime}$, for each $\bar{\mathscr{V}_{i^{\prime}}}$ to be disjoint from $\bar{\mathscr{V}_{i^{\prime}+1}}$ above it within $\mathscr{O} \cap \mathrm{Slice}_{y-1,y,y+1,y+2}$, with positive probability,

\begin{align*}
      \bar{\mathscr{V}_{i^{\prime\prime}}} \cap    \mathrm{Slice}_{y+1} \neq \emptyset \text{ } \text{ , } \text{ }    \bar{\mathscr{V}_{i^{\prime\prime}}} \cap    \mathrm{Slice}_{y+2} \equiv \emptyset  \text{ } \text{ , } \end{align*}

      \noindent there must be $i^{\prime\prime}$ so that $\bar{\mathscr{V}_{i^{\prime\prime}}} \cap \mathrm{Slice}_{y-1,y,y+1} \neq \emptyset$, with $y \equiv \text{ } 0 \text{ } \mathrm{mod}\text{ } 3$, and also, that,

      \begin{align*}
      \text{ }    \bar{\mathscr{V}_{i^{\prime\prime\prime}}} \cap \mathrm{Slice}_{y-2} \neq \emptyset \text{ } , \text{ } \bar{\mathscr{V}_{i^{\prime\prime\prime}}} \cap \mathrm{Slice}_{y-1} \equiv \emptyset        \text{ } \text{ , } 
    \end{align*}

\noindent in which there must be a domain $\bar{\mathscr{V}_{i^{\prime\prime\prime}}}$ which does not intersect any of the portions of $\mathrm{Slice}_{y-1,y,y+1}$ which intersect $\bar{\mathscr{V}_{i^{\prime\prime}}}$. To this end, for any increasing sequence of nonnegative integers tending to $i^{\prime}$, given a collection of $\bar{\mathscr{V}_{i^{\prime}}}$, fix $i_1 > i_2 >0$. Under this choice of $i_1$ and $i_2$, it is possible that either $\bar{\mathscr{V}_{i_1}} \text{ } \cap \text{ } \bar{\mathscr{V}_{i_2}} \neq \emptyset$, or that $\bar{\mathscr{V}_{i+_1}} \subset \bar{\mathscr{V}_{i_2}} \subset \big( \text{ } \mathscr{O} \text{ } \cap \text{ } \mathrm{Slice}_{y-1,y,y+1} \text{ } \big)$. By symmetry, we discard the latter scenario in which one union of \textit{symmetric} domains, indexed in $i_1$, can be contained within the union of \textit{symmetric} domains indexed in $i_2$. Additionally, given $\bar{\mathscr{V}_{i_1}} \text{ } \cap \text{ } \mathscr{V}_{i_1}$, this intersection contains the left and right boundaries, respectively introduced previously with $\widetilde{(\gamma_{i_1})}_L$ and $ \widetilde{(\gamma_{i_1})}_R$, at $i_1$. Within the restricted cylindrical finite volume,

\begin{align*}
     \bar{\mathscr{V}_{i_1}} \supset   \bigg[       \widetilde{(\gamma_{i_1})}_L        \cap  \widetilde{(\gamma_{i_1})}_R  \cap       \mathcal{T} \big(  \widetilde{\mathscr{V}_{i_1}}\big)    \cap  \mathcal{B} \big( \widetilde{\mathscr{V}_{i_1}} \big)   \bigg]        \text{ } \text{ . } 
\end{align*}

\noindent Besides the left and right boundaries provided above, for $i_2$ there also exists $\widetilde{(\gamma_{i_2})}_L$ and $ \widetilde{(\gamma_{i_2})}_R$ with the accompanying quad structure,

\begin{align*}
    \bar{\mathscr{V}_{i_2}}  \equiv  \underset{i_2 < i_1}{\bigcup}    \widetilde{\mathscr{V}_{i_2}}        \equiv  \underset{i_2 < i_1}{\bigcup} \big\{         \widetilde{(\gamma_{i_2})}_L   \cap    \widetilde{(\gamma_{i_2})}_R  \cap  \mathcal{T} \big(  \widetilde{V_{i_1}}  \big) \cap   \mathcal{B} \big(    \widetilde{V_{i_1}}   \big)   \big\}     \text{ } \text{ , } 
\end{align*}

\noindent so that, with positive probability, within $\bar{\mathscr{V}_{i_2}}$ the crossing between the left and right boundaries centered about $i_2$, which were introduced for $\bar{\mathscr{V}_{i_2}}$,

\begin{align*}
      \big\{  \widetilde{(\gamma_{i_2})}_L   \underset{( \text{ } \bar{\mathscr{V}_{i_2}} \text{ } )^c  \cap  \mathscr{O} }{\overset{h \geq ck}{\longleftrightarrow}}  \widetilde{(\gamma_{i_2})}_R    \big\}         \text{ } \text{ , } 
\end{align*}

\noindent occurs, which in turn implies that the intersection,

\begin{align*}
 \bigg[    (  \bar{\mathscr{V}_{i_2}} )^c  \cap  \mathscr{O}       \bigg]  \cap       \bigg[     \widetilde{(\gamma_{i_2})}_L                   \cap    \widetilde{(\gamma_{i_2})}_R       \text{ } \cap \text{ }            \mathcal{T} \big( \text{ }  \widetilde{\mathscr{V}_{i_2}} \text{ } \big)       \text{ } \cap \text{ }     \mathcal{B} \big( \text{ }   \widetilde{\mathscr{V}_{i_2}}   \text{ } \big)     \bigg]          \text{ }  \neq \emptyset    \text{ } \text{ , } 
\end{align*}

\noindent does not equal the empty set, and hence that the left to right crossing between $\widetilde{(\gamma_{i_2})}_L $, and $\widetilde{(\gamma_{i_2})}_R$, does not occur within $\bar{\mathscr{V}_{i_2}}$, but instead within the intersection of the complement with $\mathscr{O}$. To begin, suppose,

\begin{align*}
          \bar{\mathscr{V}_0}        \text{ } \cap \text{ }  \bar{\mathscr{V}_1}      \equiv \emptyset                        \text{ } \text{ , } \tag{\textit{HYP}}
\end{align*}

\noindent corresponding to the first \textit{symmetric} domain in $\mathscr{O} \text{ } \cap \text{ } \mathrm{Slice}{y_0,y_1,y_2}$, indexed with $i^{\prime}\equiv 0$. The next inductive step would suggest,

\begin{align*}
 \bigg[   \underset{j\equiv 1}{\bigcup} \text{ } \widetilde{\mathscr{V}_j}  \bigg]  \text{ } \cap \text{ } \widetilde{\mathscr{V}_2}   \text{ } \text{ . } 
\end{align*}

\noindent For all $i^{\prime} > 1$, therefore it suffices to show that,

\begin{align*}
    \bigg[  \underset{j < i^{\prime}}{\bigcup} \text{ }    \bar{\mathscr{V}_j}      \bigg]  \text{ } \cap \text{ }            \bar{\mathscr{V}_{i^{\prime}}} \equiv \emptyset      \text{ } \text{ , } 
\end{align*}

\noindent where in the intersection above, the union of all $\bar{\mathscr{V}_j}$, as long as $j < i^{\prime}$,  . By induction, the containment can be shown to hold for all $i^{\prime}$ for which $\mathscr{O} \text{ } \cap \text{ } \mathrm{Slice}_{i^{\prime}} \neq \emptyset$. Proceeding, if the property is to hold for the next domain in the increasing sequence of nonnegative integers tending to $i^{\prime}+1$, 

\begin{align*}
      \bigg[   \underset{j < i^{\prime}+1}{\bigcup} \text{ }    \bar{\mathscr{V}_j}       \bigg]  \text{ } \cap \text{ }  \bar{\mathscr{V}_{i^{\prime}}}  \text{ } \text{ , }
\end{align*}

\noindent yields a nonempty intersection between the two sets. Besides $\bar{\mathscr{V}_{i^{\prime}}}$ in the intersection above, the union over $j < i^{\prime}+1$ is a subset of $\mathscr{O}$, as,

\begin{align*}
    \underset{j < i^{\prime}+1}{\bigcup} \text{ } \widetilde{\mathscr{V}_j} \text{ } \equiv \text{ }  \underset{j < i^{\prime}+1}{\bigcup} \text{ } \bigg[     \widetilde{\mathscr{V}_{j-1}}    \cap          \widetilde{\mathscr{V}_{j+1}}  \bigg]     \text{ } \equiv \text{ }  \bigg[   \underset{j < i^{\prime}+1}{\bigcup}  \text{ }                \widetilde{\mathscr{V}_{j-1}}     \bigg]  \text{ }  \cap \text{ } \bigg[     \underset{j < i^{\prime}+1}{\bigcup} \text{ }           \widetilde{\mathscr{V}_{j+1}}      \bigg]  \text{ }      \text{ } \text{ , }
\end{align*}

\noindent which in turn implies that,

\begin{align*}
      \underset{j < i^{\prime}+1}{\bigcup} \text{ } \widetilde{\mathscr{V}_j} \text{ }  \neq \emptyset      \text{ } \text{ , } 
\end{align*}

\noindent because,

\begin{align*}
      \underset{j < i^{\prime}+1}{\bigcup}  \text{ }                \widetilde{\mathscr{V}_{j-1}}     \subset \mathscr{O} \neq \emptyset   \text{ } \text{ , } 
\end{align*}

\noindent and similarly, also because,

\begin{align*}
      \underset{j < i^{\prime}+1}{\bigcup} \text{ }           \widetilde{\mathscr{V}_{j+1}}     \subset \mathscr{O} \neq \emptyset      \text{ } \text{ , } 
\end{align*}

\noindent in which each of the two containments above hold, and moreover, are nonempty, due to the fact that not only does each $\widetilde{\mathscr{V}} \cap \mathrm{Slice}_{y-1,y,y+1} \neq \emptyset$, but also that $\widetilde{\mathscr{V}_{j-1}}$, and $\widetilde{\mathscr{V}_{j+1}}$ are by construction nonempty subsets of $\mathscr{O}$. Hence, these arguments demonstrate, from the induction hypothesis (\textit{HYP}), that,

\begin{align*}
          \bar{\mathscr{V}_0}        \text{ } \cap \text{ }  \bar{\mathscr{V}_1}      \equiv \emptyset   \text{ }     \overset{(\textit{HYP})}{\Rightarrow} \text{ }        \widetilde{\mathscr{V}_{j-1}} \cap \widetilde{\mathscr{V}_j} \text{ } \equiv \text{ }       \bigg[       \underset{k^{\prime\prime}}{\bigcup } \text{ } \widetilde{\mathscr{V}_{k^{\prime\prime}}}   \bigg]  \text{ } \cap \text{ }          \widetilde{\mathscr{V}_j}   \text{ } \equiv \text{ } \emptyset     \text{ } \text{ . }
\end{align*}

\noindent The above intersection is nonempty because, besides $\widetilde{\mathscr{V}_j}$, as a countable intersection of the union of \textit{symmetric} domains, being nonempty,

\begin{align*}
         \underset{k^{\prime\prime}}{\bigcup } \text{ } \widetilde{\mathscr{V}_{k^{\prime\prime}}}  \text{ }            \text{ } \text{ , } 
\end{align*}

\noindent demonstrating that (\textit{HYP}) holding for $i$ would then imply that (\textit{HYP}) would also hold for $i+1$, namely,

\begin{align*}
        \underset{1 \leq o \leq i+1}{\bigcap} \text{ }               \widetilde{\mathscr{V}_o} \equiv  \bigg[     \underset{1 \leq o \leq i}{\bigcap}             \widetilde{\mathscr{V}_o} \bigg]  \cap \widetilde{\mathscr{V}_{i+1}} \equiv   \emptyset              \text{ } \text{ , } 
\end{align*}

\noindent so that the intersection of all such sets would be nonempty hence making them disjoint. With the following series of equivalences, observe, 

\begin{align*}
   \bigg[      \underset{k^{\prime\prime}}{\bigcup } \text{ } \widetilde{\mathscr{V}_{k^{\prime\prime}}}    \bigg]   \cap         \widetilde{\mathscr{V}_j}   \text{ } \equiv \text{ } \emptyset \text{ } \overset{(\textit{1})}{\Longleftrightarrow}       \bigg[  \text{ }       \underset{k^{\prime\prime}(i)}{\bigcup}  \text{ }           \bigg[      \widetilde{(\gamma_{i})_L}          \cap           \widetilde{(\gamma_{i})_R}    \bigg]     \cap \mathcal{T}\mathcal{B} \big(  \widetilde{\mathscr{V}_i} \big)         \text{ } \bigg]           \cap \widetilde{\mathscr{V}_j}      \equiv \text{ } \emptyset \text{ } \\   \overset{(\textit{2})}{\Longleftrightarrow}         \bigg[  \text{ }    \underset{k^{\prime\prime}(i)}{\bigcup}  \text{ }           \bigg[     \widetilde{(\gamma_{i})_L}         \cap     \widetilde{(\gamma_{i})_R}        \bigg]     \cap \widetilde{\mathscr{V}_j}  \text{ }  \bigg]   \cap  \bigg[  \text{ } \mathcal{T}\mathcal{B} \big(  \widetilde{\mathscr{V}_i} \big)     \cap      \widetilde{\mathscr{V}_j} \text{ } \bigg]                     \text{ }  \equiv \text{ } \emptyset \text{ } \\   \overset{(\textit{3})}{\Longleftrightarrow}      \bigg[     \underset{k^{\prime\prime}(i)}{\bigcup}    \text{ }  \widetilde{(\gamma_{i})_L}         \cap     \widetilde{\mathscr{V}_j}      \bigg]         \cap \text{ }  \bigg[     \underset{k^{\prime\prime}(i)}{\bigcup}    \widetilde{(\gamma_{i})_R}         \cap     \widetilde{\mathscr{V}_j}      \bigg]   \cap \bigg[    \mathcal{T}\mathcal{B} \big( \widetilde{\mathscr{V}_i}   \big)     \cap      \widetilde{\mathscr{V}_j} \bigg]       \equiv \emptyset  \text{ } \text{ , } 
\end{align*}

\noindent where, 

\begin{align*}
  \mathcal{T}\mathcal{B} \big( \widetilde{\mathscr{V}_k^{\prime\prime}}   \big)  \equiv   \mathcal{T}\big(  \widetilde{\mathscr{V}_k^{\prime\prime}}  \big)  \cap   \mathcal{B}\big( \widetilde{\mathscr{V}_k^{\prime\prime}}  \big) \subset \mathscr{O} \text{ } \text{ , } 
\end{align*}

 \begin{figure}
\begin{align*}
\includegraphics[width=0.14\columnwidth]{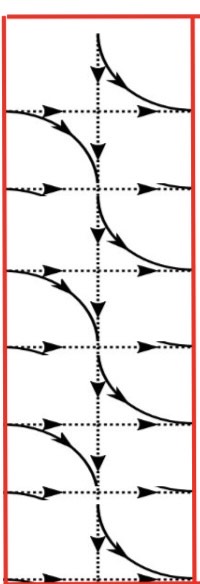}\\
\end{align*}
\caption{\textit{A depiction of the first five faces incident to the left finite volume boundary provided in the previous figure, Figure 35, before reversing the distribution of directed, oriented, non-intersecting, loops under} $\big(\mathcal{T}\big)^{-1}$.}
\end{figure}

\noindent where the rearrangements above immediately imply, for $i+1$, 

\begin{align*}
   \bigg[      \underset{k^{\prime\prime}+1}{\bigcup } \text{ } \widetilde{\mathscr{V}_{k^{\prime\prime}}}    \bigg]\cap         \widetilde{\mathscr{V}_j}   \text{ } \equiv \text{ } \emptyset \text{ } \overset{(\textit{1})^{\prime}}{\Longleftrightarrow}       \bigg[       \underset{k^{\prime\prime}(i+1)}{\bigcup}  \text{ }           \bigg[ \text{ }     \widetilde{(\gamma_{i+1})_L}          \cap           \widetilde{(\gamma_{i+1})_R}        \text{ } \bigg]     \cap \mathcal{T}\mathcal{B} \big(  \widetilde{\mathscr{V}_{i+1}}  \big)       \bigg]           \cap \widetilde{\mathscr{V}_j}      \equiv \text{ } \emptyset \text{ } \\   \overset{(\textit{2})^{\prime}}{\Longleftrightarrow}         \bigg[  \text{ }    \underset{k^{\prime\prime}(i+1)}{\bigcup}  \text{ }           \bigg[ \text{ }     \widetilde{(\gamma_{i+1})_L}         \cap     \widetilde{(\gamma_{i+1})_R}      \bigg]     \cap \widetilde{\mathscr{V}_j}  \text{ }  \bigg]   \cap  \bigg[ \ \mathcal{T}\mathcal{B} \big(  \widetilde{\mathscr{V}_{i+1}}   \big)     \cap      \widetilde{\mathscr{V}_j}\bigg]                     \text{ }  \equiv \text{ } \emptyset \text{ } \\   \overset{(\textit{3})^{\prime}}{\Longleftrightarrow}      \bigg[    \underset{k^{\prime\prime}(i+1)}{\bigcup}    \widetilde{(\gamma_{i+1})_L}         \cap     \widetilde{\mathscr{V}_j}        \bigg]  \text{ }        \cap \text{ }  \bigg[  \text{ }    \underset{k^{\prime\prime}(i+1)}{\bigcup}    \text{ }   \widetilde{(\gamma_{i+1})_R}         \cap     \widetilde{\mathscr{V}_j}         \text{ } \bigg]  \text{ } \cap \bigg[   \text{ } \mathcal{T}\mathcal{B} \big(  \widetilde{\mathscr{V}_{i+1}}   \big)     \cap      \widetilde{\mathscr{V}_j} \text{ } \bigg]      \equiv \emptyset  \text{ , } 
\end{align*}

\noindent which, upon inspection, is also empty because,

\begin{align*}
     \underset{k^{\prime\prime}(i+1)}{\bigcup}    \text{ }   \widetilde{(\gamma_{i+1})_L}      \cap     \widetilde{\mathscr{V}_j}  \equiv   \bigg[      \underset{k^{\prime\prime}(i)}{\bigcup}       \widetilde{(\gamma_{i})_L}   \bigg]     \cup   \widetilde{(\gamma_{i+1})_L}   \cap     \widetilde{\mathscr{V}_j}  \equiv        \emptyset  \text{ } \text{ , }
\end{align*}

\noindent and also because,

\begin{align*}
       \underset{k^{\prime\prime}(i+1)}{\bigcup}    \text{ }  \widetilde{(\gamma_{i+1})_R}         \cap     \widetilde{\mathscr{V}_j}  \equiv   \bigg[      \underset{k^{\prime\prime}(i)}{\bigcup}       \widetilde{(\gamma_{i})_R} \bigg]     \cup   \widetilde{(\gamma_{i+1})_R}   \cap     \widetilde{\mathscr{V}_j} \equiv \emptyset      \text{ } \text{ , } 
\end{align*}

\noindent hence demonstrating that (\textit{HYP}) holding for $i$ is equivalent to (\textit{HYP}) holding for the left and right boundaries, respectively given by $\widetilde{(\gamma_{k^{\prime\prime}})_L}$ and $\widetilde{(\gamma_{k^{\prime\prime}})_R}$. In the sequence of equivalences above, first, for (\textit{1}),$\widetilde{\mathscr{V}_j}$ has nonempty intersection with the left and right boundaries, respectively given by the two conditions,

\begin{align*}
              \widetilde{(\gamma_{k^{\prime\prime}})_L}      \text{ }   \cap \text{ }         \widetilde{\mathscr{V}_j}             \equiv \emptyset            \text{ } \text{ , } 
\end{align*}

\noindent and by,

\begin{align*}
    \widetilde{(\gamma_{k^{\prime\prime}})_R}       \text{ }   \cap \text{ }   \widetilde{\mathscr{V}_j}         \equiv \emptyset                      \text{ } \text{ , } 
\end{align*}

\noindent respectively, for (\textit{2}) the equivalence provided in (\textit{1}) was rewritten by distributing the intersection with $\widetilde{\mathscr{V}_j}$, for (\textit{3}) the equivalence provided in (\textit{2}) was rewritten by distribution the intersection signs between either $ \widetilde{(\gamma_{k^{\prime\prime}})_L}$, $ \widetilde{(\gamma_{k^{\prime\prime}})_R}$, or $\widetilde{\mathscr{V}_j}$ as needed. The same reasoning is used in $(\textit{1})^{\prime}$, $(\textit{2})^{\prime}$ and $(\textit{3})^{\prime}$.

\bigskip

\noindent Next, observe, from arguments provided in {\color{blue}[11]}, that the proportion of domains is strictly less than $\frac{11}{12}$, following a pigeonhole estimate for the total number of faces that are crossed by $\mathrm{Bottom}_{x,y}$, holding $y$ fixed while $x$ varies between $1$ and $n$, that is dependent upon $\delta$, for some $\delta >0$, which takes the form,

\begin{align*}
\lfloor  \delta r \rfloor \text{ }   \equiv 12  \lfloor \frac{N}{n} \rfloor \equiv 12       \lfloor    \eta r     \rfloor  \text{ }   >       \lfloor \delta^{\prime} r^{\prime} \rfloor         \overset{(\textit{*})}{\equiv } N^{\prime} \lfloor   \frac{N}{n^{\prime}}      \rfloor  \overset{(\textit{**})}{\equiv}     N^{\prime} \lfloor       \eta r^{\prime}      \rfloor                        \text{ } \text{ , } 
\end{align*}

\noindent where the parameters $\delta^{\prime}$ and $r^{\prime}$ were previously introduced in Properties (\textit{XII}), (\textit{XIII}), and (\textit{XIV}) of cylindrical \textit{symmetric} domains, where in (\textit{*}), pick $N^{\prime}$ sufficiently large so that $N^{\prime} > 12$, while in (\textit{**}), obtain the final expression for $\lfloor \delta^{\prime} r^{\prime} \rfloor$ upon substituting in for the parameter $n^{\prime}$ for scaling parameters introduced over the cylinder. Under such a choice of parameters, the same property holds for $\mathrm{Top}_{x,y}$.

\bigskip

\noindent With the above results, denoting the set of quads previously described above that have at most $\lfloor \delta^{\prime} r^{\prime} \rfloor$ of which are \textit{tight} with $\mathscr{T}$, observe,

\begin{align*}
 \big|    \mathscr{T}     \big|  < \frac{n}{6}  \text{ , } 
\end{align*}

\noindent from which a bound on the number of triplet quads, of the form,

\begin{align*}
  \underset{j \in \textbf{Z}: j \equiv 0 \text{ } \mathrm{mod} 3}{\underset{i-1 \leq j \leq i+1}{\bigcup}} \widetilde{\mathscr{V}_j}  \text{ } \text{ , } 
\end{align*}

\noindent which, denoted with $\mathscr{T}_{i-1,i,i+1}$, satisfies,

\begin{align*}
 \big|     \mathscr{T}_{i-1,i,i+1}  \big|  <  3 \text{ } \big|  \mathscr{T}  \big|    <   \frac{n}{2}         \text{ , } 
\end{align*}

\noindent from which we conclude the argument. \boxed{}

\bigskip

\noindent The following items below are put to use in future arguments. Arguments making use of the definition below demonstrate that cylindrical symmetric domains under sloped boundary condition introduced for quantifying the free energy of the six-vertex model do not satisfy the modulo $6$ property. Despite the fact that the sloped symmetric cylindrical domains is not valid for the annulus supported from $6r$ to $12r$, with positive probability there exists a parameter that one can take sufficiently large, which is strictly larger than $12r$, for which a similar bound to that provided in {\color{blue}[11]} holds for a sloped boundary energy functional. It continues to remain of interest to explore new applications of the analysis of the free energy of the six-vertex model to other models, in which arguments initially presented for analysis of the free energy of the six-vertex model under sufficiently flat boundary conditions. Following the remaining results in \textit{5} before introducing the Ashkin-Teller model in \textit{6}, a significant difficulty emerges from the encoding of boundary conditions for the Ashkin-Teller model. Exploring whether other models with Hamiltonians similar to that of the $\big(q_{\sigma} , q_{\tau} \big)$-spin model, which consist of one-point and two-point interactions, can classified, in the sense that there boundaries in the phase diagram can be conjectured, in addition to qualitative differences between phase diagram sectors, is of great interest.

\bigskip

\noindent \textbf{Definition} \textit{11} (\textit{sloped ridge events from vertical $\mathrm{x}$-crossings of }$\mathrm{Top}_{x,y}$ to $\mathrm{Bottom}_{x,y}$ \textit{in} $\bar{\mathscr{V}_{i^{\prime}}}$). Denote,

\begin{align*}
  \mathrm{Sloped \text{ } Ridge\text{ }  Event} \text{ } \equiv \text{ } \mathcal{R}^{\mathrm{Sloped}}_{x,y} \text{ } \equiv \text{ }  \mathcal{R}^{\mathrm{Sloped}} \text{ } \equiv \text{ }  \bigg\{ \text{ }    \mathrm{Top}_{x,y} {\underset{\bar{\mathscr{V}_{i^{\prime}}}}{\overset{h\geq ck}{\longleftrightarrow}}}_{\text{ } \mathrm{x}} \mathrm{Bottom}_{x,y}      \text{ } \bigg\}                 \text{ } \text{ , } 
\end{align*}

\noindent as the $\mathrm{x}$-crossing event within $\bar{\mathscr{V}_{i^{\prime}}}$ for which $h \geq ck$, which was previously introduced following the statement of \textit{Property XI} in \textit{4.3.1}.

\bigskip

\noindent \textbf{Definition} \textit{12} (\textit{sloped fence events from paths along the boundary of the interior of the cylinder}). Denote,

\begin{align*}
   \mathrm{Sloped \text{ } Fence\text{ }  Event} \text{ } \equiv \text{ }  \big( \text{ }    \mathcal{F}^{\text{ } \mathrm{Sloped}}_{x,y}\text{ } \big)_h \text{ } \equiv \text{ }      \mathcal{F}^{\text{ } \mathrm{Sloped}}_h \text{ } \equiv \text{ } \underset{j^{\prime} \in \{  i^{\prime}-1 , i^{\prime}+1 \} }{\bigcup} \text{ }  \bigg\{ \text{ }      \widetilde{(\gamma_{j^{\prime}})}_L     \underset{\mathscr{O}\text{ } \cap\text{ } \mathrm{Slice}_{i^{\prime}-1,i^{\prime},i^{\prime}+1} \text{ } \cap \text{ }( \bar{\mathscr{V}_{j^{\prime}}} )^c }{\overset{h \geq ck}{\not\longleftrightarrow}}    \widetilde{(\gamma_{i^{\prime}})}_R               \text{ } \bigg\}         \text{ }   \text{ , } 
\end{align*}

\noindent where the \textit{disconnectivity} event defined above relates to the condition that, with positive probability,

\begin{align*}
       \textbf{P}^{\xi^{\mathrm{Sloped}}}_{\mathscr{O}} \bigg[ \text{ }     \widetilde{(\gamma_{i^{\prime}-1})}_L      \underset{\bar{\mathscr{V}_{i^{\prime}-1}} \text{ } }{\overset{h \geq ck}{\not\longleftrightarrow}}     \widetilde{(\gamma_{i^{\prime}-1})}_R     \text{ } \bigg]             \text{ } \text{ , } 
\end{align*}

\noindent occurs, in addition to,

\begin{align*}
    \textbf{P}^{\xi^{\mathrm{Sloped}}}_{\mathscr{O}} \bigg[ \text{ }     \widetilde{(\gamma_{i^{\prime}+1})}_L  \underset{\bar{\mathscr{V}_{i^{\prime}+1}}\text{ } }{\overset{h \geq ck}{\not\longleftrightarrow}}    \widetilde{(\gamma_{i^{\prime}+1})}_R   \text{ } \bigg]         \text{ } \text{ , } 
\end{align*}

\noindent occurring with positive probability, for $y \text{ } \equiv \text{ } 2 \text{ } \mathrm{mod} \text{ } 3$. 

\bigskip

\noindent \textbf{Definition} \textit{13} (\textit{union of faces outside those contained within each} $\gamma_{i^{*}}$). Denote,

\begin{align*}
        \mathcal{X}^{\text{ } \mathrm{Sloped}} \equiv \underset{\mathscr{F}^{\prime\prime}}{\underset{\mathscr{F}^{\prime}}{\bigcup}} \text{ }  \bigg\{      \forall \text{ }  \gamma_{i^{*}} , \exists \text{ }    \mathscr{F}^{\prime} ,  \mathscr{F}^{\prime\prime} \subset  \bigg[  \big(  F \big(      \gamma_{i^{*}}  \big) \big)^c  \cap F \big(     \mathscr{O} \cap \mathrm{Slice}_{i^{\prime}-1,i^{\prime},i^{\prime}+1}  \big)     \bigg]     :     \textbf{P}^{\xi^{\mathrm{Sloped}}}_{\mathscr{O}} \big[     \mathscr{F}^{\prime}    \overset{h < ck}{\longleftrightarrow}   \mathscr{F}^{\prime\prime}   \big] \text{ } \geq \text{ }  0  \bigg\}    \text{ } \text{ , } 
\end{align*}

\noindent as the union of faces outside of each $\gamma_{i^{*}}$, where the complementary faces in $\mathscr{O} \cap \mathrm{Slice}_{y-1,y,y+1}$ to each $\gamma_{i^{*}}$ are,

\begin{align*}
     \big(  F \big(      \gamma_{i^{*}}  \big) \big)^c \equiv    \big\{   F \in F \big(  \mathscr{O} \cap \mathrm{Slice}_{i^{\prime}-1,i^{\prime},i^{\prime}+1} \big)   :   F \notin F \big(  \gamma_{i^{*}}  \big)      \big\}            \text{ } \text{ , } 
\end{align*}

\noindent where $F \big( \text{ } \gamma_{i^{*}} \text{ } \big)$ is the subset of faces intersecting $\gamma_{i^{*}}$,

\begin{align*}
    F \big( \gamma_{i^{*}}  \big) \equiv      \big\{    F \in F \big(  \mathscr{O} \cap  \mathrm{Slice}_{i^{\prime}-1,i^{\prime},i^{\prime}+1} \big) :     F \cap \gamma_{i^{*}} \neq \emptyset      \big\}      \text{ } \text{ . } 
\end{align*}        

\bigskip

\noindent From the definitions above, the result below states that the conditional event of obtaining a sloped fence, given the existence of $\mathcal{R}^{\mathrm{Sloped}}$, is bound below by a strictly positive constant.

\bigskip

\noindent \textbf{Lemma} \textit{4.5} (\textit{sloped fences from sloped ridges}). The conditional pushforward of a sloped fence event under $\textbf{P}^{\mathrm{Sloped}}_{\mathscr{O}} \big[ \text{ } \cdot \text{ } \big]$,

\begin{align*}
          \textbf{P}^{\mathrm{Sloped}}_{\mathscr{O}} \big[                   \mathcal{F}^{\text{ } \mathrm{Sloped}}\text{ }  \big| \text{ }  \mathcal{R}^{\text{ } \mathrm{Sloped}} \text{ } \cup \text{ }    \big\{  \mathcal{X}^{\text{ } \mathrm{Sloped}}   \equiv  X  \big\}      \big] \text{ } \geq  \text{ } \mathcal{C}_{\mathcal{F},\mathcal{R}}   \text{ } \text{ , } 
\end{align*}

\noindent is bound below with a strictly positive constant. In the conditioning, the other event taken into union with $\mathcal{R}^{\text{ } \mathrm{Sloped}}$, as provided in \textbf{Definition} \textit{12}, represents the union of faces, excluding those contained within $\mathcal{I} \big( \text{ } \mathscr{O}\text{ } \big)$, for which a path exists in $\big(  F \big(      \gamma_{i^{*}}  \big) \big)^c  \cap F \big( \text{ }     \mathscr{O} \cap \mathrm{Slice}_{i^{\prime}-1,i^{\prime},i^{\prime}+1}     \text{ } \big)$. $\big(  F \big(      \gamma_{i^{*}}  \big) \big)^c$ is given as the subset of faces provided in the second part of \textbf{Definition} \textit{12}.

\bigskip

\noindent \textit{Proof of Lemma 4.5}. We make use of \textbf{Definition} \textit{10} and \textbf{Definition} \textit{11} given above. In summary, to argue that the conditional event for $\mathcal{F}^{\text{ } \mathrm{Sloped}}$ is uniformly positive, a realization of the environment from the faces in $\mathscr{O} \cap \mathrm{Slice}_{i^{\prime}-1,i^{\prime},i^{\prime}+1} $ must be fixed, which is denoted with $\mathrm{Env}$. Furthermore, the sloped version of $\mathcal{R}$ and of $\mathcal{F}$ differ when a slope is imposed from boundary conditions $\xi$ imposed upon the measure $\textbf{P}^{\xi} \big[ \text{ } \cdot \text{ } \big]$. To ensure that $\mathcal{R}^{\text{ } \mathrm{Sloped}}$ occurs with positive probability, as in previous arguments for \textbf{Lemma} \textit{4.1} which made use of a decomposition of crossing probabilities in terms of a summation over all realizations of \textit{symmetric} domains, in this case we similarly introduce a decomposition of $\mathrm{Env}$ in terms of the tuple $\big( \text{ }   \mathrm{Env}_k      \text{ }  ,  \text{ }     h|_{\mathrm{Env}_k}      \text{ } \big)$, for some $k>0$, which denotes the set of faces in $F \big( \text{ } \mathscr{O} \cap \mathrm{Slice}_{i^{\prime}-1,i^{\prime},i^{\prime}+1}  \text{ } \big)$ which have been explored in $\mathrm{Env}_k$, and the corresponding union of the values of the graph homomorphism for the height function when restricted to each $\mathrm{Env}_k$. Explicitly,

\begin{align*}
   h|_{\mathrm{Env}_k}   \equiv  \big\{  h \subset F \big(  \mathscr{O} \cap \mathrm{Slice}_{i^{\prime}-1,i^{\prime},i^{\prime}+1}  \big) :  h \cap \mathrm{Env}_k \neq \emptyset      \big\}   \text{ } \text{ , } 
\end{align*}

\noindent For $\mathrm{Env}_k \subset \mathrm{Env}$, each subset indexed in $k$ has the form,

\begin{align*}
     \mathrm{Env}_k  \equiv  \bigg\{    \emptyset \neq \mathscr{F}_k \text{ } \subset \text{ }    F \big(  \mathscr{O} \cap \mathrm{Slice}_{i^{\prime}-1,i^{\prime},i^{\prime}+1}   \big) :  \text{ }    \underset{k}{\bigcup} \text{ } \mathscr{F}_k         \text{ }  \subset \text{ }  \bigg[      \underset{i^{\prime} \in \textbf{N}}{\bigcup}  \text{ } \big( \text{ } \bar{\mathscr{V}_{i^{\prime}}} \text{ } \big)^c      \bigg] \text{ } \cap \text{ } \bigg[     \underset{y^{\prime} \equiv 0 \text{ } \mathrm{mod} \text{ } 3}{\underset{y^{\prime} \in \textbf{N}}{\bigcup}}  \mathrm{Slice}_{y^{\prime}-1,y^{\prime},y^{\prime}+1}     \bigg]     \bigg\}        \text{ , }
\end{align*}

\noindent where, for each $k$,

\begin{align*}
   \mathrm{Env}_k \text{ } \subset \text{ }              \mathrm{Env}_{k+1}           \text{ } \text{ . } 
\end{align*}

\noindent The family of subset, for $k\equiv 0$, is also subject to the condition,

\begin{align*}
      \big|  \mathrm{Env}_0      \big| \equiv 1     \text{ } \text{ , } 
\end{align*}

\bigskip

\noindent which, in words, states the the beginning of the exploration process to determine $\mathrm{Env}$, from the sequence of subsets contained within $\mathrm{Env}$ for every $k$, initially begins with the exploration of a single face within $F \big( \text{ } \mathscr{O} \cap \mathrm{Slice}_{i^{\prime}-1,i^{\prime},i^{\prime}+1}  \text{ } \big)$. Hence it suffices to demonstrate that $\mathcal{C}_{\mathcal{F},\mathcal{R}}$ is of the following form, in which we begin with the probability under $\textbf{P}^{\xi^{\mathrm{Sloped}}}_{\mathscr{O}} \big[  \cdot \big]$,

\begin{align*}
    \textbf{P}^{\xi^{\mathrm{Sloped}}}_{\mathscr{O}} \bigg[                  \mathcal{F}^{\text{ } \mathrm{Sloped}}_{(1-c)k+1} \bigg|  \mathcal{R}^{\text{ } \mathrm{Sloped}}  \cup   \big\{   \mathcal{X}^{\text{ } \mathrm{Sloped}}    \equiv    X   \big\}     \bigg]   \equiv  \underset{  \text{ } \mathrm{Env}_k \subset \mathrm{Env} \text{ } \forall k \text{ } }{\underset{h|_{\mathrm{Env}_k} \subset h|_{\mathscr{O}}}{\underset{(    \mathrm{Env}_k     \text{ }  , \text{ }   h|_{\mathrm{Env}_k } )}{\sum}} }  \mathscr{P}^{\prime} \text{ } \textbf{P}^{\xi^{\mathrm{Sloped}}}_{\mathscr{O}} \bigg[     \big\{   \mathrm{Env} \equiv \big(  \mathrm{Env}_k   , h|_{\mathrm{Env}_k } \big)          \big\}    \\ \cup    \big\{   \mathcal{X}^{\text{ } \mathrm{Sloped}}     \equiv    X  \big\}   \bigg]       \text{ } \text{ , } 
\end{align*}

\noindent where the restriction of the height function over the cylinder, $ h|_{\mathscr{O}}$, strictly contains the restrict height function over each $\mathrm{Env}_k$, $\forall \text{ } k$, and,

\begin{align*}
  \mathscr{P}^{\prime} \equiv   \textbf{P}^{\xi^{\mathrm{Sloped}}}_{\mathscr{O}} \bigg[                   \mathcal{F}^{\text{ } \mathrm{Sloped}}_{(1-c)k+1} \text{ }  \bigg|  \text{ } \big\{ \text{ }   \mathrm{Env} \equiv \big( \text{ } \mathrm{Env}_k   \text{ } , \text{ } h|_{\mathrm{Env}_k } \text{ } \big)              \text{ } \big\}   \text{ } \cup \text{ }    \big\{ \text{ }  \mathcal{X}^{\text{ } \mathrm{Sloped}}    \text{ } \equiv \text{ }   X   \text{ } \big\}      \bigg]  \text{ } \text{ , } 
\end{align*}

\noindent is an abbreviation for the conditional \textit{sloped fence event}. From such a decomposition, with the definition of the \textit{sloped fence event} given in \textbf{Definition} \textit{11}, instead of considering the probability of a vertical $\mathrm{x}$-crossing of the height function in $\bar{\mathscr{V}_{i^{\prime}}}$, the conditional probability,

\begin{align*}
      \textbf{P}^{\xi^{\mathrm{Sloped}}}_{\mathscr{O}} \bigg[  \big( \text{ }    \mathcal{F}^{\text{ } \mathrm{Sloped}}_{(1-c)k+1}        \text{ } \big)^{*} \text{ }   \bigg| \text{ }     \big\{ \text{ }   \mathrm{Env} \equiv \big( \text{ } \mathrm{Env}_k   \text{ } , \text{ } h|_{\mathrm{Env}_k } \text{ } \big)              \text{ } \big\}   \cup   \big\{ \text{ }  \mathcal{X}^{ \mathrm{Sloped}}   \equiv    X   \text{ } \big\}   \bigg]     \text{ } \text{ , } 
\end{align*}

\noindent of obtaining a dual horizontal crossing, where,

\begin{align*}
  \big(    \mathcal{F}^{\text{ } \mathrm{Sloped}}_l      \big)^{*}  \equiv \text{ }     \underset{j^{\prime\prime} \in \{  i^{\prime}-1 , i^{\prime}+1 \} }{\bigcup }  \bigg\{      \widetilde{(\gamma_{i^{\prime}})}_L     \underset{\mathscr{O}\text{ } \cap\text{ } \mathrm{Slice}_{i^{\prime}-1,i^{\prime},i^{\prime}+1} \text{ } \cap \text{ } \big( \text{ } \bar{\mathscr{V}_{j^{\prime\prime}}}\text{ } \big)^c }{\overset{h < l }{\longleftrightarrow}}    \widetilde{(\gamma_{j^{\prime\prime}})}_R     \bigg\}               \text{ , } 
\end{align*}

\noindent for some $l < ck$, and strictly positive product indices $j^{\prime\prime}>0$, which is a restatement of the event $\mathcal{F}^{\text{ } \mathrm{Sloped}}$ of $\textbf{Definition}$ \textit{12}, as the two crossings,

\begin{align*}
     \textbf{P}^{\xi^{\mathrm{Sloped}}}_{\mathscr{O}} \bigg[      \widetilde{(\gamma_{i^{\prime}-1})}_L     \underset{\bar{\mathscr{V}_{i^{\prime}-1}} }{\overset{h < l }{\longleftrightarrow}}    \widetilde{(\gamma_{i^{\prime}-1})}_R     \bigg]      \text{ } \text{ , } 
\end{align*}

\noindent and,

\begin{align*}
        \textbf{P}^{\xi^{\mathrm{Sloped}}}_{\mathscr{O}} \bigg[    \widetilde{(\gamma_{i^{\prime}+1})}_L     \underset{\bar{\mathscr{V}_{i^{\prime}+1}} }{\overset{h < l }{\longleftrightarrow}}    \widetilde{(\gamma_{i^{\prime}+1})}_R    \bigg]       \text{ } \text{ , } 
\end{align*}

\noindent given $i \text{ } \neq \text{ } 2 \text{ } \mathrm{mod} \text{ } 3$. The two crossings between the left and right boundaries instead both occur with positive probability, instead of the \textit{disconnectivity} event for $h \geq ck$. As a result, first,

\begin{align*}
 \textbf{P}^{\xi^{\mathrm{Sloped}}}_{\mathscr{O}} \bigg[  \big( \text{ }    \mathcal{F}^{\text{ } \mathrm{Sloped}}_{(1-c)k+1}         \text{ } \big)^{*}   \bigg|    \big\{  \mathrm{Env} \equiv \big( \text{ } \mathrm{Env}_k   \text{ } , \text{ } h|_{\mathrm{Env}_k } \text{ } \big)              \big\}  \cup \big\{  \mathcal{X}^{\text{ } \mathrm{Sloped}} \equiv   X  \big\}     \bigg] \text{ } \\ \overset{( \textbf{Corollary}\text{ }  \textit{1.2})}{\geq} \text{ }   \textbf{P}^{\xi^{\mathrm{Sloped}}}_{\bar{\mathscr{V}_{i^{\prime}}}} \bigg[ \text{ }  \big( \text{ }    \big(    \mathcal{F}^{\text{ } \mathrm{Sloped}}        \big)^{*}  \text{ } \big)_1 \text{ }  \bigg|    \text{ }   \big\{  \mathrm{Env} \equiv \big( \mathrm{Env}_k   \text{ } , \text{ } h|_{\mathrm{Env}_k } \big)              \big\}  \cup \big\{  \mathcal{X}^{\text{ } \mathrm{Sloped}}    \text{ } \equiv \text{ }   X  \big\}               \text{ } \bigg]   \text{ } \text{ , }  \tag{\textit{Env lower bound}}
\end{align*}

\noindent because $\bar{\mathscr{V}_{i^{\prime}}} \subset \mathscr{O}$, where in the lower bound probability above, the modified \textit{sloped fence} event is,

\begin{align*}
     \big( \text{ }    \big(    \mathcal{F}^{\text{ } \mathrm{Sloped}}        \big)^{*}  \text{ } \big)_1  \text{ } \equiv \text{ }               \big(    \mathcal{F}^{\text{ } \mathrm{Sloped}}_{(1-c)k+1 - j }        \big)^{*}          \text{ } \text{ , } 
\end{align*}

\noindent implying,

\begin{align*}
  (\textit{Env lower bound}) \text{ } \geq \text{ } \textbf{P}^{\xi^{\mathrm{Sloped}}}_{\bar{\mathscr{V}_{i^{\prime}}}} \big[     \big( \text{ }    \mathcal{F}^{\text{ } \mathrm{Sloped}}_{(1-c)k+1 - j }        \text{ } \big)^{*}    \big]   \text{ } \text{ . } 
\end{align*}

\noindent From previously used methods, if an $\mathrm{x}$-crossing across the strip or cylinder is given, then the complementary crossing event can be formulated by imposing that the value of the height function along all faces in the path lie in the complementary values that $h$ can taken, in addition to the complementary event for a horizontal (resp. vertical) crossing being a vertical (resp. horizontal) crossing. Altogether, one lower bound for the pushforward of $\big\{ \text{ }        \big( \text{ }    \mathcal{F}^{\text{ } \mathrm{Sloped}}_{(1-c)k+1 - j }        \text{ } \big)^{*}     \text{ } \big\}$ is,

\begin{align*}
   \textbf{P}^{\xi^{\mathrm{Sloped}}}_{\bar{\mathscr{V}_{i^{\prime}}}} \big[       \big( \text{ }    \mathcal{F}^{\text{ } \mathrm{Sloped}}_{(1-c)k+1 - j }        \text{ } \big)^{*}          \big] \text{ } \geq  1 -       \textbf{P}^{\xi^{\mathrm{Sloped}}}_{\bar{\mathscr{V}_{i^{\prime}}}} \bigg[    \mathrm{Top}_{x,y} \underset{\bar{\mathscr{V}_{i^{\prime}}}}{\overset{h \geq (1-c)k+1-j}{\longleftrightarrow}    }     \mathrm{Bottom}_{x,y}   \bigg]         \text{ } \text{ , } 
\end{align*}

\noindent in which the probabilistic quantity in the lower bound above for a crossing between $\mathrm{Top}_{x,y}$, and between $\mathrm{Bottom}_{x,y}$, is quantified under $\xi^{\mathrm{Sloped}}$ across $\bar{\mathscr{V}_{i^{\prime}}}$. Moreover,

\begin{align*}
    \textbf{P}^{\xi^{\mathrm{Sloped}}}_{\bar{\mathscr{V}_{i^{\prime}}}} \bigg[    \mathrm{Top}_{x,y} \underset{\bar{\mathscr{V}_{i^{\prime}}}}{\overset{h \geq (1-c)k+1-j}{\longleftrightarrow}    }     \mathrm{Bottom}_{x,y}  \bigg] \overset{(\textbf{Proposition}\textit{1.1})}{\leq}    \textbf{P}^{\big( \text{ } \xi^{\mathrm{Sloped}}\text{ } \big)^{\prime}}_{\textbf{Z} \times \big[i^{\prime} ,\big(  i^{\prime}+1   \big) r\big]} \bigg[    \mathrm{Top}_{x,y} \underset{\bar{\mathscr{V}_{i^{\prime}}}}{\overset{h \geq (1-c)k+1-2j}{\longleftrightarrow}    }     \mathrm{Bottom}_{x,y}   \bigg]     \text{ } \text{ , } 
\end{align*}

\noindent in which the upper bound above to the probability measure supported over $\bar{\mathscr{V}_{i^{\prime}}}$, under $\big( \text{ } \xi^{\mathrm{Sloped}}\text{ } \big)^{\prime}$, holds because $\bar{\mathscr{V}_{i^{\prime}}} \subset \big( \text{ } \textbf{Z} \times \big[i^{\prime} ,\big(  i^{\prime}+1   \big) r\big]\text{ } \big)$. Also,

\begin{align*}
 \textbf{P}^{\big( \text{ } \xi^{\mathrm{Sloped}}\text{ } \big)^{\prime}}_{\textbf{Z} \times \big[i^{\prime} ,\big(  i^{\prime}+1   \big) r\big]} \bigg[    \mathrm{Top}_{x,y} \underset{\bar{\mathscr{V}_{i^{\prime}}}}{\overset{h \geq (1-c)k+1-2j}{\longleftrightarrow}    }     \mathrm{Bottom}_{x,y}   \bigg]     \leq     \textbf{P}^{\big( \text{ } \xi^{\mathrm{Sloped}}\text{ } \big)^{\prime}}_{\textbf{Z} \times \big[i^{\prime} ,\big(  i^{\prime}+1   \big) r\big]} \bigg[\text{ }    \mathrm{Top}_{x,y} {\overset{h \geq (1-c)k+1-2j}{\longleftrightarrow}    }     \mathrm{Bottom}_{x,y}   \text{ } \bigg]        \text{ } \text{ , } 
\end{align*}

\noindent in which removing the conditioning on the crossing event,

\begin{align*}
     \bigg\{     \mathrm{Top}_{x,y} {\overset{h \geq (1-c)k+1-2j}{\longleftrightarrow}    }     \mathrm{Bottom}_{x,y}   \bigg\}     \text{ } \text{ , } 
\end{align*}

\noindent to not depend upon the connected components being contained within $\bar{\mathscr{V}_{i^{\prime}}}$ makes the probability of crossing in the upper bound more likely to occur. $\big( \text{ } \xi^{\mathrm{Sloped}} \text{ } \big)^{\prime}$ denoted the sloped boundary conditions on the probability measure $\textbf{P}^{\xi}_{\textbf{Z} \times \big[ i^{\prime} , \big( i^{\prime}+1 \big) \text{ } r\big] } \big[ \cdot \big]$. Under this measure, the values of the height function, along the boundary segment $\partial \text{ } \textbf{Z}\text{ }  \times \text{ } \big[ \text{ } i^{\prime} , \big( \text{ } i^{\prime}+1 \text{ } \big) r \big]$, are precisely the image of the graph homomorphism over each boundary face which would equal $\big( \text{ } \xi^{\mathrm{Sloped}} \text{ } \big)^{\prime} - \mathcal{C}_{f}$, in the case of the minimal height function for which one would expect $\big( \text{ } \xi^{\mathrm{Sloped}} \text{ } \big)^{\prime}$ to be an admissible boundary condition. The quantity $\mathcal{C}_f$ subtracted from $\big( \text{ } \xi^{\mathrm{Sloped}} \text{ } \big)^{\prime}$ in the previous sentence represents some constant that is dependent upon the total number of \textit{frozen} faces within each \textit{freezing cluster} that has a nonempty intersection with $\partial \big(  \textbf{Z}  \times  \big[ \text{ } i^{\prime} , \big( \text{ } i^{\prime}+1 \text{ } \big) r \big]\big)$.

\bigskip

\noindent Recall a previous result in the strip,

\begin{align*}
    \text{ } \textbf{P}^{\xi^{\mathrm{Sloped}}}_{\textbf{Z} \times [0,n^{\prime}N]}  \big[   [0 , \delta n]                             \text{ }   \times \{ 0 \}      \text{ } \overset{ h \geq (1-c) k }{ \longleftrightarrow } \text{ } [ i , i + \delta n ]      \text{ }  \times \{ n \}                    \big] \text{ } < \text{ } 1 - c  \text{ } \text{ , } 
\end{align*}

\noindent from the upper bound that is provided in \textbf{Proposition} \textit{1.2}. With the same probability measure supported over $\big( \text{ } \textbf{Z} \times \big[ i^{\prime} , \big( i^{\prime}+1 \big) \text{ } r\big]\text{ } \big)$,

\begin{align*}
    \textbf{P}^{\big( \text{ } \xi^{\mathrm{Sloped}}\text{ } \big)^{\prime}}_{\textbf{Z} \times \big[i^{\prime} ,\big(  i^{\prime}+1   \big) r\big]} \bigg[   \mathrm{Top}_{x,y} {\overset{h \geq (1-c)k+1-2j}{\longleftrightarrow}    }     \mathrm{Bottom}_{x,y}   \bigg]        <  1 - c           \text{ } \text{ , } 
\end{align*}

\noindent under the same choice of $\delta^{\prime\prime}$, which, as stated in \textit{Property XII}, and in \textit{Property XIII}, that the number of faces between the union $\bar{\mathscr{V}_{i^{\prime},y}}$ and either the top, or the bottom, of $\bar{\mathscr{V}_{i^{\prime},y}}$ can at most be equal to $\lfloor \text{ } \delta^{\prime\prime} r \text{ } \rfloor$. With such a $\delta^{\prime\prime}$,

\begin{align*}
    \textbf{P}^{\big( \text{ } \xi^{\mathrm{Sloped}}\text{ } \big)^{\prime}}_{\textbf{Z} \times \big[i^{\prime} ,\big(  i^{\prime}+1   \big) r\big]} \bigg[     \widetilde{(\gamma_{i^{\prime}})}_L     \underset{\bar{\mathscr{V}_{i^{\prime}}} }{\overset{h < (1-c_0)k }{\longleftrightarrow}}    \widetilde{(\gamma_{i^{\prime}})}_R       \text{ }      \bigg|   \text{ }    \big\{  \mathrm{Env} \equiv \big( \mathrm{Env}_k   \text{ } , \text{ } h|_{\mathrm{Env}_k } \big)              \big\}  \cup \big\{  \mathcal{X}^{\text{ } \mathrm{Sloped}}    \text{ } \equiv \text{ }   X  \big\}         \bigg] \geq c_0  \text{ } \text{ , } 
\end{align*}

\noindent in which the strictly positive constant $c_0$ lower bounds the conditional event dependent upon the connected components of $h$ from $\bigg\{  \widetilde{(\gamma_{i^{\prime}})}_L    \underset{\bar{\mathscr{V}_{i^{\prime}}} }{\overset{h < (1-c_0)k }{\longleftrightarrow}}    \widetilde{(\gamma_{i^{\prime}})}_R      \bigg\}$.

\bigskip

\noindent Before completing the argument, to lighten notation, set,

\begin{align*}
   \Gamma_{L,R,i^{\prime}}   \equiv \text{ } \bigg\{    \widetilde{(\gamma_{i^{\prime}})}_L     \underset{\bar{\mathscr{V}_{i^{\prime}}} }{\overset{h < (1-c_0)k }{\longleftrightarrow}}    \widetilde{(\gamma_{i^{\prime}})}_R       \bigg\}    \text{ }   \text{ , } 
\end{align*}

\noindent for the crossing event under the strip probability measure above, 

\begin{align*}
    \mathcal{E}\mathrm{n}  \equiv \text{ }    \big\{  \mathrm{Env} \equiv \big( \mathrm{Env}_k   \text{ } , \text{ } h|_{\mathrm{Env}_k } \big)              \big\}  \text{ } \text{ , } 
\end{align*}

\noindent for the decomposition of $\mathrm{Env}$, and,

\begin{align*}
   \mathrm{XS}  \equiv \text{ }    \big\{  \mathcal{X}^{\text{ } \mathrm{Sloped}}    \text{ } \equiv \text{ }   X  \big\}     \text{ } \text{ , } 
\end{align*}

\noindent for the clusters of the height function excluding those induced from the crossings of each $\gamma_{i^{*}}$ and $\gamma_{L+i^{*}}$, from which we conclude, given a suitable index set $\mathcal{J}^{\prime}$,

\begin{align*}
     \textbf{P}^{\big( \text{ } \xi^{\mathrm{Sloped}}\text{ } \big)^{\prime}}_{\textbf{Z} \times \big[i^{\prime} ,\big(  i^{\prime}+1   \big) r\big]} \bigg[    \Gamma_{L,R,i^{\prime}} \text{ } \big| \text{ }   \mathcal{E}\mathrm{n}  \cap      \mathrm{XS}    \bigg] \text{ }  \equiv \text{ }    \textbf{P}^{\big( \text{ } \xi^{\mathrm{Sloped}}\text{ } \big)^{\prime}}_{\textbf{Z} \times \big[i^{\prime} ,\big(  i^{\prime}+1   \big) r\big]} \bigg[   \underset{\text{countably many } j^{\prime},\text{ }  j^{\prime} < i}{\bigcap} \text{ }    \Gamma_{L,R,j^{\prime}}  \text{ } \big| \text{ }   \mathcal{E}\mathrm{n}  \cap      \mathrm{XS}    \bigg]    \text{ } \text{ , } 
\end{align*}

\noindent where $j^{\prime} \in \mathcal{J}^{\prime}$. We bound the ultimate term above, as,

\begin{align*}
\textbf{P}^{\big( \text{ } \xi^{\mathrm{Sloped}}\text{ } \big)^{\prime}}_{\textbf{Z} \times \big[i^{\prime} ,\big(  i^{\prime}+1   \big) r\big]} \bigg[ \underset{\text{countably many } j^{\prime},\text{ }  j^{\prime} < i}{\bigcap} \text{ }    \Gamma_{L,R,j^{\prime}} \text{ }  \text{ } \big| \text{ }   \mathcal{E}\mathrm{n}  \cap      \mathrm{XS}  \bigg]  \text{ } \overset{(\mathrm{FKG})}{\geq} \text{ }  \underset{j^{\prime} \in \mathcal{J}^{\prime}}{\prod} \textbf{P}^{\big( \text{ } \xi^{\mathrm{Sloped}}\text{ } \big)^{\prime}}_{\textbf{Z} \times \big[i^{\prime} ,\big(  i^{\prime}+1   \big) r\big]} \bigg[     \Gamma_{L,R,j^{\prime}}     \text{ }  \big| \text{ }   \mathcal{E}\mathrm{n}  \cap      \mathrm{XS}     \bigg]      \text{ } \\       \geq        c^{\text{ } | \mathcal{J}^{\prime} | \text{ }}_0   \text{ }      \text{ }  \text{ , } 
\end{align*}

\noindent in which there the $(\mathrm{FKG})$ inequality is applied to a collection of countably many crossing events between left and right boundaries, to the conditionally defined measure,

\begin{align*}
         \textbf{P}^{\big( \text{ } \xi^{\mathrm{Sloped}}\text{ } \big)^{\prime}}_{\textbf{Z} \times \big[i^{\prime} ,\big(  i^{\prime}+1   \big) r\big]} \bigg[      \cdot \text{ }         \big|   \text{ } \mathcal{E}\mathrm{n}  \cap      \mathrm{XS}      \bigg]    \text{ } \text{ , } 
\end{align*}

\noindent where,

\begin{align*}
 c_0 \equiv \text{ } \underset{j^{\prime} \in \mathcal{J}^{\prime}}{\mathrm{inf}} \text{ }  \bigg\{   \textbf{P}^{\big( \text{ } \xi^{\mathrm{Sloped}}\text{ } \big)^{\prime}}_{\textbf{Z} \times \big[i^{\prime} ,\big(  i^{\prime}+1   \big) r\big]} \bigg[     \Gamma_{L,R,j^{\prime}}     \text{ }  \big| \text{ }   \mathcal{E}\mathrm{n}  \cap      \mathrm{XS}   \bigg]  \bigg\}  \text{ } \text{ , }
\end{align*}

\noindent is obtained from the restricted event over the cylinder,

\begin{align*}
    \Gamma_{L,R,j^{\prime}}      \text{ } \equiv \text{ }  \underset{j^{\prime} \in \textbf{N}}{\bigcup}   \bigg\{             \widetilde{(\gamma_{j^{\prime}})}_L     \underset{\bar{\mathscr{V}_{i^{\prime}}} }{\overset{h < (1-c_0)k }{\longleftrightarrow}}    \widetilde{(\gamma_{j^{\prime}})}_R                \bigg\}         \text{ } \subset \text{ } \bigg\{    \widetilde{(\gamma_{i^{\prime}})}_L     \underset{\bar{\mathscr{V}_{i^{\prime}}} }{\overset{h < (1-c_0)k }{\longleftrightarrow}}    \widetilde{(\gamma_{i^{\prime}})}_R       \bigg\}      \text{ } \text{ , } 
\end{align*}

\noindent $\forall \text{ } j^{\prime}$ with some fixed $i^{\prime}$, in which for the first, and last, $j^{\prime}$, which we respectively denote with $j_0$ and $j_f$, appearing in the union of connectivity events between each $\widetilde{(\gamma_{j^{\prime}})}_L$ and $   \widetilde{(\gamma_{j^{\prime}})}_R$, satisfy,

\begin{align*}
           \widetilde{(\gamma_{j_0})}_L \text{ } \cap  \text{ }   \widetilde{(\gamma_{i^{\prime}})}_L    \text{ } \neq \emptyset \text{ } \text{ , } \\ \text{ }          \widetilde{(\gamma_{j_f})}_R \text{ } \cap  \text{ }   \widetilde{(\gamma_{i^{\prime}})}_R    \text{ } \neq \emptyset        \text{ } \text{ , } 
\end{align*}

\noindent hence yielding an intermediate lower bound $c^{\text{ } | \mathcal{J}^{\prime} | \text{ }}_0$ to $\mathcal{C}_{\mathcal{F} , \mathcal{R}}$.

\bigskip

 \noindent The collection of events $\bigg\{ \widetilde{(\gamma_{j^{\prime}})}_L     \underset{\bar{\mathscr{V}_{i^{\prime}}} }{\overset{h < (1-c_0)k }{\longleftrightarrow}}    \widetilde{(\gamma_{j^{\prime}})}_R          \bigg\}$ satisfies the containment,

\begin{align*}
        \bigg\{   \widetilde{(\gamma_{j^{\prime}})}_L     \underset{\bar{\mathscr{V}_{i^{\prime}}} }{\overset{h < (1-c_0)k }{\longleftrightarrow}}    \widetilde{(\gamma_{j^{\prime}})}_R       \bigg\}  \subset    \bigg\{   \widetilde{(\gamma_{i^{\prime}})}_L     \underset{\bar{\mathscr{V}_{i^{\prime}}} }{\overset{h < (1-c_0)k }{\longleftrightarrow}}    \widetilde{(\gamma_{i^{\prime}})}_R       \bigg\}  \text{ } \text{ , } 
\end{align*}

\noindent $\forall \text{ } j^{\prime}$. From the decomposition of $ \textbf{P}^{\xi^{\mathrm{Sloped}}}_{\mathscr{O}} \bigg[                \mathcal{F}^{\text{ } \mathrm{Sloped}}_{(1-c)k+1} \text{ } \big| \text{ } \mathcal{R}^{\text{ } \mathrm{Sloped}} \text{ } \cup \text{ }    \big\{ \text{ }  \mathcal{X}^{\text{ } \mathrm{Sloped}}    \text{ } \equiv \text{ }   X   \text{ } \big\}   \bigg]$ provided at the beginning of the proof, 

\begin{align*}
     \underset{  \text{ } \mathrm{Env}_k \subset \mathrm{Env} \text{ } \forall k \text{ } }{\underset{h|_{\mathrm{Env}_k} \subset h|_{\mathscr{O}}}{\underset{(    \mathrm{Env}_k     \text{ }  , \text{ }   h|_{\mathrm{Env}_k } )}{\sum}} }  \mathscr{P}^{\prime} \text{ } \textbf{P}^{\xi^{\mathrm{Sloped}}}_{\mathscr{O}} \big[  \big\{ \text{ }   \mathrm{Env} \equiv \big( \text{ } \mathrm{Env}_k   \text{ } , \text{ } h|_{\mathrm{Env}_k } \text{ } \big)              \text{ } \big\}   \text{ } \cup  \text{ }    \big\{ \text{ }  \mathcal{X}^{\text{ } \mathrm{Sloped}}    \text{ } \equiv \text{ }   X   \text{ } \big\}    \big]         \geq    c^{\text{ } | \mathcal{J}^{\prime} | \text{ }}_0   \text{ } \text{ , } 
\end{align*}

\noindent also holds, with the identical lower bound dependent upon $c_0$ applying.

\bigskip

\noindent Finally, there exists some sufficiently large integer $X$, where $X < | \mathcal{J}^{\prime}|$, which is proportional to the number of $\bar{\mathscr{V}_{i^{\prime}}}$ that are needed for the crossings across each $\bigg\{ \widetilde{(\gamma_{j^{\prime}})}_L     \underset{\bar{\mathscr{V}_{i^{\prime}}} }{\overset{h < (1-c_0)k }{\longleftrightarrow}}    \widetilde{(\gamma_{j^{\prime}})}_R          \bigg\}$ to occur. The desired lower bound then takes the form,

\begin{align*}
      c^{\text{ } | \mathcal{J}^{\prime} | \text{ }}_0  \text{ } \geq \text{ }    c^{\text{ } X}       \text{ } \text{ , } 
\end{align*}

\noindent from which we conclude the argument, after having set $\mathcal{C}_{\mathcal{F},\mathcal{R}} \equiv c^{\text{ } X}$. \boxed{}

\bigskip

\noindent To proceed, in the following result a lower bound is also provided for a conditionally defined \textit{sloped ridge} event.

\bigskip

\noindent \textbf{Lemma} \textit{4.6} (\textit{probabilistic upper bound for conditional sloped ridge events}). For all $r>0$, and all $k$ sufficiently large, 

\begin{align*}
    \textbf{P}^{\xi^{\mathrm{Sloped}}}_{\mathscr{O}}      \bigg[                     \mathcal{R}^{\text{ } \mathrm{Sloped}} \big| \text{ }            \mathcal{F}^{\text{ } \mathrm{Sloped}}_{(1-c_0)k+1} \text{ } \cup \big\{ \text{ } \mathcal{X}^{\text{ } \mathrm{Sloped}}    \text{ } \equiv \text{ }   X   \text{ } \big\}\bigg]  \leq   \mathcal{C}_{\text{ } \mathcal{R},  \mathcal{F}\text{ } }      \text{ } \text{ , } 
\end{align*}

\noindent where the lower bound is a strictly positive constant dependent upon the probability,

\begin{align*}
       \textbf{P}^{\xi^{\mathrm{Sloped}}}_{\textbf{Z} \times [ -r + \delta_1 , 2r + \delta_1]}      \bigg[              [ 0 , \lfloor \delta^{\prime\prime} r  \rfloor ] \times \{0\} \underset{\textbf{Z} \times [0,r]}{\overset{h \geq c_0 k  - j}{\longleftrightarrow}}                  \textbf{Z} \times \{ n \}              \bigg]            \text{ } \text{ , } 
\end{align*}

\noindent given a $\delta_1$ sufficiently small.

\bigskip

\noindent \textit{Proof of Lemma 4.6}. In the conditioning for the probability defined in the statement, if $\mathcal{F}^{\text{ } \mathrm{Sloped}}_{(1-c_0)k+1}$ occurs with positive probability, then within $\bar{\mathscr{V}_{i^{\prime},y+1}}$, there exists a bottom most domain, which can be decomposed into the following four components, 

\begin{align*}
    \bar{\mathscr{V}_{i^{\prime},y+1}}  \text{ } \equiv \text{ } \underset{i^{\prime} : i^{\prime} \equiv 1 \text{ } \mathrm{mod} \text{ } 3 }{\bigcup} \text{ } {\underset{1 \text{ } \leq \text{ } j \text{ } \leq \text{ } 4}{\bigcap}}   \text{ }        \gamma_{j,i^{\prime}}   \text{ }\equiv \text{ } \underset{i^{\prime} : i^{\prime} \equiv 1 \text{ } \mathrm{mod} \text{ } 3 }{\bigcup} \text{ } \bigg\{ \text{ }     \underset{1,i^{\prime}}{\underset{\Updownarrow}{\underset{j\equiv 1}{\underbrace{\widetilde{(\gamma_{i^{\prime}})}_L}}}       }       \text{ } \cap \text{ }    \underset{2,i^{\prime}}{\underset{\Updownarrow}{\underset{j\equiv 2}{\underbrace{\widetilde{(\gamma_{i^{\prime}})}_R}  }}}     \text{ } \cap \text{ }    \underset{3,i^{\prime}}{\underset{\Updownarrow}{\underset{j\equiv 3}{\underbrace{ \mathcal{T} \big( \text{ } \widetilde{\mathscr{V}_{i^{\prime}}}} \text{ } \big)    }}}   \text{ } \cap \text{ }     \underset{4,i^{\prime}}{\underset{\Updownarrow}{\underset{j\equiv 4}{\underbrace{ \mathcal{B} \big( \text{ } \widetilde{\mathscr{V}_{i^{\prime}}} \text{ } \big)    }}}}  \text{ } \bigg\} \text{ }  \text{ , } 
\end{align*}

\noindent which exist with positive probability when $\bar{\mathscr{V}_{i^{\prime},y+1}} \cap \mathrm{Slice}_{y+1} \neq \emptyset$.

\bigskip

\noindent Proceeding, denote,

\begin{align*}
          P_m          \text{ } \equiv \text{ }   \textbf{P}^{\xi^{\mathrm{Sloped}}}_{[-m,m] \times [ -r + \delta_1 , 2r + \delta_1]}      \bigg[              [ 0 , \lfloor \delta^{\prime\prime} r  \rfloor ] \times \{0\} \underset{\textbf{Z} \times [0,r]}{\overset{h \geq c_0 k  - j}{\longleftrightarrow}}                 [-m,m] \times \{ n \}            \bigg]            \text{ } \text{ , } 
\end{align*}

\noindent as the probability of $\bigg\{ \text{ }   [ 0 , \lfloor \delta^{\prime\prime} r  \rfloor ] \times \{0\} \underset{\textbf{Z} \times [0,r]}{\overset{h \geq c_0 k  - j}{\longleftrightarrow}}                 [-m,m] \times \{ n \}   \text{ } \bigg\}$ occurring within $\textbf{Z} \times [0,r]$. To study the sequence of measures $P_m \longrightarrow P_{+\infty}$, introduce,

\begin{align*}
            \chi_B^{x,y} \equiv \text{ }  \underset{\mathscr{F}^{\prime}_2 \neq \mathscr{F}^{\prime}_1}{\underset{\mathscr{F}^{\prime}_2}{\underset{\mathscr{F}^{\prime}_1}{\bigcup}}} \text{ }  \bigg\{                 \mathscr{F}^{\prime}_1 , \mathscr{F}^{\prime}_2 \in F \big(  \mathscr{O} \cap \mathrm{Slice}_{i-1} \big) :   \textbf{P}^{\xi^{\mathrm{Sloped}}}_{\mathscr{O}} \big[ \text{ }      \mathscr{F}^{\prime}_1 \underset{\mathscr{V}_{i-1}}{\overset{h \leq (1-c_0)k}{\longleftrightarrow}} \mathscr{F}^{\prime}_2       \text{ } \big]      \text{ } > 0     \bigg\}             \text{ } \text{ , } 
\end{align*}

\noindent as the bottom, and top, boundaries,

\begin{align*}
              \chi_T^{x,y} \equiv  \underset{\mathscr{F}^{\prime\prime}_2 \neq \mathscr{F}^{\prime\prime}_1}{\underset{\mathscr{F}^{\prime\prime}_2}{\underset{\mathscr{F}^{\prime\prime}_1}{\bigcup}}}      \bigg\{         \mathscr{F}^{\prime\prime}_1 , \mathscr{F}^{\prime\prime}_2 \in F \big( \mathscr{O} \cap \mathrm{Slice}_{i+1}  \big) :    \textbf{P}^{\xi^{\mathrm{Sloped}}}_{\mathscr{O}} \big[ \text{ }      \mathscr{F}^{\prime\prime}_1 \underset{\mathscr{V}_{i+1}}{\overset{h \leq (1-c_0)k}{\longleftrightarrow}} \mathscr{F}^{\prime\prime}_2       \text{ } \big]      \text{ } > 0                   \bigg\}                 \text{ } \text{ , } 
\end{align*}

\noindent respectively, of some domain within $\bar{\mathscr{V}_{i^{\prime}}}$ so that the subdomain $\mathcal{D}_{i^{\prime},y} \subset \mathscr{V}_{i^{\prime}}$ satisfies the decomposition,

\begin{align*}
       \mathcal{D}_{i^{\prime},x,y} \text{ } \equiv \text{ }   \underset{1 \leq k^{\prime} \leq 4}{\bigcap}   \gamma_{k^{\prime}}    \text{ }  \equiv \text{ }  (\gamma_{i^{\prime}})_L         \text{ } \cap \text{ }   \chi_B^{x,y}     \text{ } \cap \text{ }  (\gamma_{i^{\prime}})_R         \text{ } \cap    \text{ }    \chi_T^{x,y}   \text{ } \text{ , } 
\end{align*}

\noindent is bound by the top and bottom boundaries $\chi^B_{x,y}$, and $\chi^T_{x,y}$, and the left and right boundaries from the definition of $\mathscr{V}_{i^{\prime}}$ which were introduced in \textit{4.3.1} following (\textit{Property V}). Irrespective of whether boundary conditions on the six-vertex probability measured are flat or sloped, the restriction $\big\{   \mathcal{D}_{i^{\prime},x,y} \text{ }           : (x,y) \text{ satisfy (\textit{Property XIV})}   \big\} \text{ }  \subset \text{ } \mathcal{D}_{i^{\prime},x,y}$, corresponding to admissible pairs $(x,y)$ satisfying the \textit{good} property, is dependent upon,

\begin{align*}
     h|_{\mathcal{D}_{i^{\prime},x,y}^c} \text{ } \equiv \text{ }   \big\{     \forall \text{ } i^{\prime} ,  \exists F_{i^{\prime}} \in F \big( \mathrm{Slice}_{y-1,y,y+1} \cap      \mathcal{D}_{i^{\prime},x,y}^c       \big)         :   h|_{F_{i^{\prime}}}   \in \textbf{Z}     \big\}         \text{ } \text{ , } 
\end{align*}

\noindent which is the value of the height function in the complementary region of $\mathcal{D}_{i^{\prime},x,y}^c \subset \mathscr{O}$, in addition to,

\begin{align*}
     h|_{\partial \mathcal{D}_{i^{\prime},x,y}} \text{ } \equiv \text{ } \big\{         \forall i^{\prime} ,  \exists F_{i^{\prime}} \in F \big( \text{ }  \mathcal{D}_{i^{\prime},x,y}^c  \cap      \mathcal{D}_{i^{\prime},x,y}       \text{ } \big)         :   h|_{F_{i^{\prime}}}   \in \textbf{Z}              \big\}         \text{ } \text{ , } 
\end{align*}

\noindent which is the value of the height function on the boundary $\partial \mathcal{D}_{i^{\prime},x,y} \subset \mathscr{O}$. Hence each realization of the height function is independent of another realization, and furthermore, within the interior of each such $\mathcal{D}_{i^{\prime},x,y}$, one realization of the height function is also independent of another realization of the height function. From the six-vertex probability measure $\textbf{P}^{\xi^{\mathrm{Sloped}}}_{\mathcal{D}_{i^{\prime},x,y}} \big[ \text{ } \cdot \text{ } \big]$ supported over ${\mathcal{D}_{i^{\prime},x,y}}$, given some \textit{good} $(x,y)$,

\begin{align*}
       \textbf{P}^{\xi^{\mathrm{Sloped}}}_{\mathcal{D}_{i^{\prime},x,y}} \bigg[ \text{ }      \mathrm{Top}_{x,y} \underset{\mathscr{V}_{i^{\prime}}}{\overset{h \geq k}{\longleftrightarrow}    }     \mathrm{Bottom}_{x,y}                  \text{ } \bigg] \text{ }    \equiv \text{ }    \textbf{P}^{\xi^{\mathrm{S}}}_{D} \bigg[ \text{ }                       \mathrm{Top}_{x,y} \underset{\mathscr{V}_{i^{\prime}}}{\overset{h \geq k}{\longleftrightarrow}    }     \mathrm{Bottom}_{x,y}     \text{ } \bigg]     \text{ } \text{ , } 
\end{align*}

\noindent after having introduced the abbreviations $\xi^{\mathrm{Sloped}}\equiv \xi^{\mathrm{S}}$, and $\mathcal{D}_{i^{\prime},x,y} \equiv D$, from which the probability of the crossing $\bigg\{ \text{ }       \mathrm{Top}_{x,y} \underset{\mathscr{V}_{i^{\prime}}}{\overset{h \geq k}{\longleftrightarrow}    }     \mathrm{Bottom}_{x,y}             \text{ } \bigg\}$ occurring can also be upper bounded with,

\begin{figure}
\begin{align*}
\includegraphics[width=1.05\columnwidth]{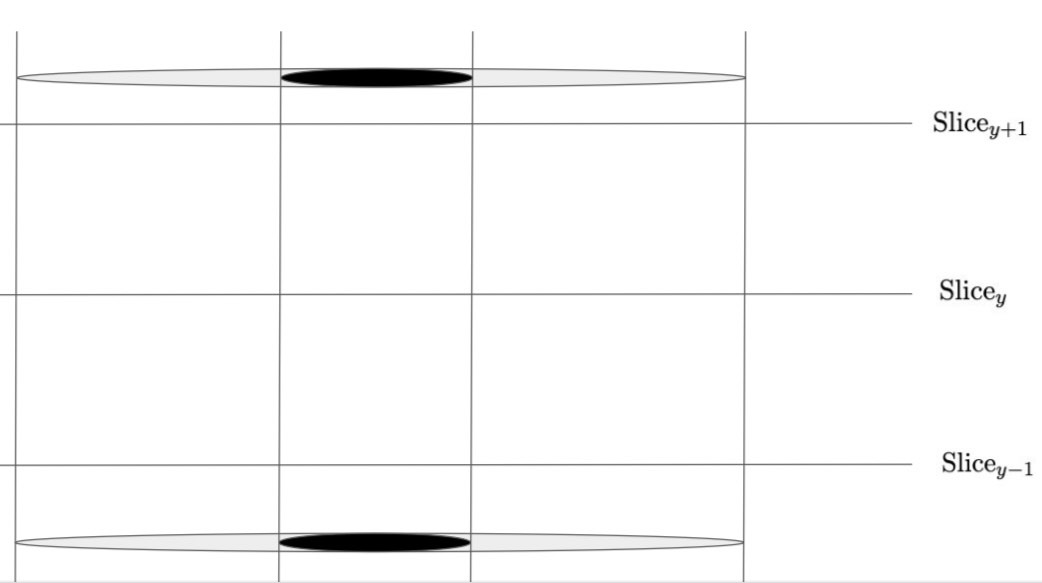}\\
\end{align*}
\caption{\textit{One cylindrical domain, with top and bottom regions depicted in black, embedded within} the top and bottom regions of $\mathscr{O}$ \textit{intersecting} $\mathrm{Slice}_{y-1,y,y+1}$.}
\end{figure}

\begin{align*}
    \textbf{P}^{\xi^{\mathrm{S}}}_{D} \bigg[ \text{ }                       \mathrm{Top}_{x,y} \underset{\mathscr{V}_{i^{\prime}}}{\overset{h \geq k}{\longleftrightarrow}    }     \mathrm{Bottom}_{x,y}     \text{ } \bigg]   \text{ } \overset{(\mathrm{CBC})}{\leq} \text{ }      \textbf{P}^{\big( \text{ } \xi^{\mathrm{S}}\text{ } \big)^{\prime}}_{D} \bigg[ \text{ }                       \mathrm{Top}_{x,y} \underset{\mathscr{V}_{i^{\prime}}}{\overset{h \geq k-j}{\longleftrightarrow}    }     \mathrm{Bottom}_{x,y}     \text{ } \bigg]  \text{ } \tag{\textit{4.5.1}} \\   \underset{\big( \text{ } \xi^{\mathrm{S}}\text{ } \big)^{\prime} \leq \big( \text{ } \xi^{\mathrm{S}}\text{ } \big)^{\prime\prime}}{\overset{(\textbf{Corollary}\text{ } \textit{1.2})}{\leq}}     \textbf{P}^{\big( \text{ } \xi^{\mathrm{S}}\text{ } \big)^{\prime\prime}}_{[-m,m] \times [ (y-1)r , (y+2) r]} \bigg[ \text{ }                       \mathrm{Top}_{x,y} \underset{\mathscr{V}_{i^{\prime}}}{\overset{h \geq k-2j}{\longleftrightarrow}    }     \mathrm{Bottom}_{x,y}     \text{ } \bigg]                \text{ } \tag{\textit{4.5.2}} \\  \text{ } \overset{(***)}{\leq} \text{ }            P_m \\ \overset{m \longrightarrow + \infty}{\leq} \text{ }              P_{+\infty}    \text{ } \text{ , }\tag{\textit{4.5.3}}
\end{align*}

\noindent where in the last term, the corresponding probability pushed forwards under $\textbf{P}^{\big( \text{ } \xi^{\mathrm{S}}\text{ } \big)^{\prime\prime}}_{[-m,m] \times [ (y-1)r , (y+2) r]} \big[    \cdot \big]$ is given by,

\begin{align*}
 \bigg\{            [ 0 , \lfloor \delta^{\prime\prime} r  \rfloor ] \times \{0\} \underset{\textbf{Z} \times [0,r]}{\overset{h \geq c_0 k  - j}{\longleftrightarrow}}              \textbf{Z} \times \{ n \}             \bigg\}            \text{ } \text{ , } 
\end{align*}

\noindent given sloped boundary conditions for which $\big( \text{ } \xi^{\mathrm{S}} \text{ } \big)^{\prime} \geq \xi^{\mathrm{S}}$, which holds by $(\mathrm{CBC})$, followed by an application of \textbf{Corollary} \textit{1.2}, given two finite volumes satisfying $\Lambda^{\prime} \equiv D \subset \big( \text{ }  \textbf{Z} \times [(y-1)r , (y+2)r] \text{ } \big) \equiv \Lambda$, and also that $\big( \text{ } \textbf{Z} \times [(y-1)r , (y+2)r] \text{ } \big) \cap D \equiv D \neq \emptyset$. Finally, in $(***)$, an upper bound for $\textbf{P}^{\big( \text{ } \xi^{\mathrm{S}}\text{ } \big)^{\prime\prime}}_{\textbf{Z} \times [ (y-1)r , (y+2) r]} \bigg[ \text{ }                       \mathrm{Top}_{x,y} \underset{\mathscr{V}_{i^{\prime}}}{\overset{h \geq k-2j}{\longleftrightarrow}    }     \mathrm{Bottom}_{x,y}     \text{ } \bigg]$ is obtained in terms of a crossing probability between $[0,\lfloor \delta^{\prime\prime} r \rfloor] \times \{ 0 \}$ and between $[-m,m] \times \{ n \}$ within $\textbf{Z} \times [0,r]$, instead with the requirement that the height along all faces in the crossing be $\geq k-2j$ in $\mathscr{V}_{i^{\prime}}$.

\bigskip

\noindent On the other hand, the probability of crossings within each $D$ for \textit{good} $(x,y)$, within $\mathcal{I} \big( \text{ } \mathscr{O} \text{ } \big)$, admits the following upper bound,

\begin{align*}
     \textbf{P}^{\xi^{\mathrm{S}}}_{\mathscr{O}}      \bigg[ \text{ }                     \mathcal{R}^{\text{ } \mathrm{Sloped}} \bigg| \text{ }            \mathcal{F}^{\text{ } \mathrm{Sloped}}_{(1-c_0)k+1} \text{ } \cup \big\{ \text{ } \mathcal{X}^{\text{ } \mathrm{Sloped}}    \text{ } \equiv \text{ }   X   \text{ } \big\}    \text{ } \bigg] \text{ }  \tag{\textit{4.5.4}}\\   \leq \text{ } \underset{\textit{good} \text{ } (x,y)}{\prod} \text{ }    \textbf{P}^{\xi^{\mathrm{S}}}_{D}      \bigg[ \text{ }     \mathrm{Top}_{x,y} \underset{\mathscr{V}_{i^{\prime}}}{\overset{h \geq k}{\longleftrightarrow}    }     \mathrm{Bottom}_{x,y}            \text{ }      \bigg| \text{ }            \mathcal{F}^{\text{ } \mathrm{Sloped}}_{(1-c_0)k+1} \text{ } \cup \big\{ \text{ } \mathcal{X}^{\text{ } \mathrm{Sloped}}    \text{ } \equiv \text{ }   X   \text{ } \big\}    \text{ } \bigg]  \text{ } \tag{\textit{4.5.5}} \\ 
     \text{ } \overset{(\textit{4.5.1})}{\leq} \text{ }   \underset{\textit{good} \text{ } (x,y)}{\prod} \text{ }   \textbf{P}^{\big( \text{ } \xi^{\mathrm{S}}\text{ } \big)^{\prime}}_{D} \bigg[ \text{ }                     \mathrm{Top}_{x,y} \underset{\mathscr{V}_{i^{\prime}}}{\overset{h \geq k-j}{\longleftrightarrow}    }     \mathrm{Bottom}_{x,y}     \text{ } \bigg| \text{ }          \mathcal{F}^{\text{ } \mathrm{Sloped}}_{(1-c_0)k+1} \text{ } \cup \big\{ \text{ } \mathcal{X}^{\text{ } \mathrm{Sloped}}    \text{ } \equiv \text{ }   X   \text{ } \big\}        \text{ } \bigg]  \tag{\textit{4.5.6}} \end{align*}

     \begin{align*}   \text{ }       \overset{(\textit{4.5.2})}{\leq} \text{ }      \underset{\textit{good} \text{ } (x,y)}{\prod} \text{ }         \textbf{P}^{\big( \text{ } \xi^{\mathrm{S}}\text{ } \big)^{\prime\prime}}_{[-m,m] \times [ (y-1)r , (y+2) r]} \bigg[ \text{ }      \mathrm{Top}_{x,y} \underset{\mathscr{V}_{i^{\prime}}}{\overset{h \geq k-2j}{\longleftrightarrow}    }     \mathrm{Bottom}_{x,y}     \text{ } \bigg| \text{ }          \mathcal{F}^{\text{ } \mathrm{Sloped}}_{(1-c_0)k+1} \text{ } \cup \big\{ \text{ } \mathcal{X}^{\text{ } \mathrm{Sloped}}    \text{ } \equiv \text{ }   X   \text{ } \big\}                \text{ } \bigg]   \tag{\textit{4.5.7}} \\ 
     \text{ } \overset{m \longrightarrow + \infty}{\leq} \text{ }       \underset{\textit{good} \text{ } (x,y)}{\prod} \text{ }    \textbf{P}^{\big( \text{ } \xi^{\mathrm{S}}\text{ } \big)^{\prime\prime}}_{\textbf{Z} \times [ (y-1)r , (y+2) r]} \bigg[ \text{ }      \mathrm{Top}_{x,y} \underset{\mathscr{V}_{i^{\prime}}}{\overset{h \geq k-2j}{\longleftrightarrow}    }     \mathrm{Bottom}_{x,y}     \text{ } \bigg| \text{ }          \mathcal{F}^{\text{ } \mathrm{Sloped}}_{(1-c_0)k+1} \text{ } \cup \big\{ \text{ } \mathcal{X}^{\text{ } \mathrm{Sloped}}    \text{ } \equiv \text{ }   X   \text{ } \big\}                \text{ } \bigg]   \tag{\textit{4.5.8}}     \\ \leq \text{ }        \underset{\textit{good} \text{ } (x,y)}{\prod} \text{ }    \textbf{P}^{\big( \text{ } \xi^{\mathrm{S}}\text{ } \big)^{\prime\prime}}_{\textbf{Z} \times [ (y-1)r , (y+2) r]} \bigg[ \text{ }      \mathrm{Top}_{x,y} \underset{\mathscr{V}_{i^{\prime}}}{\overset{h \geq k-2j}{\longleftrightarrow}    }     \mathrm{Bottom}_{x,y}     \text{ } \bigg]  \tag{\textit{4.5.9}} \\ \text{ }     \equiv \bigg[ \text{ }       \textbf{P}^{\big( \text{ } \xi^{\mathrm{S}}\text{ } \big)^{\prime\prime}}_{\textbf{Z} \times [ (y-1)r , (y+2) r]} \bigg[ \text{ }      \mathrm{Top}_{x,y} \underset{\mathscr{V}_{i^{\prime}}}{\overset{h \geq k-2j}{\longleftrightarrow}    }     \mathrm{Bottom}_{x,y}     \text{ } \bigg]       \text{ } \bigg]^{\big| \text{ } \textit{good}\text{ } (x,y) \text{ } \big|}   \tag{\textit{4.5.10}} \\ \text{ }           \leq              \big( \text{ }    c_{\mathcal{R},\mathcal{F}} \text{ } \big)^{\big| \text{ } \textit{good}\text{ } (x,y) \text{ } \big|}       \text{ } \tag{\textit{4.5.11}}    \\   \text{ } \leq \text{ }    C^{\big| \text{ } \textit{good}\text{ } (x,y) \text{ } \big|}  \equiv \mathcal{C}_{\mathcal{R}, \mathcal{F}}   \text{ } \text{ , }  \tag{\textit{4.5.12}} 
\end{align*}

\noindent where, from the sequence of inequalities above, beginning from the probability of a \textit{sloped ridge event} occurring in $(\textit{4.5.4})$, in $(\textit{4.5.5})$ the independence of the probability measures $\textbf{P}^{\xi^{\mathrm{S}}}_D \big[ \text{ } \cdot \text{ } \big]$ is employed so that a factorization over $\textit{good} \text{ } (x,y)$ can be obtained, instead supported with the measure $\textbf{P}^{\xi^{\mathrm{S}}}_D\big[ \text{ } \cdot \text{ } \big]$, instead of with $\textbf{P}^{\xi^{\mathrm{S}}}_{\mathscr{O}} \big[ \text{ } \cdot \text{ } \big]$, in $(\textit{4.5.6})$ the conditional pushforward of $\bigg\{       \mathrm{Top}_{x,y} \underset{\mathscr{V}_{i^{\prime}}}{\overset{h \geq k}{\longleftrightarrow}    }     \mathrm{Bottom}_{x,y}        \bigg\} $ is upper bounded, from $(\mathrm{CBC})$ in $(\textit{4.5.2})$, with the conditional event,

\begin{align*}
        \bigg\{  \mathrm{Top}_{x,y} \underset{\mathscr{V}_{i^{\prime}}}{\overset{h \geq k-2j}{\longleftrightarrow}    }     \mathrm{Bottom}_{x,y}     \text{ } \big| \text{ }          \mathcal{F}^{\text{ } \mathrm{Sloped}}_{(1-c_0)k+1} \text{ } \cup \big\{ \text{ } \mathcal{X}^{\text{ } \mathrm{Sloped}}    \text{ } \equiv \text{ }   X   \text{ } \big\}  \bigg\}          \text{ } \text{ , } 
\end{align*}

\noindent hence providing a pushforward under $\textbf{P}^{\big( \text{ } \xi^{\mathrm{S}}\text{ } \big)^{\prime\prime}}_{[-m,m] \times [ (y-1)r , (y+2) r]} \big[ \text{ } \cdot \text{ } \big]$, in $(\textit{4.5.7})$ the probability measure from $\textit{good} \text{ } (x,y)$ in each $D$ is upper bounded with the probability measure $\textbf{P}^{\big( \text{ } \xi^{\mathrm{S}} \text{ } \big)^{\prime}}_{[-m,m] \times [ (y-1)r , (y+2) r]} \big[  \cdot  \big]$, under the same boundary conditions $\big( \text{ } \xi^{\mathrm{S}} \text{ } \big)^{\prime}$, in $(\textit{4.5.8})$ the weak infinite volume limit of $P_m \longrightarrow P_{+\infty}$ is taken for $m \longrightarrow +\infty$, in $(\textit{4.5.9})$ the conditional probability from the event,

\begin{align*}
    \bigg\{    \mathrm{Top}_{x,y} \underset{\mathscr{V}_{i^{\prime}}}{\overset{h \geq k-2j}{\longleftrightarrow}    }     \mathrm{Bottom}_{x,y}     \text{ } \big| \text{ }          \mathcal{F}^{\text{ } \mathrm{Sloped}}_{(1-c_0)k+1} \text{ } \cup \big\{ \text{ } \mathcal{X}^{\text{ } \mathrm{Sloped}}    \text{ } \equiv X \text{ } \big\}  \bigg\} \text{ , } 
\end{align*}

\noindent is less likely to occur than the unconditioned event,

\begin{align*}
    \bigg\{    \mathrm{Top}_{x,y} \underset{\mathscr{V}_{i^{\prime}}}{\overset{h \geq k-2j}{\longleftrightarrow}    }     \mathrm{Bottom}_{x,y}      \bigg\}         \text{ } \text{ , } 
\end{align*}

\noindent when the pushforward, taken under $ \textbf{P}^{\big( \text{ } \xi^{\mathrm{S}}\text{ } \big)^{\prime\prime}}_{\textbf{Z} \times [ (y-1)r , (y+2) r]} \big[ \text{ } \cdot \text{ } \big]$, is computed, in $(\textit{4.5.10})$ the same probability appearing in the upper bound provided in $(\textit{4.5.8})$ is raised to the cardinality of \textit{good} $(x,y)$, in $(\textit{4.5.11})$ a lower bound to $(\textit{4.5.10})$ is provided with suitable $c_{\mathcal{R}, \mathcal{F}}$, where the cardinality of \textit{good} points satisfies,

\begin{align*}
 \big|  \textit{good}\text{ } (x,y)  \big|      \text{ }  < \text{ } \frac{mn}{6}   \text{ } \text{ , } 
\end{align*}

\noindent and to conclude, in $(\textit{4.5.12})$ the final upper bound, as stated in \textbf{Lemma} \textit{4.5}, is provided with $\mathcal{C}_{\mathcal{R},\mathcal{F}}$. \boxed{}

\bigskip

\noindent To provide the form of the exponent to which the free energy for the six-vertex model is raised, we also present the following arguments.

\subsection{Incorporating results from \textbf{Lemma} \textit{4.5} and \textbf{Lemma} \textit{4.6}}

\noindent See the following, pertaining to $\textbf{Definition}$ \textit{10} and $Theorem \textit{6V 2}$, respectively introduced at the beginning of \textit{4.3}, and at the beginning of \textit{4.4}.

\bigskip

\noindent \textit{Proof of Theorem $6V$ \textit{2}}. From previous arguments, a lower bound of the following form,

\begin{align*}
  P^{\prime}   \geq       \textbf{P}^{\xi^{\mathrm{Sloped}}}_{\mathscr{O}}      \bigg[                     \mathcal{F}^{\text{ } \mathrm{Sloped}} \big|  \text{ }                  \mathcal{R}^{\text{ } \mathrm{Sloped}}  \text{ } \cup \text{ }    \big\{ \text{ } \mathcal{X}^{\text{ } \mathrm{Sloped}}    \text{ } \equiv \text{ }   X   \text{ } \big\}             \bigg]  \text{ }    \textbf{P}^{\xi^{\mathrm{Sloped}}}_{\mathscr{O}}      \bigg[           \mathcal{R}^{\text{ } \mathrm{Sloped}}        \big| \mathcal{X}^{\text{ } \mathrm{Sloped}}    \text{ } \equiv \text{ }   X       \bigg]        \\ \overset{(\textbf{Lemma} \text{ } \textit{4.5})}{\geq} \text{ }  \mathcal{C}_{\mathcal{F}, \mathcal{R}}     \text{ }  \textbf{P}^{\xi^{\mathrm{Sloped}}}_{\mathscr{O}}      \bigg[           \mathcal{R}^{\text{ } \mathrm{Sloped}}        \big| \text{ } \mathcal{X}^{\text{ } \mathrm{Sloped}}    \text{ } \equiv \text{ }   X                   \bigg] \text{ }   \\ {\geq}     \mathcal{C}_{\mathcal{F}, \mathcal{R}}    \textbf{P}^{\xi^{\mathrm{Sloped}}}_{\mathscr{O}}      \big[      \mathscr{U} \mathscr{V}         \big]             \text{ } \text{ , } \tag{$\mathscr{U}\mathscr{V}$ \textit{dependent lower bound}}
\end{align*}

\noindent where,

\begin{align*}
    P^{\prime} \equiv  \textbf{P}^{\xi^{\mathrm{Sloped}}}_{\mathscr{O}}      \bigg[                   \mathcal{R}^{\text{ } \mathrm{Sloped}} \text{ } \big| \text{ }            \mathcal{F}^{\text{ } \mathrm{Sloped}}_{(1-c_0)k+1} \text{ } \cup \big\{ \text{ } \mathcal{X}^{\text{ } \mathrm{Sloped}}    \text{ } \equiv \text{ }   X   \text{ } \big\}    \bigg] \text{ } \text{ , }  \text{ } 
\end{align*}

\noindent holds because the \textit{sloped ridge event} is a subset of the union of vertical $\mathrm{x}$- crossing events contained within $\mathcal{I} \big( \text{ } \mathscr{O} \text{ } \big)$ induced by each $\gamma_{i^{*}}$, and $\gamma_{L+i^{*}}$. Recall the denomination provided for $\mathscr{U}\mathscr{V}$ in \textit{4.1},

\begin{align*}
      \mathscr{U} \mathscr{V} \text{ }   \equiv \text{ } \underset{j \text{ } :\text{ }  \mathscr{F}_1 \in \mathcal{B} (  \mathscr{O} ) \text{ } ,\text{ }  \mathscr{F}_2 \in \mathcal{T} (  \mathscr{O} )}{\bigcup} \text{ }     \big\{     \mathscr{F}_1 \underset{\mathcal{I}_j ( \mathscr{O} )}{\overset{h \geq ck}{\longleftrightarrow}}      \mathscr{F}_2  \big\}      \text{ } \text{ . } 
\end{align*}

\noindent

\noindent From \textbf{Lemma} \textit{4.6}, the probability $P^{\prime}$ admits the upper bound $ \mathcal{C}_{\mathcal{F},\mathcal{R}}$,

\begin{align*}
 \mathcal{C}_{\mathcal{F},\mathcal{R}}\text{ }  \geq \bigg[      \textbf{P}^{\big( \xi^{\mathrm{S}}\big)^{\prime\prime}}_{\textbf{Z} \times [ (y-1)r , (y+2) r]} \bigg[      \mathrm{Top}_{x,y} \underset{\mathscr{V}_{i^{\prime}}}{\overset{h \geq k-2j}{\longleftrightarrow}    }     \mathrm{Bottom}_{x,y}  \bigg]    \bigg]^{\big|  \textit{good}\text{ } (x,y)  \big|}    \geq  P^{\prime}   \text{ } \text{ , } 
\end{align*}

\noindent which can be applied once more to obtain the relation below,

\begin{align*}
 \textbf{P}^{\xi^{\mathrm{Sloped}}}_{\textbf{Z} \times [ -r + \delta_1 , 2r + \delta_1]}      \bigg[              [ 0 , \lfloor \delta^{\prime\prime} r  \rfloor ] \times \{0\} \underset{\textbf{Z} \times [0,r]}{\overset{h \geq c_0 k  - j}{\longleftrightarrow}}                 [-m,m] \times \{ n \}             \bigg]   \\
 \equiv  \textbf{P}^{\xi^{\mathrm{Sloped}}}_{\textbf{Z} \times [ -r + \delta_1 , 2r + \delta_1]}      \big[            [ 0 , \lfloor \delta^{\prime\prime} r  \rfloor ] \times \{0\} {\overset{h \geq c_0 k  - j}{\longleftrightarrow}}                 [-m,m] \times \{ n \}             \big]  \text{ } \\ \overset{(\mathscr{U}\mathscr{V}\textit{dependent lower bound})\text{ },\text{ }  (\textbf{Lemma} \text{ }   \textit{4.6})}{\geq} \text{ } \big(    \mathcal{C}_{\mathcal{F},\mathcal{R}} \text{ }           \textbf{P}^{\xi^{\mathrm{Sloped}}}_{\mathscr{O}}      \big[     \mathscr{U} \mathscr{V}      \big]               \big)^{\big| \textit{good}\text{ }  (x,y)  \big|^{-1}}  \text{ . } \tag{\textit{4.7.1}}
\end{align*}

\noindent Next, also recall that the probability of the union of $\mathrm{x}$-connectivity events, $\mathscr{U}\mathscr{V}$, consists of vertical $\mathrm{x}$-crossings occurring for every $\gamma_{i^{*}}$ and $\gamma_{L+i^{*}}$. The pushforward of $\mathscr{U} \mathscr{V}$, under the probability measure with sloped boundary conditions supported over $\mathscr{O}$, admits a lower bound dependent upon a perturbation of the sloped free energy function away from $0$, in which, 

\begin{align*}
           \textbf{P}^{\xi^{\mathrm{Sloped}}}_{\mathscr{O}}      \big[    \mathscr{U} \mathscr{V}      \big]         \geq \mathrm{exp} \big[    M   N  \big(         g_c (\beta) - g_c(0)                          \big)   + \mathcal{O}       \big] \text{ } \tag{\textit{4.7.1.1}} \\ \equiv  \mathrm{exp} \big[   M   N    \big(            g_c (\beta) - g_c(0)                            \big)               \big]   \mathrm{exp} \big[      \mathcal{O}               \big] \\ \overset{(\textit{4.7.2.1})}{\geq}       \epsilon   \mathrm{exp} \big[      M   N  \big(           g_c (\beta) - g_c(0)                        \big)         \big]       \text{ } \text{ , }  \text{ } \tag{\textit{4.7.2}}
\end{align*}

\noindent which includes a prefactor that is dependent upon the number of preimages of $\mathcal{T}^{\prime}$given by (\textit{Property Two}) of $\mathcal{T}^{\prime}$ introduced in \textit{4.3}. In (\textit{4.7.2.1}) to obtain the final lower bound, there exists a sufficiently small, strictly positive, $\epsilon$, which satisfies,

\begin{align*}
\mathrm{exp} \big[ o(1) \big] \equiv  \mathrm{exp} \big[ \mathcal{O}  \big] \text{ } \geq \text{ } \epsilon  \big( N , M \big)  \longrightarrow 0 \text{ , } \text{ as } N , M  \longrightarrow + \infty \text{ } \text{ . } 
\end{align*}

\noindent Hence, the directionality of the following inequality,

\begin{align*}
         \big(  \textbf{P}^{\xi^{\mathrm{Sloped}}}_{\mathscr{O}}      \big[      \mathscr{U} \mathscr{V}      \big]\big)^{{\big| \textit{good}\text{ }  (x,y)  \big|}^{-1}}         \geq \bigg\{  \epsilon  \text{ }  \mathrm{exp} \big[  M N \big(          g_c (\beta) - g_c(0)           \big)    \big]  \bigg\}^{{\big| \textit{good}\text{ }  (x,y)  \big|}^{-1}}  \text{ } \text{ , }  \text{ } \tag{\textit{4.7.3}}
\end{align*}

\noindent is preserved after applying a monotonic transformation dependent upon the rescaled quantities for $N$ and $M$, $n$ and $m$, respectively, by raising the inequality in (\textit{4.7.2}) to the reciprocal of the number of \textit{good} points in the finite volume cylinder. From the exponential lower bound dependent upon the sloped free energy function in (\textit{4.7.3}) above,

\begin{align*}
        \epsilon \text{ }      \mathrm{exp} \big[   M N  \big(         g_c (\beta) - g_c(0)         \big)          \big]      \text{ } \text{ , } \tag{\textit{4.7.3.0}}
\end{align*}

\noindent from the finite-volume prefactor in the possible number of preimages under the pullback $\big(\mathcal{T}^{\prime}\big)^{-1}$ given by (\textit{Property Two}) of $\mathcal{T}^{\prime}$ introduced in \textit{4.3}, therefore implying an equivalent formulation of (\textit{4.7.3}) with (\textit{4.7.3.1}),

\begin{align*}
      \big(   \textbf{P}^{\xi^{\mathrm{Sloped}}}_{\mathscr{O}}      \big[     \mathscr{U} \mathscr{V}      \big] \big)^{{\big| \textit{good}\text{ }  (x,y)  \big|}^{-1}}       \geq  \bigg\{  \epsilon  \text{ } \mathrm{exp} \big[     M N       \big(           g_c (\beta) - g_c(0)             \big)      \big] \text{ } \bigg\}^{{\big| \textit{good}\text{ }  (x,y)  \big|}^{-1}}   \text{ } \tag{\textit{4.7.3.1}}       \text{ } \text{ . } 
\end{align*}

\noindent From a previous result introduced at the beginning of Section \textit{3} for \textit{Russo-Seymour-Welsh in the strip} from \textit{Theorem 6V 1},

\begin{align*}
   \big(   c^{-1}   \textbf{P}^{\xi^{\mathrm{Sloped}}}_{\textbf{Z} \times [0,n^{\prime} N]} \big[      \mathcal{O}_{h\geq ck} ( 6n , n_R n)            \big]  \big)^{\big(  C^{\prime}  \big)^{-1}}  \geq           \textbf{P}^{\xi^{\mathrm{Sloped}}}_{\textbf{Z} \times [-n^{\prime} , 2 n^{\prime} ] }   \bigg[      [ 0 , \lfloor \delta^{\prime} n \rfloor ]           \times \{ 0 \}       \underset{ \textbf{Z} \times [0,n^{\prime} N]}{\overset{h\geq ck}{\longleftrightarrow}}          \textbf{Z}         \times \{ n \}   \bigg]                \text{ } \text{ , } \tag{\textit{4.7.4}}
\end{align*}

\noindent the lower bound that is dependent upon the probability of connectivity between $[ 0 , \lfloor \delta^{\prime} n \rfloor ]           \times \{ 0 \}$, and $\textbf{Z}         \times \{ n \}$ holds. From (\textit{4.7.4}), incorporating lower bounds obtained for either the probability of a connectivity event between $1$ and $1$ occurring, or for the probability of occurring, yields the lower bound, from (\textit{4.7.7}), below,

\begin{align*}
    \textbf{P}^{\xi^{\mathrm{Sloped}}}_{\textbf{Z} \times [0,n^{\prime} N]} \big[      \mathcal{O}_{h\geq ck} ( 6n , n_R n)         \big]  \overset{(\textit{4.7.4})}{\geq}      \bigg[ c \text{ } \textbf{P}^{\xi^{\mathrm{Sloped}}}_{\textbf{Z} \times [ -r + \delta_1 , 2r + \delta_1]}      \bigg[              [ 0 , \lfloor \delta^{\prime\prime} r  \rfloor ] \times \{0\} \underset{\textbf{Z} \times [0,r]}{\overset{h \geq c_0 k  - j}{\longleftrightarrow}}                 [-m,m] \times \{ n \}              \bigg] \bigg]^{C^{\prime}} \\  \text{ } \overset{(\textit{4.7.1})}{\geq} \text{ } c^{C^{\prime}} \bigg[   \big(    \mathcal{C}_{\mathcal{F},\mathcal{R}}        \textbf{P}^{\xi^{\mathrm{Sloped}}}_{\mathscr{O}}      \big[     \mathscr{U} \mathscr{V}       \big]               \big)^{{\big| \textit{good}\text{ }  (x,y)  \big|}^{-1}} \text{ } \bigg]^{C^{\prime}}  \\ \overset{(\textit{4.7.1.1})}{\geq}   c^{C^{\prime}} \big( \big(     \mathcal{C}_{\mathcal{F},\mathcal{R}}  \mathrm{exp} \big[  M N  \big(          g_c (\beta) - g_c(0)                       \big)   + \mathcal{O}          \big]  \text{ } \big)^{{\big| \textit{good}\text{ }  (x,y)  \big|}^{-1}}   \text{ }     \big)^{C^{\prime}}    \\ \overset{(\textit{4.7.2})\text{ } , \text{ } (\textit{4.7.3})}{\geq}  c^{C^{\prime}}   \big(     \epsilon   \mathcal{C}_{\mathcal{F},\mathcal{R}}          \mathrm{exp} \big[   MN    \big(            g_c (\beta) - g_c(0)                          \big)                \big]                         \big)^{C^{\prime}{\big| \textit{good}\text{ }  (x,y)  \big|}^{-1}} \\ \equiv     c^{C^{\prime}{\big| \textit{good}\text{ }  (x,y)  \big|}^{-1}}  \epsilon^{C^{\prime}{\big| \textit{good}\text{ }  (x,y)  \big|}^{-1}}         \big(   \mathcal{C}_{\mathcal{F},\mathcal{R}}    \big)^{C^{\prime}{\big| \textit{good}\text{ }  (x,y)  \big|}^{-1}}     \\ \times  \big(     \mathrm{exp} \big[     M   N \big(         g_c (\beta) - g_c(0)                       \big)              \big]              \big)^{C^{\prime}{\big| \textit{good}\text{ }  (x,y)  \big|}^{-1}}                   \text{ } \\
    \overset{(\textit{4.7.5})}{\equiv}  \big(   c^{C^{\prime}} \epsilon^{C^{\prime}} \mathcal{C}_{\mathcal{F}, \mathcal{R}}            \big)^{{\big| \textit{good}\text{ }  (x,y)  \big|}^{-1}} \text{ }     \big( \mathrm{exp} \big[    M   N \big(           g_c (\beta) - g_c(0)                       \big)            \big]              \big)^{C^{\prime}{\big| \textit{good}\text{ }  (x,y)  \big|}^{-1}}                  \text{ }                \\  \overset{(\textit{4.7.6})}{\geq}    \mathcal{C}_{ c ,\epsilon , \mathcal{F} , \mathcal{R} }   \big(       \mathrm{exp} \big[     M   N   \big(            g_c (\beta) - g_c(0)                            \big)               \big]                  \big)^{C^{\prime}{\big| \textit{good}\text{ }  (x,y)  \big|}^{-1}}                           \text{ } \\ \overset{(\textit{4.7.3.1})}{\geq} \text{ }      \mathcal{C}_{ c , \epsilon , \mathcal{F} , \mathcal{R}}   \big(       \mathrm{exp} \big[        M   N \big(      g_c (\beta) - g_c(0)        \big)          \big]                                             \big)^{C^{\prime}{\big| \textit{good}\text{ }  (x,y)  \big|}^{-1}}  \\ \equiv        \mathcal{C}_{c, \epsilon}       \big(    \mathrm{exp} \big[        M   N  \big(        g_c (\beta) - g_c(0)             \big)        \big] \big)^{C^{\prime}{\big| \textit{good}\text{ }  (x,y)  \big|}^{-1}}   \text{ } \text{ , }     \tag{\textit{4.7.7}}   
\end{align*}

\noindent where, in (\textit{4.7.5}),

\begin{align*}
              c^{{C^{\prime}{{\big| \textit{good}\text{ }  (x,y)  \big|}^{-1}}}} \epsilon^{{C^{\prime}{{\big| \textit{good}\text{ }  (x,y)  \big|}^{-1}}}}         \text{ }     \big(    \mathcal{C}_{c , \epsilon , \mathcal{F},\mathcal{R}}   \big)^{{\big| \textit{good}\text{ }  (x,y)  \big|}^{-1}}       \text{ }  \equiv   \big(    c  \epsilon \big)^{C^{\prime}{\big| \textit{good}\text{ }  (x,y)  \big|}^{-1}}           \big(     \mathcal{C}_{c , \epsilon , \mathcal{F},\mathcal{R}}    \big)^{{{\big| \textit{good}\text{ }  (x,y)  \big|}^{-1}}}     \text{ }   \\ \equiv         \big(    \big(  c  \epsilon \big)^{C^{\prime}}   \mathcal{C}_{c,\epsilon,\mathcal{F}, \mathcal{R}}              \big)^{{\big| \textit{good}\text{ }  (x,y)  \big|}^{-1}}  \text{ } \\ \equiv         \big(   c^{C^{\prime}} \text{ } \epsilon^{C^{\prime}}  \mathcal{C}_{\mathcal{F} , \mathcal{R}}             \big)^{{\big| \textit{good}\text{ }  (x,y)  \big|}^{-1}}  \text{ } \text{ , } 
\end{align*}

 \noindent while, in (\textit{4.7.6}), the inequality between constants,

\begin{align*}
        \big(   c^{C^{\prime}}  \epsilon^{C^{\prime}} \mathcal{C}_{\mathcal{F}, \mathcal{R}}       \big)^{{\big| \textit{good}\text{ }  (x,y)  \big|}^{-1}}       \geq  \mathcal{C}_{c, \epsilon , \mathcal{F} , \mathcal{R}}   \equiv \mathcal{C}_{c , \epsilon}     \text{ } \text{ , } 
\end{align*}

\noindent holds for a suitably chosen constant $ \mathcal{C}_{\text{ } c , \mathcal{F} , \mathcal{R}\text{ } }$ which bounds the product of constants, provided above, from below with $\mathcal{C}_{c , \epsilon}$.

\bigskip

\noindent Hence,

\begin{align*}
       \textbf{P}^{\xi^{\mathrm{Sloped}}}_{\textbf{Z} \times [0,n^{\prime} N]} \big[     \mathcal{O}_{h\geq ck} ( 6n , n_R n)     \big] \overset{(\textit{4.7.7})}{\geq}              \mathcal{C}_{c, \epsilon}               \text{ }  \big(   \mathrm{exp} \big[    M N        \big(           g_c (\beta) - g_c(0)           \big)        \big] \big)^{C^{\prime}{\big| \textit{good}\text{ }  (x,y)  \big|}^{-1}}                          \text{ } \text{ , } 
\end{align*}

\noindent from which we conclude the argument. \boxed{}

\bigskip

\noindent With these results, we progress towards concluding RSW arguments for the six-vertex model.

\section{Logarithmic delocalization of the height function under sloped boundary conditions}

\subsection{Lower bound from $v^{\xi^{\mathrm{Sloped}}}_n$}

\noindent In this section, we collect the following results. First, we modify one expectation value that is introduced in {\color{blue}[11]} for logarithmic bound on variance of the height function,

\begin{align*}
  v_n  \equiv \underset{\xi:\partial \Lambda_n \rightarrow \{ -1 , 0 , +1\}}{\mathrm{min}}\text{ }       \textbf{E}^{\xi}_{\Lambda_n} [ h(0)^2 ]   \text{ }   \text{ , } 
\end{align*}

\noindent to instead hold for sloped boundary conditions, with,

\begin{align*}
  v^{\xi^{\mathrm{Sloped}}}_{n} \equiv v^{\xi^{\mathrm{S}}}_{n} \equiv \underset{\xi^{\mathrm{Sloped}} :\text{ }  \partial \Lambda_n \rightarrow \textbf{B}\textbf{C}^{\mathrm{Sloped}}}{\mathrm{min}} \text{ } \textbf{E}^{\xi^{\mathrm{Sloped}}}_{\Lambda_n} [ h(0)^2 ]   \text{ } \text{ , }
\end{align*}

\noindent in which the boundary conditions over which the infimum is taken do not have to be over $\{ -1 , 0 , +1 \}$. With this modification to $v_n$, arguments from logarithmic bounds for the height function under flat boundary conditions can be applied to the height function under sloped boundary conditions, upon taking into account differences in the scales across the annulus

\bigskip

\noindent \textbf{Lemma} \textit{5.0} (\textit{lower bounding $v^{\xi^{\mathrm{S}}}_n$ across longer scales with $v^{\xi^{\mathrm{S}}}_n$ across shorter scales}). Fix $c \in [1,2]$. For $R \geq 1$, 

\begin{align*}
  v^{\xi^{\mathrm{S}}}_{R^{\prime}n} \geq v^{\xi^{\mathrm{S}}}_n + c_n  \text{ } \text{ , } 
\end{align*}

\noindent given a suitably chosen $R^{\prime} < R$, and a strictly positive constant $c_n$.

\bigskip

\noindent \textit{Proof of Lemma 5.0}. Consider the crossing probability for the annulus, from the absolute value of $h$,

\begin{align*}
  \textbf{P}^{\xi^{\mathrm{Sloped}}}_{\Lambda_{R^{\prime}n}} \big[     \mathcal{O}_{|h| \geq k} \big(    n , R^{\prime} n   \big)          \big]           \text{ } \text{ . } 
\end{align*}

\noindent It suffices to show that the quantity above, for the annulus, across the restricted scale from $n$ to $R^{\prime} n$, occurs with positive probability, and more specifically, some probability greater than $\frac{1}{2}$ for $n$ sufficiently large. To this end, fix boundary conditions $\xi^{\mathrm{S}}$. From our choice of $\xi^{\mathrm{S}}$, minimizing the expectation value for the square of the height function at the origin, as provided in the modification to $v_n$ with $v^{\xi^{\mathrm{S}}}_n$. Given sloped boundary conditions $(\xi^{\mathrm{S}})^{\prime}\geq \xi^{\mathrm{Sloped}}$, the aforementioned observations amount to the following lower bound,

\begin{align*}
  v^{(\xi^{\mathrm{S}})^{\prime}}_{R^{\prime} n} \equiv  \text{ }         \textbf{E}^{(\xi^{\mathrm{S}})^{\prime}}_{\Lambda_{R^{\prime}n}} [ h(0)^2 ]   \text{ } \geq \text{ }  \textbf{E}^{(\xi^{\mathrm{S}})^{\prime}}_{\Lambda_{R^{\prime}n}}  \text{ }    [ ( h(0) +  k_0 )^2 ]                    \text{ }  - \text{ } \sqrt{k_0}    \text{ }  \\   \equiv \text{ }     \textbf{E}^{(\xi^{\mathrm{S})^{\prime}+k_0}}_{\Lambda_{R^{\prime}n}}  \text{ }    [ ( h(0) )^2 ]                    \text{ }  - \text{ } \sqrt{k_0}            \\     \overset{(\mathrm{CBC})}{\geq} \textbf{E}^{\xi^{\mathrm{Sloped}}}_{\Lambda_{R^{\prime}n}} [ h(0)^2 ] - \sqrt{k_0}         \tag{\textit{5.0.1}} \text{ , } 
\end{align*}

\noindent where $k_0$ is a strictly positive parameter from a set $I_{k_0}$, which is the collection of all $k_0$ satisfying,

\begin{align*}
  I_{k_0} \equiv I \big( k_0 \big) \equiv     \big\{     k_0 \in \textbf{N} :    (\xi^{\mathrm{S}})^{\prime} - \xi^{\mathrm{Sloped}}  >  k_0    \big\} \text{ } \text{ . } 
\end{align*}

\noindent Furthermore, in the arguments above, the sloped boundary expectation of the square of the height function at the origin can be decomposed as,

\begin{align*}
    \textbf{E}^{\xi^{\mathrm{S}}}_{\Lambda_{R^{\prime}n}} \big[           h(0)^2  \big]  \text{ } \equiv \text{ }   \textbf{E}^{\xi^{\mathrm{S}}}_{\Lambda_{R^{\prime}n}} \big[           h(0)^2 \text{ } \textbf{1}_{\mathcal{O}_{|h| \geq k}(n,R^{\prime}n)}   \big]  +  \textbf{E}^{\xi^{\mathrm{S}}}_{\Lambda_{R^{\prime}n}} \big[        h(0)^2  \text{ }  \text{ } \textbf{1}_{\mathcal{O}_{|h| \geq k}(n,R^{\prime}n)^c}          \big]     \text{ } \text{ . } 
\end{align*}

\noindent In the equality above, the expected value of the square of the height function at the origin can either depend upon $\mathcal{O}_{|h| \geq k}(n,R^{\prime}n)$, or upon the complementary event, $\mathcal{O}_{|h| \geq k}(n,R^{\prime}n)^c$, occurring with positive probability. To lower bound the first expectation dependent upon $\mathcal{O}_{|h| \geq k}(n,R^{\prime}n)$, observe,

\begin{align*}
       \textbf{E}^{\xi^{\mathrm{S}}}_{\Lambda_{R^{\prime}n}} \big[   h(0)^2 \text{ } \textbf{1}_{\mathcal{O}_{|h| \geq k}(n,R^{\prime}n)}\big]  \text{ }     \equiv \text{ } \underset{\mathrm{domains}\text{ }  \mathscr{D}^{\prime}}{\sum}  \text{ }     \textbf{E}^{\xi^{\mathrm{S}}}_{\mathscr{D}^{\prime}} [ h(0)^2 ] \text{ } \textbf{P}^{\xi^{\mathrm{S}-k}}_{\Lambda_{R^{\prime}n}} [ \mathscr{D} \equiv \mathscr{D}^{\prime} ]    \text{ } \text{ , } \tag{\textit{DOM-SUM}}
\end{align*}

\noindent holds in which the sloped expectation can be expressed as a summation over admissible $\mathscr{D}^{\prime}$ from the set of all possible realizations $\mathscr{D}$, similar to arguments previously given for \textbf{Lemma} \textit{4.1}, and for \textbf{Lemma} \textit{4.5}, respectively in \textit{4.3} and in \textit{4.6}. Proceeding,

\begin{align*}
   \underset{\mathrm{domains}\text{ }  \mathscr{D}^{\prime}}{\sum}  \text{ }     \textbf{E}^{\xi^{\mathrm{S}}}_{\mathscr{D}^{\prime}} [ h(0)^2 ] \text{ } \textbf{P}^{\xi^{\mathrm{S}-k}}_{\Lambda_{R^{\prime}n}} [ \mathscr{D} {\equiv} \mathscr{D}^{\prime} ] \text{ } \overset{(\textit{DOM-SUM})}{\equiv} \text{ }  \underset{\mathrm{domains}\text{ }       \mathscr{D}^{\prime}}{\sum}  \text{ }     \textbf{E}^{\xi^{\mathrm{S}}+ k_0}_{\mathscr{D}^{\prime}} [ h(0)^2 ]        \text{ }      \textbf{P}^{\xi^{\mathrm{S}-k}}_{\Lambda_{R^{\prime}n}} [ \mathscr{D} \equiv \mathscr{D}^{\prime} ]    \text{ } \tag{\textit{5.0.2}} \\ {\equiv} \underset{\mathrm{domains}\text{ }       \mathscr{D}^{\prime}}{\sum}  \text{ }     \textbf{E}^{\xi^{\mathrm{S}}}_{\mathscr{D}^{\prime}} [ (h(0)+k_0)^2 ]        \text{ }      \textbf{P}^{\xi^{\mathrm{S}-k}}_{\Lambda_{R^{\prime}n}} [ \mathscr{D} \equiv \mathscr{D}^{\prime} ]   \tag{\textit{5.0.3}}\\ \overset{(\textit{5.0.1})\text{ } ,\text{ }  h \longleftrightarrow -h }{\geq}  \text{ }           \underset{\mathrm{domains}\text{ }       \mathscr{D}^{\prime}}{\sum}               \text{ }   \big( \text{ }        \textbf{E}^{\xi^{\mathrm{Sloped}}}_{\Lambda_{R^{\prime}n}} [     h(0)^2 ]                   \text{ } + \text{ }    \sqrt{k_0}    \text{ } \big)    \text{ }      \textbf{P}^{\xi^{\mathrm{S}-k}}_{\Lambda_{R^{\prime}n}} [ \mathscr{D} \equiv \mathscr{D}^{\prime} ]    \tag{\textit{5.0.4}} \\\ \equiv \text{ }        \underset{\mathrm{domains}\text{ }       \mathscr{D}^{\prime}}{\sum}      \textbf{E}^{\xi^{\mathrm{Sloped}}}_{\Lambda_{R^{\prime}n}} [     h(0)^2 ]   \text{ }                  \textbf{P}^{\xi^{\mathrm{S}-k}}_{\Lambda_{R^{\prime}n}} [ \mathscr{D} \equiv \mathscr{D}^{\prime} ]   + \sqrt{k_0} \text{ }    \underset{\mathrm{domains}\text{ }       \mathscr{D}^{\prime}}{\sum}                                   \textbf{P}^{\xi^{\mathrm{S}-k}}_{\Lambda_{R^{\prime}n}} [ \mathscr{D} \equiv \mathscr{D}^{\prime} ]    \text{ } \text{ . } \tag{\textit{5.0.5}}
\end{align*}

\noindent In the sequence of rearrangements above, in (\textit{5.0.2}) the (\textit{DOM-SUM}) expression, as a decomposition over $\mathscr{D}^{\prime}$ is applied, in (\textit{5.0.3}) a shift by $k_0$ to the square of the height function at the origin is incorporated, with a corresponding shift of $k_0$ units downwards to the boundary conditions $\xi^{\mathrm{S}}$ on the expectation, in (\textit{5.0.4}) the previous bound obtained in (\textit{5.0.1}) is applied, in addition to the symmetry $h \longleftrightarrow -h$, in (\textit{5.0.5}) an equivalent summation, from two terms taken over domains $\mathscr{D}^{\prime}$, is obtained from the lower bound to (\textit{5.0.3}) with (\textit{5.0.4}). 

\bigskip

\noindent We conclude our arguments for the lower bound of the first expectation value with,

\begin{align*}
  (\textit{5.0.5}) \text{ } \overset{(\textit{5.0.4})\text{ }}{\geq} \text{ }       \textbf{E}^{\xi^{\mathrm{Sloped}}}_{\Lambda_{R^{\prime}n}} [     h(0)^2 ]   \text{ }                  \textbf{P}^{\xi^{\mathrm{S}-k}}_{\Lambda_{R^{\prime}n}} [  \mathcal{O}_{|h| \geq k} (   n , R^{\prime} n )           ]   + \sqrt{k_0} \text{ }          \textbf{P}^{\xi^{\mathrm{S}-k}}_{\Lambda_{R^{\prime}n}} [    \mathcal{O}_{|h| \geq k} (   n , R^{\prime} n )        ]   \\ \equiv \text{ }   \big( \text{ }       \textbf{E}^{\xi^{\mathrm{Sloped}}}_{\Lambda_{R^{\prime}n}} [     h(0)^2 ]     +  \sqrt{k_0}  \big) \text{ }      \textbf{P}^{\xi^{\mathrm{S}-k}}_{\Lambda_{R^{\prime}n}} [    \mathcal{O}_{|h| \geq k} (  n , R^{\prime} n  )        ]   \text{ }  \\ \geq \text{ }   \big( \text{ }  v^{\xi^{\mathrm{S}}}_{n} + \sqrt{k_0}     \text{ }  \big) \text{ }    \textbf{P}^{\xi^{\mathrm{S}-k}}_{\Lambda_{R^{\prime}n}} [    \mathcal{O}_{|h| \geq k} (  n , R^{\prime} n  )        ]     \text{ } \text{ , } \tag{\textit{5.0.6}}
\end{align*}

\noindent concluding the lower bound for the first expectation.

\bigskip

\noindent To bound the second expectation dependent upon $\textbf{1}_{\mathcal{O}_{|h| \geq k}(n,R^{\prime}n)^c}$ in the decomposition of $ \textbf{E}^{\xi^{\mathrm{S}}}_{\Lambda_{R^{\prime}n}} \big[ \text{ }          h(0)^2  \text{ } \big]$, observe,

\begin{align*}  
     \textbf{E}^{\xi^{\mathrm{Sloped}}}_{\Lambda_{R^{\prime}n}}   \big[         h(0)^2\text{ }  \big| \text{ }  |h| > \chi^{\mathrm{Sloped}} \text{ over } \Lambda_{R^{\prime}n}^c   \big]      \text{ } \overset{(\mathrm{FKG-|h|})}{\geq}\text{ }     \textbf{E}^{\xi^{\mathrm{Sloped}}}_{\Lambda_{R^{\prime}n}}   \big[         h(0)^2\text{ } \big| \text{ }          |h| > \xi^{\mathrm{Sloped}} \text{ over } \Lambda_{R^{\prime}n}^c       \big] \text{ } \text{ , } \tag{\textit{5.0.7}}
\end{align*}

\noindent in which, for $\chi^{\mathrm{Sloped}} > \xi^{\mathrm{Sloped}}$, conditionally upon the absolute value of the height function on the complement of $\Lambda_{R^{\prime}n}$. Next,

 \begin{figure}
\begin{align*}
\includegraphics[width=0.98\columnwidth]{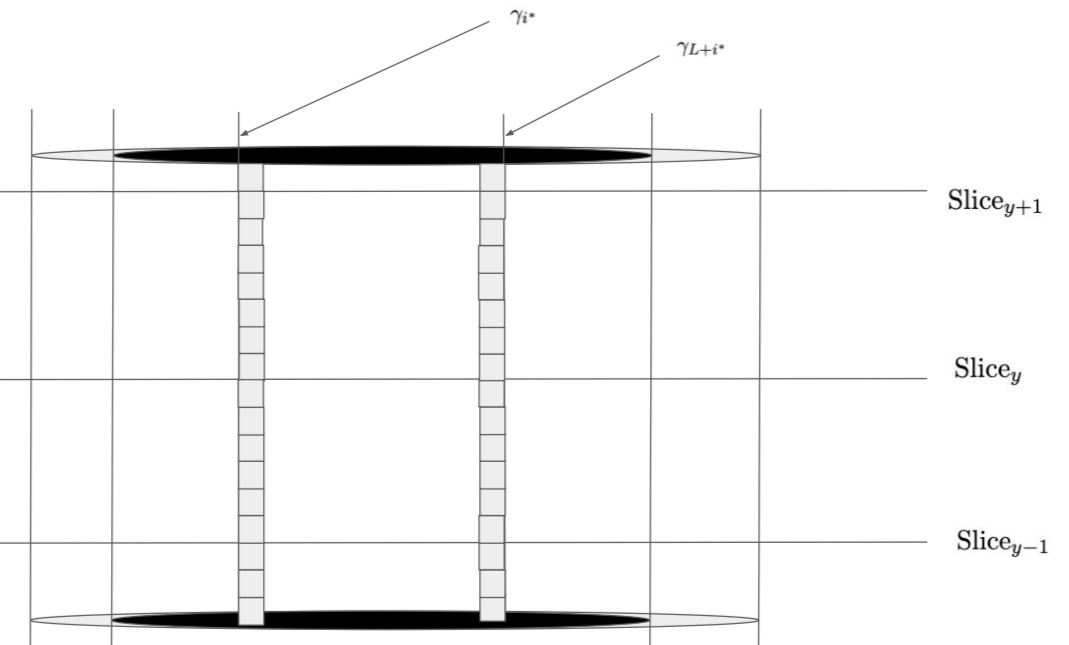}\\
\end{align*}
\caption{\textit{A depiction of incident gray faces to the boundaries induced by the existence of paths $\gamma_i$ and $\gamma_{L+i^{*}}$ in the cylinder considered in the previous section}.}
\end{figure}

\begin{align*}
    (\textit{5.0.7}) \geq \text{ } \underset{\xi^{\mathrm{Sloped}} : \partial \Lambda_{R^{\prime}n} \rightarrow \textbf{B}\textbf{C}^{\mathrm{Sloped}}}{\mathrm{min}}        \textbf{E}^{\xi^{\mathrm{Sloped}}}_{\Lambda_{R^{\prime}n}}   \big[         h(0)^2\text{ } \big|          |h| > \xi^{\mathrm{Sloped}} \text{ over } \Lambda_{R^{\prime}n}^c      \text{ } \big] \text{ } \geq \text{ }   \underset{\xi^{\mathrm{Sloped}} : \partial \Lambda_{R^{\prime}n} \rightarrow \textbf{B}\textbf{C}^{\mathrm{Sloped}}}{\mathrm{min}}  \textbf{E}^{\xi^{\mathrm{Sloped}}}_{\Lambda_{R^{\prime}n}}   [ h(0)^2 ] \\ \equiv   v^{\xi^{\mathrm{S}}}_n   \text{ } \text{ , } \tag{\textit{5.0.8}} 
\end{align*}

\noindent concluding the argument for lower bounding the second expectation, as,

\begin{align*}
    \textbf{E}^{\xi^{\mathrm{S}}}_{\Lambda_{R^{\prime}n}} \big[ \text{ }          h(0)^2  \text{ }  \text{ } \textbf{1}_{\mathcal{O}_{|h| \geq k}(n,R^{\prime}n)^c}            \text{ } \big]          \text{ }  \overset{(\textit{5.0.8})}{\geq} \text{ }  v^{\xi^{\mathrm{S}}}_n \text{ } \textbf{P}^{\xi^{\mathrm{Sloped}}}_{\Lambda_{R^{\prime}n}} [ \text{ } \mathcal{O}_{|h| \geq k} ( n \text{ } , \text{ } R^{\prime} n )^c  \text{ } ] \text{ } \text{ . } \tag{\textit{5.0.9}}
\end{align*}

\noindent Altogether,

\begin{align*}
 v^{(\xi^{\mathrm{S}})^{\prime}}_{R^{\prime} n} \equiv    \textbf{E}^{\xi^{\mathrm{S}}}_{\Lambda_{R^{\prime}n}} \big[ \text{ }          h(0)^2 \text{ } \textbf{1}_{\mathcal{O}_{|h| \geq k}(n,R^{\prime}n)}  \text{ } \big]  \text{ } + \text{ }  \textbf{E}^{\xi^{\mathrm{S}}}_{\Lambda_{R^{\prime}n}} \big[ \text{ }          h(0)^2  \text{ }  \text{ } \textbf{1}_{\mathcal{O}_{|h| \geq k}(n,R^{\prime}n)^c}            \text{ } \big]  \overset{(\textit{5.0.1}
  )\text{ } , \text{ } (\textit{5.0.6}
  )\text{ } ,\text{ } (\textit{5.0.9})}{\geq}  \text{ }    \big( \text{ }  v^{\xi^{\mathrm{S}}}_{n} + \sqrt{k_0}     \text{ }  \big) \text{ } \\ \times    \textbf{P}^{\xi^{\mathrm{S}-k}}_{\Lambda_{R^{\prime}n}} [    \mathcal{O}_{|h| \geq k} (  n , R^{\prime} n  )        ]  +    v^{\xi^{\mathrm{S}}}_n \text{ } \textbf{P}^{\xi^{\mathrm{Sloped}}}_{\Lambda_{R^{\prime}n}} [ \text{ } \mathcal{O}_{|h| \geq k} ( n ,  R^{\prime} n )^c  \text{ } ]  \text{ } - \text{ } \sqrt{k_0} \tag{\textit{5.10}} \\ \text{ }  \equiv \bigg[     \textbf{P}^{\xi^{\mathrm{S}-k}}_{\Lambda_{R^{\prime}n}} [    \mathcal{O}_{|h| \geq k} (   n , R^{\prime} n  )        ]  + \textbf{P}^{\xi^{\mathrm{Sloped}}}_{\Lambda_{R^{\prime}n}} [  \mathcal{O}_{|h| \geq k} ( n  ,  R^{\prime} n )^c   ]       \bigg] v^{\xi^{\mathrm{S}}}_n \text{ } \\ +  \text{ }  \sqrt{k_0} \text{ }         \textbf{P}^{\xi^{\mathrm{S}-k}}_{\Lambda_{R^{\prime}n}} [    \mathcal{O}_{|h| \geq k} (   n , R^{\prime} n )        ] \text{ }   - \sqrt{k_0} \text{ } \tag{\textit{5.11}} \\ \geq                 v^{\xi^{\mathrm{S}}}_n \text{ } + \text{ }       \sqrt{k_0} \text{ }         \textbf{P}^{\xi^{\mathrm{S}-k}}_{\Lambda_{R^{\prime}n}} [    \mathcal{O}_{|h| \geq k} (  n , R^{\prime} n  )    \text{ }     ] \text{ }   - \sqrt{k_0} \text{ }  \tag{\textit{5.12}} \text{ } \text{ , } 
\end{align*}

\noindent where, in (\textit{5.10}) after substituting in the lower bounds for each expectation value from (\textit{5.0.6}) and (\textit{5.0.9}), in (\textit{5.11}) common terms with $v^{\xi^{\mathrm{S}}}_n$ are grouped together, in (\textit{5.12}), the fact that,

\begin{align*}
 \bigg[      \textbf{P}^{\xi^{\mathrm{S}-k}}_{\Lambda_{R^{\prime}n}} [    \mathcal{O}_{|h| \geq k} (   n , R^{\prime} n )        ]  + \textbf{P}^{\xi^{\mathrm{Sloped}}}_{\Lambda_{R^{\prime}n}} [ \mathcal{O}_{|h| \geq k} ( n  ,  R^{\prime} n )^c   ]      \bigg]    v^{\xi^{\mathrm{S}}}_n \geq  v^{\xi^{\mathrm{S}}}_n   \text{ } \text{ , } 
\end{align*}

\noindent holds because,

\begin{align*}
 \textbf{P}^{\xi^{\mathrm{S}-k}}_{\Lambda_{R^{\prime}n}} [    \mathcal{O}_{|h| \geq k} (  n , R^{\prime} n  )        ]  + \textbf{P}^{\xi^{\mathrm{Sloped}}}_{\Lambda_{R^{\prime}n}} [  \mathcal{O}_{|h| \geq k} ( n  , R^{\prime} n )^c ]     > 0   \text{ } \text{ , } 
\end{align*}

\noindent implies that the given lower bound in (\textit{5.12}) holds for (\textit{5.11}). To ensure that the random variable,

\begin{align*}
    \sqrt{k_0}       \textbf{P}^{\xi^{\mathrm{S}-k}}_{\Lambda_{R^{\prime}n}} [    \mathcal{O}_{|h| \geq k} (   n , R^{\prime} n  )      ] -   \sqrt{k_0}       \text{ } \text{ , } 
\end{align*}

\noindent given in (\textit{5.12}) is lower bounded with some $c_n >0$, observe,

\begin{align*}
    \sqrt{k_0}        \textbf{P}^{\xi^{\mathrm{S}-k}}_{\Lambda_{R^{\prime}n}} [    \mathcal{O}_{|h| \geq k} (   n , R^{\prime} n  )      ] -   \sqrt{k_0}   \equiv  \bigg[   \textbf{P}^{\xi^{\mathrm{S}-k}}_{\Lambda_{R^{\prime}n}} [    \mathcal{O}_{|h| \geq k} (   n , R^{\prime} n )       ]  - 1       \bigg]   \sqrt{k_0}     
    \\ \geq   \bigg[    \textbf{P}^{\xi^{\mathrm{S}-k}}_{\Lambda_{R^{\prime}n}} [    \mathcal{O}_{|h| \geq k} (   n , R^{\prime} n )       ]  - 1        \bigg]   c_0  \tag{\textit{5.13}} \\ \geq       c_1 c_0 \geq c_n \tag{\textit{5.14}}              
\end{align*}

\noindent where, in (\textit{5.1.3}), there exists a strictly positive, suitably chosen, $c_0$ so that $\sqrt{k_0} \geq c_0$, while in (\textit{5.14}), there exists a strictly positive, suitably chosen, $c_1$ so that,

\begin{align*}
      \textbf{P}^{\xi^{\mathrm{S}-k}}_{\Lambda_{R^{\prime}n}} [    \mathcal{O}_{|h| \geq k} (   n , R^{\prime} n  )       ] \geq \text{ } c_1        \text{ } \text{ , } 
\end{align*}

\noindent because, as a probability,

\begin{align*}
          1 > \text{ }     \textbf{P}^{\xi^{\mathrm{S}-k}}_{\Lambda_{R^{\prime}n}} [    \mathcal{O}_{|h| \geq k} (  n , R^{\prime} n  )      ]  \text{ } > 0     \text{ } \text{ , } 
\end{align*}

\noindent hence yielding,

\begin{align*}
      v^{\xi^{\mathrm{S}}}_{R^{\prime}n} \geq v^{\xi^{\mathrm{S}}}_n + c_n       \text{ } \text{ , } 
\end{align*}

\noindent from which we conclude the argument. \boxed{}

\bigskip

\noindent In addition to the previous \textbf{Lemma}, also introduce the following result for lower bounding the probability, under sloped boundary conditions instead of flat boundary conditions as provided in \textbf{Lemma 5.2} of {\color{blue}[11]}.

\bigskip

\noindent \textbf{Lemma} \textit{5.1} (\textit{annulus crossing probability across restricted scales from the absolute value of the height function}). For every $k \geq 0$, there exists $c,C,n_0 > 0$ such that for all $N > N^{\prime}$ with $\frac{N^{\prime}}{2} \geq n > n_0$, 

\begin{align*}
 \textbf{P}^{\xi^{\mathrm{S}}}_{\Lambda_{N}} [   \mathcal{O}_{|h| \geq k} ( n , N^{\prime} )  ]   \geq  1 -  C \big(  \frac{n}{N^{\prime}}  \big)^c   \text{ } \text{ , } 
\end{align*}

\noindent for $\xi^{\mathrm{S}} \sim \textbf{B}\textbf{C}^{\mathrm{Sloped}}$.

\bigskip

\noindent \textit{Proof of Lemma 5.1}. We directly apply the inductive argument provided in {\color{blue}[11]}, introducing modifications to the inductive step with the sloped boundary conditions on the six-vertex probability measure, in addition to introducing other modifications throughout the inductive argument through different conditioning on the absolute value of the height function. First, denote the annulus with inner radius $r_1>0$, and outer radius $r_2>0$ with $\mathcal{A} \big( \text{ } r_1 , r_2 \text{ } \big)$. In particular, set $\Lambda_N \equiv N$, fix $i \ge 1$, and for $i^{\prime} < i$ so that $2^i < 2^{i^{\prime}} \leq N^{\prime} < N <  2^{i^{\prime}+1} < 2^{i+1}$,

\begin{align*}
   \textbf{P}^{\xi^{\mathrm{Sloped}}}_N [    \mathcal{O}_{|h| \geq k} ( n , N^{\prime} )    ]   \geq    \textbf{P}^{\xi^{\mathrm{Sloped}}}_N [         \mathcal{O}_{|h| \geq k} (      n     ,    2^{i^{\prime}} n         )    ]  \\ \overset{(\mathrm{FKG-|h|})}{\geq}     \textbf{P}^{\xi^{\mathrm{Sloped}}}_N [         \mathcal{O}_{|h| \geq k} (      n     ,    2^{i^{\prime}} n         )     \text{ }        |    | h | \leq \xi^{\mathrm{Sloped}}                     \text{ over } \partial \Lambda_{2^{i^{\prime}}n}         ]  \tag{\textit{5.1.1}} \\  \underset{(\xi^{\mathrm{Sloped}})^{\prime} \leq \xi^{\mathrm{Sloped}}}{\overset{(\mathrm{SMP})}{\equiv}}     \textbf{P}^{(\xi^{\mathrm{Sloped}})^{\prime}}_{2^{i^{\prime}}n} [      \mathcal{O}_{|h| \geq k + \xi^{\mathrm{Sloped}}- ( \xi^{\mathrm{Sloped})^{\prime}}} (      n     ,    2^{i^{\prime}} n         )     ]  \tag{\textit{5.1.2}} \\ \geq  \underset{(\xi^{\mathrm{Sloped}})^{\prime}: \partial \Lambda_{2^{i^{\prime}}n} \rightarrow \textbf{B}\textbf{C}^{\mathrm{Sloped}} }{\mathrm{min}}   \text{ } \textbf{P}^{(\xi^{\mathrm{Sloped}})^{\prime}} [     \mathcal{O}_{|h| \geq k + \xi^{\mathrm{Sloped}}- ( \xi^{\mathrm{Sloped})^{\prime}}} (      n     ,    2^{i^{\prime}} n         )   ]   \tag{\textit{5.1.3}}  \\ \underset{\xi^{\mathrm{S}} \leq (\xi^{\mathrm{Sloped}})^{\prime}}{\overset{(\mathrm{CBC-|h|})}{\geq}}   \textbf{P}^{\xi^{\mathrm{S}}} [        \mathcal{O}_{|h| \geq k + \xi^{\mathrm{Sloped}}- ( \xi^{\mathrm{Sloped})^{\prime}}} (      n     ,    2^{i^{\prime}} n         )   ]    \tag{\textit{5.1.4}}                \text{ } \text{ , } 
\end{align*}

\noindent where in the series of rearrangements above, beginning in (\textit{5.1.1}), ($\mathrm{FKG-|h|}$) is applied so that the probability of the annulus crossing, $\mathcal{O}_{|h| \geq k} ( n ,2^{i^{\prime}}n)$, is bound below by a conditionally defined event, under $\mathcal{O}_{|h| \geq k} ( n ,2^{i^{\prime}}n)$, is still pushed forwards under $\textbf{P}^{\xi^{\mathrm{Sloped}}} [ \text{ } \cdot \text{ } ]$, in (\textit{5.1.2}) by (SMP), given boundary conditions $\sim \textbf{B}\textbf{C}^{\mathrm{Sloped}}$ for which $(\xi^{\mathrm{Sloped}})^{\prime} \leq \xi^{\mathrm{Sloped}}$, modifying the height required by the annulus crossing $\mathcal{O}$, with height $k + \xi^{\mathrm{Sloped}} - (\xi^{\mathrm{Sloped}})^{\prime}$ instead of $k$, yields a lower bound for the crossing probability in (\textit{5.1.1}), in (\textit{5.1.3}) a lower bound to (\textit{5.1.2}) is provided through an infimum over admissible boundary conditions $\sim \textbf{B} \textbf{C}^{\mathrm{Sloped}}$, for which the annulus crossing $\mathcal{O}_{|h| \geq k + \xi^{\mathrm{Sloped}}- ( \xi^{\mathrm{Sloped})^{\prime}}} (      n     ,    2^{i^{\prime}} n         )$ occurs with positive probability, and finally, in (\textit{5.1.4}), given admissible boundary conditions $\sim \textbf{B}\textbf{C}^{\mathrm{Sloped}}$, for which $\xi^{\mathrm{S}} \leq (\xi^{\mathrm{Sloped}})^{\prime}$, a lower bound is obtained with ($\mathrm{CBC-|h|}$).

\bigskip

\noindent From the final lower bound obtained in (\textit{5.1.4}), we make use of previous arguments to obtain the lower bound,

\begin{align*}
      \textbf{P}^{\xi^{\mathrm{Sloped}}}_{2^{i^{\prime}n}} [    \mathcal{O}_{|h| \geq k} (           n , 2^{i^{\prime}} n             )^c   ] \leq      \big(    1 - \delta^{\prime}                        \big)^{i^{\prime}}                      \text{ } \text{ , }
\end{align*}

\noindent given $\delta^{\prime} \equiv \delta^{\prime} ( k)$. First, for $i \equiv 1$, by strict positivity of the crossing probability across the annulus, the annulus crossing probability across $2n$ admits a similar lower bound as in the inductive argument for flat boundary conditions,

\begin{align*}
       \textbf{P}^{\xi^{\mathrm{Sloped}}}  _{2^{i^{\prime}n}} [   \mathcal{O}_{|h| \geq k} (           n , 2^{i^{\prime}} n         \text{ }     )^c   ]  \overset{j >0}{\geq}               \textbf{P}^{\xi^{\mathrm{Sloped}}}  _{2^{i^{\prime}n}} [   \mathcal{O}_{h \geq k+j} (           n , 2^{i^{\prime}} n         )^c]    \geq    \delta^{\prime}                                     \text{ } \text{ , } 
\end{align*}

\noindent where,

\begin{align*}
       \delta^{\prime}(k) < \delta(k)      \text{ } \text{ , } 
\end{align*}

\noindent from which, directly applying inclusion of events and conditioning on $|h|$, as provided in {\color{blue}[11]}, yields, 

\begin{align*}
     \textbf{P}^{\xi^{\mathrm{Sloped}}}_{2^{i^{\prime}n}}   [ \mathcal{O}_{|h| \geq k} ( n , 2^{i^{\prime}}n )^c ]  \leq       \textbf{P}^{\xi^{\mathrm{Sloped}}}_{2^{i^{\prime}n}}   [       \mathcal{O}_{|h| \geq k}    ( 2n , 2^{i^{\prime}} n)^c                 \cap  \mathcal{O}_{|h| \geq k}    ( n , 2n )^c     ]              \text{ }  \\ \equiv      \textbf{P}^{\xi^{\mathrm{Sloped}}}_{2^{i^{\prime}n}}   [     \mathcal{O}_{|h|\geq k} ( 2n , 2^{i^{\prime}} n )^c         ]         \\ \times  \textbf{P}^{\xi^{\mathrm{Sloped}}}_{2^{i^{\prime}}n} [     \mathcal{O}_{|h| \geq k} ( n , 2n)^c  |         \mathcal{O}_{|h| \geq k} (2n ,    2^{i^{\prime}} n    )^c  ]         \text{ }  \text{ , } 
\end{align*}

\noindent Denoting,

\begin{align*}
   \mathcal{P} \equiv \text{ }     \textbf{P}^{\xi^{\mathrm{Sloped}}}_{2^{i^{\prime}}n} [    \mathcal{O}_{|h| \geq k} ( n , 2n)^c  |         \mathcal{O}_{|h| \geq k} (2n ,    2^{i^{\prime}} n    )^c ]     \text{ } \text{ , } 
\end{align*}

\noindent recall that the conditionally defined measure on the absolute value of the height function,

\begin{align*}
  \mathcal{P} \big[  \cdot  \big] \equiv  \text{ }     \textbf{P}^{\xi^{\mathrm{Sloped}}}_{2^{i^{\prime}}n} [   \mathcal{O}_{|h| \geq k} ( n , 2n)^c \text{ } |        \cdot ]   \text{ } \text{ , } 
\end{align*}

\noindent satisfies ($\mathrm{FKG-|h|}$), in which for two increasing events $\mathcal{E}_1$ and $\mathcal{E}_2$,

\begin{align*}
         \mathcal{P} \big[  \mathcal{E}_1 \cap \mathcal{E}_2  \big]  \geq  \mathcal{P} \big[  \mathcal{E}_1  \big] \text{ }     \mathcal{P} \big[  \mathcal{E}_2  \big]      \text{ } \text{ . } 
\end{align*}

\noindent With further rearrangement from previous steps above imply, again in a similar vein to arguments provided in {\color{blue}[11]}, in which, upon modifying the conditioning upon $|h|$,

\begin{align*}
    \mathcal{P} \text{ }        \overset{(\mathrm{CBC-|h|})}{\geq} \text{ }                        \textbf{P}^{\xi^{\mathrm{Sloped}}}_{2^{i^{\prime}}n} [    \mathcal{O}_{|h| \geq k} ( n , 2n)^c |        \mathcal{O}_{|h| \geq k} (2n ,    2^{i^{\prime}} n^{\prime\prime\prime}    )^c ]       \text{ } \text{ , } 
\end{align*}

\noindent for suitable $n^{\prime\prime\prime} > n$, for which,

 \begin{figure}
\begin{align*}
\includegraphics[width=0.87\columnwidth]{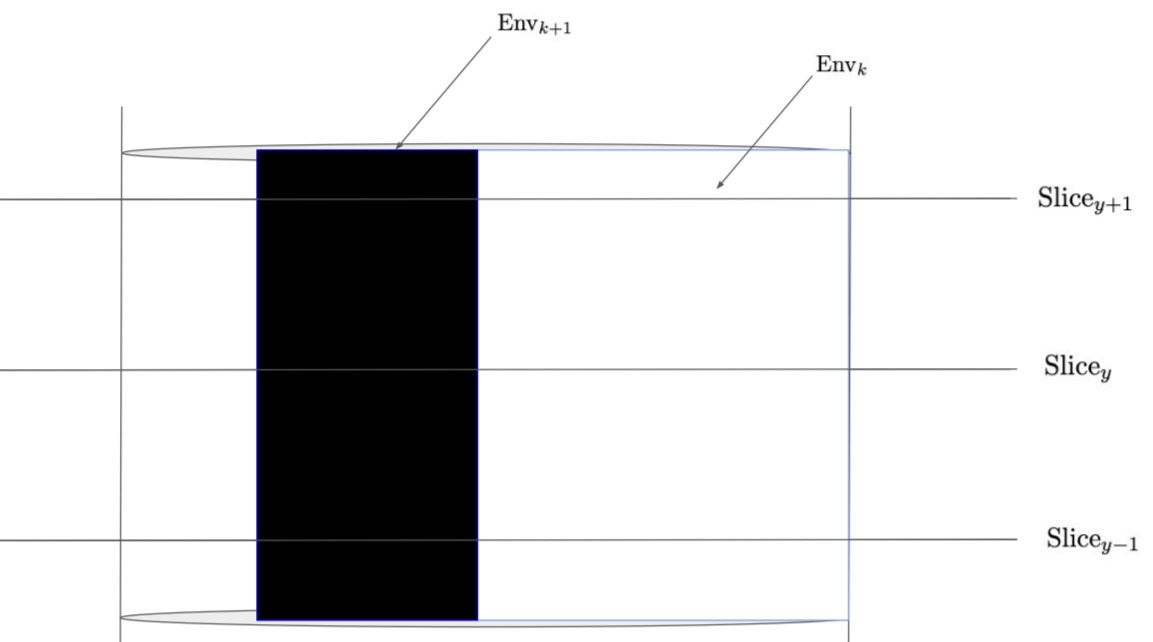}\\
\end{align*}
\caption{\textit{Two exploration sets of the environment of random faces within the interior of the cylinder from arguments used in the previous Lemma}. Depicted above are two exploration sets, $\mathrm{Env}_k$, and $\mathrm{Env}_{k+1}$, each of which contain a nonempty subset of faces in the cylinder. $\mathrm{Env}_k$ is positioned insofar as to intersect with the right cylindrical boundary, in addition to the top and bottom portions of the cylinder, while $\mathrm{Env}_{k+1}$ is positioned insofar as to intersect more faces in the bulk of the interior of the cylinder, while sharing all of the faces that were contained within $\mathrm{Env}_k$.}
\end{figure}

\begin{align*}
   \textbf{P}^{\xi^{\mathrm{Sloped}}}_{2^{i^{\prime}}n} [   \mathcal{O}_{|h| \geq k} (2n ,    2^{i^{\prime}} n^{\prime\prime\prime}    )^c         ]      \leq    \textbf{P}^{\xi^{\mathrm{Sloped}}}_{2^{i^{\prime}}n} [     \mathcal{O}_{|h| \geq k} (2n ,    2^{i^{\prime}} n    )^c   ]    \text{ } \text{ , } 
\end{align*}

\noindent because,

\begin{align*}
   \mathcal{O}_{|h| \geq k} (2n ,    2^{i^{\prime}} n^{\prime\prime\prime}    )^c   \subsetneq     \mathcal{O}_{|h| \geq k} (2n ,    2^{i^{\prime}} n    )^c            \text{ } \text{ , } 
\end{align*}

\noindent which in turn gives, beginning with an application of $(\mathrm{FKG-|h|})$,

\begin{align*}
       \mathcal{P} \text{ } \overset{(\mathrm{FKG-|h|})}{\leq}  \textbf{P}^{\xi^{\mathrm{Sloped}}}_{2^{i^{\prime}}n} \big[              \mathcal{O}_{|h| \geq k} ( n , 2n)^c \big|  \big| h \big| \leq \xi^{\mathrm{Sloped}} , \forall x \in \mathcal{A} \big( \text{ } 2n -1 , 2^{i^{\prime}} n \text{ } \big)                                           \big] \tag{\textit{5.1.5}} \\ \overset{\mathcal{O}_{|h| \geq k}(2n , 2^{i^{\prime}}n) \subsetneq \mathcal{O}_{h \geq k} (2n , 2^{i^{\prime}} n ) }{\leq}     \textbf{P}^{\xi^{\mathrm{Sloped}}}_{2^{i^{\prime}}n} \big[                 \mathcal{O}_{h \geq k} ( n,2n)^c  \big|  \big| h \big| \leq \xi^{\mathrm{Sloped}} , \forall x \in \mathcal{A} \big(  2n -1 , 2^{i^{\prime}} n \big)                                                   \big] \tag{\textit{5.1.6}} \\ \overset{(\mathrm{SMP})}{\equiv} \text{ } \textbf{P}^{(\xi^{\mathrm{Sloped}})_{\mathcal{A}}}_{2^{i^{\prime}}n} \big[                       \mathcal{O}_{h \geq k}(n,2n )                 \big]  \text{ } \\     \tag{\textit{5.1.7}} \\ \overset{(\mathrm{CBC})}{\leq} \text{ } \textbf{P}^{(\xi^{\mathrm{Sloped}})_{\mathcal{A}^{\prime}}}_{2^{i^{\prime}}n} \big[            \mathcal{O}_{h \geq k}(n,2n )                      \big]    \tag{\textit{5.1.8}}        \\  \overset{ h  \leftrightarrow - h}{\equiv}          \textbf{P}^{-(\xi^{\mathrm{Sloped}})_{\mathcal{A}^{\prime}}}_{2^{i^{\prime}}n} \big[              \mathcal{O}_{h \geq k}(n,2n )                      \big]    \tag{\textit{5.1.9}}        \\ \equiv     1 - \textbf{P}^{(\xi^{\mathrm{Sloped}})_{\mathcal{A}^{\prime}}}_{2^{i^{\prime}}n} \big[            \mathcal{O}_{h \geq k+j}(n,2n )                   \big]                     \tag{\textit{5.1.10}}   \\   \leq 1 -   \delta_F                  \tag{\textit{5.1.11}} \text{ } \text{ , } 
\end{align*}

\noindent where, in the sequence of rearrangements above, in (\textit{5.1.5}) $\mathcal{P}$ is upper bounded with the conditional probability on the absolute value of the height function in the annulus $\mathcal{A} \big( \text{ } 2n-1 , 2^{i^{\prime}} n \text{ } \big)$ with $(\mathrm{FKG-|h|})$, in which the annulus satisfies the containment $2^{i^{\prime}} n \subsetneq \mathcal{A} \big( \text{ } 2n-1 , 2^{i^{\prime}} n \text{ } \big)$, in (\textit{5.1.6}) the annulus crossing probability in (\textit{5.1.5}) is upper bounded with the same conditionally defined annulus crossing event, which is instead dependent upon $h$ rather than $|h|$, in (\textit{5.1.7}), by (SMP), the probability measure supported over $2^{i^{\prime}} n $ is equivalent to the probability measure with the same support, but instead with annulus boundary conditions $(\xi^{\mathrm{Sloped}})_{\mathcal{A}} \geq \xi^{\mathrm{Sloped}}$, where the boundary conditions in (\textit{5.1.6}) are dominated by the annulus boundary conditions, in (\textit{5.1.8}), by (CBC), an upper bound to (\textit{5.1.7}) is obtained with $(\xi^{\mathrm{Sloped}})_{\mathcal{A}^{\prime}} \geq \xi^{\mathrm{Sloped}})_{\mathcal{A}^{\prime}}$, in (\textit{5.1.9}), the symmetry $h \leftrightarrow -h$ is applied, in (\textit{5.1.10}) an equivalent random variable for the probability given in (\textit{5.1.9}) is obtained by subtracting $1$ from the probability, instead under boundary conditions $-(\xi^{\mathrm{Sloped}})_{\mathcal{A}^{\prime}}$, which is used to push forward the annulus crossing event $\mathcal{O}_{h \geq k+j}(2n , 2^{i^{\prime}} n ) $ across higher level lines of the height function than required by the annulus crossing event $\mathcal{O}_{h \geq k}(2n , 2^{i^{\prime}} n )$, in (\textit{5.1.11}) a final upper bound is obtained by making use of the fact that the annulus crossing probability, from previous arguments of \textbf{Theorem} \textit{6V 0}, can be lower bounded with a suitable $\delta_F >0$.

\bigskip

\noindent To conclude, given that the induction hypothesis demonstrates that the inequality holds for $i-1$,  

\begin{align*}
     \textbf{P}^{\xi^{\mathrm{Sloped}}}_{2^{i^{\prime}n}} [    \mathcal{O}_{|h| \geq k} (           n , 2^{i^{\prime}-1} n              )^c  ] \leq       \big(     1 - \delta^{\prime}                     \big)^{i^{\prime}-1}                    \Rightarrow  \textbf{P}^{\xi^{\mathrm{Sloped}}}_{2^{i^{\prime}n}} [   \mathcal{O}_{|h| \geq k} (           n , 2^{i^{\prime}} n            )^c  ] \leq      \big(    1 - \delta^{\prime}                     \big)^{i^{\prime}}                           \text{ } \text{ , }
\end{align*}

\noindent from which we conclude the argument. \boxed{}

\bigskip

\noindent \textit{Proof of Theorem 6V 0}. Begin with the annulus crossing probability,

\begin{align*}
    \textbf{P}^{\xi^{\mathrm{Sloped}}}_D [ \mathcal{O}_{h \geq k}(n, n^{\prime} n)] \overset{(\textit{*})}{\equiv} \text{ } \textbf{P}^{\xi^{\mathrm{Sloped}+i^{\prime}}}_D [  \mathcal{O}_{h \geq k + i^{\prime}}(n, n^{\prime} n)] \text{ } \overset{(\textit{**})}{\geq} \text{ } \textbf{P}^{\xi^{\mathrm{Sloped}}}_D [  \mathcal{O}_{h \geq k + i^{\prime}}(n, n^{\prime} n)]  \text{ , } 
\end{align*}

\noindent where, in the sequence of rearrangements above, in (\textit{*}) we obtain an equivalent expression for the first annulus crossing probability by shifting the boundary conditions, in addition to the threshold height of the crossing in the annulus crossing event $\mathcal{O}$ for some $i^{\prime}>0$, and in (\textit{**}) the fact that $\xi^{\mathrm{Sloped}}+i^{\prime} > \xi^{\mathrm{Sloped}}$ readily implies that the crossing probability obtained in (\textit{*}) dominates the crossing probability obtained in (\textit{**}). 

\bigskip

\noindent Next, from an application of a previous result, 

\begin{align*}  
  \textbf{P}^{\xi^{\mathrm{Sloped}}}_{\Lambda}      [       \mathcal{O}_{h \geq k}(n,n^{\prime}n)      ]    \text{ }  \equiv \text{ }  \textbf{P}^{\xi^{\mathrm{Sloped}}}_D      [       \mathcal{O}_{h \geq k}(n,n^{\prime}n)      ]                   \text{ }  \geq \frac{1}{\mathscr{C}} \text{ }      \textbf{P}^{\xi^{\mathrm{Sloped}}}_{\Lambda_{2n}}   [  \mathcal{O}_{h \geq k}(n,n^{\prime}n)  ]    \text{ }  \equiv \text{ } \textbf{P}^{\xi^{\mathrm{Sloped}}}_{\Lambda^{\prime}}   [  \mathcal{O}_{h \geq k}(n,n^{\prime}n)  ]      \text{ } \text{ , } \tag{\textit{***}}
\end{align*}

\noindent with $\Lambda \equiv D$, and $\Lambda^{\prime} \equiv \Lambda_{2n}$ so that $\Lambda \supset \Lambda^{\prime}$. For $1 \leq c \leq 2$, fixing sufficiently large integers so that the annulus crossing probability above, for $h \geq k$, implies an exponential lower bound as a function of the sloped free energy,

\begin{align*}
       \textbf{P}^{\xi^{\mathrm{Sloped}}}_{\Lambda_{12n}} [   \mathcal{O}(6n,n_R n)           ]       \geq       \mathrm{exp} \big[      C^{\prime}\big( r^{\prime}  \big)^2    \big(            g_c (     \beta                  )    -  g_c ( 0 )   \big)           \big]                \text{ } \text{ , } 
\end{align*}

\noindent where, as shown in previous arguments, the parameter at which the sloped free energy function is evaluated at away from $0$ is,

\begin{align*}
  \beta \equiv \frac{k}{\eta r^{\prime}}  \text{ } \text{ . } 
\end{align*}

\noindent From such an exponential lower bound, further manipulate the probability from (\textit{***}), in which, for suitable, strictly positive, $n^{\prime\prime} \geq n^{\prime} \geq n > 0$, and $j$,

\begin{align*}
    \textbf{P}^{\xi^{\mathrm{Sloped}}}_{\Lambda_{n_R n^{\prime\prime}}}  [  \mathcal{O}_{h \geq k+j}   (n^{\prime\prime},n_R n^{\prime\prime}             )  ]   \overset{(\mathrm{SMP})}{\equiv} \text{ } \textbf{P}^{(\chi^{\mathrm{Sloped}})^{\prime}}_{\Lambda_{n_R n}} \big[   \mathcal{O}_{h \geq k+j} (n^{\prime\prime} , n_R n^{\prime\prime} ) \big|         h \equiv   \xi^{\mathrm{Sloped}}        \text{ over }    \partial \Lambda_{n_R n^{\prime\prime}}   \big]  \\     \overset{(\textit{*})}{\leq}            \textbf{P}^{(\chi^{\mathrm{Sloped}})^{\prime}}_{\Lambda_{n_R n}} \big[      \mathcal{O}_{h \geq k+j} (n^{\prime\prime} , n_R n^{\prime\prime} )  \big]    \\  \overset{(\textit{**})}{\equiv} \text{ }                 \textbf{P}^{(\chi^{\mathrm{Sloped}})^{\prime}-j}_{\Lambda_{n_R n}} \big[   \mathcal{O}_{h \geq k}(n^{\prime\prime} , n_R n^{\prime\prime} )     \big]  \\ \overset{(\textit{***})}{\leq}      \textbf{P}^{(\chi^{\mathrm{Sloped}})^{\prime}-j}_{\Lambda_{n_R n}} \big[    \mathcal{O}_{h \geq k}(n^{\prime} , n_R n^{\prime\prime} )     \big]             \\ \overset{(\textit{****})}{\leq}   \textbf{P}^{(\chi^{\mathrm{Sloped}})^{\prime}-j}_{\Lambda_{n_R n}} \big[   \mathcal{O}_{h \geq k}(n , n_R n^{\prime\prime} )        \big]  \\ \overset{(\textit{*****})}{\leq}      \textbf{P}^{(\chi^{\mathrm{Sloped}})^{\prime}-j}_{\Lambda_{n_R n}} \big[    \mathcal{O}_{h \geq k}(n , n_R n^{\prime} )        \big]    \text{ } \\ \overset{(\mathrm{CBC})}{\leq}            \textbf{P}^{\xi^{\mathrm{Sloped}}}_{\Lambda_{n_R n}} \big[    \mathcal{O}_{h \geq k}(n , n_R n^{\prime} )      \big]                   \text{ } \text{ , } 
\end{align*}

\noindent where, in the sequence of rearrangements above, modifying the support of the probability, that is first taken under boundary conditions $\xi^{\mathrm{Sloped}}$, and then under $(\chi^{\mathrm{Sloped}})^{\prime}$, with $(\chi^{\mathrm{Sloped}})^{\prime} \sim \textbf{B}\textbf{C}^{\mathrm{Sloped}}$ being the boundary conditions over $\partial \Lambda_{n_R n^{\prime\prime}}$, of the annulus crossing $\mathcal{O}_{h \geq k}(n,n^{\prime}n)$ occurring can equivalently be expressed with the conditional crossing event provided in the first line. Next, in (\textit{*}), an upper bound to the probability expressed with ($\mathrm{SMP}$) can be obtained by only pushing forwards the annulus crossing event, instead of the conditional annulus crossing event. In (\textit{**}) the probability from (\textit{*}) can be equivalently expressed by shifting down the boundary conditions $(\chi^{\mathrm{Sloped}})^{\prime}$ by some strictly positive parameter $j$, with appropriate modifications in the height required by the annulus crossing. In (\textit{***}) the probability obtained in (\textit{**}) can be upper bounded because, from our choice of strictly positive parameters $n$, $n^{\prime}$ and $n^{\prime\prime}$,

\begin{align*}
 \mathcal{O}_{h \geq k}(     n^{\prime}        ,  n_R n^{\prime\prime}   )  \subset \mathcal{O}_{h \geq k}(  n^{\prime\prime}       ,       n_R n^{\prime\prime}      )  \text{ } \text{ , } 
\end{align*}

\noindent as the annulus crossing event between $n^{\prime}$ and $n_R n^{\prime\prime}$ is strictly contained with the number of faces required for the annulus crossing event between $n^{\prime\prime}$ and $ n_R n^{\prime\prime}$. Along similar lines, in (\textit{****}), 

\begin{align*}
 \mathcal{O}_{h \geq k}(     n        ,  n_R n^{\prime\prime}   )  \subset \mathcal{O}_{h \geq k}(  n^{\prime}       ,       n_R n^{\prime\prime}      )  \text{ } \text{ , } 
\end{align*}

\noindent implies that an upper bound to the probability in (\textit{***}) can be obtained because the annulus crossing event between $n$ and $n_R n^{\prime\prime}$ is strictly contained within the number of faces required for the annulus crossing between $n^{\prime}$ and $n_R n^{\prime\prime}$. Finally, for (\textit{*****}),

\begin{align*}
    \mathcal{O}_{h \geq k}(     n     ,      n_R n^{\prime}         )  \subset \mathcal{O}_{h \geq k}(   n ,       n_R n^{\prime\prime}           )       \text{ } \text{ , } 
\end{align*}

\noindent implies that the last upper bound in the sequence of inequalities above can be obtained, for the same reason as cited in previous steps, hence exhibiting that the desired lower bound holds, as,

\begin{align*}
  \textbf{P}^{\xi^{\mathrm{Sloped}}}_{\Lambda_{n_R n}} \big[   \mathcal{O}_{h \geq k}(n , n_R n^{\prime} )   \big]   \geq    \mathrm{exp} \big[       C^{\prime}  \big(  r^{\prime}\big)^2   \big(          g_c (     \beta                  )    -  g_c ( 0 )   \big)                  \big]    > 0   \text{ } \text{ , } 
\end{align*}

\noindent from which we conclude the argument, with $D \equiv \Lambda_{n_R n}$. \boxed{}

\bigskip

\noindent With the two previous results, in the result below we obtain an upper bound on the sloped variance of the height function, which, recall, was introduced following \textbf{Corollary} \textit{1} in \textit{1.3}.

\bigskip

\noindent \textit{Proof of lower bound on the sloped variance from Corollary 6V 1}. Under the sloped variance,

\begin{align*}
  \mathrm{Var}^{\xi^{\mathrm{Sloped}}}_D [    h(0)              ] \text{ }    \underset{l^{\prime\prime\prime} > 0}{\overset{\text{height function shift invariance}}{\equiv}}        \text{ }       \mathrm{Var}^{\xi^{\mathrm{Sloped}}+l^{\prime\prime\prime}}_D [    h(0)              ] \text{ }                       \text{ }       \text{ , }  \tag{\textit{SLOPE-VAR}}
\end{align*}

\noindent which can be further rearranged insofar as to obtain the following lower bound,

\begin{align*}
     (\textit{SLOPE-VAR}) \text{ } \overset{(\mathrm{FKG-|h|})}{\geq} \text{ } \textbf{E}^{\xi^{\mathrm{Sloped}}}_D [    h(0)^2                          ]      \text{ } - 4 | \xi^{\mathrm{Sloped}}|^2  \overset{(\mathrm{CBC-|h|})}{\geq} \text{ }             \textbf{E}^{\xi^{\mathrm{Sloped}}}_D [    h(0)^2           \big| \text{ } |h| \leq \xi^{\mathrm{Sloped}}     \text{ over } D \backslash \Lambda_n          ]                      \text{ }   \\ -   4 | \xi^{\mathrm{Sloped}}|^2      \text{ } \text{ , } 
\end{align*}

\noindent from which, by (SMP), the above sloped conditional expectation in the final lower bound is equivalent to,

\begin{align*}
 \textbf{E}^{\xi^{\mathrm{Sloped}}}_D [    h(0)^2           \big| \text{ } |h| \leq \xi^{\mathrm{Sloped}}     \text{ over } D \backslash \Lambda_n          ]    \text{ }  \equiv \text{ } \textbf{E}^{(\xi^{\mathrm{Sloped}})_{D \backslash \Lambda_n}} [ h(0)^2 ] \geq v^{\xi^{\mathrm{S}}}_n          \text{ } \text{ , } 
\end{align*}

\noindent which can be bounded from below with the infimum over sloped boundary conditions for which the sloped expected value of the square of the height function occurs, as provided in the definition of $v^{\xi^{\mathrm{S}}}_n$. Finally, from the choice of parameters used in arguments earlier in the section, $v^{\xi^{\mathrm{S}}}_n$ itself admits the following logarithmic lower bound,

\begin{align*}
      v^{\xi^{\mathrm{S}}}_n \geq \text{ } c_n \mathrm{log} [ n ]      \text{ } \text{ , } 
\end{align*}

\noindent for suitable $c_n$, which, upon adjusting constants if necessary, can be equivalent to a lower bound,

\begin{align*}
C_n \text{ }   \mathrm{log} \big[ \frac{d(x,y)}{2} \big]   \text{ } \text{ , } 
\end{align*}

\noindent for suitable $C_n$, from which we conclude the argument. \boxed{}

\bigskip

\noindent We also make use of previous arguments to establish the lower bound of the delocalized phase, with the following.

\bigskip

\noindent \textit{Proof of the lower bound of the delocalized phase provided in Theorem 1.1}. The argument is nearly identical to that provided in {\color{blue}[11]}. The differences in the argument result from the lower bound to the conditional balanced expectation of $h(y)^2$ with sloped boundary conditions instead of an infimum of a pushforward taken under $\textbf{E}^{\xi^{\mathrm{Sloped}}} [ \cdot ]$ for $|\xi|\leq 1$ in finite volume, which was then shown to admit a product lower bound of some constant with a logarithm in $\frac{d(x,y)}{2}$. \boxed{}

\bigskip

\noindent From the lower bounds, with a similarly defined quantity taken under the sloped expectation, we obtain the remaining upper bounds, first for simply connected domains which is then extended to the torus.

\subsection{Upper bound from $w^{\xi^{\mathrm{Sloped}}}_n$}

\noindent Modify the expectation value for $w_n$, defined in {\color{blue}[11]} as,

\begin{align*}
     w_n \equiv \text{ } \underset{\partial D \cap \Lambda_n \neq \emptyset}{\mathrm{sup}} \text{ } \textbf{E}^{0,1}_{D} [ h(0)^2 ]         \text{ } \text{ , } 
\end{align*}

\noindent with,

\begin{align*}
     w^{\xi^{\mathrm{S}}}_n  \equiv   w_n \equiv  \text{ } \underset{\partial D \cap \Lambda_n \neq \emptyset}{\mathrm{sup}} \text{ } \textbf{E}^{\xi^{\mathrm{Sloped}}}_{D} [ h(0)^2 ]         \text{ } \text{ , } 
\end{align*}

\noindent which we further manipulate in the following result.

\bigskip

\noindent \textbf{Lemma} \textit{5.2} (\textit{lower bounding the sloped expectation value}). Fix $c \in [1,2]$. $\forall n \geq 1$, there exists $C^{\prime} >0$ so that,

\begin{align*}
  w_{2n} \leq w_n + C^{\prime}  \text{ } \text{ , } 
\end{align*}

\noindent in which the lower bound for the sloped expectation value across scale $n$ is of scale $2n$.

\bigskip

\noindent \textit{Proof of Lemma 5.2}. Introduce the following random variable previously analyzed in {\color{blue}[11]},

\begin{align*}
 K \equiv \mathrm{inf} \big\{ k \geq 1 : \partial D \longleftrightarrow \Lambda_n  \big\}      \text{ } \text{ , } 
\end{align*}

\noindent from which the connected components of faces over which $|h| \leq k$ on $\partial D$, is given by,

\begin{align*}
   C_k \equiv  \underset{i>0}{\bigcup}   \big\{    \mathscr{F}_i \in F ( \text{ } D \text{ } )  : \big\{ \mathscr{F}_i \notin C_{k-1}  \big\} \cap \big\{  \mathscr{F}_i \notin C_{k+1}   \big\}      \big\}                              \text{ } \text{ , } 
\end{align*}

\noindent in addition to the exploration process about which each cluster up to $C_k$ can be explored, through the revealment procedure for an arbitrary cluster $C_k$,

\begin{align*}
    \mathrm{Reveal} \big(  C_k   \big)  \equiv   {\underset{k > k^{\prime} > 0}{\bigcup}}             \big\{           \mathscr{F} \subset F ( D  )  : \mathscr{F} \cap C_{k^{\prime}} \neq \emptyset        \big\}            \text{ } \text{ , } 
\end{align*}

\noindent in which there will exist some cluster $C_{k+1}$ for which $C_{k+1} \cap D \equiv \emptyset$, after which no more additional executions of the revealment procedure will be necessary. For the subset of faces bound within $D$ that are complementary to $C_k$, denote,

\begin{align*}
     \Omega \equiv \text{ } \underset{u>0}{\bigcup} \big\{ \forall \text{ } k \in \textbf{Z} \text{ }  ,  \text{ }    \exists  \text{ }      \mathscr{F} \subset F ( \text{ } D \text{ } ) :  \mathscr{F} \cap C_k \equiv \emptyset \big\}       \text{ } \text{ , } 
\end{align*}

\noindent as the collection of all faces in the complementary region in $D$ to each $C_k$, $D \backslash C_k$. From $\Omega$, introduce,

\begin{align*}
\textbf{E}^{\xi^{\mathrm{Sloped}}}_D [ h(0)^2 ]  \equiv      \underset{\text{ boundary conditions } \big( \xi_{\mathscr{D}} \big) }{\underset{\text{ domains } \mathscr{D}^{\prime\prime} }{\sum}}    \textbf{P}^{\xi^{\mathrm{Sloped}}}_D \big[ h(0)^2  \big| \Omega \equiv \mathscr{D}^{\prime\prime} , h \equiv \xi_{\mathscr{D}} \text{ over } C_k  \big] \\  \text{ } \times \textbf{P}^{\xi^{\mathrm{Sloped}}}_D \big[    \Omega \equiv \mathscr{D}^{\prime\prime} , h \equiv \xi_{\mathscr{D}} \text{ over } C_k       \big]         \text{ } \text{ , } 
\end{align*}

\noindent as the decomposition of the sloped expectation of the square of the height function at the origin, which is dependent upon the total number of possible realizations of domains, and corresponding sloped boundary conditions. 

\bigskip

\noindent Introduce the probability for an annulus $\mathrm{x}$-crossing,

\begin{align*}
  \textbf{P}^{\xi^{\mathrm{Sloped}}}_D \big[          \mathcal{O}^{\mathrm{x}}_{|h| \geq k}(n)   \big] \text{ }      \text{ , } \tag{\textit{Annulus-x crossing}}
\end{align*}

\noindent which admits an identical upper bound through the following series of inequalities introduced for flat boundary conditions,

\begin{align*}
  (\textit{Annulus-x-crossing}) \text{ } \overset{\mathcal{O}^{\mathrm{x}}_{|h| \geq k} (n )  \subsetneq ( \mathcal{O}^{\mathrm{x}}_{h \geq k} (n )   \cup   \mathcal{O}^{\mathrm{x}}_{|h| \leq - k} (n ) ) }{\leq} \text{ }  \textbf{P}^{\xi^{\mathrm{Sloped}}}_D \big[         \mathcal{O}^{\mathrm{x}}_{h \geq k}(n)   \big] + \textbf{P}^{\xi^{\mathrm{Sloped}}}_D \big[          \mathcal{O}^{\mathrm{x}}_{h \leq  - k}(n)    \big]                 \text{ }  \\    \text{ }        \overset{\textbf{P}^{\xi^{\mathrm{Sloped}}}_D [ \mathcal{O}^{\mathrm{x}}_{h \leq - k }(n) ] \leq \textbf{P}^{\xi^{\mathrm{Sloped}}}_D [ \mathcal{O}^{\mathrm{x}}_{h \geq  k }(n) ] }{\leq}  \text{ }  \textbf{P}^{\xi^{\mathrm{Sloped}}}_D \big[    \mathcal{O}^{\mathrm{x}}_{h \geq k}(n)    \big] + \textbf{P}^{\xi^{\mathrm{Sloped}}}_D \big[          \mathcal{O}^x_{h \geq   k}(n)  \big]  \\ \equiv \text{ } 2 \text{ } \textbf{P}^{\xi^{\mathrm{Sloped}}}_D [ \mathcal{O}^{\mathrm{x}}_{h \geq k}(n) ]    \text{ } \text{ . } 
\end{align*}

\noindent Under $\xi^{\mathrm{Sloped}} \sim \textbf{B}\textbf{C}^{\mathrm{Sloped}}$, it will be inductively demonstrated that,

\begin{align*}
 \textbf{P}^{\xi^{\mathrm{Sloped}}}_D [ \mathcal{O}^{\mathrm{x}}_{h \geq 2k^{\prime}} (n)   ] \text{ } \leq \mathrm{exp} \big(  - c  k^{\prime}     \big)   \text{ } \text{ , } 
\end{align*}

\noindent for $k > k^{\prime}$, from which it suffices to first demonstrate that the inequality holds for $k^{\prime} \equiv 1$, in which,

\begin{align*}
   \textbf{P}^{\xi^{\mathrm{Sloped}}}_D [ \mathcal{O}^{\mathrm{x}}_{h \geq 2 } (n)   ] \text{ } \leq \mathrm{exp} \big(  - c    \big) \text{ } \text{ , } 
\end{align*}

\noindent as the level lines of the height function need to exceed two in order for the annulus $\mathrm{x}$-crossing event to occur. In line with previous observations and remarks for arguments to obtain the lower bound through previous induction arguments, the arguments will demonstrate that the statement holds for $i+1$, from which it will also hold for $i$. From the $\mathrm{x}$-crossing across the annulus introduced previously, rearranging the conditional probability below yields,

\begin{align*}
  \textbf{P}^{\xi^{\mathrm{Sloped}}}_D \bigg[               \mathcal{O}_{h\geq k+j}^{\mathrm{x}}    ( n )             \big|    \mathcal{O}_{h \geq k - 2j}^{\mathrm{x}} \text{ } \cup \text{ } \big\{       \mathscr{D}^{\prime\prime} \equiv D       \big\} \bigg]  \text{ } \overset{(\mathrm{SMP})}{\equiv} \text{ }    \textbf{P}^{\xi^{\mathrm{Sloped}}}_{\mathscr{D}^{\prime\prime}}  \text{ }   \big[   \mathcal{O}_{h\geq k+j}^{\mathrm{x}}    ( n )   \big]      \text{ } \text{ , } \tag{\textit{j x-annulus crossing}}
\end{align*}

\noindent for $\mathscr{D}^{\prime\prime}$ the $\mathrm{x}$-loop that is the closest to the finite volume boundary. Upon performing a global downwards shift of the height function, is equivalent to the following probability under appropriately modified boundary conditions,

\begin{align*}
(\textit{j x-annulus crossing} ) \text{ }    \equiv      \text{ }  \textbf{P}^{\xi^{\mathrm{Sloped}}- j^{\prime}}_{\mathscr{D}^{\prime\prime}} [   \mathcal{O}^{\mathrm{x}}_{h \geq k + j - j^{\prime}} ( n ) ] \text{ } \equiv \text{ } \textbf{P}^{\xi^{\mathrm{Sloped}}_{k,j^{\prime}}}_{\mathscr{D}^{\prime\prime}} [     \mathcal{O}^{\mathrm{x}}_{h \geq k + j - j^{\prime}} ( n )               ]  \text{ }                 \text{ } \text{ . } 
\end{align*}

\noindent given that the modification to the height function crossing, $k+j-j^{\prime}$ for the boundary conditions $\xi^{\mathrm{Sloped}}_{k,j^{\prime}}$, which are $\xi^{\mathrm{Sloped}}$ shifted downwards by some strictly positive translation $j^{\prime}$, satisfy,

\begin{align*}
      k + j - j^{\prime} > 0    \text{ } \text{ . } 
\end{align*}

\noindent As a result, to demonstrate that the inductive step holds, which hence would imply that the annulus crossing probabilities from the absolute value of the height function exponentially decay, write,

\begin{align*}
    \textbf{P}^{\xi^{\mathrm{Sloped}}_{k,j^{\prime}}}_{\mathscr{D}^{\prime\prime}} [     \mathcal{O}^{\mathrm{x}}_{h \geq k + j - j^{\prime}} ( n )          ] \text{ }    \leq \text{ }    \textbf{P}^{\xi^{\mathrm{Sloped}}_{k,j^{\prime}}}_{\mathscr{D}^{\prime\prime}} [ \partial \Lambda_n   (z) \text{ }  {\overset{h \geq k + j - j^{\prime}}{\longleftrightarrow}}_{\mathrm{x}}   \partial \Lambda_{2n} (z)  ]                       \text{ } \text{ , } 
\end{align*}

\noindent in which, from the upper bound above, the $\mathrm{x}$-annulus crossing probability from the event $\mathcal{O}^{\mathrm{x}}_{h \geq k + j - j^{\prime}} ( n )$ is dominated by a crossing of height $h \geq k + j - j^{\prime}$, where the two boxes $\Lambda_n(z)$, and $\Lambda_{2n} (z)$, denote,

\begin{align*}
     \Lambda_n(z) \equiv \big\{   z \notin D  : \big|  \partial \Lambda_n (z)  \big| \equiv 4n    \big\}           \text{ } \text{ , } 
\end{align*}

\noindent and, similarly, 

\begin{align*}
     \Lambda_{2n}(z) \equiv \big\{     z \notin D          : \big|  \partial \Lambda_{2n} (z) \big| \equiv 8n        \big\}   \text{ } \text{ . } 
\end{align*}

\noindent Further analyzing the connectivity probability between $\partial \Lambda_n$ and $\partial \Lambda_{2n}$ yields an equivalent random variable,

\begin{align*}
     \textbf{P}^{\xi^{\mathrm{Sloped}}_{k,j^{\prime}}}_{\mathscr{D}^{\prime\prime}} [          \partial \Lambda_n( z)             \overset{h \geq k+j-j^{\prime}}{\longleftrightarrow }_{\mathrm{x}}            \partial \Lambda_{2n}(z)    ]     \equiv 1 - \textbf{P}^{\xi^{\mathrm{Sloped}}_{k,j^{\prime}}}_{\mathscr{R}^{\prime\prime}} \bigg[ \text{ }   \mathcal{O}^{\mathrm{x}}_{h \leq k + j - j^{\prime} -1} \big(    A(n,2n) + z   \big)    \big|  \text{ }  \cdots \\  \cdots   h \leq  \xi^{\mathrm{Sloped}}_{k,j^{\prime}} \text{ } \text{ over } \Lambda_{2n}(z)\text{ }  \backslash \text{ } \widetilde{\mathscr{D}^{\prime\prime}}        \bigg]   \text{ } \text{ . }  \tag{$\mathscr{R}^{\prime\prime}$\textit{probability}}
\end{align*}

\noindent where, in the (LHS) above, for the interior of the region $\mathscr{D}^{\prime\prime}$, which we denote with $\widetilde{\mathscr{D}^{\prime\prime}}$, the support on the conditional probability measure is given by,

\begin{align*}
 \mathscr{R}^{\prime} \equiv \text{ }   \mathscr{D}^{\prime\prime}       \cup  \Lambda_{2n}(z) \text{ } \text{ , } 
\end{align*}

\noindent for which,

\begin{align*}
    \big(     \mathscr{D}^{\prime\prime}       \cup  \Lambda_{2n}(z)   \big)      \cap       \mathscr{R}^{\prime}         \neq \emptyset  \text{ } \text{ , } 
\end{align*}

\noindent in addition to the modified annulus $\mathrm{x}$-crossing event, which is explicitly given with,

\begin{align*}
     \mathcal{O}^{\mathrm{x}}_{h \leq k + j - j^{\prime} -1} \big( \text{ }   A(n,2n) + z    \text{ } \big)  \equiv          \bigg\{ \forall z \notin D , \exists \text{ }  \mathscr{F}_1 , \mathscr{F}_2 \in F \big(  \mathscr{R}^{\prime\prime}  \big)   : \textbf{P}^{\xi^{\mathrm{Sloped}}}_{\mathscr{R}^{\prime\prime}} \big[ \mathscr{F}_1   \underset{A(n,2n)}{\overset{h \leq k + j - j^{\prime}-1}{\longleftrightarrow}}_{\mathrm{x}}  \mathscr{F}_2  \big] > 0 \bigg\}      \text{ } \text{ , }
\end{align*}

\noindent where the $\mathrm{x}$-crossing within $A(n,2n)$, 

\begin{align*}
 \big\{ \mathscr{F}_1  \underset{A(n,2n)}{\overset{h \leq k + j - j^{\prime}-1}{\longleftrightarrow}}_{\mathrm{x}}   \mathscr{F}_2  \big\}    \text{ } \text{ , } 
\end{align*}

\noindent occurs with positive probability.

\bigskip

\noindent From this conditionally dependent probability upon the value of the height function attained over the complementary region to $\mathscr{R}^{\prime}$.

\bigskip

\noindent From a previous expression, rearranging yields,

\begin{align*}
  (\mathscr{R}^{\prime\prime}\textit{probability})   \geq  \textbf{P}^{\xi_{2n}}_{\Lambda_{2n}(z)}               \bigg[        \mathcal{O}^{\mathrm{x}}_{h \leq k + j - j^{\prime} -1} \big( \text{ }   A(n,2n) + z    \text{ } \big)   \big|  h \leq  \xi^{\mathrm{Sloped}}_{k,j^{\prime}} \text{ } \text{ over } \Lambda_{2n}(z)\text{ }  \backslash \text{ } \widetilde{\mathscr{D}^{\prime\prime}}   \bigg]    \text{ } \text{ , } \tag{$\mathscr{R}^{\prime\prime}$ \textit{lower bound probability}}
\end{align*}

\noindent where, in the lower bound above, the probability measure supported over the translated box, by $z$, of side length $2n$ is taken under $\xi_{2n} \sim \textbf{B}\textbf{C}^{\mathrm{Sloped}}$. 

\bigskip

\noindent To perform the inductive step, set $\textbf{P}^{\xi^{\mathrm{Sloped}}_{k,j^{\prime}}}_{\mathscr{D}^{\prime\prime}} [ \cdot ]  \equiv \textbf{P}^{k,j^{\prime}}_{\mathscr{D}^{\prime\prime}} [ \cdot ]$ to lighten notation, and next, observe,

\begin{align*}
  \textbf{P}^{k,j^{\prime}}_{\mathscr{D}^{\prime\prime}}        \big[   \partial \Lambda_n(z) \overset{h\geq k + j - j^{\prime}}{\longleftrightarrow}_{\mathrm{x}}   \partial \Lambda_{2n}(z)  \big]   \leq \textbf{P}^{\xi_{2n}}_{\Lambda_{2n}(z)} \big[   \partial \Lambda_n(z) \overset{h\geq k + j - j^{\prime}}{\longleftrightarrow}_{\mathrm{x}} \text{ }   \partial \Lambda_{2n}(z)    \big|    \mathcal{C}^{\prime}   \big]    \tag{$\mathscr{R}^{\prime\prime}$\textit{1}}  \end{align*}

  \noindent which can be further upper bounded with the probability, 
  
  \begin{align*}   \textbf{P}^{\xi_{2n}}_{\Lambda_{2n}(z)} \big[       \partial  \Lambda_n(z) \overset{h\geq k + j - j^{\prime} - n }{\longleftrightarrow}    \partial \Lambda_{2n}(z)   \big|  \mathcal{C}^{\prime}         \big]             \tag{$\mathscr{R}^{\prime\prime}$\textit{2}}     \text{ } \text{ , } 
\end{align*}

\noindent where,

\begin{align*}
  \mathcal{C}^{\prime} \equiv  \text{ }  \big\{                 k+j  - j^{\prime}-1   \leq  h \leq k + j - j^{\prime}+1 \text{ over } \Lambda_{2n}(z)\text{ }  \backslash\text{ }  \widetilde{\mathscr{D}^{\prime\prime}} \big\} \text{ , } 
\end{align*}

\noindent and in the sequence of rearrangements above, in ($\mathscr{R}^{\prime\prime}$\textit{1}), we make use of the fact, that conditionally upon the absolute value of the height function in the region outside of $\Lambda_{2n}(z)$ over which the probability measure in the upper bound to $\textbf{P}^{k,j^{\prime}}_{\mathscr{D}^{\prime\prime}}        \big[ \text{ }   \partial \Lambda_n(z) \overset{h\geq k + j - j^{\prime}}{\longleftrightarrow}_{\mathrm{x}} \text{ }   \partial \Lambda_{2n}(z) \text{ } \big]$, that conditioning upon the value of the height function is conducive for making the connectivity event $\big\{ \text{ }    \partial \Lambda_n(z) \overset{h\geq k + j - j^{\prime}}{\longleftrightarrow}_{\mathrm{x}} \text{ }   \partial \Lambda_{2n}(z)      \text{ } \big\}$ occur, while in ($\mathscr{R}^{\prime\prime}$\textit{2}), we apply additional arguments to conclude that an upper bound to ($\mathscr{R}^{\prime\prime}$\textit{1}) exists, when instead the conditional crossing event pushed forwards under , which, as an annulus event that is not an $\mathrm{x}$-crossing, satisfies

\begin{align*}
 \big\{             \partial \Lambda_n(z) \overset{h\geq k + j - j^{\prime}}{\longleftrightarrow}_{\mathrm{x}}    \partial \Lambda_{2n}(z)    \big\}  \subsetneq   \big\{    \partial  \Lambda_n(z) \overset{h\geq k + j - j^{\prime} - n }{\longleftrightarrow}  \partial \Lambda_{2n}(z)      \big\}  \text{ } \text{ , } 
\end{align*}

\noindent for suitable, strictly positive, $n$, so that, upon the conditionally defined probability measure,

 \begin{figure}
\begin{align*}
\includegraphics[width=0.75\columnwidth]{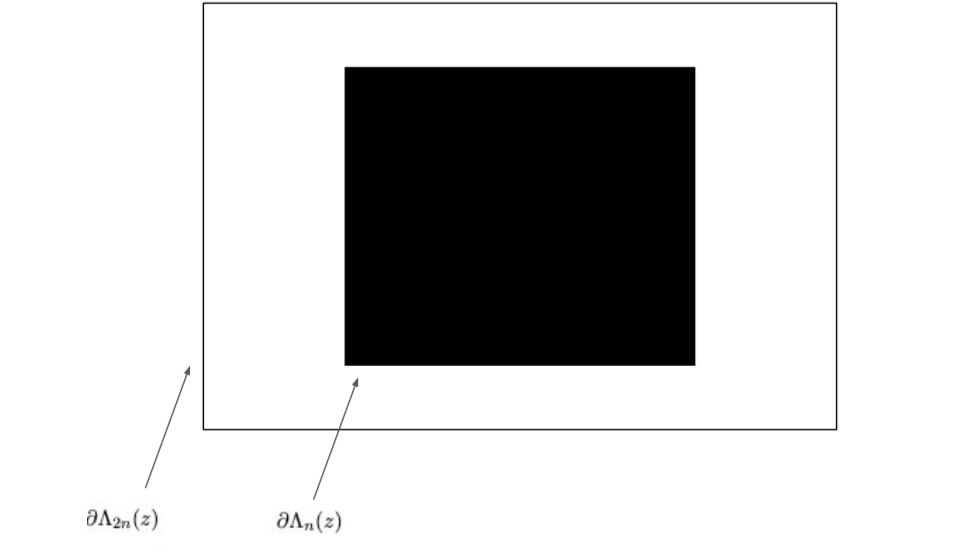}\\
\end{align*}
\caption{\textit{A depiction of the finite volumes}, $\partial \Lambda_n (z)$, and $\partial \Lambda_{2n}(z)$, which is used to form the conditional crossing event used in arguments below.}
\end{figure}

\begin{align*}
  \textbf{P}^{\xi_{2n}}_{\Lambda_{2n}(z)} \big[ \text{ } \cdot \text{ } \big|  \mathcal{C}^{\prime} \big]   \text{ } \text{ , } 
\end{align*}

\noindent $(\mathscr{R}^{\prime\prime}\textit{1}) \leq (\mathscr{R}^{\prime\prime}\textit{2})$. Proceeding further,

\begin{align*}
  (\mathscr{R}^{\prime\prime}\textit{2} ) \equiv  \textbf{P}^{\xi_{2n}}_{\Lambda_{2n}(z)} \bigg[           \partial \Lambda_n(z) \overset{h \geq k+j-j^{\prime}-n}{\longleftrightarrow}     \partial \Lambda_{2n}(z)  \big|  \xi_{2n} + 1 \geq h \geq \xi_{2n}-1 \text{ } \text{ over } \Lambda_{2n}(z)\text{ }  \backslash \text{ }  \widetilde{\mathscr{D}^{\prime\prime}} \bigg]     \text{ }  \tag{\textit{1.1}}   \\ \equiv \textbf{P}^{\xi_{2n}}_{\Lambda_{2n}(z)} \bigg[  \partial \Lambda_n(z) \text{ } {\overset{h \geq k + j - j^{\prime}}{\longleftrightarrow}}  \text{ } \big\{ \text{ }  \partial \Lambda_{2n}(z) \cap \big\{ \text{ } \xi_{2n}   \geq k + j - j^{\prime}  \text{ } \big\}  \text{ } \big\} \text{ }  \big| \text{ }        \xi_{2n} + 1 \geq h \geq \xi_{2n}-1 \text{ over }  \Lambda_{2n}(z)\text{ }  \backslash \text{ }  \widetilde{\mathscr{D}^{\prime\prime}}   \bigg] \tag{\textit{1.2}} \end{align*}

  \begin{align*}
  \overset{(\mathrm{FKG-|h|})}{\leq}      \textbf{P}^{\xi_{2n}}_{\Lambda_{2n}(z)} \bigg[    \partial \Lambda_n(z) \text{ } {\overset{|h| \geq k + j - j^{\prime}}{\longleftrightarrow}}  \text{ } \big\{ \text{ }  \partial \Lambda_{2n}(z) \cap \big\{ \text{ } \xi_{2n}   \geq k + j - j^{\prime}  \text{ } \big\}  \text{ } \big\} \text{ } \cdots \\  \big| \text{ }        \xi_{2n} + 1 \geq |h| \geq \xi_{2n}-1 \text{ over }  \Lambda_{2n}(z)\text{ }  \backslash \text{ }  \widetilde{\mathscr{D}^{\prime\prime}}   \bigg]                 \tag{\textit{1.3}} \\ 
   \leq \text{ } \textbf{P}^{\xi_{2n}}_{\Lambda_{2n}(z)} \bigg[   \partial \Lambda_n(z) {\overset{|h| \geq k + j - j^{\prime}}{\longleftrightarrow}}  \bigg\{  \partial \Lambda_{2n}(z) \cap \big\{ \text{ } \xi_{2n}   \geq k + j - j^{\prime}  \big\}  \bigg\}  \bigg]        \tag{\textit{1.4}} \\ 
     \equiv \text{ }      \textbf{P}^{\xi_{2n}}_{\Lambda_{2n}(z)} \big[   \partial \Lambda_n(z)       {\overset{h \geq k + j - j^{\prime}}{\longleftrightarrow}}  \partial \Lambda_{2n}(z)            \big]         \tag{\textit{1.5}}  \\
   \text{ } \equiv      \textbf{P}^{-\xi_{2n}}_{\Lambda_{2n}(z)}                    \big[    \partial \Lambda_n(z) \text{ }      {\overset{h < k + j - j^{\prime}}{\longleftrightarrow}} \text{ }  \partial \Lambda_{2n}(z)         \big]                  \text{ }      \text{ , }   \text{ } \tag{\textit{1.6}} 
\end{align*}

\noindent where, in the sequence or rearrangements above, beginning in (\textit{1.1}) the probability given in ($\mathscr{R}^{\prime\prime}$\textit{1}) is equivalent to the given conditionally defined event above pushed forwards under $\textbf{P}^{\xi_{2n}}_{\Lambda_{2n}(z)} \big[ \cdot \big]$, in (\textit{1.2}) the conditionally defined probability that ($\mathscr{R}^{\prime\prime}$ \textit{1}) is equivalent to is also equivalent to the probability, under the same measure $\textbf{P}^{\xi_{2n}}_{\Lambda_{2n}(z)} \big[ \text{ } \cdot \text{ } \big]$, that $\big\{ \text{ }     \partial \Lambda_n(z) \text{ } {\overset{h \geq k + j - j^{\prime}}{\longleftrightarrow}}  \text{ } \big\{ \text{ }  \partial \Lambda_{2n}(z) \cap \big\{ \text{ }  \xi_{2n}  \geq k + j - j^{\prime}   \text{ }     \big\}    \text{ } \big\}$ occurs, in (\textit{1.3}), the probability in (\textit{1.2}) is upper bounded with the same conditionally defined event in (\textit{1.2}), with the exception that $|h|$, instead of $h$, appears in the connectivity event, in (\textit{1.4)} an upper bound to (\textit{1.3}) is obtained from the connectivity event dependent upon $|h|$,

\begin{align*}
     \bigg\{              \partial \Lambda_n(z) \text{ } {\overset{|h| \geq k + j - j^{\prime}}{\longleftrightarrow}}  \bigg\{   \partial \Lambda_{2n}(z) \cap \big\{  \xi_{2n}   \geq k + j - j^{\prime}  \big\}         \bigg\}  \bigg\}   \text{ } \text{ , } 
\end{align*}

\noindent because,

\begin{align*}
      \bigg\{              \partial \Lambda_n(z) \text{ } {\overset{|h| \geq k + j - j^{\prime}}{\longleftrightarrow}}  \text{ } \bigg\{   \partial \Lambda_{2n}(z) \cap \big\{  \xi_{2n}   \geq k + j - j^{\prime}  \big\} \bigg\}           \bigg\}  \equiv    \bigg\{      \partial \Lambda_n(z) \text{ } {\overset{h \geq k + j - j^{\prime}}{\longleftrightarrow}}  \text{ }  \partial \Lambda_{2n}(z)      \bigg\}        \text{ } \text{ , } 
\end{align*}

\noindent while, for the second to last term obtained in the series of rearrangements above, in (\textit{1.5}) the equality amongst the connectivity events, namely,

\begin{align*}
     \bigg\{    \partial \Lambda_n(z) \text{ } {\overset{|h| \geq k + j - j^{\prime}}{\longleftrightarrow}}  \text{ } \bigg\{ \text{ }  \partial \Lambda_{2n}(z) \cap \big\{ \text{ } \xi_{2n}   \geq k + j - j^{\prime}  \text{ } \big\}  \text{ } \bigg\}  \bigg\}           \text{ } \text{ , } 
\end{align*}

\noindent with,

\begin{align*}
     \big\{    \partial \Lambda_n(z) \text{ }      {\overset{h \geq k + j - j^{\prime}}{\longleftrightarrow}} \text{ }  \partial \Lambda_{2n}(z)  \big\}      \text{ } \text{ , } 
\end{align*}

\noindent is applied to obtain an equivalent probability with $\textbf{P}^{\xi_{2n}}_{\Lambda_{2n}(z)} \big[ \text{ }   \partial \Lambda_n(z) \text{ }      {\overset{h \geq k + j - j^{\prime}}{\longleftrightarrow}} \text{ }  \partial \Lambda_{2n}(z)         \text{ }   \big]$, and, in (\textit{1.6}), an equivalent expression for the probability obtained in (\textit{1.5}) is provided by switching the sign of the height function, as,

\begin{align*}
 h  < k + j - j^{\prime}  \text{ } \text{ , } 
\end{align*}

\noindent is the value of the height function required for the connectivity event to occur, in comparison to,

\begin{align*}
  h  \geq k + j - j^{\prime} \text{ } \text{ . } 
\end{align*}

\noindent To conclude the inductive argument, observe,

\begin{align*}
  (\textit{1.6}) \overset{(\mathrm{CBC})}{\leq}        \textbf{P}^{-\xi_{2n}+j}_{\Lambda_{2n}(z)} \big[            \partial \Lambda_n(z)    {\overset{h < k + j - j^{\prime}}{\longleftrightarrow}}  \partial \Lambda_{2n}(z)  \big]  \tag{\textit{1.7}}\\   \leq \textbf{P}^{-\xi_{2n}+j}_{\Lambda_{2n}(z)} \big[   \mathcal{O}_{h < k + j - j^{\prime}}(n)    \big]                     \tag{\textit{1.8}}  \\  \equiv     1 - \textbf{P}^{\xi_{2n}-j}_{\Lambda_{2n}(z)} \big[      \mathcal{O}^{\mathrm{x}}_{h < k + 2j - j^{\prime} }(n)                                    \big]                \text{ }           \tag{\textit{1.9}}                    \\ \leq \text{ }      1 - \textbf{P}^{\xi_{2n}-j}_{\Lambda_{2n}(z)} \big[           \mathcal{O}^{\mathrm{x}}_{h < k + 3j - j^{\prime} }(n)          \big]               \text{ }                              \tag{\textit{1.10}}  \\      \overset{h \longleftrightarrow -h}{\equiv}  1 - \textbf{P}^{\xi_{2n}-j}_{\Lambda_{2n}(z)} \big[                  \mathcal{O}^{\mathrm{x}}_{h \geq k + 3j - j^{\prime} }(n)                 \big]                   \tag{\textit{1.11}} \\    \overset{(\mathrm{CBC})}{\leq}  1 -\textbf{P}^{k,j^{\prime}}_{\Lambda_{2n}(z)}          \big[       \mathcal{O}^{\mathrm{x}}_{h \geq k + 3j - j^{\prime} }(n)     \big]                  \tag{\textit{1.12}} \\    \text{ }                    \leq \mathrm{exp} \big( \text{ } - c^{\prime} \text{ } \big)                          \tag{\textit{1.13}} \\  \Updownarrow \\                    \textbf{P}^{\xi^{\mathrm{Sloped}}}_D \bigg[              \mathcal{O}_{h\geq k+j}^{\mathrm{x}}    ( n )                 \big|      \mathcal{O}_{h \geq k - 2j}^{\mathrm{x}}  \cup  \big\{          \mathscr{D}^{\prime\prime} \equiv D          \big\} \bigg]           \leq  \mathrm{exp} \big( \text{ } - c^{\prime}  \text{ } \big)   \text{ } \text{ , } \tag{\textit{1.14}}
\end{align*}

\noindent which demonstrates that the induction hypothesis holds for all possible $\mathscr{D}^{\prime\prime}$, in which, from the sequence of rearrangements above, in (\textit{1.7}) (CBC) is applied so that the probability given in (\textit{1.6}) can be upper bounded with a probability, under the same support, and crossing event, except under $-\xi_{2n}+j$, in (\textit{1.8}) an upper bound to (\textit{1.7}) is provided, in light of the observation that,

\begin{align*}
  \mathcal{O}_{h < k + j - j^{\prime}}(n)     \subset     \big\{  \partial \Lambda_n(z) \text{ }      {\overset{h < k + j - j^{\prime}}{\longleftrightarrow}} \text{ }  \partial \Lambda_{2n}(z) \big\}       \text{ } \text{ , } 
\end{align*}

\noindent in (\textit{1.9}) an equivalent probability to what is provided in (\textit{1.8}) is obtained from the similar observation, from the rearrangement in (\textit{1.7}), that,

\begin{align*}
    \big\{   \mathcal{O}^{\mathrm{x}}_{h < k + j - j^{\prime}}(n)           \big\}^c   \equiv         \mathcal{O}_{h < k + 2j - j^{\prime}}(n)         \text{ } \text{ , } 
\end{align*}

\noindent as the $\mathrm{x}$-annulus crossing probability is more likely to occur across the smaller height $k+j-j^{\prime}$, in (\textit{1.10}) an upper bound for (\textit{1.9}) is obtained, in light of the observation that,

\begin{align*}
\mathcal{O}^{\mathrm{x}}_{h < k + 3j - j^{\prime}}(n)  \subset \mathcal{O}^{\mathrm{x}}_{h < k + j - j^{\prime}}(n)  \text{ } \text{ , } 
\end{align*}

\noindent in (\textit{1.11}) an equivalent random variable that is provided in (\textit{1.10}) is obtained under the symmetry $h \longleftrightarrow -h$, in (\textit{1.12}), from a second application of (CBC), under boundary conditions $k,j^{\prime} \equiv \xi^{\mathrm{Sloped}}_{k,j^{\prime}}$,

\begin{align*}
 \textbf{P}^{k,j^{\prime}}_{\Lambda_{2n}(z)}          \big[       \mathcal{O}^{\mathrm{x}}_{h \geq k + 3j - j^{\prime} }(n)   \big]      \leq  \textbf{P}^{\xi_{2n}}_{\Lambda_{2n}(z)}          \big[       \mathcal{O}^{\mathrm{x}}_{h \geq k + 3j - j^{\prime} }(n)       \big]                            \text{ , } 
\end{align*}

\noindent in (\textit{1.13}), the final exponential upper bound, in $c^{\prime}$ from the strict positivity of the annulus crossing across restricted scales, is obtained, hence concluding arguments for showing that the induction hypothesis holds, and in (\textit{1.14}), the exponential upper bound in $c^{\prime}$ to the conditionally defined crossing event, from the $\mathrm{x}$-annulus crossing probability, is provided.

\bigskip

\noindent Hence,

\begin{align*}
 \textbf{P}^{\xi^{\mathrm{Sloped}}} \big[     \mathcal{O}^{\mathrm{x}}_{h \geq 2i}  (n)          \big]  \leq \mathrm{exp} \big( \text{ }  - c^{\prime}_i k \text{ } \big)   \Rightarrow    \textbf{P}^{\xi^{\mathrm{Sloped}}} \big[   \mathcal{O}^{\mathrm{x}}_{h \geq 2i + 1}   (n)   \big]   \leq \mathrm{exp} \big( \text{ }  - c^{\prime}_{i+1} k \text{ } \big)   \text{ } \text{ , } 
\end{align*}

\noindent for $c^{\prime}_{i+1} > c^{\prime}_{i} > 0$, and,

\begin{align*}
  \textbf{P}^{\xi^{\mathrm{Sloped}}}_D \big[    \partial D \overset{|h| \geq k}{\longleftrightarrow}   \Lambda_n    \big] \geq 1 -  C^{\prime}  \mathrm{exp} \big( \text{ }  - c^{\prime} k     \text{ } \big)   \text{ } \text{ , } 
\end{align*}

\noindent admits the exponential lower bound provided above, which from the induction hypothesis follows by applying the same argument to,

\begin{align*}
 \frac{1}{C^{\prime}} \text{ }    \textbf{P}^{\xi^{\mathrm{Sloped}}}_D \bigg[             \mathcal{O}_{h\geq k+j}^{\mathrm{x}}    ( n )                \big|      \mathcal{O}_{h \geq k - 2j}^{\mathrm{x}}  \cup \big\{          \mathscr{D}^{\prime\prime} \equiv D          \big\} \bigg]  \text{ } \text{ , } 
\end{align*}

\noindent for strictly positive $C^{\prime}$, which is upper bounded with,

\begin{align*}
 C^{\prime} \text{ } \mathrm{exp} \big( \text{ } - c^{\prime} \text{ } \big)          \text{ } \text{ . } 
\end{align*}

\noindent As a result of the above induction hypothesis holding for all $k$ and realizations of $\mathscr{D}^{\prime\prime}$, fix $\mathscr{D}^{\prime\prime}$ so that $0 \in \mathscr{D}^{\prime\prime}$. The sloped conditional expectation of the square of the height function at the origin, 

\begin{align*}
  \textbf{E}^{\xi^{\mathrm{Sloped}}}_D \bigg[ h(0)^2  \big|          \Omega \equiv \mathscr{D}^{\prime\prime} , h \equiv \xi_{\mathscr{D}} \text{ over } C_k      \bigg]      \equiv         \textbf{E}_{\mathscr{D}^{\prime\prime}}^{\xi^{\mathrm{Sloped}}_{i^{\prime}}, \xi^{\mathrm{Sloped}}_{i^{\prime}+1}}   \big[         h(0)^2    \big]  \\ \equiv    \textbf{E}_{\mathscr{D}^{\prime\prime}}^{\xi^{\mathrm{Sloped}}_{i^{\prime}}-k^{\prime}, \xi^{\mathrm{Sloped}}_{i^{\prime}+1}-k^{\prime}-1} \big[    \big( h(0)  +       k^{\prime}         \big)         \big]\text{ } \\ \equiv    \text{ }     \textbf{E}_{\mathscr{D}^{\prime\prime}}^{\xi^{\mathrm{Sloped}}_{i^{\prime}}-k^{\prime}, \xi^{\mathrm{Sloped}}_{i^{\prime}+1}-k^{\prime}-1} \big[   h(0)^2       \big] +    \textbf{E}_{\mathscr{D}^{\prime\prime}}^{\xi^{\mathrm{Sloped}}_{i^{\prime}}-k^{\prime}, \xi^{\mathrm{Sloped}}_{i^{\prime}+1}-k^{\prime}-1} \big[    2 k^{\prime} h(0)   \big] +   \big(  k^{\prime} \big)^2 \\      \equiv   \textbf{E}_{\mathscr{D}^{\prime\prime}}^{\xi^{\mathrm{Sloped}}_{i^{\prime}}-k^{\prime}, \xi^{\mathrm{Sloped}}_{i^{\prime}+1}-k^{\prime}-1} \big[  h(0)^2     \big]  +        2 k^{\prime}+    \big(  k^{\prime}  \big)^2  \\       {\leq} \text{ }   w^{\xi^{\mathrm{S}}}_n \text{ } + 2k^{\prime} +     \big(  k^{\prime}  \big)^2 \text{ } \text{ , }  \tag{$w^{\xi^{\mathrm{S}}}_n$ \textit{upper bound}}
\end{align*}

\noindent where, in the ultimate expression from the upper bound above, the supremum over the sloped boundary conditions for which the expectation value under boundary conditions $\xi^{\mathrm{Sloped}}_{i^{\prime}}-k^{\prime}, \xi^{\mathrm{Sloped}}_{i^{\prime}+1}-k^{\prime}-1$ occurs is bounded above by the supremum of the expectation of the square of the height function at the origin. Hence,

\begin{align*}
    ( w^{\xi^{\mathrm{S}}}_n \textit{upper bound} )        \text{ } \leq \text{ } w^{\xi^{\mathrm{S}}}_n \text{ } + \text{ }     \textbf{E}^{\xi^{\mathrm{S}}}_{\mathscr{D}^{\prime\prime}} \bigg[ \bigg[  2 K^{\prime}  +  \big(  K^{\prime} \big)^2 \bigg]  \bigg]      \text{ } \text{ , } \tag{$w^{\xi^{\mathrm{S}}}_n$ \textit{upper bound II}}
\end{align*}

\noindent which, as a superposition of a sloped expectation with $w^{\xi^{\mathrm{S}}}_n$, is also bound below some suitably chosen constant $C_{k^{\prime}}$, because,

\begin{align*}
 \textbf{E}^{\xi^{\mathrm{S}}}_{\mathscr{D}^{\prime\prime}} \bigg[ \bigg[  2 K^{\prime}  +  \big( K^{\prime}  \big)^2 \bigg]  \bigg]   \equiv \bigg[   2 K^{\prime} + (  K^{\prime}  )^2  \bigg]  \underset{ \mathscr{D}^{\prime\prime}}{\sum}     \mathrm{d} \textbf{P}^{\xi^{\mathrm{S}}}_{\mathscr{D}^{\prime\prime}} > 0  \text{ } \text{ , } 
\end{align*}

\noindent after which ($w^{\xi^{\mathrm{S}}}_n$ \textit{upper bound}) into $\textbf{E}^{\xi^{\mathrm{Sloped}}}_D[h(0)^2]$, which was previously decomposed over pairs of admissible domains and boundary conditions, with,

\begin{align*}
  \underset{\text{ boundary conditions } \big(  \xi_{\mathscr{D}}  \big) }{\underset{\text{ domains } \mathscr{D}^{\prime\prime} }{\sum}}    \textbf{P}^{\xi^{\mathrm{Sloped}}}_D \bigg[ h(0)^2  \big|  \Omega \equiv \mathscr{D}^{\prime\prime} , h \equiv \xi_{\mathscr{D}} \text{ over } C_k  \bigg]  \textbf{P}^{\xi^{\mathrm{Sloped}}}_D \big[   \Omega \equiv \mathscr{D}^{\prime\prime} , h \equiv \xi_{\mathscr{D}} \text{ over } C_k     \big]     \text{ } \text{ . } 
\end{align*}

\noindent Otherwise, there exists a large enough constant so that $h(0)^2$ can be bounded from above if $0 \notin \mathscr{D}^{\prime\prime}$, from which we conclude the argument. \boxed{}

\subsection{Simultaneously incorporating arguments from upper and lower bound results to the torus from simply connected domains}

\noindent Denote the faces of the torus $\textbf{T}$ with $F ( \textbf{T})$. For any such $F$, the sloped expectation for any face over $\textbf{T}$ is determined by the restriction of the height function to the face of contact, as,

\begin{align*}
      \textbf{E}^{\xi^{\mathrm{Sloped}}}_F [ \text{ } \cdot \text{ } ] \text{ } \equiv \text{ } \textbf{E}^{(\mathrm{bal})}_{\textbf{T}} \bigg[ \text{ } \cdot \big| \text{ } h|_F \in \{ \xi^{\mathrm{Sloped}} - 1 , \xi^{\mathrm{Sloped}} \}  \bigg]      \text{ } \text{ , } 
\end{align*}

\noindent whose behavior we will similarly characterize with arguments for two additional crossing events, each of which are dependent upon the absolute value of the height function.

\bigskip

\noindent Next, for $x,y \in F ( \textbf{T})$, introduce,

\begin{align*}
         w^{\xi^{\mathrm{Sloped}}}_n \equiv w^{\xi^{\mathrm{S}}}_n \equiv  \text{ }        \underset{F}{\mathrm{sup}} \big\{        \textbf{E}^{\xi^{\mathrm{Sloped}}}_{F^c} [ h(x)^2 ] :     F \subset F ( \textbf{T} ) \text{ connected , with, } F \cap \Lambda_{n}(x) \neq \emptyset , |F| \geq 2n      \big\}                      \text{ } \text{ , } 
\end{align*}

\noindent for the sloped expectation of the square of the height function at the origin supported over the complement of the faces $F$ in $\textbf{T}$, and,

\begin{align*}
        u^{\xi^{\mathrm{Sloped}}}_n \equiv u^{\xi^{\mathrm{S}}}_n \equiv  \text{ }          \underset{F}{\mathrm{sup}} \big\{    \textbf{E}^{\xi^{\mathrm{Sloped}}}_{F^c} [ h(x)^2 ] :     F \subset F(\textbf{T}) , F \cap \partial \Lambda_n (y) \neq \emptyset  \big\}                        \text{ } \text{ , } 
\end{align*}

\noindent for the sloped expectation of the square of the height function at the origin, which is also supported over the complement of the faces $F$ in $\textbf{T}$, which satisfies a different condition than what is given in $w^{\xi^{\mathrm{Sloped}}}_n$.

\bigskip

\noindent \textbf{Lemma} \textit{5.3} (\textit{induction arguments for two height function dependent crossing events}). There exists two constants $c$,$C^{\prime}$, coinciding with those provided in arguments for previous results, for which,

\begin{align*}
  \textbf{P}^{\xi^{\mathrm{Sloped}}}_{F^c} \big[  F \overset{|h| \leq k}{\longleftrightarrow}         \Lambda_n(x)       \big]  \geq   1 - C^{\prime} \mathrm{exp} \big( \text{ } -c^{\prime} k \text{ } \big)        \text{ }  \text{ , } 
\end{align*}

\noindent and also for which,

\begin{align*}
  \textbf{P}^{\xi^{\mathrm{Sloped}}}_{F^c} \big[  F \overset{|h| \leq k}{\longleftrightarrow}   \partial \Lambda_{\frac{3}{2}n} (y) \big] \geq 1 -  C^{\prime} \mathrm{exp} \big( \text{ } -c^{\prime} k \text{ } \big)  \text{ } \text{ , } 
\end{align*}

\noindent under the assumption that $d(x,y) \leq \frac{N}{16}$ for $x,y \in F( \textbf{T})$. Inductively, arguing that each one of the two crossing probabilities holds, implies that,

\begin{align*}
  w^{\xi^{\mathrm{S}}}_{\frac{3}{2}n} \leq w^{\xi^{\mathrm{S}}}_n + C^{\prime}  \text{ } \text{ , } \\     u^{\xi^{\mathrm{S}}}_n \leq u^{\xi^{\mathrm{S}}}_{\frac{3}{2}n} + C^{\prime} \text{ } \text{ , } \\   u^{\xi^{\mathrm{S}}}_{\frac{3}{2}n} \text{ } \leq \text{ }  w^{\xi^{\mathrm{S}}}_n     \text{ } \text{ , } 
\end{align*}

 \noindent for $1 \leq n \leq \frac{N}{8}$.

\bigskip 

\noindent \textit{Proof of Lemma 5.3}. Along the lines of two induction arguments argued for the, respectively, lower, and upper, bounds in \textit{5.1} and \textit{5.2}, if both of the crossing events under the sloped probability measure are to be exponentially bounded from below, then,

\begin{align*}
  \textbf{P}^{\xi^{\mathrm{Sloped}}}_{F^c} \big[    \mathcal{O}^{\mathrm{x}}_{|h| \geq k} (n )          \big]      \geq      1 - C^{\prime\prime} \mathrm{exp} \big( \text{ } - c^{\prime\prime} k \text{ } \big)    \text{ } \text{ , } 
\end{align*}

\noindent corresponding to the first condition of the above \textbf{Lemma} would hold, in addition to,

\begin{align*}
   \textbf{P}^{\xi^{\mathrm{Sloped}}}_{F^c} \big[        \mathcal{O}^{\mathrm{x}}_{|h| \geq k} (2n )          \big]         \geq             1 - C^{\prime\prime\prime} \mathrm{exp} \big( \text{ } - c^{\prime\prime\prime} k \text{ } \big)          \text{ }  \text{ , } 
\end{align*}

\noindent corresponding to the second condition of the above would also hold, for $c^{\prime\prime}, c^{\prime\prime\prime} , C^{\prime\prime} , C^{\prime\prime\prime} > 0$. As a brief remark, each of the induction hypotheses above are dependent upon $\mathrm{x}$-annulus crossing events, and hence the directionality of the crossing upon the absolute value of the height function is reversed. For the first claim, to demonstrate that the claim holds for $i$ first, begin with the conditional sloped probability,

\begin{align*}
      \textbf{P}^{\xi^{\mathrm{Sloped}}}_{F^c} \bigg[    \mathcal{O}^{\mathrm{x}}_{|h| \geq k+2}(n)  \big|  \mathcal{O}^{\mathrm{x}}_{|h| \geq k}(n) \cup \big\{    \mathscr{D}^{\prime\prime\prime} \equiv D^{\prime}  \big\}      \bigg] \text{ }            \text{ } \text{ , } \tag{\textit{1}}
\end{align*}

\noindent where, as in previous argument, $\mathscr{D}^{\prime\prime\prime}$ denotes the $\mathrm{x}$-loop that is closest to the finite volume boundary. Further rewriting the conditional probability above, for $k\equiv 2$, yields,

\begin{align*}
      \textbf{P}^{\xi^{\mathrm{Sloped}}}_{\mathscr{D}^{\prime\prime\prime}} \big[   \mathcal{O}^{\mathrm{x}}_{|h| \geq 4} (n)        \big]   \equiv  \textbf{P}^{\xi^{\mathrm{Sloped}}-1, \xi^{\mathrm{Sloped}}+1}_{\mathscr{D}^{\prime\prime\prime}} \big[              \mathcal{O}^{\mathrm{x}}_{|h| \geq 4}(n)     \big]    \text{ } \text{ , } \tag{\textit{2}}
\end{align*}

\noindent which, upon adapting familiarly oriented finite volumes that can be centered about some $z$ on $\partial \Lambda_{2n}(z)$, yields,

\begin{align*}
     \textbf{P}^{\xi^{\mathrm{Sloped}}-1, \xi^{\mathrm{Sloped}}+1}_{\mathscr{D}^{\prime\prime\prime}} \big[         \mathcal{O}^{\mathrm{x}}_{|h| \geq 4}(n)     \big]    \leq   \textbf{P}^{\xi^{\mathrm{Sloped}}-1, \xi^{\mathrm{Sloped}}+1}_{\mathscr{D}^{\prime\prime\prime}} \big[   \partial \Lambda_{n}(z) {\overset{|h| \geq 4}{\longleftrightarrow}}_{\mathrm{x}}   \partial \Lambda_{2n} (z)   \big] \tag{$\mathscr{D}^{\prime\prime\prime}$ \textit{probability}} \\ \equiv  1 -  \textbf{P}^{\xi^{\mathrm{Sloped}}-1, \xi^{\mathrm{Sloped}}+1}_{\mathscr{R}^{\prime\prime\prime}}  \bigg[   \mathcal{O}_{|h| \leq 4}\big(   A ( n , 2n) + z  \big)             \big|   \xi^{\mathrm{Sloped}}-2   \leq h \leq    \xi^{\mathrm{Sloped}}  \text{ over } \Lambda_{2n}(z) \backslash \mathscr{D}^{\prime\prime\prime} \bigg]   \text{ } \text{ , } \tag{$\mathscr{R}^{\prime\prime\prime}$ \textit{probability}}
\end{align*}

\noindent in which, from the sequence of rearrangements above, ($\mathscr{R}^{\prime\prime\prime}$ \textit{probability}) is obtained from ($\mathscr{D}^{\prime\prime\prime}$ \textit{probability}) upon modification of the support of the probability measure from $\mathscr{D}^{\prime\prime\prime}$ to $\mathscr{R}^{\prime\prime\prime}$. With the following, further manipulation shows that the induction hypothesis holds, in which one of the previously obtained probabilities would satisfy,

\begin{align*}
    (\mathscr{D}^{\prime\prime\prime} \textit{probability}) \text{ } \overset{(\mathrm{CBC})}{\leq} \textbf{P}^{\xi^{\mathrm{Sloped}-k^{\prime}},\xi^{\mathrm{Sloped}+k^{\prime}}}_{\mathscr{D}^{\prime\prime}} \big[    \partial \Lambda_{n}(z) {\overset{|h| \geq 4}{\longleftrightarrow}}_{\mathrm{x}}    \partial \Lambda_{2n} (z) \big]  \\  {\leq}  \textbf{P}^{\xi^{\mathrm{Sloped}-k^{\prime}},\xi^{\mathrm{Sloped}+k^{\prime}}}_{\mathscr{D}^{\prime\prime\prime}} \bigg[     \partial \Lambda_{n}(z) {\overset{|h| \geq 4}{\longleftrightarrow}}    \partial \Lambda_{2n} (z)   \big|   1  \leq | h | \leq      2    \text{ over } \Lambda_{2n}(z) \backslash    \widetilde{\mathscr{D}^{\prime\prime\prime}}    \bigg] \\ \equiv \textbf{P}^{\xi^{\mathrm{Sloped}-k^{\prime}},\xi^{\mathrm{Sloped}+k^{\prime}}}_{\mathscr{D}^{\prime\prime\prime}} \bigg[   \partial \Lambda_{n}(z) {\overset{|h| \geq 4}{\longleftrightarrow}}    \partial \Lambda_{2n} (z)  \big|  1 \leq | h | \leq         2  \text{ over } \Lambda_{2n}(z ) \backslash        \widetilde{\mathscr{D}^{\prime\prime\prime}}   \bigg] \tag{\textit{3}} \\  \equiv     \textbf{P}^{\xi^{\mathrm{Sloped}-k^{\prime},+k^{\prime}}}_{\mathscr{D}^{\prime\prime\prime}} \bigg[   \partial \Lambda_{n}(z) {\overset{|h| \geq 4}{\longleftrightarrow}}   \partial \Lambda_{2n} (z)   \big|  1 \leq | h | \leq       2  \text{ over } \Lambda_{2n}(z)   \backslash        \widetilde{\mathscr{D}^{\prime\prime\prime}}      \bigg]        \tag{\textit{3}-\textit{2}} \\  {\leq}    \textbf{P}^{\xi^{\mathrm{Sloped}-k^{\prime},+k^{\prime}}}_{\mathscr{D}^{\prime\prime\prime}} \bigg[     \partial \Lambda_{n}(z) {\overset{|h| \geq 4}{\longleftrightarrow}} \mathcal{C}^{\prime} \text{ } \big|  1 \leq | h | \leq        2  \text{ over } \Lambda_{2n}(z) \text{ }  \backslash         \widetilde{\mathscr{D}^{\prime\prime\prime}}     \bigg]     \tag{\textit{4}}\\
    \overset{(\mathrm{SMP})}{\equiv}   \textbf{P}^{\xi^{\mathrm{Sloped}-k^{\prime}_{2n} ,+k^{\prime}_{2n}}}_{\Lambda_{2n}(z)} \bigg[         \partial \Lambda_{n}(z) {\overset{|h| \geq 4}{\longleftrightarrow}}  \bigg\{   \partial \Lambda_{2n} (z) \text{ } \cap  \big\{       \xi^{\mathrm{Sloped}-k^{\prime},+k^{\prime}} -2          \leq \xi^{\mathrm{Sloped}-k^{\prime},+k^{\prime}}\leq  \xi^{\mathrm{Sloped}-k^{\prime},+k^{\prime}} + 1          \big\}   \bigg\}                \bigg]          \tag{\textit{5}}   \\      \leq    \textbf{P}^{\xi^{\mathrm{Sloped}-(k^{\prime}_{2n} )^{\prime},+(k^{\prime}_{2n})^{\prime}}}_{\Lambda_{2n}(z)} \big[         \partial \Lambda_{n}(z) {\overset{|h| \geq 4}{\longleftrightarrow}}     \partial \Lambda_{2n} (z)         \big]     \tag{\textit{6}} \text{ , } 
\end{align*}

\noindent for sufficiently chosen $k^{\prime}_{2n}$ and $\big( k^{\prime}_{2n} \big)^{\prime}$, in which, from the sequence of rearrangements above, from inductive arguments throughout \textit{5.1} and \textit{5.2}, from the probability in (\textit{1}), associating boundary conditions to the $\mathrm{x}$-annulus crossing in (\textit{2}) then allows us to obtain (\textit{3}), which is a conditionally defined crossing from the absolute value of the height function, and the values of the height function over $\Lambda_{2n}(z) \backslash \widetilde{\mathscr{D}^{\prime\prime\prime}}$. Following (\textit{3}), (\textit{3}-\textit{2}) is obtained by shortening the notation for the boundary conditions on the probability measure with $\xi^{\mathrm{Sloped}}-k^{\prime},\xi^{\mathrm{Sloped}}+k^{\prime} \equiv \xi^{\mathrm{Sloped}-k^{\prime},+k^{\prime}}$, from which an upper bound for (\textit{4}) is obtained with (\textit{5}) after removing the conditioning,

\begin{align*}
   \big\{     1 \leq | h | \leq 2 \text{ over } \Lambda_{2n}(z) \backslash \widetilde{\mathscr{D}^{\prime\prime\prime}}    \big\}  \text{ } \text{ , } 
\end{align*}

\noindent also, where, in (\textit{4}), the connectivity event is,

\begin{align*}
\mathcal{C}^{\prime} \equiv  \bigg\{    \partial \Lambda_{2n} (z)  \cap \bigg\{ \text{ }      \xi^{\mathrm{Sloped}-k^{\prime},+k^{\prime}} -2          \leq \xi^{\mathrm{Sloped}-k^{\prime},+k^{\prime}}\leq  \xi^{\mathrm{Sloped}-k^{\prime},+k^{\prime}} + 1          \bigg\} \bigg\}     \text{ } \text{ , } 
\end{align*}

\noindent from the pushforward that is computed under $\textbf{P}^{\xi^{\mathrm{Sloped}-k^{\prime},+k^{\prime}}}_{\mathscr{D}^{\prime\prime\prime}}[  \cdot  ]$, from which a final upper bound is given with (\textit{6}), which does not include conditioning upon the value of the boundary condition $\xi^{\mathrm{Sloped}-k^{\prime},+k^{\prime}}$.

\bigskip

\noindent Due to the fact that the support of the probability measure in (\textit{6}) is also supported over $\Lambda_{2n}(z)$, we are in a position to apply the final sequence of arguments as provided in (\textit{1.7})-(\textit{1.14}) for the proof of \textbf{Lemma} \textit{5.3}, as, given a similar choice of parameters,

\begin{align*}
   \textbf{P}^{\xi^{\mathrm{Sloped}-k^{\prime},+k^{\prime}}}_{\Lambda_{2n}(z)} \big[       \partial \Lambda_{n}(z) {\overset{|h| \geq 4}{\longleftrightarrow}}    \partial \Lambda_{2n} (z)       \big]  \text{ } \overset{(\mathrm{CBC})}{\leq}          \textbf{P}^{\xi^{\mathrm{Sloped}-k^{\prime}+j,+k^{\prime}+j}}_{\Lambda_{2n}(z)} \big[        \partial \Lambda_{n}(z) {\overset{|h| \geq 4}{\longleftrightarrow}} \text{ }    \partial \Lambda_{2n} (z)           \big]                   \text{ }      \\ \leq \textbf{P}^{\xi^{\mathrm{Sloped}-k^{\prime}+j,+k^{\prime}+j}}_{\Lambda_{2n}(z)} \big[     \mathcal{O}_{|h| \geq 4}   (n)  \big] \text{ } \\ \text{ } \equiv 1 - \text{ } \textbf{P}^{\xi^{\mathrm{Sloped}-k^{\prime}+j,+k^{\prime}+j}}_{\Lambda_{2n}(z)} \big[      \mathcal{O}_{|h| \geq 4}(n)    \big]   \text{ } \\ \leq   1 - \text{ }   \textbf{P}^{\xi^{\mathrm{Sloped}-k^{\prime}+j,+k^{\prime}+j}}_{\Lambda_{2n}(z)} \big[      \mathcal{O}_{|h| \geq 5}(n)        \big] \\
   \equiv \text{ }     1 - \textbf{P}^{\xi^{\mathrm{Sloped}-k^{\prime}+j,+k^{\prime}+j}}_{\Lambda_{2n}(z)} \big[    \mathcal{O}_{|h| < 5}(n)         \big]   \text{ , } \tag{\textit{1.11 II}}
\end{align*}

\noindent where, in the sequence of rearrangements above we implement the steps given in (\textit{1.7})-(\textit{1.11}), in which we begin with an application of (CBC). The remaining steps yield an exponential upper bound, for some $c^{\prime\prime}$ so that the third induction hypothesis holds. With a second application of (CBC),

\begin{align*}
    (\textit{1.11 II}) \leq \text{ } 1 -     \textbf{P}^{k^{\prime},j^{\prime}}_{\Lambda_{2n}(z)} \big[      \mathcal{O}_{|h| < 5}(n)        \big] \leq \mathrm{exp} \big( \text{ } -c^{\prime\prime} \text{ } \big)       \text{ } \text{ , } \tag{\textit{Hypothesis II}}
\end{align*}

\noindent from which the statement in the first condition is shown to hold. The same justification for the sequence of rearrangements above has been mentioned in arguments for a previous \textbf{Lemma}, with the only difference in the boundary conditions and parameters that are used to change boundary conditions, and the threshold value of the height of the annulus events, which is left to the reader.

\bigskip

\noindent For the second condition of the \textbf{Lemma}, the induction arguments are similarly applied, from the observation that repeating the same line of observations as given above, for $\mathcal{O}^{\mathrm{x}}(2n)$ instead of $\mathcal{O}^{\mathrm{x}}(n)$. That is, with the existence of another exterior most $\mathrm{x}$-loop as denoted with $\mathscr{D}^{\prime\prime\prime}$ when arguing that the previous induction statement holds, rearrangements would yield an inequality of the form,

\begin{align*}
      \textbf{P}^{\xi^{\mathrm{Sloped}-k^{\prime\prime},+k^{\prime\prime}}}_{\Lambda_{4n}(z)} \big[     \partial \Lambda_{2n}(z) {\overset{|h| \geq 4}{\longleftrightarrow}}  \partial \Lambda_{4n} (z)        \big]   \overset{(\textit{1.11 III})}{\leq}  \mathrm{exp} \big( \text{ } -c^{\prime\prime\prime} \text{ } \big)      \text{ } \text{ , } \tag{\textit{Hypothesis III}}
\end{align*}

\noindent in which (\textit{1.11 III}) indicates,

\begin{align*}
       \textbf{P}^{\xi^{\mathrm{Sloped}-k^{\prime\prime},+k^{\prime\prime}}}_{\Lambda_{4n}(z)} \big[      \partial \Lambda_{2n}(z) {\overset{|h| \geq 4}{\longleftrightarrow}} \text{ }    \partial \Lambda_{4n} (z)         \big]        \leq      1 - \textbf{P}^{\xi^{\mathrm{Sloped}-k^{\prime\prime}+j^{\prime},+k^{\prime\prime}+j^{\prime}}}_{\Lambda_{4n}(z)} \big[    \mathcal{O}_{|h| < 5}(2n)      \big]          \text{ , }   \text{ } \tag{\textit{1.11 III}}
\end{align*}

\noindent which is used to demonstrate that the exponential upper bound in $c^{\prime\prime\prime}$ holds, for other suitably chosen constants $k^{\prime\prime}$ and $j^{\prime\prime}$ appearing in the shift of boundary conditions throughout different steps of the argument. Altogether, the previous two induction statements demonstrate,

\begin{align*}
    (\textit{Hypothesis II}) \Leftrightarrow                      \textbf{P}^{\xi^{\mathrm{Sloped}}}_{F^c} \big[  F \overset{|h| \leq k}{\longleftrightarrow}         \Lambda_n(x)     \big] \geq   1 - C^{\prime\prime} \mathrm{exp} \big( \text{ } -c^{\prime\prime} k \text{ } \big)    \Leftrightarrow   w^{\xi^{\mathrm{S}}}_{\frac{3}{2}n} \leq w^{\xi^{\mathrm{S}}}_n + C^{\prime}  \text{ } \text{ , }   \\      (\textit{Hypothesis III}) \Leftrightarrow        \textbf{P}^{\xi^{\mathrm{Sloped}}}_{F^c} \big[  F \overset{|h| \leq k}{\longleftrightarrow}   \partial \Lambda_{2n} (y) \big]  \geq 1 -  C^{\prime\prime\prime} \mathrm{exp} \big( \text{ } -c^{\prime\prime\prime} k \text{ } \big)                \Leftrightarrow                    u^{\xi^{\mathrm{S}}}_n \leq u^{\xi^{\mathrm{S}}}_{\frac{3}{2}n} + C^{\prime}         \text{ } \text{ . } 
\end{align*}

\noindent To conclude the arguments, apply the induction argument as given in the second statement across a slightly larger scale, which would implicate an upper bound of the following form,

\begin{align*}
     \textbf{P}^{\xi^{\mathrm{Sloped}-k^{\prime\prime\prime},+k^{\prime\prime\prime}}}_{\Lambda_{8n}(z)} \big[        \partial \Lambda_{4n}(z) {\overset{|h| \geq 4}{\longleftrightarrow}} \text{ }    \partial \Lambda_{8n} (z)             \big]             \overset{(\textit{1.11 IV})}{\leq}    \mathrm{exp} \big( \text{ }  -c^{\prime\prime\prime\prime} \text{ } \big) \text{ } \text{ , } \tag{\textit{Hypothesis IV}}
\end{align*}

\noindent for suitable $c^{\prime\prime\prime\prime}$, $k^{\prime\prime\prime}$ and $j^{\prime\prime}$, where, as one can expect,

\begin{align*}
      \textbf{P}^{\xi^{\mathrm{Sloped}-k^{\prime\prime},+k^{\prime\prime}}}_{\Lambda_{8n}(z)} \big[     \partial \Lambda_{4n}(z) {\overset{|h| \geq 4}{\longleftrightarrow}} \text{ }    \partial \Lambda_{8n} (z)            \big]             \text{ } \leq     \text{ }   1 - \textbf{P}^{\xi^{\mathrm{Sloped}-k^{\prime\prime\prime}+j^{\prime\prime},+k^{\prime\prime}+j^{\prime\prime}}}_{\Lambda_{8n}(z)} \big[    \mathcal{O}_{|h| < 5}(2n)       \big]                        \text{ } \text{ . } \tag{\textit{1.11 IV}}
\end{align*}

\noindent In addition to the previous two statements, this readily implies, for suitable $C^{\prime\prime\prime}$,

\begin{align*}
   (\textit{Hypothesis IV}) \Leftrightarrow                      \textbf{P}^{\xi^{\mathrm{Sloped}}}_{F^c} \big[ F \overset{|h| \leq k}{\longleftrightarrow}         \Lambda_{8n}(x)   \big] \text{ } \geq   1 - C^{\prime\prime\prime} \mathrm{exp} \big( \text{ } -c^{\prime\prime} k \text{ } \big)    \Leftrightarrow           u^{\xi^{\mathrm{S}}}_{\frac{3}{2}n} \text{ } \leq \text{ }  w^{\xi^{\mathrm{S}}}_n                         \text{ } \text{ , } 
\end{align*}

\noindent from which we conclude the argument. \boxed{}

\bigskip

\noindent \textit{Proof of the upper bound for the sloped variance in Corollary 6V 1 for arbitrary domains}. We appropriate the argument for the upper bound from {\color{blue}[11]}, in light of previously obtained results concerning $w^{\xi^{\mathrm{S}}}_n$ and $u^{\xi^{\mathrm{S}}}_n$. Fix a face $x$ of a simple connected domain $D$. Then, the sloped variance,

\begin{align*}
        \mathrm{Var}^{\xi^{\mathrm{Sloped}}}_D [      h(x)   ]   \equiv \mathrm{Var}^{\xi^{\mathrm{Sloped}}+l^{\prime}}_D  [ h(x) ]  \text{ } \tag{$l^{\prime}$ \textit{Variance}}
\end{align*}

\noindent upon shifting the boundary conditions by some $l^{\prime\prime}>0$, from which,

\begin{align*}
    (l^{\prime\prime} \textit{Variance}) \leq \text{ }  \textbf{E}^{\xi^{\mathrm{Sloped}+l^{\prime}}}_D [           h(x)^2                ] \text{ } \overset{(\mathrm{CBC})}{\leq} \text{ } \textbf{E}^{\chi^{\mathrm{Sloped}}+l^{\prime}}_D [ h(x)^2] \text{ } \text{ , } \tag{\textit{CBC 1}}
\end{align*}

\noindent for $\chi^{\mathrm{Sloped}} \geq \xi^{\mathrm{Sloped}}$ minimizing the sloped expectation,

\begin{align*}
 \textbf{E}^{\chi^{\mathrm{Sloped}}+l^{\prime}}_D [ h(x)^2]   \text{ } \text{ . } 
\end{align*}

\noindent Next,

\begin{align*}
 (\textit{CBC 1}) \equiv \text{ }  \mathrm{Var}^{\chi^{\mathrm{Sloped}}+l^{\prime}}_D [  h(x)    ]  + \textbf{E}^{\chi^{\mathrm{Sloped}}+l^{\prime}}_D [h(x)  ]^2 \text{ , } 
\end{align*}

\noindent in addition to the fact that the variance term from the equivalent formulation of (\textit{CBC 1}) above admits the following upper bound for $\chi^{\mathrm{Sloped}}$ introduced previously,

\begin{align*}
     \mathrm{Var}^{\xi^{\mathrm{Sloped}}}_D [      h(x)   ]   \leq     \mathrm{Var}^{\chi^{\mathrm{Sloped}}+l^{\prime}}_D [  h(x)    ] \text{ } + 4 |\xi^{\mathrm{Sloped}}|^2  \text{ } \text{ , } 
\end{align*}

\noindent in which the constant factor appearing in the upper bound with the sloped variance that is taken under $\chi^{\mathrm{Sloped}}$, which quantifies the change in variance of the height function that are respectively taken under $\xi^{\mathrm{Sloped}}$ and $\chi^{\mathrm{Sloped}}$. As a result, for some $y$ even,

\begin{align*}
    \mathrm{Var}^{\chi^{\mathrm{Sloped}}}_D [ h(x) ]  \leq \textbf{E}^{(\mathrm{bal})}_{\textbf{T}} [ h(x)^2 | \text{ } |h| \leq \chi^{\mathrm{Sloped}} \text{ over } \partial D ] \overset{(\mathrm{CBC-|h|})}{\leq}     \textbf{E}^{(\mathrm{bal})}_{\textbf{T}} [ h(x)^2  | | h | \leq \xi^{\mathrm{Sloped}} \text{ over } \partial D ] \\ \equiv \textbf{E}^{(\mathrm{bal})}_{\textbf{T}} [  (h(x)-h(y))^2 ]                \text{ } \text{ , } \tag{\textit{Squared height function difference}}
\end{align*}

\noindent given a suitable embedding of the simply connected domain $D$ into $\textbf{T}$, in which $N > 2 |D|$. Therefore,

\begin{align*}
   (\textit{Squared height function difference}) \leq C^{\prime} \text{ } \mathrm{log} [ d_{\textbf{T}}(x,y) ]   \text{ } \leq \text{ }  C^{\prime}\text{ }  \mathrm{log} [  d(x,y)  ] \text{ } \leq \text{ } C^{\prime} \text{ } \mathrm{log} [  d(x,\partial D) + 1     ]  \leq C^{\prime\prime}  \text{ } \text{ , } 
\end{align*}

\noindent for suitable $C^{\prime\prime}$ so that,

\begin{align*}
        \mathrm{log} [ d(x,\partial D)+1] \leq \frac{C^{\prime\prime}}{C^{\prime}}       \text{ } \text{ , } 
\end{align*}

\noindent from which we conclude the argument. \boxed{}

\bigskip

\noindent Finally, we also prove the upper bound for the delocalized phase under the sloped expectation. The sequence of inequalities directly mirrors arguments for obtaining the logarithmic delocalization result for the height function of the six-vertex model under sufficiently flat boundary conditions.

\bigskip

\noindent \textit{Proof of upper bound on the delocalized phase of the sloped expectation, as given in Theorem 1.1}. From arguments in \textbf{Lemma} \textit{5.3}, upper bound the sloped expectation with,

\begin{align*}
 \textbf{E}^{(\mathrm{bal})}_{\textbf{T}} [ (h(x)-h(y))^2] \leq u^{\xi^{\mathrm{S}}}_1 \leq            u^{\xi^{\mathrm{S}}}_{\frac{3}{2}} + C^{\prime} {\leq}     u^{\xi^{\mathrm{S}}}_{\frac{3}{2}n}          + C^{\prime}      \leq   u^{\xi^{\mathrm{S}}}_{\frac{3}{2}n}          + C^{\prime}  \text{ } \mathrm{log} \big[ \frac{3}{2} d \big] \leq w^{\xi^{\mathrm{S}}}_n +   C^{\prime}  \text{ } \mathrm{log} \big[ \frac{3}{2} d \big]  \\ \leq \text{ }      w^{\xi^{\mathrm{S}}}_1 +   C^{\prime}  \text{ } \mathrm{log} \big[ \frac{3}{2} d \big]         \text{ } \text{ , } 
\end{align*}

\noindent which yields the desired bound, as,

\begin{align*}
   w^{\xi^{\mathrm{S}}}_1 <  C^{\prime\prime\prime} \text{ } \text{ , } 
\end{align*}

\noindent for suitable $C^{\prime\prime\prime}$, hence implying,

\begin{align*}
        w^{\xi^{\mathrm{S}}}_1 +   C^{\prime}  \text{ } \mathrm{log} \big[ \frac{3}{2} d \big]   \leq     C^{\prime\prime\prime}  +      C^{\prime}  \text{ } \mathrm{log} \big[ \frac{3}{2} d \big]  \leq  \sup \big\{     C^{\prime\prime\prime}  , C^{\prime}     \big\} \bigg[  1 +  \mathrm{log} \big[ \frac{3}{2} d \big]   \bigg]  \leq C^{\prime\prime\prime,\prime\prime} \equiv C \\ \Updownarrow 
        \\ \textbf{E}^{(\mathrm{bal})}_{\textbf{T}} [ (h(x)-h(y))^2]     \leq     C   \text{ } \text{ , } 
\end{align*}

\noindent for a suitable upper bound $C$, from which we conclude the argument. \boxed{}

\section{Propagating RSW results from the strip environment to the spin-representation of the six-vertex model to the self-dual parameter line of the Ashkin-Teller model}

\subsection{Introduction}

\noindent From many variants of the six-vertex model, the spin-representation captured by the Ashkin-Teller model shares in correspondence with the six-vertex model, on the self-dual curve, defined via the equality $\mathrm{sinh}(2J) \equiv \mathrm{exp}(-2U)$, for $a \equiv b \equiv 1$ and $c \equiv \mathrm{coth}(J)$, where $U$ and $J$, respectively, are two coupling constants for two coupled Ising models. This choice of parameters for the correspondence between the six-vertex and Ashkin-Teller models satisfies the properties of the previously presented arguments for sloped boundary conditions, as $a\equiv b\equiv 1$ implies that \textit{symmetric} domains in the strip are invariant under horizontal and vertical rotations, from the symmetry $\sigma \xi \geq -\xi$, in addition to $1 \leq c \leq 2$ if the coupling $J$ is properly tuned. Under this parameter choice for the mixed Ashkin-Teller model, the parameter range specified previously for the six-vertex model is also satisfied, in which $\mathrm{max} \big\{ a,b \} \leq c$. From such a representation of $+$ and $-$ face variables, properties from the six-vertex model, including ($\mathrm{CBC-|h|}$) and ($\mathrm{FKG-|h|}$), do not hold. Furthermore, there is no analog of sloped boundary conditions for the spin-representation, as the possible boundary conditions that can be imposed on the measure $\textbf{P}^{\mathrm{AT}} \big[ \cdot \big] \equiv \mathcal{P} \big[ \cdot \big] $ can either consist of completely $+$ faces distributed along the finite volume boundary, with $\mathcal{P}^{++}$, completely $-$ faces distributed along the finite volume boundary, with $\mathcal{P}^{--}$. With such conventions, the strip portions of the argument can be applied to obtain the first RSW result, from which different conditions on $\mathcal{P} \big[ \cdot \big]$ that are introduced to compensate for the lack of an analog for sloped boundary conditions, in addition to the absence of an analog for the absolute value of the height function, $|h|$.

From the set of possible boundary conditions for the Ashkin-Teller model, there does not exist an analog for flat, or sloped, boundary conditions, as well as a lack of an analog for $|h|$. Under the presence of $++$ or $--$ boundary conditions, quantifying sufficiently good probabilities so that a different form of crossing events can be obtained is poossible by directly appealing to crossing probability estimates obtained for the height function and flat and sloped boundary conditions alike. As a result, arguments for estimating crossing proobabiliteis under the Ashkin-Teller measure rely upon determining whether there are connecteted components, between even or odd faces, for which the height function has a continuous segment of $+$, or of $-$, faces. Determining the regions of square lattice for which such a property holds is essential for obtaining weakened analogs of crossing probability estimates for the height function of the six-vertex model, due to the fact that the FKG inequality for the Ashkin-Teller along the self-dual line only holds for collections of even, or for odd, faces.

\subsection{Differences in encoding boundary conditions for the Ashkin-Teller model, from encoding flat and sloped boundary conditions for the six-vertex model}

In opposition to the six-vertex model which can be studied under flat and sloped boundary conditions, boundary conditions for the Ashkin-Teller model are encoded by collections of spins that one fixes along the faces of the height function lying in collections and even, and odd, faces of the square lattice. To obtain an analog of crossing probability estimates for the six-vertex height function under flat and sloped boundary conditions alike in the strip for the Ashkin-Teller model, one can consider the same environment of strips of the square lattice over which crossing probabilities are estimated, with the exception that the configurations drawn from the Ashkin-Teller sample space solely occupy collections of even, or of odd, faces of $\textbf{Z}^2$.

\begin{figure}
\begin{align*}
\includegraphics[width=0.53\columnwidth]{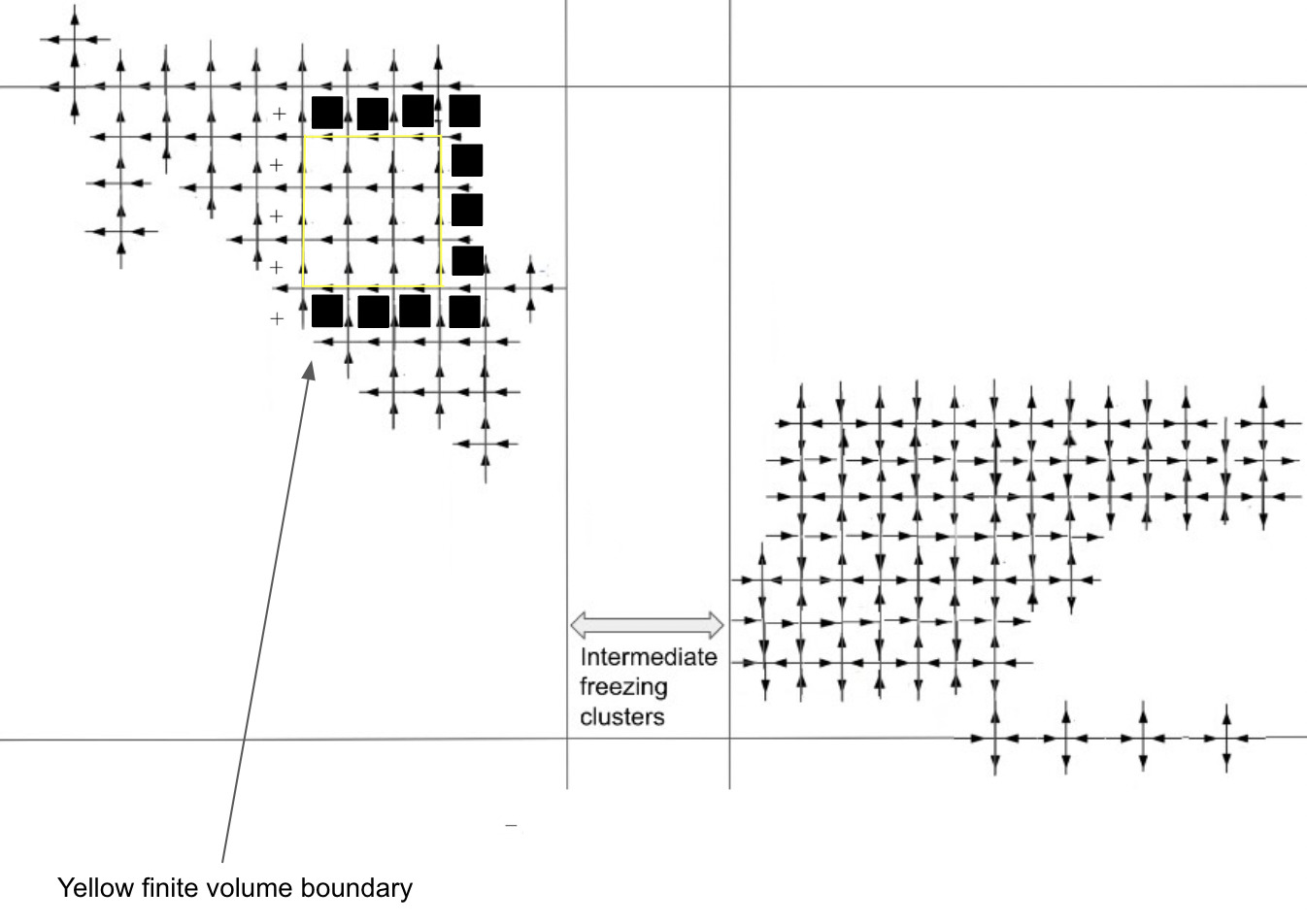}\\
\includegraphics[width=0.53\columnwidth]{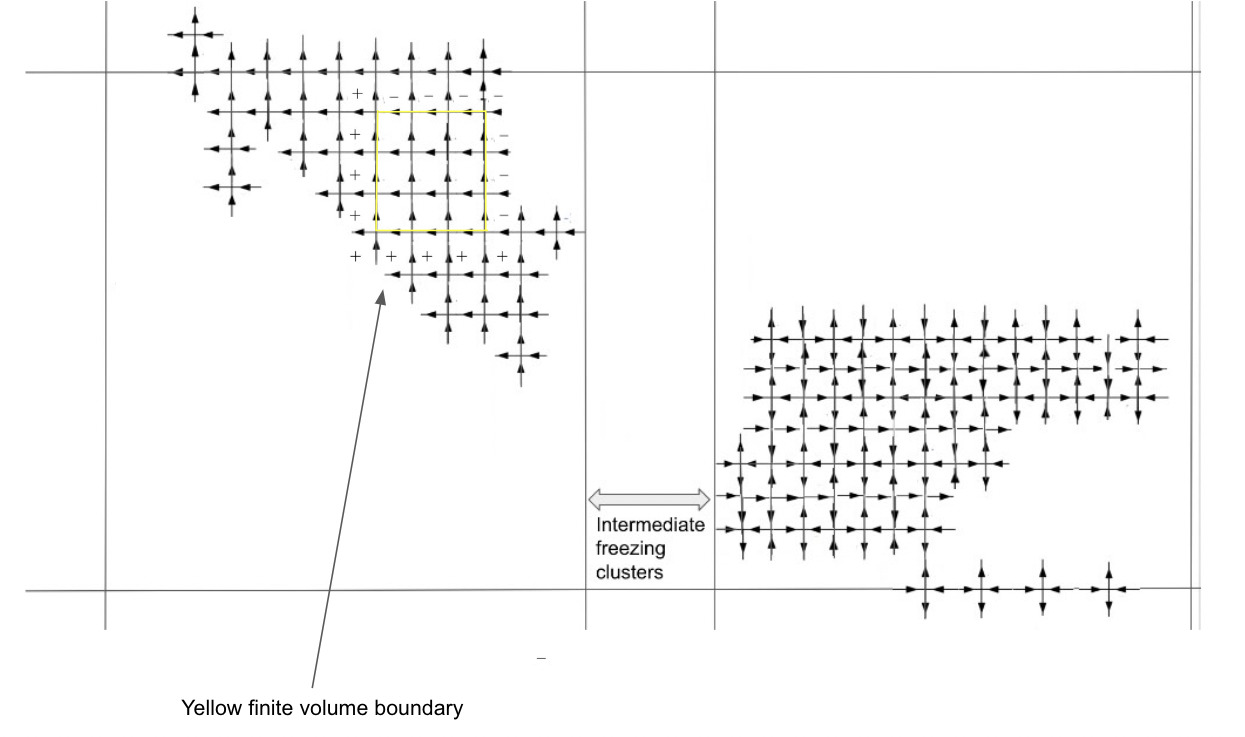}\\\includegraphics[width=0.53\columnwidth]{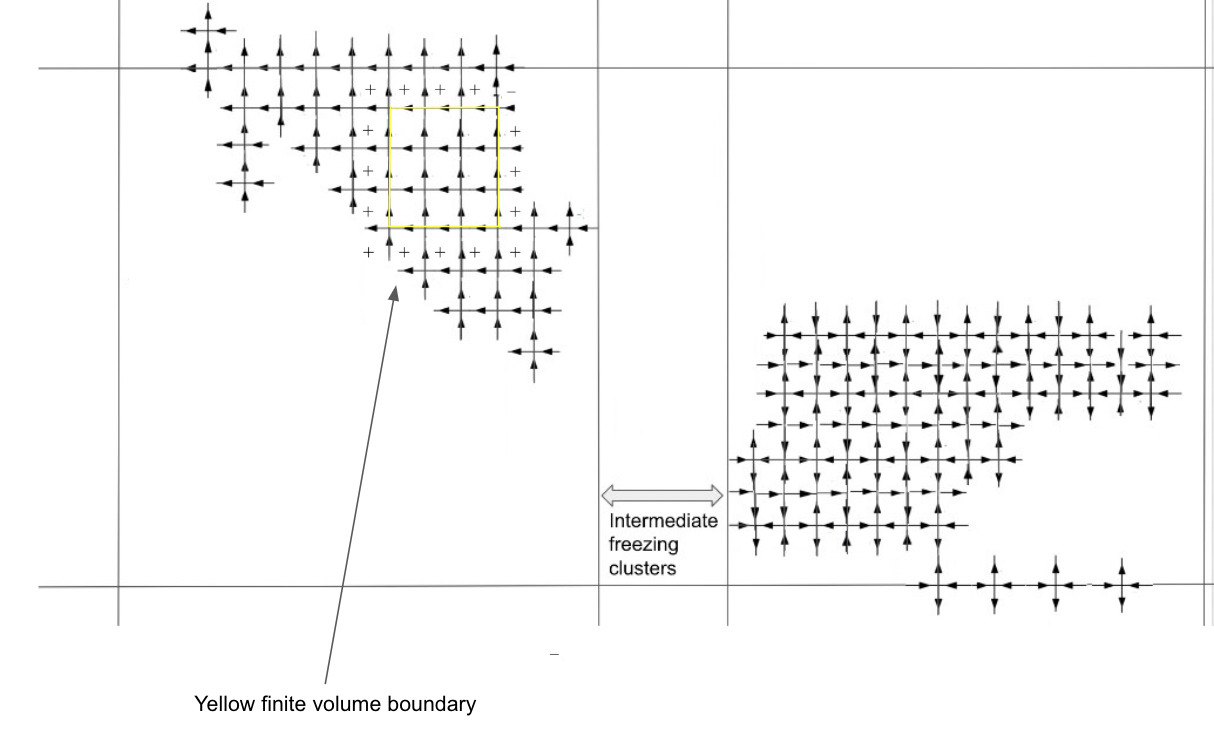}\\
\end{align*}
\caption{\textit{Ashkin-Teller spin representation of configurations from the six-vertex model in the strip}. A six-vertex configuration with \textit{freezing clusters} intersecting the strip is reproduced from \textit{Figure 4}, in which connectivity events between \textit{freezing clusters} and lines embedded in the strip were quantified. To pass to the spin-representation of the six-vertex model, $+$ and $-$ face variables are assigned to the faces enclosed in finite volume. In the same six-vertex configuration with frozen faces from \textit{Figure 4}, the corresponding Ashkin-Teller configuration is depicted, which consists of the following characteristics. Within the \textit{freezing cluster}, the boundary of a finite subvolume is depicted in yellow. To obtain the distribution of $+$ and $-$ variables within the interior, first, the distribution of $+$ face variables are assigned to all faces that do not have any overlap with any of the other surrounding incident faces to the yellow finite volume which are highlighted in black. Second, proceed to distribute more $+$ and $-$ spins on the remaining faces of $\textbf{Z}^2$ that are incident to the finite volume, in which faces in the uppermost right corner of each plaquette, a collection of four faces around each vertex, surrounding each frozen vertex are $-$.}
\end{figure}

\subsection{Mixed Ashkin-Teller objects}

\noindent In the following, denote some finite volume $\Lambda \subset ( \textbf{Z}^2 )^{*} \equiv (1/2 , 1/2) + \textbf{Z}$. Moving forwards, all quantities with a $*$ superscript live on odd faces of $\textbf{Z}^2$, and otherwise, on even faces of $\textbf{Z}^2$. In the following, fix two increasing functions $f(\sigma) \equiv f( + , - $, and $g(\sigma^{*}) \equiv g  \equiv g( + , - )$, with $f \equiv f(\sigma^{*}) $, or $f \equiv f(\sigma)$, and $g \equiv g(\sigma^{*})$, or $g \equiv g(\sigma^{*})$. From this dependency of the increasing functions only on the even or odd faces of the lattice, given by the respective distributions of spins from $\sigma$ and $\sigma^{*}$, respectively denoted with $F^{\mathrm{odd}} ( \text{ } \textbf{Z}^2 \text{ } ) \subset F( \text{ } \textbf{Z}^2 \text{ } )$, and with $F^{\mathrm{even}} ( \text{ } \textbf{Z}^2 \text{ } ) \subset F( \text{ } \textbf{Z}^2 \text{ } )$.

\bigskip

\noindent \textbf{Definition} \textit{13} (\textit{Ashkin-Teller Hamiltonian}). Define,

\begin{align*}
  \mathcal{H} \big( i , j , \tau(i ) , \tau(j) , \tau^{\prime}(i) , \tau^{\prime}(j) , J , U \big) \equiv \mathcal{H} \equiv  \underset{i \sim j}{\sum}  \bigg[   J \big( \text{  }     \tau(i) \tau(j) + \tau^{\prime} ( i ) \tau^{\prime} (j)   \text{ } \big) \text{ } + U \big( \text{ }  \tau(i) \tau(j) \tau^{\prime} ( i ) \tau^{\prime} (j)      \text{ } \big)    \bigg]          \text{ } \text{ , } 
\end{align*}

\noindent as the Ashkin-Teller Hamiltonian $\mathcal{H}$, consisting of one four-point interaction term multiplied by one Ising model coupling constant $U$, and two two-point interaction terms multipled by the other Ising model coupling constant $J$. Additionally, given real coupling constants $J< U$ of two coupled Ising models, and,

\begin{align*}
    \big( \tau , \tau^{\prime} \big) \text{ } \in \text{ }  \{ - 1 , 1 \}^{V(\Lambda)} \times \{ -1 , 1 \}^{V(\Lambda)}          \text{ } \text{ . } 
\end{align*}

\noindent Next, also introduce the probability measure for the Ashkin-Teller model, informally mentioned in the Introduction from the previous section.

\bigskip

\noindent \textbf{Definition} \textit{14} (\textit{Ashkin-Teller probability measure}). Define the probability of sampling some spin-representation configuration $\omega$ of the six-vertex model, where $ \omega \sim \mathcal{P}^{\xi}_{\Lambda} \big[ \cdot \big]$, as,

\begin{align*}
    \mathcal{P}^{\xi}_{\Lambda} [ \omega ] \equiv \frac{\mathrm{exp} \big[  - \mathcal{H} \big] }{Z^{\mathrm{AT},\xi} \big( i , j , \tau(i) , \tau(j) , \tau^{\prime}(i) , \tau^{\prime}(j) ,  \Lambda         \big) }   \equiv \frac{\mathrm{exp} \big[  - \mathcal{H} \big] }{Z^{\mathrm{AT}} \big( \Lambda         \big) }   \text{ } \text{ , } 
\end{align*}

\noindent the Ashkin-Teller probability measure, under boundary conditions $\xi \in \{ ++ , -- , + / -   \} \equiv \{ ++, -- , \mathrm{Mixed} \}$, with $Z^{\mathrm{AT}}$ the Ashkin-Teller partition function and configuration $\cdot \in \Omega^{\mathrm{AT}}$, uniformly drawn from the Ashkin-Teller sample space, that is given by the summation over all $i \sim j$ in finite volume,

\begin{align*}
  Z^{\mathrm{AT},\xi} \big(  i , j , \tau(i) , \tau(j) , \tau^{\prime}(i) , \tau^{\prime}(j) , \Lambda     \big) \equiv  Z^{\mathrm{AT}} \big( \Lambda         \big) \equiv  \underset{i,j \in \Lambda}{\underset{i \sim j}{\sum}}  \bigg[     J \big( \text{  }     \tau(i) \tau(j) + \tau^{\prime} ( i ) \tau^{\prime} (j)   \text{ } \big) \text{ } + U \big( \text{ }  \tau(i) \tau(j) \tau^{\prime} ( i ) \tau^{\prime} (j)      \text{ } \big)     \bigg]      \text{ } \text{ . } 
\end{align*}

\bigskip

\noindent \textbf{Proposition} \textit{AT-1} (\textit{existence of self-dual parameters from the Ashkin-Teller model establishing correspondence with the six-vertex model along the self-dual curve for six-vertex weights $a\equiv b\equiv 1$, and $c \equiv \mathrm{coth}(2J)$}, ({\color{blue}[16]}, \textbf{Lemma} \textit{7.1})). In the following relation, self-dual parameters of the Ashkin-Teller model, first due to Mittag and Stephen, satisfy,

\begin{align*}
       \big(  \tau , \tau^{\prime} ,       \xi      ,     \sigma^{*}     ,    \sigma           \big) \text{ }  \propto \text{ }     2^{-k( \xi^{*}) } \text{ } \big( c - 1  \big)^{|\xi| - |\omega(\sigma^{\prime})|} \text{ }        \textbf{1}_{\tau \tau^{\prime} \equiv \sigma^{*}} \text{ } \textbf{1}_{\tau \bot \xi^{*}} \text{ } \textbf{1}_{\sigma^{*} \bot \xi^{*}} \text{ } \textbf{1}_{\sigma \bot \xi}                                             \text{ } \text{ , } 
\end{align*}

\noindent with,

\begin{align*}
   (  \tau , \tau^{\prime} ) \in \text{ }  \{ - 1 , 1 \}^{V(\Lambda)} \times \{ -1 , 1 \}^{V(\Lambda)}     \text{ } \text{ . } 
\end{align*}

\noindent The above quantities represent two random variables sampled from the Hamiltonian used to define the Ashkin-Teller measure, distributed under $\mathcal{P}^{\xi} [\cdot]$,

\begin{align*}
       \xi    \sim \textbf{B}\textbf{C}^{\mathrm{AT}}_{++} \equiv  \textbf{B}\textbf{C}_{++ } \text{ } \text{ , } 
\end{align*}

\noindent or under $\textbf{BC}_{--}$, each of which represent another random variable distributed under $\mathrm{FK-IS}$, the law of an FK-Ising representation, from the set of all possible mixed $+$ and $-$ boundary conditions, distributed amongst odd and even faces of the square lattice, $\textbf{B}\textbf{C}^{\mathrm{AT}}_{\mathrm{Mixed}} \equiv\textbf{B}\textbf{C}_{\mathrm{Mixed}} \equiv \textbf{B}\textbf{C}_{+-} \equiv \textbf{B}\textbf{C}_{-+}$ for $\mathcal{P}[ \cdot ]$, and, 

\begin{align*}
    ( \sigma^{*} , \sigma) \text{ } \text{ . } 
\end{align*}

\noindent which is distributed under the spin-representation of the six-vertex model, under $\textbf{P}^{\xi} [ \cdot ]$. Finally, 

\begin{align*}
          \tau \bot \xi^{*} \equiv \big\{ \forall \text{ clusters on } \xi^{*}\text{ } ,  \text{ } \exists \text{ }  c \text{ } \in \textbf{R}   : \tau \equiv  \text{ } c   \big\}   \text{ } \text{ , } 
\end{align*}

\noindent appearing in the second indicator function indicates that the value of $\tau^{\prime}$ is constant on every cluster of $\xi^{\prime}$, with identical conditions holding for the remaining indicator functions, $\textbf{1}_{\sigma^{*} \bot \xi^{*}}$ and $\textbf{1}_{\sigma \bot \xi}$, appearing in the product from the above correspondence.

\bigskip

\noindent With the \textbf{Proposition} above, we proceed to introduce similarly oriented objects to those defined for the six-vertex model when executing RSW arguments in the strip in \textit{2.2}. Furthermore, introduce an index set for the total number of faces, $\mathcal{I}_{\mathrm{AT}}$, where the faces have $+$ or $-$ distribution, as described from configurations belonging to $\Omega^{\mathrm{AT}}$ which are pushed forwards under $\mathcal{P}^{\xi}_{\Lambda}[\cdot]$. In the case of the strip which was the first setting over which RSW arguments were provided for the six-vertex model, set $\Lambda \equiv [-m,m] \times [0,n^{\prime} N]$.

\bigskip

\noindent \textbf{Definition} \textit{15} (\textit{Ashkin-Teller equivalent of freezing clusters in the six-vertex model}). From \textbf{Definition} \textit{15}, introduce,

\begin{align*}
    \mathscr{F}\mathscr{C}^{\mathrm{even}}_{\mathrm{AT}} \equiv \mathscr{F} \mathscr{C}^{\mathrm{even}} \equiv  \underset{\{j_1 , \cdots , j_k     \} \in \mathcal{I}_{\mathrm{AT}}}{\underset{j_1 < \cdots < j_k  }{\bigcup}} \text{ } \big\{        \forall \text{ } + \backslash - \in F^{\mathrm{even}} ( \textbf{Z}^2 \cap \Lambda       ) \text{ } ,\text{ }  \exists \text{ }                 j_1 <         \cdots  <   j_k \neq i \in \mathcal{I}_{\mathrm{AT}}        \text{ }  :  \text{ }             |      + \backslash -         |    \leq k                       \big\}     \text{ } \text{ , } 
\end{align*}

\noindent for the Ashkin-Teller \textit{freezing cluster}, where the subsets of faces from $F( \textbf{Z}^2 \cap \Lambda)$, $+ \backslash -$, denote, a connected components of faces entirely colored with $+$ or $-$ face variables. Instead, if the connectivity event occurs over odd faces of the square lattice,

\begin{align*}
    \mathscr{F}\mathscr{C}^{\mathrm{odd}}_{\mathrm{AT}} \equiv \mathscr{F} \mathscr{C}^{\mathrm{odd}} \equiv  \underset{\{ j_1 , \cdots , j_k \} \in  \mathcal{I}_{\mathrm{AT}}}{\underset{j_1 < \cdots < j_k }{\bigcup}} \text{ } \big\{      \forall \text{ } + \backslash - \in F^{\mathrm{odd}} ( \textbf{Z}^2 \cap \Lambda       ) \text{ } ,\text{ }  \exists \text{ }                 j_1 <         \cdots  <   j_k \neq i \in \mathcal{I}_{\mathrm{AT}}        \text{ }  :  \text{ }             |       + \backslash -        |    \leq k                                       \big\}     \text{ } \text{ , } 
\end{align*}

\bigskip

\noindent Figures in the section will provide examples of such configurations that can be sampled from the Ashkin-Teller measure with positive probability. The set of all such clusters is denoted with $\mathcal{F}\mathcal{C}_{\mathrm{AT}}$. Over the even sublattice of the square lattice, $\mathcal{F}\mathcal{C}^{\mathrm{even}}_{\mathrm{AT}} \equiv \mathcal{F}\mathcal{C} \cap F^{\mathrm{even}} \big( \textbf{Z}^2 \big)$, while similarly, one also has that $\mathcal{F}\mathcal{C}^{\mathrm{odd}}_{\mathrm{AT}} \equiv \mathcal{F}\mathcal{C} \cap F^{\mathrm{odd}} \big( \textbf{Z}^2 \big)$, in order for the set of all Ashkin-Teller \textit{freezing clusters} to admit the decomposition, $\mathcal{F}\mathcal{C}_{\mathrm{AT}} \equiv \mathcal{F}\mathcal{C}^{\mathrm{even}}_{\mathrm{AT}} \cup \mathcal{F}\mathcal{C}^{\mathrm{odd}}_{\mathrm{AT}}$, over the even and odd sublattices of the square lattice that are dual to each other. The fact that the self-dual line Ashkin-Teller model does not satisfy the FKG inequality for the entire square lattice at the same time implies that the crossing events across strips of the square lattice must be constructed only from faces that are along a diagonal line from whichever face of the odd or even square of the square lattice that the height function previously occupied.  The restriction of the FKG inequality to the odd, or even, set of faces of the square lattice at a time also influences the random geometry from configurations in the generalized random-cluster, and $\big( q_{\sigma} , q_{\tau} \big)$-spin models, in that configurations associated with the probability measures of these two models can only reside over even, or odd, faces of the square lattice. In finite volume, the same restriction for encoding boundary conditions also exists, in which the spins along all odd and even faces incident to the boundary be set to either $+$ or $-$ uniformly. 

\bigskip

\noindent \textbf{Definition} \textit{15} and \textbf{Definition} \textit{16} below are direct translations of the requirements for the inner and outer diameters of six-vertex \textit{freezing clusters} introduced in \textbf{Definition} \textit{3}, and in \textbf{Definition} \textit{4}. 

\bigskip

\noindent \textbf{Definition} \textit{16} (\textit{Ashkin-Teller inner diameter of freezing clusters}). From \textbf{Definition} \textit{15}, introduce,

\begin{align*}
     \mathcal{D}^I_{\mathrm{AT}} \equiv  \mathcal{D}^I ( D, I )   \equiv   \mathcal{D}^I   \equiv   \underset{\mathscr{F}}{\mathrm{max}}  \text{ }   \bigg\{ \forall \mathscr{F} \text{ }  \exists \text{ }  \{   \mathscr{F}_j \}_{1 \leq j \leq D^{\prime}}  \neq \mathscr{F}        \text{ }   : \mathcal{P}^{+}_{\Lambda} \big[     \mathscr{F}  \overset{+ \backslash -}{\underset{ \mathscr{F}\mathscr{C}\cap F^{\mathrm{even}}      }{\longleftrightarrow}}    \mathscr{F}_{D^{\prime}}         \big] > 0    \text{ } \bigg\}            \text{ } \text{ , } 
\end{align*}

\noindent for the inner diameter of an Ashkin-Teller \textit{freezing cluster}, with $D^{\prime} > 0$. Instead, if the connectivity event occurs across odd faces of the square lattice,

\begin{align*}
    \mathcal{D}^I_{\mathrm{AT}} \equiv  \mathcal{D}^I ( D, I )   \equiv   \mathcal{D}^I   \equiv   \underset{\mathscr{F}}{\mathrm{max}}  \text{ }   \bigg\{ \forall \mathscr{F} \text{ }  \exists \text{ }  \{   \mathscr{F}_j \}_{1 \leq j \leq D^{\prime}}  \neq \mathscr{F}        \text{ }   : \mathcal{P}^{+}_{\Lambda} \big[     \mathscr{F}  \overset{+ \backslash -}{\underset{ \mathscr{F}\mathscr{C}\cap F^{\mathrm{odd}}      }{\longleftrightarrow}}    \mathscr{F}_{D^{\prime}}     \big] > 0    \text{ } \bigg\}            \text{ }   \text{ . } 
\end{align*}

\noindent As denoted above the connectivity event between $\mathscr{F}$ and $\mathscr{F}_{D^{\prime}}$, the superscript $+\backslash -$ indicates connectivity between connected components which entirely consist of $+$, or of $-$, faces.

\bigskip

\noindent The outer diameter of a \textit{freezing cluster} is similarly defined below.

\bigskip

\noindent \textbf{Definition} \textit{17} (\textit{Ashkin-Teller outer diameter of freezing clusters}). Along similar lines of \textbf{Definition} \textit{16}, introduce,

\begin{figure}
\begin{align*}
\includegraphics[width=0.35\columnwidth]{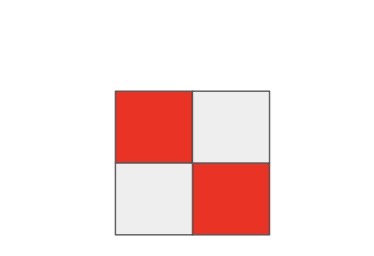}\\
\end{align*}
\caption{\textit{A depiction of even faces of the square lattice in red, and neighboring, dual odd faces of the square latice in grey}. By default, crossing events in the mixed Ashkin-Teller model for which the FKG lattice inequality is applicable pertain only to crossings between faces on the even square lattice, or between faces on the odd square lattice.}
\end{figure}

\begin{align*}
   \mathcal{D}^O_{\mathrm{AT}} \equiv  \mathcal{D}^O ( D-L, O )    \equiv \mathcal{D}^O \equiv    \underset{\mathscr{F}}{\mathrm{max}}  \text{ }   \bigg\{ \text{ } \forall \mathscr{F}\text{ }  \exists\text{ }    \{ \text{ } \mathscr{F}_k \text{ } \}_{1 \leq k \leq D^{\prime\prime}}     \neq \mathscr{F}             : \mathcal{P}^{+}_{\Lambda} \big[                  \mathscr{F}                 \overset{ + \backslash -}{\underset{ \mathscr{F}\mathscr{C}\cap F^{\mathrm{even}}  }{\longleftrightarrow}}   \mathscr{F}_{D^{\prime\prime}}     \big] > 0      \text{ } \bigg\}  \text{ } \text{ , } 
\end{align*}

\noindent for the outer diameter of an Ashkin-Teller \textit{freezing cluster}, with $D^{\prime\prime}>0$. Instead, if the connectivity event occurs across odd faces of the square lattice,

\begin{align*}
   \mathcal{D}^O_{\mathrm{AT}} \equiv  \mathcal{D}^O ( D-L, O )    \equiv \mathcal{D}^O \equiv    \underset{\mathscr{F}}{\mathrm{max}}  \text{ }   \bigg\{ \text{ } \forall \mathscr{F}\text{ }  \exists\text{ }    \{ \text{ } \mathscr{F}_k \text{ } \}_{1 \leq k \leq D^{\prime\prime}}     \neq \mathscr{F}             : \mathcal{P}^{+}_{\Lambda} \big[                 \mathscr{F}                 \overset{ + \backslash -}{\underset{\mathscr{F}\mathscr{C}\cap F^{\mathrm{odd}}  }{\longleftrightarrow}}   \mathscr{F}_{D^{\prime\prime}}    \big] > 0      \text{ } \bigg\}  \text{ } \text{ . } 
\end{align*}

\noindent As with the object introduced in \textbf{Definition} \textit{16}, the superscript $+\backslash -$ indicates connectivity between connected components which entirely consist of $+$, or of $-$, faces.

\bigskip

\noindent Due to the fact that there is no direct equivalent to sloped boundary conditions in the Ashkin-Teller model, we introduce $+$ boundary domains with the following.

\bigskip

\noindent \textbf{Definition} \textit{18} (\textit{Ashkin-Teller symmetric domains whose interior consists of connected components of $+$ or $-$ face variables that are invariant under global spin flips}). Introduce,

\[
F^{\mathrm{even}}( \textbf{Z}^2 )  \supsetneq  \Gamma^{\mathrm{even}} \equiv \Gamma^{\mathrm{even}}\big( \gamma_L , \gamma_R , [a_L , \widetilde{a_L}] , [a_R , \widetilde{a_R} ]  \big) \Longleftrightarrow  
\left\{\!\begin{array}{ll@{}>{{}}l} \gamma_L \subsetneq F^{\mathrm{even}}( \textbf{Z}^2  )  & \text{ : }  & \text{ } \mathcal{P}^{+}_{\Lambda} \big[      \mathcal{I}_L                   \underset{  \Lambda \cap F^{\mathrm{even}}}{\overset{+ \backslash -}{\longleftrightarrow }}                    \widetilde{\mathcal{I}_L}        \big] > 0   \text{ } \text{ , } \\
 \gamma_R \subsetneq F^{\mathrm{even}}(  \textbf{Z}^2 )  & \text{ : }  &    \mathcal{P}^{+}_{\Lambda} \big[        \mathcal{I}_R                        \underset{  \Lambda \cap F^{\mathrm{even}}}{\overset{+ \backslash -}{\longleftrightarrow }}                    \widetilde{\mathcal{I}_R}       \big] > 0        \text{ } \text{ , } \\
\end{array}\right.
\]

\noindent where, in the statement above, over the strip the Ashkin-Teller \textit{symmetric domain} is the union of four boundaries, $\mathcal{I}_L , \widetilde{\mathcal{I}_L}$, $\mathcal{I}_R$, $\widetilde{\mathcal{I}_R} \subset \Lambda$, for which there exists positive probability of a \textit{face contour} forming between $\mathcal{I}_L$ and $\widetilde{\mathcal{I}_L}$, and also between $\mathcal{I}_R$ and $\widetilde{\mathcal{I}_R}$, for $a_L \in \mathcal{I}_L$, $\widetilde{a_L} \in \widetilde{\mathcal{I}_L}$, $a_R \in \mathcal{I}_R$ and $\widetilde{a_R} \in \widetilde{\mathcal{I}_R}$. Instead, if the connectivity event occurs for odd faces of the square lattice, by making use of the same quantities described previously, the Ashkin-Teller domain the strip would take the form,

\[
F^{\mathrm{odd}}( \textbf{Z}^2  )  \supsetneq  \Gamma^{\mathrm{odd}} \equiv \Gamma^{\mathrm{odd}}\big( \gamma_L , \gamma_R , [a_L , \widetilde{a_L}] , [a_R , \widetilde{a_R} ] \big) \Longleftrightarrow 
\left\{\!\begin{array}{ll@{}>{{}}l} \gamma_L \subsetneq F^{\mathrm{odd}}(  \textbf{Z}^2  )  & \text{ : }  & \text{ } \mathcal{P}^{+}_{\Lambda} \big[       \mathcal{I}_L                   \underset{  \Lambda \cap F^{\mathrm{odd}}}{\overset{+ \backslash -}{\longleftrightarrow }}                    \widetilde{\mathcal{I}_L}        \big] > 0   \text{ } \text{ , } \\
 \gamma_R \subsetneq F^{\mathrm{odd}}(  \textbf{Z}^2  )  & \text{ : }  &    \mathcal{P}^{+}_{\Lambda} \big[        \mathcal{I}_R                        \underset{  \Lambda \cap F^{\mathrm{odd}}}{\overset{+ \backslash -}{\longleftrightarrow }}                    \widetilde{\mathcal{I}_R}      \big] > 0        \text{ } \text{ . } \\
\end{array}\right.
\]

\subsection{Properties of the probability measure}

     \begin{figure}
\begin{align*}
\includegraphics[width=0.73\columnwidth]{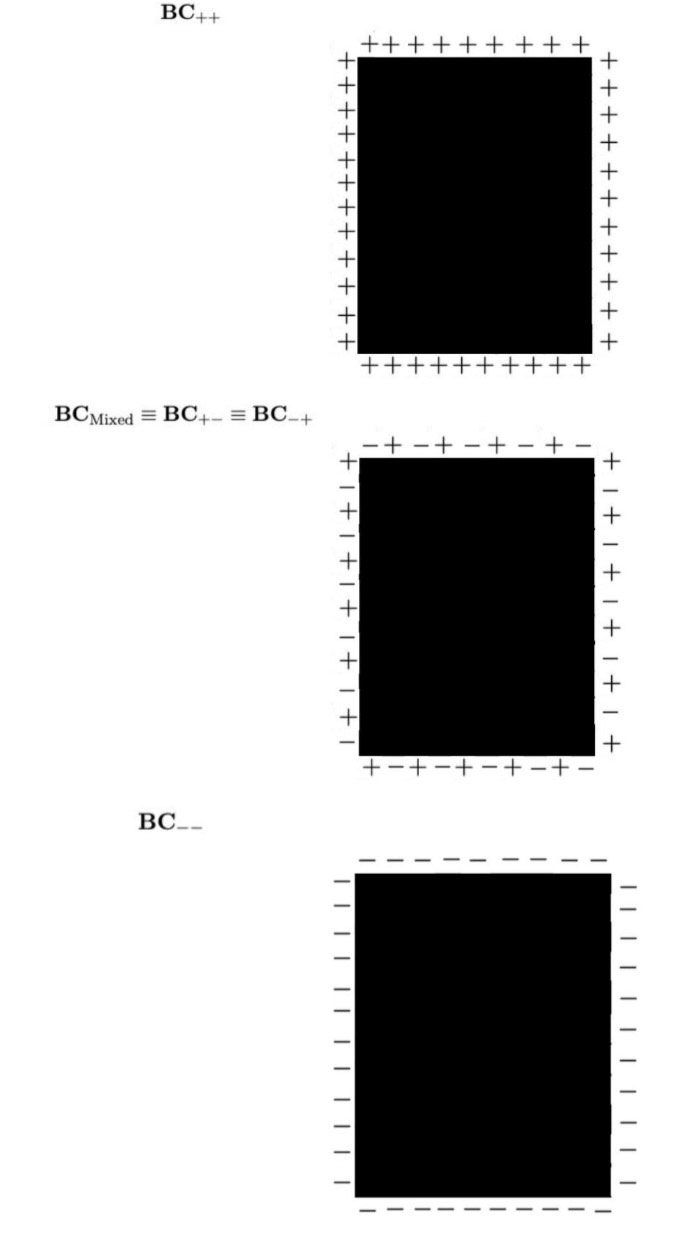}
\end{align*}
\caption{\textit{A depiction of the three possible classes of boundary conditions for the Ashkin-Teller model}.}
\end{figure}

\noindent We state the range of properties that are common to both the six-vertex and Ashkin-Teller probability measures below. In the following, fix the same quantities as defined in the beginning of \textit{6.2}, namely the increasing functions $f\equiv f( + , - )$, and $g  \equiv g( + , - )$, with $f \equiv f(\sigma^{*}) $, or $f \equiv f(\sigma)$, and $g \equiv g(\sigma^{*})$, or $g \equiv g(\sigma^{*})$. Being able to define these functions as dependent on $+$, or $-$, implies that a strategy along similar lines for the Ashkin-Teller model in the previous section can be applied. However, due to the lack of analogy for $|h|$, and for the sloped of boundary conditions in assigning $+$ or $-$ variables for the faces of the Ashkin-Teller model, a strategy for estimating arbitrarily long horizontal crossing probabilities relies upon quantifying the probability of obtaining crossings across connected components of $+$, or of $-$, faces. Given the random geometry of an Ashkin-Teller configuration outside of such $+$ and $-$ connected components of faces on $\textbf{Z}^2$, crossing probabilities can also be quantified depending on connectivity events between paths and analogues of unfrozen faces in the six-vertex model. From this dependency of the increasing functions only on the even or odd faces of the lattice, given by the respective distributions of spins from $\sigma$ and $\sigma^{*}$, respectively denoted with $F^{\mathrm{odd}} ( \text{ } \textbf{Z}^2 \text{ } ) \subset F( \text{ } \textbf{Z}^2 \text{ } )$, and $F^{\mathrm{even}} ( \text{ } \textbf{Z}^2 \text{ } ) \subset F( \text{ } \textbf{Z}^2 \text{ } )$, the properties below,

\begin{align*}
  \forall \text{ increasing functions } f, g     \text{ }  \Rightarrow \text{ }  \mathcal{P}^{\xi}_{\Lambda} [ f  g ] \geq \mathcal{P}^{\xi}_{\Lambda} [ f ] \text{ }   \mathcal{P}^{\xi}_{\Lambda} [ g ]    \text{ , }           \tag{$\mathrm{marginal \text{ } FKG}$} \\ 
  \forall \text{ } \xi^{\prime} \geq \xi \Rightarrow       \mathcal{P}^{\xi^{\prime}}_{\Lambda} [ f ]      \geq \mathcal{P}^{\xi}_{\Lambda} [ f ]     \text{ , }  \tag{$\mathrm{marginal \text{ } CBC}$} \text{ , } 
\end{align*}

\noindent hold for $\Lambda^{\prime} \subset \Lambda$. 
From the three properties above, next we introduce analogs of \textbf{Proposition} \textit{1.1} and of \textbf{Corollary} \textit{1.2}. Relatedly, for the Ashkin-Teller expectation $\mathcal{E}^{\xi}_{\Lambda} [ \cdot ]$, the following properties hold,

\begin{align*}
 \mathcal{E}^{\xi}_{\Lambda} [ f g] \text{ } \geq  \mathcal{E}^{\xi}_{\Lambda} [ f ] \text{ } \mathcal{E}^{\xi}_{\Lambda} [ g ] \text{ }  \tag{$\mathrm{marginal \text{ }FKG}$}\\
 \forall \xi^{\prime} \geq \xi \Rightarrow     \text{ }  \mathcal{E}^{\xi^{\prime}}_{\Lambda} [ f ] \geq  \mathcal{E}^{\xi}_{\Lambda} [ f  ]  \text{ } \tag{$\mathrm{marginal \text{ } CBC}$} \text{ , } 
\end{align*}

\noindent for increasing functions $X,Y$, each of which can be a function of a finite volume that entirely consist of $+$ faces, or of $-$ faces, as the boundary conditions for $F^{\mathrm{even}}$, or for $F^{\mathrm{odd}}$. The (marginal CBC) property is established in the Appendix following the conclusion of arguments for quantifying crossing probabilities in the strip for the Ashkin-Teller model, which makes use of observations for the proof of the (FKG) and (CBC) inequalities from {\color{blue}[11]}, in addition to \textbf{Theorem} \textit{4} from {\color{blue}[16]} which provides conditions on the (FKG) inequality which only holds for \textit{marginals} of mixed-spin configurations.

\bigskip

  \begin{figure}
\begin{align*}
\includegraphics[width=0.54\columnwidth]{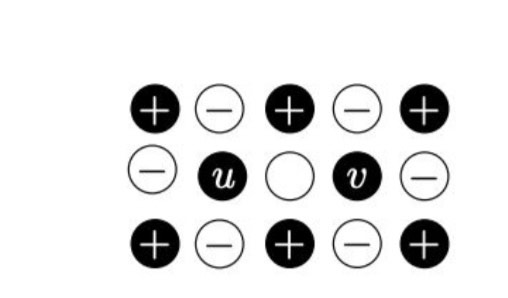}
\end{align*}
\caption{\textit{An adaptation to the spin configuration depicted in Figure 4 of {\color{blue}[16]}, with $-$ boundary conditions on odd faces of the square lattice instead of a combination of $+$ and $-$ boundary spins on odd faces, for a mixed Ashkin-Teller spin configuration so that the positive association inequality is always satisfied}.}
\end{figure}

\noindent For results that remain later in the section, we introduce the following.

\bigskip

\noindent \textbf{Definition} \textit{19} (\textit{Ashkin-Teller crossing events in the strip for combining SMP and CBC properties}). Define, for $\mathscr{F}_i,\mathscr{F}_j \in F^{\mathrm{even}}\big( \textbf{Z}^2 \big)$,

\begin{align*}
 E^{\mathrm{even}}_{\mathrm{AT}} \equiv   E_{\mathrm{AT}} \equiv E \equiv  \big\{     \mathscr{F}_i      \underset{F^{\mathrm{even}}}{ \overset{+ \backslash -}{\longleftrightarrow}}    \mathscr{F}_j   \big\}  \equiv    \big\{       \mathscr{F}_i      { \overset{+ \backslash -}{\longleftrightarrow}}    \mathscr{F}_j    \big\}    \text{ } \text{ , } 3
 \end{align*}

\noindent as the crossing event across even faces, for some $i \neq j$, in the strip. Instead, for $\mathscr{F}_i , \mathscr{F}_j \in F^{\mathrm{odd}} ( \text{ }      \textbf{Z}^2 \text{ } ) $, namely if the connectivity event is taken over odd faces, the crossing event takes the form,

\begin{align*}
 E^{\mathrm{odd}}_{\mathrm{AT}} \equiv   E_{\mathrm{AT}} \equiv E \equiv  \big\{      \mathscr{F}_i      \underset{F^{\mathrm{odd}}}{ \overset{+ \backslash -}{\longleftrightarrow}}    \mathscr{F}_j  \big\}   \equiv    \big\{   \mathscr{F}_i      { \overset{+ \backslash -}{\longleftrightarrow}}    \mathscr{F}_j   \big\}      \text{ } \text{ , } 
 \end{align*}

\bigskip

\noindent \textit{Remark}. In comparison to the theorems combining (SMP) and (CBC) properties for the six-vertex probability measure $\textbf{P}^{\xi^{\mathrm{Sloped}}}_{\Lambda}[\cdot]$, from the marginal properties of (FKG) and (CBC) that are satisfied by the mixed Ashkin-Teller spin measure $\mathcal{P}^{\xi}_{\Lambda} [ \cdot]$, there is no difference between values of the level of the height function, which appears in the height function crossing events $\mathcal{E}_{k+j}$ and $\mathcal{E}_k$, in the lower and upper bounds, respectively, of \textbf{Corollary} \textit{1.2}.

\bigskip

\noindent With the aforementioned differences between crossings of the six-vertex and Mixed Ashkin-Teller spin representation, we turn towards quantifying the pushforward of the strip connectivity event that was introduced in \textbf{Proposition} \textit{1.2}, under $\mathcal{P}^{\xi}_{\Lambda} [ \cdot]$, with the crossing probability,

\begin{align*}
\mathcal{P}^{\xi}_{\Lambda} \big[     [ 0, \lfloor \delta^{\prime\prime} n \rfloor ] \times \{ 0 \} \underset{F^{\mathrm{even}}}{\longleftrightarrow}  [i , i+ \lfloor \delta^{\prime\prime}n \rfloor ]  \times \{n \}  \big]     \text{ } \text{ , } 
\end{align*}

\noindent over even faces, in addition to,

\begin{align*}
\mathcal{P}^{\xi}_{\Lambda} \big[     [ 0, \lfloor \delta^{\prime\prime} n \rfloor ] \times \{ 0 \} \underset{F^{\mathrm{odd}}}{\longleftrightarrow}  [i , i+ \lfloor \delta^{\prime\prime}n \rfloor ]  \times \{n \}  \big]     \text{ } \text{ , }  \tag{Event 1}
\end{align*}

\noindent corresponding to the crossing probability over odd faces, with the following result.

\bigskip

\noindent \textbf{Proposition} \textit{AT 1} (\textit{Mixed Ashkin-Teller analog of upper bound for connectivity event occurring in the strip}). WLOG suppose that the connectivity event below occurs over $F^{\mathrm{even}}$. Introduce the same choice of parameters mentioned in \textbf{Proposition} \textit{1.2}, including, $\delta , \delta^{\prime\prime}, n >0$ such that $\lfloor{}\delta n\rfloor{}> \lfloor{}\delta^{\prime\prime} n\rfloor{} \in \textbf{Z}$. For any $k \geq \frac{1}{c}$ and $i  \in \textbf{Z}$, the crossing probability,

\begin{align*}
 \mathcal{P}^{\xi}_{\Lambda} \big[     [ 0, \lfloor \delta^{\prime\prime} n \rfloor ] \times \{ 0 \} \underset{F^{\mathrm{even}}}{\longleftrightarrow}  [i , i+ \lfloor \delta^{\prime\prime}n \rfloor   ]  \times \{n \}  \big]  \equiv \mathcal{P}^{\xi}_{\Lambda} \big[     [ 0, \lfloor \delta^{\prime\prime} n \rfloor ] \times \{ 0 \} {\longleftrightarrow}  [i , i+ \delta^{\prime\prime}n  ]  \times \{n \}  \big]  \leq 1 - c \text{ } \text{ , } 
\end{align*}

\noindent admits an upper bound dependent upon $c$.

\bigskip

\noindent \textit{Proof of Proposition AT 1}. To prove the upper bound provided above that is dependent upon $c$, we provide arguments for the following \textbf{Lemma}.

\bigskip

\noindent \textbf{Lemma} \textit{AT 2.1}. (\textit{upper bound for the Ashkin-Teller segment connectivity event between $\mathcal{I}_0$ and $\widetilde{\mathcal{I}_0}$}). WLOG suppose that the segment connectivity event is dependent upon faces of the even sublattice, namely that,

\begin{align*}
      \bigg\{   \mathcal{I}_0 \underset{F^{\mathrm{even}}}{\overset{+ \backslash -}{\longleftrightarrow} } \widetilde{\mathcal{I}_0}      \bigg\}  \equiv \big\{       \mathcal{I}_0 {\overset{+ \backslash -}{\longleftrightarrow} } \widetilde{\mathcal{I}_0}     \big\}  \equiv \big\{ \mathcal{I}_0     \longleftrightarrow \widetilde{\mathcal{I}_0}  \big\}      \text{ } \text{ . }  \tag{Event 2}
\end{align*}

\noindent Fix $\chi \in \textbf{B}\textbf{C}_{++}$. With regards to the segment connectivity event between two segments $\mathcal{I}_0$ and $\widetilde{\mathcal{I}_0}$, which are respectively on the top and bottom boundaries of the finite volume strip,

\begin{align*}
       \mathcal{P}^{\chi}_{\Lambda} \big[      \mathcal{I}_0 {\longleftrightarrow}  \widetilde{\mathcal{I}_0}       \big]  \leq 1 -        \widetilde{\mathscr{R}^{\prime}}     \mathcal{P}^{\chi}_{\Lambda} \bigg[  \big\{ \mathcal{I}_{-\mathscr{I}} {\longleftrightarrow} \widetilde{\mathcal{I}_{-\mathscr{I}}} \big\}  \cap    \big\{   \mathcal{I}_{-\mathscr{J}} {\longleftrightarrow}     \widetilde{\mathcal{I}_{-\mathscr{J}}}     \big\}   \cap       \big\{          \mathcal{I}_{-\mathscr{K}} {\longleftrightarrow}     \widetilde{\mathcal{I}_{-\mathscr{K}}}           \big\}            \bigg]                    \text{ } \text{ , } \tag{Event 3}
\end{align*}

\noindent for $0 \neq \mathscr{I} < \mathscr{J} < \mathscr{H}$ with a strictly positive multiplicative factor $\widetilde{\mathscr{R}^{\prime}}$ for the intersection of crossing events appearing in the upper bound.

\bigskip

\noindent \textit{Proof of Lemma AT 2.1}. First, observe, by (marginal FKG), that the upper bound intersection of crossing probabilities is bound below by the product,

\begin{align*}
     \mathcal{P}^{\chi}_{\Lambda} \big[  \mathcal{I}_{-\mathscr{I}} {\longleftrightarrow} \widetilde{\mathcal{I}_{-\mathscr{I}}} \big]  \mathcal{P}^{\chi}_{\Lambda} \big[  \mathcal{I}_{-\mathscr{J}} {\longleftrightarrow}     \widetilde{\mathcal{I}_{-\mathscr{J}}}   \big]   \mathcal{P}^{\chi}_{\Lambda} \big[              \mathcal{I}_{-\mathscr{K}} {\longleftrightarrow}     \widetilde{\mathcal{I}_{-\mathscr{K}}}                    \big]         \text{ }\text{ .  }
\end{align*}

\noindent Following the application of (marginal FKG) above, additionally observe,

\begin{align*}
     \underset{k \in \{      \mathscr{I} , \mathscr{J} , \mathscr{K} \} }{\prod}   \mathcal{P}^{\chi}_{\Lambda} \big[    \mathcal{I}_{-k} {\longleftrightarrow}    \widetilde{\mathcal{I}_{-k}}          \big]     \geq \big( \mathcal{P}^{\chi}_{\Lambda} \big[   \mathcal{I}_0         { \longleftrightarrow} \widetilde{\mathcal{I}_0}        \big] \big)^3              \text{ } \text{ , }
\end{align*}

\noindent from which, along the lines of arguments previously provided for obtaining the upper bound for the six-vertex model as given in \textbf{Lemma} \textit{2.5}, implies, for the upper bound to be strictly less than $1$, that,

\begin{align*}
      \mathcal{P}^{\chi}_{\Lambda} \big[      \mathcal{I}_0 {\longleftrightarrow}  \widetilde{\mathcal{I}_0}       \big]  \leq 1 \\ \Updownarrow \\   \mathcal{P}^{\chi}_{\Lambda} \big[      \mathcal{I}_0 {\longleftrightarrow}  \widetilde{\mathcal{I}_0}       \big]  \leq 1 -  \widetilde{\mathscr{R}^{\prime}}     \mathcal{P}^{\chi}_{\Lambda} \bigg[  \big\{ \mathcal{I}_{-\mathscr{I}} {\longleftrightarrow} \widetilde{\mathcal{I}_{-\mathscr{I}}} \big\}  \cap    \big\{   \mathcal{I}_{-\mathscr{J}} {\longleftrightarrow}     \widetilde{\mathcal{I}_{-\mathscr{J}}}     \big\}   \cap       \big\{          \mathcal{I}_{-\mathscr{K}} {\longleftrightarrow}     \widetilde{\mathcal{I}_{-\mathscr{K}}}           \big\}            \bigg]  \\    \leq  1 -  \widetilde{\mathscr{R}^{\prime}}     \mathcal{P}^{\chi}_{\Lambda} \bigg[  \big\{ \mathcal{I}_{-\mathscr{I}} {\longleftrightarrow} \widetilde{\mathcal{I}_{-\mathscr{I}}} \big\}  \cap    \big\{   \mathcal{I}_{-\mathscr{J}} {\longleftrightarrow}     \widetilde{\mathcal{I}_{-\mathscr{J}}}     \big\}   \cap       \big\{          \mathcal{I}_{-\mathscr{K}} {\longleftrightarrow}     \widetilde{\mathcal{I}_{-\mathscr{K}}}           \big\}   \bigg] \\ \equiv  1 -  \widetilde{\mathscr{R}^{\prime}}     \mathcal{P}^{\chi}_{\Lambda} \bigg[   \underset{k \in \{ \mathscr{I} , \mathscr{J} , \mathscr{K} \}}{\bigcap}     \big\{      \mathcal{I}_{-k} \longleftrightarrow   \widetilde{\mathcal{I}_{-k}}    \big\}    \bigg] \\ \overset{(\mathrm{marginal \text{  }FKG})}{\leq} 1  - \widetilde{\mathscr{R}^{\prime}}    \bigg[   \underset{k \in \{ \mathscr{I} , \mathscr{J} , \mathscr{K}    \}    }{\prod}      \mathcal{P}^{\chi}_{\Lambda} \big[   \mathcal{I}_{-k} \longleftrightarrow   \widetilde{\mathcal{I}_{-k}}  \big]     \bigg]   \\ \Updownarrow \\  1 -  \widetilde{\mathscr{R}^{\prime}}    \big(    \mathcal{P}^{\chi}_{\Lambda} \big[   \mathcal{I}_0         { \longleftrightarrow} \widetilde{\mathcal{I}_0}        \big]    \big)^3 < 1  \\ \Updownarrow  \\
      \widetilde{\mathscr{R}^{\prime}}    \big(    \mathcal{P}^{\chi}_{\Lambda} \big[   \mathcal{I}_0         { \longleftrightarrow} \widetilde{\mathcal{I}_0}        \big]    \big)^3  > 0 \\ \Updownarrow \\ \big(    \mathcal{P}^{\chi}_{\Lambda} \big[   \mathcal{I}_0         { \longleftrightarrow} \widetilde{\mathcal{I}_0}        \big]        \big)^{-1} > \sqrt[3]{\widetilde{\mathscr{R}^{\prime}}}  >  0        \text{ } \text{ , }
\end{align*}

\noindent exhibiting that the claim in the \textbf{Lemma} above holds for $i \neq 0$, given $\widetilde{\mathscr{R}^{\prime}}$ sufficiently small satisfying the conditions above. On the other hand, for the two remaining cases of the argument, introduce some parameter $\delta^{\prime\prime}$ satisfying previous arguments for the upper bound in the six-vertex model under flat and sloped boundary conditions, in which for $i$ not belonging to the first case of the argument, namely $i \in [ - 2 \lfloor \delta^{\prime\prime} n \rfloor , 2 \lfloor \delta^{\prime\prime} n \rfloor]$, observe,

\begin{align*}
    \mathcal{P}^{\chi}_{\Lambda} \big[             [ 0, \lfloor \delta^{\prime\prime} n \rfloor ] \times \{ 0 \} {\longleftrightarrow}  [i , i+ \lfloor \delta^{\prime\prime}n \rfloor   ]  \times \{n \}                     \big]   \leq        \mathcal{P}^{\chi}_{\Lambda} \big[             [ 0, \lfloor \delta n \rfloor ] \times \{ 0 \} {\longleftrightarrow}  [i , i+ \lfloor \delta  n \rfloor ]  \times \{n \}                     \big]   \leq 1 - c     \text{ } \text{ , }         \tag{Event 4}
\end{align*}

\noindent while, for the remaining case in which $|i| > 2 \delta^{\prime\prime}n$, observe,

\begin{align*}
      \mathcal{P}^{\chi}_{\Lambda} \big[               [ \lfloor \delta^{\prime\prime} n \rfloor, 2 \lfloor \delta^{\prime\prime} n \rfloor ] \times \{ 0 \} {\longleftrightarrow}  [ - i + \lfloor \delta^{\prime\prime} n \rfloor , i+ 2  \lfloor \delta^{\prime\prime}n \rfloor ]  \times \{n \}              \big] \\  \geq     \mathcal{P}^{\chi}_{\Lambda} \big[        [ 0 ,  \lfloor \delta^{\prime\prime} n \rfloor ] \times \{ 0 \} {\longleftrightarrow}  [  i + \lfloor \delta^{\prime\prime} n \rfloor , i+  2 \lfloor \delta^{\prime\prime}n \rfloor ]  \times \{n \}                   \big]       \text{ } \text{ . }        \tag{Event 5} 
\end{align*}

\noindent Altogether,

\begin{align*}
       \mathcal{P}^{\chi}_{\Lambda} \big[       [ 0 , \lfloor \delta^{\prime\prime} n \rfloor ] \times \{ 0 \} \longleftrightarrow   [ i , i + \lfloor \delta^{\prime\prime} n \rfloor ]          \times \{ n \}           \big] < \frac{1}{2}          \text{ } \text{ , } 
\end{align*}

\noindent for a suitably chosen constant $< \frac{1}{2}$, from which we conclude the argument. \boxed{}

\bigskip

\noindent With the results from the above \textbf{Lemma}, to complete the arguments for the proof of \textbf{Proposition} \textit{AT 1}, make use of the following series of \textbf{Lemmas} for completing the proof as in the case for flat, and for sloped, boundary conditions in the six-vertex model, in which the desired upper bound for the segment connectivity event provided in \textbf{Proposition} \textit{AT 1} is obtained from estimations of vertical and horizontal crossing probabilities across \textit{symmetric strip domains}. From properties of Mixed Ashkin-Teller \textit{symmetric domains} introduced in \textbf{Definition} \textit{18}, we implement steps of the argument providing lower bounds for (\textit{2.2 left boundary symmetric domain lower bound}).

\bigskip

\noindent \textbf{Lemma} \textit{6.3} (\textit{Mixed Ashkin-Teller analog of lower bound for connectivity event between the left symmetric domain and Ashkin-Teller freezing clusters}). WLOG suppose that the connectivity event below occurs over $\Lambda^{\mathrm{even}} \equiv \Lambda$, and $\mathscr{F}\mathscr{C}^{\mathrm{even}} \equiv \mathscr{F}\mathscr{C}$. Fix $\chi_1  \in \textbf{B}\textbf{C}_{++}$, $\chi_2 \in \textbf{BC}_{-+}$, and $\mathscr{F}\mathscr{C} \in \mathcal{F}\mathcal{C}$. Under the assumption that $\chi_1 \leq \chi_2$, the induced $\mathrm{L}$-$1$ absolute-value distance between the following two crossing probabilities,

\begin{align*}
\bigg|  \mathcal{P}^{\chi_2}_{\Lambda} \bigg[      \gamma_L \overset{+ \backslash -}{\underset{\Lambda^{\mathrm{even}} \cap (\mathscr{F}\mathscr{C}_1)^c}{\longleftrightarrow}}    \mathscr{F}\mathscr{C}_1          \bigg] -                \mathcal{P}^{\chi_1}_{\Lambda} \bigg[         \gamma_L \overset{+ \backslash -}{\underset{\Lambda^{\mathrm{even}}\cap (\mathscr{F}\mathscr{C}_1)^c}{\longleftrightarrow}}    \mathscr{F}\mathscr{C}_1            \bigg]      \bigg|  \equiv \bigg|  \mathcal{P}^{\chi_2}_{\Lambda} \bigg[      \gamma_L \overset{+ \backslash -}{\underset{\Lambda \cap (\mathscr{F}\mathscr{C}_1)^c}{\longleftrightarrow}}    \mathscr{F}\mathscr{C}_1          \bigg] -                \mathcal{P}^{\chi_1}_{\Lambda} \bigg[         \gamma_L \overset{+ \backslash -}{\underset{\Lambda\cap (\mathscr{F}\mathscr{C}_1)^c}{\longleftrightarrow}}    \mathscr{F}\mathscr{C}_1            \bigg]      \bigg| \text{ } \text{ , } \\  \tag{Event 6} 
\end{align*}

\noindent admits a suitably chosen, strictly positive lower bound, given by,

\begin{align*}
      C^{\chi_1 , \chi_2} \big( \Lambda \big) \equiv C^{\chi_1, \chi_2}              \text{ } \text{ , } 
\end{align*}

\noindent where the superscript $+ \backslash -$, as provided above the connectivity event between the left \textit{symmetric} domain boundary, and its intersection with the first mixed Ashkin-Teller \textit{freezing cluster}, indicates that there can be a path connecting $\gamma_L$ with $\mathscr{F}\mathscr{C}_1$ consisting of $+$ and $-$ face variables.

\bigskip

\noindent \textit{Proof of Lemma 6.3}. Fix a mixed-spin configuration $\cdot \in \Omega^{\mathrm{AT}}$. To avoid continually making use of burdensome notation, suppose that the following connectivity events occur over $F^{\mathrm{even}}(\textbf{Z}^2) \equiv F^{\mathrm{even}}$, in addition to setting $\Lambda^{\mathrm{even}} \equiv \Lambda$. The arguments for connectivity events over the odd faces of the square lattice follow by instead setting $F^{\mathrm{odd}}(\textbf{Z}^2) \equiv F^{\mathrm{even}}$, in addition to setting $\Lambda^{\mathrm{odd}} \equiv \Lambda$. In the following arguments, to emphasize the fact that \textit{Ashkin-Teller} crossings, in comparison to six-vertex crossings, can be simultaneously dependent upon $+$ and $-$ face variables, denote,

\begin{align*}
      \mathcal{P}^{\chi_1}_{\Lambda} \big[      \gamma_L \overset{+ \backslash -}{\underset{\Lambda \cap (\mathscr{F}\mathscr{C}_1)^c}{\longleftrightarrow}}    \mathscr{F}\mathscr{C}_1          \big]     \equiv      \mathcal{P}^{\chi_1}_{\Lambda} \big[      \gamma_L \overset{\mathrm{Mixed}}{\underset{\Lambda \cap (\mathscr{F}\mathscr{C}_1)^c}{\longleftrightarrow}}    \mathscr{F}\mathscr{C}_1          \big]      \text{ } \text{ , } 
\end{align*}

\noindent in addition to,

\begin{align*}
      \mathcal{P}^{\chi_2}_{\Lambda} \big[      \gamma_L \overset{+ \backslash -}{\underset{\Lambda \cap (\mathscr{F}\mathscr{C}_1)^c}{\longleftrightarrow}}    \mathscr{F}\mathscr{C}_1          \big]     \equiv      \mathcal{P}^{\chi_2}_{\Lambda} \big[      \gamma_L \overset{\mathrm{Mixed}}{\underset{\Lambda \cap (\mathscr{F}\mathscr{C}_1)^c}{\longleftrightarrow}}    \mathscr{F}\mathscr{C}_1          \big]      \text{ } \text{ . } 
\end{align*}

\noindent For boundary conditions in $\textbf{B}\textbf{C}_{++}$, as a pushforward under $\mathcal{P}^{\xi}_{\Lambda} [ \cdot]$, from the absolute value difference between connectivity event, respectively taken under $\chi_2 \in \textbf{B}\textbf{C}_{++}$ and $\chi_1 \in \textbf{B}\textbf{C}_{++}$, in the \textit{Ashkin-Teller freezing cluster},

\begin{align*}
  \bigg|  \mathcal{P}^{\chi_2}_{\Lambda} \big[      \gamma_L \overset{\mathrm{Mixed}}{\underset{\Lambda \cap (\mathscr{F}\mathscr{C}_1)^c}{\longleftrightarrow}}    \mathscr{F}\mathscr{C}_1          \big] -                \mathcal{P}^{\chi_1}_{\Lambda} \big[         \gamma_L \overset{\mathrm{Mixed}}{\underset{\Lambda\cap (\mathscr{F}\mathscr{C}_1)^c}{\longleftrightarrow}}    \mathscr{F}\mathscr{C}_1            \big]      \bigg| \text{ }  \text{ } \text{ } \text{ , } 
\end{align*}

\noindent express each probability with the following decomposition for the crossing probability first under $\chi_2$ with the linear combination,

\begin{align*}
  \mathcal{P}^{\chi_2}_{\Lambda} \big[     \gamma_L        \underset{\Lambda \cap (\mathscr{F}\mathscr{C})_1}{\overset{-}{\longleftrightarrow}}                 \gamma_{-}                   \big] \text{ }  +  \text{ }    \mathcal{P}^{\chi_2}_{\Lambda} \big[       \big(  \gamma_{-,+}   \big)_1        \underset{\Lambda \cap (\mathscr{F}\mathscr{C}_1 )^c}{\overset{\mathrm{Mixed}}{\longleftrightarrow}}                        \gamma_{L^{\prime}}             \big] \text{ } +   \text{ }    \mathcal{P}^{\chi_2}_{\Lambda} \big[     \big(   \gamma_{-,+}  \big)_2          \underset{\Lambda \cap (\mathscr{F}\mathscr{C})_1}{\overset{+}{\longleftrightarrow}}                     \mathscr{F}\mathscr{C}_1                \big] \text{ }        \text{ , }   \tag{Event 7} 
\end{align*}

\noindent as well as the following expression for the second probability under $\chi_2$ with the linear combination,

\begin{align*}
    \mathcal{P}^{\chi_1}_{\Lambda} \big[     \gamma_L          \underset{\Lambda \cap (\mathscr{F}\mathscr{C} )_1}{\overset{-}{\longleftrightarrow}}                 \gamma_{-}                   \big] \text{ }  +  \text{ }    \mathcal{P}^{\chi_1}_{\Lambda} \big[       \big(  \gamma_{-,+}   \big)_1        \underset{\Lambda \cap (\mathscr{F}\mathscr{C}_1)^c}{\overset{\mathrm{Mixed}}{\longleftrightarrow}}                        \gamma_{L^{\prime}}             \big] \text{ } +   \text{ }    \mathcal{P}^{\chi_1}_{\Lambda} \big[     \big(   \gamma_{-,+}  \big)_2          \underset{\Lambda \cap (\mathscr{F}\mathscr{C} )_1}{\overset{+}{\longleftrightarrow}}                     \mathscr{F}\mathscr{C}_1                \big] \text{ }                \text{ } \text{ , } 
\end{align*}

\noindent in which crossings of the height function, from the distribution of $+$ and $-$ face variables from the Ashkin-Teller model, consist of connected components entirely consisting of $+$ faces, $-$ faces, in addition to $+ \cup -$ faces, in which the first $-$, or $+$, face of the $+ \cup -$ path which is connected to a connected component of the opposite sign within $\Lambda \cap \big(\mathscr{F} \mathscr{C}\big)^c$; within each component of the connectivity events decomposed above respectively under $\chi_2$ and $\chi_1$, the faces $\gamma_{-}$, $\big(\gamma_{+,-}\big)_1$, $\gamma_{L^{\prime}}$, and $\big(\gamma_{-,+}\big)_2$, satisfy,

\begin{align*}
d \big( \gamma_L , \mathscr{F}\mathscr{C}_1 \big)   > d \big( \gamma_L ,  \big( \gamma_{-,+} \big)_2      \big) > d \big( \gamma_L , \gamma_{L^{\prime}}       \big) > d \big(   \gamma_L , \big( \gamma_{-,+} \big)_1       \big) >      d \big(  \gamma_L , \gamma_{-}        \big)       \text{ } \text{ , } 
\end{align*}

\noindent in which the $\gamma_{-}$, followed by $\big(\gamma_{-,+} \big)_1$, $\gamma_{L^{\prime}}$, and $\big(\gamma_{-,+} \big)_2$, successively are more distance from the first face $\gamma_L$ at which the connectivity event begins.

\bigskip

    \begin{figure}
\begin{align*}
\includegraphics[width=0.81\columnwidth]{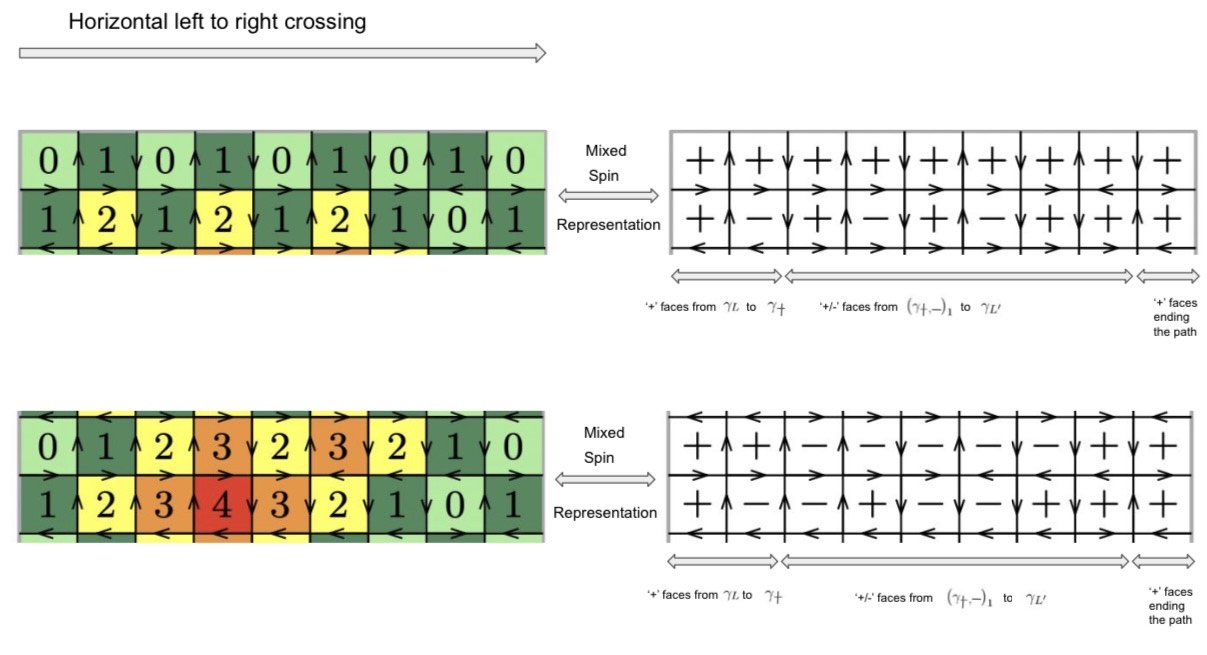}\\
\end{align*}
\caption{\textit{A second adaptation of mixed Ashkin-Teller configurations obtained over the square lattice as provided in Figure 3 of} {\color{blue}[16]}. From finite subvolumes of the square lattice as depicted above, passing to the Mixed Spin-Representation demonstrates that horizontal crossing events from the height function, which were originally considered earlier in the paper on the left throughout the strip and for directed, oriented, non-intersecting loops over the cylinder in the Six-vertex model, correspond to varying connected components of $+$ and $-$ faces in the Ashkin-Teller model.}
\end{figure}

\noindent Each crossing probability, respectively taken under $\chi_2$, and $\chi_1$, under the global symmetric flip sending $+ \longrightarrow -$ is equivalent to the superposition of three terms,

\begin{align*}
       \mathcal{P}^{-\chi_2}_{\Lambda} \big[     \gamma_L        \underset{\Lambda \cap (\mathscr{F}\mathscr{C})_1}{\overset{+}{\longleftrightarrow}}                 \gamma_{-}                   \big] \text{ }  +  \text{ }    \mathcal{P}^{-\chi_2}_{\Lambda} \big[       \big(  \gamma_{-,+}   \big)_1        \underset{\Lambda \cap (\mathscr{F}\mathscr{C}_1)^c}{\overset{\mathrm{Mixed}}{\longleftrightarrow}}                      \gamma_{L^{\prime}}             \big] \text{ } +   \text{ }    \mathcal{P}^{-\chi_2}_{\Lambda} \big[     \big(   \gamma_{-,+}  \big)_2          \underset{\Lambda \cap (\mathscr{F}\mathscr{C})_1}{\overset{-}{\longleftrightarrow}}                     \mathscr{F}\mathscr{C}_1                \big]          \text{ } \text{ , } 
\end{align*}

\noindent and also to,

\begin{align*}
          \mathcal{P}^{-\chi_1}_{\Lambda} \big[     \gamma_L        \underset{\Lambda \cap (\mathscr{F}\mathscr{C})_1}{\overset{+}{\longleftrightarrow}}                 \gamma_{-}                   \big] \text{ }  +  \text{ }    \mathcal{P}^{-\chi_1}_{\Lambda} \big[       \big(  \gamma_{-,+}   \big)_1         \underset{\Lambda \cap (\mathscr{F}\mathscr{C}_1)^c}{\overset{\mathrm{Mixed}}{\longleftrightarrow}}                       \gamma_{L^{\prime}}             \big] \text{ } +   \text{ }    \mathcal{P}^{-\chi_1}_{\Lambda} \big[     \big(   \gamma_{-,+}  \big)_2          \underset{\Lambda \cap (\mathscr{F}\mathscr{C})_1}{\overset{-}{\longleftrightarrow}}                     \mathscr{F}\mathscr{C}_1                \big]            \text{ } \text{ , } 
\end{align*}

\noindent in which, appearing in both of the two random variables above, the face variable that is used to condition on the first component of the crossing occurring is instead dependent upon $+$ face variables instead of $-$ face variables, which is then followed by a change of sign to a connected component of $-$ face variables between $\gamma_{-}$ and $\big( \gamma_{-,+}\big)_1$.

\bigskip

\noindent Denoting,

\begin{align*}
 \mathcal{P}^{\chi_2}_{-,+\cup - , -}  \big( \Lambda \cap (\mathscr{F}\mathscr{C})_1 ) \equiv  \mathcal{P}^{\chi_2}_{-,+\cup - , -}  \equiv  \text{ }   \mathcal{P}^{\chi_2}_{\Lambda} \big[     \gamma_L        \underset{\Lambda \cap (\mathscr{F}\mathscr{C})_1}{\overset{-}{\longleftrightarrow}}                 \gamma_{-}                   \big] \text{ }  +  \text{ }    \mathcal{P}^{\chi_2}_{\Lambda} \big[       \big(  \gamma_{-,+}   \big)_1          \underset{\Lambda \cap (\mathscr{F}\mathscr{C}_1)^c}{\overset{\mathrm{Mixed}}{\longleftrightarrow}}                 \gamma_{L^{\prime}}             \big] \text{ } \\ +   \text{ }    \mathcal{P}^{\chi_2}_{\Lambda} \big[     \big(   \gamma_{-,+}  \big)_2          \underset{\Lambda \cap (\mathscr{F}\mathscr{C})_1}{\overset{+}{\longleftrightarrow}}                     \mathscr{F}\mathscr{C}_1                \big]         \text{ } \text{ , } 
\end{align*}

\noindent in addition to the quantity,

\begin{align*}
  \mathcal{P}^{\chi_1}_{-,+\cup - , -}  \big( \Lambda \cap \big(\mathscr{F}\mathscr{C}\big)_1 \big) \equiv  \mathcal{P}^{\chi_1}_{-,+\cup - , -} \equiv \text{ }      \mathcal{P}^{\chi_1}_{\Lambda} \big[     \gamma_L        \underset{\Lambda \cap (\mathscr{F}\mathscr{C})_1}{\overset{-}{\longleftrightarrow}}                 \gamma_{-}                   \big] \text{ }  +  \text{ }    \mathcal{P}^{\chi_1}_{\Lambda} \big[       \big(  \gamma_{-,+}   \big)_1      \underset{\Lambda \cap (\mathscr{F}\mathscr{C}_1)^c}{\overset{\mathrm{Mixed}}{\longleftrightarrow}}                     \gamma_{L^{\prime}}             \big] \text{ } \\ +   \text{ }    \mathcal{P}^{\chi_1}_{\Lambda} \big[     \big(   \gamma_{-,+}  \big)_2          \underset{\Lambda \cap (\mathscr{F}\mathscr{C})_1}{\overset{+}{\longleftrightarrow}}                     \mathscr{F}\mathscr{C}_1                \big]        \text{ } \text{ , } 
\end{align*}

\noindent allows for grouping together of terms from $+$, $-$, and $- \cup +$, components under the absolute value,

\begin{align*}
    \big| \mathcal{P}^{\chi_2}_{-,+\cup - , -} -      \mathcal{P}^{\chi_1}_{-,+\cup - , -}  \big| \overset{\text{global face variable flip}}{\equiv}  \big| \mathcal{P}^{\chi_2}_{+,-\cup + , +} -      \mathcal{P}^{\chi_1}_{+,-\cup + , +}  \big|   \text{ } \text{ , } 
\end{align*}

\noindent implying that it suffices to individually bound each of the three following contributions from $-$, $- \cup +$ and $+$ face variables,

\begin{align*}
\big| \mathcal{P} \textit{1} \big|  \equiv  \bigg|\mathcal{P}^{\chi_2}_{\Lambda} \big[     \gamma_L        \underset{\Lambda \cap (\mathscr{F}\mathscr{C})_1}{\overset{-}{\longleftrightarrow}}                 \gamma_{-}                   \big]   -  \mathcal{P}^{\chi_1}_{\Lambda} \big[     \gamma_L        \underset{\Lambda \cap (\mathscr{F}\mathscr{C})_1}{\overset{-}{\longleftrightarrow}}                 \gamma_{-}                   \big] \bigg| \equiv    \bigg|\mathcal{P}^{-\chi_2}_{\Lambda} \big[     \gamma_L        \underset{\Lambda \cap (\mathscr{F}\mathscr{C})_1}{\overset{+}{\longleftrightarrow}}                 \gamma_{-}                   \big]  \\ -  \mathcal{P}^{-\chi_1}_{\Lambda} \big[     \gamma_L        \underset{\Lambda \cap (\mathscr{F}\mathscr{C})_1}{\overset{+}{\longleftrightarrow}}                 \gamma_{-}                   \big] \bigg|    \text{ } \text{ , }      \\ \big| \mathcal{P} \textit{2} \big|  \equiv  \bigg|   \mathcal{P}^{\chi_2}_{\Lambda} \big[       \big(  \gamma_{-,+}   \big)_1        \underset{\Lambda \cap (\mathscr{F}\mathscr{C}_1)^c}{\overset{\mathrm{Mixed}}{\longleftrightarrow}}                     \gamma_{L^{\prime}}             \big]               -    \mathcal{P}^{\chi_1}_{\Lambda} \big[       \big(  \gamma_{-,+}   \big)_1        \underset{\Lambda \cap (\mathscr{F}\mathscr{C})_1}{\overset{\mathrm{Mixed}}{\longleftrightarrow}}                        \gamma_{L^{\prime}}             \big]        \bigg| \equiv    \bigg|   \mathcal{P}^{-\chi_2}_{\Lambda} \big[       \big(  \gamma_{-,+}   \big)_1        \underset{\Lambda \cap (\mathscr{F}\mathscr{C}_1)^c}{\overset{\mathrm{Mixed}}{\longleftrightarrow}}                      \gamma_{L^{\prime}}             \big]             \\ -  \mathcal{P}^{-\chi_1}_{\Lambda} \big[       \big(  \gamma_{-,+}   \big)_1        \underset{\Lambda \cap (\mathscr{F}\mathscr{C})_1}{\overset{\mathrm{Mixed}}{\longleftrightarrow}}                        \gamma_{L^{\prime}}             \big]        \bigg|  \text{ } \text{ , } \\
\big| \mathcal{P} \textit{3} \big|  \equiv  \bigg| \mathcal{P}^{\chi_2}_{\Lambda} \big[     \big(   \gamma_{-,+}  \big)_2          \underset{\Lambda \cap (\mathscr{F}\mathscr{C})_1}{\overset{+}{\longleftrightarrow}}                     \mathscr{F}\mathscr{C}_1                \big]   - \mathcal{P}^{\chi_1}_{\Lambda} \big[     \big(   \gamma_{-,+}  \big)_2          \underset{\Lambda \cap (\mathscr{F}\mathscr{C})_1}{\overset{+}{\longleftrightarrow}}                     \mathscr{F}\mathscr{C}_1                \big]   \bigg| \equiv \bigg| \mathcal{P}^{\chi_2}_{\Lambda} \big[     \big(   \gamma_{-,+}  \big)_2          \underset{\Lambda \cap  ( \mathscr{F}\mathscr{C})_1}{\overset{+}{\longleftrightarrow}}                     \mathscr{F}\mathscr{C}_1                \big]  \\ -   \mathcal{P}^{\chi_1}_{\Lambda} \big[     \big(   \gamma_{-,+}  \big)_2          \underset{\Lambda \cap (\mathscr{F}\mathscr{C})_1}{\overset{+}{\longleftrightarrow}}                     \mathscr{F}\mathscr{C}_1                \big]   \bigg|  \text{ , } 
\end{align*}

\noindent which, from the denominations introduced above, is upper bounded by the difference below as a result of the triangle inequality,

\begin{align*}
    \bigg|  \mathcal{P} \textit{1}  +     \mathcal{P} \textit{2}  +  \mathcal{P} \textit{3} \bigg|   \leq \big| \mathcal{P} \textit{1} \big| + \big| \mathcal{P} \textit{2} \big| +  \big| \mathcal{P} \textit{3} \big| \text{ } \text{ . } \tag{$\Delta$}
\end{align*}

\noindent To achieve such an upper bound, of the form given above, which will be further analyzed, beginning with the random variable given in $\mathcal{P} \textit{1}$,

\begin{align*}
\mathcal{P} \textit{1} \equiv \mathcal{P}^{\chi_2}_{\Lambda} \big[     \gamma_L        \underset{\Lambda \cap (\mathscr{F}\mathscr{C})_1}{\overset{-}{\longleftrightarrow}}                 \gamma_{-}                   \big]   -  \mathcal{P}^{\chi_1}_{\Lambda} \big[     \gamma_L        \underset{\Lambda \cap (\mathscr{F}\mathscr{C})_1}{\overset{-}{\longleftrightarrow}}                 \gamma_{-}                   \big]     \equiv   \mathcal{P}^{\chi_2}_{\Lambda} \big[     \gamma_L        \underset{\Lambda \cap (\mathscr{F}\mathscr{C})_1}{\overset{-}{\longleftrightarrow}}                 \gamma_{-}                   \big] \bigg[ 1 - \frac{\mathcal{P}^{\chi_1}_{\Lambda} \big[     \gamma_L        \underset{\Lambda \cap (\mathscr{F}\mathscr{C})_1}{\overset{-}{\longleftrightarrow}}                 \gamma_{-}                   \big]}{\mathcal{P}^{\chi_2}_{\Lambda} \big[     \gamma_L        \underset{\Lambda \cap  ( \mathscr{F}\mathscr{C} )_1}{\overset{-}{\longleftrightarrow}}                 \gamma_{-}                   \big]} \bigg]     \text{ } \text{ , } 
\end{align*}

\noindent from which it suffices to demonstrate that the ratio of crossing probabilities between $\gamma_L$ and $\gamma_{-}$ is strictly less than $1$, as a result of the following observation,

\begin{align*}
  \frac{\mathcal{P}^{\chi_1}_{\Lambda} \big[     \gamma_L        \underset{\Lambda \cap (\mathscr{F}\mathscr{C} )_1}{\overset{-}{\longleftrightarrow}}                 \gamma_{-}                   \big]}{\mathcal{P}^{\chi_2}_{\Lambda} \big[     \gamma_L        \underset{\Lambda \cap (\mathscr{F}\mathscr{C})_1}{\overset{-}{\longleftrightarrow}}                 \gamma_{-}                   \big]} \equiv         \frac{\underset{\chi_1 \in \{ + , + \}}{\sum}             \mathcal{P}^{\chi_1}_{\Lambda} \bigg[    \gamma_L        \underset{\Lambda \cap (\mathscr{F}\mathscr{C})_1}{\overset{-}{\longleftrightarrow}}                 \gamma_{-}        \big|                       \mathcal{C}_{-}   \cup  \mathscr{D}_{\chi,-}   \bigg]      \mathcal{P}^{\chi_1}_{\Lambda} \big[      \mathcal{C}_{-}   \cup  \mathscr{D}_{\chi,-}        \big]          }{\underset{\chi_2 \in \{ - , + \}}{\sum}     \mathcal{P}^{\chi_2}_{\Lambda} \bigg[    \gamma_L        \underset{\Lambda \cap (\mathscr{F}\mathscr{C} )_1}{\overset{-}{\longleftrightarrow}}                 \gamma_{-}        \big|      \mathcal{C}_{-}   \cup  \mathscr{D}_{\chi,-}        \bigg]      \mathcal{P}^{\chi_2}_{\Lambda} \big[       \mathcal{C}_{-}   \cup  \mathscr{D}_{\chi,-}         \big]   }   \tag{$\chi_1$-$\chi_2$ \textit{ratio}}   \\ \leq    \frac{\underset{\chi_1 \in \{ + , + \}}{\sum}             \mathcal{P}^{\chi_1}_{\Lambda} \bigg[    \gamma_L        \underset{\Lambda \cap (\mathscr{F}\mathscr{C})_1}{\overset{-}{\longleftrightarrow}}                 \gamma_{-}        \big|                       \mathcal{C}_{-}   \cup  \mathscr{D}_{\chi,-}   \bigg]      \mathcal{P}^{\chi_1}_{\Lambda} \big[       \mathcal{C}_{-}   \cup  \mathscr{D}_{\chi,-}        \big]          }{\underset{\chi_2 \in \{ + , + \}}{\sum}     \mathcal{P}^{\chi_2}_{\Lambda} \bigg[    \gamma_L        \underset{\Lambda \cap (\mathscr{F}\mathscr{C})_1}{\overset{-}{\longleftrightarrow}}                 \gamma_{-}        \big|      \mathcal{C}_{-}   \cup  \mathscr{D}_{\chi,-}        \bigg]      \mathcal{P}^{\chi_2}_{\Lambda} \big[       \mathcal{C}_{-}   \cup  \mathscr{D}_{\chi,-}      \big]                }      \text{ } \text{ , }  \tag{Event $8$}
\end{align*}

\noindent where the conditioning in the crossing probabilities taken under $\chi_1$ and $\chi_2$ is,

\begin{align*}
 \mathcal{C}_{-} \big(    \Lambda , \mathscr{L}_{j,-}   , \gamma_{-}   \big) \equiv    \mathcal{C}_{-} \equiv   \underset{j \in \textbf{N}}{\bigcap}    \bigg\{         \gamma_{j-1,-}                            \underset{\Lambda \cap ( ( \mathscr{F}\mathscr{C} )_1 \cap     \mathscr{L}_{j,-}             )}{\overset{-}{\longleftrightarrow}}         \gamma_{-}                     \bigg\}           \text{ } \text{ , } 
\end{align*}

\noindent where in the intersection given above, over $j \in \textbf{N}$, the paths $\gamma_{j-1,-}$ consist of faces after the first face $\gamma_L$ with which the path begins, which connect to $\gamma_{-}$ with a path of strictly $-$ face variables, each of which intersects the intersection of the \textit{freezing cluster} with a line $\mathscr{L}_{j,-}$ intersecting the top and bottom boundaries of the strip; the remaining item in the conditioning that is included in the union with $\mathcal{C}_{-}$, respectively under $\mathcal{P}^{\chi_1}_{\Lambda} \big[ \cdot \big]$, and $\mathcal{P}^{\chi_2}_{\Lambda} \big[ \cdot \big]$, 

\begin{align*}
   \mathcal{D}_{\chi,-} \big(    \Lambda , \mathscr{L}_{j,-}  , \gamma_{-}    \big) \equiv    \mathscr{D}_{\chi,-} \equiv    \big\{                F \big( \gamma_L \cap               \gamma_{j-1,-}      \big)  \in -     \text{ over }  \Lambda \cap \big( \big( \mathscr{F}\mathscr{C}\big)_1 \cap     \mathscr{L}_{j,-}          \big)      \big\}     \text{ } \text{ , } 
\end{align*}

\noindent indicates the existence of suitable domain for which,

\begin{align*}
        \mathcal{P}^{\chi_2 \text{ or } \chi_1}_{\Lambda} \bigg[                  \gamma_L        \underset{\Lambda \cap (\mathscr{F}\mathscr{C} )_1}{\overset{-}{\longleftrightarrow}}                 \gamma_{-}  \text{ }       \big|      \text{ }                  \mathcal{C}_{-}   \cup  \mathscr{D}_{\chi,-}        \bigg]           \text{ } \text{ , } \tag{Event $9$}
\end{align*}

\noindent occurs with positive probability. Therefore, the summation over $\mathscr{F}$, 

\begin{align*}
     (\chi_1 - \chi_2 \textit{ ratio}) \leq   \underset{\chi_3 \in \{ + , + \}}{\underset{\emptyset \neq \mathscr{F} \in F ( \Lambda \cap ( \mathscr{F} \mathscr{C})_1 )}{\sum}}         \mathcal{P}^{\chi_3}_{\Lambda} \bigg[     \gamma_L \underset{\Lambda \cap ( \mathscr{F}\mathscr{C} )_1}{\overset{-}{\longleftrightarrow}}         \mathscr{F}  \big| \mathcal{C}_{-} \big( \Lambda , \mathscr{L}_{j,-} , \mathscr{F} \big) \cup   \mathscr{D}_{\chi_3 , -} \big( \Lambda , \mathscr{L}_{j,-} , \mathscr{F} \big)     \bigg]  \mathcal{P}^{\chi_3}_{\Lambda} \big[      \mathcal{C}_{-}   \cup  \mathscr{D}_{\chi,-}       \big]  \\ \equiv   \underset{\chi_3 \in \{ + , + \}}{\underset{\emptyset \neq\mathscr{F} \in F ( \Lambda \cap ( \mathscr{F} \mathscr{C} )_1 )}{\sum}}         \mathcal{P}^{\chi_3}_{\Lambda} \bigg[     \gamma_L \underset{\Lambda \cap ( \mathscr{F}\mathscr{C})_1}{\overset{-}{\longleftrightarrow}}         \mathscr{F}  \big| \mathcal{C}_{-}  \cup   \mathscr{D}_{\chi_3 , -}     \bigg]   \mathcal{P}^{\chi_3}_{\Lambda} \big[     \mathcal{C}_{-}   \cup  \mathscr{D}_{\chi,-}         \big]       \text{ } \text{ , } \tag{$\chi_1$-$\chi_2$ \textit{ratio 2}} \\ \tag{Event $10$}
\end{align*}

\noindent the final expression above can be upper bounded with,

\begin{align*}
 \delta^{\prime\prime}_2  , \delta^{\prime\prime}_1 \equiv \underset{k>0}{\mathrm{inf}} \big\{    \big( \delta^{\prime\prime}_1 \big)_{k}     \big\}  > 0         \Longrightarrow  ( \chi_1 - \chi_2 \textit{ ratio 2}  ) \geq \underset{k>0}{\sum}    \big( \delta^{\prime\prime}_1 \big)_k  \delta^{\prime\prime}_2            \\ \equiv  \big( \delta^{\prime\prime}_1 \underset{k-2\text{ copies}}{\underbrace{+ \cdots +}}  \delta^{\prime\prime}_1  \big) \delta^{\prime\prime}_2   \equiv k \delta^{\prime\prime}_1  \delta^{\prime\prime}_2     \text{ } \text{ , } 
\end{align*}

\noindent under the assumption that $\delta^{\prime\prime}_1$ above satisfies,

\begin{align*}
\mathcal{P}^{\chi_3}_{\Lambda} \bigg[     \gamma_L \underset{\Lambda \cap ( \mathscr{F}\mathscr{C} )_1}{\overset{-}{\longleftrightarrow}}         \mathscr{F}_k  \big| \mathcal{C}_{-}  \cup   \mathscr{D}_{\chi_3 , -}     \bigg] \geq \big(\delta^{\prime\prime}_1\big)_k \text{ } \forall k \Longleftrightarrow  \underset{k>0}{\mathrm{inf}} \bigg\{   \mathcal{P}^{\chi_3}_{\Lambda} \bigg[     \gamma_L \underset{\Lambda \cap ( \mathscr{F}\mathscr{C})_1}{\overset{-}{\longleftrightarrow}}         \mathscr{F}_k  \big| \mathcal{C}_{-}  \cup   \mathscr{D}_{\chi_3 , -}     \bigg]   \bigg\}   \equiv \big(\delta^{\prime\prime}_1\big)_k        \text{ } \text{ , } \tag{Event $11$}
\end{align*}

\noindent for faces satisfying,

\begin{align*}
  \mathscr{F}_k \cap F \big( \Lambda \cap \big( \mathscr{F} \mathscr{C} \big)_1 \big) \neq \emptyset  \text{ } \text{ , }
\end{align*}

\noindent in addition to the requirement that the probability below, without conditioning on the occurrence of, 

\begin{align*}
 \mathcal{P}^{\chi_3}_{\Lambda} \big[  \mathcal{C}_{-} \cup \mathscr{D}_{\chi_3, -}  \big]  \text{ } \text{ , } 
\end{align*}

\noindent admit the following estimate,

\begin{align*}
  \mathcal{P}^{\chi_1}_{\Lambda} \big[     \gamma_L        \underset{\Lambda \cap (\mathscr{F}\mathscr{C})_1}{\overset{-}{\longleftrightarrow}}                 \gamma_{-}        \big]   \geq \delta^{\prime\prime}_2    \text{ } \text{ , }
\end{align*}

\noindent for a suitably chosen, absolute constant. Collecting estimates from previous steps hence provides the lower bound for,

\begin{align*}
     \mathcal{P} \textit{1} \equiv \mathcal{P}^{\chi_2}_{\Lambda} \big[     \gamma_L        \underset{\Lambda \cap (\mathscr{F}\mathscr{C})_1}{\overset{-}{\longleftrightarrow}}                 \gamma_{-}                   \big]   -  \mathcal{P}^{\chi_1}_{\Lambda} \big[     \gamma_L        \underset{\Lambda \cap (\mathscr{F}\mathscr{C})_1}{\overset{-}{\longleftrightarrow}}                 \gamma_{-}                   \big]            \text{ } \text{ , } \tag{Event $12$}
\end{align*}

    \begin{figure}
\begin{align*}
\includegraphics[width=0.95\columnwidth]{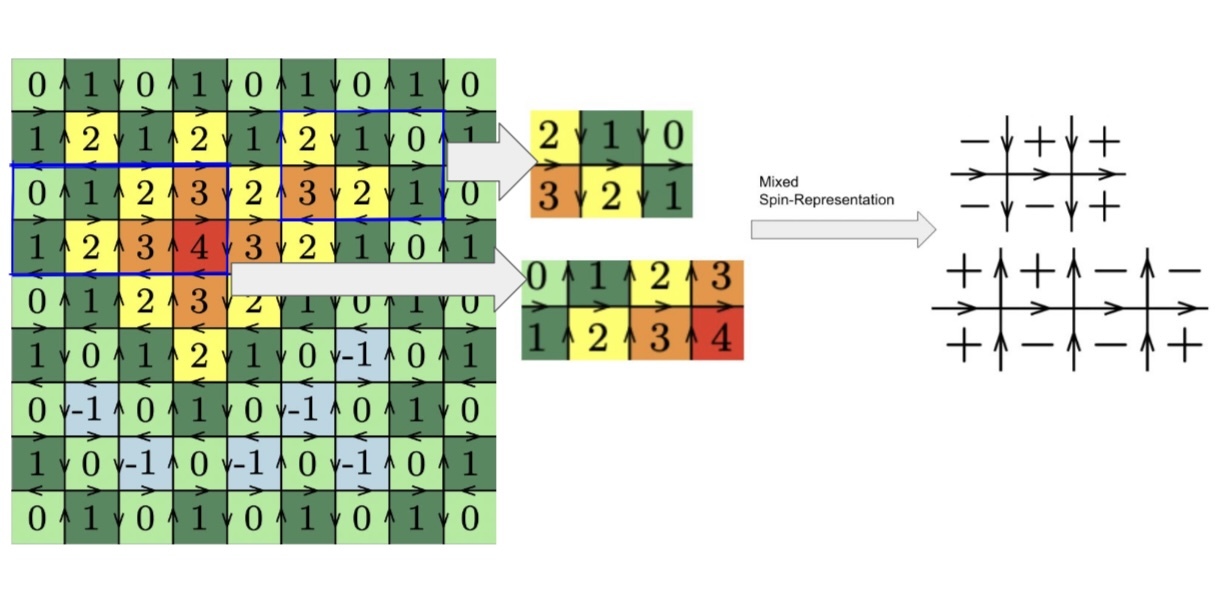}\\
\end{align*}
\caption{\textit{An adaptation of mixed Ashkin-Teller configurations obtained over the square lattice as provided in Figure 3 of} {\color{blue}[16]}. Above, for a subset of faces enclosed for the given finite volume, from the example configuration provided in {\color{blue}[16]} one can observe subregions over which faces from the graph homomorphism of the six-vertex model are frozen, which are precursors to larger collections of faces that can become frozen under sloped boundary conditions. When passing to the Mixed Ashkin-Teller spin representation from the six-vertex configuration, the collections of faces over the square lattice enclosed by the two {\color{blue}blue} finite volumes correspond to distribution of $+$ and $-$ face variables as shown, in which diagonal faces with opposite spins are assigned in correspondence to differences in the image of the height function from the six-vertex configuration. }
\end{figure}

\noindent in which,

\begin{align*}
        \mathcal{P} \textit{1}  \geq k \delta^{\prime\prime}_1   \delta^{\prime\prime}_2      \text{ } \text{ , } 
\end{align*}

\noindent as desired.

\bigskip

\noindent Similar rearrangements, as above, imply the existence of another suitable, strictly positive, lower bound for $\mathcal{P} \textit{2}$, as,

\begin{align*}
  \mathcal{P} \textit{2} \equiv      \mathcal{P}^{\chi_2}_{\Lambda} \big[       \big(  \gamma_{-,+}   \big)_1      \underset{\Lambda \cap (\mathscr{F}\mathscr{C}_1)^c}{\overset{\mathrm{Mixed}}{\longleftrightarrow}}                      \gamma_{L^{\prime}}             \big]               -    \mathcal{P}^{\chi_1}_{\Lambda} \big[       \big(  \gamma_{-,+}   \big)_1          \underset{\Lambda \cap (\mathscr{F}\mathscr{C}_1)^c}{\overset{\mathrm{Mixed}}{\longleftrightarrow}}                       \gamma_{L^{\prime}}             \big]       \equiv   \mathcal{P}^{\chi_2}_{\Lambda} \big[       \big(  \gamma_{-,+}   \big)_1        \underset{\Lambda \cap (\mathscr{F}\mathscr{C}_1)^c}{\overset{\mathrm{Mixed}}{\longleftrightarrow}}                     \gamma_{L^{\prime}}             \big]\\  \times \bigg[  1  - \frac{\mathcal{P}^{\chi_1}_{\Lambda} \big[       \big(  \gamma_{-,+}   \big)_1          \underset{\Lambda \cap (\mathscr{F}\mathscr{C}_1)^c}{\overset{\mathrm{Mixed}}{\longleftrightarrow}}                    \gamma_{L^{\prime}}             \big]   }{\mathcal{P}^{\chi_2}_{\Lambda} \big[       \big(  \gamma_{-,+}   \big)_1      \underset{\Lambda \cap (\mathscr{F}\mathscr{C}_1)^c}{\overset{\mathrm{Mixed}}{\longleftrightarrow}}                     \gamma_{L^{\prime}}             \big]   }  \bigg]        \text{ , }   \tag{Event $13$}
\end{align*}

\noindent from which implementing previous steps to obtain the lower bound yields, 

\begin{align*}
    \frac{\mathcal{P}^{\chi_1}_{\Lambda} \big[       \big(  \gamma_{-,+}   \big)_1        \underset{\Lambda \cap (\mathscr{F}\mathscr{C}_1)^c}{\overset{\mathrm{Mixed}}{\longleftrightarrow}}                      \gamma_{L^{\prime}}             \big]   }{\mathcal{P}^{\chi_2}_{\Lambda} \big[       \big(  \gamma_{-,+}   \big)_1       \underset{\Lambda \cap (\mathscr{F}\mathscr{C}_1)^c}{\overset{\mathrm{Mixed}}{\longleftrightarrow}}                   \gamma_{L^{\prime}}             \big]   } \leq \underset{\chi_4 \in \{ + , + \}}{\underset{\emptyset\neq\mathscr{F} \neq \mathscr{F}^{\prime} \in F ( \Lambda \cap ( \mathscr{F} \mathscr{C} )_1 )}{\sum}} \mathcal{P}^{\chi_4}_{\Lambda} \bigg[              \big(  \gamma_{-,+}   \big)_1         \underset{\Lambda \cap (\mathscr{F}\mathscr{C}_1)^c}{\overset{\mathrm{Mixed}}{\longleftrightarrow}}                    \mathscr{F}^{\prime} \big|     \mathcal{C}_{-} \big( \Lambda , \mathscr{L}_{j,-} , \mathscr{F}^{\prime} \big) \\   \cup \mathscr{D}_{\chi_3 , -} \big( \Lambda , \mathscr{L}_{j,-} , \mathscr{F}^{\prime} \big)                           \bigg]  \mathcal{P}^{\chi_4}_{\Lambda} \big[              \big(  \gamma_{-,+}   \big)_1          \underset{\Lambda \cap (\mathscr{F}\mathscr{C}_1)^c}{\overset{\mathrm{Mixed}}{\longleftrightarrow}}                   \mathscr{F}^{\prime} \big]    \\ \equiv  \underset{\emptyset\neq\mathscr{F} \neq \mathscr{F}^{\prime} \in F ( \Lambda \cap ( \mathscr{F} \mathscr{C} )_1 )}{\sum} \mathcal{P}^{\chi_4}_{\Lambda} \bigg[              \big(  \gamma_{-,+}   \big)_1       \underset{\Lambda \cap (\mathscr{F}\mathscr{C}_1)^c}{\overset{\mathrm{Mixed}}{\longleftrightarrow}}                    \mathscr{F}^{\prime} \big|     \mathcal{C}_{-}  \cup     \mathscr{D}_{\chi_3 , -}                       \bigg] \\ \times 
    \mathcal{P}^{\chi_4}_{\Lambda} \big[              \big(  \gamma_{-,+}   \big)_1       \underset{\Lambda \cap (\mathscr{F}\mathscr{C}_1)^c}{\overset{\mathrm{Mixed}}{\longleftrightarrow}}                     \mathscr{F}^{\prime} \big]   \text{ } \text{ . } 
\end{align*}

\noindent Concluding for $\mathcal{P} \textit{2}$, from the summation over $\mathscr{F}^{\prime}$ above,

\begin{align*}
    \underset{\emptyset\neq\mathscr{F} \neq \mathscr{F}^{\prime} \in F ( \Lambda \cap ( \mathscr{F} \mathscr{C} )_1 )}{\sum} \mathcal{P}^{\chi_4}_{\Lambda} \bigg[              \big(  \gamma_{-,+}   \big)_1        \underset{\Lambda \cap (\mathscr{F}\mathscr{C}_1)^c}{\overset{\mathrm{Mixed}}{\longleftrightarrow}}                       \mathscr{F}^{\prime} \big|     \mathcal{C}_{-}  \cup   \mathscr{D}_{\chi_3 , -}                       \bigg]  \geq       \underset{\emptyset\neq\mathscr{F} \neq \mathscr{F}^{\prime} \in F ( \Lambda \cap ( \mathscr{F} \mathscr{C} )_1 )}{\sum} \delta_4 \big( \mathscr{F}^{\prime} \big)                   \delta^{\prime}_4 \big( \mathscr{F}^{\prime} \big)  \\ \equiv        \underset{\emptyset\neq\mathscr{F} \neq \mathscr{F}^{\prime} \in F ( \Lambda \cap ( \mathscr{F} \mathscr{C} )_1 )}{\sum} \Delta_{4} \big( \mathscr{F} \big)  \\ \geq  \big( \Delta_{4}     \big)^{|\mathscr{F}^{\prime} | } \\ \geq \big( \Delta_{4} \big)^{N}      \text{ } \text{ , } 
\end{align*}

\noindent given the existence of the following several quantities, the first of which is,

\begin{align*}
         \delta_4 \big( \mathscr{F}^{\prime} \big)   \equiv     \mathrm{inf} \bigg\{                 \mathcal{P}^{\chi_4}_{\Lambda} \bigg[              \big(  \gamma_{-,+}   \big)_1          \underset{\Lambda \cap (\mathscr{F}\mathscr{C}_1)^c}{\overset{\mathrm{Mixed}}{\longleftrightarrow}}                      \mathscr{F}^{\prime} \big|     \mathcal{C}_{-}  \cup   \mathscr{D}_{\chi_3 , -}                       \bigg]             \bigg\} \leq  \mathcal{P}^{\chi_4}_{\Lambda} \bigg[              \big(  \gamma_{-,+}   \big)_1          \underset{\Lambda \cap (\mathscr{F}\mathscr{C}_1)^c}{\overset{\mathrm{Mixed}}{\longleftrightarrow}}                       \mathscr{F}^{\prime} \big|     \mathcal{C}_{-}  \cup   \mathscr{D}_{\chi_3 , -}                       \bigg]                    \text{ } \text{ , } 
\end{align*}

\noindent for faces satisfying,

\begin{align*}
     \mathscr{F}^{\prime} \cap F \big( \Lambda \cap \big( \mathscr{F} \mathscr{C} \big)_1 \big) \neq \emptyset        \text{ , } 
\end{align*}

\noindent the second of which is,

\begin{align*}
      \delta^{\prime}_4 \big( \mathscr{F}^{\prime} \big)  \leq      \mathcal{P}^{\chi_4}_{\Lambda} \big[              \mathcal{C}_{-}  \cup   \mathscr{D}_{\chi_3 , -}                      \big]      \text{ } \text{ , } 
\end{align*}

\noindent the third of which is,

\begin{align*}
  \Delta_{4} \equiv \delta_4 \big( \mathscr{F}^{\prime} \big) \delta^{\prime}_4 \big( \mathscr{F}^{\prime} \big)   \text{ }\text{ , } 
\end{align*}

\noindent and the fourth of which is,

\begin{align*}
  N < | \mathscr{F}^{\prime} |  \text{ }\text{ , } 
\end{align*}

\noindent which altogether hence imply the lower bound,

\begin{align*}
  \mathcal{P} \textit{2} \geq   \big( \Delta_4 \big)^N \text{ }   \equiv \text{ } \Delta^N_4   \text{ } \text{ . } 
\end{align*}

\noindent Lastly, concluding the argument with the lower bound for $\mathcal{P} \textit{3}$ yields,

\begin{align*}
    \text{ }   \mathcal{P} \textit{3} \equiv      \mathcal{P}^{\chi_2}_{\Lambda} \big[     \big(   \gamma_{-,+}  \big)_2          \underset{\Lambda \cap (\mathscr{F}\mathscr{C})_1}{\overset{+}{\longleftrightarrow}}                     \mathscr{F}\mathscr{C}_1                \big]   - \mathcal{P}^{\chi_1}_{\Lambda} \big[     \big(   \gamma_{-,+}  \big)_2          \underset{\Lambda \cap (\mathscr{F}\mathscr{C})_1}{\overset{+}{\longleftrightarrow}}                     \mathscr{F}\mathscr{C}_1                \big]                  \tag{Event $14$}      \\ \equiv  \mathcal{P}^{\chi_2}_{\Lambda} \big[     \big(   \gamma_{-,+}  \big)_2          \underset{\Lambda \cap (\mathscr{F}\mathscr{C} )_1}{\overset{+}{\longleftrightarrow}}                     \mathscr{F}\mathscr{C}_1                \big]   \bigg[ 1 - \frac{ \mathcal{P}^{\chi_1}_{\Lambda} \big[     \big(   \gamma_{-,+}  \big)_2          \underset{\Lambda \cap (\mathscr{F}\mathscr{C})_1}{\overset{+}{\longleftrightarrow}}                     \mathscr{F}\mathscr{C}_1                \big] }{ \mathcal{P}^{\chi_2}_{\Lambda} \big[     \big(   \gamma_{-,+}  \big)_2          \underset{\Lambda \cap (\mathscr{F}\mathscr{C})_1}{\overset{+}{\longleftrightarrow}}                     \mathscr{F}\mathscr{C}_1                \big] } \bigg] \text{ } \text{ . }  \end{align*}

\noindent The condition above holds iff, given the existence of two suitable constants for which,

    \begin{align*}
 \forall \mathscr{F}\mathscr{C}_1 \in F^{\mathrm{odd}} \text{ } \text{ or } \text{ } F^{\mathrm{even}} \text{ } , \exists \text{ }   \delta_5 > 0 :   \mathcal{P}^{\chi_2}_{\Lambda} \big[     \big(   \gamma_{-,+}  \big)_2          \underset{\Lambda \cap (\mathscr{F}\mathscr{C} )_1}{\overset{+}{\longleftrightarrow}}                     \mathscr{F}\mathscr{C}_1                \big] \geq \delta_5   \text{ }  \text{ , } \\  \forall \mathscr{F}\mathscr{C}_1 \in F^{\mathrm{odd}} \text{ } \text{ or } \text{ } F^{\mathrm{even}} \text{ } ,  \exists    \text{ 
} \big( \delta_5 \big)^{\prime} > 0  :    \frac{ \mathcal{P}^{\chi_1}_{\Lambda} \big[     \big(   \gamma_{-,+}  \big)_2          \underset{\Lambda \cap(\mathscr{F}\mathscr{C})_1}{\overset{+}{\longleftrightarrow}}                     \mathscr{F}\mathscr{C}_1                \big] }{ \mathcal{P}^{\chi_2}_{\Lambda} \big[     \big(   \gamma_{-,+}  \big)_2          \underset{\Lambda \cap (\mathscr{F}\mathscr{C})_1}{\overset{+}{\longleftrightarrow}}                     \mathscr{F}\mathscr{C}_1                \big] }  \leq  - \big( \delta_5 \big)^{\prime}   + 1          \text{ }  \text{ . } 
\end{align*}

\noindent To obtain the desired constants $\delta_5$ and $\big( \delta_5 \big)^{\prime}$ shown above, observe, first, that $\big( \delta_5 \big)^{\prime}$ takes the form,

\begin{align*}
      \frac{ \mathcal{P}^{\chi_1}_{\Lambda} \big[     \big(   \gamma_{-,+}  \big)_2          \underset{\Lambda \cap (\mathscr{F}\mathscr{C})_1}{\overset{+}{\longleftrightarrow}}                     \mathscr{F}\mathscr{C}_1                \big] }{ \mathcal{P}^{\chi_2}_{\Lambda} \big[     \big(   \gamma_{-,+}  \big)_2          \underset{\Lambda \cap (\mathscr{F}\mathscr{C})_1}{\overset{+}{\longleftrightarrow}}                     \mathscr{F}\mathscr{C}_1                \big] }     \\ \leq           \underset{\chi_5 \in \{ + , + \}}{\underset{\emptyset \neq \mathscr{F} \neq \mathscr{F}^{\prime} \neq \mathscr{F}^{\prime\prime} \in F ( \Lambda \cap (\mathscr{F} \mathscr{C})_1 )}{\sum}}          \mathcal{P}^{\chi_5}_{\Lambda} \bigg[          \big(   \gamma_{-,+}  \big)_2                       \underset{\Lambda \cap (\mathscr{F}\mathscr{C})_1}{\overset{+}{\longleftrightarrow}}    \mathscr{F}^{\prime\prime}          \big|         \mathcal{C}_{-} \big( \Lambda , \mathscr{L}_{j,-} , \mathscr{F}^{\prime\prime} \big) \cup     \mathscr{D}_{\chi_3 , -}    \big( \Lambda , \mathscr{L}_{j,-} , \mathscr{F}^{\prime\prime} \big)           \bigg]  \\ \times  \mathcal{P}^{\chi_5}_{\Lambda} \big[           \mathcal{C}_{-} \big( \Lambda , \mathscr{L}_{j,-} , \mathscr{F}^{\prime} \big) \cup    \mathscr{D}_{\chi_3 , -}                 \big] \\ \equiv       {\underset{\emptyset \neq \mathscr{F} \neq \mathscr{F}^{\prime} \neq \mathscr{F}^{\prime\prime} \in F ( \Lambda \cap (\mathscr{F} \mathscr{C})_1 )}{\sum}}          \mathcal{P}^{\chi_5}_{\Lambda} \bigg[          \big(   \gamma_{-,+}  \big)_2                       \underset{\Lambda \cap (\mathscr{F}\mathscr{C})_1}{\overset{+}{\longleftrightarrow}}    \mathscr{F}^{\prime\prime}          \big|         \mathcal{C}_{-} \cup    \mathscr{D}_{\chi_3 , -}              \bigg]   \mathcal{P}^{\chi_5}_{\Lambda} \big[           \mathcal{C}_{-}  \cup    \mathscr{D}_{\chi_3 , -}                 \big]      \text{ } \text{ . } \end{align*}
      
     \noindent The condition above holds iff, given suitable $\delta_5$ for which,

      \begin{align*}
      \big( \delta_5 \big)^{\prime} \overset{(*)}{\leq}     \big( \Delta_5 \big)^{N^{\prime}}    \\ \overset{(**)}{\leq}   \big( \Delta_5 \big)^{|\mathscr{F}^{\prime\prime}|}                  \\ \overset{(***)}{\leq}        {\underset{\emptyset \neq \mathscr{F} \neq \mathscr{F}^{\prime} \neq \mathscr{F}^{\prime\prime} \in F ( \Lambda \cap (\mathscr{F} \mathscr{C})_1 )}{\sum}}           \delta_5 \big(   \mathscr{F}^{\prime\prime}     \big)     \delta^{\prime}_5 \big( \mathscr{F}^{\prime\prime} \big) \\ 
      \overset{(****)}{\leq} {\underset{\emptyset \neq \mathscr{F} \neq \mathscr{F}^{\prime} \neq \mathscr{F}^{\prime\prime} \in F ( \Lambda \cap (\mathscr{F} \mathscr{C} )_1 )}{\sum}}         \mathcal{P}_{\Lambda} \bigg[              \big(   \gamma_{-,+}  \big)_2                       \underset{\Lambda \cap (\mathscr{F}\mathscr{C})_1}{\overset{+}{\longleftrightarrow}}    \mathscr{F}^{\prime\prime}          \big|         \mathcal{C}_{-} \cup    \mathscr{D}_{\chi_3 , -}              \bigg]  \mathcal{P}^{\chi_2}_{\Lambda} \big[           \mathcal{C}_{-}  \cup    \mathscr{D}_{\chi_3 , -}                 \big]                    \text{ } \text{ , } 
\end{align*}

\noindent where, in the sequence of rearrangements above, in the upper bound provided from the final expression in $(\textit{****})$, lower bounds for,

\begin{align*}
     \mathcal{P}_{\Lambda} \bigg[              \big(   \gamma_{-,+}  \big)_2                       \underset{\Lambda \cap (\mathscr{F}\mathscr{C} )_1}{\overset{+}{\longleftrightarrow}}    \mathscr{F}^{\prime\prime}          \big|         \mathcal{C}_{-} \cup    \mathscr{D}_{\chi_3 , -}              \bigg]    \text{ } \text{ , }              \tag{Event $15$}
\end{align*}

\noindent and for,

\begin{align*}
  \mathcal{P}^{\chi_2}_{\Lambda} \big[           \mathcal{C}_{-}  \cup    \mathscr{D}_{\chi_3 , -}                 \big]       \text{ } \text{ , } 
\end{align*}

\noindent are respectively provided with,

\begin{align*}
      \delta_5 \big( \mathscr{F}^{\prime\prime} \big)          \text{ } \text{ , } 
\end{align*}

\noindent which is obtained in light of the sequence of observations below,

\begin{align*}
        \mathcal{P}^{\chi_2}_{\Lambda} \bigg[              \big(   \gamma_{-,+}  \big)_2                       \underset{\Lambda \cap (\mathscr{F}\mathscr{C})_1}{\overset{+}{\longleftrightarrow}}    \mathscr{F}^{\prime\prime}          \big|         \mathcal{C}_{-} \cup    \mathscr{D}_{\chi_3 , -}              \bigg]  \equiv \mathcal{P}^{\chi_2}_{\Lambda} \bigg[        \underset{\emptyset \neq \mathscr{F}^{\prime\prime} \in F ( \Lambda \cap ( \mathscr{F} \mathscr{C} )_1 )}{\underset{\underset{i}{\mathrm{max}} ( \mathscr{F}^{\prime\prime}_i )  \cap   \mathscr{F}^{\prime\prime} \neq \emptyset    }{\underset{i \in \textbf{N}}{\bigcap}}}    \big\{       \big(   \gamma_{-,+}  \big)_2                       \underset{\Lambda \cap (\mathscr{F}\mathscr{C})_1}{\overset{+}{\longleftrightarrow}}    \mathscr{F}^{\prime\prime}_i        \big\}       \big|         \mathcal{C}_{-} \cup    \mathscr{D}_{\chi_3 , -}          \bigg] \\  \overset{(\mathrm{FKG})}{\geq}   \underset{\emptyset \neq \mathscr{F}^{\prime\prime} \in F ( \Lambda \cap ( \mathscr{F} \mathscr{C} )_1 )}{\underset{\underset{i}{\mathrm{max}} ( \mathscr{F}^{\prime\prime}_i )  \cap   \mathscr{F}^{\prime\prime} \neq \emptyset    }{\underset{ i \in \textbf{N}}{\prod}}}             \mathcal{P}^{\chi_2}_{\Lambda} \bigg[              \big(   \gamma_{-,+}  \big)_2                       \underset{\Lambda \cap (\mathscr{F}\mathscr{C})_1}{\overset{+}{\longleftrightarrow}}    \mathscr{F}^{\prime\prime}_i          \big|         \mathcal{C}_{-} \cup    \mathscr{D}_{\chi_3 , -}              \bigg] \\ \geq    \underset{\emptyset \neq \mathscr{F}^{\prime\prime} \in F ( \Lambda \cap ( \mathscr{F} \mathscr{C} )_1 )}{\underset{\underset{i}{\mathrm{max}} ( \mathscr{F}^{\prime\prime}_i )  \cap   \mathscr{F}^{\prime\prime} \neq \emptyset    }{\underset{ i \in \textbf{N}}{\prod}}}              \gamma_i   \equiv \big( \gamma_i \big)^{|\mathscr{F}^{\prime\prime}|}   \\   \equiv      \delta_5 \big( \mathscr{F}^{\prime\prime} \big)                              \text{ } \text{ , } 
\end{align*}

\noindent where, in the sequence of rearrangements above, following the application of (FKG), which is satisfied for the conditionally defined \textit{Ashkin-Teller measure}, from the conditioning enforced on the crossing event,

\begin{align*}
        \big\{    \big(   \gamma_{-,+}  \big)_2                       \underset{\Lambda \cap (\mathscr{F}\mathscr{C})_1}{\overset{+}{\longleftrightarrow}}    \mathscr{F}^{\prime\prime}_i            \big\}                     \text{ } \text{ , }
\end{align*}

\noindent provided above, for the following conditionally defined measure below, that is dependent upon $\big\{  \mathcal{C}_{-} \cup    \mathscr{D}_{\chi_3 , -}          \big\}$,

\begin{align*}
          \mathcal{P}^{\chi_2}_{\Lambda} \bigg[        \cdot          \big|   \text{ }       \mathcal{C}_{-} \cup    \mathscr{D}_{\chi_3 , -}              \bigg]                        \text{ } \text{ , }
\end{align*}

\noindent given some realization of a mixed \textit{Ashkin-Teller} configurations $\cdot \in \Omega^{\mathrm{AT}}$. Hence, each of the crossing probabilities in the intersection admits the lower bound from the crossing from the intersection over $i \in \textbf{N}$, in which, as argued previously in several other instances,

\begin{align*}
   \delta^{i^{\prime}}_5    \big(        \mathscr{F}^{\prime\prime}                             \big)       \equiv         \delta^{i^{\prime}}_5  \equiv  \underset{\emptyset \neq \mathscr{F}^{\prime\prime}_i}{\mathrm{inf}} \bigg\{  \mathcal{P}^{\chi_2}_{\Lambda} \bigg[              \big(   \gamma_{-,+}  \big)_2                       \underset{\Lambda \cap (\mathscr{F}\mathscr{C})_1}{\overset{+}{\longleftrightarrow}}    \mathscr{F}^{\prime\prime}_i          \big|         \mathcal{C}_{-} \cup    \mathscr{D}_{\chi_3 , -}              \bigg]          \bigg\}    \leq    \mathcal{P}^{\chi_2}_{\Lambda} \bigg[              \big(   \gamma_{-,+}  \big)_2                       \underset{\Lambda \cap (\mathscr{F}\mathscr{C})_1}{\overset{+}{\longleftrightarrow}}    \mathscr{F}^{\prime\prime}_i          \big|         \mathcal{C}_{-} \cup    \mathscr{D}_{\chi_3 , -}              \bigg]                                         \text{ } \text{ , }      \tag{Event $16$}
\end{align*}

\noindent due to the fact that,

\begin{align*}
     \mathcal{P}^{\chi_2}_{\Lambda} \bigg[              \big(   \gamma_{-,+}  \big)_2                       \underset{\Lambda \cap (\mathscr{F}\mathscr{C})_1}{\overset{+}{\longleftrightarrow}}    \mathscr{F}^{\prime\prime}_i          \big|         \mathcal{C}_{-} \cup    \mathscr{D}_{\chi_3 , -}              \bigg]           \subsetneq \underset{\mathscr{F}^{\prime\prime}_i}{\bigcup} \text{ }  \mathcal{P}^{\chi_2}_{\Lambda}  \bigg[              \big(   \gamma_{-,+}  \big)_2                       \underset{\Lambda \cap (\mathscr{F}\mathscr{C})_1}{\overset{+}{\longleftrightarrow}}    \mathscr{F}^{\prime\prime}_i          \big|         \mathcal{C}_{-} \cup    \mathscr{D}_{\chi_3 , -}              \bigg]          \text{ } \text{ , } 
\end{align*}

\noindent taken in the product over $i \in \textbf{N}$, $\delta_5 \big( \mathscr{F}^{\prime\prime} \big)$ is upper bounded with,

\begin{align*}
        \delta_5 \big( \mathscr{F}^{\prime\prime} \big)                        \text{ } \text{ , }
\end{align*}

\noindent and also with,

\begin{align*}
        \delta^{\prime}_5 \big( \mathscr{F}^{\prime\prime} \big)             \text{ } \text{ . } 
\end{align*}

\noindent Proceeding from $(\textit{****})$, in the lower bound provided in $(\textit{***})$, there exists a suitable constant $\Delta_5$ for which,

\begin{align*}
   \big(  \Delta_5 \big)^{|\mathscr{F}^{\prime\prime}|} \leq \delta_5 \big( \mathscr{F}^{\prime\prime} \big) \delta^{\prime}_5 \big( \mathscr{F}^{\prime\prime} \big)            \text{ } \text{ , } 
\end{align*}

\noindent and finally, in (\textit{*}), the existence of a suitable constant for which,

\begin{align*}
   \big(  \Delta_5 \big)^{|\mathscr{F}^{\prime\prime}|} \equiv  \Delta^{|\mathscr{F}^{\prime\prime}|}_5 \geq        \big( \Delta_5 \big)^{N^{\prime}}   \equiv \text{ }           \Delta^{N^{\prime}}_5 \text{ } \text{ . } 
\end{align*}

\noindent Altogether, incorporating the estimates for $\mathcal{P} \textit{1}$, $\mathcal{P} \textit{2}$, and $\mathcal{P}\textit{3}$, yields the cumulative lower bound,

\begin{align*}
     (\Delta ) \equiv \big|      \underset{\mathcal{P} \textit{1}}{\underbrace{k \delta^{\prime\prime}_1   \delta^{\prime\prime}_2 }}                   \big| +     \big|                     \underset{\mathcal{P} \textit{2}}{\underbrace{ \big( \Delta_4 \big)^N   }}          \big| +  \big|       \underset{\mathcal{P} \textit{3}}{\underbrace{ \big( \Delta_5 \big)^{N^{\prime}} }}   \big|  \geq \big|  k \delta^{\prime\prime}_1   \delta^{\prime\prime}_2                   \big| + \big|    \big( \Delta_4 \big)^{N} + \big( \Delta_5 \big)^{N^{\prime}}      \big|  \equiv   \big| k \delta^{\prime\prime}_1 \delta^{\prime\prime}_2    \big| + \big|  (\Delta_5)^{N^{\prime}}     \big(       \frac{(\Delta_4)^N}{( \Delta_5)^{N^{\prime}}}     +     1  \big)            \big| \\ \equiv     \big| k \delta^{\prime\prime}_1 \delta^{\prime\prime}_2    \big| + \big|  (\Delta_5)^{N^{\prime}}     \big| \big|     \frac{(\Delta_4)^N}{( \Delta_5)^{N^{\prime}}}     +     1             \big|          \\ \geq     \big| k \delta^{\prime\prime}_1 \delta^{\prime\prime}_2    \big| + \big|  \underset{N^{\prime}}{\mathrm{inf}}      \big( \Delta_5 \big)^{N^{\prime}}                 \big|   \text{ }       \big|         \frac{(\Delta_4)^N}{( \Delta_5)^{N^{\prime}}}     +     1               \big|       \\ \geq     \big| k \delta^{\prime\prime}_1 \delta^{\prime\prime}_2    \big| + \big|  \underset{N^{\prime}}{\mathrm{inf}}      \big( \Delta_5 \big)^{N^{\prime}}                 \big|   \text{ }       \big|    \underset{N,N^{\prime}}{\mathrm{inf}}   \big( \frac{(\Delta_4)^N}{( \Delta_5)^{N^{\prime}}}     +     1      \big)                  \big|  \\ \equiv     k \big| \delta^{\prime\prime}_1 \delta^{\prime\prime}_2    \big| + \big|  \underset{N^{\prime}}{\mathrm{inf}}      \big( \Delta_5 \big)^{N^{\prime}}                 \big|   \text{ }       \big|    \underset{N,N^{\prime}}{\mathrm{inf}}   \big( \frac{(\Delta_4)^N}{( \Delta_5)^{N^{\prime}}}     +     1      \big)                  \big| \text{ , } \tag{\textit{1}}  \end{align*}
     
     \noindent which is further rearranged with, 
     
     \begin{align*}
   (\textit{1})   \overset{(\textit{*})}{\geq} k \text{ } \underset{\delta^{\prime\prime}_1 , \delta^{\prime\prime}_2 , N^{\prime},k}{{\mathrm{inf}}} \big\{   \big| \delta^{\prime\prime}_1 \delta^{\prime\prime}_2 \big| ,  k^{-1} \big|    \underset{N^{\prime}}{\mathrm{inf}}     \big( \Delta_5 \big)^{N^{\prime}}       \big|          \big\}           \big( 1 +  \text{ }        \big|    \underset{N,N^{\prime}}{\mathrm{inf}}   \big( \frac{(\Delta_4)^N}{( \Delta_5)^{N^{\prime}}}     +     1      \big)                  \big| \text{ }  \big)  \\ \overset{\textit{(**)}}{\geq} K   \big( 1 +  \text{ }       \big|    \underset{N,N^{\prime}}{\mathrm{inf}}   \big( \frac{(\Delta_4)^N}{( \Delta_5)^{N^{\prime}}}     +     1      \big)                  \big| \text{ }  \big) \\ \geq   K^{\prime} + K^{\prime}  \big|    \underset{N,N^{\prime}}{\mathrm{inf}}   \big( \frac{(\Delta_4)^N}{( \Delta_5)^{N^{\prime}}}     +     1      \big)                  \big| \\ \overset{(\textit{***})}{\geq} K^{\prime} +  K^{\prime}  \big|    \epsilon_{N,N^{\prime}}    +     1                   \big|    \\                    \overset{(\textit{****})}{\geq}    K^{\prime} + K^{\prime} \epsilon_{N,N^{\prime},1}        \\ \geq K^{\prime}_{N,N^{\prime},1} \big( \epsilon \big) \equiv K^{\prime}_{N,N^{\prime},1} \text{ } \text{ , } 
\end{align*}

\noindent where, in the sequence of rearrangements above, beginning with the lower bound provided in (\textit{*}), there exists suitable $K$ for which,

\begin{align*}
      k \text{ } \underset{\delta^{\prime\prime}_1 , \delta^{\prime\prime}_2 , N^{\prime},k}{{\mathrm{inf}}} \big\{   \big| \delta^{\prime\prime}_1 \delta^{\prime\prime}_2 \big| ,  k^{-1} \big|    \underset{N^{\prime}}{\mathrm{inf}}      \big( \Delta_5 \big)^{N^{\prime}}       \big|          \big\}   \geq K      \text{ } \text{ , } 
\end{align*}

\noindent hence implying that the lower bound following (\textit{**}) holds. Proceeding, for (\textit{***}), there exists suitable $K^{\prime}$ so that,

\begin{align*}
    K \geq K^{\prime}          \text{ } \text{ , } 
\end{align*}

\noindent hence implying that the lower bound following $(\textit{***})$ holds, while finally, for $(\textit{****})$, observe, that there exists a suitable constant for which,

\begin{align*}
   \big| \epsilon_{N,N^{\prime}} + 1 \big| \geq \epsilon_{N,N^{\prime},1}            \text{ } \text{ , } 
\end{align*}

\noindent hence providing the desired lower bound upon setting $K^{\prime}_{N,N^{\prime},1} \equiv C^{\chi_1, \chi_2}\big( \Lambda \big) \equiv  C^{\chi_1, \chi_2}$, from which we conclude the argument. \boxed{}

\subsection{Results for the spin-representation measure $\mathcal{P}$ over the strip}

\noindent In the following, we introduce additional objects that are counterparts to the argument executed in previous sections for RSW results in the strip. Namely, introduce the vertical crossing event, which can consist of $+$ and $-$ faces, $\mathscr{F}_1 , \cdots , \mathscr{F}_n$ and finite strip volume FV, with $\mathcal{V}^{\mathrm{Mixed}}_{[-m,m] \times [0,n^{\prime}N]} \big( \mathscr{F}_1 , \cdots , \mathscr{F}_n , FV \big) \equiv \mathcal{V}^{+ \backslash -}_{[-m,m] \times [0,n^{\prime}N]} \big( \mathscr{F}_1 , \cdots , \mathscr{F}_n , FV \big) \equiv \mathcal{V}^{+ \backslash -}_{\Lambda} \big( \mathscr{F}_1 , \cdots , \mathscr{F}_n , FV \big)  \equiv \mathcal{V}^{+ \backslash -}_{\Lambda} \big(  \mathscr{F}_1 , \cdots , \mathscr{F}_n , FV   \big) \equiv \mathcal{V}^{+ \backslash -}_{\Lambda}$. With this quantity, we turn to arguments presented for lower bounding the probability of connectivity between $\mathcal{I}_j$ and $\widetilde{\mathcal{I}_j}$ under $\textbf{P}^{\xi^{\mathrm{Sloped}}}_{\Lambda} [ \cdot ]$ as considered throughout Section \textit{2}, instead for the spin-representation law $\mathcal{P}^{\xi}_{\Lambda}[ \cdot ]$, with the following.

\bigskip

\noindent \textbf{Lemma} \textit{6.2} (\textit{analog of Lemma 2.1 for possible arrangements of blocking freezing cluster interfaces in the strip}). WLOG denote $F^{\mathrm{even}}\big( \Lambda \big) \equiv F \big( \Lambda \big)$. Define $\mathscr{F}\mathscr{C} \equiv \underset{\text{countably many } u}{\bigcup} \mathscr{F}\mathscr{C}_u \in \big( F ( \Lambda) \cap \mathcal{F}\mathcal{C}_{\mathrm{AT}} \big)$, for $\mathscr{F}\mathscr{C} \cap \mathcal{F}\mathcal{C}^{\mathrm{even}}_{\mathrm{AT}} \neq \emptyset$, in addition to the following vertical crossing probability, 

\begin{align*}
 \mathcal{V}^{\mathrm{even}, + \backslash -}_{\Lambda} \bigg( \underset{u}{\bigcup} \text{ }  \mathscr{F}\mathscr{C}_u  \text{ } , FV  \bigg) \equiv  \mathcal{V}^{+ \backslash -}_{\Lambda} \bigg( \underset{u}{\bigcup} \text{ }  \mathscr{F}\mathscr{C}_u  \text{ } , FV  \bigg) \equiv  \mathcal{V}^{+ \backslash -}_{\Lambda} \big( \mathscr{F}\mathscr{C}_1 , \cdots , FV \big) \equiv \mathcal{V}^{+ \backslash -}_{\Lambda}   \text{ } \text{ , }
\end{align*}

\noindent where $\mathcal{V}^{+ \backslash -} \cap F \neq \emptyset$, across the strip admits the following lower bound,

\begin{align*}
C^{\emptyset}_{\mathcal{V}} \big( \Lambda \big) \equiv  C^{\emptyset}_{\mathcal{V}} \geq  \mathcal{P}^{\chi}_{\mathcal{D}_{\mathrm{AT}}} [  \mathcal{V}^{\mathrm{Mixed}}_{\Lambda}     ] \equiv \mathcal{P}^{\chi}_{\mathcal{D}_{\mathrm{AT}}} [  \mathcal{V}^{+ \backslash -}_{\Lambda}     ] \equiv  \mathcal{P}^{\chi}_{\mathcal{D}} [       \mathcal{V}^{+ \backslash -}_{\Lambda}    ]  \geq  C^{1}_{\mathcal{V}} \text{ } \text{ , } 
\end{align*}

\noindent for non-intersecting left and right boundaries of $\mathcal{D}$, in addition to the vertical crossing probability lower bound,

\begin{align*}
  C_{\mathcal{V}} \big( \Lambda \big) \equiv  C_{\mathcal{V}}   \geq \mathcal{P}^{\chi}_{\mathcal{D}} [       \mathcal{V}^{+ \backslash -}_{\Lambda}    ]   \text{ } \text{ , } 
\end{align*}

\noindent for intersecting left and right boundaries of $\mathcal{D}$, and boundary conditions $\chi \sim \textbf{B}\textbf{C}_{++}$.

\bigskip
   
\noindent \textit{Proof of Lemma 6.2}. For the first case, to obtain the desired lower bound $\mathcal{C}^{\emptyset}_{\mathcal{V}}$, write,

\begin{align*}
  \mathcal{P}^{\chi}_{\Lambda} [    \mathcal{V}^{+ \backslash -}_{\Lambda}    ] {\equiv}  \mathcal{P}^{\chi}_{\Lambda} \bigg[  \overset{| \mathcal{I}^{+ \backslash -}| }{\underset{j=1}{\bigcup}}         \mathcal{V}^{+ \backslash -}_j             \bigg]          \overset{(\mathrm{FKG})}{\geq}    \prod_{j=1}^{| \mathcal{I}^{+ \backslash -}|}   \mathcal{P}^{\chi}_{\Lambda} \big[         \mathcal{V}^{+ \backslash -}_j        \big]    \text{ } \text{ , }   \end{align*}

  \noindent for $\mathcal{I}^{\pm}_j$ such that,

\begin{align*}
   \mathcal{I}^{+ \backslash -}_j \equiv \{  j :   \mathcal{V}^{+ \backslash -}_j \cap \mathcal{V}^{+ \backslash -}_{\Lambda} \neq \emptyset    \} \text{ } \text{ , } \\     \big| \mathcal{I}^{+ \backslash - }_j \big| < + \infty        \text{ } \text{ . } 
\end{align*}

 \noindent where the final term after the application of $(\mathrm{FKG})$ is equivalent to the following intersection of crossing events, for $\gamma_L \big( \Lambda \big) \equiv \gamma_L$, and $\gamma_R \big( \Lambda \big) \equiv \gamma_R$,
  
  \begin{align*}
 \underset{R_1 \cap \gamma_R \neq \emptyset}{ \underset{L_1 \cap \gamma_L \neq \emptyset}{\underset{1\leq j \leq {\big| \mathcal{I}^{+ \backslash -}\big|}}{\prod}}}  \mathcal{P}^{\chi}_{\Lambda} \big[   \big\{ \gamma_{L_j} \overset{+ \backslash -}{\longleftrightarrow} \mathscr{F}\mathscr{C}_1       \big\}  \cap    \big\{   \mathscr{F}\mathscr{C}_2 \overset{+\backslash -}{\not\longleftrightarrow} \mathscr{F}\mathscr{C}_{N-1}         \big\}              \cap       \big\{    \mathscr{F}\mathscr{C}_N              \overset{+\backslash -}{\longleftrightarrow}     \gamma_{R_j}          \big\}            \big]    \overset{(\mathrm{FKG})}{\geq}      {\prod_{j=1}^{ \big| \mathcal{I}^{+ \backslash -}\big|}}   \mathcal{P}^{\chi}_{\Lambda} \big[   \gamma_L \overset{+ \backslash -}{\longleftrightarrow} \mathscr{F}\mathscr{C}_1      \big]  \tag{Event $17$} \\  \times  \mathcal{P}^{\chi}_{\Lambda} \big[    \mathscr{F}\mathscr{C}_2 \overset{+\backslash -}{\not\longleftrightarrow} \mathscr{F}\mathscr{C}_{N-1}      \big]   \mathcal{P}^{\chi}_{\Lambda} \big[    \mathscr{F}\mathscr{C}_N              \overset{+\backslash -}{\longleftrightarrow}     \gamma_R       \big]                 \text{ }                           \tag{\textit{AT 2.6.1}} \text{ , }    
\end{align*}

\noindent where, as demonstrated in previous arguments for six-vertex \textit{freezing clusters} in the strip, the middle disconnectivity event between \textit{Ashkin-Teller freezing clusters} over the finite strip $\Lambda$,

\begin{align*}
    \mathcal{P}^{\chi}_{\Lambda} \big[    \mathscr{F}\mathscr{C}_2 \overset{+\backslash -}{\not\longleftrightarrow} \mathscr{F}\mathscr{C}_{N-1}      \big]            \text{ } \text{ , } 
\end{align*}

\noindent can be lower bounded with, by $(\mathrm{FKG})$,

\begin{align*}
       \underset{2 \leq i \leq i+1 \leq N-1}{\prod}  \mathcal{P}^{\chi}_{\Lambda} \big[    \mathscr{F}\mathscr{C}_i \overset{+\backslash -}{\not\longleftrightarrow} \mathscr{F}\mathscr{C}_{i+1}      \big]              \text{ } \text{ , } \tag{\textit{i - (i+1)}} \\ \tag{Event $18$}
\end{align*}

\noindent which itself can be lower bounded with, for $F_1 , \cdots , F_n \in F^{\mathrm{even}} ( \textbf{Z}^2 \backslash \mathscr{F} \mathscr{C})$, or, for $ F_1 , \cdots , F_n \in F^{\mathrm{odd}} ( \textbf{Z}^2 \backslash \mathscr{F} \mathscr{C})$,

\begin{align*}
    \mathcal{P}^{\chi}_{\Lambda} \big[    \mathscr{F}\mathscr{C}_i \overset{+\backslash -}{\not\longleftrightarrow} \mathscr{F}\mathscr{C}_{i+1}      \big]   \geq       \mathcal{P}^{\chi}_{\Lambda} \big[    \big\{   F_1         \underset{(\mathscr{F} \mathscr{C})^c}{\overset{+ \backslash -}{\longleftrightarrow}}  F_n  \big\}  \cap  \big\{ \mathscr{F}\mathscr{C}_i \overset{+\backslash -}{\not\longleftrightarrow}  \mathscr{F}\mathscr{C}_{i+1}   \big\}  \big]  \overset{(\mathrm{FKG})}{\geq}     \text{ } \underset{1 \leq i^{\prime} \leq i^{\prime}+1 < + \infty }{\prod}      \mathcal{P}^{\chi}_{\Lambda}  \big[         F_i  \underset{(\mathscr{F} \mathscr{C})^c}{\overset{+ \backslash -}{\longleftrightarrow}}      F_{i+1}          \big]    \\ \times \mathcal{P}^{\chi}_{\Lambda}     \big[  \mathscr{F}\mathscr{C}_i \overset{+\backslash -}{\not\longleftrightarrow}  \mathscr{F}\mathscr{C}_{i+1}  \big]     \text{ } \text{ , } \tag{\textit{i - (i+1) II}}
\end{align*}

\noindent implying,

\begin{align*}
  (\textit{i - (i+1)}) \overset{(\textit{i-(i+1) II})}{\geq}     \underset{2 \leq i \leq i+1 \leq N-1}{\prod}    \mathcal{P}^{\chi}_{\Lambda} \big[          \gamma_L \overset{+ \backslash -}{\longleftrightarrow} \mathscr{F}\mathscr{C}_1     \big] \text{ } \bigg[  \underset{1 \leq i^{\prime} \leq i^{\prime}+1 < + \infty }{\prod}  \text{ }      \mathcal{P}^{\chi}_{\Lambda}  \big[         F_i  \underset{(\mathscr{F} \mathscr{C})^c}{\overset{+ \backslash -}{\longleftrightarrow}}      F_{i+1}          \big]      \text{ } \mathcal{P}^{\chi}_{\Lambda}     \big[  \mathscr{F}\mathscr{C}_i \overset{+\backslash -}{\not\longleftrightarrow}  \mathscr{F}\mathscr{C}_{i+1}  \big] \bigg]   \\ \times   \mathcal{P}^{\chi}_{\Lambda} \big[        \mathscr{F}\mathscr{C}_N              \overset{+\backslash -}{\longleftrightarrow}     \gamma_R     \big]    \text{ } \text{ , }
\end{align*}

\noindent and hence that a previously obtained product of probabilities can be lower bound with,

\begin{align*}
 (\textit{AT 2.6.1}) \geq   \underset{1 \leq i^{\prime} \leq i^{\prime}+1 < + \infty}{\underset{2 \leq i \leq i+1 \leq N-1}{\underset{1 \leq j \leq | \mathcal{I}^{+ \backslash -}|}{\prod}}}   \mathcal{P}^{\chi}_{\Lambda} \big[     \gamma_L \overset{+ \backslash -}{\longleftrightarrow} \mathscr{F}\mathscr{C}_1      \big]          \mathcal{P}^{\chi}_{\Lambda}  \big[         F_i  \underset{(\mathscr{F} \mathscr{C})^c}{\overset{+ \backslash -}{\longleftrightarrow}}      F_{i+1}          \big]    \mathcal{P}^{\chi}_{\Lambda}     \big[  \mathscr{F}\mathscr{C}_i \overset{+\backslash -}{\not\longleftrightarrow}  \mathscr{F}\mathscr{C}_{i+1}  \big] \mathcal{P}^{\chi}_{\Lambda} \big[        \mathscr{F}\mathscr{C}_N              \overset{+\backslash -}{\longleftrightarrow}     \gamma_R     \big]                            \text{ } \text{ , }      \tag{\textit{AT 2.6.1 II}}
\end{align*}

\noindent Next, continue to execute arguments along the lines of \textbf{Lemma} \textit{2.6.1.1}, \textbf{Lemma} \textit{2.6.1.2}, and \textbf{Lemma} \textit{2.6.1.3}, with the following arguments.

\bigskip

\noindent \textbf{Lemma} \textit{AT 2.6.1.1} (\textit{a lower bound for the first crossing probability of (\textit{AT 2.6.1})}). For the line $\mathscr{L}_1$ appearing before the first \textit{Ashkin-Teller freezing cluster} in the finite-volume strip, one has the strictly positive lower bound,

\begin{align*}
   \mathcal{P}^{\chi}_{\Lambda} \big[   \gamma_L \underset{F^{\mathrm{even}}}{\overset{+ \backslash -}{\longleftrightarrow}} \mathscr{F}\mathscr{C}_1      \big]  \equiv  \mathcal{P}^{\chi}_{\Lambda} \big[   \gamma_L \overset{+ \backslash -}{\longleftrightarrow} \mathscr{F}\mathscr{C}_1      \big]  \geq \mathscr{C}^{\mathrm{AT}}_1 \equiv \mathscr{C}_1          \text{ } \text{ , } 
\end{align*}

\noindent for $\chi \in \textbf{B}\textbf{C}_{++}$.

\bigskip

\noindent \textit{Proof of Lemma AT 2.6.1.1}. The probability that the connectivity event between the line $\mathscr{L}_1$ and $\mathscr{F}\mathscr{C}_1$ will be analyzed through the following sequence of inequalities. Express the crossing event above with the decomposition,

\begin{align*}
    \mathcal{P}^{\chi}_{\Lambda}  [  \gamma_L \overset{+ \backslash -}{\longleftrightarrow} \mathscr{F}\mathscr{C}_1  ]        \overset{ \underset{L_i}{\cup} \{      \gamma_L \overset{+ \backslash -}{\longleftrightarrow}   \gamma_{L_i}       \}    \subsetneq    \{    \gamma_L \overset{+ \backslash -}{\longleftrightarrow} \mathscr{F}\mathscr{C}_1        \}         }{\equiv}    \mathcal{P}^{\chi}_{\Lambda}  \big[    \text{ }   \underset{L_i}{\bigcap}  \big\{     \gamma_L \overset{+ \backslash -}{\longleftrightarrow}                   \gamma_{L_i}    \big\} \text{ }  \big] \overset{(\mathrm{FKG})}{\geq}      \underset{L_i}{\prod}         \mathcal{P}^{\chi}_{\Lambda}  \big[      \gamma_L \overset{+ \backslash -}{\longleftrightarrow}                   \gamma_{L_i}  \big]        \text{ } \text{ . }          \tag{Event $19$}
\end{align*}

\noindent For further rearrangements below, introduce,

\begin{align*}
   \mathcal{L}_i \equiv  \big| L_i \big| \equiv \big|  \big\{ L_i : \mathcal{P}^{\chi}_{\Lambda}    \big[      \gamma_L \overset{+ \backslash -}{\longleftrightarrow}                   \gamma_{L_i}  \big] > 0   \big\} \big| < + \infty   \text{ } \text{ . } 
\end{align*}    

      \begin{figure}
\begin{align*}
\includegraphics[width=0.86\columnwidth]{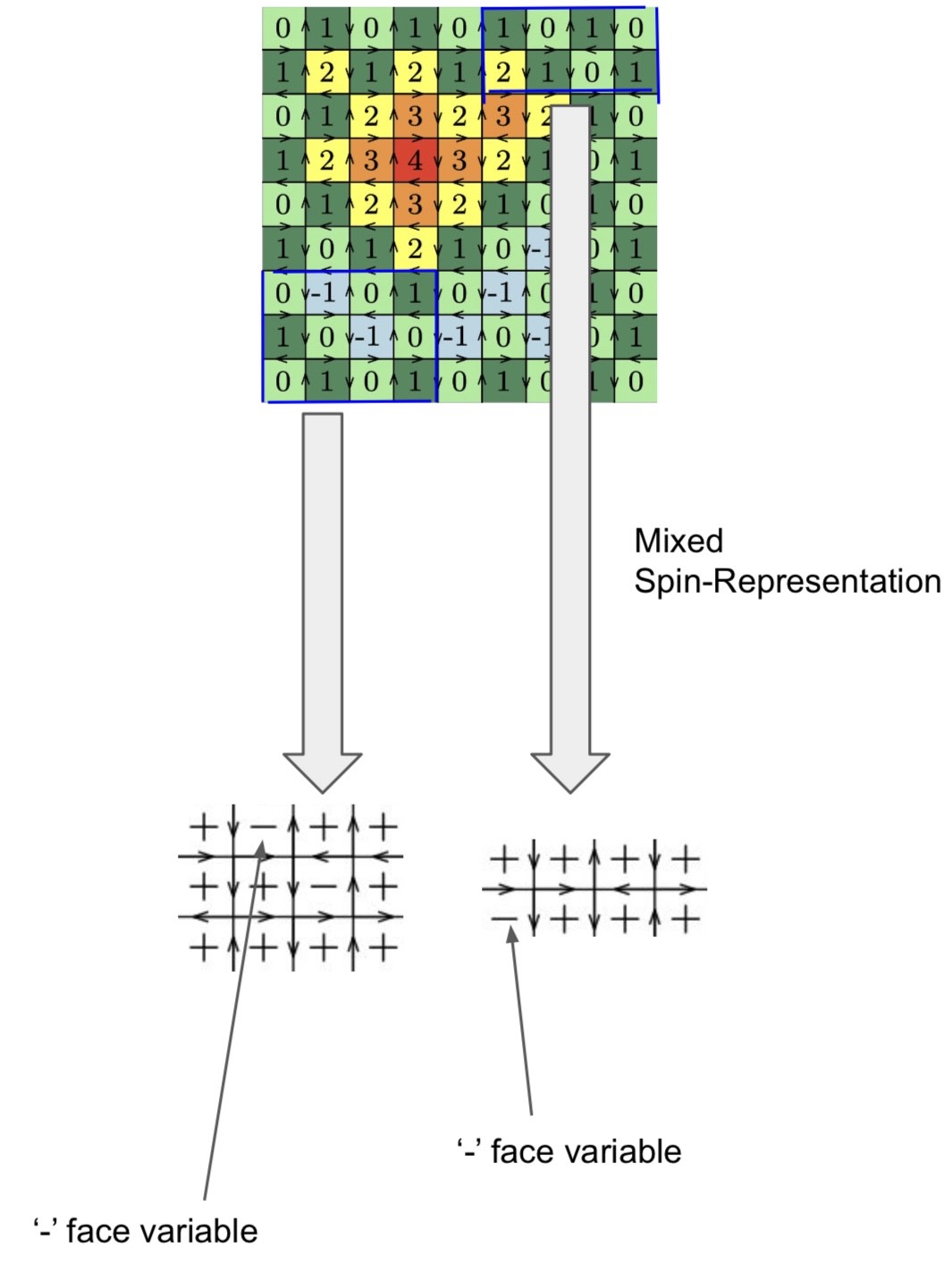}
\end{align*}
\caption{\textit{Other finite volumes within the example provided by Glazman and Peled for mixed Ashkin-Teller spin representations}. From \textit{Figure 32}, the restriction of the faces to the two finite volumes given above, each of which are segmented in blue, also exhibit that $-$, or $+$ face variables together contribute to crossing events pushed forwards under the Ashkin-Teller probability measure $\mathcal{P}_{\Lambda} [\cdot]$.}
\end{figure}

\noindent Therefore, one can deduce the proportionality, 

\begin{align*}
    \mathcal{P}^{\chi}_{\Lambda}  \big[      \gamma_L \overset{+ \backslash -}{\longleftrightarrow}                   \gamma_{L_i}  \big]    \overset{(\textit{*})}{\propto}    \mathcal{P}^{\chi}_{\Lambda}  \big[ \forall  \mathscr{D}_1 , \mathscr{D}_2 \subsetneq \Lambda , \exists \chi_1 , \chi_2 \in \textbf{B}\textbf{C}_{++}:  \mathcal{P}^{\chi_1}_{\mathscr{D}_1 } \big[ \big(  \gamma_{L} \cap L \big( \mathscr{D}_1 \big)  \big) \overset{+ \backslash -}{\longleftrightarrow }  \big( \gamma_{L_i} \cap R \big( \mathscr{D}_1 \big) \big)      \big] \\ ,         \mathcal{P}^{\chi_2}_{\mathscr{D}_2 } \big[ \big( \gamma_{L}   \cap   L \big( \mathscr{D}_2 \big)    \big) \overset{+ \backslash -}{\longleftrightarrow }  \big( \gamma_{L_i} \cap  R \big( \mathscr{D}_2 \big)    \big) \big]  > 0  \big]              \overset{(\textit{**})}{ \geq}          C \big( \big|  \frac{\Lambda \cap L_i \cap L }{\Lambda}   \big|      \big) \end{align*}

\noindent where $L \big( \mathscr{D}_1 \big)$ and $L \big( \mathscr{D}_2 \big)$ respectively denote the left sides of the domains $\mathscr{D}_1$ and $\mathscr{D}_2$, and $L \big( \mathscr{D}_1 \big)$ and $L \big( \mathscr{D}_2 \big)$ respectively denote the right sides of the domains $\mathscr{D}_1$ and $\mathscr{D}_2$, in which the statement captures the probability of two crossings occurring, as long as there exists $++$ boundary conditions for which both of the crossing probabilities between the left and right sides of the symmetric domain are strictly positive. The constant $C$ appearing in the lower bound above can further be lower bound from a product of crossing probabilities between $\gamma_L$ and $\gamma_{L_i}$, in which,

    \begin{align*}
   C \big( \big|  \frac{\Lambda \cap L_i \cap L }{\Lambda}   \big|      \big)  \equiv   \mathscr{C}^{\frac{1}{\mathcal{L}_i}}_1                             \Longleftrightarrow \underset{L_i}{\prod}         \mathcal{P}^{\chi}_{\Lambda}  \big[      \gamma_L \overset{+ \backslash -}{\longleftrightarrow}                   \gamma_{L_i}  \big]  \geq  \underset{\mathrm{lines \text{ }} L_i}{\prod}   \sqrt[{\frac{1}{\mathcal{L}_i}}]{\mathscr{C}_1}               \geq   \mathscr{C}_1 \big( \mathcal{L}_i \big) \equiv \mathscr{C}_1  \text{ } \text{ , } 
\end{align*}

\noindent as desired, where in the proportionality in $(\textit{*})$, with positive probability a mixed path consisting of $+$ and $-$ face variables occurs if the faces traversed in the crossing event $\big\{ \gamma_L    \overset{+ \backslash -}{\longleftrightarrow}     \gamma_{L_i} \big\}$, which is expressed through the product of indicators under stipulations $\mathcal{C}_1$ and $\mathcal{C}_2$, where $C$ is dependent upon the total number of faces bound by the finite volume,

\begin{align*}
    \big|     \frac{\Lambda \cap L_i \cap L }{\Lambda}   \big|      \text{ } \text{ , } 
\end{align*}

\noindent that is given by,

\begin{align*}
  \mathcal{P}^{\chi}_{\Lambda} \big[ \Lambda \subsetneq  \textbf{Z}^2 :   \big| F_{+} \big| \geq  \big|      F ( \Lambda \cap L_i \cap L )           \big|                        \big]  \equiv \mathcal{P}^{\chi}_{\Lambda} \bigg[    \Lambda \subsetneq  \textbf{Z}^2 :  \underset{\Lambda^{\prime} \subsetneq \Lambda}{\sum} \textbf{1}_{\{ F_{+} ( \Lambda^{\prime} ) \geq F ( \Lambda^{\prime} \cap L_i \cap L )  \}}  \approx \big| \Lambda^{\prime} \big|  \bigg] \\ \propto  \underset{\Lambda^{\prime} \subsetneq \Lambda}{\sum}  \mathcal{P}^{\chi}_{\Lambda} \big[    \textbf{1}_{\{ F_{+} ( \Lambda^{\prime} ) \geq F ( \Lambda^{\prime} \cap L_i \cap L )  \}}  \approx \big| \Lambda^{\prime} \big|^{- | \Lambda^{\prime} \cap \Lambda |}                    \big]  \\ \propto   C \big(  \big(   F \big( \Lambda \cap L_i \cap L \big)     \big)^{-1}  \big)  C \big( \big|  \Lambda^{\prime} \big|  ,  \big| \Lambda  \big| ,   \big| F \big( \Lambda \big) \big| \big)     \text{ , } \tag{$\mathscr{F}_1$-$\mathscr{F}_2$ \textit{bound}}
\end{align*}

\noindent for some $\alpha >0$, in which the disjoint union of the faces with $+$ variables, with the faces with $-$ variables is equal to the total number of faces over $\Lambda \cap L_i \cap L$, as,

\begin{align*}
     \big| F_{+} \big| \text{ } \dot{\cup} \text{ } \big| F_{-} \big| \equiv  \text{ } \big| F \big( \Lambda \cap L_i \cap L \big)  \big| < \big| F ( \textbf{Z}^2 ) \big|  \text{ } \text{ , }
\end{align*}

\noindent in addition to the total number of faces between $L$ and $L_i$,

\begin{align*}
               F \big( \Lambda \cap L_i \cap L \big)           \subset  F ( \Lambda)        \subset F( \textbf{Z}^2 ) \text{ } \text{ , } 
\end{align*}

\noindent and also that,

\begin{align*}
\underset{\Lambda^{\prime} \subsetneq \Lambda}{\sum} \frac{1}{\big| \Lambda^{\prime} \big|^{| \Lambda^{\prime} \cap \Lambda |}}  =    \frac{1}{\big| \Lambda^{\prime} \big|^{| \Lambda^{\prime} \cap \Lambda |}}  + \underset{| \Lambda^{\prime} \cap \Lambda | - 2 \text{ } \mathrm{terms}}{\cdots} +  \frac{1}{\big| \Lambda^{\prime} \big|^{| \Lambda^{\prime} \cap \Lambda |}}  \approx  \frac{| \Lambda^{\prime} \cap \Lambda |}{|\Lambda^{\prime} |^{| \Lambda^{\prime} \cap \Lambda |}}     \approx  \frac{| \Lambda^{\prime} |}{|\Lambda^{\prime} |^{| \Lambda^{\prime} \cap \Lambda |}} \approx   | \Lambda^{\prime} |^{1 - \epsilon}  \approx | \Lambda^{\prime} |  \text{ } \text{ , } 
\end{align*}

\noindent for $\epsilon > 0$ sufficiently small, with,

\begin{align*}
     \epsilon < | \Lambda^{\prime} \cap \Lambda |   \text{ } \text{ . } 
\end{align*}

\noindent The lower bound following $(\textit{**})$ results from the fact that,

\begin{align*}
      (\mathscr{F}_1 - \mathscr{F}_2 \text{ } \textit{bound}) \equiv     C \big(    \big( F \big( \Lambda \big)  F \big( \Lambda \cap L_i \cap L \big)     \big)^{-1} \big)      C \big( \big|  \Lambda^{\prime} \big|  ,  \big| \Lambda  \big| ,   \big| F \big( \Lambda \big) \big| \big)     \propto     \mathcal{P}^{\chi}_{\Lambda} \bigg[ \forall \Lambda \subsetneq \textbf{Z}^2 ,  \exists \mathscr{D}_1 , \mathscr{D}_2  \subsetneq F ( \textbf{Z}^2 ) \\  , \chi_1 , \chi_2 \in \textbf{B}\textbf{C}^{++}:    \mathcal{P}^{\chi_1}_{\mathscr{D}_1} \big[                 \big(  \gamma_{L_i} \cap L \big( \mathscr{D}_1 \big)   \big)   \overset{+ \backslash -}{\longleftrightarrow}    \big( \gamma_L \cap R \big( \mathscr{D}_1 \big) \big)  \big] \\  ,  \mathcal{P}^{\chi_2}_{\mathscr{D}_2} \big[     \big( \gamma_{L_i} \cap L \big( \mathscr{D}_2 \big) \big)    \overset{+ \backslash -}{\longleftrightarrow}       \big( \gamma_L \cap R \big( \mathscr{D}_2 \big) \big)    \big] \\ ,   \underset{\Lambda^{\prime} \subsetneq \Lambda}{\sum} \textbf{1}_{\{ F_{+} ( \Lambda^{\prime} ) \geq F ( \Lambda^{\prime} \cap L_i \cap L )  \}}  \approx \big| \Lambda^{\prime} \big|   \bigg] \end{align*}
      
     \noindent which is proportional to the following product of probabilities,

   \begin{align*}    \mathcal{P}^{\chi}_{\Lambda} \big[       \mathcal{P}^{\chi_1}_{\mathscr{D}_1} \big[                 \big(  \gamma_{L_i} \cap L \big( \mathscr{D}_1 \big)   \big)   \overset{+ \backslash -}{\longleftrightarrow}    \big( \gamma_L \cap R \big( \mathscr{D}_1 \big) \big)  \big]  ,  \mathcal{P}^{\chi_2}_{\mathscr{D}_2} \big[     \big( \gamma_{L_i} \cap L \big( \mathscr{D}_2 \big) \big)    \overset{+ \backslash -}{\longleftrightarrow}       \big( \gamma_L \cap R \big( \mathscr{D}_2 \big) \big)    \big] \big]  \\  \times   \mathcal{P}^{\chi}_{\Lambda} \bigg[      \underset{\Lambda^{\prime} \subsetneq \Lambda}{\sum} \textbf{1}_{\{ F_{+} ( \Lambda^{\prime} ) \geq F ( \Lambda^{\prime} \cap L_i \cap L )  \}}  \approx \big| \Lambda^{\prime} \big|      \bigg]        \\  \overset{(\mathrm{FKG})}{\geq} \mathcal{P}^{\chi}_{\Lambda} \big[    \mathcal{P}^{\chi_1}_{\mathscr{D}_1} \big[                 \big(  \gamma_{L_i} \cap L \big( \mathscr{D}_1 \big)   \big)   \overset{+ \backslash -}{\longleftrightarrow}    \big( \gamma_L \cap R \big( \mathscr{D}_1 \big) \big)  \big]      \big] \mathcal{P}^{\chi}_{\Lambda} \big[  \mathcal{P}^{\chi_2}_{\mathscr{D}_2} \big[     \big( \gamma_{L_i} \cap L \big( \mathscr{D}_2 \big) \big)    \overset{+ \backslash -}{\longleftrightarrow}       \big( \gamma_L \cap R \big( \mathscr{D}_2 \big) \big)    \big] \big]  \\ \times \mathcal{P}^{\chi}_{\Lambda} \bigg[ \underset{\Lambda^{\prime} \subsetneq \Lambda}{\sum} \textbf{1}_{\{ F_{+} ( \Lambda^{\prime} ) \geq F ( \Lambda^{\prime} \cap L_i \cap L )  \}}  \approx \big| \Lambda^{\prime} \big|   \bigg]  \\  =   \mathcal{P}^{\chi}_{\Lambda} \big[    \mathcal{P}^{\chi_1}_{\mathscr{D}_1} \big[                 \big(  \gamma_{L_i} \cap L \big( \mathscr{D}_1 \big)   \big)   \overset{+ \backslash -}{\longleftrightarrow}    \big( \gamma_L \cap R \big( \mathscr{D}_1 \big) \big)  \big]      \big] \mathcal{P}^{\chi}_{\Lambda} \big[  \mathcal{P}^{\chi_2}_{\mathscr{D}_2} \big[     \big( \gamma_{L_i} \cap L \big( \mathscr{D}_2 \big) \big)    \overset{+ \backslash -}{\longleftrightarrow}       \big( \gamma_L \cap R \big( \mathscr{D}_2 \big) \big)    \big] \big]    \\ \times     \bigg[ \underset{\Lambda^{\prime} \subsetneq \Lambda}{\sum} \mathcal{P}^{\chi}_{\Lambda} \big[        \textbf{1}_{\{ F_{+} ( \Lambda^{\prime} ) \geq F ( \Lambda^{\prime} \cap L_i \cap L )  \}}  \approx \big| \Lambda^{\prime} \big|^{-|\Lambda \cap \Lambda^{\prime}|}             \big] \bigg]        \\  \propto    \bigg[  C \big( \big|  \frac{\Lambda \cap L_i \cap L }{\Lambda}   \big|      \big) \bigg] \bigg[      C \big(  \big(   F \big( \Lambda \cap L_i \cap L \big)     \big)^{-1}  \big)  C \big( \big|  \Lambda^{\prime} \big|  ,  \big| \Lambda  \big| ,   \big| F \big( \Lambda \big) \big| \big)          \bigg]     \text{ } \text{ , } 
\end{align*}

\noindent for suitable $C_1$ and $C_2$, which as with the previously provided constant $C$, are also dependent upon,

\begin{align*}
    \big|     \frac{\Lambda \cap L_i \cap L }{\Lambda}   \big|  \text{ } \text{ , } 
\end{align*}

\noindent from which we conclude the argument. \boxed{}

\bigskip

\noindent \textbf{Lemma} \textit{AT 2.6.1.2} (\textit{a lower bound for the product of crossing probabilities appearing in (\textit{i - (i+1)})}). WLOG suppose that $\mathscr{F}_i , \mathscr{F}_{i+1} \in F^{\mathrm{even}} \big( \textbf{Z}^2 \big) \equiv F \big( \textbf{Z}^2 \big), \forall \text{ } i >0$. For the family of disconnectivity probabilities, $\big\{  \mathscr{F}\mathscr{C}_i \overset{+\backslash -}{\not\longleftrightarrow}  \mathscr{F}\mathscr{C}_{i+1}     \big\}$, indexed in $i^{\prime}$ for $1 \leq i^{\prime} \leq i^{\prime}_1 < +\infty$, pushed forwards under $\mathcal{P}^{\chi}_{\Lambda} [ \cdot ]$ satisfies,

\begin{align*}
       \mathcal{P}^{\chi}_{\Lambda} \big[     \mathscr{F}\mathscr{C}_i \overset{+\backslash -}{\not\longleftrightarrow}  \mathscr{F}\mathscr{C}_{i+1}                    \big]   \geq          c^{\mathrm{AT}}_2 \equiv c_2                               \text{ } \text{ , } 
\end{align*}

\noindent for $\chi \in \textbf{B}\textbf{C}_{++}$, implying that the desired lower bound for the product over $i$, given above, takes the form,

\begin{align*}
   \underset{\text{ countably many } i^{\prime} }{ \prod} \mathcal{P}^{\chi}_{\Lambda} \big[      \mathscr{F}\mathscr{C}_i  \overset{+\backslash -}{\not\longleftrightarrow}  \gamma_{i^{\prime}}                \big]         \geq   \underset{\text{ countably many } i^{\prime} }{ \prod} \big(    c^{\mathrm{AT}}_2      \big)_{i^{\prime}}        \geq \mathscr{C}^{\mathrm{AT}}_2                           \text{ } \text{ , } 
\end{align*}

\noindent for paths of faces $\gamma_{i^{\prime}}$ occurring to the right of $\mathscr{F}\mathscr{C}_i$.

\bigskip

\noindent \textit{Proof of Lemma AT 2.6.1.2}. From the family of  disconnectivity events given in the statement, express,

\begin{align*}
     \mathcal{P}^{\chi}_{\Lambda} \big[     \mathscr{F}\mathscr{C}_i \overset{+\backslash -}{\not\longleftrightarrow}  \mathscr{F}\mathscr{C}_{i+1}                    \big]  \equiv   \mathcal{P}^{\chi}_{\Lambda} \big[   \underset{\text{ countably many } i^{\prime} }{\bigcap} \text{ } \big\{ \mathscr{F}\mathscr{C}_i   \overset{+\backslash -}{\not\longleftrightarrow}         \gamma_{i^{\prime}}           \big\}      \big]  \tag{Event $20$} \\ \geq      \underset{\text{ countably many } i^{\prime} }{\prod}     \mathcal{P}^{\chi}_{\Lambda} \big[      \mathscr{F}\mathscr{C}_i   \overset{+\backslash -}{\not\longleftrightarrow}         \gamma_{i^{\prime}}              \big]   \\     \overset{(\textit{*})}{\geq}    \underset{\text{ countably many } i^{\prime} }{\prod}         C_{i^{\prime}}        \geq \mathscr{C}^{\mathrm{AT}}_2      \text{ } \text{ , } 
\end{align*}

\noindent where, in $(\textit{*})$, apply,

\begin{align*}
      C_{i^{\prime}}         \leq \underset{i^{\prime}}{\mathrm{inf}} \big\{    \mathcal{P}^{\chi}_{\Lambda} \big[      \mathscr{F}\mathscr{C}_i   \overset{+\backslash -}{\not\longleftrightarrow}         \gamma_{i^{\prime}}              \big]         \big\}  \text{ } \text{ , } 
\end{align*}

\noindent for suitably chosen, strictly positive $C$, reflecting arguments for obtaining the lower bound in \textbf{Lemma} \textit{2.6.1.2}, from which we conclude the argument. \boxed{}

\bigskip

\noindent \textbf{Lemma} \textit{AT 2.6.1.3} (\textit{a lower bound for all other crossing probabilities appearing in the estimate for (AT 2.6.1)}). WLOG, suppose that the two crossing events,

\begin{align*}
  \bigg\{     F_i  \underset{(\mathscr{F} \mathscr{C})^c \cap F^{\mathrm{even}} }{\overset{+ \backslash -}{\longleftrightarrow}}      F_{i+1}     \bigg\} \equiv \big\{     F_i  \underset{(\mathscr{F} \mathscr{C})^c}{\overset{+ \backslash -}{\longleftrightarrow}}      F_{i+1}     \big\}  \text{ } \text{ , } 
\end{align*}

\noindent and,

\begin{align*}
 \bigg\{       \mathscr{F}\mathscr{C}_N             \underset{F^{\mathrm{even}}}{\overset{+\backslash -}{\longleftrightarrow}}     \gamma_R      \bigg\}   \equiv \big\{    \mathscr{F}\mathscr{C}_N              \overset{+\backslash -}{\longleftrightarrow}     \gamma_R        \big\}  \text{ } \text{ , } 
\end{align*}

\noindent each are only dependent upon faces that have nonempty intersection with $F^{\mathrm{even}} \big( \textbf{Z}^2 \big)$. The product of crossing random variables, 

\begin{align*}
    \mathcal{P}^{\chi}_{\Lambda}  \bigg[         F_i  \underset{(\mathscr{F} \mathscr{C})^c}{\overset{+ \backslash -}{\longleftrightarrow}}      F_{i+1}          \bigg]           \text{ }       \mathcal{P}^{\chi}_{\Lambda} \big[        \mathscr{F}\mathscr{C}_N              \overset{+\backslash -}{\longleftrightarrow}     \gamma_R     \big]                \text{ } \text{ , } 
\end{align*}

\noindent admits the lower bound,

\begin{align*}
       \mathscr{C}^{\mathrm{AT}}_3 \big( \mathscr{F}\mathscr{C}_N , \gamma_R \big) \equiv    \mathscr{C}^{\mathrm{AT}}_3        \text{ } \text{ , } 
\end{align*}

\noindent for $\mathcal{I}(i^{\prime\prime}, j^{\prime\prime}) \equiv \mathcal{I}$ with $i^{\prime\prime} \neq j^{\prime\prime}>0$, and for $\chi \in \textbf{B}\textbf{C}_{++}$.

\bigskip

\noindent \textit{Proof of Lemma AT 2.6.1.3}. For the final lower bound, express the lower bound from the rearrangements,

\begin{align*}
      \mathcal{P}^{\chi}_{\Lambda}  \big[         F_i  \underset{(\mathscr{F} \mathscr{C})^c}{\overset{+ \backslash -}{\longleftrightarrow}}      F_{i+1}          \big]    \geq     \mathcal{P}^{\chi}_{\Lambda}  \big[    \underset{F^{\prime}_j \in F ( \Lambda) }{\underset{\text{ countably many } j^{\prime}}{\bigcap}} \big\{      F_i  \underset{(\mathscr{F} \mathscr{C})^c}{\overset{+ \backslash -}{\longleftrightarrow}}    F^{\prime}_j              \big\}       \big]    \overset{(\mathrm{FKG})}{\geq}   \underset{F^{\prime}_j \in F ( \Lambda) }{\underset{\text{ countably many } j^{\prime}}{\prod}}      \mathcal{P}^{\chi}_{\Lambda}  \big[       F_i  \underset{(\mathscr{F} \mathscr{C})^c}{\overset{+ \backslash -}{\longleftrightarrow}}    F^{\prime}_j                   \big] \geq C^{\mathrm{AT}}_1        \text{ } \text{ , }  \\ \tag{Event $21$}
\end{align*}

\noindent corresponding to the first probability with connectivity over $\big( \mathscr{F} \mathscr{C} \big)^c \subset \Lambda$, and, for $\mathscr{L}_i \big( \gamma_R \big) \equiv \mathscr{L}_i$,

\begin{align*}
      \mathcal{P}^{\chi}_{\Lambda} \big[              \mathscr{F}\mathscr{C}_N \overset{+ \backslash - }{\longleftrightarrow}     \gamma_R                     \big]   \geq         \mathcal{P}^{\chi}_{\Lambda} \big[   \text{ }    \underset{\exists \text{ }  i > i^{\prime} \text{ } : \text{ } \mathscr{L}_i \cap  \gamma_R \neq \emptyset}{\underset{\text{countably many } i^{\prime}}{\bigcap}}  \text{ }  \big\{              \mathscr{F}\mathscr{C}_N \overset{+ \backslash - }{\longleftrightarrow}   \mathscr{L}_{i^{\prime}}  \big\}   \big]   \overset{(\mathrm{FKG})}{\geq}       {\underset{\text{countably many } i^{\prime}}{\prod}}    \mathcal{P}^{\chi}_{\Lambda} \big[     \mathscr{F}\mathscr{C}_N \overset{+ \backslash - }{\longleftrightarrow}   \mathscr{L}_{i^{\prime}} \big]   \geq C^{\mathrm{AT}}_2                              \text{ } \text{ , } \\ \tag{Event $22$} 
\end{align*}

\noindent corresponding to the second probability, from which we conclude the argument upon setting $\mathscr{C}^{\mathrm{AT}}_3 \big( i^{\prime} , j^{\prime} \big) \equiv \mathscr{C}^{\mathrm{AT}}_3 \leq \bigg[ C^{\mathrm{AT}}_1 \big( i^{\prime} \big) C^{\mathrm{AT}}_2 \big( j^{\prime} \big)       \bigg]^{|\mathcal{I}|} \equiv \bigg[ C^{\mathrm{AT}}_1  C^{\mathrm{AT}}_2       \bigg]^{|\mathcal{I}|}$. \boxed{}

\bigskip

\noindent \textit{Proof of Lemma 6.2}. Incorporating the previously obtained estimates simultaneously yields, for the lower bound on the vertical crossing probability over $F^{\mathrm{even}} \big( \textbf{Z}^2 \big)$, 

\begin{align*}
 (\textit{AT 2.6.1 II}) \overset{(\textbf{Lemma}\text{ }  \textit{AT 2.6.1.1}), (\textbf{Lemma}\text{ }  \textit{AT 2.6.1.2}), (\textbf{Lemma}\text{ }  \textit{AT 2.6.1.3})}{\geq}          \text{ } \underset{1 \leq i^{\prime} \leq i^{\prime}+1 < + \infty}{\underset{2 \leq i \leq i+1 \leq N-1}{\underset{1 \leq j \leq | \mathcal{I}^{+ \backslash -}|}{\prod}}} \big( \mathscr{C}^{\mathrm{AT}}_1 \big)_j \text{ }             \big(  \mathscr{C}^{\mathrm{AT}}_1 \big)_{i} \text{ } \big( \mathscr{C}^{\mathrm{AT}}_3 \big)_{i^{\prime}} \\ \equiv     \bigg[ \underset{1 \leq j \leq | \mathcal{I}^{+ \backslash -}|}{ \prod}    \big( \mathscr{C}^{\mathrm{AT}}_1 \big)_j \bigg] \bigg[   \underset{2 \leq i \leq i+1 \leq N-1}{ \prod}    \big(  \mathscr{C}^{\mathrm{AT}}_1 \big)_{i} \bigg]     \bigg[   \underset{1 \leq i^{\prime} \leq i^{\prime}+1 < + \infty}{ \prod}      \big( \mathscr{C}^{\mathrm{AT}}_3 \big)_{i^{\prime}}  \bigg]    \\  \geq \mathscr{C}^{\mathrm{AT}}_1 \text{ } \mathscr{C}^{\mathrm{AT}}_2 \text{ } \mathscr{C}^{\mathrm{AT}}_3 \\     \equiv   C^{1}_{\mathcal{V}}             \text{ } \text{ , } 
\end{align*}

\noindent where the product over $j$, $i$ and $i^{\prime}$ can be decomposed as a product over $\mathscr{C}_1$, $\mathscr{C}_2$, and $\mathscr{C}_3$, respectively. The lower bound above concludes the argument for non-intersecting left and right boundaries of the strip symmetric domain. To obtain the upper bound on $\mathcal{P}^{\chi}_{\mathcal{D}} [       \mathcal{V}^{+ \backslash -}_{\Lambda}    ] $, observe,

\begin{align*}
  \mathcal{P}^{\chi}_{\Lambda} \big[  \mathcal{I}_0 \longleftrightarrow \widetilde{\mathcal{I}_0}  \big] \leq \mathcal{P}^{\chi}_{\Lambda} \big[ \mathcal{V}^{+ \backslash -}_{\Lambda} \big]         \text{ } \text{ , } 
\end{align*}

\noindent and furthermore, that,

\begin{align*}
  \mathcal{P}^{\chi}_{\Lambda} \big[ \mathcal{V}^{+ \backslash -}_{\Lambda} \big] \leq 1 - \mathcal{P}^{\chi}_{\Lambda} \big[ \big(\mathcal{H}^{+ \backslash -}_{\Lambda} \big)^{\mathrm{odd}} \big] < 1 - \mathcal{P}^{\chi}_{\mathscr{D}} \big[ \big(\mathcal{H}^{+ \backslash -}_{\mathscr{D}} \big)^{\mathrm{odd}} \big] \text{ } 
\end{align*}

\bigskip

\noindent On the other hand, if the left and right boundaries of the strip \textit{symmetric domain intersect}, consider the following. If such an intersection occurs, then there exists some line situated within the square lattice that the left boundary of the \textit{strip symmetric domain} intersects with, which we denote as $\mathscr{L}_{\mathrm{int}}$. To obtain a smaller \textit{strip symmetric domain} about which the probability of a vertical crossing will occur, perform reflections about the intersection line $\mathscr{L}_{\mathrm{int}}$, in which one obtains a mixed-spin configuration with the same distribution of $+$, or of $-$, face variables that were present in the path before the left boundary of the \textit{strip symmetric domain} intersects $\mathscr{L}_{\mathrm{int}}$. Furthermore, observe,

\begin{align*}
 \mathcal{P}^{\chi}_{\mathscr{D}} \big[     \mathcal{I}_0 \longleftrightarrow \widetilde{\mathcal{I}_0}         \big]   \leq   \mathcal{P}^{\chi}_{\mathscr{D}^{\prime\prime}} \big[       \mathcal{I}_0 \longleftrightarrow \widetilde{\mathcal{I}_0}     \big]       \text{ } \text{ , } 
\end{align*}

\noindent where in the upper bound probability, the support $\mathscr{D}^{\prime\prime} \subsetneq \mathscr{D}$ on the Ashkin-Teller measure denotes a larger \textit{strip symmetric domain} consisting of the original \textit{strip symmetric domain}, $\mathscr{D}$, with identical boundary conditions on $\partial \mathscr{D} \cap \partial \mathscr{D}^{\prime\prime} \neq \emptyset$, and on $\partial \mathscr{D} \cap \partial \mathscr{D}^{\prime} \neq \emptyset$, in addition to the \textit{strip symmetric domain} $\mathscr{D}^{\prime}$, introduce another \textit{strip symmetric domain}, $\mathscr{D}^{\prime\prime}$, satisfying $\mathscr{D} \supsetneq \mathscr{D}^{\prime} \cap \mathscr{D}^{\prime\prime} \neq \emptyset$ which satisfies, WLOG,

\[
F^{\mathrm{even}} \big( \textbf{Z}^2 \big) \supsetneq \mathscr{D}^{\prime\prime} \equiv \text{ } 
\left\{\!\begin{array}{ll@{}>{{}}l}       \gamma_L     \text{ } \text{ , }   \\  \gamma^{\prime}_L
          \text{ } \text{ , } \\        \mathscr{B} \text{ } \text{ , } \\  \mathscr{B}^{\prime} \text{ } \text{ , } 
\end{array}\right.
\]

\noindent where from the definition of $\mathscr{D}^{\prime} \equiv \gamma_L \cup \gamma^{\prime}_L \cup \mathscr{B} \cup \mathscr{B}^{\prime}$ above, one may furthermore conclude that,

\begin{align*}
 \mathcal{P}^{\chi}_{\mathscr{D}^{\prime\prime}}\big[     \mathcal{I}_0 \longleftrightarrow \widetilde{\mathcal{I}_0}         \big] \leq 1 -  \mathcal{P}^{\chi}_{\mathscr{D}^{\prime\prime}} \big[     \gamma_L \longleftrightarrow \mathscr{B}^{\prime}  \big] \text{ } \text{ , } 
\end{align*}

\noindent in which from the inequality above, the random variable $ \mathcal{P}^{\chi}_{\mathscr{D}^{\prime\prime}}\big[     \mathcal{I}_0 \longleftrightarrow \widetilde{\mathcal{I}_0}         \big]$ can be bound above with another random variable dependent upon the crossing probability, under the same support by the new \textit{strip symmetric domain} $\mathscr{D}^{\prime\prime}$ and identical boundary conditions $\chi$. By finite energy, there exists a strictly positive constant $C_{\mathrm{temp}}$ for which,

\begin{align*}
  \mathcal{P}^{\chi}_{\mathscr{D}^{\prime\prime}} \big[     \gamma_L \longleftrightarrow \mathscr{B}^{\prime}  \big]  \geq C_{\mathrm{temp}} \text{ }\text{ , }    \tag{Event $23$}
\end{align*}

\noindent hence implying that there exists some suitable upper bound for,

\begin{align*}
 \textbf{P}^{\chi}_{\mathscr{D}^{\prime\prime}}\big[        \mathcal{I}_0 \longleftrightarrow \widetilde{\mathcal{I}_0}       \big]   \text{ } 
\end{align*}

\noindent concludes the argument for intersecting left and right boundaries of the strip symmetric domain. \boxed{}

\bigskip

\noindent \textit{Completing the proof of Proposition AT 1 in light of upper bounds for $\textbf{P}^{\chi}_{\Lambda} [ \mathcal{I}_0 \longleftrightarrow \widetilde{\mathcal{I}_0} ]$ for non-intersecting, and intersecting, boundaries}. To conclude with an upper bound for the segment connectivity event under $\chi \sim \textbf{B}\textbf{C}_{++}$, 

\begin{align*}
 \mathcal{P}^{\chi}_{\Lambda} \big[[ 0 , \lfloor \delta^{\prime\prime} n \rfloor ] \times \{ 0 \}     \longleftrightarrow [ i , i + \lfloor \delta^{\prime\prime} n \rfloor ] \times \{ n \}   \big]   \text{ } \text{ , } 
\end{align*}

\noindent observe the following inequality,

\begin{align*}
1 - \mathcal{P}^{\chi}_{\mathscr{D}} \big[ \big(\mathcal{H}^{+ \backslash -}_{\mathscr{D}} \big)^{\mathrm{odd}} \big] \geq \mathcal{P}^{\chi}_{\Lambda} \big[     \mathcal{V}^{+ \backslash -}_{\Lambda}  \big]      \geq \textbf{P}^{\chi}_{\Lambda}\big[        \mathcal{I}_0 \longleftrightarrow \widetilde{\mathcal{I}_0}       \big]  \geq  \mathcal{P}^{\chi}_{\Lambda} \big[[ 0 , \lfloor \delta^{\prime\prime} n \rfloor ] \times \{ 0 \}     \longleftrightarrow [ i , i + \lfloor \delta^{\prime\prime} n \rfloor ] \times \{ n \}   \big]   \text{ } \text{ , } 
\end{align*}

\noindent holds for non-intersecting domain boundaries, while the following inequality,

\begin{align*}
    \mathcal{P}^{\chi}_{\mathscr{D}^{\prime\prime}} \big[   \mathcal{I}_0 \longleftrightarrow \widetilde{\mathcal{I}_0}    \big]  \geq \textbf{P}^{\chi}_{\Lambda}\big[        \mathcal{I}_0 \longleftrightarrow \widetilde{\mathcal{I}_0}       \big]  \geq  \mathcal{P}^{\chi}_{\Lambda} \big[[ 0 , \lfloor \delta^{\prime\prime} n \rfloor ] \times \{ 0 \}     \longleftrightarrow [ i , i + \lfloor \delta^{\prime\prime} n \rfloor ] \times \{ n \}   \big]   \text{ } \text{ , } 
\end{align*}

\noindent holds for intersecting domain boundaries, from which we conclude the argument. \boxed{}


\subsection{Ashkin-Teller bridging events in the strip}

\noindent From the estimate above, now we turn towards establishing weakened RSW results with the following. WLOG suppose that $\mathscr{F}\mathscr{C}^{\mathrm{even}} \equiv \mathscr{F}\mathscr{C}$, and that $\Lambda^{\mathrm{even}} \equiv \Lambda$. Consider the family of crossing events,

\begin{align*}
     \mathcal{B}^{\mathrm{AT}}_{+ \backslash -} (j) \equiv \mathcal{B}_{+ \backslash -}(j) \equiv \big\{    \mathcal{I}_j               {\underset{\mathscr{F}\mathscr{C}\cap \Lambda}{\longleftrightarrow}}             \mathcal{I}_{j+1}     \big\}           \text{ } \text{ , } 
\end{align*}

\noindent which correspond to the series of modified bridge events that were introduced in {\color{blue}[11]}, and also adapted in \textit{2.4}, for concluding bounds on crossing probabilities in strips. To this end, introduce the following series of results below.

\bigskip

\noindent \textbf{Lemma} \textit{AT 2.7} (\textit{Ashkin-Teller equivalent of} \textbf{Lemma} \textit{2.7.4} \textit{for a lower-bound of bridging events for the six-vertex model}). WLOG, suppose,

\begin{align*}
 \bigg\{          \mathcal{I}_j \underset{F^{\mathrm{even}} \cap \mathscr{F}\mathscr{C} \cap \Lambda}{\longleftrightarrow}   \widetilde{\mathcal{I}_j}   \bigg\}       \equiv   \bigg\{          \mathcal{I}_j \underset{ \mathscr{F}\mathscr{C} \cap \Lambda}{\longleftrightarrow}   \widetilde{\mathcal{I}_j}   \bigg\}     \equiv  \big\{        \mathcal{I}_j {\longleftrightarrow}   \widetilde{\mathcal{I}_j}    \big\} \text{ } \text{ , } 
\end{align*}

\noindent $\forall j \geq 0$, from which the lower-bound for \textit{Ashkin-Teller} bridging events takes the form,

\begin{align*}
      \mathcal{P}^{\chi}_{\Lambda} \big[      \mathcal{B}_{+ \backslash -}(j)                      \big]    \geq                C_{\mathrm{AT}} \text{ }  \big(  \mathcal{P}^{\chi}_{\Lambda} \big[         \mathcal{I}_0 {\longleftrightarrow}   \widetilde{\mathcal{I}_0}    \big] \big)^2                                           \text{ } \text{ , } 
\end{align*}

\noindent for a suitably chosen, strictly positive, constant $C_{\mathrm{AT}}$ appearing in the lower-bound, and $++$ boundary conditions $\chi$.

\bigskip

\noindent \textit{Proof of Lemma AT 2.7}. We follow the strategy that was implemented for the six-vertex model in the case of flat and sloped boundary conditions adapted from {\color{blue}[11]}, in which we introduce the intersection of crossing events from the top to bottom of the strip,

\begin{align*}
  \mathscr{C}^{\mathrm{AT}} \equiv \mathscr{C} \equiv \big\{    \mathcal{I}_{-j} \underset{\mathscr{F}\mathscr{C}\cap \Lambda}{\overset{+ \backslash -}{\longleftrightarrow}}   \widetilde{\mathcal{I}_{-j}} \big\} \cap  \big\{   \mathcal{I}_{-k} \underset{\mathscr{F}\mathscr{C}\cap \Lambda}{\overset{+ \backslash -}{\longleftrightarrow}}   \widetilde{\mathcal{I}_{-k}}      \big\} \text{ }\text{ , } 
\end{align*}

\noindent for $k>j$, which, can further be made use of to conclude,

\begin{align*}
  \mathcal{P}_{\Lambda}^{\chi} \big[ \mathscr{C} \big] \equiv  \mathcal{P}_{\Lambda}^{\chi} \big[         \big\{    \mathcal{I}_{-j} \underset{\mathscr{F}\mathscr{C} \cap \Lambda}{\overset{+ \backslash -}{\longleftrightarrow}}   \widetilde{\mathcal{I}_{-j}} \big\} \cap  \big\{   \mathcal{I}_{-k} \underset{\mathscr{F}\mathscr{C} \cap \Lambda}{\overset{+ \backslash -}{\longleftrightarrow}}   \widetilde{\mathcal{I}_{-k}}      \big\}      \big]  \overset{(\mathrm{FKG})}{\geq} \text{ }    \mathcal{P}_{\Lambda}^{\chi} \big[ \mathcal{I}_{-j} \underset{\mathscr{F}\mathscr{C} \cap \Lambda}{\overset{+ \backslash -}{\longleftrightarrow}}   \widetilde{\mathcal{I}_{-j}}   \big] \text{ }  \mathcal{P}_{\Lambda}^{\chi} \big[  \mathcal{I}_{-k} \underset{\mathscr{F}\mathscr{C} \cap \Lambda}{\overset{+ \backslash -}{\longleftrightarrow}}   \widetilde{\mathcal{I}_{-k}}   \big] \tag{Event $24$} \\ \Updownarrow \\      \mathcal{P}^{\chi}_{\Lambda} \bigg[    \mathcal{B}_{+ \backslash -}(j) \text{ }    \big|    \text{ }       \big\{    \mathcal{I}_{-j} \underset{\mathscr{F}\mathscr{C} \cap \Lambda}{\overset{+ \backslash -}{\longleftrightarrow}}   \widetilde{\mathcal{I}_{-j}} \big\} \cap  \big\{   \mathcal{I}_{-k} \underset{\mathscr{F}\mathscr{C} \cap \Lambda}{\overset{+ \backslash -}{\longleftrightarrow}}   \widetilde{\mathcal{I}_{-k}}      \big\}        \bigg]  \geq C_{\mathrm{AT}}     \text{ , } 
\end{align*}

\noindent implying that it suffices to provide the lower bound of $C_{\mathrm{AT}}$ for the bridging event occurring conditionally upon the occurrence of,

  \begin{figure}
\begin{align*}
\includegraphics[width=0.53\columnwidth]{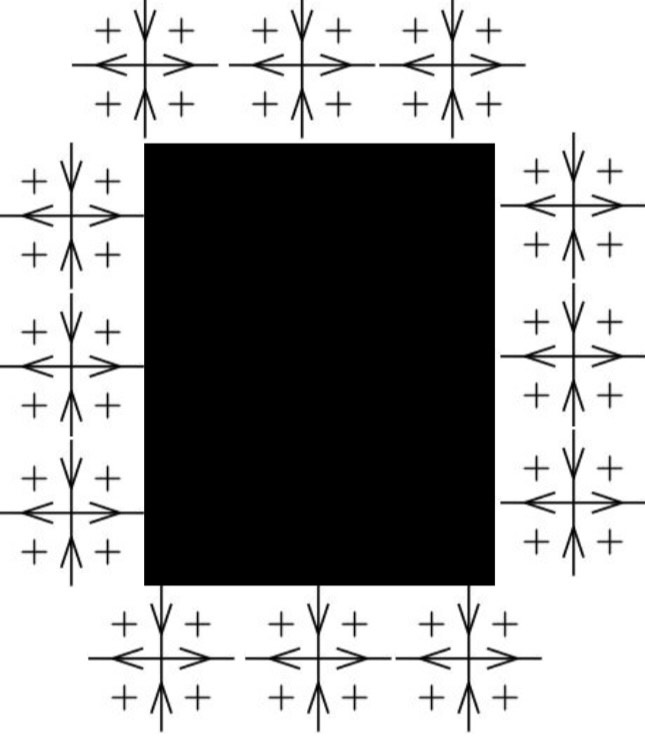}
\end{align*}
\caption{\textit{A depiction of the distribution of the distribution of $+$ spins on even, or on odd faces, outside of the black finite volume, as used in arguments for lower bounding the occurrence of the conditional bridging event with the occurrence of the conditional horizontal event}.}
\end{figure}

\begin{align*}
   \big\{    \mathcal{I}_{-j} \underset{\mathscr{F}\mathscr{C} \cap \Lambda}{\overset{+ \backslash -}{\longleftrightarrow}}   \widetilde{\mathcal{I}_{-j}} \big\} \cap  \big\{   \mathcal{I}_{-k} \underset{\mathscr{F}\mathscr{C} \cap \Lambda}{\overset{+ \backslash -}{\longleftrightarrow}}   \widetilde{\mathcal{I}_{-k}}      \big\}     \text{ } \text{ , } 
\end{align*}

\noindent under $\mathcal{P} \big[ \cdot \big]$. Proceeding, consider the conditional horizontal crossing event,

\begin{align*}
  \mathcal{P}^{\chi}_{\Lambda} \big[   \mathcal{H}_{+ \backslash -} ( \mathscr{D} )  \text{ }     \big|  \text{ }     \mathrm{sign} \big(   \mathscr{F}|_{\mathscr{D}^c \cap \Lambda} \big) \equiv +                          \big]   \text{ } \text{ , } 
\end{align*}

\noindent dependent upon the distribution of $+$ faces $\mathscr{F}|_{\mathscr{D}^c \cap \Lambda}$ of the complementary region to the \textit{Ashkin-Teller symmetric domain} $\mathscr{D}$ in the strip, which satisfies,

\begin{align*}
 \mathcal{P}^{\chi}_{\Lambda} \big[   \mathcal{B}_{+ \backslash -}(j)                              \big|     \mathrm{sign} \big(   \mathscr{F}|_{\mathscr{D}^c \cap \Lambda} \big) \equiv +                    \big] \overset{\mathrm{Bayes}}{\equiv}  \frac{\mathcal{P}^{\chi}_{\Lambda} \big[   \mathcal{B}_{+ \backslash -}(j)  \cap \big\{      \mathrm{sign} \big(   \mathscr{F}|_{\mathscr{D}^c \cap \Lambda} \big) \equiv +                  \big\}     \big]   }{\mathcal{P}^{\chi}_{\Lambda} \big[       \mathrm{sign} \big(   \mathscr{F}|_{\mathscr{D}^c \cap \Lambda} \big) \equiv +                     \big]  }  \tag{Event $25$} \\ \overset{(\mathrm{FKG})}{\geq}    \frac{\mathcal{P}^{\chi}_{\Lambda} \big[   \mathcal{B}_{+ \backslash -}(j)  \big] \mathcal{P}^{\chi}_{\Lambda} \big[     \mathrm{sign} \big(   \mathscr{F}|_{\mathscr{D}^c \cap \Lambda} \big) \equiv +              \big]   }{\mathcal{P}^{\chi}_{\Lambda} \big[       \mathrm{sign} \big(   \mathscr{F}|_{\mathscr{D}^c \cap \Lambda} \big) \equiv +                     \big]  }   \\ \overset{\mathcal{H}_{+ \backslash -} ( \mathscr{D} ) \subsetneq     \mathcal{B}_{+ \backslash -}(j)  }{\geq}  \frac{\mathcal{P}^{\chi}_{\Lambda} \big[   \mathcal{H}_{+ \backslash -} ( \mathscr{D} )   \big] \mathcal{P}^{\chi}_{\Lambda} \big[     \mathrm{sign} \big(   \mathscr{F}|_{\mathscr{D}^c \cap \Lambda} \big) \equiv +              \big]    }{\mathcal{P}^{\chi}_{\Lambda} \big[       \mathrm{sign} \big(   \mathscr{F}|_{\mathscr{D}^c \cap \Lambda} \big) \equiv +                     \big]  } \\ \equiv \mathcal{P}^{\chi}_{\Lambda} \big[   \mathcal{H}_{+ \backslash -} ( \mathscr{D} )  |   \mathrm{sign} \big(   \mathscr{F}|_{\mathscr{D}^c \cap \Lambda} \big) \equiv +              \big]    \\ \Updownarrow \\ 
 \mathcal{P}^{\chi}_{\Lambda} \big[    \mathcal{B}_{+ \backslash -}(j)      \big|      \mathrm{sign} \big(   \mathscr{F}|_{\mathscr{D}^c \cap \Lambda} \big) \equiv +        \big]  \geq \mathcal{P}^{\chi}_{\Lambda} \big[    \mathcal{H}_{+ \backslash -} ( \mathscr{D} )       \big|      \mathrm{sign} \big(   \mathscr{F}|_{\mathscr{D}^c \cap \Lambda} \big) \equiv +        \big]            \text{ } \text{ , } 
\end{align*}

\noindent where $\mathscr{D}$ satisfies the same set of properties that were introduced in the argument for the six-vertex bridging events, namely that the existence of a crossing between $\mathcal{I}_{-j}$ and $\widetilde{\mathcal{I}_{-j}}$, and between $\mathcal{I}_{-k}$ and $\widetilde{\mathcal{I}_{-k}}$, respectively induce left and right boundaries of $\mathscr{D}$, while the top and bottom boundaries of $\mathscr{D}$. From such a symmetric domain, the conditional horizontal crossing probability across $\mathscr{D}$, conditionally upon the $+$ distribution of faces, satisfies,

\begin{align*}
     \mathcal{P}^{\chi}_{\Lambda} \big[      \mathcal{H}_{+ \backslash -} \big( \mathscr{D} \big) |             \mathrm{sign} \big(   \mathscr{F}|_{\mathscr{D}^c \cap \Lambda} \big) \equiv +                     \big] \geq    \mathscr{C}_{+}   \mathcal{P}^{\chi}_{\Lambda} \big[      \mathcal{H}_{+ \backslash -} \big( \mathscr{D} \big)  \big]       \tag{Event $26$}  \\ \geq     \mathscr{C}_{+}  \big(      1 -  \mathcal{P}^{\chi}_{\Lambda} \big[    \mathcal{V}^{\mathrm{odd}}_{+ \backslash -} \big( \mathscr{D} \big)            \big]  \big)         \text{ } \text{ , } 
\end{align*}

\noindent where in the lower bound for the conditional horizontal crossing event with a vertical crossing event, we make use of the fact that \textit{duality} for the mixed-spin Ashkin-Teller representation encapsulates that an odd (resp. even) horizontal crossing event on $F^{\mathrm{odd}}$ (resp. $F^{\mathrm{even}}$) has a dual counterpart corresponding to a vertical crossing event on even (resp. odd) sublatices of $\textbf{Z}^2$, which are respectively given by $F^{\mathrm{even}}$, and by $F^{\mathrm{odd}}$. Otherwise, the remaining component in the lower bound is given by a suitably chosen, strictly positive $\mathscr{C}_{+}$, for which,

\begin{align*}
    \mathcal{P}^{\chi}_{\Lambda} \big[   \mathrm{sign} \big(   \mathscr{F}|_{\mathscr{D}^c \cap \Lambda} \big) \equiv +             \big]         \geq     \mathscr{C}_{+}   \text{ } \text{ , } 
\end{align*}

\noindent is satisfied by finite energy. Proceeding, from additional properties of the \textit{symmetric} domain, the following vertical crossing probability,

\begin{align*}
 \mathcal{P}^{-\chi}_{\Lambda} \big[   \mathcal{V}^{\mathrm{odd}}_{+ \backslash -}\big( \mathscr{D} \big)      \big]     \overset{(\mathrm{marginal\text{ } CBC})}{\leq} \mathcal{P}^{\chi}_{\Lambda} \big[   \mathcal{V}^{\mathrm{odd}}_{+ \backslash -}\big( \mathscr{D} \big)      \big]  \leq \mathcal{P}^{\chi}_{\Lambda} \big[      [ 0, \lfloor \delta^{\prime\prime} n \rfloor ] \times \{ 0 \} {\longleftrightarrow}  [i , i+ \delta^{\prime\prime}n  ]  \times \{n \}       \big]      \text{ } \text{ , } 
\end{align*}

\noindent in which the probability, under $\chi$, of the segment connectivity occurring upper bounds the probability of the vertical event happening, across \textit{odd} faces of the square lattice, under $\chi$. Relatedly, previous arguments demonstrate,

\begin{align*}
 \mathcal{P}^{\chi}_{\Lambda} \big[  [- \lfloor \delta^{\prime\prime} n \rfloor   , 2 \lfloor   \delta^{\prime\prime} n \rfloor          \big]  \times \{ 0 \} {\longleftrightarrow}  [- \lfloor \delta^{\prime\prime}  n \rfloor   , 2 \lfloor   \delta^{\prime\prime} n  \rfloor          \big]  \times \{ n \} \big] \leq 1 - c_0   \text{ } \text{ , }      \tag{Event $27$} \end{align*}

\noindent as it has been show, from previous arguments for the upper bound provided in , that,

\begin{align*}
 \mathcal{P}^{\chi}_{\Lambda} \big[     [ 0, \lfloor \delta^{\prime\prime} n \rfloor ] \times \{ 0 \} {\longleftrightarrow}  [i , i+ \delta^{\prime\prime}n  ]  \times \{n \}  \big]  \leq 1 - c  \text{ } \text{ , }    \tag{Event $28$} 
\end{align*}

\noindent from which we conclude the argument. \boxed{}

\bigskip

\noindent For the final result below, we introduce an analog of \textbf{Proposition} \textit{2.7.5} which provided a lower bound for the segment connectivity event in the six-vertex model.

\bigskip

\noindent \textbf{Proposition} \textit{AT 2} (\textit{adaptation of crossing probability estimate between two segment connectivity events from the six-vertex model onto the symmetric mixed-spin representation}). WLOG, suppose $F^{\mathrm{even}} \big( \Lambda \cap \textbf{Z}^2\big)\equiv F$. To quantify a lower bound for a segment connectivity event, fix $\chi \in \textbf{B}\textbf{C}_{++}$, $n^{\prime} > n > 0$, and $\rho > 0$, from which the lower bound takes the form, given strictly positive $\delta^{\prime\prime}$ from previous arguments,

\begin{align*}
  \mathcal{P}^{\chi}_{\Lambda} \big[  [-n^{\prime} , 2 n^{\prime} ] \times \{ 0 \} \longleftrightarrow                \times \{ \rho \lfloor \delta^{\prime\prime} n \rfloor  \}   \big] \geq \bigg[     \sqrt[\mathscr{N}^{\prime}]{\mathscr{O}^{\prime}}   \mathcal{P}^{\chi}_{\Lambda} \big[      [ 0 , \lfloor \delta^{\prime\prime} n \rfloor  ] \times \{ 0 \} \longleftrightarrow \textbf{Z} \times \{ n \}            \big]        \bigg]^{\mathscr{N}^{\prime}}  \text{ } \text{ , } 
\end{align*}

\noindent where in the inequality above, the multiplicative constant to the lower bound probability for quantifying connectivity with $\textbf{Z} \times \{ n \}$, $\mathscr{O}^{\prime}$, is strictly less than $1$, while the power to which the probability in the lower bound is raised, $\mathscr{N}^{\prime}$, is another strictly positive integer.

\bigskip

\noindent \textit{Proof of Proposition AT 2}. WLOG, suppose $\rho \in \textbf{Z}$. From previous estimates on the probability of a bridging event occurring, introduce a lower bound for the segment connectivity event with the bridging event, in which,

\begin{align*}
   \mathcal{P}^{\chi}_{\Lambda} \big[    \{ 0 \} \times [ - n^{\prime} , 2 n^{\prime} ]     \longleftrightarrow      \{  \rho \lfloor \delta^{\prime\prime} n \rfloor      \} \times [ - n^{\prime} , 2 n^{\prime} ]      \big]    \tag{Event $29$}            \\ \Updownarrow \\   \mathcal{P}^{\chi}_{\Lambda} \bigg[ \underset{0 \leq i \leq f}{\underset{\rho_f \equiv      \rho \lfloor \delta^{\prime\prime} n \rfloor  }{\underset{\rho_{i+1} \equiv \rho_i + \delta^{\prime\prime\prime}}{\underset{\rho_0 \equiv 0 }{\bigcap}}}} \bigg\{   \{ \rho_i \lfloor \delta^{\prime\prime} n \rfloor  \} \times [ - n^{\prime} , 2 n^{\prime} ]         \longleftrightarrow      \{  \rho_{i+1} \lfloor \delta^{\prime\prime} n \rfloor      \} \times [ - n^{\prime} , 2 n^{\prime} ]              \bigg\}  \bigg] \\ \overset{(\mathrm{marginal \text{ } FKG})}{\geq} \underset{0 \leq i \leq f}{\underset{\rho_f \equiv      \rho \lfloor \delta^{\prime\prime} n \rfloor  }{\underset{\rho_{i+1} \equiv \rho_i + \delta^{\prime\prime\prime}}{\underset{\rho_0 \equiv 0 }{\prod}}}} \mathcal{P}^{\chi}_{\Lambda} \big[   \{ \rho_i \lfloor \delta^{\prime\prime} n \rfloor  \} \times [ - n^{\prime} , 2 n^{\prime} ]         \longleftrightarrow      \{  \rho_{i+1} \lfloor \delta^{\prime\prime} n \rfloor      \} \times [ - n^{\prime} , 2 n^{\prime} ]        \big]   \\ \geq \big( \mathcal{P}^{\chi}_{\Lambda} \big[    \mathcal{B}_{+ \backslash - }(j)          \big]    \big)^f        \\ \overset{(\textbf{Lemma}\text{ } \textit{AT 2.7})}{\geq}          C^f_{\mathrm{AT}} \text{ }  \big(  \mathcal{P}^{\chi}_{\Lambda} \big[         \mathcal{I}_0 {\longleftrightarrow}   \widetilde{\mathcal{I}_0}    \big] \big)^{2f}    \text{ } \text{ , } 
\end{align*}

\noindent where, in the first line, the probability under $++$ boundary conditions $\chi$ of the segment connectivity event occurring is equivalent to the union of 'smaller' segment connectivity events occurring across smaller scales, as, for suitable $\delta^{\prime\prime\prime} < \delta^{\prime\prime}$,

\begin{align*}
 \big\{      \{ 0 \} \times [ - n^{\prime} , 2 n^{\prime} ]     \longleftrightarrow      \{  \rho \lfloor \delta^{\prime\prime} n \rfloor      \} \times [ - n^{\prime} , 2 n^{\prime} ]          \big\} \equiv \underset{0 \leq i \leq f}{\underset{\rho_f \equiv      \rho \lfloor \delta^{\prime\prime} n \rfloor  }{\underset{\rho_{i+1} \equiv \rho_i + \delta^{\prime\prime\prime}}{\underset{\rho_0 \equiv 0 }{\bigcap}}}} \bigg\{   \{ \rho_i \lfloor \delta^{\prime\prime} n \rfloor  \} \times [ - n^{\prime} , 2 n^{\prime} ]         \longleftrightarrow      \{  \rho_{i+1} \lfloor \delta^{\prime\prime} n \rfloor      \} \times [ - n^{\prime} , 2 n^{\prime} ]              \bigg\}   \text{ } \text{ . } 
\end{align*}

\noindent Next, from the lower bound obtained from \textbf{Lemma} \textit{AT 2.7}, observe that it suffices to show that a bound of the following form holds, in which,

\begin{align*}
      \mathcal{P}^{\chi}_{\Lambda} \big[    \mathcal{I}_0         \longleftrightarrow     \widetilde{\mathcal{I}_0}   \big] \geq        \text{ } \mathcal{C}_{\mathrm{AT}}     \text{ }   \big( \mathcal{P}^{\chi}_{\Lambda} \big[  \mathcal{I}_0                \longleftrightarrow    \textbf{Z} \times \{ n \}   \big]\big)^2 { , }           \tag{Event $30$} 
\end{align*}

\noindent for a suitably chosen, strictly positive, prefactor $\mathcal{C}_{\mathrm{AT}}$. Beginning with the probability in the upper bound for the quantifying the connectivity between $\mathcal{I}_0$ and $\widetilde{\mathcal{I}_0}$, it is possible that,

\begin{align*}
   \mathcal{P}^{\chi}_{\Lambda} \big[ \mathcal{I}_0 \longleftrightarrow        (- \infty , 0 ] \times \{ n \}      \big]           \text{ } \text{ , } \tag{Event $31$}
\end{align*}

\noindent which represents a mixed-spin configuration for which a connectivity event occurs between $\mathcal{I}_0$ and $ (- \infty , 0 ] \times \{ n \}  $, or,

\begin{align*}
   \mathcal{P}^{\chi}_{\Lambda} \big[ \mathcal{I}_0 \longleftrightarrow          [ \lfloor \delta^{\prime\prime} n \rfloor  , + \infty ) \times \{ n \}       \big]    \text{ } \text{ , } 
\end{align*}

\noindent which represents a mixed-spin configuration for which a connectivity event occurs between $\mathcal{I}_0$ and $\big\{ \delta^{\prime\prime} n \rfloor  , + \infty ) \times \{ n \}  \big\}$, while finally, 

\begin{align*}
   \mathcal{P}^{\chi}_{\Lambda} \big[ \widetilde{\mathcal{I}_0} \longleftrightarrow        [ \lfloor \delta^{\prime\prime} n \rfloor , + \infty ) \times \{ 0 \}       \big]     \text{ } \text{ , } 
\end{align*}

\noindent which represents a mixed-spin configuration for which a connectivity event occurs between $\mathcal{I}_0$ and . From each of the three previously mentioned possibilities, in the second case, for the event $\big\{   \mathcal{I}_0 \longleftrightarrow          [ \lfloor \delta^{\prime\prime} n \rfloor  , + \infty ) \times \{ n \}      \big\}$ to occur,

\begin{align*}
\mathcal{P}^{\chi}_{\Lambda} \big[ \mathcal{I}_0 \longleftrightarrow \widetilde{\mathcal{I}_0} \big] \geq   \mathcal{P}^{\chi}_{\Lambda} \bigg[  \big\{ \mathcal{I}_0 \longleftrightarrow \widetilde{\mathcal{I}_0}     \big\} \cap   \big\{ \mathcal{I}_0 \longleftrightarrow          [ \lfloor \delta^{\prime\prime} n \rfloor  , + \infty ) \times \{ n \}   \big\}     \bigg]    \overset{(\mathrm{Marginal \text{ } (FKG)})}{\geq}         \mathcal{P}^{\chi}_{\Lambda} \big[  \mathcal{I}_0 \longleftrightarrow \widetilde{\mathcal{I}_0}     \big] \text{ } \\ \times   \mathcal{P}^{\chi}_{\Lambda} \big[ \mathcal{I}_0 \longleftrightarrow          [ \lfloor \delta^{\prime\prime} n \rfloor  , + \infty ) \times \{ n \}     \big]        \text{ } \text{ , 
    } \tag{\textit{m-lower bound}}
\end{align*}

\noindent following one application of (marginal FKG), from which performing further manipulation yields the final lower bound,

\begin{align*}
  (\textit{m-lower bound}) \geq \text{ }  \mathcal{C}_{\mathrm {AT}} \big(  \mathcal{P}^{\chi}_{\Lambda} \big[  \mathcal{I}_0 \longleftrightarrow \widetilde{\mathcal{I}_0}      \big]     \big)^2     \text{ } \text{ . } 
\end{align*}

\noindent Besides the second case, the first case follows from the fact that,

\begin{align*}
   \mathcal{P}^{\chi}_{\Lambda} \big[ \mathcal{I}_0 \longleftrightarrow \widetilde{\mathcal{I}_0} \big]  \geq   \mathcal{P}^{\chi}_{\Lambda} \bigg[  \big\{ \mathcal{I}_0 \longleftrightarrow \widetilde{\mathcal{I}_0}     \big\} \cap   \big\{ \mathcal{I}_0 \longleftrightarrow          (- \infty , 0 ] \times \{ n \}  \big\}     \bigg]    \overset{(\mathrm{Marginal \text{ } (FKG)})}{\geq}         \mathcal{P}^{\chi}_{\Lambda} \big[  \mathcal{I}_0 \longleftrightarrow \widetilde{\mathcal{I}_0}     \big] \text{ } \\ \times  \mathcal{P}^{\chi}_{\Lambda} \big[ \mathcal{I}_0 \longleftrightarrow          (- \infty , 0 ] \times \{ n \}    \big]     \text{ } \text{ , } 
\end{align*}

\noindent in which, by making use similar arguments for obtaining the lower bound for (\textit{m-lower bound}), concludes the argument for the first case. Finally, for the third case, as in the argument for the six-vertex model under flat and sloped boundary conditions, the fact that,

\begin{align*}
    \mathcal{P}^{\chi}_{\Lambda} \big[ \mathcal{I}_0 \longleftrightarrow \widetilde{\mathcal{I}_0} \big]          \text{ } \text{ , }
\end{align*}

\noindent occurs immediately provides the desired lower bound, from which we conclude the argument. \boxed{}

\section{The generalized random-cluster model}

\noindent Besides obtaining weakened crossing probability estimates in strips of the square lattice for the Ashkin-Teller model from the arguments introduced for RSW results under sloped boundary conditions for the six-vertex model, a final area of interest relates to obtaining crossing probability estimates for the generalized random cluster model, described in {\color{blue}[23]} through the random-cluster representation of the Ashkin-Teller model. As expected, as does the Ashkin-Teller model, the generalized random-cluster model satisfies a similar range of properties, which can be applied over strips of the square lattice as the finite volume approaches the infinite volume. Boundary conditions for this model of interest are encoded through bond configurations that are defined over $\textbf{Z}^2$. The third model that we are considering below can be defined with a product probability measure over the set of bonds of the square lattice. It is of interest to determine whether cluster representations exist for other models, and furthermore, whether there exists any connection with the Ashkin-Teller model, which could yield information on crossing probabilities. Below, we gather necessary items from this reference, where the symbols $\sigma$ and $\tau$ indicate the two components of the configuration $\omega$, from the \textit{generalized random-cluster model} sample space,

\begin{align*}
   \omega \in \Omega^{\mathrm{g-r-cm}} \equiv  \Omega \equiv   [ \{-1 , 1 \}     \times \{-1,1\} ]^{\textbf{Z}^2} \text{ } \text{ , } 
\end{align*}

\noindent which is also a product space, over $\textbf{Z}^2$. 

\subsection{Generalized random-cluster objects}

\noindent \textbf{Definition} \textit{20} (\textit{percolation bond product measure for defining the probability measure of the generalized random-cluster model}, {\color{blue}[23]}). For a configuration over the set of all bonds $b$ from $\textbf{Z}^2$, denoted with $\mathcal{B} ( \textbf{Z}^2) \equiv \mathcal{B}$, introduce the following four probability measures,

\begin{align*}
   \lambda_b (( 0,0)) \equiv a_0(b)         \text{ } \text{ , }  \\  \lambda_b((1,0)) \equiv  a_{\sigma}(b)     \text{ } \text{ , }  \\      \lambda_b((1,1)) \equiv       a_{\sigma \tau}(b)         \text{ } \text{ , }  \\           \lambda_b((0,1)) \equiv a_{\tau}(b)    \text{ } \text{ , } 
\end{align*}

\noindent over bonds, which yields the product, \textit{generalized percolation measure}, on the set of all bonds in $\mathcal{B}$, which is expressed with,

\begin{align*}
    \lambda_{\mathcal{B}}\big(  \underline{n}  \big) \equiv    \lambda_{\mathcal{B}}\big( ( \underline{n}_{\sigma}   , \underline{n}_{\tau} ) \big)   \equiv     \lambda_{\mathcal{B}}\big( ( n_{\sigma} (b)   , n_{\tau} (b) ) \big)      \equiv     \lambda_{\mathcal{B}}\big(   (   n_{\sigma,b} , n_{\tau,b} )    \big)            \text{ } \text{ , } 
\end{align*}

\noindent which, as the following product measure, admits the decomposition,

\begin{align*}
   \lambda_{\mathcal{B}}(\underline{n}) \equiv \underset{n_b \equiv (0,0)}{\underset{ b \in \mathcal{B}}{\prod}}              \lambda_b((0,0))      \underset{n_b \equiv (1,0)}{\underset{b \in \mathcal{B}}{ \prod }}            \lambda_b((1,0))     \underset{n_b \equiv (0,1)}{\underset{b \in \mathcal{B}}{\prod} }          \lambda_b((1,1))      \underset{n_b \equiv (1,1)}{\underset{b \in \mathcal{B}}{\prod }} \lambda_b((0,1))     \equiv   \underset{n_b \equiv (0,0)}{\underset{ b \in \mathcal{B}}{\prod}}          a_0(b)  \\ \times   \underset{n_b \equiv (1,0)}{\underset{b \in \mathcal{B}}{ \prod }}          a_{\sigma}(b)    \underset{n_b \equiv (0,1)}{\underset{b \in \mathcal{B}}{\prod} }          a_{\sigma\tau}(b)     \underset{n_b \equiv (1,1)}{\underset{b \in \mathcal{B}}{\prod }} a_{\tau}(b)     \text{ } \text{ , } 
\end{align*}

\noindent for a configuration on bonds, $\underline{n}$.

\bigskip

\noindent \textbf{Definition} \textit{21} (\textit{notion of set connectedness for the generalized random-cluster model}, {\color{blue}[23]}). For a bond configuration $\underline{n} \in \{ 0 , 1 \} \times \{ 0 , 1 \}$, denote,

\begin{align*}
    \big\{    i       \overset{\sigma}{\longleftrightarrow}      j        \big\}           \text{ } \text{ , } 
\end{align*}

\noindent as the $\sigma$-connectivity event between sites $i$ and $j$, where the two sites are said to be \textit{connected}, if,

\begin{align*}
         t_0 \equiv i , t_1 , \cdots , t_{k-1} , t_k \equiv j                              \text{ } \text{ , } 
\end{align*}

\noindent such that from the sequence of $t$ variables above,

\begin{align*}
         n_{\sigma} \big(   < t_i , t_{i+1} >     \big)      \equiv 1         \text{ } \text{ , } 
\end{align*}

\noindent $\forall i = 0 , \cdots , k-1$.

\bigskip

\noindent Next, we provide a definition pertaining to the cell-complex \textit{lattice structure}.

\bigskip

\noindent \textbf{Definition} \textit{22} (\textit{cell-complex lattice structure}, {\color{blue}[25]}). From the set of all bonds $\mathcal{B} ( \textbf{Z}^2)$ over the square lattice, the \textit{cell-complex lattice}, $\mathscr{L}$, satisfies,

\begin{itemize}
    \item[$\bullet$]   \underline{Property one} (\textit{boundary of the bond configuration}): Given a bond $b$, the \textit{boundary} of $b$, $\{ t , t^{\prime} \}$, for nearest-neighbors $t$ and $t^{\prime}$, is given by $\delta \mathcal{B} \equiv \{   t \in \textbf{Z}^2 : t \in \delta b  \}$, for $\delta b \equiv \{ t , t^{\prime} \}$. 
    
    \item[$\bullet$]    \underline{Property two} (\textit{plaquettes from the bond configuration}): A \textit{plaquette} in the generalized random-cluster model is a subset of $\textbf{R}^2$, unit squares, whose corners are sites.                    
    
    \item[$\bullet$]   \underline{Property three} (\textit{boundaries of plaquettes from}, \underline{Property 1}, \textit{and,} \underline{Property 2}): From \textit{generalized random-cluster model plaquettes} introduced in \underline{Property two}, the \textit{boundary} of a \textit{plaquette} is given by their boundary, which is another subset of $\textbf{R}^2$.            
\end{itemize}

\bigskip

\noindent From the \textit{generalized percolation bond measure}, below introduce the \textit{generalized random-cluster probability measure}, under $++$, and under \textit{free} boundary conditions.

\bigskip

\noindent \textbf{Definition} \textit{23} (\textit{generalized random-cluster probability measure with $++$ boundary conditions}, {\color{blue}[25]}). From the \textit{generalized percolation bond measure} provided in \textbf{Definition} \textit{20}, denote,

\begin{align*}
    \nu^{+}_{\Lambda} \big( \underline{n} | q_{\sigma} , q_{\tau} \big)   \equiv  \nu^{+} \big( \underline{n} | q_{\sigma} , q_{\tau} \big)  \equiv  \frac{\lambda_{\mathcal{B}^{+}(\Lambda)} (\underline{n}) q_{\sigma}^{N_{\sigma}( \underline{n} | \Lambda)} q_{\tau}^{N_{\tau}(\underline{n} | \Lambda )}  } {\underset{\underline{n} : \underline{n}|_{\partial \Lambda} ++}{\underset{\underline{n} \equiv (\underline{n}_{\sigma} , \underline{n}_{\tau})}{\sum} } \lambda_{\mathcal{B}^{+}(\Lambda)} (\underline{n}) q_{\sigma}^{N_{\sigma}( \underline{n} | \Lambda)} q_{\tau}^{N_{\tau}(\underline{n} | \Lambda )}  }                                 \text{ } \text{ , } 
\end{align*}

\noindent for $\underline{n}$ satisfying $++$ boundary conditions on $\Lambda$, as the \textit{generalized random-cluster probability measure} with $++$ boundary conditions.

\bigskip

\noindent \textbf{Definition} \textit{24} (\textit{generalized random-cluster probability measure with free boundary conditions}, {\color{blue}[25]}). From the \textit{generalized percolation bond measure} provided in \textbf{Definition} \textit{20}, denote,

\begin{align*}
         \nu^{\mathrm{f}}_{\Lambda}  \big( \underline{n} | q_{\sigma} , q_{\tau} \big) \equiv   \nu^{\mathrm{f}}  \big( \underline{n} | q_{\sigma} , q_{\tau} \big)  \equiv   \frac{\lambda_{\mathcal{B}(\Lambda)} (\underline{n}) q_{\sigma}^{N_{\sigma}( \underline{n} | \Lambda)} q_{\tau}^{N_{\tau}(\underline{n} | \Lambda )}}{\underset{\underline{n} : \underline{n}|_{\partial \Lambda} \mathrm{f}}{\underset{\underline{n} \equiv ( \underline{n}_{\sigma} , \underline{n}_{\tau})}{\sum}  } \lambda_{\mathcal{B}(\Lambda)} (\underline{n}) q_{\sigma}^{N_{\sigma}( \underline{n} | \Lambda)} q_{\tau}^{N_{\tau}(\underline{n} | \Lambda )} }               \text{ } \text{ , } 
\end{align*}

\noindent for $\underline{n}$ satisfying \textit{free} boundary conditions on $\Lambda$, as the \textit{generalized random-cluster probability measure} with \textit{free} boundary conditions.

\bigskip

\noindent From the probability measures introduced above for $+$ and free boundary conditions, we also turn to the statement of the FKG lattice condition, in addition to the statement of the comparison inequalities. To this end, denote two bond configurations with $\underline{n}$ and $\underline{n}^{\prime}$, as well as the following partitions of the bond space $\mathcal{B}$.

\bigskip

\noindent \textbf{Definition} \textit{24} (\textit{ordering of two bond configurations}, ({\color{blue}[25]}, \textbf{Definition} \textit{4.1})). For bond configurations $\underline{n}$ and $\underline{n}^{\prime}$, $\underline{n}$ is said to \textit{dominate} $\underline{n}^{\prime}$ if $n_b \geq n^{\prime}_b$, $\forall b \in \mathcal{B}$.

\bigskip

\noindent \textbf{Definition} \textit{25} (\textit{from ordering of bond configurations provided in the previous definition to increasing functions}, ({\color{blue}[25]]}, \textbf{Definition} \textit{4.2})). From the ordering on bond configurations introduced in \textbf{Definition} \textit{24}, a function over the bond space is said to be \textit{increasing}, if,

\begin{align*}
      \underline{n} \geq \underline{n^{\prime}}   \Longleftrightarrow             f (   \underline{n} ) \geq f( \underline{n^{\prime}} )    \text{ } \text{ , } 
\end{align*}

\noindent $\forall b \in \mathcal{B}$, and otherwise, is said to be \textit{decreasing} if,

\begin{align*}
        f \text{ } \mathrm{decreasing} \Longleftrightarrow      - f \text{ } \mathrm{increasing}       \text{ }  \text{ , } 
\end{align*}

\noindent $\forall b \in \mathcal{B}$.

\bigskip

\noindent \textbf{Definition} \textit{26} (\textit{two partitions of the bond space utilized in the proof of the FKG lattice condition for the generalized random-cluster measure}, {\color{blue}[25]}). Partition $\mathcal{B}$ into,

\begin{align*}
\mathcal{B} \supsetneq \mathcal{B}_{>}  \equiv \big\{   b \in \mathcal{B} :  a_{\sigma \tau} ( b ) a_0(b) \geq a_{\sigma} (b) a_{\tau}(b)                \big\}  \text{ } \text{ , } 
\end{align*}

\noindent and into,

\begin{align*}
  \mathcal{B} \supsetneq \mathcal{B}_{<}   \equiv \big\{   b \in \mathcal{B} :   a_{\sigma \tau}(b) a_0(b) < a_{\sigma}(b) a_{\tau} (b)              \big\}  \text{ } \text{ . } 
\end{align*}

\noindent \textbf{Lemma} \textit{GRCM 1} (\textit{FKG lattice condition for the generalized random-cluster measure for increasing functions}, ({\color{blue}[25]}, \textbf{Lemma} \textit{4.1})). The \textit{generalized percolation bond measure} provided in \textbf{Defintion} \textit{20} is FKG, in which,

\begin{align*}
\lambda_{\mathcal{B}} \big( \underline{n} \vee \underline{n}^{\prime} \big) \lambda_{\mathcal{B}} \big(   \underline{n} \wedge \underline{n}^{\prime} \big) \geq \lambda_{\mathcal{B}} \big(   \underline{n}    \big)   \lambda_{\mathcal{B}} \big(    \underline{n}^{\prime}   \big)            \text{ } \text{ , } 
\end{align*}

\noindent where $\vee$ denotes the least upper bound of $\underline{n}$ and $\underline{n}^{\prime}$, while $\wedge$ denotes the greatest lower bound of $\underline{n}$ and $\underline{n}^{\prime}$.

\bigskip

\noindent Also, the comparison inequalities take the form below, which serve as a generalization of the comparison inequalities for the random-cluster model, in which one can establish comparisons between the value of the probability measure for different values of $q$, with the value of the probability measure for occupation. Because the comparison inequalities for the generalized random-cluster model are related to comparison inequalities of the nongeneralized random-cluster model, we briefly introduce the probability measu8re for the random cluster model below as well.

\bigskip

\noindent \textbf{Definition} \textit{27} (\textit{random-cluster probability measure}). Fix a graph $G$, with edge-weight $p$ and strictly positive cluster weight $q$. For wired boundary conditions $\mathrm{w}$, the random-cluster probability measure for the random-cluster model is given by,

\begin{align*}
   \phi^{\mathrm{w}}_{G , p , q} \equiv  \phi^{\mathrm{w}}_{G} \equiv      \big[ \mathcal{Z}^{\mathrm{w}} \big( G , p , q \big) \big]^{-1}      \bigg( \frac{p}{1-p} \bigg)^{|(\omega^{\prime})^{\mathrm{w}}|} q^{k(\omega^{\prime})}                         \text{ } \text{ , } 
\end{align*}

\noindent for a suitable parameter $\omega^{\prime}$ satisfying,

\begin{align*}
  \omega^{\prime} \equiv \underset{e \in E}{\sum}  \omega^{\prime}_e  \text{ } \text{ , }
\end{align*}

\noindent $k(\omega^{\prime})$ denotes the number of connected components of the configuration $\omega^{\prime}$ appearing in the power of the cluster weight $q$,

\begin{align*}
G \supsetneq \big( \omega^{\prime} \big)^{\mathrm{w}} \equiv   \big\{        v , w \in (\omega^{\prime})^{\mathrm{w}} : \textbf{1}_{\{v, w \text{ } \mathrm{wired} \text{ } : \text{ }  v \overset{(\omega^{\prime})^{\mathrm{w}}}{\longleftrightarrow} w \}} \equiv 1   \big\}      \text{ } \text{ , }
\end{align*}

\noindent corresponding to \textit{wired}, or to \textit{free}, boundary conditions, with $\phi_G^{\mathrm{f}}[\cdot]$ instead of $\phi_G^{\mathrm{w}}[ \cdot ]$, and,

\begin{align*}
  \mathcal{Z}^{\mathrm{w}} \big( G , p , q \big) \equiv \mathcal{Z}^{\mathrm{w}} \equiv  \underset{\omega^{\prime}}{\sum} \bigg( \frac{p}{1-p} \bigg)^{|(\omega^{\prime})^{\mathrm{w}}|} q^{k(\omega^{\prime})}   \text{ } \text{ , } 
\end{align*}

\noindent denotes the normalizing constant, the partition function, under wired boundary conditions.

\bigskip

\noindent From the random-cluster probability measure $\phi$ introduced above, one can equivalently define the measure by specifying the following collection of parameters,

\[
\text{ } 
\left\{\!\begin{array}{ll@{}>{{}}l}       a_0       \text{ } \text{ , }   \\ 
       a_{\sigma}   \text{ } \text{ , } \\     a_{\tau}    \text{ } \text{ , } \\ a_{\sigma\tau}  \text{ } \text{ , } \\ q_{\sigma} \text{ } \text{ , } \\  q_{\tau}  \text{ } \text{ , }     
\end{array}\right.
\]

\noindent which will be used in the statement of the comparison inequality following the next \textbf{Definition}.

\bigskip

\noindent \textbf{Definition} \textit{28} (\textit{rescaling of random-cluster model parameters for comparison inequalities of the generalized random-cluster model} {\color{blue}[25]}). Introduce,

\[
\text{ } 
\left\{\!\begin{array}{ll@{}>{{}}l}            \alpha_0 \equiv \frac{a_0}{\widetilde{a_0}}       \text{ } \text{ , }   \\ 
   \alpha_{\sigma}   \equiv   \frac{a_{\sigma}  }{\widetilde{a_{\sigma}}} \text{ } \text{ , } \\   \alpha_{\tau} \equiv   \frac{a_{\tau} }{\widetilde{a_{\tau}}}   \text{ } \text{ , } \\ \alpha_{\sigma\tau} \equiv \frac{a_{\sigma\tau} }{\widetilde{a_{\sigma\tau}}} \text{ } \text{ , } \\    \rho_{\sigma}     \equiv \frac{q_{\sigma}}{\widetilde{q_{\sigma}}} \text{ } \text{ , } \\ \rho_{\tau} \equiv \frac{q_{\tau} }{\widetilde{q_{\tau}}} \text{ } \text{ , }     
\end{array}\right.
\]

\noindent as the collection of parameters $\alpha_0$, $\alpha_{\sigma}$, $\alpha_{\tau}$, $\alpha_{\sigma\tau}$, obtained from the  collection of parameters for the \textit{random-cluster model probability measure} $\widetilde{\phi}$.

\bigskip

\noindent \textbf{Lemma} \textit{GRCM 2} (\textit{comparison inequalities for the generalized random-cluster measure}, ({\color{blue}[25]}, \textbf{Lemma} \textit{4.2})). Suppose,

\[
\text{ } 
\left\{\!\begin{array}{ll@{}>{{}}l}           q_{\sigma} \leq \widetilde{q_{\sigma}}     \text{ } \text{ , }   \\    \alpha_{\sigma\tau}  \geq  \mathrm{max} \big\{  \alpha_{\sigma}  , \alpha_{\tau}  \big\} \geq \mathrm{min} \big\{  \alpha_{\sigma} ,  \alpha_{\tau}  \big\} \text{ } \text{ , } \text{ } \forall b \in \mathcal{B}_{>} \text{ } \text{ , }   \\       \rho_{\tau} \alpha_{\sigma} \geq \mathrm{max} \big\{       \alpha_{\sigma\tau}      ,  \rho_{\tau} \alpha_0 \big\}   \geq \mathrm{min} \big\{       \alpha_{\sigma\tau}   , \rho_{\tau} \alpha_0   \big\} \geq \alpha_{\tau}     \text{ } \text{ , } \text{ }   \forall b \in  \mathcal{B}_{<}   \text{ } \text{ , }  
\end{array}\right.
\]

\noindent hence, which implies that a comparison between the random cluster measures $\phi$ and $\widetilde{\phi}$ holds, in which,

\begin{align*}
  \phi [ A ] \geq \widetilde{\phi}[A]                    \text{ } \text{ , } 
\end{align*}

\noindent for any increasing function $A$.

\bigskip

\noindent With one example of a comparison inequality for the generalized random-cluster model, we turn towards the statement for establishing equivalence between the generalized random-cluster model and Ashkin-Teller model. 

\bigskip

\noindent \textbf{Proposition} \textit{$\mathrm{GRCM}$ 1} (\textit{parameter choice for equivalence between the generalized random-cluster model and the Ashkin-Teller model}, ({\color{blue}[25]}, \textbf{Proposition} \textit{3.1})). Let,

\[
\text{ } 
\left\{\!\begin{array}{ll@{}>{{}}l}            \alpha_0 \equiv           \mathrm{exp} \big(   -2 \big( J_{\sigma} + J_{\tau} \big)      \big)      \text{ } \text{ , }   \\ 
   \alpha_{\sigma}   \equiv   \mathrm{exp} \big( - 2   J_{\tau}    \big) \big( \mathrm{exp} \big(  - 2 J_{\sigma \tau}        \big) -    \mathrm{exp} \big(    - 2 J_{\sigma}      \big)      \big)  \text{ } \text{ , } \\   \alpha_{\tau} \equiv \mathrm{exp} \big(  - 2 J_{\sigma}     \big) \big( \mathrm{exp} \big(       -2 J_{\sigma\tau}           \big) -  \mathrm{exp} \big(      - 2 J_{\tau}         \big) \big)   \text{ } \text{ , } \\  \alpha_{\sigma\tau}     \equiv 1 - \mathrm{exp} \big( - 2 \big( J_{\sigma} + J_{\sigma \tau} \big) \big) - \mathrm{exp} \big( - 2 \big( J_{\tau} + J_{\sigma\tau} \big) \big) + \mathrm{exp} \big( -2 \big( J_{\sigma} + J_{\tau} \big) \big)   \text{ } \text{ , }  
\end{array}\right.
\]

\noindent denote the collection of parameters for which the \textit{generalized bond percolation measure},

\begin{align*}
   \lambda_{\mathcal{B}} \big( \underline{n}_{\sigma\tau} \big)  \equiv \lambda_{\mathcal{B}} \big( \sigma_A  \tau_B   \big)  \text{ } \text{ , } 
\end{align*}

\noindent provided in \textbf{Definition} \textit{20}, yields the following two probability measures, the first of which is of the form,

\begin{align*}
     \mathscr{P}^{++}_{\Lambda} \big[ \underline{n}_{\sigma\tau} \big|      \underline{p} , q      \big] \equiv  \mathscr{P}^{++} \big[ \underline{n}_{\sigma\tau} \big|       \underline{p} , q      \big] \equiv \frac{\lambda_{\mathcal{B}^{+}} \big( \underline{n}_{\sigma\tau} \big) }{\mathcal{Z}^{++}  \big( \Lambda , a_0 , a_{\sigma} , a_{\tau} , a_{\sigma \tau}  \big)  }  \equiv \frac{\underset{n_b \equiv (0,0)}{\underset{ b \in \mathcal{B}^{+}}{\prod}}          a_0(b)  \underset{n_b \equiv (1,0)}{\underset{b \in \mathcal{B}^{+}}{ \prod }}          a_{\sigma}(b)    \underset{n_b \equiv (0,1)}{\underset{b \in \mathcal{B}^{+}}{\prod} }          a_{\sigma\tau}(b)     \underset{n_b \equiv (1,1)}{\underset{b \in \mathcal{B}^{+}}{\prod }} a_{\tau}(b)    }{\mathcal{Z}^{++} \big( \Lambda , a_0 , a_{\sigma} , a_{\tau} , a_{\sigma \tau}  \big)  } \\  \equiv \frac{\lambda_{\mathcal{B}^{+} ( \underline{n}_{\sigma\tau} ) }           2^{N_{\sigma} (\underline{n}_{\sigma\tau} | \Lambda)} 2^{N_{\tau} ( \underline{n}_{\sigma\tau} | \Lambda )}     }{  \mathscr{C}_2 \underset{\underline{n}_{\sigma\tau}: {\underline{n}_{\sigma\tau}}|_{\partial \Lambda} ++ }{\underset{\underline{n} = (\underline{n}_{\sigma} , \underline{n}_{\tau} ) }{\sum}} \lambda_{\mathcal{B}^{+} ( \underline{n}_{\sigma\tau} ) }           2^{N_{\sigma} (\underline{n}_{\sigma\tau} | \Lambda)} 2^{N_{\tau} ( \underline{n}_{\sigma\tau} | \Lambda )}     }  \\ \equiv    \nu^{++}_{\Lambda}  \big[ \underline{n}_{\sigma\tau} | 2 , 2 \big]       \\ \equiv  \nu^{++} \big[ \kappa^{\sigma}_A  \kappa^{\tau}_B | 2, 2 \big]    \text{ } \text{ , } 
\end{align*}


\noindent for suitable, strictly positive $\mathscr{C}_2$, corresponding to the \textit{generalized random-cluster probability measure} under $++$ boundary conditions, in addition to, the following probability measure, the second of which is of the form,

\begin{align*}
   \mathscr{P}^{\mathrm{f}}_{\Lambda} \big[ \underline{n}_{\sigma\tau} \big|      \underline{p} , q      \big] \equiv   \mathscr{P}^{\mathrm{f}} \big[ \underline{n}_{\sigma\tau} \big|      \underline{p} , q      \big]  \equiv \frac{\lambda_{\mathcal{B}} \big( \underline{n}_{\sigma\tau} \big)}{\mathcal{Z}^{\mathrm{f}} \big( \Lambda , a_0 , a_{\sigma} , a_{\tau} , a_{\sigma\tau} \big) }  \equiv  \frac{\underset{n_b \equiv (0,0)}{\underset{ b \in \mathcal{B}}{\prod}}          a_0(b)  \underset{n_b \equiv (1,0)}{\underset{b \in \mathcal{B}}{ \prod }}          a_{\sigma}(b)    \underset{n_b \equiv (0,1)}{\underset{b \in \mathcal{B}}{\prod} }          a_{\sigma\tau}(b)     \underset{n_b \equiv (1,1)}{\underset{b \in \mathcal{B}}{\prod }} a_{\tau}(b)     }{  \mathcal{Z}^{\mathrm{f}} \big( \Lambda , a_0 , a_{\sigma} , a_{\tau} , a_{\sigma\tau} \big)   }      \\ \equiv \frac{\lambda_{\mathcal{B}(\underline{n}_{\sigma\tau})} 2^{N_{\sigma} (\underline{n}_{\sigma\tau} | \Lambda )} 2^{N_{\tau} ( \underline{n}_{\sigma\tau} | \Lambda )}   }{\mathscr{C}_3 \underset{\underline{n}_{\sigma\tau}: {\underline{n}_{\sigma\tau}}|_{\partial \Lambda} \mathrm{f} }{\underset{\underline{n} = (\underline{n}_{\sigma} , \underline{n}_{\tau} ) }{\sum}} \lambda_{\mathcal{B}(\underline{n}_{\sigma\tau})} 2^{N_{\sigma} (\underline{n}_{\sigma\tau} | \Lambda )} 2^{N_{\tau} ( \underline{n}_{\sigma\tau} | \Lambda )}   }   \\ \equiv \nu^{\mathrm{f}}_{\Lambda} \big[ \underline{n}_{\sigma\tau} | 2, 2\big] \\ \equiv  \nu^{\mathrm{f}} \big[  \underline{n}_{\sigma\tau} | 2, 2\big]  \text{ } \text{ , } 
\end{align*}

\noindent for suitable, strictly positive $\mathscr{C}_3$, corresponding to the \textit{generalized random-cluster probability measure} under free boundary conditions.

\bigskip

\noindent With the identification of parameters in the Ashkin-Teller model that coincide with parameters that need to be specified for the generalized random-cluster model, also introduce the corresponding connectivity event in the strip that is analyzed under the generalized random-cluster measure. From the \textbf{Proposition} above, as a related fact, one obtains the following equalities directly relating the \textit{generalized bond percolation} measure appearing in the \textit{generalized random-cluster model measure} with the \textit{Ashkin-Teller measure} of \textbf{Definition} \textit{14}, in which,

\begin{align*}
    {\nu}^{+}_{\Lambda} \big[ i \overset{\sigma}{\longleftrightarrow}    \Lambda^c \big] \equiv  \mathcal{P}^{++} \big[      \sigma_i          \big]       \text{ } \text{ , } \\    {\nu}^{+}_{\Lambda} \big[ i \overset{\sigma}{\longleftrightarrow} j  \big]              \equiv       \mathcal{P}^{++} \big[      \sigma_i        \sigma_j   \big]     \text{ } \text{ , } 
\end{align*}

\noindent in which, from each equality provided above, the probability, under the \textit{Ashkin-Teller measure}, of sampling $\sigma_i$, for $i \in \Lambda$ under $++$ boundary conditions, is equal to the probability of the connectivity event $\big\{ i \longleftrightarrow \Lambda^c \big\}$ occurring under the \textit{generalized random-cluster measure}. In the second equality, the spin-spin correlation $\sigma_i \sigma_j$, for $j \in \Lambda$ under $++$ boundary conditions, is equal to the probability of the connectivity event $\big\{ i \longleftrightarrow j \big\}$ occurring under the \textit{generalized random-cluster measure}. 

\bigskip

\noindent From the two special cases above of \textbf{Proposition} \textit{$\mathrm{GRCM}$ 1}, it is also necessary to introduce reformulations of the crossing events previously analyzed under the \textit{Ashkin-Teller measure}, in order to establish equality with connectivity events under the \textit{generalized random-cluster measure} under the special cases of  \textbf{Proposition} \textit{GRCM 1} above, with the following items:

\begin{itemize}

    \item[$\bullet$] \textit{Segment connectivity event within the strip}. From (Event 1),
    
    \begin{align*}
          \mathcal{P}^{\xi}_{\Lambda} \big[     [ 0, \lfloor \delta^{\prime\prime} n \rfloor ] \times \{ 0 \} \underset{F^{\mathrm{odd}}}{\longleftrightarrow}  [i , i+ \lfloor \delta^{\prime\prime}n \rfloor ]  \times \{n \}  \big]    \approx        \mathcal{P}^{\xi}_{\Lambda} \big[  \sigma_0  \sigma_{\lfloor \delta^{\prime\prime} n \rfloor} \sigma_i  \sigma_{i + \lfloor \delta^{\prime\prime} n \rfloor}      \big]              \tag{Event $1^{\prime}$} \\         \overset{(\textbf{Proposition} \text{ } \textit{$\mathrm{GRCM}$ 1})}{=}                \nu^{+}_{\Lambda}  \big[              [ 0, \lfloor \delta^{\prime\prime} n \rfloor ] \times \{ 0 \} \overset{\sigma}{\longleftrightarrow}   [i , i+ \lfloor \delta^{\prime\prime}n \rfloor ]  \times \{n \}             \big]      \text{ } \text{ , }  
    \end{align*}

    \noindent with,
  
  \begin{align*}
   \sigma_0 \cap \sigma_{\lfloor \delta^{\prime\prime} n \rfloor} \cap [0 ,  \lfloor \delta^{\prime\prime} n \rfloor ] \neq \emptyset     \text{ } \text{ , } 
  \end{align*}

    \noindent in addition to,
    
     \begin{align*}
       \sigma_i \cap \sigma_{i+\lfloor \delta^{\prime\prime} n \rfloor} \cap [i , i +  \lfloor \delta^{\prime\prime} n \rfloor ] \neq \emptyset         \text{ } \text{ . } 
  \end{align*}

    \item[$\bullet$] \textit{Segment connectivity event between the top and bottom strip boundaries}. From (Event 2),
    
    \begin{align*}
     \bigg\{   \mathcal{I}_0 \underset{F^{\mathrm{even}}}{\overset{+ \backslash -}{\longleftrightarrow} } \widetilde{\mathcal{I}_0}      \bigg\}  \equiv  \big\{ \mathcal{I}_0     \longleftrightarrow \widetilde{\mathcal{I}_0}  \big\}  \approx \big\{   \sigma_0 \widetilde{\sigma_0}                             \big\}  \text{ } \text{ , } \tag{Event $2^{\prime}$} \\         \overset{(\textbf{Proposition} \text{ } \textit{$\mathrm{GRCM}$ 1})}{=}                \nu^{+}_{\Lambda}  \big[   \mathcal{I}_0  \overset{\sigma}{\longleftrightarrow}  \widetilde{\mathcal{I}_0}      \big]   \text{ } \text{ , }
    \end{align*}
    
    \noindent with,

    \begin{align*}
        \sigma_0 \cap \mathcal{I}_0 \neq \emptyset      \text{ } \text{ , } 
    \end{align*}

    \noindent in addition to,
    
        \begin{align*}
          \widetilde{\sigma_0} \cap \widetilde{\mathcal{I}_0} \neq \emptyset     \text{ } \text{ . } 
    \end{align*}

    \item[$\bullet$] \textit{Intersection of three segment connectivity events, \textbf{Lemma} \textit{AT 2.1}}. From (Event 3),
    
    \begin{align*}
           \mathcal{P}^{\chi}_{\Lambda} \bigg[  \big\{ \mathcal{I}_{-\mathscr{I}} {\longleftrightarrow} \widetilde{\mathcal{I}_{-\mathscr{I}}} \big\}  \cap    \big\{   \mathcal{I}_{-\mathscr{J}} {\longleftrightarrow}     \widetilde{\mathcal{I}_{-\mathscr{J}}}     \big\}   \cap       \big\{          \mathcal{I}_{-\mathscr{K}} {\longleftrightarrow}     \widetilde{\mathcal{I}_{-\mathscr{K}}}           \big\}            \bigg]  \approx     \mathcal{P}^{\chi}_{\Lambda}  \big[  \sigma_{-\mathscr{I}} \widetilde{\sigma_{-\mathscr{I}}} \sigma_{-\mathscr{J}}   \widetilde{\sigma_{-\mathscr{J}}}  \sigma_{-\mathscr{K}}   \widetilde{\sigma_{-\mathscr{K}}}   \big]              \text{ } \text{ , }     \tag{Event $3^{\prime}$} \\         \overset{(\textbf{Proposition} \text{ } \textit{$\mathrm{GRCM}$ 1})}{=}                \nu^{+}_{\Lambda}  \bigg[   \big\{ \mathcal{I}_{-\mathscr{I}} \overset{\sigma}{\longleftrightarrow}   \widetilde{\mathcal{I}_{-\mathscr{I}}} \big\}  \cap    \big\{   \mathcal{I}_{-\mathscr{J}} \overset{\sigma}{\longleftrightarrow}   \widetilde{\mathcal{I}_{-\mathscr{J}}}     \big\}   \cap       \big\{          \mathcal{I}_{-\mathscr{K}} \overset{\sigma}{\longleftrightarrow}   \widetilde{\mathcal{I}_{-\mathscr{K}}}           \big\}             \bigg]   \text{ } \text{ , }
    \end{align*}

    \noindent with,

    \begin{align*}
        \sigma_{-\mathscr{I}} \cap \mathcal{I}_{-\mathscr{I}} \neq \emptyset      \text{ } \text{ , } \\   \sigma_{-\mathscr{I}} \cap \widetilde{\mathcal{I}_{-\mathscr{I}}} \neq \emptyset          \text{ } \text{ , } \\       \sigma_{-\mathscr{J}} \cap \mathcal{I}_{-\mathscr{J}}      \neq \emptyset             \text{ } \text{ , } \\            \sigma_{-\mathscr{J}} \cap \widetilde{\mathcal{I}_{-\mathscr{J}} }      \neq \emptyset       \text{ } \text{ , } \\   \sigma_{-\mathscr{K}} \cap \mathcal{I}_{-\mathscr{K}}                \neq \emptyset       \text{ } \text{ , }  \\   \sigma_{-\mathscr{K}} \cap \widetilde{\mathcal{I}_{-\mathscr{K}} }  \neq \emptyset       \text{ } \text{ . } 
    \end{align*}

    \item[$\bullet$] \textit{Two segment connectivity events upper bounded by $1-c$, \textbf{Lemma} \textit{AT 2.1}}. From (Event 4),

    \begin{align*}
       \underset{ \overset{(\textbf{Proposition} \text{ } \textit{$\mathrm{GRCM}$ 1})}{=}                \nu^{+}_{\Lambda}  \big[                   [ 0, \lfloor \delta^{\prime\prime} n \rfloor ] \times \{ 0 \} \overset{\sigma}{\longleftrightarrow}   [i , i+ \lfloor \delta^{\prime\prime}n \rfloor  ]  \times \{n \}       \big]  }{\underset{\approx \mathcal{P}^{\chi}_{\Lambda} \big[ \sigma_0  \sigma_{\lfloor \delta^{\prime\prime} n \rfloor} \sigma_i \sigma_{i+ \lfloor \delta^{\prime\prime} n \rfloor}    \big] }{{\underbrace{\mathcal{P}^{\chi}_{\Lambda} \big[             [ 0, \lfloor \delta^{\prime\prime} n \rfloor ] \times \{ 0 \} {\longleftrightarrow}  [i , i+ \lfloor \delta^{\prime\prime}n \rfloor  ]  \times \{n \}                     \big] }}}}  \leq            \underset{ \overset{(\textbf{Proposition} \text{ } \textit{$\mathrm{GRCM}$ 1})}{=}                \nu^{+}_{\Lambda}  \big[              [ 0, \lfloor \delta^{\prime\prime} n \rfloor ] \times \{ 0 \}    \overset{\sigma}{\longleftrightarrow}   [i , i+ \lfloor \delta^{\prime\prime}n \rfloor  ]  \times \{n \}                      \big]  }{\underset{\approx \mathcal{P}^{\chi}_{\Lambda} \big[ \sigma_0   \sigma_{\lfloor \delta n \rfloor} \sigma_i   \sigma_{i+ \lfloor \delta^{\prime\prime}n \rfloor} \big] }{{\underbrace{\mathcal{P}^{\chi}_{\Lambda} \big[             [ 0, \lfloor \delta n \rfloor ] \times \{ 0 \} {\longleftrightarrow}  [i , i+ \lfloor \delta  n \rfloor ]  \times \{n \}                     \big] }}}}  \text{ } \text{ , } \tag{Event $4^{\prime}$}          
    \end{align*}
    
    \noindent with,
    
    \begin{align*}
        \sigma_0 \cap \sigma_{\lfloor \delta^{\prime\prime} n \rfloor}    \cap   [0 , \lfloor \delta^{\prime\prime} n \rfloor ]     \neq \emptyset   \text{ } \text{ , } \\   \sigma_i \cap \sigma_{i+\lfloor \delta^{\prime\prime} n \rfloor}    \cap   [i , i+ \lfloor \delta^{\prime\prime} n \rfloor ]     \neq \emptyset   \text{ } \text{ , } \\          \sigma_0 \cap \sigma_{\lfloor \delta n \rfloor} \cap [0 , \lfloor \delta n \rfloor ]       \neq \emptyset   \text{ } \text{ , }  \\     \sigma_i \cap \sigma_{i+\lfloor \delta n \rfloor} \cap [i , i +  \lfloor \delta n \rfloor ]          \neq \emptyset    \text{ } \text{ . } 
    \end{align*}

       \item[$\bullet$] \textit{Two other segment connectivity events, \textbf{Lemma} \textit{AT 2.1}}. From (Event 5),

       \begin{align*}
      \underset{ \overset{(\textbf{Proposition} \text{ } \textit{$\mathrm{GRCM}$ 1})}{=}                \nu^{+}_{\Lambda}  \big[               [ \lfloor \delta^{\prime\prime} n \rfloor, 2 \lfloor \delta^{\prime\prime} n \rfloor ] \times \{ 0 \}        \overset{\sigma}{\longleftrightarrow}   [ - i + \lfloor \delta^{\prime\prime} n \rfloor , i+ 2  \lfloor \delta^{\prime\prime}n \rfloor ]  \times \{n \}                       \big]  }{\underset{\approx \mathcal{P}^{\chi}_{\Lambda} \big[  \sigma_{\lfloor \delta^{\prime\prime} n \rfloor} \sigma_{2 \lfloor \delta^{\prime\prime} n \rfloor }   \sigma_{-i +  \lfloor \delta^{\prime\prime} n \rfloor }   \sigma_{i +  2 \lfloor \delta^{\prime\prime} n \rfloor }          \big] }{{\underbrace{\mathcal{P}^{\chi}_{\Lambda} \big[               [ \lfloor \delta^{\prime\prime} n \rfloor, 2 \lfloor \delta^{\prime\prime} n \rfloor ] \times \{ 0 \} {\longleftrightarrow}  [ - i + \lfloor \delta^{\prime\prime} n \rfloor , i+ 2  \lfloor \delta^{\prime\prime}n \rfloor ]  \times \{n \}              \big] }}}}  \\ \geq       \underset{ \overset{(\textbf{Proposition} \text{ } \textit{$\mathrm{GRCM}$ 1})}{=}                \nu^{+}_{\Lambda}  \big[          [ 0 ,  \lfloor \delta^{\prime\prime} n \rfloor ] \times \{ 0 \} \overset{\sigma}{\longleftrightarrow}     [  i + \lfloor \delta^{\prime\prime} n \rfloor , i+   2\lfloor \delta^{\prime\prime}n \rfloor ]  \times \{n \}         \big]  }{\underset{\approx \mathcal{P}^{\chi}_{\Lambda} \big[  \sigma_0 \sigma_{\lfloor \delta^{\prime\prime} n \rfloor} \sigma_{i + \lfloor \delta^{\prime\prime} n \rfloor}  \sigma_{i + 2 \lfloor \delta^{\prime\prime} n \rfloor}     \big] }{{\underbrace{\mathcal{P}^{\chi}_{\Lambda} \big[        [ 0 ,  \lfloor \delta^{\prime\prime} n \rfloor ] \times \{ 0 \} {\longleftrightarrow}  [  i + \lfloor \delta^{\prime\prime} n \rfloor , i+   2\lfloor \delta^{\prime\prime}n \rfloor ]  \times \{n \}                   \big]}}}} \text{ } \text{ , }    \tag{Event $5^{\prime}$}
       \end{align*}

    \noindent with,

    \begin{align*}
       \sigma_{\lfloor \delta^{\prime\prime} n \rfloor}    \cap    \sigma_{2 \lfloor \delta^{\prime\prime} n \rfloor}   \cap      [ \lfloor \delta^{\prime\prime} n \rfloor, 2 \lfloor \delta^{\prime\prime} n \rfloor ]          \neq \emptyset                \text{ } \text{ , }  \\                         \sigma_{i+\lfloor \delta^{\prime\prime} n \rfloor}    \cap    \sigma_{i+2 \lfloor \delta^{\prime\prime} n \rfloor}   \cap      [i+ \lfloor \delta^{\prime\prime} n \rfloor, i+ 2 \lfloor \delta^{\prime\prime} n \rfloor ]          \neq \emptyset   \text{ } \text{ , } \\                  \sigma_0 \cap \sigma_{\lfloor \delta^{\prime\prime} n \rfloor}        \cap [0 , \lfloor \delta^{\prime\prime} n \rfloor ]        \neq \emptyset    \text{ } \text{ , } \\          \sigma_{i + \lfloor \delta^{\prime\prime} n \rfloor} \cap \sigma_{i + 2 \lfloor \delta^{\prime\prime} n \rfloor} \cap     [ i + \lfloor \delta^{\prime\prime} n \rfloor , i + 2 \lfloor \delta^{\prime\prime} n \rfloor]         \neq \emptyset    \text{ } \text{ . }
    \end{align*}

    \item[$\bullet$] \textit{Connectivity events from the left symmetric domain boundary, \textbf{Lemma} \textit{6.3}}. From (Event $6$),
    
    \begin{align*}
     \mathcal{P}^{\chi_2}_{\Lambda} \big[      \gamma_L \overset{+ \backslash -}{\underset{\Lambda \cap (\mathscr{F}\mathscr{C}_1)^c}{\longleftrightarrow}}    \mathscr{F}\mathscr{C}_1          \big] -                \mathcal{P}^{\chi_1}_{\Lambda} \big[         \gamma_L \overset{+ \backslash -}{\underset{\Lambda\cap (\mathscr{F}\mathscr{C}_1)^c}{\longleftrightarrow}}    \mathscr{F}\mathscr{C}_1            \big] \approx     \mathcal{P}^{\chi_2}_{\Lambda} \big[                  \sigma_{\gamma_L} \sigma_{\mathscr{F}\mathscr{C}_1}     \big] -       \mathcal{P}^{\chi_1}_{\Lambda} \big[          \sigma_{\gamma_L} \sigma_{\mathscr{F}\mathscr{C}_1}             \big]         \text{ } \text{ , }     \tag{Event $6^{\prime}$} \\  \overset{(\textbf{Proposition} \text{ } \textit{$\mathrm{GRCM}$ 1})}{\approx }                \nu^{\chi_1}_{\Lambda}  \big[  \gamma_L \overset{\sigma}{\longleftrightarrow}     \mathscr{F}\mathscr{C}_1      \big]  - \nu^{\chi_2}_{\Lambda} \big[\gamma_L \overset{\sigma}{\longleftrightarrow}    \mathscr{F}\mathscr{C}_1 \big]   \text{ } \text{ , }
    \end{align*}

    \noindent for $\chi_1 , \chi_2 \in \textbf{B} \textbf{C}_{++}$, with,

   \begin{align*}
       \sigma_{\gamma_L} \cap \gamma_L \neq \emptyset     \text{ } \text{ , }       \\    \sigma_{\mathscr{F}\mathscr{C}_1} \cap   \mathscr{F}\mathscr{C}_1           \neq \emptyset \text{ } \text{ . } 
   \end{align*}

     \item[$\bullet$] \textit{Decomposition of crossing events from (Event $6$), \textbf{Lemma} \textit{6.3}}. From (Event $7$),

    \begin{align*}
      \mathcal{P}^{\chi_2}_{\Lambda} \big[     \gamma_L        \underset{\Lambda \cap (\mathscr{F}\mathscr{C} )_1}{\overset{-}{\longleftrightarrow}}                 \gamma_{-}                   \big] \text{ }  +  \text{ }    \mathcal{P}^{\chi_2}_{\Lambda} \big[       \big(  \gamma_{-,+}   \big)_1        \underset{\Lambda \cap (\mathscr{F}\mathscr{C})_1}{\overset{+ \cup -}{\longleftrightarrow}}                        \gamma_{L^{\prime}}             \big] \text{ } +   \text{ }    \mathcal{P}^{\chi_2}_{\Lambda} \big[     \big(   \gamma_{-,+}  \big)_2          \underset{\Lambda \cap (\mathscr{F}\mathscr{C})_1}{\overset{+}{\longleftrightarrow}}                     \mathscr{F}\mathscr{C}_1                \big] \\ \approx              \mathcal{P}^{\chi_2}_{\Lambda} \big[            \sigma_L \sigma_{-}               \big]  +    \mathcal{P}^{\chi_2}_{\Lambda} \big[     \big( \sigma_{-,+}\big)_1  \sigma_{L^{\prime}}                          \big]  +  \mathcal{P}^{\chi_2}_{\Lambda} \big[                 \big( \sigma_{-,+} \big)_2 \sigma_{\mathscr{F}\mathscr{C}_1 }            \big]    \text{ } \tag{Event $7^{\prime}$} \\   \overset{(\textbf{Proposition} \text{ } \textit{$\mathrm{GRCM}$ 1})}{=}                \nu^{+}_{\Lambda}  \big[   \gamma_L  \overset{\sigma}{\longleftrightarrow}   \gamma_{-}  \big] + \nu^{+}_{\Lambda}  \big[   \big( \gamma_{-,+} \big)_1 \overset{\sigma}{\longleftrightarrow}     \gamma_{L^{\prime}}    \big] + \nu^{+}_{\Lambda}  \big[    \big( \gamma_{-,+} \big)_2 \overset{\sigma}{\longleftrightarrow}    \mathscr{F} \mathscr{C}_2    \big]  \text{ } \text{ , }
    \end{align*}

    \noindent with,
  
    \begin{align*}
          \sigma_L   \cap \gamma_L \neq \emptyset                              \text{ } \text{ , }      \\  \sigma_{-}  \cap         \gamma_{-} \neq \emptyset    \text{ } \text{ , }      \\      \big( \sigma_{-,+} \big)_1            \cap    \big( \gamma_{-,+} \big)_1      \neq \emptyset  \text{ } \text{ , }      \\ \sigma_{L^{\prime}}           \cap     \gamma_{L^{\prime}}       \neq \emptyset         \text{ } \text{ , } \\     \big( \sigma_{+ , - }\big)_2   \cap \big( \gamma_{-,+} \big)_2   \neq \emptyset         \text{ } \text{ , } \\  \sigma_{\mathscr{F}\mathscr{C}_1} \cap \mathscr{F}\mathscr{C}_1    \neq \emptyset      \text{ } \text{ . }
    \end{align*}

    \item[$\bullet$] \textit{$\chi_1$-$\chi_2$ \textit{ratio}}. From (Event $8$),
    
    \begin{align*}
        \frac{\underset{\chi_1 \in \{ + , + \}}{\sum}            \mathcal{P}^{\chi_1}_{\Lambda} \bigg[    \gamma_L        \underset{\Lambda \cap (\mathscr{F}\mathscr{C} )_1}{\overset{-}{\longleftrightarrow}}                 \gamma_{-}        \big|                       \mathcal{C}_{-}   \cup  \mathscr{D}_{\chi,-}   \bigg]      \mathcal{P}^{\chi_1}_{\Lambda} \big[                   \mathcal{C}_{-}   \cup  \mathscr{D}_{\chi,-}      \big]          }{\underset{\chi_2 \in \{ + , + \}}{\sum}     \mathcal{P}^{\chi_2}_{\Lambda} \bigg[    \gamma_L        \underset{\Lambda \cap (\mathscr{F}\mathscr{C} )_1}{\overset{-}{\longleftrightarrow}}                 \gamma_{-}        \big|      \mathcal{C}_{-}   \cup  \mathscr{D}_{\chi,-}        \bigg]      \mathcal{P}^{\chi_2}_{\Lambda} \big[                    \mathcal{C}_{-}   \cup  \mathscr{D}_{\chi,-}     \big]   }        \\ \approx \frac{\underset{\chi_1 \in \{ + , + \}}{\sum}      \mathcal{P}^{\chi_1}_{\Lambda}  \big[       \sigma_L \sigma_{-}        |  \sigma_{j-1,-} \sigma_{\mathscr{L}_j} \sigma_{\chi , -}    \big]    \mathcal{P}^{\chi_1}_{\Lambda}  \big[     \sigma_{j-1,-}    \big]   }{\underset{\chi_2 \in \{ + , + \}}{\sum} \mathcal{P}^{\chi_2}_{\Lambda} \big[ \sigma_L \sigma_{-}    |   \sigma_{j-1,-} \sigma_{\mathscr{L}_j}  \sigma_{\chi , -}     \big] \mathcal{P}^{\chi_2}_{\Lambda} \big[   \sigma_{j-1,-}  \big]  }                         \text{ } \text{ , }   \tag{Event $8^{\prime}$}     \\   \approx \frac{\underset{\chi_1 \in \{ +, + \}}{\sum}      \mathcal{P}^{\chi_1}_{\Lambda}  \big[       \sigma_L \sigma_{-}        |  \sigma_{j-1,-} \sigma_{\mathscr{L}_j} \sigma_{\chi , -}    \big]    \mathcal{P}^{\chi_1}_{\Lambda}  \big[    \sigma_{j-1,-}     \big]   }{\underset{\chi_2 \in \{ +, + \}}{\sum} \mathcal{P}^{\chi_2}_{\Lambda} \big[ \sigma_L \sigma_{-}    |   \sigma_{j-1,-} \sigma_{\mathscr{L}_j}  \sigma_{\chi , -}     \big] \mathcal{P}^{\chi_2}_{\Lambda} \big[   \sigma_{j-1,-}  \big]  }   \\ \overset{(\textbf{Proposition} \text{ } \textit{$\mathrm{GRCM}$ 1})}{=}                \frac{\underset{\chi \in \{ +, + \}}{\sum} \nu^{\chi_1}_{\Lambda}  \bigg[ \gamma_L \overset{\sigma}{\longleftrightarrow}  \gamma_{-}    \big|    \big\{     \gamma_{j-1,-}     \overset{\sigma}{\longleftrightarrow} \gamma_{-}          \big\} \cup \big\{     \mathscr{L}_{j,-} \overset{\sigma}{\longleftrightarrow}       \big( \Lambda \cap \big( \mathscr{F} \mathscr{C} \big)^c \big)                \big\}       \bigg] }{\underset{\chi \in \{ +, + \}}{\sum} \nu^{\chi_2}_{\Lambda}  \bigg[ \gamma_L \overset{\sigma}{\longleftrightarrow}  \gamma_{-}  \big|     \big\{     \gamma_{j-1,-}     \overset{\sigma}{\longleftrightarrow} \gamma_{-}          \big\} \cup \big\{     \mathscr{L}_{j,-} \overset{\sigma}{\longleftrightarrow}       \big( \Lambda \cap \big( \mathscr{F} \mathscr{C} \big)^c \big)                \big\}       \bigg]   } \\ \times \frac{\nu^{\chi_1}_{\Lambda}  \big[   \big\{     \gamma_{j-1,-}     \overset{\sigma}{\longleftrightarrow} \gamma_{-}          \big\} \cup \big\{     \mathscr{L}_{j,-} \overset{\sigma}{\longleftrightarrow}       \big( \Lambda \cap \big( \mathscr{F} \mathscr{C} \big)^c \big)                \big\}       \big] }{\nu^{\chi_2}_{\Lambda}  \big[   \big\{     \gamma_{j-1,-}     \overset{\sigma}{\longleftrightarrow} \gamma_{-}          \big\} \cup \big\{     \mathscr{L}_{j,-} \overset{\sigma}{\longleftrightarrow}       \big( \Lambda \cap \big( \mathscr{F} \mathscr{C} \big)^c \big)                \big\}       \big] }  \text{ } \text{ , } 
    \end{align*}

    \noindent with,

    \begin{align*}
     \sigma_L     \cap     \gamma_L     \neq \emptyset    \text{ } \text{ , } \\           \sigma_{-} \cap   \gamma_{-}     \neq \emptyset \text{ } \text{ , } \\  \sigma_{j-1}   \cap \sigma_{-}  \cap  \mathcal{C}_{-} \neq \emptyset  \text{ } \text{ , }  \\ \sigma_{\chi , -} \cap \sigma_{-} \cap \mathscr{D}_{\chi,-} \neq \emptyset   \text{ } \text{ . } 
    \end{align*}

\item[$\bullet$] \textit{Conditional probability of crossing between $\gamma_L$ and $\gamma_{-}$}. From (Event $9$),

    \begin{align*}
    \mathcal{P}^{\chi_2 \backslash \chi_1}_{\Lambda} \bigg[                  \gamma_L        \underset{\Lambda \cap (\mathscr{F}\mathscr{C})_1}{\overset{-}{\longleftrightarrow}}                 \gamma_{-}  \text{ }       \big|      \text{ }                  \mathcal{C}_{-}   \cup  \mathscr{D}_{\chi,-}        \bigg]   \approx   \mathcal{P}^{\chi_2 \backslash  \chi_1}_{\Lambda} \big[  \sigma_L \sigma_{-} |  \sigma_{j-1,-} \sigma_{-}    \sigma_{\mathscr{L}_j} \sigma_{-}   \big] \text{ } \text{ , }  \tag{Event $9^{\prime}$} \\  \overset{(\textbf{Proposition} \text{ } \textit{$\mathrm{GRCM}$ 1})}{=}                \nu^{+ \backslash   -}_{\Lambda}  \big[     \gamma_L        \overset{\sigma}{\longleftrightarrow}  \gamma_{-}     \big| \big\{     \gamma_{j-1,-}     \overset{\sigma}{\longleftrightarrow} \gamma_{-}          \big\} \cup \big\{     \mathscr{L}_{j,-} \overset{\sigma}{\longleftrightarrow}       \big( \Lambda \cap \big( \mathscr{F} \mathscr{C} \big)^c \big)                \big\}        \big]   \text{ } \text{ , }
    \end{align*}

    \noindent with,
   
   \begin{align*}
        \sigma_L \cap \gamma_L \neq \emptyset     \text{ } \text{ , } \\ \sigma_{-}   \cap \gamma_{-} \neq \emptyset     \text{ } \text{ , }  \\   \sigma_{j-1,-} \cap \sigma_{-} \cap \mathcal{C}_{-}  \neq \emptyset \text{ } \text{ , } \\         \sigma_{\chi,-} \cap \sigma_{-}  \cap   \mathscr{D}_{\chi,-}   \neq \emptyset       \text{ } \text{ . }   
   \end{align*}

    \item[$\bullet$] \textit{$\chi_1$-$\chi_2$ \textit{ratio 2}}. From (Event $10$),
    
    \begin{align*}
   \underset{\chi_1 \in \{ + , + \}}{\underset{\emptyset \neq\mathscr{F} \in F ( \Lambda \cap ( \mathscr{F} \mathscr{C})_1 )}{\sum}}         \mathcal{P}^{\chi_1}_{\Lambda} \bigg[     \gamma_L \underset{\Lambda \cap ( \mathscr{F}\mathscr{C})_1}{\overset{-}{\longleftrightarrow}}         \mathscr{F}  \big| \mathcal{C}_{-}  \cup   \mathscr{D}_{\chi_3 , -}     \bigg]   \mathcal{P}^{\chi_1}_{\Lambda} \big[  \mathcal{C}_{-}  \cup   \mathscr{D}_{\chi_3 , -}           \big] \approx   \underset{ \sigma \cap \mathscr{F} \neq \emptyset}{\underset{\chi_1 \in \{ + , + \}}{\underset{\emptyset \neq\mathscr{F} \in F ( \Lambda \cap ( \mathscr{F} \mathscr{C})_1 )}{\sum}}}  \mathcal{P}^{\chi_1}_{\Lambda} \big[    \sigma_L \sigma     | \cdots \\ \sigma_{j-1,-} \sigma_{-}    \sigma_{\mathscr{L}_j} \sigma_{-}    \big]        \mathcal{P}^{\chi_1}_{\Lambda} \big[    \sigma_{j-1,-} \sigma_{-}    \sigma_{\mathscr{L}_j} \sigma_{-}      \big]              \text{ } \text{ , } \tag{Event $10^{\prime}$} \\  \overset{(\textbf{Proposition} \text{ } \textit{$\mathrm{GRCM}$ 1})}{=}                \nu^{+}_{\Lambda}  \bigg[   \gamma_L \overset{\sigma}{\longleftrightarrow}     \mathscr{F}         \big| \big\{     \gamma_{j-1,-}     \overset{\sigma}{\longleftrightarrow} \gamma_{-}          \big\} \cup \big\{     \mathscr{L}_{j,-} \overset{\sigma}{\longleftrightarrow}       \big( \Lambda \cap \big( \mathscr{F} \mathscr{C} \big)^c \big)                \big\}         \bigg]  \\ \times     \nu^{+}_{\Lambda}  \bigg[  \big\{     \gamma_{j-1,-}     \overset{\sigma}{\longleftrightarrow} \gamma_{-}          \big\} \cup \big\{     \mathscr{L}_{j,-} \overset{\sigma}{\longleftrightarrow}       \big( \Lambda \cap \big( \mathscr{F} \mathscr{C} \big)^c \big)                \big\}     \bigg]             \text{ } \text{ , } 
    \end{align*}

    \noindent with,

   \begin{align*}
           \sigma_L \cap \gamma_L \neq \emptyset \text{ } \text{ , } \\      \sigma_{j-1,-} \cap \sigma_{-}  \cap    \mathcal{C}_{-}   \neq \emptyset \text{ } \text{ , } \\     \sigma_{\mathscr{L}_j} \cap \sigma_{-}        \cap    \mathscr{D}_{\chi_3,-}    \neq \emptyset   \text{ } \text{ . }
   \end{align*}

       \item[$\bullet$] \textit{Conditional probability from crossing beginning at the left symmetric domain boundary}. From (Event $11$),
       
       \begin{align*}
      \mathcal{P}^{\chi_3}_{\Lambda} \bigg[     \gamma_L \underset{\Lambda \cap ( \mathscr{F}\mathscr{C})_1}{\overset{-}{\longleftrightarrow}}         \mathscr{F}_k  \big| \mathcal{C}_{-}  \cup   \mathscr{D}_{\chi_3 , -}     \bigg]   \approx      \mathcal{P}^{\chi_3}_{\Lambda} \big[              \sigma_L \sigma_{k}                     |   \sigma_{j-1,-} \sigma_{-}    \sigma_{\mathscr{L}_j} \sigma_{-}   \big]          \text{ } \text{ , }         \tag{Event $11^{\prime}$} \\ \overset{(\textbf{Proposition} \text{ } \textit{$\mathrm{GRCM}$ 1})}{=}                \nu^{+}_{\Lambda}  \bigg[  \gamma_L        \overset{\sigma}{\longleftrightarrow}        \mathscr{F}_k      \big| \big\{     \gamma_{j-1,-}     \overset{\sigma}{\longleftrightarrow} \gamma_{-}          \big\} \cup \big\{     \mathscr{L}_{j,-} \overset{\sigma}{\longleftrightarrow}       \big( \Lambda \cap \big( \mathscr{F} \mathscr{C} \big)^c \big)                \big\}        \bigg]    \text{ } \text{ , }
       \end{align*}
       
       \noindent with,
       
         \begin{align*}
         \sigma_L \cap \gamma_L \neq \emptyset   \text{ } \text{ , } \\    \sigma_k     \cap   \mathscr{F}_k  \neq \emptyset \text{ } \text{ , } \\         \sigma_{j-1,-}          \cap \sigma_{-} \cap \mathcal{C}_{-}       \neq \emptyset                \text{ } \text{ , } \\        \sigma_{\mathscr{L}_j}           \cap       \sigma_{-} \cap    \mathscr{D}_{\chi_3,-}      \neq \emptyset                \text{ } \text{ . }
    \end{align*}

    \item[$\bullet$] \textit{$\mathcal{P}_1$}. From (Event $12$),

    \begin{align*}
    \mathcal{P}^{\chi_2}_{\Lambda} \big[     \gamma_L        \underset{\Lambda \cap (\mathscr{F}\mathscr{C})_1}{\overset{-}{\longleftrightarrow}}                 \gamma_{-}                   \big]   -  \mathcal{P}^{\chi_1}_{\Lambda} \big[     \gamma_L        \underset{\Lambda \cap (\mathscr{F}\mathscr{C})_1}{\overset{-}{\longleftrightarrow}}                 \gamma_{-}                   \big] \approx  \mathcal{P}^{\chi_2}_{\Lambda} \big[   \sigma_{L} \sigma_{-}      \big] -    \mathcal{P}^{\chi_1}_{\Lambda} \big[  \sigma_L \sigma_{-}  \big] \\ \equiv \mathcal{P}^{\chi_1}_{\Lambda} \big[  \sigma_L \sigma_{-}  \big] \bigg[  1 -     \frac{\mathcal{P}^{\chi_2}_{\Lambda} \big[   \sigma_{L} \sigma_{-}      \big]  }{\mathcal{P}^{\chi_1}_{\Lambda} \big[  \sigma_L \sigma_{-}  \big]}        \bigg]   \text{ } \text{ , }       \tag{Event $12^{\prime}$}\\ \overset{(\textbf{Proposition} \text{ } \textit{$\mathrm{GRCM}$ 1})}{=}                \nu^{+}_{\Lambda}  \big[    \sigma_L \sigma_{-}      \big]  \bigg[  1 -      \frac{\nu^{+}_{\Lambda}  \big[     \sigma_L \sigma_{-}        \big] }{\nu^{-}_{\Lambda}  \big[     \sigma_L \sigma_{-}        \big] } \bigg]    \text{ } \text{ , }
    \end{align*}
    
    \noindent with,

   \begin{align*}
       \sigma_L \cap \gamma_L \neq \emptyset        \text{ } \text{ , } \\ \sigma_{-} \cap \gamma_{-} \neq \emptyset  \text{ } \text{ . } 
   \end{align*}
   
   \item[$\bullet$] \textit{$\mathcal{P}_2$}. From (Event $13$),
   
   \begin{align*}
    \mathcal{P}^{\chi_2}_{\Lambda} \big[       \big(  \gamma_{-,+}   \big)_1        \underset{\Lambda \cap (\mathscr{F}\mathscr{C})_1}{\overset{+ \cup -}{\longleftrightarrow}}                        \gamma_{L^{\prime}}             \big]               -    \mathcal{P}^{\chi_1}_{\Lambda} \big[       \big(  \gamma_{-,+}   \big)_1        \underset{\Lambda \cap (\mathscr{F}\mathscr{C})_1}{\overset{+ \cup -}{\longleftrightarrow}}                        \gamma_{L^{\prime}}             \big]     \approx  \mathcal{P}^{\chi_2}_{\Lambda} \big[ \big(  \sigma_{-,+}\big)_1  \sigma_{L^{\prime}}         \big] -  \mathcal{P}^{\chi_1}_{\Lambda} \big[   \big(  \sigma_{-,+}\big)_1  \sigma_{L^{\prime}}           \big] \\ \equiv         \mathcal{P}^{\chi_2}_{\Lambda} \big[ \big(  \sigma_{-,+}\big)_1  \sigma_{L^{\prime}}         \big]    \bigg[  1  -  \frac{\mathcal{P}^{\chi_1}_{\Lambda} \big[   \big(  \sigma_{-,+}\big)_1  \sigma_{L^{\prime}}           \big]}{\mathcal{P}^{\chi_2}_{\Lambda} \big[ \big(  \sigma_{-,+}\big)_1  \sigma_{L^{\prime}}         \big]}           \bigg]         \text{ } \text{ , }    \tag{Event $13^{\prime}$} \\ \overset{(\textbf{Proposition} \text{ } \textit{$\mathrm{GRCM}$ 1})}{=}                \nu^{+}_{\Lambda}  \big[      \big(  \sigma_{-,+}\big)_1  \sigma_{L^{\prime}}       \big] \bigg[  1 - \frac{v^{+}_{\Lambda} \big[ \big(  \sigma_{-,+}\big)_1  \sigma_{L^{\prime}}    \big]}{v^{-}_{\Lambda} \big[ \big(  \sigma_{-,+}\big)_1  \sigma_{L^{\prime}}    \big] }  \bigg]  \text{ } \text{ , }
   \end{align*}
   
    \noindent with,
   
   \begin{align*}
         \big(  \gamma_{-,+}   \big)_1  \cap    \big(  \sigma_{-,+}\big)_1  \neq \emptyset     \text{ } \text{ , }  \\     \gamma_{L^{\prime}}    \cap \sigma_{L^{\prime}} \neq \emptyset  \text{ } \text{ . } 
   \end{align*}
   
   \item[$\bullet$] \textit{$\mathcal{P}_3$}. From (Event $14$),
   
   \begin{align*}
   \mathcal{P}^{\chi_2}_{\Lambda} \big[     \big(   \gamma_{-,+}  \big)_2          \underset{\Lambda \cap (\mathscr{F}\mathscr{C})_1}{\overset{+}{\longleftrightarrow}}                     \mathscr{F}\mathscr{C}_1                \big]   - \mathcal{P}^{\chi_1}_{\Lambda} \big[     \big(   \gamma_{-,+}  \big)_2          \underset{\Lambda \cap (\mathscr{F}\mathscr{C})_1}{\overset{+}{\longleftrightarrow}}                     \mathscr{F}\mathscr{C}_1                \big] \approx  \mathcal{P}^{\chi_2}_{\Lambda} \big[         \big( \sigma_{-,+}\big)_2 \sigma_{\mathscr{F}\mathscr{C}_1}     \big] -  \mathcal{P}^{\chi_1}_{\Lambda} \big[    \big( \sigma_{-,+}\big)_2 \sigma_{\mathscr{F}\mathscr{C}_1}   \big]                \\ \equiv \mathcal{P}^{\chi_2}_{\Lambda} \big[         \big( \sigma_{-,+}\big)_2 \sigma_{\mathscr{F}\mathscr{C}_1}     \big]    \bigg[  1 -   \frac{\mathcal{P}^{\chi_1}_{\Lambda} \big[    \big( \sigma_{-,+}\big)_2 \sigma_{\mathscr{F}\mathscr{C}_1}   \big] }{\mathcal{P}^{\chi_2}_{\Lambda} \big[         \big( \sigma_{-,+}\big)_2 \sigma_{\mathscr{F}\mathscr{C}_1}     \big] }  \bigg]                 \text{ } \text{ , } \tag{Event $14^{\prime}$} \\ \overset{(\textbf{Proposition} \text{ } \textit{$\mathrm{GRCM}$ 1})}{=}                \nu^{+}_{\Lambda}  \big[        \big( \sigma_{-,+}\big)_2 \sigma_{\mathscr{F}\mathscr{C}_1}         \big]  \bigg[  1 -    \frac{\nu^{+}_{\Lambda}  \big[        \big( \sigma_{-,+}\big)_2 \sigma_{\mathscr{F}\mathscr{C}_1}         \big] }{\nu^{-}_{\Lambda}  \big[        \big( \sigma_{-,+}\big)_2 \sigma_{\mathscr{F}\mathscr{C}_1}         \big] }    \bigg]   \text{ } \text{ , }
   \end{align*}

   \noindent with,
   
   \begin{align*}
     \big( \gamma_{-,+} \big)_2 \cap   \big( \sigma_{-,+}\big)_2   \neq \emptyset       \text{ } \text{ , } \\ \mathscr{F}\mathscr{C}_1    \cap    \sigma_{\mathscr{F}\mathscr{C}_1}         \neq \emptyset  \text{ } \text{ . } 
   \end{align*}
   
    \item[$\bullet$] \textit{First conditional crossing event for finishing arguments of \textbf{Lemma} \textit{6.3}}. From (Event $15$),
    
    \begin{align*}
      \mathcal{P}_{\Lambda} \bigg[              \big(   \gamma_{-,+}  \big)_2                       \underset{\Lambda \cap (\mathscr{F}\mathscr{C})_1}{\overset{+}{\longleftrightarrow}}    \mathscr{F}^{\prime\prime}          \big|         \mathcal{C}_{-} \cup    \mathscr{D}_{\chi_3 , -}              \bigg]            \approx   \mathcal{P}_{\Lambda} \big[   \big( \sigma_{-,+} \big)_2 \sigma_{\mathscr{F}^{\prime\prime}}                             |            \sigma_{j-1,-} \sigma_{-}    \sigma_{\mathscr{L}_j} \sigma_{-}             \big]  \text{ } \text{ , }  \tag{Event $15^{\prime}$} \\ \overset{(\textbf{Proposition} \text{ } \textit{$\mathrm{GRCM}$ 1})}{=}                \nu^{+}_{\Lambda}  \bigg[  \big(   \gamma_{-,+}  \big)_2                   \overset{\sigma}{\longleftrightarrow}    \mathscr{F}^{\prime\prime}       \big| \big\{     \gamma_{j-1,-}     \overset{\sigma}{\longleftrightarrow} \gamma_{-}          \big\} \cup \big\{     \mathscr{L}_{j,-} \overset{\sigma}{\longleftrightarrow}       \big( \Lambda \cap \big( \mathscr{F} \mathscr{C} \big)^c \big)                \big\}    \bigg]  \text{ } \text{ , }  
    \end{align*}
    
    \noindent with,
    
     \begin{align*}
             \big( \sigma_{-,+} \big)_2    \cap    \big( \gamma_{-,+} \big)_2             \neq \emptyset \text{ } \text{ , }      \\  \sigma_{\mathscr{F}^{\prime\prime}}  \cap     \mathscr{F}^{\prime\prime}           \neq \emptyset \text{ } \text{ , } \\ \sigma_{j-1,-}          \cap \sigma_{-} \cap \mathcal{C}_{-}       \neq \emptyset                \text{ } \text{ , } \\        \sigma_{\mathscr{L}_j}           \cap       \sigma_{-} \cap    \mathscr{D}_{\chi_3,-}      \neq \emptyset    \text{ } \text{ . } 
    \end{align*}

     \item[$\bullet$] \textit{Second conditional crossing event for finishing arguments of \textbf{Lemma} \textit{6.3}}. From (Event $16$),
    
    \begin{align*}
        \mathcal{P}^{\chi_2}_{\Lambda} \bigg[              \big(   \gamma_{-,+}  \big)_2                       \underset{\Lambda \cap (\mathscr{F}\mathscr{C})_1}{\overset{+}{\longleftrightarrow}}    \mathscr{F}^{\prime\prime}_i          \big|         \mathcal{C}_{-} \cup    \mathscr{D}_{\chi_3 , -}              \bigg] \approx  \mathcal{P}^{\chi_2}_{\Lambda} \big[    \big( \sigma_{-,+} \big)_2 \sigma_{\mathscr{F}^{\prime\prime}_i}                                |  \sigma_{j-1,-} \sigma_{-}    \sigma_{\mathscr{L}_j} \sigma_{-}           \big] \text{ } \text{ , }  \tag{Event $16^{\prime}$} \\ \overset{(\textbf{Proposition} \text{ } \textit{$\mathrm{GRCM}$ 1})}{=}                \nu^{+}_{\Lambda}  \bigg[  \big( \gamma_{-,+} \big)_2 \overset{\sigma}{\longleftrightarrow}     \mathscr{F}^{\prime\prime}_i           \big| \big\{     \gamma_{j-1,-}     \overset{\sigma}{\longleftrightarrow} \gamma_{-}          \big\} \cup \big\{     \mathscr{L}_{j,-} \overset{\sigma}{\longleftrightarrow}       \big( \Lambda \cap \big( \mathscr{F} \mathscr{C} \big)^c \big)                \big\}  \bigg]  \text{ } \text{ , }  
    \end{align*}

    \noindent with,

    \begin{align*}
      \big( \sigma_{-,+} \big)_2 \cap \big(\gamma_{-,+} \big)_2        \neq \emptyset   \text{ } \text{ , } \\  \sigma_{j-1,-}          \cap \sigma_{-} \cap \mathcal{C}_{-}       \neq \emptyset                \text{ } \text{ , } \\        \sigma_{\mathscr{L}_j}           \cap       \sigma_{-} \cap    \mathscr{D}_{\chi_3,-}      \neq \emptyset        \text{ } \text{ . } 
    \end{align*}

    \item[$\bullet$] \textit{Intersection of three crossing events, \textbf{Lemma} \textit{6.2}}. From (Event $17$), for $\big| \mathcal{I}^{+ \backslash -} \big| < + \infty$,
    
    \begin{align*}
       \underset{R_1 \cap \gamma_R \neq \emptyset}{ \underset{L_1 \cap \gamma_L \neq \emptyset}{\underset{1\leq j \leq {\big| \mathcal{I}^{+ \backslash -}\big|}}{\prod}}}    \mathcal{P}^{\chi}_{\Lambda} \bigg[   \big\{ \gamma_{L_j} \overset{+ \backslash -}{\longleftrightarrow} \mathscr{F}\mathscr{C}_1       \big\}  \cap    \big\{   \mathscr{F}\mathscr{C}_2 \overset{+\backslash -}{\not\longleftrightarrow} \mathscr{F}\mathscr{C}_{N-1}         \big\}              \cap       \big\{    \mathscr{F}\mathscr{C}_N              \overset{+\backslash -}{\longleftrightarrow}     \gamma_{R_j}          \big\}            \bigg]   \approx \mathcal{P}^{\chi}_{\Lambda} \big[      \sigma_{\not\in \mathscr{F}\mathscr{C}_2} \sigma_{\not\in \mathscr{F}\mathscr{C}_{N-1}}            \big]      \\ \times      \underset{R_1 \cap \gamma_R \neq \emptyset}{ \underset{L_1 \cap \gamma_L \neq \emptyset}{\underset{1\leq j \leq {\big| \mathcal{I}^{+ \backslash -}\big|}}{\prod}}}  \mathcal{P}^{\chi}_{\Lambda} \big[     \sigma_{L_j} \sigma_{\mathscr{F}\mathscr{C}_1} \sigma_{\mathscr{F} \mathscr{C}_N}      \big]  \text{ } \text{ , } \tag{Event $17^{\prime}$} \\ \overset{(\textbf{Proposition} \text{ } \textit{$\mathrm{GRCM}$ 1})}{=}                \nu^{+}_{\Lambda}  \big[              \mathscr{F}\mathscr{C}_2    \overset{\sigma}{\not\longleftrightarrow}       \mathscr{F}\mathscr{C}_N        \big]     \underset{R_1 \cap \gamma_R \neq \emptyset}{ \underset{L_1 \cap \gamma_L \neq \emptyset}{\underset{1\leq j \leq {\big| \mathcal{I}^{+ \backslash -}\big|}}{\prod}}}           \nu^{+}_{\Lambda}  \bigg[    \big\{       \gamma_{L_j}               \overset{\sigma}{\longleftrightarrow}     \mathscr{F}\mathscr{C}_1      \big\}               \cap \big\{ \mathscr{F}\mathscr{C}_N  \overset{\sigma}{\longleftrightarrow}  \gamma_{R_j}   \big\}        \bigg]   \text{ } \text{ , }  
    \end{align*}
    
    \noindent with,
   
    \begin{align*}
          \sigma_{L_j} \cap \gamma_{L_j} \neq \emptyset    \text{ } \text{ , } \\ \sigma_{R_j}  \cap   \gamma_{R_j}     \neq \emptyset  \text{ } \text{ , } \\ \sigma_{\not\in\mathscr{F}\mathscr{C}_2}   \cap   \mathscr{F}\mathscr{C}_2   \equiv \emptyset   \text{ } \text{ , } \\  \sigma_{\not\in\mathscr{F}\mathscr{C}_{N-1}}   \cap \mathscr{F}\mathscr{C}_{N-1}          \equiv \emptyset   \text{ } \text{ . } 
    \end{align*}
    
   \item[$\bullet$] \textit{Second connectivity event from \textbf{Lemma} \textit{6.2}}. From (Event $18$),
   
   \begin{align*}
        \underset{2 \leq i \leq i+1 \leq N-1}{\prod}  \mathcal{P}^{\chi}_{\Lambda} \big[    \mathscr{F}\mathscr{C}_i \overset{+\backslash -}{\not\longleftrightarrow} \mathscr{F}\mathscr{C}_{i+1}      \big]  \approx  \underset{2 \leq i \leq i+1 \leq N-1}{\prod}  \mathcal{P}^{\chi}_{\Lambda} \big[            \sigma_i    {\not\longleftrightarrow} \sigma_{i+1}              \big]  \text{ } \text{ , } \tag{Event $18^{\prime}$} \\ \overset{(\textbf{Proposition} \text{ } \textit{$\mathrm{GRCM}$ 1})}{=}               \underset{2 \leq i \leq i+1 \leq N-1}{\prod}  \nu^{+}_{\Lambda}  \big[  \mathscr{F}\mathscr{C}_i  \not\longleftrightarrow   \mathscr{F}\mathscr{C}_{i+1}    \big]  \text{ } \text{ , }  
   \end{align*}

   \noindent with,
    
    \begin{align*}
     \mathscr{F}\mathscr{C}_i   \cap   \sigma_i     \neq \emptyset       \text{ } \text{ , }  \forall \text{ }  2 \leq i  \leq N - 1 \text{ } \text{ , } 
        \\   \mathscr{F}\mathscr{C}_{i+1}  \cap   \sigma_{i+1}     \neq \emptyset \text{ } \text{ , } \forall \text{ }  2  \leq i+1 \leq N - 1 \text{ }  \text{ . }
    \end{align*}
   
   \item[$\bullet$] \textit{Crossing event from \textbf{Lemma} \textit{AT 2.6.1.1}}. From (Event $19$),

   \begin{align*}
   \underset{L_i}{\prod}         \mathcal{P}^{\chi}_{\Lambda}  \big[      \gamma_L \overset{+ \backslash -}{\longleftrightarrow}                   \gamma_{L_i}  \big]  \approx  \underset{L_i}{\prod}         \mathcal{P}^{\chi}_{\Lambda}  \big[    \sigma_L \sigma_{L_i}                    \big]   \text{ } \text{ } \tag{Event $19^{\prime}$} \\ \overset{(\textbf{Proposition} \text{ } \textit{$\mathrm{GRCM}$ 1})}{=}      \underset{L_i}{\prod}            \nu^{+}_{\Lambda}  \big[    \gamma_L   \overset{\sigma}{\longleftrightarrow}    \gamma_{L_i}       \big]  \text{ } \text{ , }   
   \end{align*}

   \noindent with,
   
   \begin{align*}
     \gamma_L   \cap   \sigma_L     \neq \emptyset      \text{ } \text{ , } \\         \gamma_{L_i}   \cap   \sigma_{L_i}     \neq \emptyset   \text{ } \text{ . } 
   \end{align*}

   \item[$\bullet$] \textit{Disconnectivity event from \textbf{Lemma} \textit{AT 2.6.1.2}}. From (Event $20$),
   
    \begin{align*}
       \mathcal{P}^{\chi}_{\Lambda} \big[   \underset{\text{ countably many } i^{\prime} }{\bigcap} \text{ } \big\{ \mathscr{F}\mathscr{C}_i   \overset{+\backslash -}{\not\longleftrightarrow}         \gamma_{i^{\prime}}           \big\}      \big]  \approx  \mathcal{P}^{\chi}_{\Lambda} \big[    \underset{\text{ countably many } i^{\prime} }{\bigcap} \text{ } \big\{            \sigma_{\mathscr{F}\mathscr{C}_i} {\not\longleftrightarrow}     \sigma_{i^{\prime}}   \big\}     \big]  \text{ } \text{ , }  \tag{Event $20^{\prime}$} \\ \overset{(\textbf{Proposition} \text{ } \textit{$\mathrm{GRCM}$ 1})}{=}                \nu^{+}_{\Lambda}  \big[   \text{ }   \underset{\text{ countably many } i^{\prime} }{\bigcap}   \big\{  \mathscr{F}\mathscr{C}_i \overset{\sigma}{\not\longleftrightarrow}  \gamma_{i^{\prime}} \big\}         \text{ }        \big]  \text{ } \text{ , }  
    \end{align*}

   \noindent with,
    
        \begin{align*}
             \sigma_{\mathscr{F}\mathscr{C}_i}   \cap           \mathscr{F}\mathscr{C}_i      \neq \emptyset   \text{ } \text{ , } \\  \sigma_{i^{\prime}}   \cap        \gamma_{i^{\prime}}     \neq  \emptyset      \text{ } \text{ . } 
    \end{align*}
    
    \item[$\bullet$] \textit{First connectivity event from \textbf{Lemma} \textit{AT 2.6.1.3}}. From (Event $21$),
    
    \begin{align*}
        \mathcal{P}^{\chi}_{\Lambda}  \big[    \underset{F^{\prime}_j \in F ( \Lambda) }{\underset{\text{ countably many } j^{\prime}}{\bigcap}} \big\{      F_i  \underset{(\mathscr{F} \mathscr{C})^c}{\overset{+ \backslash -}{\longleftrightarrow}}    F^{\prime}_j              \big\}       \big]    \approx    \mathcal{P}^{\chi}_{\Lambda}  \big[           \underset{F^{\prime}_j \in F ( \Lambda) }{\underset{\text{ countably many } j^{\prime}}{\bigcap}} \big\{    \sigma_i          {\longleftrightarrow} \sigma^{\prime}_j  \big\} \big]    \text{ } \text{ , } \tag{Event $21^{\prime}$} \\ \overset{(\textbf{Proposition} \text{ } \textit{$\mathrm{GRCM}$ 1})}{=}                \nu^{+}_{\Lambda}  \big[           \underset{F^{\prime}_j \in F ( \Lambda) }{\underset{\text{ countably many } j^{\prime}}{\bigcap}}      \big\{  F_i \overset{\sigma}{\longleftrightarrow} F^{\prime}_j  \big\}        \big]  \text{ } \text{ , }  
    \end{align*}
    
    \noindent with,
   
    \begin{align*}
               \sigma_i \cap F_i \neq \emptyset        \text{ } \text{ , }   \\ \sigma^{\prime}_j  \cap F^{\prime}_j   \neq \emptyset \text{ } \text{countably many, strictly positive } j^{\prime}   \text{ } \text{ . } 
    \end{align*}

   \item[$\bullet$] \textit{Second connectivity event from \textbf{Lemma} \textit{AT 2.6.1.3}}. From (Event $22$),

   \begin{align*}
   \mathcal{P}^{\chi}_{\Lambda} \big[   \text{ }    \underset{\exists \text{ }  i > i^{\prime} \text{ } : \text{ } \mathscr{L}_i \cap  \gamma_R \neq \emptyset}{\underset{\text{countably many } i^{\prime}}{\bigcap}}  \text{ }  \big\{              \mathscr{F}\mathscr{C}_N \overset{+ \backslash - }{\longleftrightarrow}   \mathscr{L}_{i^{\prime}}  \big\}   \big]  \approx   \mathcal{P}^{\chi}_{\Lambda} \big[    \text{ }    \underset{\exists \text{ }  i > i^{\prime} \text{ } : \text{ } \mathscr{L}_i \cap  \gamma_R \neq \emptyset}{\underset{\text{countably many } i^{\prime}}{\bigcap}}  \text{ }  \big\{                    \sigma_N    {\longleftrightarrow}        \sigma_{\mathscr{L}_{i^{\prime}}}    \big\}   \text{ }        \big]    \text{ } \text{ , } \tag{Event $22^{\prime}$} \\ \overset{(\textbf{Proposition} \text{ } \textit{$\mathrm{GRCM}$ 1})}{=}                \nu^{+}_{\Lambda}  \big[      \text{ }    \underset{\exists \text{ }  i > i^{\prime} \text{ } : \text{ } \mathscr{L}_i \cap  \gamma_R \neq \emptyset}{\underset{\text{countably many } i^{\prime}}{\bigcap}}       \big\{    \mathscr{F}\mathscr{C}_N   \overset{\sigma}{\longleftrightarrow}   \mathscr{L}^{\prime}_i   \big\}             \text{ }       \big]  \text{ } \text{ , }  
   \end{align*}

   \noindent with,
   
   \begin{align*}
       \sigma_N      \cap   \mathscr{F}\mathscr{C}_N \neq \emptyset   \text{ } \text{ , } \\  \sigma_{\mathscr{L}_{i^{\prime}}}   \cap \mathscr{L}_{i^{\prime}} \neq \emptyset \text{ } \text{ . } 
   \end{align*}

   \item[$\bullet$] \textit{Connectivity event from \textbf{Lemma} \textit{6.2}}. From (Event $23$),
   
   \begin{align*}
    \mathcal{P}^{\chi}_{\mathscr{D}^{\prime\prime}} \big[     \gamma_L \longleftrightarrow \mathscr{B}^{\prime}  \big] \approx  \mathcal{P}^{\chi}_{\mathscr{D}^{\prime\prime}} \big[     \sigma_L    \longleftrightarrow \mathscr{B}^{\prime} \big]      \tag{Event $23^{\prime}$}  \overset{(\textbf{Proposition} \text{ } \textit{$\mathrm{GRCM}$ 1})}{=}                \nu^{+}_{\Lambda}  \big[  \gamma_L   \overset{\sigma}{\longleftrightarrow}  \mathscr{B}^{\prime}   \big]  \text{ } \text{ , }  
   \end{align*}
   
   \noindent with,
   
   \begin{align*}
    \sigma_L \cap \gamma_L \neq \emptyset   \text{ } \text{ . } 
   \end{align*}
   
  \item[$\bullet$] \textit{Connectivity event for bridging event lower bound \textbf{Lemma} \textit{AT 2.7}}. From (Event $24$),
  
  \begin{align*}
      \mathcal{P}_{\Lambda}^{\chi} \bigg[ \mathcal{I}_{-j} \underset{\mathscr{F}\mathscr{C} \cap \Lambda}{\overset{+ \backslash -}{\longleftrightarrow}}   \widetilde{\mathcal{I}_{-j}}   \bigg] \text{ }  \mathcal{P}_{\Lambda}^{\chi} \bigg[  \mathcal{I}_{-k} \underset{\mathscr{F}\mathscr{C} \cap \Lambda}{\overset{+ \backslash -}{\longleftrightarrow}}   \widetilde{\mathcal{I}_{-k}}   \bigg] \approx     \mathcal{P}_{\Lambda}^{\chi} \big[  \sigma_j \widetilde{\sigma_j}               \big]  \mathcal{P}_{\Lambda}^{\chi} \big[    \sigma_k \widetilde{\sigma_k}                   \big]          \text{ } \text{ , }    \tag{Event $24^{\prime}$}
 \\ \overset{(\textbf{Proposition} \text{ } \textit{$\mathrm{GRCM}$ 1})}{=}                \nu^{+}_{\Lambda}  \big[   \mathcal{I}_{-j}    \overset{\sigma}{\longleftrightarrow} \widetilde{\mathcal{I}_{-j}}  \big] \nu^{+}_{\Lambda}  \big[   \mathcal{I}_{-k}     \overset{\sigma}{\longleftrightarrow} \widetilde{\mathcal{I}_{-k}}    \big]  \text{ } \text{ , }  
\end{align*}

\noindent with,

 \begin{align*}
 \sigma_{-j}  \cap       \mathcal{I}_{-j}      \neq \emptyset    \text{ } \text{ , }     \\  \widetilde{\sigma_{-j}} \cap  \widetilde{\mathcal{I}_{-j}}           \neq \emptyset    \text{ } \text{ , }     \\   \sigma_{-k}  \cap       \mathcal{I}_{-k}          \neq \emptyset    \text{ } \text{ , }     \\ \widetilde{\sigma_{-k}} \cap  \widetilde{\mathcal{I}_{-k}}            \neq \emptyset  \text{ } \text{ . } 
  \end{align*}

  \item[$\bullet$] \textit{Condtional bridging event from \textbf{Lemma} \textit{AT 2.7}}. From (Event $25$),
  
  \begin{align*}
 \mathcal{P}^{\chi}_{\Lambda} \big[   \mathcal{B}_{+ \backslash -}(j)                              \big|     \mathrm{sign} \big(   \mathscr{F}|_{\mathscr{D}^c \cap \Lambda} \big) \equiv +                    \big] \approx      \mathcal{P}^{\chi}_{\Lambda} \big[ \sigma_j \widetilde{\sigma_j}   |   \sigma_{+}    \big]       \tag{Event $25^{\prime}$}  \overset{(\textbf{Proposition} \text{ } \textit{$\mathrm{GRCM}$ 1})}{=}                \nu^{+}_{\Lambda}  \big[     \mathcal{I}_j  \overset{\sigma}{\longleftrightarrow}  \widetilde{\mathcal{I}_j}   | \mathscr{D}^c \overset{\sigma}{\longleftrightarrow}  \partial \Lambda    \big]  \text{ } \text{ , }  
  \end{align*}

  \noindent with,
 
 \begin{align*}
 \sigma_{-j} \cap \mathcal{I}_{-j} \neq \emptyset  \text{ } \text{ , }    \\  \widetilde{\sigma_{-j}} \cap \widetilde{\mathcal{I}_{-j}} \neq \emptyset   \text{ } \text{ , }  \\  \sigma_{+}    \cap \big( \mathscr{D}^c \cap \Lambda \big) \neq \emptyset   \text{ } \text{ . } 
 \end{align*}
 
 \item[$\bullet$] \textit{Conditional horizontal crossing event from \textbf{Lemma} \textit{AT 2.7}}. From (Event $26$),
 
 \begin{align*}
    \mathcal{P}^{\chi}_{\Lambda} \big[      \mathcal{H}_{+ \backslash -} \big( \mathscr{D} \big) |             \mathrm{sign} \big(   \mathscr{F}|_{\mathscr{D}^c \cap \Lambda} \big) \equiv +                     \big] \approx    \mathcal{P}^{\chi}_{\Lambda} \big[  \sigma_{L_{\mathscr{D}}}  \sigma_{R_{\mathscr{D}}}  |    \sigma_{+}    \big]  \text{ } \text{ , } \tag{Event $26^{\prime}$} \\ \overset{(\textbf{Proposition} \text{ } \textit{$\mathrm{GRCM}$ 1})}{=}                \nu^{+}_{\Lambda}  \big[   L_{\mathscr{D}}  \overset{\sigma}{\longleftrightarrow}      R_{\mathscr{D}}        \big|      \mathscr{D}^c \overset{\sigma}{\longleftrightarrow}  \partial \Lambda          \big]  \text{ } \text{ , }  
 \end{align*}
 
 \noindent with,

\begin{align*}
    \sigma_{L_{\mathscr{D}}} \cap \mathscr{D} \neq \emptyset  \text{ } \text{ , }   \\ \sigma_{R_{\mathscr{D}}} \cap \mathscr{D} \neq \emptyset    \text{ } \text{ , }     \\ \sigma_{+} \cap \big( \mathscr{D}^c \cap \Lambda ) \neq \emptyset    \text{ } \text{ . }  
\end{align*}

  \item[$\bullet$] \textit{Connectivity event upper bounded by $1-c_0$ in \textbf{Lemma} \textit{AT 2.7}}. From (Event $27$),

  \begin{align*}
   \mathcal{P}^{\chi}_{\Lambda} \big[  [- \lfloor \delta^{\prime\prime} n \rfloor   , 2 \lfloor   \delta^{\prime\prime} n \rfloor          \big]  \times \{ 0 \} {\longleftrightarrow}  [- \lfloor \delta^{\prime\prime}  n \rfloor   , 2 \lfloor   \delta^{\prime\prime} n  \rfloor          \big]  \times \{ n \} \big] \approx            \mathcal{P}^{\chi}_{\Lambda} \big[   \sigma_{-\lfloor \delta^{\prime\prime} n \rfloor}       \sigma_{2 \lfloor \delta^{\prime\prime} n \rfloor}      \big]       \text{ } \text{ , }        \tag{Event $27^{\prime}$} \\ \overset{(\textbf{Proposition} \text{ } \textit{$\mathrm{GRCM}$ 1})}{=}                \nu^{+}_{\Lambda}  \big[  [- \lfloor \delta^{\prime\prime} n \rfloor   , 2 \lfloor   \delta^{\prime\prime} n \rfloor          \big]  \times \{ 0 \} \overset{\sigma}{\longleftrightarrow}     [- \lfloor \delta^{\prime\prime}  n \rfloor   , 2 \lfloor   \delta^{\prime\prime} n  \rfloor          \big]  \times \{ n \}    \big]  \text{ } \text{ , }  
  \end{align*}
  
  \noindent with,
  
  \begin{align*}
     \sigma_{-\lfloor \delta^{\prime\prime} n \rfloor }  \cap  \big( [ - \lfloor \delta^{\prime\prime} n \rfloor     , 2 \lfloor \delta^{\prime\prime} n \rfloor ] \times \{ 0 \} \big)  \neq \emptyset      \text{ } \text{ , } \\    \sigma_{-2 \lfloor \delta^{\prime\prime} n \rfloor }  \cap  \big( [ - \lfloor \delta^{\prime\prime} n \rfloor     , 2 \lfloor \delta^{\prime\prime} n \rfloor ] \times \{ n \} \big)       \neq \emptyset  \text{ } \text{ . } 
  \end{align*}

    \item[$\bullet$] \textit{Connectivity event upper bounded by $1-c$ in \textbf{Lemma} \textit{AT 2.7}}. From (Event $28$),
  
  \begin{align*}
    \mathcal{P}^{\xi}_{\Lambda} \big[     [ 0, \lfloor \delta^{\prime\prime} n \rfloor ] \times \{ 0 \} {\longleftrightarrow}  [i , i+ \delta^{\prime\prime}n  ]  \times \{n \}  \big] \approx    \mathcal{P}^{\xi}_{\Lambda} \big[      \sigma_0 \sigma_{\lfloor \delta^{\prime\prime} n \rfloor}        \sigma_i \sigma_{i+\lfloor \delta^{\prime\prime} n \rfloor}               \big]  \text{ } \text{ , } \tag{Event $28^{\prime}$} \\ \overset{(\textbf{Proposition} \text{ } \textit{$\mathrm{GRCM}$ 1})}{=}                \nu^{+}_{\Lambda}  \big[   [0 , \lfloor \delta^{\prime\prime} n \rfloor ] \times \{ 0 \} \overset{\sigma}{\longleftrightarrow}             [i , i+ \delta^{\prime\prime}n  ]  \times \{n \}    \big]  \text{ } \text{ , }  
  \end{align*}

  \noindent with,
 
  \begin{align*}
      \sigma_0 \cap \sigma_{\lfloor \delta^{\prime\prime} n \rfloor } \cap [0 ,  \lfloor \delta^{\prime\prime} n \rfloor  ] \neq \emptyset       \text{ } \text{ , } \\           \sigma_i \cap  \sigma_{i+\lfloor \delta^{\prime\prime} n \rfloor }   \cap [ i , i + \lfloor \delta^{\prime\prime} n \rfloor ]  \neq \emptyset      \text{ } \text{ , } \\    \sigma_0 \cap \sigma_{\lfloor \delta^{\prime\prime} n \rfloor } \cap \big( [0 ,  \lfloor \delta^{\prime\prime} n \rfloor  ] \times \{ 0 \} \big) \neq \emptyset       \text{ } \text{ , } \\      \sigma_i \cap  \sigma_{i+\lfloor \delta^{\prime\prime} n \rfloor }   \cap  \big( [ i , i + \lfloor \delta^{\prime\prime} n \rfloor ]        \times \{ n \} \big)   \neq \emptyset    \text{ } \text{ . }
  \end{align*}

 \item[$\bullet$] \textit{First connectivity event from \textbf{Proposition} \textit{AT 2}}. From (Event $29$),

 \begin{align*}
       \mathcal{P}^{\chi}_{\Lambda} \big[    \{ 0 \} \times [ - n^{\prime} , 2 n^{\prime} ]     \longleftrightarrow      \{  \rho \lfloor \delta^{\prime\prime} n \rfloor      \} \times [ - n^{\prime} , 2 n^{\prime} ]      \big] \approx        \mathcal{P}^{\chi}_{\Lambda} \big[    \sigma_0 \sigma_{-n^{\prime}} \sigma_{2n^{\prime}}     \sigma_{\rho \lfloor \delta^{\prime\prime} n \rfloor}                        \big]   \text{ } \text{ , } \tag{Event $29^{\prime}$} \\ \overset{(\textbf{Proposition} \text{ } \textit{$\mathrm{GRCM}$ 1})}{=}                \nu^{+}_{\Lambda}  \big[      \{ 0 \} \times [ - n^{\prime} , 2 n^{\prime} ]     \overset{\sigma}{\longleftrightarrow}       \{  \rho \lfloor \delta^{\prime\prime} n \rfloor      \} \times [ - n^{\prime} , 2 n^{\prime} ]      \big]  \text{ } \text{ , }   
 \end{align*}
 
  \noindent with,

  \begin{align*}
        \sigma_0 \cap \sigma_{-n^{\prime}} \cap [-n^{\prime} , 2n^{\prime} ] \neq \emptyset       \text{ } \text{ , } \\     \sigma_0 \cap \sigma_{-n^{\prime}} \cap \big( \{ 0 \} \times [-n^{\prime} , 2n^{\prime} ] \big) \neq \emptyset    \text{ } \text{ , }    \\        \sigma_{2n^{\prime}} \cap   \sigma_{n^{\prime}} \cap [ - n^{\prime} , 2 n^{\prime} ]      \neq \emptyset  \text{ } \text{ , } \\  \sigma_{2n^{\prime}} \cap   \sigma_{n^{\prime}} \cap  \big(      \{  \rho \lfloor \delta^{\prime\prime} n \rfloor      \}   \times [ - n^{\prime} , 2 n^{\prime} ]     \big) \neq \emptyset \text{ } \text{ . } 
  \end{align*}

  \item[$\bullet$] \textit{Second connectivity event for \textbf{Proposition} \textit{AT 2}}. From (Event $30$),

  \begin{align*}
    \big( \mathcal{P}^{\chi}_{\Lambda} \big[  \mathcal{I}_0                \longleftrightarrow    \textbf{Z} \times \{ n \}   \big]\big)^2  \overset{(\textbf{Proposition} \text{ } \textit{$\mathrm{GRCM}$ 1})}{\approx}       \nu^{+}_{\Lambda} \bigg[ \underset{\mathcal{I}_{i^{\prime}} \cap \textbf{Z} \neq \emptyset \text{ } , \text{ } \forall i^{\prime} }{\underset{\mathrm{countably\text{ }  many \text{ } }\text{ } i^{\prime}}{\bigcap}} \bigg\{    \mathcal{I}_0  \overset{\sigma}{\longleftrightarrow }    \big\{    \mathcal{I}_{i^{\prime}} \times \{ n \}   \big\}       \bigg\}       \bigg]^2  \tag{Event $30^{\prime}$}
  \end{align*}
  
  \noindent with,

\begin{align*}
       \sigma^{\prime}_0    \cap     \mathcal{I}_0    \neq \emptyset                     \text{ } \text{ , } \\  \sigma_{i^{\prime}}   \cap      \big( \textbf{Z} \times \{ n \} \big)    \neq \emptyset   \text{ } \text{ , } \forall i^{\prime}   \text{ } \text{ . } 
\end{align*}

\item[$\bullet$] \textit{Third connectivity event for \textbf{Proposition} \textit{AT 2}}. From (Event $31$),

\begin{align*}
     \mathcal{P}^{\chi}_{\Lambda} \big[ \mathcal{I}_0 \longleftrightarrow        (- \infty , 0 ] \times \{ n \}      \big] \approx   \mathcal{P}^{\chi}_{\Lambda} \big[            \sigma_0 \sigma_{-\infty}                \big] \text{ } \text{ , }            \tag{Event $31^{\prime}$} \\ \overset{(\textbf{Proposition} \text{ } \textit{$\mathrm{GRCM}$ 1})}{=}                \nu^{+}_{\Lambda}  \big[  \mathcal{I}_0 \overset{\sigma}{\longleftrightarrow}     \big\{     (- \infty , 0 ] \times \{ n \}  \big\}         \big]  \text{ } \text{ , }  
\end{align*}    

\noindent with,

\begin{align*}
 \sigma_0 \cap  \mathcal{I}_0    \neq \emptyset  \text{ } \text{ , } \\ \sigma_{-\infty}  \cap \big( ( - \infty , 0 ] \times \{ n \} \big)      \neq \emptyset \text{ } \text{ . } 
\end{align*}

\end{itemize}

\bigskip

\noindent To analyze crossing probabilities such as the one provided above, we provide the following statements for quantifying crossing estimates in \textit{generalized random-cluster model} setting, building upon crossing probability estimates for the Ashkin-Teller model. In doing so, we directly adapt the \textit{symmetric domains} for the Ashkin-Teller model, which correspond to \textit{symmetric domains} for the generalized random-cluster model.

\subsection{Generalized random-cluster crossing events}

\noindent We introduce equivalent formulations of each crossing event obtained from \textbf{Proposition} \textit{GRCM 1}. Each proof is a direct application of previous arguments for weakened crossing probability estimates of the Ashkin-Teller model, and are hence omitted. 

\bigskip

\noindent \textbf{Proposition} \textit{GRCM 1} (\textit{GRCM equivalent of \textbf{Proposition} AT 1}). Under equivalent choice of parameters, (Event $1^{\prime}$) admits an upper bound of $1-c$.

\bigskip

\noindent \textit{Proof sketch of Proposition GRCM 1}. Apply the same sequence of arguments as given for the proof of \textbf{Proposition} \textit{AT 1}, in which arguments for upper bounding (Event $2^{\prime}$) are given in the statement of \textbf{Lemma} \textit{GRCM 1} below.

\bigskip

\noindent \textbf{Lemma} \textit{GRCM 1} (\textit{GRCM equivalent of \textbf{Lemma} \textit{AT 2.1}}). Under equivalent choice of parameters, from (Event $2^{\prime}$), an upper bound for,

\begin{align*}
  \nu^{+}_{\Lambda}  \big[   \mathcal{I}_0 \overset{\sigma}{\longleftrightarrow} \widetilde{\mathcal{I}_0}  \big]       \text{ } \text{ , } 
\end{align*}

\noindent can be obtained for sufficiently small $\widetilde{\widetilde{\mathscr{R}}}$, from the expression,

\begin{align*}
 1 - \widetilde{\widetilde{\mathscr{R}}} ( Event 3^{\prime} )    \text{ } \text{ . }  \tag{$\widetilde{\widetilde{\mathscr{R}}}$-bound}
\end{align*}

\noindent \textit{Proof sketch of Lemma GRCM 1}. Obtain the desired upper bound, from ($\widetilde{\widetilde{\mathscr{R}}}$-bound) above by following the identical sequence of steps provided in , and in , respectively for , and for . After obtaining the appropriate range of parameters for $\widetilde{\widetilde{\mathscr{R}}}$, upper bound (Event $4^{\prime}$) with $1-c$, in addition to obtaining the final lower bound for (Event $5^{\prime}$), hence concluding the argument.

\bigskip

\noindent \textbf{Lemma} \textit{GRCM 2} (\textit{GRCM equivalent of \textbf{Lemma} \textit{6.3}, one crossing probability lower bound}). Under equivalent choice of parameters, a lower bound for (Event $6^{\prime}$) can be obtained with a suitably chosen, strictly positive, parameter.

\bigskip

\noindent \textit{Proof sketch of Lemma GRCM 2}. Introduce (Event $7^{\prime}$), from which, under similar assumptions provided for arguments in \textbf{Lemma} \textit{6.3}, yields a lower bound, obtained from the upper bound provided for,

\begin{align*}
  \bigg|   \mathcal{P}_1  +  \mathcal{P}_2  +  \mathcal{P}_3 \bigg|   \text{ } \text{ , } 
\end{align*}

\noindent in ($\Delta$). Proceeding, from following the previously implemented arguments relating to (Event $8^{\prime}$), (Event $9^{\prime}$), (Event $10^{\prime}$), (Event $11^{\prime}$), (Event $12^{\prime}$), (Event $13^{\prime}$), (Event $14^{\prime}$), (Event $15^{\prime}$), and (Event $16^{\prime}$), the suitably chosen, strictly positive, parameter described in the statement of \textbf{Lemma} \textit{GRCM 2} can be obtained, from which we conclude the argument.

\bigskip 

\noindent \textbf{Lemma} \textit{GRCM 3} (\textit{GRCM equivalent of \textbf{Lemma} \textit{6.2}}). From (Event $17^{\prime}$), the vertical crossing probability considered under the Ashkin-Teller measure,

\begin{align*}
 \mathcal{P}^{\chi}_{\mathcal{D}_{\mathrm{AT}}} \big[ \mathcal{V}^{\mathrm{Mixed}}_{\Lambda}  \big]    \text{ } \text{ , } 
\end{align*}

\noindent can be upper bounded, and lower bounded, with suitably chosen, strictly positive, constants similar to assumptions provided on the upper and lower bounds given in \textbf{Lemma} \textit{6.2}.

\bigskip

\noindent \textit{Proof sketch of Lemma GRCM 3}. Apply previous arguments, relating to (Event $17^{\prime}$), for obtaining a lower bound of similar nature to that provided in \textit{(AT 2.6.1)}. From (Event $18^{\prime}$), another lower bound can be obtained, from applications of (FKG) described in {\color{blue}[23]}. Incorporating the two aforementioned lower bounds, for (Event $17^{\prime}$), and for (Event $18^{\prime}$), readily yields a lower bound similar in nature to (\textit{AT 2.6.1 II}), from which we conclude the argument.

\bigskip

\noindent \textbf{Lemma} \textit{GRCM 4} (\textit{GRCM equivalent of \textbf{Lemma} \textit{AT 2.6.1.1}}). (Event $19^{\prime}$) can be lower bounded with a suitably chosen, strictly positive, constant.

\bigskip

\noindent \textit{Proof sketch of Lemma GRCM 4}. Lower bound (Event $19^{\prime}$) with the same properties of (FKG) mentioned in the proof sketch of \textit{Lemma GRCM 3} above. Following the lower bound for this \textit{generalized random-cluster} event, applying similar arguments in the remainder of the proof yields the desired lower bound, from which we conclude the argument.

\bigskip

\noindent \textbf{Lemma} \textit{GRCM 5} (\textit{GRCM equivalent of \textbf{Lemma} \textit{AT 2.6.1.2}}). (Event $20^{\prime}$) can be lower bounded with a suitably chosen, strictly positive constant.

\bigskip

\noindent \textit{Proof sketch of Lemma GRCM 5}. Lower bound (Event $20^{\prime}$) in a similar fashion that (Event $19^{\prime}$) is as described in the proof for the previous result, from which we conclude the argument.

\bigskip

\noindent \textbf{Lemma} \textit{GRCM 6} (\textit{GRCM equivalent of \textbf{Lemma} \textit{AT 2.6.1.3}}). (Event $21^{\prime}$), and (Event $22^{\prime}$), can be lower bounded with suitably chosen, strictly positive constants.

\bigskip

\noindent \textit{Proof sketch of Lemma GRCM 6}. Lower bound (Event $21^{\prime}$), and (Event $22^{\prime}$), in similar fashions that (Event $20^{\prime}$), and (Event $19^{\prime}$), in the two previous items, \textbf{Lemma} \textit{GRCM 4}, and \textbf{Lemma} \textit{GRCM 5}, from which we conclude the argument.

\bigskip

\noindent \textbf{Lemma} \textit{GRCM 7} (\textit{GRCM equivalent of Lemma AT 2.7}). (Event $27^{\prime}$), and (Event $28^{\prime}$), can be lower bounded with suitably chosen, strictly positive constants. 

\bigskip

\noindent \textit{Proof sketch of Lemma GRCM 7}. Obtain a lower bound for (Event $25^{\prime}$) as the lower bound that was obtained for (Event $25$). Following this lower bound, the conditional crossing event provided in (Event $26^{\prime}$) can be lower bound with a similarly related conditional crossing event that is instead dependent upon horizontal crossings across symmetric domains rather than the bridging event. Following this lower bound, (Event $27^{\prime}$) can be upper bounded with a strictly positive, suitably chosen constant. Finally, (Event $28^{\prime}$) can be upper bounded with another strictly positive, suitably chosen constant, from which we conclude the argument.

\bigskip

\noindent \textbf{Proposition} \textit{GRCM 2} (\textit{GRCM equivalent of \textbf{Proposition} AT 2}). A lower bound for (Event $31^{\prime}$) with a suitably chosen, strictly positive constant yields a lower bound for (Event $2^{\prime}$).

\bigskip

\noindent \textit{Proof sketch of Proposition GRCM 2}. Lower bound (Event $29^{\prime}$) with (FKG) properties described for the \textit{generalized random-cluster} measure in {\color{blue}[23]}. Following this lower bound, obtain another lower bound for (Event $30^{\prime}$), in addition to a lower bound for (Event $31^{\prime}$). Equipped with these two lower bounds, (Event $2^{\prime}$) can be bounded from below with a final suitably chosen, strictly positive constant, dependent upon the lower bounds for (Event $30^{\prime}$) and (Event $31^{\prime}$), from which we conclude the argument. 

\section{$(q_{\sigma} , q_{\tau})$-cubic models}

\subsection{Introduction}

\noindent In the final section, in light of arguments provided in the previous section, we briefly demonstrate how arguments provided for the \textit{generalized random-cluster} model can also be applied to the $(q_{\sigma}, q_{\tau} )$-cubic models. Below, we demonstrate how the probability measure for such models is defined from the Hamiltonian provided in {\color{blue}[23]}. In comparison to the Ashkin-Teller model Hamiltonian which consists of two two-point interactions and two four-point interactions, the Hamiltonian governing interactions for the $\big( q_{\sigma} , q_{\tau} \big)$-spin model consists of interaction terms, two of which are one-point functions, and three of which are two-point functions. From the formulation of Ashkin-Teller model, despite the fact that there is a difference in the number of two-point versus four-point interactions, by introducing mappings between the Ashkin-Teller and $\big( q_{\sigma} , q_{\tau}\big)$-spin Hamiltonians, one can imagine connections between interaction terms in a Hamiltonian and configurations can be sampled with high probability under the $\big( q_{\sigma} , q_{\tau} \big)$-spin measure.

\subsection{$\big( q_{\sigma} , q_{\tau} \big)$ objects}

\noindent \textbf{Definition} \textit{29} (\textit{cubic model Hamiltonian}, {\color{blue}[23]}). Introduce,

\begin{align*}
  \mathcal{H}^{\mathrm{cubic}} \equiv \mathcal{H} \big( i ,j , \sigma(i) , \sigma(j) , \tau(i) , \tau(j)  \big) \equiv \mathcal{H} \equiv  - 2 \big( J_{\sigma} - J_{\sigma\tau} \big) \delta_{\sigma_i,\sigma_j} - 2 \big( J_{\tau} - J_{\sigma\tau} \big) \delta_{\tau_i , \tau_j} - 4 J_{\sigma\tau} \delta_{\sigma_i , \sigma_j} \delta_{\tau_i , \tau_j}  \text{ } \text{ , } 
\end{align*}

\noindent corresponding to the \textit{cubic model Hamiltonian}, for $q_{\sigma} , q_{\tau} \in \textbf{N}$, $\sigma_i \in \big\{    1 , \cdots , q_{\sigma}     \big\}$, $\tau_i \in \big\{   1 , \cdots , q_{\tau}   \big\}$, and coupling constants which can either satisfy $J_{\tau} \equiv J_{\sigma\tau}$, $J_{\tau} \geq J_{\sigma\tau}$, or $J_{\tau} \leq J_{\sigma\tau}$.

\bigskip

\noindent Next, introduce the probability measure.

\bigskip

\noindent \textbf{Definition} \textit{30} (\textit{cubic model probability measure from the Hamiltonian}, {\color{blue}[23]}). Introduce,

\begin{align*}
 \mathscr{P}^{\mathrm{cubic}}\big[ \cdot \big]  \equiv \mathscr{P}^{\xi}_{\Lambda}  \big[ \cdot \big] \equiv \frac{\mathrm{exp} \big[ - \mathcal{H} \big]}{\mathcal{Z}^{\mathrm{cubic}} \big(\Lambda  ,    i , j , \sigma(i )  , \sigma(j) , \tau(i) , \tau(j)          \big) }   \equiv    \frac{\mathrm{exp} \big[ - \mathcal{H} \big]}{\mathcal{Z}^{\mathrm{cubic}} }           \text{ } \text{ , } 
\end{align*}

\noindent corresponding to the \textit{cubic model probability measure}, with the normalizing constant given by the \textit{cubic model} partition function $\mathcal{Z}$, for $\Lambda \subsetneq \textbf{Z}^d$ and boundary conditions $\xi$.

\bigskip

Following the definition of the Hamiltonian above and corresponding probability measure, arguments similar to those provided in the previous section also apply. That is, from the approximate correspondence provided in {\color{blue}[23]} between the measures $\mathcal{P}^{++} \big[ \cdot \big]$ and $\nu^{+} \big[ \cdot \big]$ that was applied in previous sections to reformulate weakened crossing probability estimates for the generalized random-cluster model from weakened crossing probability estimates for the Ashkin-Teller model, one can execute nearly identical steps of the argument to conclude that crossing probabilities occur with sufficiently good probability over strips of the square lattice.

\bigskip

With regards to boundary conditions, from previous remarks, encoding boundary conditions for the generalized random-cluster, and $\big(q_{\sigma} , q_{\tau} \big)$, models is similar to encoding boundary conditions for the Ashkin-Teller model, in the sense that versions of the SMP, CBC, and positive association inequalities are expected to hold. However, as mentioned in the first section giving an overview, the main obstacle for applying arguments for estimating long horizontal crossings over the strip is made far more difficult by the manner in which boundary conditions are encoded for the Ashkin-Teller model. It continues to remain of great interest to determine the ways in which RSW arguments can yield information about other models of Statistical Mechanics, beginning with information obtained from classes of symmetric domains as was done earlier for the six-vertex model in the planar random-cluster model {\color{blue}[15]}, which can likely be studied for many other models.

\section{Appendix}

\noindent In the final section of the paper, we establish a connection between the (FKG) criterionm, which was introduced in \textbf{Lemma} \textit{A.2} of {\color{blue}[11]} for the six-vertex model under flat boundary conditions, with the positive assocation on the \textit{marginals} of mixed-spin representations of the Ashkin-Teller model provided in \textbf{Theorem} \textit{4} of {\color{blue}[16]}. First, we state each of the two results below.

\bigskip

\noindent \textbf{Appendix Lemma} \textit{1} (\textit{FKG crtierion for the six-vertex model}, {\color{blue}[11]}). Suppose that $\mu$ is irreducible. If for every face $x$ in some finite volume $D$, every integer $k$, and $\mu$-almost every $\chi \in \mathcal{H}_{D \backslash \{x \}}$ - the set of all admissible height functions in the finite volume $D$ excluding the face $x$ - and $\mu$-almost every $\chi^{\prime} \in \mathcal{H}_{D \backslash \{x \}}$, with $\chi \leq \chi^{\prime}$,

\begin{align*}
    \mu \bigg[ h(x) \geq k \big|     h|_{D \backslash \{ x\} } \equiv \chi  \bigg]    \leq    \mu \bigg[ h(x) \geq k \big|      h|_{D \backslash \{ x\} } \equiv \chi^{\prime}      \bigg]        \text{ } \text{ , } 
\end{align*}

\noindent then for all increasing functions $F,G : \mathcal{H}_{D} \longrightarrow \textbf{R}$,

\begin{align*}
 \mu \big[ F(h) G(h) \big] \geq \mu \big[ F( h) \big] \mu \big[ G( h) \big]    \text{ } \text{ , }
\end{align*}

\noindent holds.

\bigskip

\noindent Next, also introduce the result for positive associativity of the mixed-spin representation below.

\bigskip

\noindent \textbf{Appendix Theorem} \textit{1} (\textit{positive association for the mixed-spin representation}, {\color{blue}[16]}). Let $D$ be a domain and $\tau \in \mathcal{E}_{\mathrm{spin}} ( \textbf{Z}^2)$ - the set of all possible $+$ and $-$ face variables assigned to spin configurations over $D$ - that is equal to $1$ at all odd faces outside of $D$. Suppose that $a,b,c > 0$ satisfy $a,b \leq c$. Then the marginal of the spin probability measure, $\mathrm{Spin}^{\tau}_{D ,a,b,c}$  on $\sigma^{\mathrm{even}}$ satisfies the FKG lattice condition. In particular, for any increasing functions $f,g : \{ -1 , +1 \}^{F^{\mathrm{even}}}$, one has,

\begin{align*}
      \mathrm{Spin}^{\tau}_{D,a,b,c} \big[ f(\sigma^{\mathrm{even}}) g( \sigma^{\mathrm{even}}) \big] \geq    \mathrm{Spin}^{\tau}_{D,a,b,c} \big[ f(\sigma^{\mathrm{even}})  \big]  \mathrm{Spin}^{\tau}_{D,a,b,c} \big[ g(\sigma^{\mathrm{even}})  \big]       \text{ } \text{ . } 
\end{align*}

\noindent To simultaneously incorporate the results together, from the (FKG) criterion with the positive association that the Ashkin-Teller measure satisfies on the marginals of spin configurations, we provide the following statement.

\bigskip

\noindent \textbf{Appendix Lemma} \textit{2} (\textbf{Appendix Lemma} \textit{1} + \textbf{Appendix Theorem} \textit{1}). The spin probability measure, $\mathrm{Spin}^{\tau}_{D,a,b,c}$, introduced above for the specified parameter range, satisfies a variant of the FKG criterion provided in \textbf{Appendix Lemma} \textit{1}, on the marginals of spin-representations.

\bigskip

\noindent \textit{Proof of Appendix Lemma 2}. The proof is a result of argumentation provided in \textbf{Appendix Theorem} \textit{1} of {\color{blue}[16]}. \boxed{}

\section{Acknowledgments}

\noindent The author would like to thank Philippe Sosoe for several comments on the previous submission of V1.

\newpage

\section{References}

\noindent [1] Baxter, R. Partition function of the eight-vertex lattice model, Ann. of Phys. \textbf{70} 193-228 (1972).

\bigskip

\noindent [2] Belov, P.A. The limit shape of the height function in the six-vertex model with domain-wall boundary conditions. IOP Publishing Ltd, Journal of Physics 1697, 012086 (2020).

\bigskip

\noindent [3] Bogoliubov, N. $\&$ Malyshev, C. The partition function of the four-vertex model in inhomogeneous external field and trace statistics. arXiv: 2011.10200v2 (2020).

\bigskip

\noindent [4] Burenev, I.N. $\&$ Pronko, A.G. Determinant formulas for the five-vertex model. arXiv: 2011.01972v2 (2020).

\bigskip

\noindent [5] Cai, J-Y., Fu, Z. $\&$ Xia, M. Complexity classification of the six-vertex model. arXiv: 1702.02863v1 (2017).

\bigskip

\noindent [6] Cai, J-Y., Liu, T. $\&$ Lu, P. Approximability of the Six-vertex model. arXiv: 1712.055880v1 (2017).

\bigskip

\noindent [7] Cai, J-Y., Liu, T., Lu, P. $\&$ Yu, J. Approximability of the Eight-vertex model. arXiv: 1811.03126v1 (2018).

\bigskip

\noindent [8] Corwin, I., Ghosal, P., Shen, H., $\&$ Tsai, L-C. Stochastic PDE limit of the six vertex model. arXiv: 1803.08210v3 (2019).

\bigskip

\noindent [9] Di Francesco, P. Twenty vertex model and domino tilings of the aztec triangle. arXiv: 2102.02920v1 (2021).

\bigskip

\noindent [10] Dimitrov, E. Six-vertex Models and the GUE-corners Process. International Mathematics Researhc Notices \textbf{6} 1794-1881 (2020).

\bigskip

\noindent [11] Duminil-Copin, H., Karrila, A., Manolescu, I., $\&$ Oulamara, M. Delocalization of the height function of the six-vertex model. arXiv: 2012.13750v1 (2020).

\bigskip

\noindent [12] Duminil-Copin, H., Krachun, D., Kozlowski, K., Manolescu, I., \& Tikhonovskaya, T. On the free energy of the six-vetex model. \textit{arXiv: 2012.11675}, 2020.

\bigskip

\noindent[13] Duminil-Copin, H., Hongler, C. \& Nolin, P. Connection probabilities and RSW-type bounds for the FK Ising Model. \textit{Communications on Pure and Applied Mathematics} \textbf{64}(9) (2011).

\bigskip

\noindent [14] Duminil-Copin, H., Sidoravicius, V. \& Tassion, V. Continuity of the phase transition for planar random-cluster and Potts models for $1 \leq q \leq 4$. \textit{Communications in Mathematical Physics} \textbf{349} 47-107 (2017).

\bigskip

\noindent[15] Duminil-Copin, H. \& Tassion, V. Renormalization of crossing probabilities in the planar random-cluster model. \textit{Moscow Mathematical Journal} \textbf{20}(4):711-740 (2020). 

\bigskip

\noindent [16] Glazman, A., Peled, R. On the transition between the disordered and antiferroelectric phases of the 6-vertex model. arXiv: 1909.03436 (2019).

\bigskip 
\noindent [17] Gorbunov, V.G., Korff, C. $\&$ Stroppel, C. Yang-Baxter algebras, convolutional algebras, and Grassmannians. Russian Math Surveys \textbf{75}:5 791-842 (2020).

\bigskip

\noindent [18] Kapitonov, V.S. $\&$ Pronko, A.G. Six-vertex model as Grasmann integral, one-point function, and Arctic ellipse. POMI \textbf{494} 168-218 (2020).

\bigskip

\noindent [19] Kenyon, R., Miller, J., Sheffield, S. $\&$ Wilson, D.B. The six-vertex model and Schramm-Loewner evolution. arXiv: 1605.0647v2 (2017).

 \bigskip

 \noindent [20] Lieb, E. H. The Residual Entropy of Square Ice, Phys Rev. \textbf{162} 162-172 (1967).

\bigskip

\noindent [21] Liu, T. Torpid mixing of Markov chains for the six-vertex model on $\textbf{Z}^2$. Approximation, Randomization, and Combinatorial Optimization \textbf{52} 1-15 (2018).

\bigskip

\noindent [22] Pauling, L. J. Am. Chem. Soc. \textbf{57}, 2680 (1935).

\bigskip

\noindent [23] Pfister, C.-E., Velenik, Y. Random-Cluster Representation of the Ashkin-Teller model. \textit{arXiv, 9704017 v1} (1997).

\bigskip

\noindent [24] Ray, G. $\&$ Spinka, Y. A short proof of the discontinuity of the phase transition of the random-cluster mode with $q>4$. arXiv: 1904.105571v1 (2019).

\bigskip

\noindent [25] \noindent Rigas, P. Renormalization of crossing probabilities in the dilute Potts model. \textit{arXiv preprint}.

\bigskip

\noindent [26] Reshetikhin, N. Lectures in the integrability of the 6-vertex model. arXiv: 1010.5031v1 (2010).

\bigskip

\noindent [27] Russo, L. A note on percolation. \textit{Zeitschrift fur Wahrscheinlichkeitstheorie und Verwandte Gebiete} \textbf{43} 39-48 (1978).

\bigskip

\noindent [28] Seymour, P. \& Welsh, D. Percolation probabilities on the square lattice. \textit{Annals
Discrete Math}, \textbf{3} 227–245 (1978).

\bigskip

\noindent [29] Sridhar, A. Limit Shapes in the Six Vetex Model. PhD Thesis.

\bigskip

\noindent [30] Vijay, S. Studying the six-vertex model with the Yang-Baxter equation. \textit{Notes}.

\bigskip

\noindent [31] Zinn-Justin, P. Six-vertex, loop and tiling models: integrability and combinatorics. arXiv: 0901.0665 (2009).

\bigskip

\noindent [32] Zhong, C. Stochastic symplectic ice. arXiv: 2102.00660v1 (2021).

\end{document}